\documentstyle[12pt]{article}
\newcommand{\numero}[1]{
\newpage
\addtocounter{section}{1}
\begin{center}{\Large \bf \thesection .\
#1\vspace{-.1in}}\end{center}
\setcounter{subsection}{0}
\setcounter{lemma}{0}\indent}

\newcommand{\subnumero}[1]{
\pagebreak[3]\begin{center}{\bf #1}\nopagebreak\end{center}
}

\newcommand{\eop}{\hfill $\Box$\vspace{.1in}}

\newtheorem{lemma}{Lemma}[section]

\newtheorem{conjecture}[lemma]{Conjecture}
\newtheorem{proposition}[lemma]{Proposition}

\newtheorem{lemme}[lemma]{Lemme}
\newtheorem{theoreme}[lemma]{Th\'eor\`eme}
\newtheorem{definition}[lemma]{D\'efinition}
\newtheorem{corollaire}[lemma]{Corollaire}
\newtheorem{probleme}[lemma]{Probl\`eme}
\newtheorem{convention}[lemma]{Convention}

\newtheorem{lemmedefinition}[lemma]{Lemme-D\'efinition}
\newtheorem{remarque}[lemma]{Remarque}

\newcommand{\rr}{{\bf R}}
\newcommand{\qq}{{\bf Q}}

\newcommand{\zz}{{\bf Z}}

\newcommand{\mm}{{\bf M}}
\newcommand{\nn}{{\bf N}}

\newcommand{\Cc}{{\cal C}}

\newcommand{\Ff}{{\cal F}}
\newcommand{\Gg}{{\cal G}}

\newcommand{\Oo}{{\cal O}}

\newcommand{\Bb}{{\cal B}}

\newcommand{\Uu}{{\cal U}}

\newcommand{\Xx}{{\cal X}}
\newcommand{\Yy}{{\cal Y}}
\newcommand{\Zz}{{\cal Z}}

\newcommand{\facon}{fa\c{c}on\,\, }

\newcommand{\tworightarrows}{\stackrel{\displaystyle \rightarrow}{\rightarrow}}

\newcommand{\direct}[1]{\stackrel{\rightarrow}{#1}}
\newcommand{\inverse}[1]{\stackrel{\leftarrow}{#1}}

\begin{document}

\section*{Descente pour les $n$-champs}

Andr\'e Hirschowitz
\footnote{
Universit\'e de Nice-Sophia Antipolis, Parc Valrose,
06108 Nice cedex 2, France}
et
Carlos Simpson
\footnote{
CNRS, Laboratoire Emile Picard, Universit\'e Toulouse 3, 31062 Toulouse cedex,
France}

\begin{center}
{\bf Abstract}
\end{center}

We develop the theory of $n$-stacks (or more generally Segal $n$-stacks
which are $\infty$-stacks such that the morphisms are invertible above
degree $n$). This is done by systematically using the theory of
closed  model categories (cmc). Our main results are: a definition of
$n$-stacks in terms of limits,  which should be perfectly general for stacks
of any type of objects;  several other characterizations of $n$-stacks in
terms of ``effectivity of descent data''; construction of the stack
associated to an $n$-prestack; a strictification result saying that
any ``weak'' $n$-stack is equivalent to a (strict) $n$-stack; and a
descent result saying that the $(n+1)$-prestack of $n$-stacks (on a
site) is an $(n+1)$-stack. As for other examples, we start from a
``left Quillen presheaf'' of cmc's and introduce the associated
Segal $1$-prestack. For this situation, we prove a general descent result,
giving sufficient conditions for this prestack to be a stack.
This applies to the case of complexes, saying how complexes of
sheaves of $\Oo$-modules can be glued together via quasi-isomorphisms.
This was the problem that originally motivated us.

\begin{center}
{\bf R\'esum\'e}
\end{center}

On d\'eveloppe une th\'eorie des $n$-champs (plus exactement celle des
$n$-champs de Segal, qui sont des $\infty$-champs o\`u les morphismes
sont inversibles en degr\'e $\geq n$). Pour cela on utilise
syst\'ematiquement la th\'eorie des cat\'egories de mod\`eles
ferm\'ees (cmf). Nos contributions principales sont: une d\'efinition
de $n$-champ en termes de limites, qui est parfaitement
g\'en\'eralisable \`a toutes sortes d'autres champs; plusieurs autres
caract\'erisations des $n$-champs en termes d'``effectivit\'e'' des
donn\'ees de descente; la construction du champ associ\'e \`a un
$n$-pr\'echamp; un r\'esultat de strictification assurant que tout
$n$-champ ``faible'' est \'equivalent \`a un
$n$-champ (strict); et un r\'esultat de descente affirmant que le
$(n+1)$-pr\'echamp des $n$-champs (sur un site) est un $(n+1)$-champ.
Pour d'autres exemples, nous partons d'un pr\'efaisceau de cmf ``de
Quillen \`a gauche'' et introduisons le $1$-pr\'echamp de Segal associ\'e.
Dans ce cadre, nous prouvons un r\'esultat de descente g\'en\'eral
donnant des conditions suffisantes pour que ce pr\'echamp soit un
champ. Ceci s'applique au cas des complexes, et dit comment on
peut recoller des complexes de $\Oo$-modules \`a l'aide de
quasi-isomorphismes. C'est ce probl\`eme qui \'etait la motivation
initiale du pr\'esent travail.

\newpage

\begin{center}
{\Large \underline{Sommaire}}
\end{center}

\noindent
{\bf 1. Introduction}---p. \pageref{intropage}
\newline
{\small
Motivation, description des r\'esultats.
}

\noindent
{\bf 2. Les $n$-cat\'egories de Segal}---p. \pageref{ncatsegpage}
\newline
{\small
D\'efinition des $n$-cat\'egories de Segal,
troncation, int\'erieur. La m\'ethode
g\'en\'erale pour obtenir une structure de cmf, et son
application aux $n$-cat\'egories de Segal.
Les classes d'homotopie de morphismes.}

\noindent
{\bf 3. Les $n$-pr\'echamps de Segal}---p. \pageref{nchampsegpage}
\newline
{\small
La structure de cmf pour les $n$-pr\'echamps de Segal (notions d'\'equivalence
faible, de
cofibration, de fibration), comparaison avec la topologie grossi\`ere,
compatibilit\'e avec la troncation, propret\'e. Donn\'ees de descente.
}

\noindent
{\bf 4. Les fonctorialit\'es $p^{\ast}$, $p_{\ast}$, $p_!$}---p.
\pageref{fonctorialitepage}
\newline
{\small
Relation entre $n$-pr\'echamps de Segal sur $\Xx$ et $\Zz$ pour un foncteur
$p:\Zz
\rightarrow \Xx$, crit\`eres de pr\'eservation des objets fibrants.
}

\noindent
{\bf 5. La structure de type Bousfield-Kan d'apr\`es Hirschhorn}---p.
\pageref{bousfieldkanpage}
\newline
{\small
La structure HBKQ de cmf pour les $n$-pr\'echamps de Segal, o\`u les
cofibrations
sont obtenues par libre addition de cellules; application \`a la
notion de donn\'ee de descente. }

\noindent
{\bf 6. Le point de vue de la localisation}---p.
\pageref{localisationpage}
\newline
{\small
Les cmf de $n$-pr\'echamps de Segal pour la topologie $\Gg$ peuvent \^etre
obtenues
\`a partir des cmf pour la topologie grossi\`ere par inversion de fl\`eches
correspondant aux cribles de $\Gg$.
}

\noindent
{\bf 7. Cat\'egories de Segal et cat\'egories simpliciales}---p.
\pageref{segsimplpage}
\newline
{\small
Relation entre les $1$-cat\'egories de Segal et les cat\'egories simpliciales;
produits et coproduits homotopiques et l'argument de Gabriel-Zisman.
}

\noindent
{\bf 8. Localisation de Dwyer-Kan}---p.
\pageref{dwyerkanpage}
\newline
{\small
La localisation de Dwyer-Kan fournit beaucoup d'exemples de cat\'egories
simpliciales, et permet en outre de montrer que toute $1$-cat\'egorie de Segal
est \'equivalente \`a une cat\'egorie simpliciale.
}

\noindent
{\bf 9. Premi\`ere d\'efinition de champ}---p.
\pageref{notionchamppage}
\newline
{\small
D\'efinition de $n$-champ de Segal comme $n$-pr\'echamp de Segal $A$ dont les
$A(X)$ sont des $n$-cat\'egories de Segal et pour lequel le morphisme
$A\rightarrow A'$ vers le remplacement $\Gg$-fibrant est une \'equivalence
objet-par-objet. Notion de champ associ\'e, compatibilit\'es.
}

\noindent
{\bf 10. Crit\`eres pour qu'un pr\'echamp soit un champ}---p.
\pageref{caracterisationpage}
\newline
{\small
$A$ est un $n$-champ de Segal si et seulement si les $A_{1/}(x,y)$ sont des
$n-1$-champs de Segal sur $\Xx /X$ pour tout $x,y\in A_0(X)$, et si les
donn\'ees de descente pour un crible $\Bb \subset \Xx /X$ sont effectives.
On a plusieurs versions de ce crit\`ere.
}

\noindent
{\bf 11. Cat\'egories de mod\`eles internes}---p.
\pageref{cmipage}
\newline
{\small
On definit la notion de cat\'egorie de mod\`eles interne et on l'utilise pour
d\'efinir le $n+1$-pr\'echamp $nSe\underline{CHAMP}(\Xx )$ des $n$-champs de
Segal au-dessus de $\Xx$, ainsi que sa $n+1$-cat\'egorie de Segal
$nSeCHAMP(\Xx )$ de sections globales.
Comparaison avec la localis\'ee de Dwyer-Kan.}

\noindent
{\bf 12. La famille universelle}---p.
\pageref{universellepage}
\newline
{\small
On d\'efinit un morphisme $nSeCHAMP(\Xx )\rightarrow \underline{Hom}(\Xx ^o,
nSeCAT')$, on montre qu'il est pleinement fid\`ele (et on \'enonce le
r\'esultat \ref{correlation}
qui identifie son image essentielle). Ceci nous conduit \`a une nouvelle version
de la notion de champ.}

\noindent
{\bf 13. Le champ associ\'e \`a un pr\'echamp}---p.
\pageref{associepage}
\newline
{\small
On d\'efinit le foncteur
``champ associ\'e \`a un pr\'echamp'' de $nSeCHAMP(\Xx ^{\rm gro})$ vers
\newline
$nSeCHAMP (\Xx ^{\Gg})$, qui est adjoint au foncteur d'inclusion dans l'autre
sens.
}

\noindent
{\bf 14. Limites}---p.
\pageref{limitespage}
\newline
{\small
Calcul des limites \`a l'aide des structures de cmf;
d\'efinition des champs en termes de limites. }

\noindent
{\bf 15. Un peu plus sur la condition de descente}---p.
\pageref{peupluspage}
\newline
{\small
On refait les crit\`eres de \S 10 en termes de limites, et on explicite le
calcul de ces limites dans le cas d'un crible d\'efini par une famille
couvrante, en termes de limites prises au-dessus de $\Delta$. }

\noindent
{\bf 16. La construction de Grothendieck}---p.
\pageref{gropage}
\newline
{\small
On construit la ``$n$-cat\'egorie fibr\'ee'' associ\'ee \`a un $n$-pr\'echamp
au-dessus d'une cat\'egorie, o\-p\'e\-ra\-tion not\'ee comme une
int\'egrale suivant
Thomason. On compare les sections du $n$-pr\'echamp et celles de la
$n$-cat\'egorie
fibr\'ee associ\'ee.
}

\noindent
{\bf 17. Pr\'efaisceaux de Quillen}---p.
\pageref{quipage}
\newline
{\small
On introduit la notion de pr\'efaisceau de Quillen \`a gauche, une sorte de
pr\'efaisceau de cmf
o\`u les foncteurs de restriction sont des foncteurs de Quillen \`a gauche. On
donne diverses structures de cmf pour les sections de l'int\'egrale d'un
pr\'efaisceau de Quillen. }

\noindent
{\bf 18. Strictification}---p.
\pageref{strictpage}
\newline
{\small
On strictifie les sections faibles des $1$-pr\'echamps de Segal qui
proviennent de pr\'efaisceaux de Quillen par localisation de Dwyer-Kan. Ceci
permet de strictifier des familles faibles de $n$-cat\'egories de Segal, en
$n$-pr\'echamps de Segal. Fin de la preuve de \ref{correlation}.}

\noindent
{\bf 19. La descente pour les pr\'efaisceaux de Quillen \`a gauche}---p.
\pageref{quidescentepage}
\newline
{\small
Pour un pr\'efaisceau de Quillen \`a gauche $\mm$ on donne un crit\`ere pour que
le $1$-pr\'echamp de Segal associ\'e $L(\mm )$ soit un champ. }

\noindent
{\bf 20. Exemple: la descente pour les
$n$-champs}---p.
\pageref{exemplepage}
\newline
{\small
On utilise le crit\`ere du \S 19 pour prouver que $nSe\underline{CHAMP}(\Xx )$
est un $n+1$-champ de Segal. Ceci g\'en\'eralise le th\'eor\`eme classique de
recollement des faisceaux. On donne aussi une preuve directe. }

\noindent
{\bf 21. Exemple: la descente pour les complexes}---p.
\pageref{complexepage}
\newline
{\small
On utilise le crit\`ere du \S 19 pour prouver que $L(Cpx _{\Oo},qis)$ est un
$1$-champ de Segal, pour un site annel\'e raisonnable $(\Xx , \Oo )$. On
\'enonce un r\'esultat de g\'eometricit\'e pour le champ de modules des
complexes
parfaits. }

\numero{Introduction}

\label{intropage}

\subnumero{La descente: des modules aux complexes}
Pour donner une id\'ee du contenu du pr\'esent travail,
on peut l'aborder, en premi\`ere approximation,
comme une contribution \`a la th\'eorie des complexes de faisceaux de
modules. Pour le mettre en perspective, nous commen\c{c}ons par
r\'esumer
la th\'eorie des faisceaux de modules localement libres
(de rang fix\'e $r$):

Il existe un champ $BG$,
muni d'un module universel $U$ et d'un fibr\'e universel dont $U$
est le module des sections. Ce champ ``classifie'' les
modules localement libres de rang $r$ sur les sch\'emas
en ce sens que le champ de modules
$Fib_{X,r}$ de tels modules sur un sch\'ema
$X$ s'identifie au champ des morphismes de $X$ vers $BG$.
Le champ $BG$ est alg\'ebrique et si $X$ est projectif, il en est de m\^eme
pour $Fib_{X,r}$. Le champ $BG$ est un ouvert dans
(au moins) deux champs plus g\'en\'eraux $Coh$ et
$Qcoh$ qui classifient respectivement les faisceaux
coh\'erents et quasi-coh\'erents (on n'entre pas ici dans les d\'etails
concernant en particulier la topologie choisie).
Le fait que $Qcoh$ soit un champ signifie plus concr\'etement que les
donn\'ees de descente pour un module quasi-coh\'erent (ou pour un
morphisme entre deux modules quasi-coh\'erents) sont effectives.
On peut observer que $Coh$ n'est pas alg\'ebrique, ce
qui est un peu contrariant.

On se propose de g\'en\'eraliser cette th\'eorie au cas des complexes.
Le point crucial est qu'on veut parler de recollement de complexes
\`a l'aide non pas d'isomorphismes mais de {\em quasi-isomorphismes}
(compatibles en un sens ad\'equat), de sorte que ce travail rel\`eve
dans une large mesure de la th\'eorie de l'homotopie.

Les champs dont on vient de parler sont des champs de cat\'egories
et lorsqu'on se pr\'eoccupe d'\'etendre ce qui pr\'ec\`ede au cas des
complexes, il faut t\^ot ou tard organiser ces complexes en cat\'egories
munies de structures ad\'equates.

La solution qui vient d'abord \`a l'esprit consiste \`a consid\'erer
des {\em cat\'egories d\'eriv\'ees}. On peut en effet former par exemple le
pr\'echamp $D_{Qcoh}$ qui \`a un sch\'ema affine $Spec A$ associe
la cat\'egorie d\'eriv\'ee $D_{Qcoh}(A)$ de celle des $A$-modules
quasi-coh\'erents. Demander si ce pr\'echamp est un champ
est la fa\c{c}on savante
de demander si on peut recoller des complexes \`a l'aide de quasi-isomorphismes
``compatibles'' (i.e. v\'erifiant une condition de cocycle par ailleurs
assez technique). Cette question a \'et\'e tr\`es t\^ot reconnue comme
impertinente, les objets de ces cat\'egories d\'eriv\'ees \'etant
``de nature essentiellement non-recollables'' \footnote{
Pour le lecteur qui ne serait pas convaincu que cette voie des
cat\'egories d\'eriv\'ees est sans issue, signalons que
le sous-pr\'echamp (plein) $Parf^{\{0,1\}}$ de $D_{Qcoh}$
des complexes parfaits \`a support cohomologique dans $[0,1]$ est bien
un champ, mais qu'il n'est pas localement alg\'ebrique (au sens d'Artin):
en fait, comme on le verra plus loin, le ``bon'' objet est un $2$-champ
localement alg\'ebrique (au sens de \cite{geometricN}).}
(\cite {SGA6} Expos\'e 0, p.11).

\subnumero{Les homotopies sup\'erieures}

Pour formuler les probl\`emes de recollement des complexes
d\'efinis \`a quasi-isomorphisme pr\`es,
il faut consid\'erer des homotopies sup\'erieures, ce qui
complique singuli\`erement le tableau. En topologie ordinaire,
le prototype d'un tel probl\`eme consiste par exemple \`a recoller sur
$U_1 \cup U_2 \cup U_3 \cup U_4$ des donn\'ees
du genre suivant: pour $1 \leq i\leq 4$, $C_i$ est un complexe sur $U_i$;
pour $1 \leq i\leq j \leq 4$, $f_{ij}$ est une \'equivalence d'homotopie
entre les restrictions $C_{ij}$ et $C_{ji}$
de $C_i$ et $C_j$ \`a $U_{ij}:=U_i \cap U_j$;
pour $1 \leq i\leq j \leq k\leq 4$, $h_{ijk}$ est une \'equivalence d'homotopie
entre les restrictions de $f_{jk}o f_{ij}$ et $f_{ik}$ \`a
$U_{ijk}:=U_i \cap U_j \cap U_k$; et les $h_{ijk}$ doivent encore v\'erifier
une condition de compatibilit\'e sur $U_1 \cap U_2 \cap U_3 \cap U_4$,
condition dont la formulation m\^eme n'est pas imm\'ediate. On con\c{c}oit
facilement
comment, pour des recouvrements plus complexes (en topologie \'etale par
exemple),
la combinatoire de ce genre de donn\'ees de descente peut devenir inextricable.
La d\'efinition que nous donnons des donn\'ees de descente d\'epasse (ou
\'evite)
les consid\'erations combinatoires.

Dans la situation pr\'ec\'edente les
$f_{ij}$ sont des fl\^eches entre complexes, les
$h_{ijk}$ doivent
\^etre consid\'er\'ees comme des $2$-fl\`eches, et les donn\'ees sur
$U_1 \cap U_2 \cap U_3 \cap U_4$ seront des $3$-fl\`eches.
Concr\'etement, $h_{ijk}$ par exemple est un morphisme d'objet
gradu\'e de degr\'e $-1$ entre $C_i$ et $C_k$ sur $U_{ijk}$, \'etablissant une
homotopie entre $f_{jk}o f_{ij}$ et $f_{ik}$. De m\^eme,
l'information cruciale dans une $3$-fl\^eche est un morphisme d'objet
gradu\'e de degr\'e $-2$ \'etablissant une homotopie entre morphismes
de degr\'e $-1$, etc.
 On a donc bien
besoin d'une notion d'{\em $\infty$-cat\'egorie} comme pr\'econis\'ee par
Grothendieck \cite{Grothendieck}, avec des
$n$-fl\`eches pour tout $n$ et les lois de composition ad\'equates.

\subnumero{Cat\'egories simpliciales ou de Segal}

En fait, une forme rudimentaire d'$\infty$-cat\'egorie adapt\'ee
\`a nos complexes
est connue depuis trente ans, c'est la notion de
cat\'egorie simpliciale. Introduite en th\'eorie d'homotopie par Kan
et Quillen
(voir \cite {Quillen}), cette notion a \'et\'e
reprise par Dwyer et Kan \cite {DwyerKan1}
\cite {DwyerKan2} \cite {DwyerKan3}: en particulier,
\`a toute cat\'egorie $M$ munie d'une sous-cat\'egorie $W$, ces auteurs
associent
une cat\'egorie
simpliciale ``localis\'ee'' $L(M,W)$. Celle-ci capture bien l'information
homotopique concernant le couple $(M,W)$ dans la mesure o\`u si
$M$ est une cat\'egorie de mod\`eles ferm\'ee
(cmf) simpliciale au sens de Quillen,
et $W$ sa sous-cat\'egorie des \'equivalences, alors $L(M,W)$ est \'equivalente
\`a la cat\'egorie simpliciale des objets cofibrants et fibrants de $M$. Ceci
montre par exemple que la structure simpliciale sur une cmf est unique \`a
\'equivalence pr\`es, et m\^eme ne d\'epend (toujours
\`a \'equivalence pr\`es) que de la sous-cat\'egorie des \'equivalences.

Plus particuli\`erement, dans la situation des complexes si $Ch$ d\'esigne
la cat\'egorie des complexes et $qis$ la sous-cat\'egorie des \'equivalences
faibles, alors la cat\'egorie simpliciale $L(Ch, qis)$ repr\'esente
ad\'equatement la th\'eorie homotopiquement des complexes.
\footnote{
Ceci veut dire que si $A^{\cdot}$ et $B^{\cdot}$ sont des complexes alors
l'ensemble simplicial de morphismes dans $L$ entre $A^{\cdot}$ et $B^{\cdot}$
est \'equivalente \`a l'ensemble obtenu par application de la construction de
Dold-Puppe au complexe tronqu\'e $\tau ^{\leq 0}{\rm Hom}^{\cdot}(A^{\cdot},
B^{\cdot})$.
On peut noter qu'une cat\'egorie simpliciale $L$ donne lieu
(suivant l'id\'ee de Grothendieck
\cite{Grothendieck}) \`a une $\infty$-cat\'egorie $\Pi _{\infty}\circ L$
par application de la construction ``$\infty$-groupoide de Poincar\'e''
$\Pi _{\infty}$ \`a chacune des ensembles simpliciales $Hom _L(x,y)$.
Dans cette $\infty$-cat\'egorie pour les complexes, l'$\infty$-cat\'egorie des
morphismes entre deux complexes est \'equivalente \`a
$\gamma {\rm Hom}^{\cdot}(A^{\cdot},
B^{\cdot})$, o\`u $\gamma$ est la construction
bien connue
qui \`a un complexe de groupes ab\'eliens $F^{\cdot}$ associe
l'$\infty$-cat\'egorie stricte $\gamma (F^{\cdot})$ dont les objets sont les
$x\in F^0$ avec
$d(x)=0$, les morphismes entre $x,y$ sont les $f\in F^{-1}$ avec $d(f)=y-x$,
les $2$-morphismes entre $f$ et $g$ sont les $u\in F^{-2}$ avec $d(u)=g-f$
ainsi de suite.}

Pour des raisons techniques, nous
utilisons syst\'ematiquement la notion un peu plus g\'en\'erale
de cat\'egorie de Segal: une cat\'egorie de Segal est une cat\'egorie
simpliciale ``large'', en ce sens qu'on n'impose pas que la composition
soit strictement associative. Toute cat\'egorie simpliciale est une
cat\'egorie de Segal et, inversement, on montre que toute cat\'egorie de
Segal est \'equivalente \`a une cat\'egorie simpliciale (\S 7).

\subnumero{Les probl\`emes pos\'es}

Ainsi les complexes de faisceaux de modules s'organisent (par exemple
sur le site \'etale)
en pr\'efaisceaux de cat\'egories de Segal. Entre autres, on peut d\'efinir
${\cal D}^+$ par ${\cal D}^+(Spec A)= L(Ch^+(A), qis)$, o\`u
$Ch^+(A) $ est la cat\'egorie des complexes de $A$-modules born\'es \`a gauche
et $qis$ sa sous-cat\'egorie des quasi-isomorphismes.
On peut
maintenant
exprimer de la fa\c{c}on suivante les probl\`emes auxquels le pr\'esent travail
est consacr\'e:
\begin{itemize}
\item
identifier les probl\`emes de descente pertinents concernant de tels
pr\'efaisceaux;
\item
formuler ces probl\`emes dans le cadre d'une th\'eorie des
champs (g\'en\'eralis\'es; dans cette introduction on dira $\infty$-champs)
adapt\'ee aux cat\'egories de Segal; dans cette th\'eorie, il faudra
bien entendu que les champs soient les pr\'echamps dans lesquels les donn\'ees
de descente se recollent;
\item
identifier une classe de tels pr\'efaisceaux, contenant nos pr\'efaisceaux
de complexes, et dans lesquels ces probl\`emes de descente admettent
une solution.
\end {itemize}

Au nombre des probl\`emes de descente qu'on veut traiter doivent figurer
les probl\`emes ``concrets'' de recollement
concernant les complexes, dont nous donnons maintenant des exemples.
Signalons toutefois que nous n'abordons pas directement ces
probl\`emes dans le pr\'esent travail---notre but \'etant dans un premier
temps de mettre en place le
cadre dans lequel de telles questions doivent naturellement
\^etre consid\'er\'ees---et nous ne les formulons
donc qu'en guise de
motivation.

\begin{itemize}
\item
Soit $k$ un corps, $X$ un $k$-sch\'ema r\'egulier, et $\Ff$ un faisceau
coh\'erent
sur $X$. On sait que $\Ff$ admet localement des r\'esolutions projectives
finies. Peut-on recoller ces r\'esolutions en un objet global?
\item
Dans \cite {SGA6} 0.4.4, pour un morphisme $f:X \rightarrow Y$ de sch\'emas,
on construit localement dans $X$ un complexe cotangent bien d\'efini \`a
quasi-isomorphisme pr\`es. Peut-on recoller ces constructions locales?
\footnote { Dans loc. cit. on consid\`ere qu'
``il n'est pas possible de proc\'eder par simple recollement'' et on
s'en tire autrement.}

\item

Toujours dans \cite {SGA6} (e.g. 0.4.2), on pose le probl\`eme
(ult\'erieurement r\'esolu par Deligne et Illusie, voir \cite {Illusie} I.4.2)
de la d\'efinition des puissances ext\'erieures d'un complexe parfait,
que nous reformulons \`a notre convenance: sachant
qu'on peut d\'efinir, sur les complexes born\'es de
modules  libres de type fini, un foncteur puissance sym\'etrique $\Lambda^i$ qui
transforme quasi-isomorphismes en quasi-isomorphismes, et sachant
que les complexes parfaits sont ceux qui sont localement quasi-isomorphes
\`a des complexes born\'es de faisceaux libres de type fini, comment
d\'efinir la puissance sym\'etrique d'un complexe parfait?

\item

Dans les travaux de O'Brian, Toledo et Tong
\cite{ObrianToledoTong1}, \cite{ObrianToledoTong2},
\cite{Toledo1}, \cite{Toledo2}, \cite{Toledo3}
consacr\'es \`a une autre question issue de SGA 6,
celle des {\em formules de Riemann-Roch}, on trouve des calculs de \v{C}ech
qui  sont certainement un exemple de situation de descente pour les
complexes. Un meilleur cadre g\'en\'eral pour ces calculs
pourrait
contribuer \`a notre compr\'ehension des formules de Riemann-Roch.

\item

Il y a certainement d'autres applications potentielles que celles
qui concernent les complexes, voir par exemple
le travail de Hinich \cite{Hinich1} qui traite un probl\`eme de descente
pour des $1$-groupo\"{\i}des en relation avec un probl\`eme de descente
correspondant pour des pr\'efaisceaux d'alg\`ebres de Lie
diff\'erentielles gradu\'es.

\end{itemize}

\subnumero {Les $\infty$-champs}

Avant de pr\'esenter notre vision des donn\'ees de descente, il nous faut
expliquer notre th\'eorie des $\infty$-champs.
On a compris plus haut qu'il nous faut---au minimum---une
th\'eorie des champs de cat\'egories simpliciales (ou de Segal), comme on a
une th\'eorie des champs d'ensembles (les faisceaux) et une th\'eorie des
champs de cat\'egories (les champs ``classiques''). Et dans une th\'eorie des
champs de cat\'egories simpliciales, le champ des morphismes entre
deux objets (i.e. sections sur un ouvert) devra \^etre un champ d'ensembles
simpliciaux. On voudrait donc bien d'une th\'eorie des $\Gamma$-champs,
o\`u $\Gamma$ pourrait \^etre
indiff\'eremment
la cat\'egorie des ensembles ou celle des
cat\'egories, ou celle des cat\'egories simpliciales, ou encore celle des
ensembles simpliciaux. Un $\Gamma$-pr\'echamp sur un site $\Xx$ est alors
un pr\'efaisceau sur $\Xx$ \`a valeurs dans $\Gamma$ (on prend des
pr\'efaisceaux stricts pour simplifier, et parce qu'on esp\`ere qu'au bout
du compte
\c{c}a revient au m\^eme). Et l'id\'ee (na\"{\i}ve) pour
exprimer qu'un $\Gamma$-pr\'echamp $F$ est un $\Gamma$-champ consisterait
\`a demander que pour tout $X \in \Xx$ et tout crible $\Bb$ couvrant $X$,
$F(X)$ soit la limite (dans $\Gamma$) de la restriction de $F$ \`a la
cat\'egorie
$\Bb$. Si $\Gamma$ est la cat\'egorie des ensembles on retrouve
\'evidemment la notion de faisceau, mais si $\Gamma$ est la cat\'egorie des
cat\'egories, \c{c}a ne va plus. Pour rectifier le tir et
retomber dans ce cas sur les champs usuels, il suffit
de  munir $\Gamma$ de sa structure de $2$-cat\'egorie et de demander
que $F(X)$ soit la $2$-limite (dans $\Gamma$) de la
restriction de $F$ \`a $\Bb$. Pour g\'en\'eraliser cette situation, on va donc
consid\'erer que $\Gamma$ est une $\infty$-cat\'egorie
\footnote{On ne prend pas la peine ici de donner une d\'efinition
d'$\infty$-cat\'egorie: on consid\`ere simplement que les $n$-cat\'egories
de Tamsamani \cite {Tamsamani} ainsi que leurs variantes introduites ci-dessous,
les $n$-cat\'egories de Segal, sont des $\infty$-cat\'egories.}
dans laquelle on sait sp\'ecifier les limites homotopiques
(il n'est pas n\'ecessaire qu'elles existent). On peut alors, suivant une
suggestion
formul\'ee par le second auteur dans \cite {limits} 6.4, dire
qu'un $\Gamma$-champ est un $\Gamma$-pr\'echamp $F$ tel
que pour tout $X \in \Xx$ et tout crible $\Bb$ couvrant $X$,
$F(X)$ soit la limite (dans $\Gamma$) de la restriction de $F$ \`a la
cat\'egorie
$\Bb$.

\subnumero{Les $0$-champs de Segal}

Dans le cas o\`u $\Gamma$ est l'$\infty$-cat\'egorie
(techniquement: la $1$-cat\'egorie de Segal) des ensembles simpliciaux, nous
confrontons avec succ\`es
\footnote{
Quelque peu temp\'er\'e dans la pr\'esente version v3\ldots}
ce point de vue avec la th\'eorie existante
des pr\'efaisceaux simpliciaux,  d\'evelopp\'ee par K. Brown \cite{KBrown},
Illusie
\cite{Illusie},  Joyal dans une lettre \`a Grothendieck \cite{Joyal}, Jardine
\cite{Jardine}, et Thomason \cite{Thomason}.
\footnote{On a m\^eme la
$\infty$-cat\'egorie (cat\'egorie simpliciale ici)
des foncteurs faibles de $\Xx ^o$ vers $\Gamma$ gr\^ace
au travail de Cordier et Porter \cite{CordierPorter}.}
Dans ce cas on dispose
de la  notion de limite homotopique
introduite par Bousfield et Kan
\cite{BousfieldKan} et nos $\Gamma$-champs, que nous appelons
$0$-champs de Segal, constituent une classe de pr\'efaisceaux
simpliciaux d\'ej\`a consid\'er\'ee par Jardine (les pr\'efaisceaux flasques
par rapport \`a tout objet du site \cite {Jardine}) et, dans le cadre
tr\`es voisin
des pr\'efaisceaux de spectres, par Thomason \cite {Thomason}. A tout
pr\'efaisceau simplicial (ou $0$-pr\'echamp de Segal), on sait associer
naturellement un $0$-champ de Segal ayant les m\^emes faisceaux (et non
pr\'efaisceaux) d'homotopie.

\subnumero{Les $n$-champs de Segal}

Il est donc naturel d'\'etendre cette d\'emarche au cas des cat\'egories
simpliciales ou de Segal. Pour cela, il nous
faut munir la classe $\Gamma$ des cat\'egories de Segal d'une structure
d'$\infty$-cat\'egorie (techniquement, c'est une $2$-cat\'egorie
de Segal) o\`u l'on puisse
d\'efinir la notion de limite (homotopique). Pour prendre de la marge,
on d\'efinit carr\'ement les $n$-cat\'egories de Segal (\S 2): leur construction
suit de pr\`es celle des $n$-cat\'egories de Tamsamani \cite {Tamsamani}.
Les $0$-cat\'egories de Segal sont les ensembles simpliciaux,
et les $(n+1)$-cat\'egories de Segal sont des objets simpliciaux dans la
cat\'egorie des $n$-cat\'egories de Segal v\'erifiant les conditions
(de Segal) qui expriment que la composition, \`a d\'efaut d'\^etre
bien d\'efinie et associative, est d\'efinie et associative \`a
homotopie pr\`es. L'avantage des $n$-cat\'egories de Segal sur
les ensembles simpliciaux r\'eside dans le fait qu'on peut y trouver
des $i$-morphismes non-inversibles, du moins pour $i \leq n$.
Gr\^ace \`a la notion de
cmf {\em interne} (\S 11),
on montre qu'\`a des consid\'erations de nature ensembliste pr\`es,
les $n$-cat\'egories
de Segal s'organisent en une $(n+1)$-cat\'egorie de Segal $nSeCAT$ (\S 11),
o\`u l'on sait d\'efinir les limites homotopiques (\S 14). On peut donc
introduire les ($\Gamma$-)champs de $n$-cat\'egories de Segal, que nous
appelons $n$-champs de Segal. Techniquement, on commence par donner une
d\'efinition \`a la Jardine des $n$-champs de Segal (\S 9) et ce n'est
qu'au \S 14 qu'on retrouve la belle d\'efinition mentionn\'ee plus haut.
Le probl\`eme de la descente des complexes prend alors la forme suivante:
le $1$-pr\'echamp de Segal des complexes est-il un $1$-champ de Segal?

\subnumero {Donn\'ees de descente g\'en\'eralis\'ees}

Cette formulation n'est pas totalement satisfaisante, dans la mesure o\`u
l'on ne voit pas suffisamment bien comment elle recouvre les exemples concrets
mentionn\'es plus haut. Heureusement, nous savons introduire une notion de
donn\'ee de descente (et aussi, bien s\^ur, de donn\'ee de descente effective)
\`a valeurs dans un $n$-pr\'echamp de Segal, qui joue le r\^ole qu'on attend
d'elle. En fait, on explique ici ce qu'on appelle
{\em donn\'ee de descente g\'en\'eralis\'ee} (dddg) et on commence par le
cas plus familier des pr\'efaisceaux simpliciaux:
une dddg
\`a valeurs dans le pr\'efaisceau simplicial
$A$ sur le site $\Xx$ est un morphisme de pr\'efaisceaux
simpliciaux $\delta: D \rightarrow A$ o\`u $D$ est {\em contractile}
en ce sens que le $0$-champ de Segal associ\'e \`a $D$
est
\'equivalent \`a
l'objet final $*_{\Xx}$.
(Un tel $D$ est essentiellement la m\^eme chose qu'un {\em
hyper-recouvrement} au sens de Verdier \cite{SGA4b}.)
Une telle donn\'ee $\delta$ est
dite {\em effective} si elle
est \'equivalente \`a une donn\'ee $\delta': D \rightarrow A$ qui se factorise
\`a travers $*_{\Xx}$. Ces d\'efinitions s'\'etendent sans changement au cas
o\`u $A$ est un $n$-pr\'echamp de Segal, les pr\'efaisceaux simpliciaux
pouvant aussi
\^etre vus comme des $n$-pr\'echamps de Segal (\S 2).

Il convient d'observer que la question de l'effectivit\'e des
dddg dans un $n$-pr\'echamp de
Segal $A$ ne concerne que l'{\em int\'erieur}
$A^{int}$ de $A$: l'int\'erieur d'une $n$-cat\'egorie de Segal
$S$ est l'ensemble
simplicial obtenu \`a partir de $S$ grosso modo en ne conservant que les
morphismes (de tous ordres) qui sont inversibles \`a homotopie pr\`es, et cette
construction s'\'etend \'evidemment aux $n$-pr\'echamps de Segal.

La d\'efinition pr\'ec\'edente rend bien compte des exemples ``concrets''
mentionn\'es plus haut.
Ainsi, pour le premier exemple, on peut prendre le pr\'efaisceau simplicial
(contractile) $D$
tel que $D(U)$ soit l'ensemble
\footnote {Ici comme ailleurs dans ce travail, on ne pr\^ete pas trop
d'attention
aux questions li\'ees \`a la diff\'erence entre ensemble et classe. }
simplicial des
r\'esolutions projectives finies de la restriction de $\Ff$ \`a $U$:
les $1$-simplexes
sont les quasi-isomorphismes entre r\'esolutions
(compatibles \`a l'augmentation)\ldots .

\subnumero{Caract\'erisation des $n$-champs de Segal par la descente}

Pour expliquer en quoi
nos notions de champs sont \`a leur tour adapt\'ees \`a
ces notions de donn\'ees de descente, il nous faut introduire
quelques notations.
Si $A$ est une $n$-cat\'egorie de Segal, les ($1$-)morphismes de $A$
s'organisent en une nouvelle $(n-1)$-cat\'egorie de Segal qu'on note $Fl(A)$.
Si $x$ et $y$ sont deux objets de $A$, on note $Fl(A)(x,y)$ la sous-
$(n-1)$-cat\'egorie de Segal de $Fl(A)$ des morphismes de source $x$ et but
$y$. Plus g\'en\'eralement, si $x$ et $y$ sont deux objets de $Fl^{i-1}(A)$,
on note $Fl^{i}(A)(x,y)$ la sous-
$(n-i)$-cat\'egorie de Segal de $Fl^i(A)$ des morphismes de source $x$ et
but $y$.
\footnote{
On fait la convention que pour un $i$ donn\'e on remplace nos
$n$-cat\'egories de Segal par des $n'$-cat\'egories de Segal \'equivalentes,
avec $n'>i$.}
Ces notations s'\'etendent au cas o\`u $A$ est un $n$-pr\'echamp de
Segal sur un site $\Xx$: $Fl^i(A)$ est alors un nouveau $(n-i)$-pr\'echamp de
Segal  sur $\Xx$, et si $x$ et $y$ sont des sections de $Fl^{i-1}(A)$ sur un
objet $U$ de $\Xx$, $Fl^{i}(A)(x,y)$ est un $(n-i)$-pr\'echamp de Segal sur le
site $\Xx/U$.
Avec ces notations,
on montre le r\'esultat suivant:
\begin{theoreme}
Soit $A$ est un $n$-pr\'echamp de Segal sur le site $\Xx$.
On suppose que les valeurs de $A$ sont
des $n$-cat\'egories de Segal (cette condition est vide pour $n=0$).
Alors les conditions suivantes
sont \'equivalentes:
\newline
$(i)$ $A$ est un $n$-champ de Segal;
\newline
$(ii)$
(a)\,\, pour tout $X\in \Xx$ et tout $x,y\in A_0(X)$ le $n$-pr\'echamp
de Segal $Fl(A)(x,y)$ sur $\Xx /X$ est un $n$-champ de Segal; et
\newline
(b) \,\, les donn\'ees de descente pour $A$
sont effectives.
\newline
$(iii)$ (valable seulement pour les $n$-champs non de Segal)\newline
(a')\,\, pour tout $i \ge 1$, pour tout
$X\in \Xx$ et pour tous $x,y\in Fl^{i-1}(A)(X)$ les donn\'ees de descente
\`a valeurs dans le $(n-i)$-pr\'echamp $Fl^{i}(A)(x,y)$ sont effectives; et
\newline
(b) \,\, les donn\'ees de descente \`a valeurs dans $A$
sont effectives.
\end {theoreme}

Ce th\'eor\`eme formalise l'id\'ee (qui ressort entre autres de Giraud
\cite{GiraudThese} et Laumon-Moret Bailly \cite{LaumonMB}) qu'un champ est un
pr\'echamp o\`u les donn\'ees de descente sont effectives, et comme
l'effectivit\'e des donn\'ees de descente ne concernent que l'int\'erieur,
il caract\'erise nos
$n$-champs de Segal en des termes propres \`a la th\'eorie des
pr\'efaisceaux simpliciaux.

\subnumero{Donn\'ees de descente et topologie grossi\`ere}

Cependant, la d\'efinition m\^eme des donn\'ees de descente
(g\'en\'eralis\'ees) qu'on a formul\'ee plus haut faisait
intervenir la notion de $0$-champ de Segal et il est donc temps de
donner notre d\'efinition des (vraies) donn\'ees de descente \`a valeur
dans un pr\'efaisceau simplicial, puisqu'il reste vrai que ce cas suffit.
Soit $\Xx$ un site muni d'un crible $\Bb$ couvrant et
d'un pr\'efaisceau simplicial $A$. Une donn\'ee de descente sur $\Bb$
\`a valeurs dans $A$ est un diagramme $\delta$ de pr\'efaisceaux
simpliciaux:
$$
*_{\Bb} \leftarrow D \rightarrow A
$$
o\`u la fl\`eche de gauche est ce que nous appelons une
{\em \'equivalence objet-par-objet}, c'est-\`a-dire qu'elle induit
une \'equivalence faible d'ensembles simpliciaux sur chaque objet
$U \in \Xx$ (ici, bien s\^ur, seuls les objets de $\Bb$ sont vraiment
concern\'es). Ceci implique que $D$ est contractile et on est bien en
pr\'esence d'une
donn\'ee de descente dans le sens g\'en\'eralis\'e donn\'e plus
haut. On peut observer que la donn\'ee du diagramme $\delta$ \'equivaut
\`a celle d'une donn\'ee de descente (g\'en\'eralis\'ee ou non,
c'est pareil) \`a valeurs dans la restriction de $A$ \`a la cat\'egorie
$\Bb$ munie de la topologie grossi\`ere. Il s'agit l\`a d'un exemple
tout-\`a-fait typique de la fa\c{c}on dont la topologie grossi\`ere
intervient syst\'ematiquement tout au long de notre travail. Un autre
exemple typique est la fa\c{c}on dont nous d\'ecrivons en deux temps
la structure
de cat\'egories de mod\`eles pour nos pr\'echamps, par localisation de
Bousfield \`a partir de la structure correspondant \`a
la topologie grossi\`ere (\S 6).

On dit qu'une donn\'ee de descente est effective si elle est \'equivalente
\`a une donn\'ee de descente qui se factorise \`a travers $*_{\Xx}$. Et
lorsque dans le th\'eor\`eme ci-dessus on dit
que les donn\'ees de descente \`a valeurs dans $A$ sont effectives,
on veut dire que pour
tout objet $X$ de $\Xx$ et tout crible $\Bb$ couvrant $\Xx/X$, toute
donn\'ee de descente sur $\Bb$ \`a valeurs dans $A*{|\Xx/X}$ est effective.

\subnumero{Les deux temps de la descente et les constructions de Grothendieck}

Le probl\`eme de l'effectivit\'e des donn\'ees de descente dans un
$n$-pr\'echamp de Segal $A$ sur un site $\Xx$
(admettant un objet final) se d\'ecompose naturellement
en deux \'etapes. En effet, une donn\'ee de descente $\delta$ comme plus haut
s'ins\`ere dans un diagramme:
$$
*_{\Xx} \leftarrow *_{\Bb} \leftarrow D \rightarrow A,
$$
Puisqu'il s'agit de factoriser, \`a \'equivalence pr\`es, le morphisme
donn\'e $D \rightarrow A$ \`a travers $*_{\Xx}$, on peut d'abord
chercher une telle factorisation \`a travers $*_{\Bb}$.
Si donc on qualifie de strictes les donn\'ees de descente avec $D=*_{\Bb}$,
on peut formuler d'une part le probl\`eme de la strictification des
donn\'ees de descente, qui consiste \`a montrer que toute donn\'ee de
descente est \'equivalente \`a une donn\'ee de descente stricte, et d'autre
part celui de l'effectivit\'e des donn\'ees de descente strictes.

Techniquement, les choses se pr\'esentent un peu diff\'eremment. Pour montrer
qu'un pr\'echamp de Segal $A$ sur $\Xx$ est un champ de Segal, on utilise
la d\'efinition par limite. On
doit donc comparer
$A(X)$ avec la limite homotopique mentionn\'ee plus haut, disons
$\Gamma(\Bb, A)$. Pour identifier cette limite, on suit le chemin trac\'e dans
\cite{SGA1}: on construit une cat\'egorie de Segal
$\int _{\Bb} A$ au-dessus de $\Bb^o$ avec l'espoir d'identifier $\Gamma(\Bb, A)$
\`a une cat\'egorie de sections de $\int _{\Bb} A \rightarrow \Bb^o$. En fait,
on trouve (\S 16)
que $\Gamma(\Bb, A)$ est \'equivalente \`a la cat\'egorie de Segal
des sections ``\'equilibr\'ees'' (ou eq-sections) non pas de
$\int _{\Bb} A \rightarrow \Bb^o$ mais de son remplacement fibrant
$\int '_{\Bb} A \rightarrow \Bb^o$. Ce remplacement fait r\'ef\'erence \`a
une structure de cmf pour les cat\'egories de Segal (\S 2). On veut donc montrer
qu'un morphisme compos\'e
$$
A(X) \rightarrow
Sect^{eq}(\int _{\Bb} A) \rightarrow
Sect^{eq}(\int' _{\Bb} A)
$$
est une \'equivalence.
On peut voir le terme de droite comme une cat\'egorie
de donn\'ees de descente (\`a droite ?)
et celui du milieu comme la sous-cat\'egorie
des donn\'ees de descente strictes (on n'a pas
cherch\'e \`a \'etablir de correspondance entre ces nouvelles notions de
donn\'ees de descente et les anciennes), et le probl\`eme de l'effectivit\'e
se d\'ecompose \`a nouveau en deux \'etapes.

Ces deux \'etapes seront trait\'es chacun \`a son tour dans nos \S 18 et \S
19. Pour rester strictement en conformit\'e avec ce qui se passe dans le
papier, il faudrait dire que cette situation de ``strictification'' des
donn\'ees de descente a lieu non au-dessus de $\Bb$ mais au-dessus de la
cat\'egorie $\Delta$ via un foncteur $\rho (U):\Delta \rightarrow \Bb$,
dans le cas o\`u $\Bb$ est le crible engendr\'e par un objet $U$ recouvrant
$X$. En quelque sort $\Delta$ approxime $\Bb$ via le foncteur
$\rho (U)$ et pour des raisons techniques il nous est plus commode de
travailler au-dessus de $\Delta$.

Il se trouve que le probl\`eme de l'effectivit\'e des donn\'ees de descente
strictes a d\'ej\`a \'et\'e r\'esolu
pour nos complexes de modules sur un topos
annel\'e  dans
l'expos\'e Vbis de SGA 4 II \cite{SGA4} de Saint-Donat (d'apr\`es des
``notes
succinctes'' de Deligne). Au \S 19, nous donnons un r\'esultat plus g\'en\'eral
pour cette \'etape ``stricte''.

\subnumero{Pr\'efaisceaux de Quillen}

Nous avons donc mis en place ce formalisme des champs de Segal qui
constitue un cadre naturel pour la
th\'eorie de l'homotopie relative (au-dessus d'un site), dans lequel nous
pouvons reformuler le probl\`eme de la descente: trouver des conditions
pour qu'un pr\'echamp de Segal soit un champ (de Segal).
Comme l'\'ecrit Jardine dans \cite{Jardine}, la structure organisationnelle
pour la th\'eorie de l'homotopie est celle de cat\'egorie de
mod\`eles ferm\'ee (cmf) de Quillen \cite{Quillen}. La th\'eorie des champs
de Segal que nous avons mise en place est le cadre naturel de la
th\'eorie de l'homotopie relative (au-dessus d'un site) et la structure
organisationnelle correspondante est tout naturellement une structure de
pr\'efaisceau de cmf: nous appelons pr\'efaisceau de Quillen tout
pr\'efaisceau de cmf dans lequel les restrictions sont
ce qu'on appelle des foncteurs de
Quillen (\`a gauche). Si ${\bf M}$ est un tel pr\'efaisceau de Quillen sur
un site $\Xx$, on peut lui associer le $1$-pr\'echamp de Segal
$L({\bf M})$ sur $\Xx$ obtenu en prenant pour $L({\bf M})(X)$ la cat\'egorie
simpliciale obtenue en inversant (au sens de Dwyer-Kan) les
\'equivalences  dans ${\bf M}(X)$. On cherche donc des conditions
suffisantes
pour que le pr\'echamp de Segal ainsi associ\'e \`a un pr\'efaisceau de
Quillen soit un champ de Segal.

\subnumero {Un th\'eor\`eme de descente}

Nous proposons donc un th\'eor\`eme (\ref {dansleschamps}) donnant
des conditions suffisantes  pour que le pr\'echamp $L({\bf M})$
soit un champ de Segal. C'est bien s\^ur \`a partir du cas des complexes
que nous avons trouv\'e une telle formulation g\'en\'erale.
Les conditions que nous imposons sont de trois ordres.

D'abord, nous nous restreignons \`a des sites o\`u, grosso modo,
on peut se r\'eduire \`a ne consid\'erer que les donn\'ees de descente
d\'efinies sur les recouvrements \`a un seul objet (et v\'erifiant
les conditions de compatibilit\'e ad\'equates sur le complexe de
Cech du recouvrement).

La principale restriction que
nous imposons au pr\'efaisceau de Quillen est l'existence,
dans les valeurs prises (ce sont des cmf), de petites limites et
colimites arbitraires. Ainsi par exemple notre th\'eor\`eme ne concerne
pas les complexes de faisceaux de type fini.
Nous lui imposons
\'egalement d'autres conditions plus ou moins
naturelles \`a savoir: le fait que, dans les valeurs prises,
tous les objets soient cofibrants; l'existence de factorisations
fonctorielles;
la compatibilit\'e aux produits fibr\'es; le
caract\`ere local des \'equivalences; et la compatibilit\'e aux sommes
disjointes.

Enfin nous imposons une condition tr\`es technique autorisant
l'utilisation de l'argument dit ``du petit objet'' de Quillen
pour la construction d'une cmf auxiliaire.

La partie la plus d\'elicate de la d\'emonstration
consiste en la strictification des donn\'ees de descente (\S 18).

\noindent
--------
\newline
Dans la version 2 du papier nous avons d\^u ajouter l'hypoth\`ese (4) que
les morphismes de restriction pour $L(\mm )$ sont exacts des
deux cot\'es (et non seulement du cot\'e correspondant au fait que $\mm$ est
de Quillen \'a gauche). En effet, la premi\`ere version comportait
une faute dans une ligne
de la preuve, o\`u nous faisions commuter un morphisme de restriction avec une
limite; nous ajoutons cette commutation comme hypoth\`ese puisqu'elle n'est pas
une cons\'equence 
des autres.
Cette condition semble \^etre primordiale, 
heureusement
elle se v\'erifie dans les applications aux \S 20 et \S 21.

\subnumero {Champs de champs}

C'est en \'etudiant notre probl\`eme de descente des complexes,
que nous nous sommes rendu compte qu'il y avait un r\'esultat analogue pour les
$n$-champs eux-m\^emes: le $n+1$-pr\'echamp $n\underline{CHAMP}$ d\'efini par
$$
n\underline{CHAMP}(X):= nCHAMP (\Xx /X)
$$
est un $n$-champ (voir les th\'eor\`emes \ref{descente} et \ref{descenteNdS}).
C'est un resultat de ``recollement'' pour les $n$-champs. Par
exemple, \'etant donn\'e des champs $\Ff _U$ et $\Ff _V$ sur deux ouverts
$U,V$ d'un  espace topologique, avec une \'equivalence
$$
\Ff _U |_{U\cap V} \cong
\Ff _V |_{U\cap V}
$$
alors il existe un champ $\Ff $ sur $U\cup V$ avec $\Ff |_U \cong \Ff _U$
et $\Ff |_V \cong \Ff
_V$. On peut voir \ref{descente} et \ref{descenteNdS}
comme des g\'en\'eralisations aux $n$-champs, du r\'esultat classique qui
dit que le
$1$-prechamp qui \`a
chaque ouvert $U$ associe la cat\'egorie des faisceaux d'ensembles (i.e.
$0$-champs) sur $U$, est un
$1$-champ.  Le cas $n=1$ de cette g\'en\'eralisation se trouve d\'ej\`a
dans SGA 1
\cite{SGA1} \cite{GiraudThese},
o\`u l'on montre que le $2$-pr\'echamp $1\underline{CHAMP}$ des $1$-champs est
un $2$-champ.

Ce r\'esultat pour $n=2$ (qui dit que $2\underline{CHAMP}$ est un $3$-champ) a
\'et\'e utilis\'e, sans d\'emonstration, par Breen \cite{Breen}. On peut
dire que notre r\'esultat g\'en\'eralise approximativement la moiti\'e de ce
que fait Breen dans
\cite{Breen} (l'autre moiti\'e \'etant d'exprimer les donn\'ees de descente
en termes de cocycles,
une chose que nous ne sommes pas arriv\'ees \`a faire
de fa\c{c}on enti\`erement satisfaisante pour $n$ quelconque).

En fait, au d\'ebut de ce travail, la seule chose que nous pouvions
d\'emontrer \`a propos des complexes
\'etait que le pr\'echamp des complexes de $\qq$-espaces vectoriels, de
longueur $n$ fix\'ee, \'etait un $n+1$-champ. Pour prouver ceci, nous avions
utilis\'e le r\'esultat disant que le $n+1$-pr\'echamp des $n$-champs
est un champ: avec ceci et une version
\'el\'ementaire de la notion de ``spectre'' (consistant juste \`a faire une
translation
pour que tout se passe entre les degr\'es $N$ et $N+n$ avec $N\geq n+2$) on
pouvait obtenir (gr\^ace au fait que l'homotopie rationelle stable est
triviale) le fait que les complexes de $\qq$-espaces vectoriels de longueur
$n$, est un $n+1$-champ.

Ce n'est qu'avec l'apport des techniques de pointe en th\'eorie des
cmf (\cite{DwyerKan1} \cite{DwyerKan2} \cite{DwyerKan3}
\cite{DwyerKanDiags}, \cite{Hirschhorn},\cite{DHK}, \cite{JardineGoerssBook},
\cite{HoveyMonMod}
\cite{HoveyBook}) que nous avons pu formuler le r\'esultat plus g\'en\'eral
\ref{dansleschamps} qui permet de traiter le cas des complexes de
$\Oo$-modules.

\subnumero {Alg\'ebricit\'e}

L'\'etude plus d\'etaill\'ee des champs de complexes sur les sites
de la g\'eom\'etrie alg\'ebrique fera l'objet d'un autre travail. Nous
\'enon\c{c}ons simplement au \S 21
notre principal r\'esultat d'alg\'ebricit\'e.

\subnumero {Historique, travaux connexes et techniques utilis\'ees}

La notion de $1$-champ est aujourd'hui bien enracin\'ee  dans la
g\'eom\'etrie alg\'ebrique, depuis notamment les travaux  d'Artin
\cite{ArtinInventiones}, de Deligne-Mumford \cite{Deligne-Mumford},
de Laumon-Moret Bailly \cite{LaumonMB}, prolong\'es par de nombreux travaux
actuels. Une notion semblable joue un r\^ole important en
g\'eom\'etrie diff\'erentielle (e.g. en g\'eom\'etrie des feuilletages),
voir par
exemple \cite{Ehresmann} \cite{Morita} \cite{Pradines} \cite{Haefliger}
\cite{Haefliger2}
\cite{Molino} \cite{Moerdijk}.
\footnote{
On peut citer Moerdijk (\cite{Moerdijk} Example 1.1 (f)):
``Effective etale groupoids are well-known to be basically equivalent
to pseudogroups\ldots''
suivi par une r\'ef\'erence \`a Molino \cite{Molino}. Au vu de cette remarque,
les travaux
d'Ehresmann {\em et al} concernent essentiellement des groupo\"{\i}des
\'etales i.e. des $1$-champs
de Deligne-Mumford dans le cadre diff\'erentiel.
Il n'est pas clair \`a nos yeux
si cette connexion
\'etait connue d'Ehresmann \`a l'epoque; en tout cas, ce point de vue est
d\'ej\`a present dans les
papiers de Pradines \cite{Pradines}; et selon Pradines
(communication orale au deuxi\`eme
auteur),
ce point de vue a
\'et\'e present\'e par Ehresmann dans son expos\'e au congr\`es international
de 1956.}
En
s'appuyant sur la notion de champ, Giraud a d\'evelopp\'e en bas degr\'es la
cohomologie non-abelienne \cite{Giraud}.

Un travail pr\'ecurseur de la descente pour l'$\infty$-champ des complexes est
celui sur
la ``descente cohomologique'' de B. Saint Donat et Deligne dans SGA 4
\cite{SGA4b}.

Dans \cite
{Grothendieck}, A. Grothendieck propose comme
g\'en\'eralisation la notion de {\em $n$-champ} ou m\^eme de {\em
$\infty$-champ}, sans toutefois donner une d\'efinition
pr\'ecise.  En fait, on peut dire aue
la th\'eorie des $\infty$-champs de groupo\"{\i}des
est d\'ej\`a disponible
\`a l'epoque de \cite{Grothendieck}
par le biais de la
th\'eorie des pr\'efaisceaux simpliciaux, gr\^ace aux travaux de K. Brown
\cite{KBrown}, Illusie \cite{Illusie},  Joyal dans une lettre \`a Grothendieck
\cite{Joyal}, Jardine \cite{Jardine}, et Thomason \cite{Thomason}. Le biais
en question passe par
l'\'equivalence entre la th\'eorie des
$\infty$-groupoides et celle des espaces topologiques ou ensembles
simpliciaux---\'equivalence pour laquelle ne manquait que la d\'efinition de
$\infty$-groupo\"{\i}de!

Plus r\'ecemment, la th\'eorie des $n$-cat\'egories pour les petites valeurs de
$n$ a \'et\'e reconnue comme importante dans une large gamme de situations
li\'ees  par exemple \`a de nouveaux invariants des noeuds, aux groupes
quantiques,
\`a la th\'eorie quantique des champs, et m\^eme \`a l'informatique\ldots .
Ceci a conduit vers une meilleure compr\'ehension notamment du r\^ole jou\'e
par la non-associativit\'e de la composition. Le premier papier de Baez-Dolan
\cite{BaezDolan} fournit un r\'esum\'e int\'eressant de ces d\'eveloppements.

La th\`ese de Tamsamani \cite{Tamsamani} fournit une
d\'efinition de $n$-cat\'egorie,
ainsi que l'\'e\-qui\-va\-len\-ce entre les $n$-groupo\"{\i}des et
les espaces topologiques $n$-tronqu\'es, pour tout $n$. Cette d\'efinition
permet d'aborder la th\'eorie des $n$-champs, non n\'ecessairement de
groupo\"{\i}des,  en utilisant les
techniques dues \`a Quillen qui sont \`a  la base des travaux de Brown, Illusie,
Joyal, Jardine et Thomason. Signalons \'egalement que
d'autres d\'efinitions de $n$-cat\'egorie
ont \'et\'e propos\'ees, voir Baez-Dolan (\cite{BaezDolanLetter}
\cite{BaezDolanIII} qui \'etait ind\'ependant de \cite{Tamsamani} et
simultan\'e) ou Batanin (\cite{Batanin}, un peu plus tard).
Pour l'heure il manque une
comparaison entre ces diff\'erentes d\'efinitions, et comme expliqu\'e plus
haut, nous utilisons une l\'eg\`ere variante de la d\'efinition de
Tamsamani.

Notre technique de base est celle des cmf, et notre travail
a \'et\'e profond\'ement
influenc\'e par
des travaux r\'ecents en  th\'eorie des cmf, ceux de Hirschhorn
\cite{Hirschhorn}, de Dwyer-Hirschhorn-Kan \cite{DHK}, aussi ceux de
de Goerss-Jardine \cite{JardineGoerssBook} et de Hovey \cite{HoveyMonMod}
\cite{HoveyBook}, ainsi que par les travaux moins r\'ecents de
Dwyer-Kan sur la
localisation \cite{DwyerKan1} \cite{DwyerKan2} \cite{DwyerKan3}
\cite{DwyerKanDiags}.

\subnumero{Remerciements}

Nous remercions C. Walter, N. Mestrano, Z. Tamsamani, P. Hirschhorn, B. Toen,
R. Brown, A. Brugui\`eres, T. Goodwillie, L. Katzarkov, M. Kontsevich,
G. Maltsiniotis,  H. Miller, S. Mochizuki,  T. Pantev, J. Pradines, C. Rezk,
J.Tapia, C. Teleman, D. Toledo
pour discussions, correspondances, cours,
th\`eses, livres qui ont
influenc\'e le pr\'esent travail. Le deuxieme auteur remercie l'universit\'e de
Californie \`a Irvine pour son hospitalit\'e \`a la fin de la preparation de ce
papier.

Nous remercions nos familles pour leur soutien pendant ce travail.

\numero{Les $n$-cat\'egories de Segal}

\label{ncatsegpage}

Commen\c{c}ons par fixer quelques notations: les cat\'egories d'ensembles,
d'ensembles
simpliciaux, et d'espaces topologiques sont not\'ees respectivement
$Ens$, $EnsSpl$ et $Top$. Un ensemble simplicial est un foncteur
$\Delta ^o\rightarrow Ens$ o\`u $\Delta$ est la cat\'egorie dont les objets
sont les ensembles ordonn\'es $p:= \{ 0', \ldots , p'\}$ et les morphismes
sont les
applications croissantes. Si $X$ est un ensemble simplicial on notera
$X_p$ sa valeur sur $p\in \Delta$. Si $A: \Delta
^o\rightarrow \Cc$ est un foncteur, on notera $A_p$ ou parfois $A_{p/}$ la
valeur
qu'il prend sur $p\in \Delta ^o$.
On notera $|A|$ la r\'ealisation topologique de l'
ensemble simplicial $A$. Si $Y$ est une cat\'egorie (i.e. $1$-cat\'egorie) on
notera son nerf $\nu Y \in EnsSpl$. On utilisera l'abr\'eviation
``cmf'' pour ``cat\'egorie de mod\`eles ferm\'ee'' de \cite{Quillen}.

Les d\'efinitions qui suivent sont des g\'en\'eralisations de d\'efinitions de
Tamsamani \cite{Tamsamani} et Dunn \cite{Dunn}, qui \`a leur tour reprenaient
des id\'ees de Segal \cite{SegalTopology} \cite{Adams}. 
La notion de $1$-cat\'egorie de Segal apparait pour la premi\`ere fois
\`a la fin de Dwyer-Kan-Smith \cite{DKS}.
Voir la
comparaison ci-dessous.

On d\'efinit la notion de {\em $n$-pr\'ecat de
Segal}, et la cat\'egorie $nSePC$ de ces objets, par r\'ecurrence sur $n$. Un
ensemble pourra \^etre consid\'er\'e comme $n$-pr\'ecat de Segal: il y aura un
foncteur (pleinement fid\`ele) $Ens\rightarrow nSePC$.
Un objet dans l'image de ce foncteur sera appel\'e ``un ensemble discret''.
Pour $n=0$, une {\em $0$-pr\'ecat de
Segal} est par d\'efinition un ensemble simplicial; et le foncteur
$Ens\rightarrow 0SePC$ associe \`a un ensemble $S$ le foncteur
$\underline{S}:\Delta ^o\rightarrow Ens$ constant \`a valeurs $S$. Pour
$n\geq 1$, une {\em
$n$-pr\'ecat de Segal} est un foncteur
$$
\Delta ^o\rightarrow (n-1)SePC
$$
not\'e $p\mapsto A_{p/}$, tel que $A_{0/}$ soit un ensemble discret qu'on
notera aussi
$A_0$. Un {\em morphisme} de $n$-pr\'ecats de Segal est une transformation
naturelle de foncteurs $\Delta ^o\rightarrow (n-1)SePC$, ce qui d\'efinit la
cat\'egorie $nSePC$. Le foncteur $Ens\rightarrow nSePC$ est celui qui \`a un
ensemble $S$ associe le foncteur compos\'e
$$
\Delta ^o\stackrel{\underline{S}}{\rightarrow} Ens \rightarrow (n-1)SePC.
$$

On peut d\'evisser la r\'ecurrence dans cette d\'efinition. Si on veut
suivre les
notations de Tamsamani \cite{Tamsamani}, une $n$-pr\'ecat de Segal est donc un
ensemble $n+1$-simplicial, i.e. un foncteur
$$
(\Delta ^o)^{n+1}\rightarrow Ens
$$
not\'e $(m_1,\ldots , m_{n+1})\mapsto A_{m_1,\ldots , m_{n+1}}$,
qui satisfait les {\em  conditions de constance} selon lesquelles
pour $m_i=0$ le
foncteur
est constant en les variables $m_{i+1},\ldots , m_{n+1}$.

Ces conditions de
constance peuvent \^etre prises en compte une fois pour toutes en
introduisant la cat\'egorie $\Theta ^{n+1}$, quotient de $\Delta ^{n+1}$
\cite{nCAT}.
On rappelle  que $\Theta ^{n+1}$ est la cat\'egorie dont les objets sont les
suites $(m_1,\ldots , m_i)$ pour $i\leq n+1$, avec $m_i\in \Delta$ et avec
les identifications
$$
(m_1,\ldots , m_j=0, \ldots , m_i) =: (m_1,\ldots , m_{j-1}).
$$
L'objet $(0,\ldots , 0)$ est aussi \'equivalent \`a la suite de longueur
$0$ qui est not\'ee $0$. Les morphismes proviennent des morphismes de
$\Delta$ avec le changement induit par ces identifications; cf \cite{nCAT}
ou \cite{limits}.

Pour $k\leq n$ et pour $M=(m_1,\ldots , m_k)\in \Delta ^k$ on note $A_{M/}$
la $n-k$-pr\'ecat de Segal
$$
(m'_1,\ldots , m'_{n+1-k})\mapsto A_{m_1,\ldots , m_k, m'_1,\ldots ,
m'_{n+1-k}}.
$$
Pour $k=1$ ceci concorde avec la notation $A_{p/}$ introduite plus haut.

Une $n$-pr\'ecat
de Segal peut \'etre vue comme foncteur
$$
(\Theta ^{n+1})^o\rightarrow Ens,
$$
sans autre condition (les conditions de constance \'etant prises en compte par
le fait que $\Theta ^{n+1}$ est le quotient ad\'equat de $\Delta ^{n+1}$).
De ce fait, $nSePC$ appara\^{\i}t comme la cat\'egorie des pr\'efaisceaux
d'ensembles sur
$\Theta ^{n+1}$, et donc admet toutes
les limites et colimites avec petite indexation,
ainsi qu'un $\underline{Hom}$  interne.

En fait, le point de vue le plus utile est un m\'elange des deux points de
vue du
paragraphe pr\'ec\'edent: une $n$-pr\'ecat de Segal peut (voire doit) \^etre
consid\'er\'ee comme un foncteur
$$
A: (\Theta ^n)^o\rightarrow EnsSpl
$$
not\'e $M\mapsto A_M$, tel que, pour $|M| < n$, l'ensemble simplicial $A_M$
est un ensemble discret. Il est parfois utile de confondre les ensembles
simpliciaux et les espaces topologiques, et on obtient essentiellement la
m\^eme notion en consid\'erant les foncteurs $(\Theta ^n)^o\rightarrow Top$.

On dira que  $A_0$ est {\em l'ensemble des objets de $A$}. Pour deux objets
$x,y\in A_0$
on notera $A_{1/}(x,y)$ la pr\'e-image de $(x,y)\in A_0\times A_0$ par le
morphisme
$$
A_{1/}\rightarrow A_0\times A_0
$$
induit par les deux applications ``faces'' de la structure simpliciale.
On dira que $A_{1/}(x,y)$ est {\em la $n-1$-pr\'ecat de Segal des
morphismes entre $x$
et $y$}. On utilisera parfois cette m\^eme notation (qui est assez compacte)
pour l'ensemble des morphismes entre deux objets d'une cat\'egorie ou pour
l'ensemble
simplicial des morphismes entre deux objets d'une cat\'egorie simpliciale.

On d\'efinit maintenant les {\em
$n$-cat\'egories de Segal}, qui sont les $n$-pr\'ecats de Segal satisfaisant
certaines conditions de ``sp\'ecialit\'e'' g\'en\'eralisant la condition
originelle
de Segal \cite{Segal}. Pour formuler ces conditions on a besoin d'une grande
r\'ecurrence comme dans \cite{Tamsamani}, cf aussi le r\'esum\'e dans
\cite{limits}. On va d\'efinir simultan\'ement par r\'ecurrence sur $n$: la
condition pour qu'une $n$-pr\'ecat de Segal soit une $n$-cat\'egorie de
Segal; la condition pour qu'un morphisme entre deux $n$-cat\'egories de
Segal soit une \'equivalence; et l'ensemble ``tronqu\'e''
$\tau _{\leq 0}(A)$ associ\'e \`a une
$n$-cat\'egorie de Segal $A$. D'abord pour $n=0$, et par d\'efinition:
 toute $0$-pr\'ecat de Segal i.e.
ensemble simplicial, est une $0$-cat\'egorie de Segal; un morphisme est une
\'equivalence si c'est une \'equivalence faible d'ensembles simpliciaux, i.e.
s'il induit des isomorphismes entre les groupes d'homotopie (et les $\pi
_0$) des
r\'ealisations topologiques; et le tronqu\'e $\tau _{\leq 0}(A)$ d'un ensemble
simplicial $A$ est l'ensemble $\pi _0(|A|)$. Pour le reste, la r\'ecurrence est
identique \`a celle de Tamsamani \cite{Tamsamani}:  soit $n\geq 1$ et on suppose
que les d\'efinitions sont disponibles au-niveau $n-1$. Soit $A$ une
$n$-pr\'ecat de
Segal, vue comme foncteur
$$
A: \Delta ^o\rightarrow (n-1)SePC,\;\;\; p\mapsto A_{p/}.
$$
On dira que $A$ est une {\em $n$-cat\'egorie de Segal} si les
$A_{p/}$ sont des $n-1$-cat\'egories de Segal, et si pour tout $p$ le morphisme
de $n-1$-pr\'ecats de Segal
$$
A_{p/}\rightarrow A_{1/}\times _{A_0} \ldots \times _{A_0}A_{1/}
$$
est une \'equivalence de $n-1$-cat\'egories de Segal. Ce dernier
``morphisme de Segal'' \'etant celui dont les composantes sont les applications
duaux aux inclusions $\{ i,i+1\} \subset \{ 0,\ldots , p\}$ dans $\Delta$,
voir \cite{SegalInventiones}, \cite{Adams}, \cite{Dunn}, \cite{Tamsamani}.

Si $A$ est une $n$-cat\'egorie de Segal, alors le foncteur
$$
p\mapsto \tau _{\leq 0} (A_{p/}) \;\;\;\; \Delta ^o \rightarrow Ens
$$
est le nerf d'une cat\'egorie qu'on note $\tau _{\leq 1}(A)$. On
d\'efinit $\tau _{\leq 0}(A)$ comme l'ensemble des classes d'isomorphisme de la
cat\'egorie $\tau _{\leq 1}(A)$. Si $A$ et $B$ sont des $n$-cat\'egories de
Segal
et $f:A\rightarrow B$ est un morphisme (i.e. morphisme de $nSePC$) on dira que
$f$ est une {\em \'equivalence} si le morphisme induit
$$
\tau _{\leq 0}(A)\rightarrow \tau _{\leq 0}(B)
$$
est surjectif (on dira alors que ``$f$ est essentiellement surjectif''), et
si pour
tout couple d'objets $x,y\in A_0$, le morphisme
$$
A_{1/}(x,y)\rightarrow B_{1/}(f(x), f(y))
$$
est une \'equivalence de $n-1$-cat\'egories de Segal (on dira alors que
``$f$ est
pleinement fid\`ele''). Comparer avec la notion
d'\'equivalence entre cat\'egories simpliciales de
\cite{DwyerKanDiags}.

C'est seulement pour une $n$-cat\'egorie de Segal $A$ qu'il est raisonnable de
consid\'erer
la $n-1$-cat\'egorie des morphismes entre deux objets $A_{1/}(x,y)$.
Une {\em $1$-fl\`eche de $A$} est un objet de $A_{1/}(x,y)$.
Par it\'eration on obtient la notion de $i$-fl\`eche de $A$ \cite{Tamsamani}.
Soit $1^i\in \Theta ^n$ l'objet $(1, \ldots , 1)$ ($i$ fois), image de l'objet
$(1, \ldots , 1, 0,\ldots , 0)\in \Delta ^n$. Une {\em $i$-fl\`eche de $A$}
($i<n$) est un \'el\'ement de l'ensemble $A_{1^i}$. Pour $i=n$, $A_{1^n}$ est un
ensemble simplicial et une {\em $n$-fl\`eche de $A$} est un sommet (i.e.
un \'el\'ement de la partie de degr\'e $0$) de l'ensemble simplicial
$A_{1^n}$. Les morphismes faces en la $i$-i\`eme variable simpliciale
fournissent les applications {\em source} et {\em but} qui \`a une $i$-fl\`eche
associent des $i-1$-fl\`eches (une $0$-fl\`eche est un objet i.e. \'el\'ement
de $A_0$).

Soit $A$ une $n$-cat\'egorie de Segal.
Suivant \cite{Tamsamani}, on dira qu'une $i$-fl\`eche ($i\geq 1$) est {\em
inversible \`a \'equivalence pr\`es} ou {\em \`a homotopie pr\`es} si son image
comme $1$-fl\`eche dans la $1$-cat\'egorie $\tau _{\leq 1}(A_{1^{i-1}/})$ est
inversible. Voir \cite{Tamsamani} pour plus de d\'etails. On dira (toujours
suivant \cite{Tamsamani}) que $A$ est un {\em $n$-groupo\"{\i}de de Segal}
si pour
tout $1\leq i\leq n$, toutes les $i$-fl\`eches sont inversibles \`a
\'equivalence
pr\`es. Plus g\'en\'eralement, pour $0\leq k \leq n$ on dira que $A$ est
{\em $k$-groupique} si les $i$-fl\`eches sont inversibles \`a \'equivalence
pr\`es pour tout $i>k$. En particulier, $A$ est toujours $n$-groupique, et
$A$ est $0$-groupique si et seulement si c'est un $n$-groupo\"{\i}de de Segal.

\subnumero{Comparaison avec Tamsamani et Dunn}
On compare avec \cite{Tamsamani} et \cite{Dunn}. Tamsamani
(\cite{Tamsamani}, 1995)
consid\`ere les foncteurs $(\Delta ^o)^n\rightarrow Ens$ et leur impose les
conditions de constance ainsi que les conditions de Segal cf ci-dessous,
obtenant la notion de {\em $n$-cat\'egorie (faible)}. Dans \cite{nCAT},
on introduit la notion de {\em $n$-pr\'ecat} en relaxant toutes sauf les
conditions de constance.
Une $n$-pr\'ecat (resp. $n$-cat\'egorie de Tamsamani),
qui est un foncteur $\Theta ^n\rightarrow Ens$, peut \^etre vue comme une
$n$-pr\'ecat
de Segal en composant avec l'inclusion $Ens\rightarrow EnsSpl$.

Plus g\'en\'eralement, les $n$-pr\'ecats de Segal
(resp. $n$-cat\'egories de Segal)
telles que les valeurs du foncteur $\Theta ^n \rightarrow EnsSpl$
soient des ensembles simpliciaux $0$-tronqu\'es (i.e. r\'eeunion disjointe de
composantes contractiles) correspondent \`a des $n$-pr\'ecats
(resp. $n$-cat\'egories) par composition avec $\pi _0: EnsSpl
\rightarrow Ens$. Avec cette traduction, la plupart de nos \'enonc\'es sur les
$n$-cat\'egories de Segal restent vrais {\em mutatis mutandis} pour les
$n$-cat\'egories (la seule exception concerne la propri\'et\'e de
descente---non pr\'eserv\'ee par troncation---pour des localis\'ees de
Dwyer-Kan, qui ne sont pas tronqu\'ees en g\'en\'eral). Nous formulons
souvent nos
\'enonc\'es seulement dans le cadre des $n$-cat\'egories de Segal, laissant au
lecteur le soin d'\'ecrire les \'enonc\'es correspondants pour les
$n$-cat\'egories.
En revanche, on reviendra plus loin sur les compatibilit\'es entre les
notions de
$n$-cat\'egorie et $n$-cat\'egorie de Segal (et aussi sur les compatibilit\'es
entre diff\'erentes valeurs de $n$).

Dunn, dans \cite{Dunn} (soumis en 1993, paru en 1996 mais r\'esumant sa th\`ese
de 1984, et qui fait
r\'ef\'erence \`a Cobb \cite{Cobb})
consid\`ere---entre autres ``delooping machines''---le {\em $n$-i\`eme
it\'er\'e de la machine de Segal} dont la cat\'egorie sous-jacente est celle des
foncteurs
$$
A: (\Delta ^o) ^n\rightarrow Top
$$
tels que $A_{m_1,\ldots , m_n}=\ast$ si $m_i=0$ (cf. \cite{Dunn}, d\'ebut du
\S 3). On pourra appeler ces objets ``$n$-pr\'ecats de Dunn''.
En rempla\c{c}ant $Top$ par $EnsSpl$ et en passant au quotient, on obtient
les $n$-pr\'ecats de Segal $A$ v\'erifiant $A_M=\ast$ pour
$|M|<n$.

Dunn impose aussi \`a ses objets la
condition de sp\'ecialit\'e en chaque variable simpliciale
(\cite{Dunn} Definition 3.1) et on pourrait appeler les objets qui satisfont
\`a cette condition, les ``$n$-cat\'egories de Dunn''. On peut v\'erifier que
cette
condition co\"{\i}ncide avec la condition d'\^etre une $n$-cat\'egorie de
Segal, dans
le cas qu'il consid\`ere o\`u $A_M= \ast$ pour $|M|< n$.
En somme, Dunn avait regard\'e les $n$-cat\'egories de Segal ayant une
seule $i$-fl\`eche pour tout $i< n$.

Les d\'efinitions de Dunn et de Tamsamani
sont bas\'ees sur l'id\'ee d'it\'erer la construction de Segal
\cite{Segal}. Elles sont compl\'ementaires  comme le sont la
tour de Postnikov et la tour de Whitehead: chez Dunn (comme chez Whitehead)
il s'agit
de l'homotopie en degr\'es $\geq n$ tandis que chez Tamsamani (comme chez
Postnikov)
il s'agit de l'homotopie en degr\'es $\leq n$. La notion de {\em $n$-pr\'ecat de
Segal} permet d'int\'egrer ces deux points de vue. En outre, dans les variables
``$n$-cat\'egoriques'' on peut avoir des fl\`eches non (homotopiquement)
inversibles. Il s'agit donc d'une d\'efinition qui permet d'avoir de l'homotopie
en tout degr\'e, et des fl\`eches non-inversibles en degr\'es $\leq n$. C'est
un avatar de la notion (non encore totalement \'elucid\'ee) de {\em
$\infty$-pr\'ecat} qui
permettrait d'avoir des fl\`eches non-inversibles en tout degr\'e. Nous avons
d\'ecid\'e de faire la pr\'esente r\'edaction dans le cadre des $n$-pr\'ecats de
Segal, par souci de simplicit\'e et aussi parce que, sur les exemples que
nous avons en vue,
la non-inversiblit\'e des fl\`eches n'entre en jeu que pour un nombre fini
de degr\'es.

\subnumero{L'op\'eration $SeCat$}

L'une des id\'ees de base qui permettaient d'obtenir la structure de cat\'egorie
de mod\`eles ferm\'ee pour les $n$-pr\'ecats dans \cite{nCAT} \'etait
l'op\'eration $Cat$, qu'on peut r\'esumer tr\`es simplement en disant
qu'elle force
la condition d'\^etre une $n$-cat\'egorie. Ce proc\'ed\'e s'applique de la
m\^eme
fa\c{c}on dans
notre cadre des $n$-pr\'ecats de Segal: on aura une op\'eration qu'on notera
$A\mapsto SeCat(A)$ avec une transformation naturelle $i_A: A\rightarrow
SeCat(A)$, qui transformera toute $n$-pr\'ecat de Segal en
$n$-cat\'egorie de Segal.

En fait au lieu de fixer pr\'ecis\'ement l'op\'eration $SeCat$, on
donne la d\'efinition suivante.

\begin{definition}
\label{secat}
On appellera {\em op\'eration de type $SeCat$} tout couple $(F, i)$ o\`u
$F:nSePC\rightarrow nSePC$ est un foncteur et $i_A: A\rightarrow F(A)$ une
transformation naturelle, poss\'edant les propri\'et\'es suivantes:
\newline
(a)---pour tout $n$-pr\'ecat de Segal $A$, $F(A)$ est une $n$-cat\'egorie de
Segal;
\newline
(b)---$i_A$ est un isomorphisme sur les ensembles d'objets, et si $A$ est une
$n$-cat\'egorie de Segal alors $i_A$ est une \'equivalence de $n$-cat\'egories
de Segal; et
\newline
(c)---pour toute $n$-pr\'ecat de Segal $A$, le morphisme $F(i_A)$ est une
\'equivalence de $n$-cat\'egories de Segal.
\end{definition}

\begin{lemme}
\label{secatunique}
Il existe une op\'eration de type $SeCat$. Si $(F^1, i^1)$ et $(F^2, i^2)$
sont deux op\'erations de type $SeCat$ alors elles sont \'equivalentes en ce
sens que
$$
F^1(i^2_A): F^1(A) \rightarrow F^1(F^2(A))
$$
est une \'equivalence de $n$-cat\'egories de Segal pour toute $n$-pr\'ecat de
Segal $A$. En particulier, si $f:A\rightarrow B$ est un morphisme  de
$n$-pr\'ecats de Segal alors $F^1(f)$ est une \'equivalence de $n$-cat\'egories
de Segal si et seulement si $F^2(f)$ l'est.
\end{lemme}
{\em Preuve:}
Pour l'existence, on peut utiliser une d\'efinition analogue \`a celle de
\cite{nCAT} qui traite le cas des $n$-pr\'ecats (non de Segal) (pour une
autre approche, voir la discussion
plus loin).
Pour les ensembles
simpliciaux du dernier degr\'e,
il faut remplir toutes les "cornes" (extensions anodines
de Kan). Pour l'unicit\'e, on utilise la m\^eme d\'emonstration que celle de
\cite{nCAT} Proposition 4.2. En fait cette d\'emonstration doit \^etre coupl\'ee
avec la d\'emonstration du th\'eor\`eme \ref{cmf} dans une r\'ecurrence sur $n$
comme c'est fait dans \cite{nCAT}. On note que le cas $n=0$ est trivialement
vrai.

Si $f: A\rightarrow B$ et si $F^2(f)$ est une \'equivalence alors dans le
diagramme
$$
\begin{array}{ccc}
F^1(A) & \stackrel{F^1(i^2_A)}{\rightarrow} &  F^1(F^2(A))\\
{\scriptstyle F^1(f)}\downarrow && \downarrow {\scriptstyle F^1(F^2(f))}\\
F^1(A) & \stackrel{F^1(i^2_A)}{\rightarrow} &  F^1(F^2(A))
\end{array}
$$
la fl\`eche $F^1(F^2(f))$ est une \'equivalence car $F^1$ transforme
\'equivalences de $n$-cat\'egories de Segal en \'equivalences d'apr\`es
\ref{secat} (b). Les fl\`eches horizontales sont des \'equivalences par le
r\'esultat d'unicit\'e ci-dessus. Donc $F^1(f)$ est une \'equivalence.
\eop

On appellera $SeCat$ toute op\'eration qui satisfait les crit\`eres de la
d\'efinition \ref{secat}. Nous donnons maintenant une variante commode, d\'ej\`a
pr\'esente (pour les $n$-pr\'ecats) dans \cite{nCAT} et explor\'ee
dans \cite{effective} pour les
$1$-pr\'ecats de Segal. Pour $n=0$ on peut prendre le
foncteur identit\'e avec la transformation naturelle identit\'e. Soit $n\geq 1$
et supposons construite l'op\'eration $SeCat$ pour les $n-1$-pr\'ecats de Segal.
Soit $A$ une $n$-pr\'ecat de Segal. On d\'efinit une $n$-pr\'ecat de Segal
$Fix(A)$
avec les m\^emes objets que $A$ en posant  (pour $p\geq 1$) $$
Fix(A)_{p/}(x_0,\ldots , x_p):= SeCat(A_{p/}(x_0,\ldots , x_p)).
$$
Ceci a pour effet de transformer les $A_{p/}$ en $n$-cat\'egories de Segal.
On a un morphisme naturel $A\rightarrow Fix(A)$.
D'autre part, pour tout $m\geq 2$ on d\'efinit une $n$-pr\'ecat de
Segal $Gen[m](A)$ munie d'un morphisme $A\rightarrow Gen[m](A)$ induisant un
isomorphisme sur les ensembles d'objets, par
$$
Gen[m](A)_{q/} := A_{q/} \cup ^{\coprod_{q\rightarrow m} A_{m/}}
\coprod _{q\rightarrow m}\Gg [m]
$$
o\`u
$$
A_{m/}\stackrel{a}{\rightarrow} \Gg [m]
\stackrel{b}{\rightarrow} A_{1/}\times _{A_0}\ldots
\times _{A_0} A_{1/}
$$
est une factorisation du morphisme de Segal, o\`u $a$ est une
cofibration et $b$ une \'equivalence faible, (qui existe par le th\'eoreme
\ref{SeCmf} qu'on suppose connu pour $n-1$). Les coproduits dans la d\'efinition
de $Gen[m]$ sont pris avec tous les morphismes $q\rightarrow m$ dans
$\Delta$ qui
ne se
factorisent pas \`a travers une ar\^ete principale $1\rightarrow m$. Voir
\cite{nCAT} ou \cite{effective} pour comment d\'efinir
les morphismes de fonctorialit\'e $Gen[m](A)_{q/}\rightarrow
Gen[m](A)_{p/}$ pour $p\rightarrow q$.
Maintenant on d\'efinit $Gen$ comme le compos\'e pour tout $m\geq 2$ des
$Gen[m]$ et $SeCat$ s'obtient en it\'erant une infinit\'e
d\'enombrable $\omega$ de fois le compos\'e $Fix \circ Gen$.
Cette op\'eration d'apparence tres ineffective a en r\'ealit\'e d'assez
bonnes propri\'et\'es d'effectivit\'e, voir \cite{effective}.

On peut voir une $n$-pr\'ecat $A$ comme un {\em syst\`eme de g\'en\'erateurs
et relations} pour d\'efinir une $n$-cat\'egorie, et $SeCat(A)$ comme la
$n$-cat\'egorie ainsi engendr\'ee.

\subnumero{La structure de cat\'egorie de mod\`eles ferm\'ee}
On d\'efinit maintenant la structure de cat\'egorie de mod\`eles ferm\'ee sur
$nSePC$.

Un morphisme $A\rightarrow B$ de $n$-pr\'ecats de Segal sera une {\em
cofibration} si $A_M\rightarrow B_M$ est injectif pour tout $M\in \Theta ^{n+1}$
(cette condition est diff\'erente de la condition pour les
$n+1$-pr\'ecats de \cite{nCAT} car ici nous exigeons aussi
l'injectivit\'e pour $|M|=n+1$). Une cofibration est
donc simplement une injection de pr\'efaisceaux
d'ensembles.

On dira qu'un morphisme $A\rightarrow B$ de $n$-pr\'ecats de Segal est
une {\em
\'equivalence faible} si le morphisme
$SeCat(A)\rightarrow SeCat(B)$ est une \'equivalence
de $n$-cat\'egories de Segal. Cette propri\'et\'e est ind\'ependante du choix de
la construction $SeCat$ gr\^ace au lemme \ref{secatunique}.

Et on dira qu'un morphisme $A\rightarrow B$ de $n$-pr\'ecats de Segal
est une {\em
fibration} s'il poss\`ede la propri\'et\'e de rel\`evement pour les
cofibrations
triviales (celles qui sont des \'equivalences faibles).
Rappelons (\cite{Quillen}) que, par d\'efinition,
la propri\'et\'e de rel\`evement assure que
pour toute cofibration triviale $E\rightarrow E'$, et tout paire de morphismes
$E\rightarrow A$ et $E'\rightarrow B$ qui forment un carr\'e commutatif, il
existe un morphisme $E'\rightarrow A$ rendant commutatif les deux triangles
dans le carr\'e.

\begin{theoreme}
\label{SeCmf}
Les trois classes de morphismes
ci-dessus
font de la cat\'egorie des $n$-pr\'ecats de Segal
une cat\'egorie de mod\`eles ferm\'ee.
\end{theoreme}

La d\'emonstration entre dans le cadre plus g\'en\'eral de la comparaison
entre notions de $n$-cat\'egorie de Segal et/ou $n$-cat\'egorie, pour
diff\'erentes valeurs de $n$. On aborde cette discussion avant
la d\'emonstration de \ref{SeCmf}.

\subnumero{Changement de $n$}

La cat\'egorie de mod\`eles ferm\'ee des $n$-pr\'ecats de Segal de
\ref{SeCmf} est une sorte de substitut pour la notion d'$\infty$-cat\'egorie,
dans laquelle les $i$-fl\`eches sont inversibles \`a \'equivalence pr\`es, pour
$i>n$. Si on admet l'existence d'une th\'eorie des $\infty$-cat\'egories
convenable, l'\'equivalence
passera par la construction $A\mapsto \Pi _{\infty} \circ A$ qui transforme une
$n$-pr\'ecat de Segal en une $\infty$-pr\'ecat, avec un $\Pi _{\infty}$ qui fait
correspondre un
$\infty$-groupo\"{\i}de \`a tout espace topologique ou ensemble simplicial
de Kan.
Nous ne voulons pas entrer dans les d\'etails de la notion de
$\infty$-cat\'egorie
ici, car la version ``$n$-cat\'egories de Segal'' convient pour nos
applications.
Ce choix est d'autant plus raisonnable qu'on peut construire des
foncteurs $A\mapsto \Pi
_{m,Se}\circ A$ pour passer de $n$ \`a $n+m$.

Rappelons que Tamsamani \cite{Tamsamani} a d\'efini un foncteur
``$n$-groupo\"{\i}de
de Poincar\'e''
$$
\Pi _n : Top \rightarrow nCat
$$
o\`u $nCat$ est la sous-cat\'egorie pleine de la cat\'egorie $nPC$ des
$n$-pr\'ecats, form\'ee des $n$-cat\'egories. Ce foncteur g\'en\'eralise
le ``groupo\"{\i}de de Poincar\'e'' $\Pi _1(X)$ classique. Tamsamani a
montr\'e que
ce foncteur \'etablit une \'equivalence entre les th\'eories homotopiques
des espaces topologiques $n$-tronqu\'es (i.e. o\`u les $\pi _i$ s'annulent pour
$i>n$) et des $n$-groupo\"{\i}des (i.e. $n$-cat\'egories o\`u les fl\`eches sont
inversibles \`a \'equivalence pr\`es). Plus pr\'ecis\'ement, il a construit un
foncteur ``r\'ealisation topologique''
$$
\Re : nCat \rightarrow Top
$$
et des \'equivalences  $\Pi _n \circ \Re\cong 1$ et $\Re \circ \Pi _n \cong 1$
pour les restrictions de ces foncteurs aux cat\'egories $Top _n$ des espaces
topologiques $n$-tronqu\'es et $nGpd$ des $n$-groupo\"{\i}des. On adoptera
donc le
point de vue que les $n$-groupo\"{\i}des peuvent \'etre librement
remplac\'es par
les espaces topologiques (ou ensembles simpliciaux) $n$-tronqu\'es.

On fait ici quelques remarques qui utilisent des notions qui seront
trait\'ees plus loin dans le pr\'esent travail, pour offrir un autre \'eclairage
sur ce r\'esultat d'\'equivalence de Tamsamani (cependant,
le lecteur peut s'en passer jusqu'au sigle ``$\oslash$'').
En appliquant la
th\'eorie de Dwyer-Kan \cite{DwyerKan3} qui sera rappel\'ee au \S 8
ci-dessous, on peut inverser les \'equivalences et obtenir
des cat\'egories simpliciales "localis\'ees"
$L(Top_n)$ et  $L(nGpd)$.
Maintenant les foncteurs $\Pi _n$ et $\Re$  et
les \'equivalences de Tamsamani, en conjonction avec la proposition 3.5 et le
corollaire 3.6 de \cite{DwyerKan2} (cf aussi notre lemme \ref{nonadjoint}),
donnent une \'equivalence
$$
L(Top_n)\cong L(nGpd )
$$
de cat\'egories simpliciales. On peut pousser un peu plus loin en utilisant
nos notations du \S 11: on note $nCAT$ la $n+1$-cat\'egorie des
$n$-cat\'egories (voir \cite{nCAT}).
\footnote{
Le lecteur prendra soin de distinguer entre la typographie en majuscules $nCAT$
et la typographie qui comporte des minuscules $nCat$.
En effet, $nCAT$ d\'esigne la $n+1$-cat\'egorie des $n$-cat\'egories,
d\'efinie \`a
l'aide de la structure interne cf \S 11 et (\cite{nCAT} \S 7); tandis que
$nCat$ d\'esigne
la $1$-cat\'egorie stricte avec les m\^emes objets et les morphismes qui
respectent
strictement la structure (en particulier, c'est celle-ci qu'avait consid\'er\'ee
Tamsamani). On utilise cette m\^eme convention pour les champs plus loin.}
Notons $nGPD\subset nCAT$  la sous-$n+1$-cat\'egorie pleine des
$n$-groupo\"{\i}des. Alors $nGPD$ est $1$-groupique (i.e. ses
$n$-cat\'egories de
$Hom$ sont en fait des $n$-groupo\"{\i}des) donc \'egale \`a son int\'erieur
$1$-groupique: $nGPD^{int, 1}= nGPD$. Par le th\'eor\`eme \ref{intereqloc}
on a
$$
L(nPC)\cong L(nPC_f)
\cong nCAT ^{int, 1}.
$$
D'autre part, comme on a des remplacements fibrants fonctoriels, toute
cat\'egorie comprise entre $nPC$ et $nPC_f$  donne la m\^eme localis\'ee de
Dwyer-Kan; ceci s'applique notamment \`a $nPC_f \subset nCat \subset nPC$.
On obtient ainsi l'\'equivalence des $n+1$-cat\'egories
$$
L(nCat)\cong nCAT^{int, 1}
$$
(plus techniquement, il s'agit ici d'une \'equivalence qui passe
\'eventuellement par des
morphismes dans les deux sens, car aucune des deux $n+1$-pr\'ecats en
question n'est
fibrante).
La sous-cat\'egorie $nGpd\subset nCat$ est d\'efinie par une condition
invariante par \'equivalence, donc (par \ref{stabilite} ci-dessous) $L(nGpd)$
est la sous-cat\'egorie pleine de $L(nCat)$ form\'ee des objets qui sont des
$n$-groupo\"{\i}des. Par d\'efinition il en est de m\^eme pour $nGPD$ et donc
aussi pour $nGPD^{int, 1}$. On obtient ainsi l'\'equivalence
$$
L(nGpd) \cong nGPD.
$$
Si on applique maintenant la variante (imm\'ediate)
du th\'eor\`eme de comparaison de Tamsamani
qui donne l'\'equivalence entre localis\'es de Dwyer-Kan,
on obtient une \'equivalence de $n+1$-cat\'egories
$$
L(Top_n)\cong nGPD.
$$
D'autre part, on peut observer que les
ensembles simpliciaux
$L(Top _n)_{1/}(X,Y)$ ont le bon type d'homotopie, \`a savoir celui
de l'espace d'applications
$\underline{Hom}(X,Y)$, d\`es que, par exemple, $X$ est un CW-complexe
(voir \cite{DwyerKan3} etc).
Pour dire
les choses autrement, $L(Top_n)$ est \'equivalente \`a la cat\'egorie
simpliciale
des objets
($n$-tronqu\'es) fibrants et cofibrants pour n'importe quelle structure de
cat\'egorie de mod\`eles ferm\'ee simpliciale pour les espaces
topologiques. Donc
$L(Top_n)$ est ``la'' $n+1$-cat\'egorie des espaces $n$-tronqu\'es. Cette
\'equivalence donne donc des pr\'ecisions sur ce qu'on entend par
``les $n$-groupo\"{\i}des s'identifient aux espaces
topologiques
$n$-tronqu\'es''.

\noindent $\oslash$

\medskip

Au vu des remarques pr\'ec\'edentes (la version courte aurait suffi), on va
{\em d\'eclarer} qu'un $\infty$-groupo\"{\i}de n'est rien d'autre qu'un espace
topologique ou ensemble simplicial. Ce point de vue qui nous guide d\'ej\`a
dans la
d\'efinition des $n$-cat\'egories de Segal nous guide aussi pour
comparer ces notions pour diff\'erentes valeurs de $n$.

Comme expliqu\'e au d\'ebut de ce chapitre, on peut voir toute $n$-pr\'ecat de
Segal comme un pr\'efaisceau simplicial au-dessus de $\Theta ^n$, i.e. un
foncteur vers la cat\'egorie des ensembles simpliciaux,
$$
A:\Theta ^n \rightarrow EnsSpl,\;\;\;\; M\mapsto A_M
$$
tel que pour $|M| < n$ l'ensemble simplicial
$A_M$ est un ensemble discret.

De ce point de vue, on peut
composer avec le foncteur
$| \; |$ de r\'ealisation topologique des ensembles simpliciaux, puis avec
le foncteur $\Pi _m$ de Tamsamani \cite{Tamsamani}. De fa\c{c}on \'equivalente
on peut d'abord appliquer le foncteur
$Ex^{\infty}$ de Kan pour arriver dans les ensembles
simpliciaux de Kan, puis la variante adapt\'ee aux ensembles
simpliciaux du foncteur $\Pi _m$ de \cite{Tamsamani}.  Par ces deux
proc\'ed\'es (essentiellement \'equivalents), on obtient une
$n+m$-pr\'ecat qu'on note  $\Pi _m \circ A$.

Dans l'autre sens, si on a une $n+m$-pr\'ecat $B$, on peut lui appliquer
l'op\'eration $\Re _{\geq n}$ qui consiste \`a prendre la r\'ealisation
topologique (ou plut\^ot comme ensemble simplicial) de \cite{Tamsamani} des
$m$-pr\'ecats $B_{M/}$ pour $|M|=n$ (on garde les ensembles discrets $B_M$ pour
$|M|<n$). Le r\'esultat est une $n$-pr\'ecat de Segal
$\Re _{\geq n}(B)$. Si $A$ est une $n$-cat\'egorie de Segal n'ayant que des
ensembles simpliciaux $m$-tronqu\'es dans l'image du foncteur
$$
(M\mapsto A_{M/}): (\Theta ^n )^o\rightarrow EnsSpl
$$
consid\'er\'e ci-dessus, alors $\Re _{\geq n}\Pi _m\circ A$ est
\'equivalente \`a
$A$. Si $A$
ne v\'erifie pas cette condition de troncation alors cette construction a
pour effet
d'op\'erer la troncation de Postnikov sur les ensembles simpliciaux.
D'autre part,
si $A$ est une $n+m$-cat\'egorie $n$-groupique (i.e. o\`u les $i$-fl\`eches sont
inversibles \`a homotopie pr\`es pour $i>n$), alors $\Pi _m\circ \Re _{\geq
n}(A)$
est \'equivalente \`a $A$. Dans le cas contraire on obtient ainsi le
{\em compl\'et\'e $n$-groupique de $A$} (i.e. le morphisme universel
de $A$ vers une $n+m$-cat\'egorie $n$-groupique),
analogue du ``groupe-compl\'et\'e'' d'un
mono\"{\i}de simplicial en topologie alg\'ebrique (qui est le cas $n=0$).

La construction $A\mapsto \Pi _{m}\circ A$ est \`a homotopie pr\`es compatible
avec les coproduits. Si $A\rightarrow B$ est une cofibration (i.e. injection,
cf. la d\'efinition ci-dessous) et $A\rightarrow C$ un morphisme de
$n$-pr\'ecats
de Segal, le coproduit $B\cup ^AC$ est un coproduit d'ensembles simpliciaux
objet-par-objet au-dessus de $\Theta ^n$. On peut donc appliquer le th\'eor\`eme
de Seifert-Van Kampen pour $\Pi _m$ de \cite{nCAT}, objet-par-objet au-dessus de
$\Theta ^n$, pour obtenir l'\'equivalence
$$
(\Pi _m \circ B)\cup ^{(\Pi _m \circ A)}(\Pi _m \circ C)
\stackrel{\cong}{\rightarrow}
\Pi _m \circ (B\cup ^AC).
$$
Dans le m\^eme ordre d'id\'ees on remarque que si $A$ est un $n$-pr\'ecat
de Segal
alors il y a une \'equivalence naturelle de $m$-cat\'egories de Segal
$$
\Pi _m \circ SeCat(A) \stackrel{\cong}{\rightarrow} Cat(\Pi _m \circ A).
$$

On peut modifier le
paragraphe sur l'op\'eration $\Pi _m\circ A$ pour obtenir un
passage des $n$-cat\'egories de Segal aux $n+m$-cat\'egories de
Segal. Pour cela, on doit
d'abord d\'efinir un $\Pi _{m,Se}(X)$ pour tout espace topologique $X$,
qui sera une $m$-cat\'egorie de Segal.  Pour $m=0$ c'est tout simplement
l'ensemble
simplicial ``complexe singulier'' $Sing(X)$. Pour $m$ g\'en\'eral nous
appliquons d'abord le foncteur $\underline{\Omega}$ d\'efini par
Tamsamani
\cite{Tamsamani}, $m$
fois it\'er\'e comme dans
\cite{Tamsamani}, ce qui donne un foncteur
$$
\underline{\Omega}^m(X): \Theta ^m\rightarrow Top.
$$
Ensuite on compose avec le foncteur $Sing$  pour obtenir
$$
\Pi _{m,Se}(X):= Sing \, \underline{\Omega}^m(X).
$$
Maintenant si $A$
est une $n$-pr\'ecat de Segal on la consid\`ere comme
un foncteur
$$
\Theta ^n\rightarrow EnsSpl
$$
qu'on compose d'abord avec la r\'ealisation $| \; |$, puis avec le foncteur
$\Pi _{m,Se}$. On notera  $\Pi _{m,Se}\circ A$ le r\'esultat de
cette construction. C'est une
$n+m$-pr\'ecat de Segal. De m\^eme si $B$ est une $n+m$-pr\'ecat de Segal
on peut
d\'efinir sa {\em r\'ealisation en degr\'es $\geq n$} $\Re _{\geq n}(B)$.
On a l'\'equivalence
$$
A \cong \Re _{\geq n}(\Pi_{m,Se}\circ A)
$$
(sans condition de troncation sur $A$ cette fois-ci).
Si de plus $B$ est $n$-groupique, alors
on a $B\cong \Pi_{m,Se}\circ \Re _{\geq n}(B)$. Et de nouveau, dans le cas
g\'en\'eral, $\Pi_{m,Se}\circ \Re _{\geq n}(B)$
sera le compl\'et\'e $n$-groupique de $B$.

On a enfin la formule
$$
(\Pi _{m,Se} \circ B)\cup ^{(\Pi _{m,Se} \circ A)}(\Pi _{m,Se} \circ C)
\cong
\Pi _{m,Se} \circ (B\cup ^AC).
$$

\subnumero{Induction}

Une $n$-pr\'ecat (non de Segal) $A$ peut \^etre consid\'er\'ee comme
une $m$-pr\'ecat de
Segal pour tout $m\geq n$, qu'on peut noter $Ind _n^{m,Se}(A)$ avec
$$
Ind _n^{m,Se}(A)_{k_1, \ldots , k_{m+1}}:= A_{k_1, \ldots , k_n}
$$
pour $(k_1,\ldots , k_{m+1})\in \Delta ^{m+1}$.

{\em Avertissement:} Si $A$ est une $n$-pr\'ecat alors $SeCat(Ind _n^{m,Se}(A))$
n'est pas en g\'en\'eral \'equivalente \`a $Ind _n^{m,Se}(Cat(A))$. On a cette
\'equivalence seulement si $A$ est d\'ej\`a une $n$-cat\'egorie.

D'autre part il faut remarquer que si  $A$ est une $n$-cat\'egorie (non de
Segal)
alors $Ind _n^{m,Se}(A)$ est d'embl\'ee une $m$-cat\'egorie de Segal.

Pour toutes ces
raisons, on utilisera  cette op\'eration ``d'induction''
le plus souvent avec un argument $A$ qui est d\'ej\`a une $n$-cat\'egorie.

Le cas qu'on rencontrera le plus est celui o\`u $Y$ est une
$1$-cat\'egorie, qu'on
consid\`ere comme $1$-pr\'ecat en prenant son nerf $\nu Y$. Alors
pour tout $n\geq 1$ on dispose de la $n$-cat\'egorie de Segal
$Ind _n^{m,Se}(\nu Y)$.

En pratique on confondra $Y$ avec la
$n$-cat\'egorie de Segal induite.

\subnumero{Troncation}

Dans le cadre des comparaisons ci-dessus,
on d\'efinit maintenant une {\em troncation} $A\mapsto \tau _{\leq n}(A)$
pour les
$m$-cat\'egories de Segal. Le r\'esultat est
une $n$-cat\'egorie (non de Segal).
Si $n\geq m$ on pose
$$
\tau _{\leq n} (A) := \Pi _{n-m} \circ A
$$
avec les notations introduites plus haut. En particulier on a
$$
\tau _{\leq m}(A)= \pi _0\circ A.
$$
Pour $n<m$ on pose
$$
\tau _{\leq n}(A):= \tau _{\leq n}(\tau _{\leq m}(A)),
$$
o\`u la troncation $\tau _{\leq n}$ de droite est la troncation
qui transforme une $m$-cat\'egorie en $n$-cat\'egorie (voir \cite{Tamsamani}).
En termes des d\'efinitions pr\'ec\'edentes, on a
$$
\tau _{\leq n}(A)_M = A_M \;\; \mbox{pour}\;\; |M|<n, \;\;\;\;
\tau _{\leq n}(A)_M = \tau _{\leq 0}(A_{M/}) \;\; \mbox{pour}\;\; |M|=n.
$$

Cette troncation est essentiellement la ``troncation de Postnikov''. Si on veut
consid\'erer que nos $m$-cat\'egories de Segal $A$ sont en fait des
$\infty$-cat\'egories, on peut pr\'eferer la notation $\tau _{\leq n}(A)$
\`a la notation $\Pi _{n-m}\circ A$.

\begin{lemme}
Un morphisme $A\rightarrow B$ de $m$-cat\'egories de Segal est une \'equivalence
si et seulement si $\tau _{\leq n}A\rightarrow \tau _{\leq n}B$ est une
\'equivalence de $n$-cat\'egories pour tout $n$.
\end{lemme}
{\em Preuve:}
C'est vrai pour $m=0$ i.e. pour les ensembles simpliciaux, et on obtient le cas
g\'en\'eral par une r\'ecurrence \'evidente.
\eop

Soit $A$ une $m$-cat\'egorie de Segal. On dira que $A$ est
{\em contractile} si le morphisme canonique $A\rightarrow \ast$ est une
\'equivalence. On dira que $A$ est {\em $0$-tronqu\'ee} si $A$ est
une r\'eunion disjointe de $m$-cat\'egories de Segal contractiles (ici la
r\'eunion peut \'eventuellement \^etre vide). Soit $n\leq m$; on dira que
$A$ est {\em $n$-tronqu\'ee} si pour tout
$M$ avec $|M| = n$, la  $m-n$-cat\'egorie de Segal $A_{M/}$ est $0$-tronqu\'ee.
Enfin, soit $n\geq m$; on dira $A$ est {\em $n$-tronqu\'ee} si pour tout
$M$ avec $|M|=m$, l'ensemble simplicial $A_{M/}$ est $n-m$-tronqu\'e (i.e. ses
groupes d'homotopie s'annulent en tout degr\'e $i>n-m$).

Pour $n\geq m$ on a un morphisme naturel
$$
\Pi _{n-m, Se}\circ A\rightarrow
\tau _{\leq
n}(A).
$$
En revanche, pour $n<m$ il n'existe pas en g\'en\'eral de morphisme
$A\rightarrow Ind ^{m,Se}_n(\tau _{\leq n}(A))$ (ni dans l'autre sens). Ceci
r\'esulte de
l'existence possible de $i$-fl\`eches non-inversibles (la notion d'int\'erieur
introduite ci-dessous
permettra de contourner cette difficult\'e).  Cependant, si $A$ est
$n$-tronqu\'ee alors on a une \'equivalence naturelle
$$
A\stackrel{\cong}{\rightarrow} Ind ^{m,Se}_n(\tau _{\leq n}(A)).
$$
En fait, on peut v\'erifier qu'une $m$-cat\'egorie de
Segal $A$  est \'equivalente \`a une de la forme $Ind ^{m,Se}_n(B)$,
o\`u $B$ est une $n$-cat\'egorie, si et seulement si $A$ est $n$-tronqu\'ee.

\subnumero{L'int\'erieur}

La r\'ealisation fournit un moyen de passer d'une $n$-pr\'ecat de Segal \`a un
ensemble simplicial; \c{c}a correspond essentiellement \`a inverser toutes les
fl\`eches pour obtenir un $\infty$-groupo\"{\i}de qu'on peut voir comme un
ensemble
simplicial. Cette op\'eration est assez brutale
et n'est gu\`ere compatible
avec les questions de descente. On introduit ici une construction
duale (en un sens assez faible), {\em l'int\'erieur} d'une $n$-pr\'ecat
de Segal, qui est bien mieux adapt\'ee aux questions de descente et qui va nous
permettre
de d\'evisser de fa\c{c}on r\'ecurrente la notion de champ en termes d'
une notion d\'ej\`a bien connue pour les pr\'efaisceaux simpliciaux.

Si $A$ est une $n$-cat\'egorie de Segal, on d\'efinit {\em l'int\'erieur de $A$}
not\'ee $A^{int, 0'}$ d'abord comme la sous-$n$-cat\'egorie de Segal de $A$
contenant seulement les $i$-fl\`eches inversibles \`a homotopie pr\`es.
Plus pr\'ecis\'ement on dira qu'une $i$-fl\`eche $a\in A_{1^i}$ est {\em
inversible
\`a homotopie pr\`es} si l'image de $a$ dans $\tau _{\leq i}A$ est une
$i$-fl\`eche inversible (rappelons que la composition des $i$-fl\`eches dans la
$i$-cat\'egorie $\tau _{\leq i}A$ est bien d\'efinie, associative, etc.).
L'int\'erieur $A^{int, 0'}$ de $A$ est la sous-$n$-pr\'ecat de Segal
d\'efinie par la
condition que $a\in A_M$ est dans $A^{int,0'}_M$ si et seulement si
$f^{\ast}(a)$ est une $i$-fl\`eche inversible \`a  homotopie pr\`es, pour
tout morphisme $f: 1^i\rightarrow M$ de $\Theta ^{n+1}$, pour $i\leq n $.
On voit que $A^{int, 0'}$ est une $n$-cat\'egorie de Segal, et
plus pr\'ecis\'ement un
$n$-groupo\"{\i}de de Segal.  On introduit maintenant
$$
A^{int, 0}:= \Re _{\geq 0} A^{int, 0'},
$$
qui est un ensemble simplicial.
Cette derni\`ere transformation ne change pas grand-chose car $A^{int, 0'}$
est d\'ej\`a
un $n$-groupo\"{\i}de de Segal, et en particulier $\Pi _{n,Se}\circ A^{int,
0}$ est
\'equivalente \`a $A^{int, 0'}$.

En termes d'$\infty$-cat\'egories,
$A^{int, 0}$ est l'$\infty$-groupo\"{\i}de universel muni d'un morphisme vers
$A$.

Plus g\'en\'eralement on d\'efinit {\em l'int\'erieur $k$-groupique}
$A^{int, k}$ d'une $n$-cat\'egorie de Segal $A$, pour tout $k<n$. D'abord on
d\'efinit $A^{int, k'}\subset A$ comme la sous-$n$-pr\'ecat de Segal avec
$a\in A_M$ contenu dans $A^{int, k'}_M$ si et seulement si $f^{\ast}a$
est une $i$-fl\`eche inversible \`a homotopie pr\`es pour tout morphisme $f:
1^i\rightarrow M$ avec $k<i\leq n$. Autrement dit, c'est seulement
pour $i>k$ qu' on ne conserve que les
$i$-fl\`eches inversibles \`a homotopie pr\`es. Ici encore, on v\'erifie que
$A^{int, k'}$ est une $n$-cat\'egorie de Segal $k$-groupique (i.e. dont les
$i$-fl\`eches sont inversibles \`a homotopie pr\`es pour $i>k$). Comme
pr\'ec\'edemment, on introduit
$$
A^{int, k}:= \Re _{\geq k} A^{int, k'},
$$
qui est une $k$-cat\'egorie de Segal. On a
$$
\Pi _{n-k,Se}\circ A^{int, k} \cong A^{int, k'}.
$$
En termes d'$\infty$-cat\'egories,  $A^{int, k}$ est l'
$\infty$-cat\'egorie $k$-groupique universelle munie d'un morphisme vers $A$.

Si on veut appliquer ces constructions \`a une $n$-pr\'ecat de Segal il faut
d'abord appliquer l'op\'eration $SeCat$ pour avoir une $n$-cat\'egorie de Segal
(le cas \'ech\'eant, ceci sera sous-entendu dans la notation:
$A^{int, k}:= SeCat(A)^{int, k}$).

Pour $k<l<n$ on a $(A^{int, l})^{int, k}=A^{int , k}$ \`a une
\'equivalence---qu'on omettra de mentionner---pr\`es. D'autre part on a
$$
\tau _{\leq m}(A^{int, k})=\tau _{\leq m} A
$$
pour $m\leq k$, et
$$
\tau _{\leq m} (A^{int, k})=(\tau _{\leq m})^{int, k}
$$
pour $k<m$.
Si $A\rightarrow B$ est une \'equivalence de $n$-cat\'egories de Segal alors il
en est de m\^eme de $A^{int, k}\rightarrow B^{int , k}$. De plus, $A$ est
$k$-groupique si et seulement si $A^{int, k}\rightarrow A$ est une
\'equivalence.

Enfin, on peut utiliser ces constructions pour relier une $n$-cat\'egorie de
Segal $A$ \`a sa troncation $\tau _{\leq k}(A)$. Rappelons que pour $k\geq n$
on a d\'ej\`a un morphisme $\Pi _{k-n,Se}\circ A\rightarrow \tau _{\leq
k}(A)$. Si de plus $k<n$, on a un couple de morphismes
$$
\Pi _{n-k,Se} \circ A^{int, k} \rightarrow A
$$
et
$$
A^{int, k} \rightarrow \tau _{\leq k}(A).
$$

\subnumero{Comment obtenir une structure de cat\'egorie de mod\`eles ferm\'ee}

Avant d'aborder la preuve du th\'eor\`eme \ref{SeCmf},
nous extrayons de \cite{Jardine}, \cite{JardineGoerssBook},
\cite{Hirschhorn},  et plus sp\'ecialement de Dwyer-Hirschhorn-Kan \cite{DHK} un
\'enonc\'e qui exprime exactement ce qu'il faut d\'emontrer pour obtenir une
structure de cat\'egorie de mod\`eles ferm\'ee comme dans  \ref{SeCmf}. Un
\'enonc\'e de m\^eme nature se trouve dans le livre de Hovey (\cite{HoveyBook}
Theorem 2.1.11).

Soit $M$ une cat\'egorie et $C\subset M$ une sous-cat\'egorie (i.e. classe de
morphismes). Soit encore $I\subset C$ un sous-ensemble de morphismes. On dit
que {\em $I$ engendre $C$} si, pour tout morphisme $f:B\rightarrow B'$ de $M$,
la propri\'et\'e de rel\`evement \cite{Quillen} vis-\`a-vis des morphismes de
$I$ assure la propri\'et\'e de rel\`evement vis-\`a-vis de tous les
morphismes de
$C$  ---i.e.  pour tout diagramme
$$
\begin{array}{ccc}
E&\rightarrow & B \\
\downarrow && \downarrow \\
E'&\rightarrow & B'
\end{array}
$$
avec  $E\rightarrow E'$ dans $C$,
le morphisme diagonal existe.

On va mentionner, sans donner de d\'efinition, la propri\'et\'e pour $I$ de {\em
permettre l'ar\-gu\-ment du petit objet} (``small objet argument''), cf
\cite{Hirschhorn}  \cite{JardineGoerssBook} et particuli\`erement \cite{DHK},
7.3. On note que si $M$ admet un foncteur ``ensemble sous-jacent''
$M\rightarrow Ens$ convenable, alors tout sous-ensemble $I$ de morphismes
permet l'argument du petit objet (\cite{DHK} 8.5, \cite{Hirschhorn}). Ce
sera le cas dans tous les exemples qu'on va consid\'erer. En fait, dans ce
cas on
pourra dire que toute sous-cat\'egorie $C\subset M$ qui est d\'efinie par des
conditions ensemblistement raisonnables, aura un sous-ensemble g\'en\'erateur
$I\subset C$ permettant l'argument du petit objet---voir par exemple
l'argument de Jardine \cite{Jardine}. De fa\c{c}on g\'en\'erale nous ne
donnerons pas de preuve pour ce type de condition.

Expliquons maintenant comment nous allons parler de ``propri\'et\'e de
rel\`evement''.
Il y a deux cas: soit on se donne le morphisme $B\rightarrow B'$ et
on veut avoir la propri\'et\'e de rel\`evement dans le diagramme
ci-dessus pour tout
morphisme $E\rightarrow E'$ dans une certaine classe; soit on se donne
$E\rightarrow E'$ et on veut avoir la propri\'et\'e pour tout $B\rightarrow B'$
dans une certaine classe. Au lieu de faire la distinction en parlant de
``rel\`evement \`a gauche'' et ``rel\`evement \`a droite'' (comme le fait
Quillen
\cite{Quillen}---mais nous n'arrivons pas \`a retenir o\`u est la gauche dans
cette affaire)
nous laisserons au lecteur le soin de d\'eduire du contexte ce qu'on veut dire,
sachant que les ``cofibrations'' sont toujours \`a gauche et les
``fibrations'' \`a droite dans ce genre de diagramme.

On rappelle la notion de {\em cat\'egorie de mod\`eles ferm\'ee engendr\'ee par
cofibrations} de Hirschhorn \cite{Hirschhorn} \cite{DHK}: c'est une cmf $M$,
admettant toutes les limites et colimites, et telle que les
sous-cat\'egories $cof\subset
M$ et $cof \cap W \subset M$ de cofibrations et de cofibrations triviales
admettent des ensembles g\'en\'erateurs permettant l'argument du petit objet.

On donne ici une version du {\em lemme de reconnaissance}
(\cite{DHK} 8.1, qu' Hirschhorn attribue \`a Kan cf \cite{Hirschhorn}). Notre
version n'est pas la plus g\'en\'erale, mais en revanche int\`egre quelques
remarques-cl\'es qui se trouvent seulement plus loin dans \cite{DHK}, au
moment de la d\'emonstration de  la structure de cmf pour les cat\'egories
simpliciales (\cite{DHK} 48.7). Le lemme et sa preuve sont contenus
implicitement
dans l'argument de Jardine \cite{Jardine} et c'est par ailleurs
plus ou moins sous cette forme que le deuxi\`eme auteur a recopi\'e
l'argument de Jardine dans \cite{nCAT}. D'autre part, l'argument se trouve
(\'eparpill\'e) dans \cite{JardineGoerssBook}, ou encore (avec
une liste analogue d'axiomes) dans \cite{JardineGoerssLoc}. Voir \'egalement
l'\'enonc\'e de Hovey \cite{HoveyBook} 2.1.11. Nous conseillons au lecteur qui
veut comprendre cette d\'emonstration, de se reporter aux r\'ef\'erences
\cite{Hirschhorn} \cite{DHK} \cite{Jardine} \cite{JardineGoerssBook}
\cite{JardineGoerssLoc} \cite{HoveyBook}.

Le lemme suivant se trouve donc
enti\`erement dans les r\'ef\'erences ci-dessus.
\begin{lemme}
\label{dhklemme}
Soit $M$ une cat\'egorie avec deux sous-cat\'egories $W\subset M$ et
$cof\subset M$
dites d'equivalences faibles et de cofibrations (et on dit que $cof \cap W$
est la
sous-cat\'egorie des cofibrations triviales). Supposons que:
\newline
(0) \,\, $M$ admet toutes les limites et colimites index\'ees par des petites
cat\'egories;
\newline
(1) \,\, tout r\'etracte d'un morphisme de $W$ (resp. de $cof$)
est dans $W$ (resp. dans $cof$);
\newline
(2)\,\, $W$ satisfait ``trois pour le prix de deux'': si deux parmi $f$, $g$,
$fg$ sont dans $W$ alors le troisi\`eme aussi;
\newline
(3)\,\, tout morphisme qui poss\`ede la propri\'et\'e de
rel\`evement par rapport aux morphismes de
$cof$, est dans $W$;
\newline
(4)\,\, $cof$ admet un sous-ensemble g\'en\'erateur $I\subset cof$ permettant
l'argument du petit objet;
\newline
(5)\,\, de m\^eme, $cof \cap W$ admet un sous-ensemble g\'en\'erateur
$J\subset cof\cap W$
permettant l'argument du petit objet;
\newline
(6)\,\, $cof$ est stable par coproduit (i.e. si $A\rightarrow B$ est
dans $cof$ et $A\rightarrow A'$ est quelconque, alors $A'\rightarrow B\cup ^AA'$
est dans $cof$) et par colimite s\'equentielle transfinie (voir
\cite{Hirschhorn} \cite{DHK} pour la d\'efinition); et
\newline
(7)\,\, de m\^eme, $cof \cap W$ est stable par coproduit et colimite
s\'equentielle
transfinie.

Alors, en notant $fib\subset M$ la sous-cat\'egorie dite ``des fibrations''
\'egale \`a la classe des morphismes qui poss\`edent la propri\'et\'e
de rel\`evement par rapport aux
cofibrations triviales, $(M,W,cof, fib)$ est une cat\'egorie de mod\`eles
ferm\'ee engendr\'ee par cofibrations.
\end{lemme}
{\em Preuve:}
On applique le crit\`ere de reconnaissance (\cite{DHK} 8.1) (dont on adopte les
notations pour la pr\'esente preuve). Le fait que $I$ engendre $cof$ en notre
sens (4), implique $cof \subset I$-${\bf cof}$.  Dans l'autre sens, d'apr\`es
(\cite{DHK} 7.3 (ii)), les morphismes de  $I$-${\bf cof}$ sont retractes de
limites s\'equentielles de coproduits par des morphismes de $I$; donc par nos
conditions (1) pour $cof$ et (6), on a $I$-${\bf cof}\subset cof$ et donc
$
I$-${\bf cof}= cof.
$
De m\^eme en utilisant nos propri\'et\'es (5) et (1) pour $cof$ $+$ (7), on a
$
J$-${\bf cof}= cof \cap W.
$
Ceci nous donne les premi\`eres conditions de \cite{DHK} 8.1 (ii) et (iii),
i.e. que $J$-${\bf cof}=I$-${\bf cof}\cap W$.  Les conditions du d\'ebut de
(\cite{DHK}
8.1) sont notre (1) pour $W$ et notre (2). La condition  (\cite{DHK} 8.1
(i)) r\'esulte de
(4) et (5).  Enfin notre (3) donne $I$-${\bf inj} \subset W$,
et comme le rel\`evement pour $I$ (resp. $J$) est \'equivalent au rel\`evement
pour $cof$ (resp. $cof \cap W$), il est clair qu'on a
$I$-${\bf inj} \subset J$-${\bf inj}$. Ceci donne la deuxi\`eme partie de
(\cite{DHK} 8.1 (ii)), \`a savoir $I$-${\bf inj} \subset J$-${\bf inj} \cap W$.
\eop

{\em Remarque:} Il r\'esulte imm\'ediatement de la d\'efinition de
\cite{Hirschhorn}
\cite{DHK} que si $(M, W, cof)$ est une cat\'egorie de mod\`eles
ferm\'ee engendr\'ee par cofibrations, alors les propri\'et\'es (0)--(7) sont
satisfaites.
Il peut \^etre utile de noter que la partie de la  d\'efinition de cmf
engendr\'ee par cofibrations
qui va au-del\`a
des axiomes de  Quillen, consiste en: (0) pour les limites et colimites
infinies; (4) et (5); et
(6) et (7) pour les limites s\'equentielles.

{\em Remarque:} Lorsque nous appliquerons ce lemme, nous ne v\'erifierons
pas les
conditions (4) et (5), qui seront cons\'equences de l'existence d'un
foncteur ``ensemble sous-jacent'' et du fait que $cof$ et $W$ seront
(ensemblistement)
raisonnables.
Voir dans Jardine \cite{Jardine} la technique n\'ecessaire pour obtenir une
d\'emonstration
rigoureuse de (4) et (5) dans notre type de situation (il traite le cas des
pr\'efaisceaux
simpliciaux mais sa technique est tr\`es g\'en\'erale).
D'autre part, les cofibrations seront en
g\'en\'eral simplement les morphismes qui induisent des injections sur les
ensembles sous-jacents; donc (1) pour $cof$ et (6) seront \'evidents.
L'existence de limites et colimites sera une cons\'equence imm\'ediate de la
d\'efinition de la cat\'egorie $M$ (g\'en\'eralement comme cat\'egorie de
pr\'efaisceaux d'ensembles sur une certaine cat\'egorie).  En somme, nous
v\'erifierons seulement (1) pour $W$, (2), (3) et (7). La principale
difficult\'e sera
constitu\'ee par (7).

\subnumero{D\'emonstration du th\'eor\`eme \ref{SeCmf}}
On pourrait simplement recopier la d\'emonstration de \cite{nCAT}
en rempla\c{c}ant ``$n$-pr\'ecat'' par ``$n$-pr\'ecat de Segal''. Nous allons
indiquer cependant comment obtenir la d\'emonstration en utilisant le lemme
\ref{dhklemme} et en utilisant directement une partie des r\'esultats de
\cite{nCAT}, sans tout recopier. Au passage, la pr\'esente
d\'emonstration am\'eliore celle du \S 6 de \cite{nCAT} en rempla\c{c}ant
l'argument {\em ad hoc} pour la condition  CM5(2) par la condition (3) du lemme
\ref{dhklemme} (plus, implicitement, via ce lemme, l'argument du petit objet).

Comme dans la remarque pr\'ec\'edente, les conditions (0), (1) pour $cof$,
et (6) sont imm\'ediates. Nous n'allons pas traiter les conditions (4) et (5)
qui sont de nature ensembliste. Il reste donc \`a traiter (1) pour
$W$, (2), (3), et (7). Fixons $n$ et traitons la cat\'egorie $nSePC$ des
$n$-pr\'ecats de Segal.

Soit $f:A\rightarrow B$ un morphisme de $n$-pr\'ecats de Segal.
Supposons que $f$ est r\'etracte d'une \'equivalence faible $g$.  Alors
$SeCat(f)$ est r\'etracte de $SeCat(g)$ et pour tout $m$, $\Pi _m \circ
SeCat(f)$ est r\'etracte de $\Pi _m \circ SeCat(g)$. Ces derniers sont
des morphismes de $n+m$-cat\'egories, donc la fermeture sous r\'etractes de la
classe des \'equivalences entre $n+m$-cat\'egories  (\cite{nCAT} \S 6) implique
que $\Pi _m \circ SeCat(f)$ est une \'equivalence. Ce fait \'etant vrai
pour tout $m$, il
s'ensuit que $SeCat(f)$ est une \'equivalence donc, par d\'efinition, $f$ est
une \'equivalence faible. Ceci prouve (1) pour $W$.

De la m\^eme fa\c{c}on la propri\'et\'e (2) ``trois pour le prix de deux''
est une cons\'equence imm\'ediate (via la construction $\Pi _m\circ SeCat$)
de la m\^eme propri\'et\'e pour les $n+m$-cat\'egories pour tout $m$
(\cite{nCAT} Lemma 3.8).

Pour la propri\'et\'e (3),
si $n=0$, la propri\'et\'e (3) est bien connue pour
les ensembles simpliciaux; on suppose donc $n\geq 1$, et on suppose que (3)
est connu pour $n-1$.  Supposons qu'un morphisme $f:A\rightarrow B$
dans $nSePC$ satisfait la propri\'et\'e de rel\`evement par rapport \`a toute
cofibration $E\hookrightarrow E'$. Avec la cofibration $\emptyset
\hookrightarrow \ast$ on
obtient que $f$ est surjectif sur les objets. On rappelle maintenant une
construction de \cite{limits} qui sera utilis\'ee souvent au travers du
pr\'esent papier. Si $E$ est une $n-1$-pr\'ecat on d\'efinit la $n$-pr\'ecat
$\Upsilon (E)$ avec deux objets $\Upsilon (E)_0 = \{ 0,1\}$ et avec
$$
\Upsilon (E)_{p/}(x_0,\ldots , x_p)
$$
\'egal \`a $E$
si $x_0=\ldots = x_i=0$ et $x_{i+1}=\ldots = x_p = 1$ ($0\leq i<p$); \'egal
\`a $\ast$
si $x_0=\ldots = x_p=0$ ou $x_0=\ldots = x_p=1$; et \'egal \`a $\emptyset$
sinon. Un morphisme $\Upsilon (E)\rightarrow A$
 n'est rien d'autre
qu'une paire d'objets $(x,y)\in A_0\times A_0$ munie d'un morphisme
$E\rightarrow
A_{1/}(x,y)$. On revient \`a la d\'emonstration de la propri\'et\'e (3). La
propri\'et\'e de rel\`evement pour $f$ vis-\`a-vis de toute cofibration,
s'applique en particulier \`a $\Upsilon (E)\rightarrow \Upsilon (E')$ pour
toute cofibration $E\rightarrow E'$ de $n-1$-pr\'ecats de Segal. Ceci implique
que si $x,y\in A_0$ alors le morphisme
$$
A_{1/}(x,y)\rightarrow B_{1/}(x,y)
$$
de $n-1$-pr\'ecats de Segal, a la propri\'et\'e de rel\`evement pour toute
cofibration de $n-1$-pr\'ecats de Segal. La propri\'et\'e (3) en degr\'e $n-1$
implique que $A_{1/} (x,y)\rightarrow B_{1/}(x,y)$ est une \'equivalence
faible. Le m\^eme type d'argument (en utilisant la variante de
la construction $\Upsilon (E)$ not\'ee $[p]E$
dans \cite{limits} 2.4.5) montre que
$$
A_{p/}(x_0,\ldots , x_p)\rightarrow B_{p/}(f(x_0),\ldots , f(x_p))
$$
est une
\'equivalence faible. Il r\'esulte clairement de la description de $SeCat$ comme
compos\'ee d'op\'erations de la forme $Fix$ et $Gen[m]$ (voir les remarques
suivant le lemme \ref{secatunique}), que $SeCat(A)_{1/}(x,y)\rightarrow
SeCat(B)_{1/}(f(x), f(y))$ est une \'equivalence. Donc $SeCat(f)$ est pleinement
fid\`ele, donc c'est une \'equivalence, et $f$ est une \'equivalence faible.
Ceci donne la propri\'et\'e (3).

Pour la propri\'et\'e (7), soit
$$
B\leftarrow     A\rightarrow C
$$
un diagramme o\`u $A\rightarrow B$ est une cofibration triviale. Comme on l'a
remarqu\'e ci-dessus,  le th\'eor\`eme de Seifert-Van Kampen, Theorem 9.1 de
\cite{nCAT}, donne une \'equivalence
$$
(\Pi _m\circ B) \cup ^{(\Pi _m \circ A)}(\Pi _m \circ C)
\stackrel{\cong}{\rightarrow} \Pi _m\circ (B\cup ^AC).
$$
La propri\'et\'e de pr\'eservation des cofibrations triviales par coproduit
dans le cadre
des $n+m$-pr\'ecats non de Segal (\cite{nCAT} Lemma 3.2) implique donc que
$$
\Pi _m \circ C\rightarrow
\Pi _m \circ (B\cup ^AC)
$$
est une \'equivalence faible. Ce fait \'etant vrai pour tout $m$, il s'ensuit
que $C\rightarrow
B\cup ^AC$ est une \'equivalence faible (il convient de rappeler ici que par
convention nous avons pr\'ecompos\'e l'op\'eration $\Pi_m$ par $SeCat$). Ceci
donne (7).
\eop

\subnumero{Classes d'homotopie de morphismes}

On rappelle (Quillen \cite{Quillen}) que si $M$ est une cat\'egorie de mod\`eles
ferm\'ee et si $W\subset M$ est la sous-cat\'egorie dont les morphismes sont
les \'equivalences faibles, alors la cat\'egorie localis\'ee
$$
Ho(M):= W^{-1}M
$$
au sens de Gabriel-Zisman \cite{GabrielZisman} poss\`ede d'autres descriptions
utiles. Soient $M_c$, $M_f$ et $M_{cf}$ les cat\'egories d'objets cofibrants,
fibrants et cofibrants et fibrants, respectivement. Soient $W_c$, $W_f$ et
$W_{cf}$ les intersections de $W$ avec celles-ci. Alors les localis\'es
$W_c^{-1} M_c$, $W_f^{-1}M_f$ et $W_{cf}^{-1}M_{cf}$ sont toutes trois
\'equivalentes (via
les foncteurs induits par les inclusions) \`a $Ho(M)$.
Si on note $\pi M_{cf}$ la cat\'egorie dont les objets sont les objets de
$M_{cf}$ et les ensembles de morphismes sont les ensembles de {\em classes
d'homotopie} de morphismes de $M_{cf}$ \cite{Quillen}, alors $\pi M_{cf}$
est aussi
\'equivalente \`a $Ho(M)$. Plus g\'en\'eralement, pour $A,B\in M$, si on choisit
des \'equivalences faibles $A'\cong A$ et $B\cong B'$ avec $A'$ cofibrant et
$B'$ fibrant, alors $Hom _{Ho(M)}(A,B)$ est l'ensemble de classes d'homotopie
de morphismes $A'\rightarrow B'$. Suivant Quillen \cite{Quillen}, on va
adopter la notation
$$
[A,B] := Hom _{Ho(M)}(A,B)
$$
pour l'ensemble des morphismes de $A$ vers $B$ dans $Ho(M)$, avec $A,B\in
Ob(M)$.
En r\'esum\'e, si $A$ est cofibrant et $B$ fibrant alors tout \'el\'ement de
$[A,B]$ provient d'un morphisme $A\rightarrow B$ de $M$, et deux tels morphismes
donnent le m\^eme \'el\'ement de $[A,B]$ si et seulement s'ils sont homotopes
au sens de \cite{Quillen}.

On applique ceci avec $M= nSePC$. Dans ce cas, on a un objet
$\overline{I}$, la cat\'egorie avec deux objets $0,1$ li\'es par un
isomorphisme, consid\'er\'e comme $n$-cat\'egorie de Segal (voir la rubrique
``Induction'' ci-dessus). Cet objet sert d'{\em intervalle}, i.e. pour tout
objet
$A\in nSePC$, $A\times \overline{I}$ est un objet ``$A\times I$'' au sens de
Quillen \cite{Quillen} et pour $B$ fibrant, deux morphismes $f,g: A\rightarrow
B$ sont homotopes si et seulement s'il existe un morphisme $h: A\times
\overline{I}\rightarrow B$ avec $h|_{A\times \{ 0\} }=f$ et $h|_{A\times \{ 1\}
}=g$. Si $A$ et $B$ sont deux objets de $nSePC$ alors on peut calculer
$$
[A,B]= Hom _{Ho(nSePC)}(A,B)
$$
en choisissant une \'equivalence
faible $B\rightarrow B'$ vers un objet fibrant et (gr\^ace au fait que $A$ est
automatiquement cofibrant) en prenant l'ensemble des morphismes $A\rightarrow
B'$ modulo la relation d'homotopie explicit\'ee plus haut \`a l'aide de
$\overline{I}$.
Voir Cordier \cite{Cordier} qui fait cette observation dans le cadre des 
diagrammes homotopie-coh\'erents.

\numero{Les $n$-pr\'echamps de Segal}

\label{nchampsegpage}

Soit maintenant $\Xx$ un site, i.e. une cat\'egorie munie d'une topologie de
Grothendieck. Rappelons que ceci veut dire une collection de cribles $\Bb
\subset \Xx /X$, avec $X$ variable,  v\'erifiant certains axiomes
(voir \cite{SGA4} \cite{MacLaneMoerdijk} \cite{Rezk}).

Un {\em $n$-pr\'echamp de Segal au-dessus de $\Xx$} est un foncteur de
$\Xx ^o$ dans la cat\'egorie $nSePC$ des $n$-pr\'ecats de Segal, autrement dit
c'est un pr\'efaisceau de $n$-pr\'ecats de Segal sur $\Xx$. Par exemple pour
$n=0$ c'est simplement un
pr\'efaisceau d'ensembles simpliciaux sur $\Xx$.

On note $nSePCh(\Xx )$ ou simplement $nSePCh$ la cat\'egorie des
$n$-pr\'echamps de Segal au-dessus de $\Xx$ avec les morphismes respectant
strictement la structure (i.e. les morphismes de pr\'efaisceaux \`a valeurs dans
$nSePC$).

Insistons sur le fait qu'il s'agit de foncteurs stricts de $\Xx ^o$ vers
la $1$-cat\'egorie $nSePC$ o\`u les morphismes sont les morphismes de
$n$-pr\'ecats de Segal respectant strictement la structure; en particulier
on {\em
ne consid\`ere pas pour le moment} de ``pr\'efaisceaux faibles'' ou
``l\^aches''.
Ceci sera justifi\'e ult\'erieurement par le th\'eor\`eme \ref{correlation}
et en
particulier sa partie ``strictification'' qui est, en fait, le corollaire
\ref{strictif3}.

On peut appliquer l'op\'eration $\tau _{\leq m}$ \`a un $n$-pr\'echamp de
Segal $A$:
d'abord
on remplace chaque $A(X)$ par la $n$-cat\'egorie de Segal $SeCat(A(X))$ et
ensuite on applique objet-par-objet la troncation $\tau _{\leq m}$ d\'ecrite
dans la section pr\'ec\'edente.

Pour $n=0$ on dispose de la cat\'egorie de mod\`eles ferm\'ee de Jardine
\cite{Jardine} des pr\'e\-faisceaux simpliciaux: les cofibrations sont les
injections, les \'equivalences faibles sont celles d'Illusie \cite{Illusie},
et les fibrations sont d\'etermin\'ees par la propri\'et\'e de rel\`evement.
Jardine montre dans \cite{Jardine} que ces trois classes donnent une
structure de
cat\'egorie de mod\`eles ferm\'ee.

Pour $n>0$ on va
suivre le fil de l'argument de Jardine \cite{Jardine} (qui est en quelque sorte
r\'esum\'e  dans le lemme \ref{dhklemme}).
Pour commencer, on impose---comme Jardine---que les {\em cofibrations} de
$n$-pr\'echamps de Segal soient les injections,
i.e. les morphismes $A\rightarrow B$ tels que $A(X)\rightarrow B(X)$ soit une
cofibration pour tout $X\in \Xx$.

La diff\'erence principale avec \cite{Jardine}
est qu'il faut donner une d\'efinition r\'ecursive (sur $n>0$) de
l'\'equivalence
faible. On dira qu'un morphisme $A\rightarrow B$ de $n$-pr\'echamps de Segal,
o\`u tous
les $A(X)$ et $B(X)$ sont des $n$-cat\'egories de Segal, est
une {\em \'equivalence faible} si:
\newline
---pour tout $X\in \Xx$ et tout $x,y\in A_0(X)$ le morphisme
$$
A_{1/}(x,y)\rightarrow B_{1/}(x,y)
$$
est une \'equivalence faible de $n-1$-pr\'echamps de Segal au-dessus de
$\Xx /X$;
et
\newline
---le morphisme de pr\'efaisceaux d'ensembles
$\tau _{\leq 0} A\rightarrow \tau _{\leq 0}B$ sur $\Xx$
induit une surjection sur les faisceaux associ\'es.

On dira qu'un morphisme $A\rightarrow B$ de $n$-pr\'echamps de Segal est
une {\em
\'equivalence faible} si le morphisme $SeCat(A)\rightarrow SeCat(B)$ est une
\'equivalence faible au sens de la d\'efinition ci-dessus.

Rappelons que pour $n=0$ on a pris l'\'equivalence d'Illusie comme notion
d'\'equivalence
faible (les d\'efinitions ci-dessus n'ayant pas vraiment de sens).
On peut observer que ce choix est compatible avec la construction $\Pi
_{m,Se}\circ A$ et l'op\'eration de troncation, voir le lemme \ref{compatible}
et le corollaire \ref{compPin} ci-dessous.

Un morphisme de $n$-pr\'echamps de Segal est une {\em fibration} s'il poss\`ede
la  propri\'et\'e de rel\`evement vis-\`a-vis
des cofibrations triviales.

{\em Comparaison avec la topologie grossi\`ere:}
Notons $\Gg$ la topologie de Grothendieck donn\'ee sur $\Xx$ pour
la distinguer de la
topologie grossi\`ere.
On rappel que la topologie grossi\`ere est celle o\`u la seule crible
recouvrant $X\in \Xx$ est
$\Bb = X/\Xx$ (autrement dit pour qu'une famille recouvre $X$ il faut
qu'elle contient $X$).
On utilisera souvent cette topologie donc il convient de comparer les
notions principales pour
$\Gg$ et pour la topologie grossi\`ere.  La classe des cofibrations ne
d\'epend pas de
$\Gg$. Soit $f:A\rightarrow
B$ un morphisme de $n$-pr\'echamps de Segal. Alors:
\newline
---$f$ est une
\'equivalence faible pour la topologie grossi\`ere si et seulement si pour tout
$X\in \Xx$, $A(X)\rightarrow B(X)$ est une \'equivalence faible de $n$-pr\'ecats
de Segal;
\newline
---si $f$ est une \'equivalence faible pour la topologie grossi\`ere alors $f$
est une $\Gg$-\'equivalence faible;
\newline
---si $f$ est une $\Gg$-fibration alors $f$ est une fibration pour la topologie
grossi\`ere;
\newline
---si $f$ est une fibration pour la topologie grossi\`ere alors chaque $f_X:
A(X)\rightarrow B(X)$ est une fibration de $n$-pr\'ecats de Segal;
\newline
---(ici on utilise le Th\'eor\`eme \ref{cmf} ci-dessous) $f$ est une
$\Gg$-fibration
triviale si et seulement si $f$ est une fibration triviale pour la topologie
grossi\`ere.

On revient maintenant \`a la consid\'eration de notre site $\Xx$ avec sa
topologie $\Gg$.

\begin{theoreme}
\label{cmf}
La cat\'egorie $nSePCh(\Xx )$ des $n$-pr\'echamps de Segal avec les trois
classes de morphismes
(cofibrations, $\Gg$-\'equivalences et $\Gg$-fibrations) est une
cat\'egorie de mod\`eles ferm\'ee.
\end{theoreme}

Pour $n=0$ c'est le th\'eor\`eme de
Jardine \cite{Jardine}, Joyal \cite{Joyal}, K.
Brown \cite{KBrown}. On
peut donc supposer $n\geq 1$ et supposer par r\'ecurrence
que le r\'esultat est vrai pour $n-1$.

Pour la d\'emonstration du Th\'eor\`eme \ref{cmf} on aura besoin de
la notion suivante. Soit $f: A\rightarrow B$ un morphisme de
$n$-pr\'echamps de
Segal. On dira que $f$ est une {\em \'equivalence pr\'eliminaire} si:
\newline
---le morphisme de pr\'efaisceaux d'ensembles $A_0\rightarrow B_0$ induit un
isomorphisme entre les faisceaux associ\'es; et
\newline
---pour tout $X\in \Xx$ et pour chaque suite d'objets $x_0,\ldots , x_m\in
A_0(X)$
le morphisme
$$
A_{m/} (x_0,\ldots , x_m)\rightarrow B_{m/}(f(x_0),\ldots ,f(x_m))
$$
est une \'equivalence faible de $n-1$-pr\'echamps de Segal sur $\Xx /X$.

\begin{lemme}
Soit
$$
\begin{array}{ccc}
A & \stackrel{p}{\rightarrow} & B \\
\downarrow && \downarrow \\
A' & \stackrel{q}{\rightarrow} &B'
\end{array}
$$
un diagramme commutatif de morphismes de $n-1$-pr\'echamps de Segal, o\`u $B$ et
$B'$ sont des pr\'efaisceaux d'ensembles. Supposons que le morphisme vertical de
droite (qu'on note $g$) induit un isomorphisme entre les faisceaux
associ\'es. Alors
le morphisme $A\rightarrow A'$ est une
\'equivalence faible de $n-1$-pr\'echamps de Segal si et seulement si
pour tout $X\in \Xx$ et $x\in B(X)$ le morphisme
$p^{-1}(x)\rightarrow q^{-1}(g(x))$ est une \'equivalence faible de
$n-1$-pr\'echamps de Segal.
\end{lemme}
C'est une cons\'equence directe de la d\'efinition de l'\'equivalence faible.
\eop

\begin{corollaire}
\label{etape}
Un morphisme $A\rightarrow B$ de $n$-pr\'echamps de Segal est une \'equivalence
pr\'eliminaire si et seulement si pour tout $m\in \Delta$ le morphisme
$A_{m/} \rightarrow B_{m/}$ de $n-1$-pr\'echamps de Segal est une \'equivalence
faible. Dans ce cas les morphismes
$$
A_{1/} \times _{A_0}\ldots \times _{A_0} A_{1/}\rightarrow
B_{1/} \times _{B_0}\ldots \times _{B_0} B_{1/}
$$
sont aussi des \'equivalences faibles.
\end{corollaire}
On note que pour $m=0$ c'est la premi\`ere condition pour une
\'equivalence pr\'eliminaire;  pour la deuxi\`eme condition on applique le lemme
pr\'ec\'edent.
\eop

\begin{lemme}
\label{etape2}
Soit $f: A\rightarrow B$ une \'equivalence pr\'eliminaire de $n$-pr\'echamps de
Segal.
Alors $f$ est une \'equivalence faible de $n$-pr\'echamps de Segal.
\end{lemme}
{\em Preuve:}
Le cas o\`u les $A(X)$ et $B(X)$ sont des
$n$-cat\'egories de Segal est imm\'ediat. Il
suffit donc de prouver que si $f$ est une \'equivalence pr\'eliminaire alors
il en est de m\^eme de $SeCat(f): SeCat(A)\rightarrow SeCat(B)$. Pour ceci on
analysera les \'etapes d\'ecrites dans \cite{nCAT} pour passer de $A$ \`a
$SeCat(A)$ (la discussion de \cite{nCAT} pour l'op\`eration $Cat$ s'adapte
directement pour l'op\`eration $SeCat$, comme on l'a rappel\'e au \S 2).
En outre,  il s'agit ici
de pr\'efaisceaux de $n$-pr\'ecats de Segal sur $\Xx$;
l'op\'eration $SeCat$ est faite objet-par-objet (elle est fonctorielle). La
construction $SeCat$ est obtenue en it\'erant les constructions not\'ees
$Fix$ et
$Gen[m]$ dans \cite{nCAT} et au \S 2 ci-dessus.
L'op\'eration $Fix$ pr\'eserve le
type d'\'equivalence
faible des $A_{m/}$ donc il est clair qu'elle pr\'eserve la condition
d'\^etre une
\'equivalence pr\'eliminaire. Il s'agit donc de traiter de $Gen[m]$.

L'op\'eration $Gen[m]$ comporte d'abord une \'etape not\'ee $A\mapsto A'$ dans
\cite{nCAT}, qui pr\'eserve le type d'\'equivalence faible des $A_{p/}$. Nous
pouvons donc ignorer cette op\'eration et noter encore  $A$ (resp. $B$) son
r\'esultat.

Ensuite on consid\`ere un diagramme de la forme
$$
A_{m/} \rightarrow \Gg [m](A) \stackrel{g}{\rightarrow }
A_{1/}\times _{A_0} \ldots \times _{A_0} A_{1/}
$$
o\`u $g$ est une cofibration triviale de $n$-pr\'echamps de Segal
objet-par-objet (i.e. pour la topologie grossi\`ere). (NB la notation
$\Gg [m]$ n'a pas de lien avec la notation $\Gg$ pour la topologie de
Grothendieck, on garde $\Gg [m]$ seulement pour respecter les
notations de \cite{nCAT}.) Maintenant $Gen[m](A)_{p/}$ est le coproduit de
$A_{p/}$ et $A_{m/}\rightarrow \Gg [m](A)$ pour divers morphismes
$A_{m/}\rightarrow A_{p/}$ induits par divers morphismes $p\rightarrow
m$.
La stabilit\'e des cofibrations $\Gg$-triviales par
coproduit avec des objets cofibrants, pour les $n-1$-pr\'echamps de Segal,
est une
cons\'equence du Th\'eor\`eme \ref{cmf} pour $n-1$, via le lemme de Reedy
\cite{Reedy}.
Donc le fait que les morphismes
$$
\Gg [m](A)\rightarrow \Gg [m](B)
$$
et
$$
A_{p/} \rightarrow B_{p/}
$$
soient des \'equivalences faibles de $n-1$-pr\'echamps de Segal implique
qu'il en
est de m\^eme de  $Gen[m](A)_{p/}\rightarrow Gen[m](B)_{p/}$. On obtient par le
corollaire \ref{etape} que $Gen[m](A)\rightarrow Gen[m](B)$ est une
\'equivalence pr\'eliminaire. En observant que l'\'equivalence faible des
$n-1$-pr\'echamps de Segal et donc l'\'equivalence pr\'eliminaire des
$n$-pr\'echamps
de Segal est pr\'eserv\'ee par  colimite filtrante,
on obtient que
$SeCat(A)\rightarrow SeCat(B)$ est une \'equivalence pr\'eliminaire,
donc une
\'equivalence faible.
\eop

\begin{lemme}
\label{locale}
La propri\'et\'e d'\^etre une \'equivalence faible est locale: si
$f:A\rightarrow
B$ est un morphisme de $n$-pr\'echamps sur $\Xx /X$ et si $\Bb \subset
\Xx /X$
est un crible couvrant $X$
alors $f$ est une \'equivalence faible si et seulement si pour tout $Y\in
\Bb$, $f|_{\Xx /Y}$ est une \'equivalence faible.
\end{lemme}
{\em Preuve:}
On raisonne par r\'ecurrence sur $n$. Pour $n=0$ c'est une cons\'equence
directe de la d\'efinition d'Illusie. Si on suppose que c'est vrai pour $n-1$
alors il suffit (au vu de la d\'efinition d'\'equivalence faible) d'observer
que le fait que $\tau _{\leq 0}SeCat(A)\rightarrow \tau
_{\leq 0} SeCat(B)$ induise
une surjection sur les faisceaux associ\'es est une propri\'et\'e locale;
ce qui est \'evident
(on utilise ici les axiomes d'une topologie de Grothendieck).
\eop

{\em D\'emonstration du th\'eor\`eme \ref{cmf}:}
Notre preuve est mod\'el\'ee sur l'argument de Jardine \cite{Jardine}
et nous appliquons donc
notre lemme \ref{dhklemme}.

La propri\'et\'e (0) est imm\'ediate.

Il est facile de voir que les \'equivalences faibles sont stables sous
r\'etractes, par r\'e\-cur\-ren\-ce sur $n$ (le cas $n=0$ est
cons\'equence directe de la d\'efinition d'\'equivalence faible d'Illusie).
La stabilit\'e des cofibrations sous r\'etractes est imm\'ediate, ce qui donne
la condition (1) de
\ref{dhklemme}.

Pour la propri\'et\'e ``trois pour le prix de deux'' ((2) du \ref{dhklemme}),
soient $f: A\rightarrow B$ et $g: B\rightarrow C$ deux morphismes. Il est
imm\'ediat que si $f$ et $g$ sont
des \'equivalences faibles alors $gf$ aussi; et que si $gf$ et $g$ sont des
\'equivalences faibles alors $f$ aussi (on remarque que la premi\`ere
condition pour $g$ entraine que le morphisme $\tau _{\leq 0}B\rightarrow \tau
_{\leq 0}C$ induit une injection entre les faisceaux associ\'es, donc la
surjectivit\'e essentielle de $gf$ implique celle de $f$). Il faut montrer
que si $gf$
et $f$
sont des \'equivalences faibles alors $g$ est une \'equivalence faible.
La condition d'essentielle surjectivit\'e est \'evidente; le probl\`eme est la
premi\`ere condition (voir la d\'emonstration analogue dans \cite{nCAT}, Lemma
3.8). Cette condition est facile \`a v\'erifier pour des objets
$x,y\in B_0(X)$ de la forme
$x=f(u)$ et $y=f(v)$ avec $u,v\in A_0(X)$. Pour le cas g\'en\'eral, soient
$x,y\in
B_0(X)$. Alors pour tout $Y$ dans une famille couvrante de $X$ on peut
trouver des
objets $u_Y$ et $v_Y$ dans $A_0(Y)$ avec des \'equivalences
$$
i: f(u_Y)\cong x|_Y, \;\;\; j: f(v_Y)\cong y|_Y
$$
(ce sont des $1$-morphismes dans $B(Y)$). On peut dire que ``la composition avec
$i$ et $j$'' induit une \'equivalence objet-par-objet entre
$B_{1/}(f(u_Y), f(v_Y))|_{\Xx /Y}$ et
$B_{1/}(x,y)|_{\Xx /Y}$ (pour pr\'eciser ceci il faut une discussion
analogue \`a celle pr\'ec\'edant le   Lemma 3.8 dans \cite{nCAT}). Le fait que
$$
B_{1/}(f(u_Y), f(v_Y))|_{\Xx /Y}\stackrel{\cong}{\rightarrow}
C_{1/}(gf(u_Y), gf(v_Y))|_{\Xx /Y}
$$
soit une \'equivalence
implique que
$$
B_{1/}(x, y)|_{\Xx /Y}\rightarrow
C_{1/}(g(x), g(y))|_{\Xx /Y}
$$
est une \'equivalence. Ceci pour tout $Y$ dans la famille couvrant $X$.
Comme la condition d'\^etre une \'equivalence faible est locale
(lemme \ref{locale})
on obtient que
$$
B_{1/}(x, y)|_{\Xx /X}\rightarrow
C_{1/}(g(x), g(y))|_{\Xx /X}
$$
est une \'equivalence, ce qui est la condition cherch\'ee.

Si $B'\rightarrow B$ est un morphisme qui poss\`ede la propri\'et\'e de
rel\`evement pour toute cofibration $A\rightarrow A'$, alors on obtient en
particulier la propri\'et\'e de rel\`evement vis-\`a-vis des cofibrations
de la forme
$C_X\rightarrow C'_X$ o\`u $C_X$ est le $n$-pr\'echamp de Segal engendr\'e
librement par une $n$-pr\'ecat de Segal $C$ au-dessus d'un objet $X\in \Xx$;
et ce morphisme provient par la m\^eme construction d'une cofibration
$C\rightarrow C'$. Un morphisme  $C_X\rightarrow B'$ s'identifie \`a un
morphisme $C\rightarrow B'(X)$, et donc pour tout $X\in \Xx$, $B'(X)\rightarrow
B(X)$ poss\`ede la propri\'et\'e de rel\`evement vis-\`a-vis des cofibrations
de $n$-pr\'ecats de Segal. Ceci implique (par \ref{SeCmf}) que $B'(X)\rightarrow
B(X)$ est une \'equivalence faible, donc $B'\rightarrow B$ est une
\'equivalence faible objet-par-objet i.e. pour la topologie grossi\`ere. On
obtient ainsi la condition (3) de \ref{dhklemme}.

Comme d'habitude nous laissons au lecteur les questions ensemblistes
(4) et (5) de \ref{dhklemme}.

La stabilit\'e des cofibrations par
coproduit ((6) du \ref{dhklemme}) est imm\'ediate.
Il ne reste donc qu'\`a d\'emontrer la stabilit\'e des cofibrations
triviales par coproduit, la condition (7) de \ref{dhklemme}, ce qui
occupera le reste de la
d\'emonstration (nous laissons au lecteur le soin de v\'erifier la
stabilit\'e des
cofibrations triviales par limite s\'equentielle).

Soient $A\rightarrow B$ une cofibration $\Gg$-triviale
et $A\rightarrow C$ un morphisme de
$n$-pr\'echamps de Segal. Soit $P:= B\cup ^AC$. On veut montrer que
$C\rightarrow P$ est une cofibration $\Gg$-triviale.  Il suffit de montrer
que $SeCat(C)\rightarrow SeCat(P)$ est une $\Gg$-\'equivalence faible, or on
sait (Th\'eor\`eme \ref{SeCmf}) que le morphisme
$$
SeCat(B)\cup ^{SeCat(A)}SeCat(C) \rightarrow SeCat(P)
$$
est une \'equivalence faible objet-par-objet. Le morphisme $SeCat(A)\rightarrow
SeCat(B)$ est encore une cofibration $\Gg$-triviale. Quitte \`a remplacer
$A$, $B$ et $C$ par $SeCat(A)$, $SeCat(B)$ et $SeCat(C)$ respectivement, on peut
donc supposer que, pour tout $X$,  $A(X)$, $B(X)$ et $C(X)$ sont des
$n$-cat\'egories de
Segal.

On traite d'abord quelques cas particuliers.

\noindent
{\em Cas 1.}\,  On suppose que $A_0\rightarrow B_0$ induit un isomorphisme
sur les
faisceaux associ\'es. Dans ce cas, il est facile de voir que $A\rightarrow
B$ est
une \'equivalence pr\'eliminaire. Le r\'esultat qu'on est en train de
d\'emontrer,
appliqu\'e aux $n-1$-pr\'echamps de Segal et combin\'e avec le Corollaire
\ref{etape},
implique que $C\rightarrow B\cup ^AC$ est une \'equivalence pr\'eliminaire.
Le Lemme \ref{etape2} implique alors que c'est une \'equivalence faible.

\noindent
{\em Cas 2.}\, On suppose que $A\rightarrow B$ est une \'equivalence
faible et que pour tout
$X\in \Xx$ le morphisme $A(X)\rightarrow B(X)$ est essentiellement
surjectif.
On commence par remarquer que le pr\'esent Th\'eor\`eme \ref{cmf} est
imm\'ediat
dans le cas de la topologie
grossi\`ere.
Il en r\'esulte qu'il existe une factorisation
$$
A\rightarrow B' \rightarrow B
$$
dans laquelle, pour la topologie grossi\`ere,
le premier morphisme est une cofibration et le deuxi\`eme
est une fibration triviale.
Soit $B'' \subset B'$ d\'efini par la condition que pour tout
$X\in \Xx$, $B''(X) \subset B'(X)$ est la sous-$n$-cat\'egorie de Segal pleine
ayant pour objets ceux qui sont dans l'image de $A_0$. Le morphisme
$B'' (X)\rightarrow B(X)$ est encore essentiellement surjectif et
pleinement fid\`ele (car $B'(X)\rightarrow B(X)$ est une
\'equivalence de $n$-cat\'egories de Segal). Donc $B'' \rightarrow B$ est
une \'equivalence pour la topologie grossi\`ere.
Notre morphisme se factorise
\`a travers une cofibration $A\rightarrow B''$.
Le lemme de Reedy \cite{Reedy}
(cf \cite{Hirschhorn} \cite{DHK} \cite{HoveyBook} \cite{JardineGoerssBook}),
qui s'applique ici car tous les objets sont cofibrants, implique que le
morphisme
$$
B'' \cup ^A C \rightarrow B\cup ^AC
$$
est une \'equivalence pour la topologie grossi\`ere. Il suffit alors de savoir
que
$$
C \rightarrow B'' \cup ^A C
$$
est une \'equivalence faible. Le morphisme $A_0\rightarrow B''_0$ est un
isomorphisme (et $A\rightarrow B''$ est une \'equivalence faible car $B''$ est
objet-par-objet \'equivalent \`a $B$), donc notre Cas 1 s'applique.

\medskip

\noindent
{\em Fin de la d\'emonstration:}\,
On traite maintenant le cas g\'en\'eral (avec les m\^emes notations
$B\leftarrow A\rightarrow C$). Soit $U_0\subset B_0$ l'image de $A_0$. Soit
$U\subset B$ le sous-$n$-pr\'echamp de Segal plein ayant $U_0$ comme
ensemble d'objets. Par ailleurs,
pour $X\in \Xx$, soit $V_0(X)$ le sous-ensemble des objets de $B_0(X)$
qui sont \'equivalents \`a un objet de $U_0(X)$ et soit $V\subset B$ le
sous-$n$-pr\'echamp de Segal plein correspondant. Objet-par-objet, le morphisme
$U\rightarrow V$ est essentiellement surjectif et pleinement
fid\`ele, donc notre Cas 2 s'applique.

Le morphisme $A\rightarrow U$ est une \'equivalence faible et il est bijectif
sur les objets; donc le Cas 1 s'applique.

Enfin le morphisme de pr\'efaisceaux d'ensembles $V_0\rightarrow B_0$ induit un
isomorphisme entre les faisceaux associ\'es: l'injectivit\'e d\'ecoule du
fait que le
morphisme de pr\'efaisceaux est injectif; et la surjectivit\'e d\'ecoule du
fait que
tout objet de $B_0(X)$ est localement \'equivalent \`a un objet qui provient de
$A_0$, donc aussi de $V_0$. Et donc le Cas 1 s'applique \`a nouveau.

On a ainsi une suite de cofibrations
$$
A\rightarrow U \rightarrow V\rightarrow B
$$
avec qui les coproduits induisent des \'equivalences faibles
d'apr\`es les cas 1,  2 et 1 respectivement. Ceci prouve que le
coproduit le long
de $A\rightarrow B$ induit lui aussi une \'equivalence faible.
\eop

\subnumero{Compatibilit\'e avec les troncations}

Soit $A$ un $m$-pr\'echamp de Segal, et $A\rightarrow A'$
un remplacement \'equivalent objet-par-objet tel que les $A'(X)$ soient des
$m$-cat\'egories de Segal (e.g. $A'$ est un  remplacement fibrant pour la
topologie
grossi\`ere). On consid\`ere le $n$-pr\'echamp sur $\Xx$:
$$
\tau _{\leq n} A := X\mapsto \tau _{\leq n}(A'(X)).
$$
On utilise ici le mot ``$n$-pr\'ecat'' dans son sens de
\cite{nCAT} et ``$n$-pr\'echamp'' pour pr\'efaisceau de $n$-pr\'ecats.

{\em Remarque:} cette construction a \'et\'e not\'e $\tau _{\leq n}^{\rm
pre}$ par le
deuxi\`eme auteur dans d'autres papiers; en effet, dans ces papiers on avait
r\'eserv\'e
la notation
$\tau _{\leq n}$ pour le champ associ\'e (cf \S 9 et \S 13 ci-dessous) \`a $\tau
_{\leq n}^{\rm
pre}$; mais du point de vue qu'on adopte ici, un pr\'echamp est faiblement
\'equivalent \`a son champ associ\'e, et on n'a pas besoin de faire cette
distinction de notation.
Il faudra cependant faire attention que si $A$ est un champ (cf \S 9) alors
$\tau _{\leq n}A$
n'est plus forcement un champ.

\begin{lemme}
\label{compatible}
Soit $f:A\rightarrow B$ un morphisme de $m$-pr\'echamps de Segal. Alors
$f$ est une \'equivalence faible pour la topologie $\Gg$ si et seulement si,
pour tout $n$,
$$
\tau _{\leq n} (A)\rightarrow
\tau _{\leq n} (B)
$$
est une $\Gg$-\'equivalence faible de $n$-pr\'echamps.
\end{lemme}
{\em Preuve:}
On v\'erifie d'abord que si $A$ et $B$ sont des  $n$-pr\'echamps de
$n$-groupo\"{\i}des
alors $f:A\rightarrow B$ est une \'equivalence faible (en utilisant notre
d\'efinition pour les $n$-pr\'echamps de Segal) si et seulement si $\Re
_{\geq 0} A
\rightarrow \Re _{\geq 0}B$ est une \'equivalence faible d'Illusie de
pr\'efaisceaux simpliciaux. On montre ceci par r\'ecurrence sur $n$, on peut
donc supposer que c'est vrai pour $n-1$. Pour $x,y\in A_0(X)$ on a
$$
Path ^{x,y} \Re _{\geq 0} A \cong \Re _{\geq 0} (A_{1/}(x,y))
$$
(objet-par-objet on d\'eduit ceci de \cite{Tamsamani}), et la m\^eme chose pour
$B$. Gr\^ace \`a l'hypoth\`ese de r\'ecurrence on obtient que
$$
A_{1/}(x,y) \rightarrow B_{1/}(fx,fy)
$$
est une \'equivalence faible si et seulement si
$$
Path ^{x,y} \Re _{\geq 0} A \rightarrow Path ^{fx,fy} \Re _{\geq 0} B
$$
est une \'equivalence faible d'Illusie. D'autre part le pr\'efaisceau
d'ensembles
$\tau _{\leq 0} A$ est isomorphe au pr\'efaisceau d'ensembles $\pi _0\circ
\Re _{\geq 0} A$, donc on obtient que
$$
\tau _{\leq 0} A\rightarrow \tau _{\leq 0} B
$$
induit une surjection entre les faisceaux associ\'es, si et seulement si
$$
\pi _0\circ
\Re _{\geq 0} A\rightarrow \pi _0\circ
\Re _{\geq 0} B
$$
induit une surjection entre les faisceaux associ\'es.
On conclut cette partie en remarquant qu'un morphisme $f:U\rightarrow V$ de
pr\'efaisceaux simpliciaux est une \'equivalence faible d'Illusie, si et
seulement
si (a) pour tout $x,y\in U_0(X)$, $Path ^{x,y}(U)\rightarrow Path^{fx,fy}(V)$
est une \'equivalence faible d'Illusie, et  (b) le morphisme $\pi _0\circ
U\rightarrow \pi _0\circ V$ induit une surjection entre les faisceaux associ\'es
(l'argument pour justifier cette remarque est imm\'ediat).

Maintenant on traite le lemme pour $m=0$: soit $f: U\rightarrow V$ un morphisme
de pr\'efaisceaux simpliciaux. On sait que la troncation de Postnikov
$F \rightarrow\Re \tau _{\leq n}F$ respecte
les $\pi _i$ pour $i\leq n$. Il s'ensuit que $f$ est une \'equivalence faible
d'Illusie si et seulement si
$$
\Re \tau _{\leq n}U\rightarrow \Re \tau _{\leq n}V
$$
est une \'equivalence faible d'Illusie pour tout $n$. En appliquant le
r\'esultat du paragraphe pr\'ec\'edent, on obtient que $f$ est
une \'equivalence faible d'Illusie
si et seulement si
$$
\tau _{\leq n}(f):\tau _{\leq n}U\rightarrow  \tau _{\leq
n}V
$$
est une \'equivalence faible pour tout $n$, ce qui donne le lemme pour $m=0$.

On traite maintenant le cas g\'en\'eral par r\'ecurrence sur $m\geq 1$; on
suppose
que le r\'esultat est vrai pour $m-1$.
Pour $f: A\rightarrow B$ et $x,y\in A(X)$ on a
$$
(\tau _{\leq n}A)_{1/}(x,y) \cong \tau _{\leq n-1}(A_{1/}(x,y))
$$
et de m\^eme pour $B$ (ces \'equivalences entrant dans un diagramme
commutatif de naturalit\'e pour $f$). Les $A_{1/}(x,y)$ sont des
$m-1$-pr\'echamps de Segal.  En utilisant l'hypoth\`ese de
r\'ecurrence, on obtient
que
$$
A_{1/}(x,y)\rightarrow B_{1/}(fx,fy)
$$
est une \'equivalence si et seulement si
$$
(\tau _{\leq n}A)_{1/}(x,y)\rightarrow (\tau _{\leq n}B)_{1/}(fx,fy)
$$
est une \'equivalence pour tout $n$.
D'autre part on a l'isomorphisme naturel
$$
\tau _{\leq 0}A \cong \tau _{\leq 0}(\tau _{\leq n}A).
$$
Donc $\tau _{\leq 0}A\rightarrow \tau _{\leq 0}B$ induit une
surjection entre les
faisceaux associ\'es si et seulement si
$$
\tau _{\leq 0}(\tau _{\leq n}A)\rightarrow \tau _{\leq 0}(\tau _{\leq n}B)
$$
induit une surjection entre les faisceaux associ\'es pour tout $n$.
Ceci termine la preuve du lemme.
\eop

\begin{corollaire}
\label{compPin}
Soient $m$ et $n$ fix\'es. Un morphisme $f: A\rightarrow B$ de
$m$-pr\'echamps de
Segal est
une $\Gg$-\'equivalence si et seulement si
$$
\Pi _{n,Se}\circ f: \Pi _{n,Se}\circ A\rightarrow \Pi _{n,Se}\circ A
$$
est une $\Gg$-\'equivalence de $n+m$-pr\'echamps de Segal.
\end{corollaire}
\eop

\subnumero{Compatibilit\'e avec le produit fibr\'e homotopique}

La cat\'egorie de mod\`eles ferm\'ee ci-dessus est presque propre. Le seul
probl\`eme est que $nSePC$ elle-m\^eme n'est pas propre (pour un contre-exemple,
voir \cite{nCAT}). Les $\Gg$-\'equivalences
sont pr\'eserv\'ees
par produit fibr\'e avec une fibration \`a condition que les valeurs
des $n$-pr\'echamps de Segal consid\'er\'es soient des $n$-cat\'egories de
Segal.
Pour $n=0$ (auquel cas il n'y a pas de condition) le fait que la cmf des
pr\'efaisceaux simpliciaux est propre est prouv\'e par Jardine dans
\cite{JardineBool}.

En fait on a un \'enonc\'e un peu plus fort, \`a savoir la stabilit\'e des
$\Gg$-\'equivalences par produit fibr\'e avec une fibration pour la topologie
grossi\`ere. Autrement dit, le produit fibr\'e homotopique pour la
topologie grossi\`ere pr\'eserve les $\Gg$-\'equivalences. Cet \'enonc\'e a
\'et\'e conjectur\'e pour $n=0$ par C. Rezk dans
un courrier \'el\'ectronique \`a P. Hirschhorn \cite{Rezk}. Rezk
montre la propri\'et\'e analogue
pour les pr\'efaisceaux d'ensembles (i.e. le cas des
$0$-champs non de Segal). Encore plus int\'eressant, Rezk montre, dans le cadre
de la localisation par une pr\'e-topologie (i.e. collection de cribles)
quelconque, que cette propri\'et\'e pour les ``\'equivalences locales'' est
\'equivalente au fait que la pr\'e-topologie soit une topologie i.e.
satisfasse aux
axiomes pour une topologie de Grothendieck. Le lemme suivant a donc \'et\'e
inspir\'e par cette lettre de Rezk et nous remercions P. Hirschhorn de nous
l'avoir fait suivre.

\begin{lemme}
\label{proper}
Si $A\rightarrow B$ est une cofibration et $A\rightarrow C$ une
$\Gg$-\'equivalence faible alors $B\rightarrow B\cup ^AC$ est une
$\Gg$-\'equivalence
faible. D'autre part, si $f:F\rightarrow E$ est une fibration (m\^eme seulement
pour la topologie grossi\`ere) et $g:G\rightarrow E$ est une $\Gg$-\'equivalence
faible, et si $G(X)$ et $E(X)$ sont des $n$-cat\'egories de Segal pour tout
$X\in
\Xx$ alors $F\times _EG\rightarrow F$ est une $\Gg$-\'equivalence faible.
\end{lemme}
{\em Preuve:}
La premi\`ere partie est imm\'ediate (c'est le ``lemme de Reedy''
\cite{Reedy} \cite{DHK} \cite{Hirschhorn}) car tous les objets sont cofibrants.

Pour la deuxi\`eme partie, de fa\c{c}on g\'en\'erale, on notera
$F\times ^h_EG$ le produit fibr\'e homotopique objet-par-objet; il
s'identifie (objet-par-objet)
au produit fibr\'e ordinaire dans le cas o\`u
$F\rightarrow E$ est une fibration pour la topologie
grossi\`ere et les $G(X)$ et $E(X)$ sont des $n$-cat\'egories de Segal (voir
\cite{nCAT} Theorem 6.7).
L'op\'eration de troncation pr\'eserve---presque---le
produit fibr\'e homotopique en ce sens qu'on a une \'equivalence
objet-par-objet
$$
\tau _{\leq n-1} (F\times ^h _EG) \cong
\tau _{\leq n-1}\left( \tau _{\leq n}F\times ^h _{\tau _{\leq n}E}\tau
_{\leq n}G
\right) .
$$
Cette formule se d\'emontre pour les $0$-pr\'echamps de Segal
\`a l'aide de la suite exacte
des groupes d'homotopie pour le produit fibr\'e homotopique d'espaces. On
l'obtient ensuite pour les $n$-pr\'echamps de Segal par r\'ecurrence sur
$n$ (c'est le
m\^eme argument que ci-dessous).

En utilisant cette formule et le lemme
\ref{compatible}, on se ram\`ene au cas o\`u $E$, $F$
et $G$ sont des $n$-pr\'echamps de Segal $n$-tronqu\'es dont les valeurs
sont des
$n$-cat\'egories de Segal (ils correspondent en fait \`a des
$n$-pr\'echamps non de
Segal). Et on va traiter ce cas par r\'ecurrence sur $n$. Par exemple, s'il
s'agissait au d\'epart du cas des pr\'efaisceaux simpliciaux, on se ram\`ene au
cas des pr\'efaisceaux simpliciaux qui sont $n$-tronqu\'es objet-par-objet.

On suppose donc que $E,F,G$ sont des $n$-pr\'echamps de Segal $n$-tronqu\'es,
et que les valeurs $E(X)$, $F(X)$ et $G(X)$ sont des $n$-cat\'egories de Segal;
avec un morphisme $f: F\rightarrow E$ fibrant pour la topologie grossi\`ere,
et un morphisme $g: G\rightarrow E$ qui est une $\Gg$-\'equivalence faible.
Le produit fibr\'e homotopique objet-par-objet se calcule comme le produit
fibr\'e habituel $F\times _EG$. (Ici on utilise le th\'eor\`eme 6.7 de
\cite{nCAT}
dans sa version adapt\'ee aux $n$-cat\'egories de Segal).

Si $(u,v)$ et $(u',v')$ sont dans $(F\times _EG)_0(X)$ alors on a
(en notant par $w$ et $w'$ leurs images dans $E_0(X)$):
$$
(F\times _EG)_{1/}((u,v), (u',v')) =
F_{1/}(u,u')\times _{E_{1/} (w,w')}G_{1/}(v,v').
$$
Dans cette formule, \`a droite, les trois termes sont des
$n-1$-pr\'echamps de Segal $n-1$-tronqu\'es sur $\Xx /X$.
Les morphismes structurels sont, d'une part une fibration pour la topologie
grossi\`ere, et d'autre part une $\Gg$-\'equivalence faible.

 Par r\'ecurrence
sur $n$ on peut supposer que la premi\`ere projection
$$
(F\times _EG)_{1/}((u,v), (u',v'))
\rightarrow F_{1/}(u,u')
$$
est une $\Gg$-\'equivalence faible. D'autre part
$$
\tau _{\leq 0} (F\times ^h_EG) \rightarrow
(\tau _{\leq 0} F)\times _{\tau _{\leq 0} E} (\tau _{\leq 0} G)
$$
est une surjection de pr\'efaisceaux d'ensembles.
La propri\'et\'e en question pour
les pr\'e\-faisceaux d'ensembles (voir \cite{Rezk}) implique que  le morphisme
$$
(\tau _{\leq 0} F)\times _{\tau _{\leq 0} E} (\tau _{\leq 0} G)\rightarrow
\tau _{\leq 0} F
$$
induit une surjection sur les faisceaux associ\'es; donc
$$
\tau _{\leq 0} (F\times ^h_EG)\rightarrow \tau _{\leq 0} F
$$
induit \'egalement
une surjection sur les faisceaux associ\'es. On a v\'erifi\'e tout
ce qu'il faut
pour affirmer que $F\times ^h_EG\rightarrow F$ est une $\Gg$-\'equivalence
faible. \eop

\begin{corollaire}
\label{cordereedy}
Soit
$$
\begin{array}{ccccc}
F&\rightarrow & E & \leftarrow & G\\
\downarrow &&\downarrow && \downarrow \\
F'&\rightarrow & E' & \leftarrow & G'
\end{array}
$$
un diagramme o\`u les morphismes verticaux sont des $\Gg$-\'equivalences
faibles.
Alors le morphisme induit sur le produit homotopique objet-par-objet
$$
F\times ^h_EG\rightarrow F'\times ^h_{E'}G'
$$
est une $\Gg$-\'equivalence faible.
\end{corollaire}
{\em  Preuve:}
Ceci se d\'emontre par un argument standard de Reedy.
\eop

Au \S 9 (corollaire \ref{rezkQ}) on verra que le corollaire pr\'ec\'edent
implique que le
produit fibr\'e homotopique objet-par-objet est stable par l'op\'eration de
passage au ``champ associ\'e''. La question de Rezk \'etait en fait pos\'ee
sous cette
forme.

\subnumero{Donn\'ees de descente}

Notons $Ho(nSePCh(\Xx ^{\rm gro}))$ la cat\'egorie homotopique
de la cat\'egorie de mod\`eles fer\-m\'ee $nSePCh$ pour $\Xx^{\rm gro}$.
Rappelons (\cite{Quillen}) qu'elle peut \^etre d\'efinie comme la localis\'ee
\`a la Gabriel-Zisman \cite{GabrielZisman} de $nSePCh$ par rapport aux
\'equivalences faibles (pour la topologie grossi\`ere) ou \'egalement comme la
cat\'egorie dont les objets sont les $n$-pr\'echamps de Segal fibrants pour la
topologie grossi\`ere (et automatiquement cofibrants) et dont les morphismes
sont les classes d'\'equivalence de morphismes. En fait on dispose de
l'intervalle
$\overline{I}$ (la cat\'egorie avec deux objets $0,1$ isomorphes) et comme tous
les objets sont cofibrants, deux morphismes $f,g:A\rightarrow B$ avec $B$
fibrant, sont homotopes si et seulement s'il existe $h: A\times
\overline{I}\rightarrow B$ avec $h|_{A\times \{ 0\} }=f$ et
$h|_{A\times \{ 1\} }=g$.

Soient $A, B\in nSePCh$. On note
$$
[A,B]_{\Xx ^{\rm gro}}
$$
l'ensemble de morphismes de $A$ vers $B$ dans $Ho(nSePCh(\Xx ^{\rm gro}))$.

De la m\^eme fa\c{c}on on d\'efinit la cat\'egorie homotopique
 $Ho(nSePCh(\Xx ^{\Gg}))$ pour la
topologie $\Gg$,  soit comme la localis\'ee de Gabriel-Zisman de $nSePCh$ par
les $\Gg$-\'equivalences faibles, soit comme la cat\'egorie des objets
$\Gg$-fibrants et classes d'homotopie de morphismes. On note
$$
[A,B]_{\Xx ^{\Gg}}
$$
l'ensemble des morphismes de $A$ vers $B$ dans $Ho(nSePCh(\Xx ^{\Gg}))$.

Soit $X\in \Xx$ et $\Bb \subset \Xx /X$ un crible couvrant $X$.
Soit $A$ un $n$-pr\'echamp de Segal sur $\Xx$. Soit $\ast _{\Bb}$ le
$n$-pr\'echamp sur $\Bb$ \`a valeurs la $n$-pr\'ecat de Segal finale
$\ast$. Une {\em
donn\'ee de descente pour $A$ par rapport au crible $\Bb$} est un \'el\'ement de
l'ensemble
$$
[\ast _{\Bb},A|_{\Bb}]_{\Bb ^{\rm gro}}.
$$
On peut calculer cet ensemble de la fa\c{c}on suivante. Soit $A\rightarrow A'$
un remplacement fibrant pour la topologie grossi\`ere. On verra dans la
prochaine section que $A|_{\Bb}\rightarrow A'|_{\Bb}$ est un remplacement
fibrant pour $\Bb ^{\rm gro}$. Donc une donn\'ee de descente i.e.
\'el\'ement de
$[\ast _{\Bb},A|_{\Bb}]_{\Bb ^{\rm gro}}$ est un morphisme $\ast
_{\Bb}\rightarrow A'|_{\Bb}$, \`a \'equivalence pr\`es (l'\'equivalence \'etant
la relation d'homotopie avec l'intervalle $\overline{I}$ par exemple).
Par leger abus de notation on dira qu'une donn\'ee de descente est un
morphisme $\ast _{\Bb}\rightarrow A'|_{\Bb}$.

Avec le foncteur $p: \Bb \rightarrow \Xx$ et la construction $p_!$ (voir \S 4
ci-dessous) on peut identifier un morphisme
$\ast _{\Bb}\rightarrow A'|_{\Bb}$ avec un morphisme
$p_!\ast _{\Bb}\rightarrow A'$. La $n$-pr\'ecat de Segal $p_!\ast
_{\Bb}$
sur $\Xx$ (qui est en fait un pr\'efaisceau d'ensembles sur $\Xx$) sera
aussi souvent
not\'ee
$\ast _{\Bb}$, et avec cette notation une donn\'ee de descente pour $A$ par
rapport \`a $\Bb$ est un morphisme  $\ast _{\Bb} \rightarrow A'$ \`a
\'equivalence d'homotopie pr\`es.

Pour la question de ``l'effectivit\'e des donn\'ees de descente'' on
s'int\'eresse
aux donn\'ees de descente \`a homotopie pr\`es; c'est pour cela qu'on utilise
la cat\'egorie homotopique ici pour donner la d\'efinition de
``l'ensemble des donn\'ees de descente''. On pourrait aussi
envisager de d\'efinir la {\em $n$-cat\'egorie de Segal des donn\'ees de
descente}
dont le pr\'esent ensemble est le
$\tau _{\leq 0}$
\footnote{
La $n$-cat\'egorie de Segal des donn\'ees de
descente serait $\underline{Hom}(\ast _{\Bb}, A')$ avec les
notations du \S 11, ce qui est \'equivalent
\`a  $\lim _{\leftarrow , \Bb}
A|_{\Bb}$ d'apr\`es \ref{calclim}.
L'ensemble des donn\'ees de descente
peut alors \^etre exprim\'e comme
$$
[\ast _{\Bb},A|_{\Bb}]_{\Bb ^{\rm gro}}=\tau _{\leq 0} \lim _{\leftarrow , \Bb}
A|_{\Bb}.
$$};
dans ce cas la terme ``donn\'ee de descente''
introduit ici devrait \^etre remplac\'ee par ``donn\'ee de descente \`a
homotopie pr\`es''. Cependant, pour l'utilisation qu'on en fera
dans le pr\'esent papier, la terminologie qu'on donne ici convient
bien.

On a maintenant l'\'egalit\'e
$$
[\ast _{\Xx /X},A|_{\Xx /X}]_{(\Xx /X) ^{\rm gro}}
= [\ast , A(X)],
$$
o\`u la notation \`a droite signifie l'ensemble des classes d'homotopie de
morphismes $\ast
\rightarrow A(X)$ dans $nSePC$; et la restriction de $\Xx /X$ \`a $\Bb$
induit donc un morphisme
$$
[\ast , A(X)]\rightarrow
[\ast _{\Bb},A|_{\Bb}]_{\Bb ^{\rm gro}}.
$$
On dira qu'une donn\'ee de descente pour $A$ par rapport \`a $\Bb$
est {\em effective} si elle est dans l'image de ce morphisme, i.e. si elle
provient d'un \'el\'ement de $[\ast , A(X)]$.
La notion {\em d'effectivit\'e des donn\'ees de descente} entrera dans notre
crit\`ere pour \^etre un champ \ref{critere} ci-dessous.

En principe, les donn\'ees de descente peuvent arriver sous forme d'un diagramme
$$
\ast _{\Bb} \leftarrow D \rightarrow A'\leftarrow A
$$
o\`u la premi\`ere et la derni\`ere fl\`eche sont des \'equivalences faibles
objet-par-objet. En fait,  les donn\'ees de descente proviendront le plus
souvent de diagrammes de la forme
$$
\ast _{\Bb} \leftarrow D \rightarrow A.
$$
Au \S 5 ci-dessous on construira une seule fl\`eche $D\rightarrow \ast
_{\Bb}$ ind\'ependamment de $A$, tel que pour tout $n$-pr\'echamp $A$ dont les
valeurs $A(X)$ sont des $n$-cat\'egories de Segal fibrantes, les donn\'ees de
descente pour $A$ proviennent uniform\'ement de diagrammes
$\ast _{\Bb} \leftarrow D \rightarrow
A$.

On pourrait aussi consid\'erer des {\em donn\'ees de descente
g\'en\'eralis\'ees}, qu'on d\'efinirait comme des \'el\'ements de
$[\ast _{\Xx
/X}, A|_{\Xx /X}]_{\Xx ^{\Gg}}$.
Dans cette optique une telle donn\'ee peut provenir par exemple d'un diagramme
$$
\ast _{\Xx /X} \leftarrow D \rightarrow A
$$
o\`u la premi\`ere fl\`eche est une $\Gg$-\'equivalence faible. Tel est
certainement le cas si $D$ est \'equivalent objet-par-objet \`a $\ast_{\Bb}$
pour un crible $\Bb \subset \Xx /X$ couvrant $X$; mais un tel diagramme
pourrait aussi arriver \`a partir d'un {\em hyper-recouvrement} de $X$
(cf l'expos\'e de Verdier dans \cite{SGA4b}).
On laisse au lecteur le soin d'expliciter la construction d'un $D$
ad\'equat \`a partir
d'un hyper-recouvrement.

\numero{Les fonctorialit\'es $p^{\ast}$, $p_{\ast}$, $p_!$}

\label{fonctorialitepage}

Nous traiterons des fonctorialit\'es dans le cas des cat\'egories avec la
topologie grossi\`ere. La compatibilit\'e avec une topologie de Grothendieck
sera trait\'ee ult\'erieurement (voir Corollaire \ref{Gfibrestrict}).

Il convient (du moins pour les auteurs!) de rappeler la distinction
entre adjoint \`a gauche et adjoint \`a droite (il s'agit de deux
foncteurs $L$ et $R$
de sens contraire
entre deux cat\'egories): on aura une formule de la forme
$$
\{ LX \rightarrow Y \} \cong \{ X \rightarrow RY\}, et
$$
{\bf l'adjoint \`a gauche} est le foncteur $L$ qui appara\^{\i}t dans la
formule \`a
gauche de la fl\`eche; et
\newline
{\bf l'adjoint \`a droite} est le foncteur $R$ qui appara\^{\i}t dans la
formule \`a
droite de la fl\`eche.
\newline
On a les morphismes d'adjonction
$$
X \rightarrow RL X, \;\;\; \mbox{et} \;\;\; LRY \rightarrow Y.
$$

On rappelle la notion de {\em foncteur de Quillen} \cite{Quillen}
\cite{DHK} \cite{HoveyBook} \cite{JardineGoerssBook}. Si $M,M'$ sont deux
cat\'egories de mod\`eles
ferm\'ees, un foncteur $f: M\rightarrow M'$ est dit {\em foncteur de Quillen
\`a gauche} s'il pr\'eserve les cofibrations et les cofibrations triviales (en
particulier, on ne demande pas forcement que $f$ pr\'eserve toutes les
\'equivalences faibles, mais c'est parfois le cas tout de m\^eme); et si $f$
admet un adjoint \`a droite. La notion de {\em
foncteur de Quillen \`a droite} s'en d\'eduit par dualit\'e: c'est un
foncteur qui
pr\'eserve les fibrations et les fibrations triviales, et qui admet un
adjoint \`a gauche.
C'est un fait (pas compl\'etement trivial)
\cite{Quillen} que si $f$ est un foncteur de Quillen \`a gauche alors son
adjoint \`a droite est un foncteur de Quillen \`a droite et vice-versa.
On rappelle aussi qu' un foncteur de Quillen \`a gauche  pr\'eserve les
\'equivalences faibles
entre objets cofibrants \cite{JardineGoerssBook} \cite{DHK} \cite{HoveyBook}.

Soit $p: \Zz \rightarrow \Xx$ un foncteur entre $1$-cat\'egories munies de
la topologie
grossi\`ere. Les $n$-pr\'echamps de Segal sont
des pr\'efaisceaux d'objets dans une cat\'egorie admettant
toutes les limites avec toutes les commutations n\'ecessaires
(en effet la cat\'egorie $nSePC$  des $n$-pr\'ecats de
Segal en question,  est elle-m\^eme une cat\'egorie de pr\'efaisceaux
d'ensembles).
En cons\'equence on dispose des op\'erations  $p^{\ast}$, $p_{\ast}$, $p_!$ (et
$p^!$ dont on n'aura pas besoin). On les explicite pour plus de
commodit\'e.

D'abord si $B$ est un $n$-pr\'echamp de Segal
sur $\Xx$ alors $p^{\ast}(B)$ est le
remont\'e de $B$ sur $\Zz$ i.e. le compos\'e
$$
\Zz ^o \rightarrow \Xx ^o \rightarrow nSePC.
$$

Ce foncteur admet des adjoints \`a droite et \`a gauche. Si $A$ est un
$n$-pr\'echamp de Segal sur $\Zz$ on pose
$$
p_{\ast}(A)(X) := \lim _{\leftarrow , p(Z)\rightarrow X}A(Z)
$$
la limite \'etant prise sur la cat\'egorie des paires $(Z, f)$ avec
$f: p(Z)\rightarrow X$. D'autre part on pose
$$
p_!(A)(X) := \lim _{\rightarrow , X\rightarrow p(Z)}A(Z),
$$
la limite \'etant prise sur la cat\'egorie des paires $(Z, i)$ avec
$i: X\rightarrow p(Z)$. On a les formules d'adjonction:
$$
\{ p_!A\rightarrow B\} = \{ A\rightarrow p^{\ast} B\}
$$
et
$$
\{ B\rightarrow p_{\ast}A\} = \{ p^{\ast}B \rightarrow A\} .
$$

Il est clair que le foncteur $p^{\ast}$ pr\'eserve les cofibrations et les
cofibrations triviales. C'est donc un {\em foncteur de Quillen \`a gauche}
et son adjoint \`a droite $p_{\ast}$ envoie les fibrations (resp. fibrations
triviales) sur des fibrations (resp. fibrations triviales).

On note  $A\mapsto \Gamma (\Zz , A)$ le cas particulier
du foncteur $p_{\ast}$ obtenu en faisant $\Xx = \ast$. On dit que
c'est le foncteur des
``sections
globales''.
Il
transforme les fibrations (resp. triviales) en fibrations (resp. triviales).
Il convient parfois de remarquer que $\Gamma (\Zz , A)$ peut \'etre d\'efinie
comme la $n$-pr\'ecat de Segal qui repr\'esente le foncteur
$E\mapsto Hom _{nPC}(\underline{E}, A)$ de $nSePC$ vers $Ens$, o\`u
$\underline{E} = p^{\ast}(E)$ est le $n$-pr\'echamp de Segal constant \`a
valeurs $E$.

On a la formule de composition $(pq)_{\ast} = p_{\ast} q_{\ast}$ (comme
toujours, nous n\'egligeons les \'eventuels probl\`emes de coh\'erence entre
isomorphismes naturels dans des $1$-cat\'egories strictes comme $nSePCh$). En
termes de $\Gamma$ on obtient
$$
\Gamma (\Xx , p_{\ast}(A))=\Gamma (\Zz , A).
$$
Rappelons pour donner un sens homotopique \`a cette formule, que si $A$ est
fibrant sur $\Xx$ alors $p_{\ast}(A)$ est fibrant sur $\Xx$.

Dans l'autre sens la situation n'est pas aussi facile. On aimerait pouvoir dire
que si $B$ est fibrant sur $\Xx$ alors $p^{\ast}B$
est fibrant sur $\Zz$. Ceci n'est
pas toujours le cas (m\^eme pour $p: \Zz \rightarrow \ast$ par exemple).
On a le lemme suivant.

\begin{lemme}
\label{shriek}
Supposons que $p: \Zz \rightarrow \Xx$ a la propri\'et\'e que pour tout $X\in
\Xx$, l'op\'eration
$$
A\mapsto {\rm colim}_{ X\rightarrow p(Z)}A(Z)
$$
pr\'eserve les cofibrations et cofibrations triviales. Alors
$p^{\ast}$
pr\'eserve
les objets fibrants, les fibrations et les fibrations triviales (i.e. c'est
un foncteur
de Quillen \`a droite).
\end{lemme}
{\em Proof:}
L'hypoth\`ese revient exactement \`a dire que le foncteur $p_!$
pr\'eserve les
cofibrations et les cofibrations triviales; dans ce cas son adjoint \`a droite
$p^{\ast}$ pr\'eserve les fibrations et les fibrations triviales.
\eop

{\em Remarque:} Pour $X\in \Xx$ on note habituellement $X/p$ la cat\'egorie
des paires  $(Z,f)$ o\`u $Z\in \Zz$ et $f: X\rightarrow p(Z)$. La limite
dans l'hypoth\`ese du lemme est prise sur $X/p$. Comme on parle de
pr\'efaisceaux, tout est contravariant. En particulier, si $X/p$  est une
cat\'egorie cofiltrante (et notamment si elle a un objet initial) alors
l'hypoth\`ese est satisfaite.

Le cas qui nous int\'eresse est celui o\`u
$\Zz$ est la cat\'egorie $\Xx /Y$ pour un certain objet  $Y\in \Xx$.
Dans ce cas, pour $X\in \Xx$, la cat\'egorie
$$
X/p = \{  X\rightarrow Z \rightarrow Y\}
$$
est la cat\'egorie des diagrammes. Elle se d\'ecompose en r\'eunion disjointe
de cat\'egories $(X/p)_{\varphi}$ index\'ees par les fl\'eches $\varphi :
X\rightarrow Y$. Et la cat\'egorie
$$
(X/p)_{\varphi} = (X\stackrel{\varphi}{\rightarrow} Y) / (\Xx / Y)
$$
est la cat\'egorie des objets au-dessous de $(X,\varphi )$ dans la cat\'egorie
des objets au-dessus de $Y$. En particulier, $(X/p)_{\varphi}$ admet
$(X,\varphi )$ comme objet initial. Ainsi l'op\'eration
$$
{\rm co}\lim _{X/p} A(Z) = \coprod _{\varphi : X\rightarrow Y} A(X)
$$
vue comme foncteur en $A$
pr\'eserve les cofibrations et les cofibrations triviales,
donc pour
$$
p: \Xx /Y\rightarrow \Xx ,
$$
$p_!$ est un foncteur de Quillen \`a gauche et son adjoint $p^{\ast}$ est un
foncteur de Quillen \`a droite. On a la formule
$$
p_!(A) (X) = \coprod _{\varphi : X\rightarrow Y} A(X).
$$

\begin{corollaire}
\label{fibrestrict}
Si $\Xx$ est une cat\'egorie munie de la topologie grossi\`ere alors pour
tout $X\in
\Xx$ la restriction de $n$-pr\'echamps de Segal
$$
A\mapsto A|_{\Xx /X}
$$
pr\'eserve les fibrations (donc les objets fibrants), les
cofibrations, et les
\'equivalences faibles; c'est un foncteur de Quillen \`a
la fois \`a gauche et
\`a droite pour la topologie grossi\`ere.
\end{corollaire}
\eop

On verra plus loin (Lemma \ref{grest}) que la restriction sur $\Xx /X$
pr\'eserve
aussi les $\Gg$-fibrations (il est \'evident sur la d\'efinition qu'elle
pr\'eserve les  $\Gg$-\'equivalences faibles).

Un autre exemple qui nous int\'eresse est celui o\`u $\Zz \subset \Xx$ est
un {\em
crible}, i.e. une sous-cat\'egorie pleine telle que si
$X\rightarrow Y$ est un morphisme de $\Xx$ et si $Y\in \Zz$, alors $X\in \Zz$.
Dans ce cas, la colimite dans \ref{shriek} est, ou bien vide si $X\not \in
\Zz$, ou bien \'egale \`a $A(X)$ si $X\in \Zz$. Cette
construction pr\'eserve les cofibrations
et cofibrations triviales. On obtient le

\begin{corollaire}
\label{crible}
Soit $i:\Zz \hookrightarrow \Xx$ un crible (les deux cat\'egories
consid\'er\'ees \'etant munies de leurs topologies grossi\`eres). Alors
$i_{!}$ est un
foncteur de Quillen \`a gauche et $i^{\ast}$ un foncteur de Quillen \`a droite,
i.e. $i^{\ast}$ pr\'eserve les objets fibrants et \'equivalences faibles
entre eux.
\end{corollaire}
\eop

On peut aussi expliciter $i_{!}$ dans ce corollaire: si $A$ est un
$n$-pr\'echamp
de Segal sur $\Zz$ alors $i_!(A)(X) = A(X)$ pour $X\in \Zz$, et
$i_!(A)(X) = \emptyset$ pour $X\not \in \Zz$.

\numero{La structure de type Bousfield-Kan d'apr\`es Hirschhorn}
\label{bousfieldkanpage}

On veut indiquer ici une autre fa\c{c}on d'obtenir une cat\'egorie de
mod\`eles ferm\'ee pour les $n$-champs de Segal, bas\'ee sur le livre de
Hirschhorn \cite{Hirschhorn} et qui trouve ses origines dans une construction
de Bousfield-Kan \cite{BousfieldKan} (qui \`a son tour remonte \`a Quillen
\cite{Quillen}) pour le cas des pr\'efaisceaux simpliciaux.

Dans \cite{Hirschhorn}, Hirschhorn d\'efinit d'abord la notion de {\em
cat\'egorie de mod\`eles ferm\'ee engendr\'ee par cofibrations}. Nous avons
repris cette notion dans l'\'enonc\'e \ref{dhklemme} ci-dessus, qui peut donc
servir de d\'efinition. Ensuite,
Hirschhorn montre que, pour toute cat\'egorie $\Xx$
petite et toute cmf $M$ engendr\'ee par
cofibrations, la cat\'egorie des
pr\'efaisceaux sur $\Xx$ \`a valeurs dans $M$ admet une structure de cat\'egorie
de mod\`eles ferm\'ee not\'ee
$$
M^{\Xx ^o}.
$$
Pour $M=EnsSpl$ la construction est due \`a Bousfield-Kan (\cite{BousfieldKan},
p. 314). A son tour, leur construction \'etait une application directe de la
construction de Quillen \cite{Quillen} d'une cmf pour les objets simpliciaux
dans une cat\'egorie convenable, en l'occurrence  celle des pr\'efaisceaux
d'ensembles. Pour cette raison nous appellerons {\em structure de HBKQ (i.e.
Hirschhorn-Bousfield-Kan-Quillen)} cette structure de cmf
d\'efinie par Hirschhorn sur $M^{\Xx ^o}$.

Pour $M=nSePC$
\'egale \`a la cat\'egorie des $n$-pr\'ecats de Segal, la structure de HBKQ est
diff\'erente de celle du theor\`eme \ref{cmf}, notre
cat\'egorie de mod\`eles des $n$-pr\'echamps de Segal
pour la topologie grossi\`ere. Les
cat\'egories sous-jacentes sont les m\^emes, ainsi que les \'equivalences
faibles (qui sont les \'equivalences faibles objet-par-objet).
En revanche, les
fibrations de HBKQ sont les morphismes $f:A\rightarrow B$
induisant pour tout $X$ une
fibration $f(X): A(X)\rightarrow B(X)$ ---c'est donc une classe plus
grande que celle des fibrations du th\'eor\`eme \ref{cmf}.
Et l\`a o\`u nous
prenons toutes les injections comme cofibrations, la classe des cofibrations
de HBKQ est plus petite: ce sont les r\'etractions de
morphismes qui
sont obtenus  par une suite (transfinie) d'
additions libres de cellules, o\`u une {\em addition libre de cellules} est une
cofibration  $A\rightarrow A'$ engendr\'ee librement par une cofibration
$A(X)\rightarrow C'$ dans $M$ (pour un seul objet $X$). Plus pr\'ecis\'ement si
on note $p_X: \Xx /X\rightarrow \Xx$ le foncteur naturel, et par $\underline{C}$
le diagramme constant sur $\Xx /X$ \`a valeurs un objet $C$ dans $M$, alors une
libre addition de cellules est un coproduit de la forme
$$
A' = A \cup ^{p_{X,!}(C)}p_{X,!}(C')
$$
pour une cofibration $C\rightarrow C'$ dans $M$ (dans l'expression
pr\'ec\'edente
on avait $C=A(X)$).

L'int\'er\^et des morphismes obtenus par libre addition de cellules est le fait
suivant: si $A\rightarrow A'$ est obtenu par l'addition d'une cellule
$A(X)\rightarrow C'$ au-dessus de $X$, et si  $B$ est un autre diagramme,
alors le prolongement d'un morphisme $f:A\rightarrow B$ en $f':A'\rightarrow B$
n'est rien d'autre que le prolongement du morphisme
$f(X): A(X)\rightarrow B(X)$ en un morphisme $C'\rightarrow B(X)$.

Comme les deux cat\'egories et leurs sous-cat\'egories d'\'equivalences faibles
sont les m\^emes, les deux th\'eories homotopiques sont les m\^emes (d'apr\`es
la th\'eorie de la localisation de Dwyer-Kan
\cite{DwyerKan1} \cite{DwyerKan2} \cite{DwyerKan3}).

{\em Conclusion:} On aurait pu prendre, comme cat\'egorie de mod\`eles ferm\'ee
pour les $n$-pr\'echamps de Segal sur $\Xx$ avec la topologie grossi\`ere, la
structure de HBKQ sur $M^{\Xx ^o}$ avec pour $M$ la cat\'egorie de
mod\`eles ferm\'ee $nSePC$ des $n$-pr\'ecats de Segal (Th\'eor\`eme
\ref{SeCmf}).
Dans la structure de HBKQ, les fibrations sont simplement les morphismes
qui sont objet-par-objet des fibrations de $n$-pr\'ecats de Segal.

{\em Remarque:} Si $\Xx$ admet des produits fibr\'es alors les cofibrations de
HBKQ sont stables par restriction aux sites $\Xx /X$ et aussi stables par
produit direct. Ceci permet d'appliquer la th\'eorie des $Hom$
internes
qui sera expliqu\'ee au \S 11 ci-dessous.

Pour l\'egitimer cet autre point de vue possible, il reste \`a voir
comment int\'egrer une topologie de Grothendieck dans cette th\'eorie.
Soit encore $M$ la cat\'egorie des $n$-pr\'ecats de Segal du Th\'eor\`eme
\ref{SeCmf}. Sur $M^{\Xx ^o}$ qui est la cat\'egorie des $n$-pr\'echamps de
Segal, on garde la notion de cofibration de HBKQ; on utilise la
notion de $\Gg$-\'equivalence faible d\'efinie ci-dessus; et on d\'efinit les
{\em $\Gg$-fibrations de HBKQ} comme les morphismes qui poss\`edent la
propri\'et\'e de rel\`evement pour les cofibrations de HBKQ qui sont
des $\Gg$-\'equivalences faibles.

\begin{theoreme}
\label{cmfHirschho}
La cat\'egorie des $n$-pr\'echamps de Segal munie des cofibrations de
HBKQ, des
$\Gg$-\'equivalences faibles, et des $\Gg$-fibrations de HBKQ, est une
cat\'egorie de mod\`eles ferm\'ee.
\end{theoreme}
{\em Preuve:}
Un morphisme est une cofibration $\Gg$-triviale de HBKQ, si et seulement
si c'est une cofibration de HBKQ et aussi une cofibration $\Gg$-triviale
pour la structure \ref{cmf}. Hirschhorn prouve que ses cofibrations (pour la
topologie grossi\`ere) sont stables par coproduit \cite{Hirschhorn}.
Donc la pr\'eservation des cofibrations $\Gg$-triviales de HBKQ par
coproduit est une
con\-s\'e\-quen\-ce imm\'ediate du m\^eme r\'esultat pour
la structure du Th\'eor\`eme
\ref{cmf}.
Les autres propri\'et\'es sont aussi imm\'ediates.
\eop

Le lecteur trouvera peut-\^etre cette structure de cat\'egorie de mod\`eles
ferm\'ee plus conviviale que celle de \ref{cmf}. Nous avons gard\'e la
structure de \ref{cmf} comme notion de base parce qu'elle est
compatible avec la structure de Jardine (devenue standard notamment
en  K-th\'eorie) pour les pr\'efaisceaux simpliciaux i.e. pour le cas des
$0$-pr\'echamps de Segal.

\subnumero{Calcul des classes d'homotopie d'applications}

La structure de HBKQ permet de calculer plus facilement les classes
d'homotopie d'applications. Fixons une cat\'egorie $\Xx$ munie de la topologie
grossi\`ere. La cat\'egorie homotopique $Ho(nSePCh(\Xx ^{\rm gro}))$
est d\'efinie comme la localis\'ee de Gabriel-Zisman de $nSePCh(\Xx )$ par la
sous-cat\'egorie des \'equivalences faibles objet-par-objet. En particulier,
$Ho(nSePCh(\Xx ^{\rm gro}))$ ne d\'epend pas du choix de structure de cmf
compatible avec ces \'equivalences faibles; en particulier on peut utiliser la
structure de HBKQ pour calculer les ensembles $Hom$ dans la cat\'egorie
homotopique, i.e. les ensembles de classes d'homotopie d'applications.

Commen\c{c}ons par une remarque d'ordre g\'en\'eral: si $E$ est une
$n$-pr\'ecat de Segal et si $\underline{E}$ est le $n$-pr\'echamp de Segal
constant
sur $\Xx$ \`a valeurs $E$, alors le produit direct avec $\underline{E}$
respecte les cofibrations de HBKQ (sans hypoth\`ese sur l'existence de produits
fibr\'es dans $\Xx$, et sans supposer que $\underline{E}$ soit lui-m\^eme
cofibrant). En effet, si $A\rightarrow A'$ est obtenu par addition d'une
cellule $A(X)\rightarrow C$ au-dessus de l'objet $X$, alors
$A\times \underline{E}\rightarrow A'\times \underline{E}$ est obtenu par
addition de la cellule $A(X)\times E \rightarrow C\times E$. Plus
g\'en\'eralement si $E\rightarrow E'$ est une cofibration de $n$-pr\'ecats de
Segal et $A\rightarrow A'$ est une cofibration de HBKQ alors
$$
A\times \underline{E}' \cup ^{A\times \underline{E}}A'\times \underline{E}
\rightarrow A'\times \underline{E}'
$$
est encore une cofibration de HBKQ.

Soient $A$ et $B$ deux
$n$-pr\'echamps de Segal sur $\Xx$ et supposons que
$B(X)$ est fibrant pour tout $X\in \Xx$. Alors $B$ est un objet fibrant pour la
structure HBKQ (pour la topologie grossi\`ere). Si on choisit un remplacement
HBKQ-cofibrant $C\rightarrow A$ (\'equivalence faible objet-par-objet avec
$C$ cofibrant de HBKQ) alors tout \'el\'ement de $[A,B]_{\Xx ^{\rm gro}}$
provient d'un morphisme $C\rightarrow B$. Notons $\overline{I}$ la
cat\'egorie avec deux objets $0,1$ et un isomorphisme entre eux
(consid\'er\'ee comme $n$-cat\'egorie de Segal, voir ``Induction'' au \S 2).
Et notons aussi par $\overline{I}$ le
$n$-pr\'echamp de Segal constant \`a valeurs $\overline{I}$.  Le diagramme
$$
C\times \{ 0,1\} \rightarrow C\times \overline{I}\rightarrow C
$$
constitue un ``objet $C\times I$'' dans la terminologie de Quillen
\cite{Quillen} pour la structure de HBKQ, et en particulier la premi\`ere
application
est une cofibration. Donc, deux applications $C\rightarrow B$ sont homotopes,
i.e. induisent le m\^eme \'el\'ement de $[A,B]_{\Xx ^{\rm gro}}$, si et
seulement si elles sont li\'ees par une homotopie de la forme $C\times
\overline{I}\rightarrow B$. On a ainsi obtenu l'\'enonc\'e suivant.

\begin{lemme}
\label{howtocalc}
Soit $A$ un $n$-pr\'echamp de Segal sur $\Xx$ et fixons un remplacement
$C\rightarrow A$ avec $C$ cofibrant
de HBKQ. Alors pour tout $n$-pr\'echamp de Segal $B$ tel que
$B(X)$ soit une $n$-cat\'egorie de Segal fibrante pour tout $X\in \Xx$,
l'ensemble de classes d'homotopie d'applications
$$
[A,B]_{\Xx ^{\rm gro}}
$$
se calcule comme l'ensemble de  morphismes
de $n$-pr\'echamps de Segal $C\rightarrow B$, modulo la relation d'\'equivalence
qui identifie deux morphismes s'ils sont li\'es par une homotopie de
la forme $C\times \overline{I}\rightarrow B$.
\end{lemme}
\eop

\subnumero{Expression uniforme des donn\'ees de descente}

On peut utiliser ce lemme \ref{howtocalc} pour donner une expression uniforme
des donn\'ees de descente \`a valeurs dans un $n$-pr\'echamp de Segal $A$
tel que tous les $A(X)$ soient des $n$-cat\'egories de Segal fibrante. Pour
cette discussion  il faut avoir $n\geq 1$ (dans le
cas $n=0$, il faut remplacer $A$ par $\Pi _{1,Se}\circ A$
voir \S 10 par exemple).

Soit $X\in \Xx$ et soit $\Bb \subset \Xx /X$ un crible couvrant $X$.
Choisissons un remplacement HBKQ-cofibrant
$$
D\rightarrow \ast _{\Bb}
$$
(disons, pour fixer les id\'ees, qu'on parle ici de $n$-pr\'echamps de Segal
sur $\Bb$---mais on peut aussi prendre le $p_!$ pour se retrouver sur $\Xx /X$
ou sur $\Xx$). Maintenant, d'apr\`es le lemme \ref{howtocalc}, on sait
que pour tout $n$-pr\'echamp de Segal $A$ tel que tous les $A(X)$ soient
fibrants,
l'ensemble des donn\'ees de descente
$$
[\ast _{\Bb}, A|_{\Bb}]_{\Bb ^{\rm gro}}
$$
s'exprime comme l'ensemble de morphismes $D\rightarrow A|_{\Bb}$
modulo la relation qui identifie deux morphismes s'ils sont li\'es par
une homotopie $D\times \overline{I}\rightarrow A|_{\Bb}$.

Supposons maintenant que le site $\Xx$ admet des produits fibr\'es et que
le crible $\Bb$ provient d'une famille couvrante $\Uu = \{ U_{\alpha}
\rightarrow X\}$. Dans ce cas on peut donner une formule explicite pour un
remplacement HBKQ-cofibrant $D\rightarrow \ast _{\Bb}$.

Pour cela, soit $\overline{I}^{(p)}$ la cat\'egorie avec $p+1$ objets
not\'es $0',\ldots , p'$ et un isomorphisme entre chaque paire d'objets
(c'est en quelque sorte le ``produit sym\'etrique'' de $p$ exemplaires de
$\overline{I}$). Consid\'erons-la comme $n$-pr\'ecat de Segal
(par ``induction'' voir \S 2). On obtient un foncteur not\'e
$$
R: \Delta \rightarrow nSePC,\;\;\;\;
R(p):= \overline{I}^{(p)}.
$$
Ce foncteur munit $nSePC$ d'une structure de cmf simpliciale \cite{Quillen}.
En particulier, si $K$ est un ensemble simplicial (on n'a pas besoin dans ce
cas que $K$ soit fini) et $E$ une $n$-pr\'ecat de Segal, on obtient
la $n$-pr\'ecat de Segal $K\otimes ^RE$; dans notre cas on peut \'ecrire
$$
K\otimes ^R E = (K\otimes ^R\ast )\times E.
$$
Ici $K\otimes ^R\ast$ est juste le coproduit des $R(p)$ sur les $p$-simplexes
de $K$ (pris suivant le proc\'ed\'e habituel de r\'ealisation d'un ensemble
simplicial). On dira que $K\otimes ^R\ast$ est la {\em r\'ealisation de $K$ par
rapport \`a $R$}.

Cette construction \'etant fonctorielle, s'\'etend \`a des pr\'efaisceaux
d'ensembles simpliciaux: si $K$ est un pr\'efaisceau d'ensembles simpliciaux sur
$\Xx$ on obtient le $n$-pr\'echamp de Segal $K\otimes ^R\ast$.

Maintenant notre famille couvrante $\Uu$ repr\'esente un pr\'efaisceau
d'ensembles: c'est tout simplement le pr\'efaisceau r\'eunion disjointe
de ceux repr\'esent\'es
par les $U_{\alpha}$. On obtient le pr\'efaisceau simplicial
$$
P(\Uu ):= \ldots \Uu \times _X \Uu \tworightarrows \Uu
$$
i.e. $P(\Uu )_k=\Uu \times _X \ldots \times _X\Uu $ ($k+1$ fois).
C'est le {\em nerf du recouvrement $\Uu$}. En appliquant la construction
``r\'ealisation par rapport \`a $R$'' on obtient le $n$-pr\'echamp de Segal
$$
D:= P(\Uu )\otimes ^R \ast .
$$
Le morphisme $D\rightarrow \ast _{\Bb}$ est une
\'equivalence objet-par-objet.  En effet $P(\Uu )(Y)$ est contractile si
$Y\in \Bb$
et vide sinon. Observons que l'objet $D$ est d\'efini
comme un $n$-pr\'echamp sur $\Xx /X$, de sorte que
pour d\'efinir le morphisme $D\rightarrow \ast
_{\Bb}$ il faut interpr\^eter $\ast _{\Bb}$ comme $n$-pr\'echamp sur $\Xx /X$.

Nous pr\'etendons que $D$ est HBKQ-cofibrant. En fait il est int\'eressant de voir
explicitement comment obtenir $D$ \`a partir de $\emptyset$ par une suite
d'additions libres de cellules. D'abord on ajoute $R(0)=\ast$ au-dessus de
chaque
\'el\'ement $U_{\alpha}$ de la famille couvrante. Appelons $D^0$ le
r\'esultat obtenu. Ensuite on ajoute  $R(1)=\overline{I}$ au-dessus de chaque
$U_{\alpha}
\times _XU_{\beta}$, via le diagramme
$$
D^0(U_{\alpha} \times _XU_{\beta}) \leftarrow p_1^{\ast}R(0)\sqcup
p_2^{\ast}R(0) = \{ 0,1\} \hookrightarrow \overline{I}.
$$
Notons $D^1$ le r\'esultat obtenu. Plus g\'en\'eralement notons
$\partial R(p)$ le
``bord'' de $R(p)$, i.e. la r\'eunion dans $R(p)$ des $p+1$ faces de la forme
$R(p-1)$. On ajoute $R(2)$ au-dessus de chaque $U_{\alpha}\times _XU_{\beta}
\times _XU_{\gamma}$ via le diagramme
$$
D^1(U_{\alpha}\times _XU_{\beta}
\times _XU_{\gamma})\leftarrow \partial R(2) \hookrightarrow  R(2)
$$
pour obtenir $D^2$. A chaque \'etape on ajoute \`a $D^{p-1}$ un exemplaire de
$R(p)$ au-dessus de chaque
$$
U_{\alpha _0,\ldots , \alpha _p}:=U_{\alpha _0}
\times _X\ldots \times _X U_{\alpha_p},
$$
via le diagramme
$$
D^{p-1}(U_{\alpha _0} \times _X\ldots \times _X
U_{\alpha_p})\leftarrow \partial R(p) \hookrightarrow  R(p).
$$
Enfin $D$ est la colimite (reunion croissante) des $D^p$.

On note que $D$ peut \^etre consid\'er\'e comme un $n$-pr\'echamp de Segal
au-dessus
de $\Bb$ et la description ci-dessus montre qu'en tant que tel (aussi bien
qu'en tant que $n$-pr\'echamp de Segal au-dessus de $\Xx /X$ ou $\Xx$), il est
HBKQ-cofibrant. On peut donc l'utiliser dans le lemme \ref{howtocalc}.
On obtient l'\'enonc\'e suivant.

\begin{lemme}
\label{calcdonnees}
Supposons que $\Xx$ admet des produits fibr\'es.
Soit $\Bb \subset \Xx /X$ un crible engendr\'e par une famille $\Uu$
couvrant $X$. Soit $D\rightarrow \ast _{\Bb}$ le $n$-pr\'echamp de Segal
construit ci-dessus. Alors pour tout $n$-pr\'echamp de Segal $A$ sur $\Xx$ tel
que tous les $A(X)$ soient fibrants, l'ensemble des donn\'ees de descente
$$
[\ast _{\Bb}, A|_{\Bb}]_{\Bb ^{\rm gro}}
$$
se calcule comme l'ensemble de morphismes $D\rightarrow A$ modulo la relation
qui identifie deux morphismes s'ils sont li\'es par un morphisme
$D\times \overline{I}\rightarrow A$.
\end{lemme}
\eop

On peut donner une description explicite de ce que c'est qu'un morphisme
$D\rightarrow A$. On proc\`ede par r\'ecurrence sur le $p$ de $D^p$. Si on a un
morphisme $D^{p-1}\rightarrow A$, alors on obtient pour chaque
$\alpha _0,\ldots , \alpha _p$ le morphisme compos\'e
$$
\partial R(p)\rightarrow D^{p-1}(U_{\alpha _0,\ldots ,\alpha _p})\rightarrow
A(U_{\alpha _0,\ldots ,\alpha _p}).
$$
Une extension en morphisme $D^p\rightarrow A$ n'est rien d'autre qu'une
extension pour chaque $\alpha _0,\ldots , \alpha _p$, du morphisme ci-dessus en
un morphisme
$$
R(p)\rightarrow
A(U_{\alpha _0,\ldots ,\alpha _p}).
$$
Un morphisme $D\rightarrow A$ consiste en une telle collection
d'extensions successives.

En r\'esum\'e: un morphisme $\eta : D\rightarrow A$ est un ensemble de
donn\'ees comme suit:
pour chaque $\alpha$, un objet
$\eta (\alpha )\in A(U_{\alpha})$; pour chaque $\alpha , \beta$ une
\'equivalence
$$
\eta (\alpha , \beta ): \overline{I}\rightarrow A(U_{\alpha \beta})
$$
entre les restrictions des objets $\eta (\alpha )$ et $\eta (\beta )$;
pour chaque $\alpha , \beta , \gamma$ un morphisme
$$
\eta (\alpha , \beta , \gamma ):
\overline{I}^{(2)}\rightarrow A(U_{\alpha \beta\gamma})
$$
qui, restreint au triangle $\partial R(2) \subset \overline{I}^{(2)}$,
est \'egal au morphisme obtenu par restriction des $\eta (\alpha ,\beta )$,
$\eta (\beta , \gamma )$ et $\eta (\alpha , \gamma )$; etc.

\subnumero{Commentaire}
Dans cette description des donn\'ees de descente, on a
utilis\'e de fa\c{c}on essentielle l'hypoth\`ese que les $A(X)$ sont des
$n$-cat\'egories de Segal fibrantes. C'est pour cette
raison qu'on ne voit pas appara\^{\i}tre de $2$-morphismes, $3$-morphismes
etc
dans les donn\'ees. Par exemple, comme $A(U_{\alpha \beta \gamma})$ est
fibrant, la donn\'ee comme dans \cite{Breen} de trois $1$-fl\`eches
avec une $2$-\'equivalence entre $fg$ et $h$, peut \^etre remplac\'ee par
un vrai morphisme de $R(2)$ dans $A(U_{\alpha \beta \gamma})$
(i.e. en quelque sorte un objet de $A(U_{\alpha \beta \gamma})_{2/}$).

Ceci explique pourquoi nous arrivons \`a traiter le cas des $n$-cat\'egories
sans trop nous emm\^eler dans des ``formules de cocycles'' mais en fait
c'est un d\'esavantage si on cherche une compr\'ehension explicite de ce qui se
passe, car la notion de $n$-cat\'egorie de Segal fibrante n'est pas tr\`es
explicitement calculable.
\footnote{
D'un point de vue plus optimiste on peut cependant observer
que si $T$ est un espace topologique
(ou ensemble simplicial de Kan) alors $\Pi _{n,Se}(T)$ est fibrant; donc les
$n$-pr\'echamps de Segal qui proviennent par application de $\Pi _{n,Se}$ aux
pr\'efaisceaux d'espaces topologiques (ou ensembles simpliciaux de Kan)
v\'erifient
la condition  en question.}
En particulier, si on n'arrive pas \`a donner
une forme totalement ``en termes de cocycles'' de la notion de donn\'ees de
descente, \`a
la fa\c{c}on dont Breen l'a fait pour les $2$-champs dans
\cite{Breen},
c'est parce que ce probl\`eme est plus difficile que ceux qu'on peut
actuellement r\'esoudre. On peut formuler un

\begin{probleme}
Trouver une structure de cmf analogue \`a la structure de HBKQ, pour les
$n$-pr\'ecats (de Segal) dont les objets fibrants
soient les $n$-cat\'egories de Segal.
\end{probleme}

Avec une solution affirmative \`a ce probl\`eme, on pourrait reprendre la
discussion ci-dessus des donn\'ees de descente: si $D'\rightarrow \ast _{\Bb}$
est un remplacement cofibrant pour la structure recherch\'ee dans ce
probl\`eme (qui donnerait aussi une structure de cmf pour les
$\Xx^o$-diagrammes par la
m\'ethode de HBKQ), alors on pourrait dire que pour tout $A$ tel que les $A(X)$
soient des $n$-cat\'egories de Segal, une donn\'ee de descente pour $A$ par
rapport \`a $\Bb$ est la m\^eme chose qu'un morphisme $D'\rightarrow A$.
Dans ce cas, on verrait appara\^{\i}tre (beaucoup) des $i$-fl\`eches
pour tout $i$,
et on obtiendrait en fait une description des donn\'ees de
descente en termes de cocycles au sens
classique du terme (qui g\'en\'eraliserait \cite{Breen}).

\numero{Le point de vue de la localisation}
\label{localisationpage}


Le passage de la structure de cat\'egorie de mod\`eles
ferm\'ee pour la topologie grossi\`ere \`a celle pour la topologie $\Gg$,
est une {\em localisation de Bousfield}, op\'eration qui est
tr\`es bien expliqu\'ee dans Hirschhorn \cite{Hirschhorn}; voir aussi
Goerss-Jardine \cite{JardineGoerssLoc} et qui remonte bien s\^ur \`a
Bousfield-Kan \cite{BousfieldKan}. L'op\'eration de passage de la topologie
grossi\`ere \`a la topologie $\Gg$ pr\'eserve la notion
de cofibration (aussi bien dans le cadre de \ref{cmf} que dans le cadre de
\ref{cmfHirschho}), donc  c'est une localisation de Bousfield \`a gauche
dans la terminologie de \cite{Hirschhorn}.
Il y avait une erreur dans notre d\'emonstration (de {\tt v1, v2}) 
du lemme \ref{pourloc} ci-dessous. Nous remercions D. Dugger d'avoir trouv\'e 
cette faute, et B. Toen de nous l'avoir communiqu\'ee. Le lemme \ref{pourloc},
et donc la proposition \ref{naqualocaliser}, ne sont pas vrais dans
le cadre des $n$-pr\'echamps de Segal (et m\^eme pas pour les ensembles
simpliciaux i.e. les $0$-pr\'echamps de Segal). Le probl\`eme
provient de l'interaction entre l'homotopie en degr\'e arbitrairement grand,
et la cohomologie en degr\'e arbitrairement grand. Comme montre Dugger dans son
papier r\'ecent \cite{Dugger}, on obtient des \'enonc\'es corrects en 
rempla\c{c}ant
partout ``crible couvrant'' par ``hyperrecouvrement''. Nous ne souhaitons  pas
entrer dans les d\'etails ici et nous renvoyons \`a \cite{Dugger} pour cela.
On se restreint donc, pour cette version {\tt v3}, \`a \'enoncer la proposition
\ref{naqualocaliser} et le lemme \ref{pourloc} dans le cadre des $n$-pr\'echamps
non de Segal (autrement dit, au cas ``$n$-tronqu\'e'' dans lequel il n'y a de
l'homotopie qu'en degr\'e $\leq n$).

Cette modification se propage aux r\'esultats suivants: 
\newline
---on doit refaire la preuve du corollaire \ref{Gfibrestrict};
\newline
---on doit modifier les \'enonc\'es des crit\`eres 
\ref{newcritere}, \ref{critereaveclim} et \ref{critere7};
\ref{newcritere} et les parties (c) de ces deux derniers 
ne s'appliquent dor\'enavant qu'au cas $n$-tronqu\'e;
\newline
---on doit modifier la caract\'erisation de l'image essentielle
dans \ref{correlation};
\newline
---on doit rajouter une hypoth\`ese (5) dans l'\'enonc\'e du th\'eor\`eme
\ref{dansleschamps};
\newline
---la preuve de \ref{descente} est l\'eg\`erement modifi\'ee mais le r\'esultat
reste le m\^eme;
\newline
---et enfin, le r\'esultat de descente pour les complexes \ref{cestunchamp}
n'est dor\'enavant \'enonc\'e que pour les complexes born\'es inf\'erieurement.

Pour revenir aux  r\'esultats \ref{naqualocaliser} et \ref{pourloc} de cette version
{\tt v3} nous
nous restreignons donc \`a la consid\'eration des cmf des $n$-pr\'ecats ou
$n$-pr\'echamps non de Segal. On note qu'on d\'efinit de fa\c{c}on analogue
au Th\'eor\`eme \ref{cmf} (resp. \ref{cmfHirschho}) une cmf des $n$-pr\'echamps
(non de Segal); pour le cas de \ref{cmf} on le notera $nPCh$. 

Si $p:\Yy \rightarrow \Xx$ est un foncteur et $U$ une $n$-pr\'ecat 
alors on note $U_{\Yy}$ le $n$-pr\'echamp
$p_!(\underline{U})$ o\`u
$\underline{U}$ est le $n$-pr\'echamp
constant sur $\Yy$ de valeur $U$.

On va consid\'erer des cofibrations du type suivant: soit $X\in \Xx$ et
$\Bb \subset \Xx /X$ un crible; et soit $U\hookrightarrow U'$ une cofibration de
$n$-pr\'ecats. On note
$$
cof(U\hookrightarrow U'; \Bb \subset
\Xx /X)
$$
la cofibration
$$
U_{\Xx /X} \cup ^{U_{\Bb}}U'_{\Bb} \rightarrow U'_{\Xx /X}.
$$

\begin{proposition}
\label{naqualocaliser}
La
cat\'egorie de mod\`eles ferm\'ee $nPCh$ de \ref{cmf} (resp. de \ref{cmfHirschho}) pour
la topologie $\Gg$, mais pour les $n$-pr\'echamps non de Segal,
est la  {\em localis\'ee de Bousfield \`a gauche}
(\cite{Hirschhorn} D\'efinition 3.3.1) de la cat\'egorie de mod\`eles
ferm\'ee pour
la topologie grossi\`ere, par rapport aux morphismes
$cof(U\hookrightarrow U'; \Bb \subset \Xx /X)$, avec $\Bb$
crible couvrant $X$ pour $\Gg$.
\end{proposition}
{\em Preuve:} On traite la structure de \ref{cmf}, l'autre cas \'etant
analogue.
Nous allons juste montrer l'\'etape-cl\'e qui est le lemme suivant.
Ce lemme
impliquera, par la th\'eorie de Hirschhorn \cite{Hirschhorn} \S 3, la
proposition.
\eop
\begin{lemme}
\label{pourloc}
Si $E\hookrightarrow E'$ est une cofibration $\Gg$-triviale de 
$n$-pr\'echamps (non de Segal),
alors  il
existe une
suite (\'e\-ven\-tu\-elle\-ment transfinie voir la discusion de
\cite{Hirschhorn}
\cite{DHK}
\cite{HoveyBook}) de cofibrations
$$
E=E_{-1} \hookrightarrow E_0
\hookrightarrow \ldots \hookrightarrow E_k \hookrightarrow \ldots
$$
de colimite not\'ee $E_{\infty}$, et une cofibration triviale pour
la topologie grossi\`ere $E_{\infty} \hookrightarrow E''$ telle que
$E\rightarrow E'$ soit une r\'etraction de $E\rightarrow E''$; et telle que les
cofibrations $E_i \hookrightarrow E_{i+1}$ soient  ou bien des coproduits de
cofibrations de la forme $cof(U\hookrightarrow U'; \Bb \subset \Xx /X)$,
ou bien
des cofibrations triviales pour la topologie grossi\`ere.
\end{lemme}
{\em  Preuve:}
On appellera {\em suite ad\'equate} toute suite de cofibrations comme dans
l'\'enonc\'e du lemme, et {\em compos\'e} de cette suite, la colimite
correspondante
$E_{\infty}$. On va prouver que si $f:E\rightarrow E'$ est une fibration
$\Gg$-triviale alors il existe une suite ad\'equate, avec un morphisme
d\'efini sur le compos\'e
$E_{\infty}\rightarrow E'$ qui \'etend $f$ et qui est une \'equivalence
pour la topologie grossi\`ere. Ceci impliquera le lemme (exercice laiss\'e au
lecteur).

On prouve ce r\'esultat par r\'ecurrence sur $n$. On montre
d'abord le pas de
r\'ecurrence
puis le cas $n=0$. Soit $n\geq 1$ et on suppose donc le r\'esultat vrai
pour $n-1$. On peut sans perte supposer que $E$ et $E'$ sont fibrants pour la
topologie grossi\`ere. Pour tout $X\in \Xx$ et toute
paire d'objets $x,y\in E(X)$,
le morphisme $E_{1/} (x,y)\rightarrow E'_{1/}(x,y)$ est une fibration
$\Gg$-triviale de $n-1$-pr\'echamps au-dessus de $\Xx /X$; il
est donc
\'equivalent (pour la topologie grossi\`ere) au compos\'e d'une suite
ad\'equate. On remarque que pour tout coproduit de $E_{1/}(x,y)$ avec une
cofibration de $n-1$-pr\'echamps  $V\rightarrow V'$ au-dessus de $\Xx
/X$, on peut d\'efinir un coproduit correspondant de $E$ avec $p_!\Upsilon
(V)\rightarrow p_!\Upsilon (V')$ o\`u $p: \Xx /X \rightarrow \Xx$
(ici on utilise la notation $\Upsilon$ introduite dans \cite{limits}).

D'autre part, si $\Bb \subset \Xx /Y$ est un crible couvrant $Y\in \Xx /X$
et $U\hookrightarrow U'$ une cofibration des $n-1$-pr\'ecats alors
$$
p_! \Upsilon (cof(U\hookrightarrow U'; \Bb \subset (\Xx /X) / Y))
$$
est identique \`a
$$
cof(\Upsilon (U)\hookrightarrow \Upsilon (U'); \Bb \subset \Xx /Y);
$$
donc en appliquant $p_!\Upsilon $ \`a une suite ad\'equate
de $n-1$-pr\'echamps sur $\Xx /X$, on obtient une suite ad\'equate
de
$n$-pr\'echamps sur $\Xx$.

Donc, \`a partir de la suite de
cofibrations pour $E_{1/} (x,y)\rightarrow E'_{1/}(x,y)$ on obtiendra une suite
ad\'equate de cofibrations pour
$$
E\cup ^{p_! \Upsilon (E_{1/} (x,y))}p_! \Upsilon (E'_{1/} (x,y)).
$$
On peut ensuite ajouter
les coproduits avec une cofibration $\ast _{\Bb} \rightarrow \ast
_{\Xx}$ (qui a la bonne forme
$cof (\emptyset \hookrightarrow \ast , \Bb \subset \Xx /X )$)
pour chaque morphisme $u:\ast _{\Bb} \rightarrow E$ se
prolongeant en
$u':\ast _{\Xx /X}\rightarrow E'$.
Notons que le r\'esultat obtenu admet un morphisme vers $E'$ qui est
(gr\^ace \`a ``trois pour le prix de deux'') une $\Gg$-\'equivalence faible.

Si  on it\`ere cette construction pour
tout $(X,x,y)$ et tout $(X,\Bb ,u,u')$  et cela une infinit\'e d\'enombrable de
fois, on
obtient une factorisation
$$
E\rightarrow E'' \rightarrow E'
$$
o\`u le premier morphisme est un compos\'e de cofibrations de la
forme voulue; et
le deuxi\`eme est une $\Gg$ \'equivalence faible et qui induit une
\'equivalence pour la topologie grossi\`ere sur tous
les $E'' _{1/}(x,y)\rightarrow
E'_{1/}(x,y)$; de plus $E''\rightarrow E'$ a la propri\'et\'e de rel\`evement
pour les morphismes de la forme $\ast _{\Bb} \rightarrow \ast _{\Xx /X}$. Ceci
implique que $E'' \rightarrow E'$ est une \'equivalence pour la topologie
grossi\`ere, ce qui donne le r\'esultat voulu.

Il reste \`a traiter le cas $n=0$, i.e. le cas des pr\'efaisceaux d'ensembles.
Pour cela, on applique essentiellement deux fois de plus le pas de recurrence
ci-dessus. Dans ce cas on pourrait l'\'ecrire explicitement de fa\c{c}on plus
\'el\'ementaire. On laisse cela aux soins du lecteur.
\eop

Rezk montre dans \cite{Rezk} que
la cmf des pr\'efaisceaux pour la th\'eorie des faisceaux, s'obtient,
\`a partir de la cmf
des pr\'efaisceaux sans topologie, par la localisation ci-dessus
(i.e. Rezk traite le cas des $0$-champs non de Segal). En outre il sugg\`ere
\`a la fin de sa lettre que la m\^eme chose devrait \^etre vraie pour la
structure de Jardine.

La proposition reste \'egalement vraie pour la structure de cat\'egorie de
mod\`eles ferm\'e de type HBKQ. En fait, les r\'esultats de
\cite{Hirschhorn} permettraient de d\'efinir une notion de {\em $M$-champ}
et une
cat\'egorie de mod\`eles ferm\'ee correspondante, pour toute cat\'egorie de
mod\`eles ferm\'ee $M$, cellulaire et propre \`a gauche (cf la terminologie de
\cite{Hirschhorn}). On prendrait la cat\'egorie de mod\`eles ferm\'ee de
diagrammes $M^{\Xx ^o}$ et on localiserait par rapport \`a tous les morphismes
$cof(U\hookrightarrow U', \Bb \subset \Xx /X)$ (qui ne sont malheureusement plus
des cofibrations pour cette structure) obtenus en faisant varier
$U\hookrightarrow U'$ dans la classe des
cofibrations de $M$.
On peut maintenant faire l'analogue de
\ref{fibrestrict}, pour la
topologie $\Gg$.
On voudrait obtenir cela m\^eme pour les $n$-pr\'echamps de Segal,
en d\'epit de la restriction au cas $n$-tronqu\'e pour le lemme
\ref{pourloc} dans cette {\tt v3}. 

On note $\Upsilon ^{\langle k \rangle }$ le 
$k$-i\`eme it\'er\'e de l'op\'eration $\Upsilon$ introduite dans
\cite{limits} (mais appliqu\'ee cette fois-ci aux $n$-pr\'ecats de Segal); 
si $A$ est une $n-k$-(pr\'e)cat\'egorie
(de Segal) alors $\Upsilon ^{\langle k\rangle } A$ est une $n$-(pr\'e)cat\'egorie
(de Segal) et ceci s'\'etend de fa\c{c}on \'evidente aux pr\'echamps. 

Soit $p: \Yy \rightarrow \Xx$ un foncteur, et $A\rightarrow B$ une cofibration
de $n$-pr\'echamps de Segal au-dessus de $\Yy$. On obtient alors un carr\'e
cart\'esien
$$
\begin{array}{ccc}
p_!\Upsilon A & \rightarrow & p_!\Upsilon B \\
\downarrow && \downarrow \\
\Upsilon p_!A & \rightarrow & \Upsilon p_!B
\end{array}
$$
et par r\'ecurrence on obtient pour tout $k$ le carr\'e cart\'esien
$$
\begin{array}{ccc}
p_!\Upsilon ^{\langle k \rangle }
A & \rightarrow & p_!\Upsilon ^{\langle k \rangle }B \\
\downarrow && \downarrow \\
\Upsilon ^{\langle k \rangle }p_!A & \rightarrow 
& \Upsilon ^{\langle k \rangle } p_!B
\end{array} .
$$

En particulier, tout coproduit avec une cofibration de la forme
$$
\Upsilon ^{\langle k \rangle }p_!A  \rightarrow 
\Upsilon ^{\langle k \rangle } p_!B
$$
peut aussi \^etre interpret\'e comme coproduit avec une cofibration
de la forme 
$$
p_!\Upsilon ^{\langle k \rangle}
A \rightarrow  p_!\Upsilon ^{\langle k \rangle}B .
$$

\begin{corollaire}
\label{Gfibrestrict}
Soit $\Xx$ un site et $X\in \Xx$ avec $p: \Xx /X\rightarrow \Xx$ le foncteur
d'oubli. Soit $\Gg$ la topologie sur $\Xx$ qui induit une topologie (not\'ee
aussi $\Gg$) sur $\Xx /X$. Alors $p_!$ transforme les cofibrations
$\Gg$-triviales de $n$-pr\'echamps de Segal
sur $\Xx$ en cofibrations $\Gg$-triviales sur $\Xx $.
Autrement dit, $p_!$ est un foncteur de Quillen \`a gauche et $p^{\ast}$ un
foncteur de Quillen \`a droite par rapport aux structures de \ref{SeCmf} pour la
topologie $\Gg$. Donc $p^{\ast}$ pr\'eserve les objets fibrants et les
fibrations.
\end{corollaire}
{\em Preuve:}
On traite d'abord le cas $n=0$. Soit $f:A\rightarrow B$ une \'equivalence faible
d'Illusie de pr\'efaisceaux simpliciaux sur $\Xx /X$. On affirme que $p_!f$ est
une \'equivalence faible d'Illusie sur $\Xx$. Soit $Y\in \Xx$ et
$a\in (p_!A)Y$. On rappelle que cela veut dire qu'on a un morphisme
$g: Y\rightarrow X$ et $a\in A(Y,g)$. Soit (pour $i\geq 1$)
$$
\eta \in \pi _i((p_!B )Y, fa) = \pi _i(B(Y,g), fa).
$$
Par la propri\'et\'e d'Illusie de $f$ sur $\Xx / X$,
il existe un recouvrement ouvert de $(Y,g)$ par des $(Y_j,g_j)$ tel que 
les restrictions de $\eta$ sur les $Y_j$ proviennent d'\'el\'ements de
$\pi _i (A(Y_j,g_j), fa|_{Y_j})$. Ceci donne la m\^eme chose pour quand on
consid\`ere la propri\'et\'e d'Illusie sur $\Xx$. De la m\^eme fa\c{c}on
si 
$$
\phi ,\phi  '\in \pi _i((p_!A)Y, a) = \pi _i(A(Y,g), a)
$$
et si $f\phi = f\phi '$ dans 
$\pi _i((p_!B )Y, fa) = \pi _i(B(Y,g), fa)$ alors
on a un recouvrement de $(Y,g)$ par des $(Y_j,g_j)$ tel que  
$f\phi |_{Y_j} = f\phi '|_{Y_j}$.  En outre ce recouvrement pour $\Xx / X$ est
\'egalement un recouvrement de $Y$ dans $\Xx$ et on obtient cette propri\'et\'e
pour $f$ sur $\Xx$. Il ne reste qu'\`a traiter le cas de $\pi _0$. Soit 
$$
\eta \in \pi _0((p_!B)Y) = \coprod _{Y\rightarrow X}\pi _0(B Y).
$$
Il existe donc $g: Y\rightarrow X$ tel que $\eta $ provient d'un \'el\'ement
$u$ de $\pi _0(B(Y,g))$. Il existe donc un recouvrement de $(Y,g)$ tel que les
$u|_{Y_j}$ proviennent de $A(Y_j,g_j))$ et en particulier $\eta|_{Y_j}$ provient
de $p_!A(Y_j)$ puisque provenant de la  composante $A(Y_j,g_j)$.
Si 
$$
\phi , \phi '\in \pi _0((p_!A)Y) = \coprod _{Y\rightarrow X}\pi _0(A Y)
$$
et $f\phi =f\phi '$ alors en particulier $f\phi$ et $f\phi '$ sont dans le m\^eme
composant $\pi _0(B Y)$ de $(p_!B)Y$ ce qui implique que $\phi$ et $\phi'$
sont dans le m\^eme composant $\pi _0(AY)$ de $\pi _0(p_!A)(Y)$. Le m\^eme argument
que ci-dessus permet encore de conclure.

Maintenant notons que le proc\'ed\'e de la preuve de \ref{pourloc} permet, dans le
cas des $n$-pr\'echamps de Segal, de d\'ecomposer un morphisme $E\rightarrow E'$
qui est une $\Gg$-\'equivalence faible, en une suite de coproduits avec des
cofibrations de la forme 
$$
q_!\Upsilon ^{\langle k \rangle } (F\hookrightarrow F')
$$
avec $q: \Xx /Y\rightarrow X$ et 
o\`u les $F\hookrightarrow F'$ sont ou bien de la forme 
$$
F = \ast _{\Bb} \hookrightarrow \ast _{\Xx /Y} = F'
$$ 
o\`u $\Bb$ est un crible recouvrant $Y$, ou bien
 des cofibrations triviales d'Illusie
de pr\'efaisceaux simpliciaux dans le cas $k=n$. Pour cela 
on utilise le principe de
commutation de $q_!$ avec $\Upsilon$ pour les coproduits \'etabli avant le
pr\'esent \'enonc\'e. 
D'autre part, toute cofibration de cette forme 
$q_!\Upsilon ^{\langle k \rangle } (F\hookrightarrow F')$ est une
$\Gg$-\'equivalence sur $\Xx$. 

Si $f:E\rightarrow E'$ est une $\Gg$-\'equivalence au-dessus de $\Xx / X$ alors
d\'ecomposons-le comme dans le paragraphe pr\'ec\'edent (sur le site $\Xx /X$). 
Si $p: \Xx /X\rightarrow \Xx$ est la projection alors pour $q: Y\rightarrow X$
on obtient $p_! q_! = (pq)_!$. L'image de la d\'ecomposition de
$f$ par $p_!$ est donc une d\'ecomposition de la m\^eme forme pour $p_!f$.
Un morphisme qui se d\'ecompose ainsi est une $\Gg$-\'equivalence, donc 
$p_!f$ est une $\Gg$-\'equivalence sur $\Xx$.

Ceci prouve que $p_!$ est un foncteur de
Quillen \`a gauche. Il s'ensuit que  son adjoint $p^{\ast}$ est un foncteur de
Quillen \`a droite, i.e. pr\'eserve les fibrations.
\eop

{\em Remarque:} Ce corollaire est valable aussi pour la structure de HBKQ.

Dans la structure de HBKQ on a aussi une r\'eciproque 
(nous pensons que ce r\'esultat
reste vrai pour le cas des $n$-pr\'echamps de Segal mais ({\tt v3})
nous ne l'\'enon\c{c}ons
que pour le cas $n$-tronqu\'e):

\begin{corollaire}
\label{restrpourhbkq}
Un morphisme $f: A\rightarrow B$ de $n$-pr\'echamps non de Segal 
est $\Gg$-fibrant pour la structure de HBKQ
si et seulement si $A|_{\Xx /X} \rightarrow B|_{\Xx /X}$ est $\Gg$-fibrant pour
HBKQ pour chaque $X\in \Xx$.
\end{corollaire}
{\em Preuve:}
Une direction est fournie par le corollaire pr\'ec\'edent (pour la structure
HBKQ). Pour l'autre direction, supposons que $f|_{\Xx /X}$ est une
$\Gg$-fibration de HBKQ pour tout $X$. Alors $f$ poss\`ede la propri\'et\'e de
rel\`evement pour les cofibrations qui sont triviales objet-par-objet
(en fait d'apr\`es la construction de Hirschhorn
\cite{Hirschhorn}, les cofibrations triviales
pour la topologie grossi\`ere sont des r\'etractes de compos\'es de suites
de cofibrations obtenues par additions libres de cellules correspondant
\`a des
cofibrations triviales dans la cmf de base $M$); et $f$ satisfait la
propri\'et\'e de rel\`evement pour les cofibrations
$cof_{\Xx }(U\hookrightarrow U'; \Bb \subset \Xx /Y)$. Par l'analogue du
lemme \ref{pourloc} pour la structure de HBKQ, on obtient que $f$ est une
$\Gg$-fibration de HBKQ.
\eop

{\bf Contre-exemple:} Le corollaire \ref{restrpourhbkq} n'est pas vrai pour la
structure de \ref{cmf} (pour le cas de Segal au moins), 
si la cat\'egorie $\Xx$ n'a pas d'objet final; et ceci
m\^eme pour la topologie grossi\`ere. On donne un contre-exemple avec $n=0$,
i.e. pour les pr\'efaisceaux simpliciaux. La cat\'egorie sous-jacente $\Xx$
sera la
cat\'egorie avec $1$ seul objet $x$, et un groupe $\zz$ d'automorphismes de
$x$. Un pr\'efaisceau simplicial au-dessus de $X$ est juste un ensemble
simplicial avec action de $\zz$, i.e. un couple $(B,t)$ o\`u $B$ est un
ensemble simplicial et $t$ un automorphisme de $B$. Supposons que $t$ agit sans
point fixe \`a chaque niveau $B_k$. Si $(A,s)$ est un autre objet, et si
l'automorphisme $s$ a un point fixe dans un $A_k$, alors il n'existe pas
de morphisme $(A,s)\rightarrow (B,t)$. Ceci conduit \`a une contradiction avec
la propri\'et\'e du corollaire \ref{restrpourhbkq}. En effet,
la cat\'egorie $\Xx /x$ est la r\'eunion disjointe d'exemplaires
de la cat\'egorie
triviale, donc $(B,t)|_{\Xx /x}$ est fibrant si et seulement si $B$ est un
ensemble simplicial fibrant de Kan. Or il existe un tel $B$
avec $t$ agissant sans
point fixe---par exemple on prend le complexe singulier total $Sing(\rr )$
avec la translation qui agit sur $\rr$ par $t: b\mapsto b+1$. Maintenant on peut
fabriquer une cofibration triviale
$$
g:(B,t)\hookrightarrow (A,s)
$$
tel que $(A,s)$ ait des points fixes, pour compl\'eter l'exemple. On prend
$A=Sing (\rr ^2)$ avec l'automorphisme $s: (a,b)\mapsto (a,b+a)$ et le morphisme
$B\rightarrow A$ d\'efini par $g:b\mapsto (1,b)$. Notons que $A$ et $B$ sont
contractiles, donc $g$ est une cofibration triviale pour la structure \ref{cmf}
(la condition pour \^etre une cofibration est juste d'\^etre injectif).
Or $A$ a pour lieu fixe le sous-complexe $Sing(\{ 0\} \times \rr
)\subset Sing
(\rr ^2)$. En particulier il n'existe pas d'application $(A,s)\rightarrow (B,t)$
qui induise l'identit\'e sur $(B,t)$. Ceci montre que $(B,t)$ n'est pas fibrant
pour la structure de \ref{cmf}.

\subnumero{La localisation de Boole}

Signalons l'existence du papier de Jardine sur la
``localisation de Boole'' \cite{JardineBool}: il est probable qu'on pourrait
donner un traitement de la cmf de \ref{cmf} avec la technique de
\cite{JardineBool} mais nous n'avons pas explor\'e cette voie.


\numero{Cat\'egories de Segal et cat\'egories simpliciales}

\label{segsimplpage}

La plupart des exemples de $n$-champs de Segal concernent le cas $n=1$.
Il convient de souligner encore une fois qu'un $1$-champ de Segal est en
fait un certain type de $\infty$-champ, o\`u les morphismes en degr\'e $\geq 2$
sont inversibles. En particulier, un $1$-champ de Segal peut \^etre un $n$-champ
pour $n>1$ (c'est le cas quand les ensembles simpliciaux de morphismes sont
$n-1$-tronqu\'es).

Une {\em cat\'egorie simpliciale}
\footnote{On fait ici un l\'eger abus de notation: en principe
une ``cat\'egorie simpliciale'' devrait plut\^ot
\^etre un objet simplicial dans $Cat$,
en particulier les objets devraient former
un ensemble simplicial. Nous
imposons la condition suppl\'ementaire
que l'ensemble simplicial d'objets soit constant i.e. un
ensemble
discret.}
est une cat\'egorie enrichie sur les ensembles
simpliciaux (voir \cite{Kelly}). C'est donc un ensemble d'objets
$Ob(C)$, muni, pour toute paire d'objets $x,y\in Ob(C)$, d' un ensemble
simplicial
$Mor_C(x,y)$ dit des morphismes de $x$ vers $y$, avec une loi de composition
$$
Mor_C(x,y)\times Mor_C(y,z)\rightarrow Mor_C(x,z)
$$
qui est strictement associative, et des identit\'es $1_x\in
(Mor_C(x,x))_0$.

Les cat\'egories simpliciales
ont jou\'e, depuis longtemps d\'ej\`a, un r\^ole
cl\'e en th\'eorie de l'homotopie, \`a commencer par Kan et Quillen; voir
ensuite
\cite{DwyerKan1} \cite{DwyerKan2} \cite{DwyerKan3}, particuli\`erement
l'introduction de
\cite{DwyerKanDiags}, et plus r\'ecemment \cite{DHK}; dans la m\^eme volume que
\cite{DwyerKanDiags} on rel\`eve imm\'ediatement d'autres papiers qui utilisent
de fa\c{c}on essentielle cette notion, par exemple celui de Waldhausen.

Si $C$ est une cat\'egorie simpliciale, on d\'efinit la {\em $1$-cat\'egorie de
Segal associ\'ee}, qu'on note encore $C$, comme la cat\'egorie de Segal
d\'efinie par
$$
C_0:= Ob(C),
$$
et  pour $x_0,\ldots , x_p\in C_0$,
$$
C_{p/}(x_0,\ldots ,x_p):= Mor_C(x_0,x_1)\times Mor _C(x_1,x_2)\times \ldots
\times Mor_C(x_{p-1},x_p).
$$
Les morphismes de la structure simpliciale sont d\'efinis (de fa\c{c}on
\'evidente)
en utilisant la composition en cas de besoin. Par exemple, la composition
elle-m\^eme devient le morphisme de face $02$
$$
C_{1/}(x,y)\times C_{1/}(y,z)=C_{2/}(x,y,z)\rightarrow C_{1/}(x,z).
$$

Un foncteur entre deux cat\'egories simpliciales est une \'equivalence si et
seulement si le morphisme correspondant entre $1$-cat\'egories de Segal est une
\'equivalence. On va montrer, en utilisant la th\'eorie de la localisation de
Dwyer-Kan qui sera trait\'ee ci-dessous, que toute $1$-cat\'egorie de Segal est
\'equivalente \`a une cat\'egorie simpliciale. Plus
pr\'ecis\'ement, si $A$ est une $1$-cat\'egorie de Segal et $A'$ son
remplacement fibrant alors  il existe une cat\'egorie simpliciale
$C$ et une
\'equivalence $C\rightarrow A'$.

Pour notre traitement des morphismes, champs, etc. nous utilisons
syst\'ematiquement le point de vue des $1$-cat\'egories de Segal. Signalons
ici que
les techniques que nous utilisons
sont aussi disponibles, ou en voie de d\'eveloppement, pour les
cat\'egories simpliciales. Dwyer, Hirschhorn et Kan d\'ecrivent une
cat\'egorie de
mod\`eles ferm\'ee pour les cat\'egories simpliciales, o\`u les cofibrations
sont les r\'etractions d'extensions libres (it\'er\'ees) \cite{DHK}.
D'autre part, Cordier et Porter ont d\'evelopp\'e une th\'eorie des morphismes
``ho\-mo\-to\-pi\-que\-ment-coh\'erents'' entre deux cat\'egories simpliciales
\cite{CordierPorter}. Si $C$ et $D$ sont deux cat\'egories simpliciales, ils
d\'efinissent une cat\'egorie simpliciale $Coh(C,D)$ de foncteurs et
transformations naturelles ``coh\'erents'' entre $C$ et $D$. La notion de
composition pour de tels foncteurs est un peu probl\'ematique (car la
$2$-cat\'egorie de Segal $1SeCAT$ n'est pas \'equivalente \`a une
$2$-cat\'egorie stricte simpliciale), et c'est le sujet
d'une grande partie du papier \cite{CordierPorter}.

Il reste notamment \`a faire le lien entre ces diff\'erents
points de vue. On peut isoler deux probl\`emes:

\begin{probleme}
\label{probleme1}
Donner une \'equivalence (\`a la Quillen) entre la cmf des cat\'egories
simpliciales de Dwyer-Hirschhorn-Kan \cite{DHK}, et la cmf des $1$-pr\'ecats de
Segal du th\'eor\`eme \ref{SeCmf}; \end{probleme}

Le probl\`eme ci-dessus est en grande partie 
trait\'e dans Dwyer-Kan-Smith \cite{DKS}:
ils donnent une \'equivalence entre les cat\'egories homotopiques; il reste \`a
v\'erifier que leur construction donne un foncteur de Quillen. 

\begin{probleme}
\label{probleme2}
Etant donn\'ees deux cat\'egories
simpliciales $C$ et $D$, construire une \'e\-qui\-va\-len\-ce de
$1$-cat\'egories de Segal
$$
Coh (C,D) \cong \underline{Hom}(C',D')
$$
o\`u $Coh(C,D)$ est la cat\'egorie simpliciale d\'efinie par Porter et Cordier
\cite{CordierPorter},
$C'$ et $D'$ sont des remplacements fibrants des
cat\'egories de Segal correspondant \`a $C$ et $D$
et $\underline{Hom}(C',D')$ est le $Hom$ interne des $1$-pr\'ecats de Segal (qui
sera trait\'e au \S 11 ci-dessous).
\end{probleme}

On fait un pas vers le premier de ces probl\`emes dans le corollaire
\ref{catsimpl} ci-dessous, o\`u on montre que toute $1$-cat\'egorie de Segal est
\'equivalente \`a une $1$-cat\'egorie simpliciale.

{\em Remarque:} La cat\'egorie de mod\`eles ferm\'ee de Dwyer-Hirschhorn-Kan
ne peut pas \^etre ``interne'' (voir \S 11 ci-dessous) car les cofibrations
sont des extensions libres (it\'er\'ees) de cat\'egories et ne sont
donc pas stables
par produit direct. En particulier, la construction de \cite{nCAT}
(voir aussi \S 11 ci-dessous) d'une cat\'egorie stricte enrichie pour les
objets en question, ne peut pas s'adapter \`a
cette cat\'egorie de mod\`eles ferm\'ee.

En fait, il ne
peut pas exister une cat\'egorie stricte $SplCAT$ enrichie sur les
$1$-cat\'egories simpliciales et ayant le bon type d'homotopie (i.e.
\'equivalente \`a $1SeCAT$). Car si $SplCAT$ existait, on pourrait effectuer
la troncation $\tau _{\leq 3}SplCAT$, c'est-\`a-dire une troncation $\leq 2$
pour les cat\'egories simpliciales  de morphismes, ou encore une troncation
$\leq 1$ pour les
ensembles simpliciaux de $2$-morphismes. Or la troncation $\tau_{\leq 1}$ donne
une $1$-cat\'egorie stricte, et est compatible aux produits directs. Ceci
entrainerait que $\tau _{\leq 3}SplCAT$ serait une $3$-cat\'egorie stricte.
Mais $1SeCAT$ a des produits de Whitehead non triviaux, ce qui implique
que $\tau _{\leq 3}(1SeCAT)$
ne peut pas \'etre \'equivalente \`a une $3$-cat\'egorie stricte.

\subnumero{Produits et coproduits homotopiques}

Dans une cat\'egorie simpliciale on a une notion de {\em produit fibr\'e
homotopique} et la notion duale de {\em coproduit homotopique}. En fait, on
dispose de ces notions dans une $n$-cat\'egorie de Segal cf \cite{limits}
mais on se borne ici au cas $n=1$. C'est un cas particulier
de la notion de limite (ou de
colimite), voir \cite{limits} et \S 14 ci-dessous.

Soit $A$ une cat\'egorie simpliciale, avec $X,Y,Z\in A_0$ et des
morphismes (i.e. sommets de l'ensemble simplicial $Hom _A(\cdot , \cdot )$)
$f: X\rightarrow Y$ et $g: Z\rightarrow Y$. Si $U\in A_0$ est un objet
muni de morphismes $p: U\rightarrow X$ et $q: U\rightarrow Z$ et d'une homotopie
$\gamma : fp\sim gq$ on dira que $(U,p,q,\gamma )$ est un {\em produit fibr\'e
homotopique} et on \'ecrira
$$
U= X\times ^h_YZ
$$
si, pour tout objet $V\in A_0$ le morphisme $\mu _{V,U}$ induit par $(p,q,
\gamma )$
$$
Hom _A(V,U)\rightarrow Hom _A(V,X)\times ^h_{Hom _A(V,Y)}Hom _A(V,Z)
$$
est une \'equivalence faible; \`a droite il s'agit du produit fibr\'e
homotopique d'ensembles simpliciaux (obtenu par remplacement d'un des morphismes
par une fibration de Kan), et le morphisme est d\'efini \`a homotopie pr\`es
apr\`es remplacement \`a droite par un ensemble simplicial de Kan. En fait
le lecteur pourrait imaginer qu'on a $fp=gq$ et que $\gamma$ est l'homotopie
constante---ce sera le cas notamment dans nos exemples provenant des
cat\'egories
de mod\`eles ferm\'ees---mais si on imposait cette condition alors
la notion ne serait pas invariante par \'equivalence $A\cong A'$ de
cat\'egories simpliciales. Si
$fp$ et $gq$ sont \'egaux, le morphisme ci-dessus est bien d\'efini
\`a valeurs dans le produit fibr\'e strict habituel (non-homotopique).

 Si $(V, p',q',\gamma ')$ se trouve \^etre un autre  produit fibr\'e homotopique
(pour les m\^emes donn\'ees)
alors la fibre homotopique de $\mu _{V,U}$  au-dessus de $(p',q',\gamma ')$ est
contractile, tout comme celle de $\mu _{U,V}$
au-dessus de $(p,q,\gamma)$, ce qui
fournit des ensembles
simpliciaux contractiles canoniques de morphismes $U\rightarrow V$ et
$V\rightarrow U$, et assure l'unicit\'e essentielle du produit fibr\'e
homotopique $(U,p,q,\gamma )$ s'il en existe un.

On dira que {\em $A$ admet des produits fibr\'es homotopiques} s'il
existe un
produit fibr\'e homotopique $(U,p,q,\gamma )= X\times ^h_YZ$ pour tout
couple de morphismes $X\rightarrow Y\leftarrow Z$ de $A$.

On a une notion duale de {\em coproduit homotopique}: si
$X\stackrel{f}{\leftarrow} Y\stackrel{g}{\rightarrow} Z$ est un
diagramme dans
$A$, le coproduit homotopique correspondant est un
quadruplet$(U, i,j,\delta )$ avec $i: X\rightarrow U$,
$j: Z\rightarrow U$ et $\delta$ une homotopie $if\sim jg$. On laisse au lecteur
le soin de donner des pr\'ecisions analogues \`a ce qui est dit plus
haut pour le produit.
On notera le coproduit homotopique
$$
(U,i,j,\delta ) = X\cup _h^YZ,
$$
et on dira que $A$ {\em admet des coproduits homotopiques} s'il existe
un tel coproduit
homotopique pour tout couple de morphismes
$X\stackrel{f}{\leftarrow} Y\stackrel{g}{\rightarrow} Z$.

Ces notions sont les g\'en\'eralisations aux objets d'une cat\'egorie
simpliciale quelconque, des $holim$ et $hocolim$ de Bousfield-Kan
\cite{BousfieldKan} {\em et al.}.

Ces notions peuvent \^etre d\'efinies de la m\^eme fa\c{c}on
pour une $1$-cat\'egorie de Segal $A$ \'a condition de
choisir un scindage $A_{1/}\times
_{A_0}A_{1/}\rightarrow A_{2/}$ (un tel scindage existe par exemple
si $A$ est fibrante)
fournissant un morphisme de composition
$A_{1/}\times
_{A_0}A_{1/}\rightarrow A_{1/}$; et on n'a pas besoin des relations de
coh\'erence
sup\'erieures. On laisse ce d\'eveloppement
au lecteur, ainsi que le soin de v\'erifier que ces
notions de limites et colimites sont compatibles avec les d\'efinitions de
\cite{limits}.

Dans \cite{GabrielZisman}, Gabriel et Zisman ont introduit une notion de produit
fibr\'e homotopique dans une $2$-cat\'egorie stricte, et ils l'ont utilis\'ee
pour construire le foncteur $\Omega$ sur la cat\'egorie homotopique $Ho(Top)$
donnant lieu \`a la suite exacte longue d'une fibration. Dualement ils
ont introduit une
notion de coproduit homotopique et l'ont utilis\'ee pour construire le
foncteur $\Sigma$
sur $Ho(Top)$. Dans les deux cas, en utilisant la structure de
$2$-cat\'egorie, ils
produisent des limites ou colimites bien d\'efinies dans la
$1$-cat\'egorie obtenue par troncation (en l'occurrence la cat\'egorie
homotopique). On peut voir la notion de produit fibr\'e homotopique ci-dessus
comme une g\'en\'eralisation de leur notion au cas des
$n$-cat\'egories ($1$-groupiques) quand
$n$ tend vers l'infini. Expliquons en quoi les deux constructions sont
compatibles:
si $A$ est une cat\'egorie simpliciale alors la construction de
\cite{GabrielZisman}
pour $Top$ s'\'etend en une version de $\tau _{\leq 2}A$ qui est une
$2$-cat\'egorie stricte; et on peut alors appliquer l'argument de
\cite{GabrielZisman} \`a
$\tau _{\leq 2}A$. Si $A$ admet des produits fibr\'es homotopiques au sens
ci-dessus, alors $\tau _{\leq 2}A$ v\'erifie les conditions (A), (B),
(C)
de \cite{GabrielZisman},
et le produit qu'ils d\'efinissent est l'image de notre produit
fibr\'e homotopique dans le tronqu\'e $\tau _{\leq 1}A$
(en fait ils regardent
surtout la fibre homotopique au-dessus d'un point de base mais pour le produit
fibr\'e l'argument est le m\^eme).

L'argument de Gabriel-Zisman permet---en passant par $\tau _{\leq 2}A$---
de munir $\tau _{\leq 1}A$
d'une structure suppl\'ementaire avec les
op\'erations $\Omega$ (quand $A$ est point\'e et admet des
produits fibr\'es homotopiques) et $\Sigma$ (quand $A$ est point\'e et admet des
coproduits homotopiques). Alternativement
(mais on souligne quand-m\^eme que l'argument de Gabriel-Zisman
date de 1964) on peut induire directement vers
le tronqu\'e $\tau _{\leq 1}A$, ces op\'erations qu'on peut d\'efinir
sur $A$ par
$$
\Omega X := \ast \times ^h_X \ast
$$
et
$$
\Sigma Y:= \ast \cup _h^Y \ast .
$$
Une version de l'argument de Gabriel-Zisman adapt\'ee aux cat\'egories
d.g. de complexes, est donn\'ee par Bondal-Kapranov
\cite{BondalKapranov}. Nous ne nous sommes aper\c{c}u de l'existence de
l'argument de Gabriel-Zisman qu'apr\`es que M. Kontsevich nous a
indiqu\'e
le r\'esultat de Bondal-Kapranov.

\numero{Localisation de Dwyer-Kan}
\label{dwyerkanpage}

Une bonne source d'exemples de cat\'egories simpliciales, et donc de
$1$-cat\'egories de Segal, est la th\'eorie des localis\'ees simpliciales
$L(C,W)$
(et leurs variantes $L^H(C,W)$ qui sont \'equivalentes) d\'evelopp\'ee dans
les
papiers de Dwyer et Kan \cite{DwyerKan1}, \cite{DwyerKan2}, \cite{DwyerKan3}.
Si $C$ est une cat\'egorie et $W$ une sous-cat\'egorie on obtient une
cat\'egorie simpliciale
$L(C,W)$ ``en inversant homotopiquement les fleches de $W$''
(souvent la sous-cat\'egorie ``des \'equivalences'' $W$ est sous-entendu,
dans ce cas on utilisera
la notation $L(C)$ voir Convention \ref{convenir} ci-dessous).  Nous
renvoyons le lecteur aux
travaux \cite{DwyerKan1}, \cite{DwyerKan2}, \cite{DwyerKan3} pour plus sur
les notations et
d\'efinitions.

Ici, le terme ``localisation'' n'a pas exactement
la m\^eme signification qu'au  \S 6 (bien qu'il y ait un rapport entre les
deux notions). Pour le reste du papier nous utiliserons le mot ``localisation''
pour parler de la localisation de Dwyer-Kan (ou de son analogue pour les
$n$-cat\'egories
de Segal).

Pour \^etre historiquement correct il convient de pr\'eciser que Quillen
\cite{Quillen} avait d\'ej\`a, avec
sa notion de cat\'egorie de mod\`eles ferm\'ee simpliciale,  trouv\'e
la plupart des cat\'egories simpliciales qu'on construit avec la localisation de
Dwyer-Kan; et en particulier, dans tous nos exemples, la cmf $M$
qu'on localise admet d\'ej\`a une
structure simpliciale et on pourrait donc faire r\'ef\'erence \`a
\cite{Quillen}.
Cependant, l'observation de Dwyer-Kan selon laquelle
la structure simpliciale ne d\'epend,
\`a \'equivalence pr\`es, que de $M$ et de sa sous-cat\'egorie d'\'equivalences
faibles est une am\'elioration tr\`es importante par rapport \cite{Quillen}
qui permet
d'obtenir facilement des fonctorialit\'es sans avoir \`a v\'erifier qu'un
foncteur pr\'eserve la structure simpliciale. C'est pour cette raison
que nous utilisons
syst\'ematiquement cette construction de Dwyer-Kan.

On commence par une remarque plus ou moins triviale mais qui
n'appara\^{\i}t pas explicitement dans
\cite{DwyerKan1}, \cite{DwyerKan2}, \cite{DwyerKan3}. C'est
une variante du Corollaire 3.6 de \cite{DwyerKan2} dans laquelle on
consid\`ere des foncteurs non n\'ecessairement
adjoints et des transformations naturelles de direction quelconque.
Pour exprimer cette observation, convenons d'appeler {\em cha\^{\i}ne de
transformations naturelles} entre deux foncteurs $A$ et $B$ la donn\'ee
d'une suite $A_0=A, A_1, \dots, A_n=B$ de foncteurs et d'une suite
$u_1, \dots, u_n$ o\`u $u_i$ est une  transformation naturelle
soit de $A_{i-1}$ vers $A_i$ soit de $A_i$ vers $A_{i-1}$. Convenons
aussi
de dire qu'une transformation naturelle $u$ entre deux
foncteurs $F, G: C \rightarrow D$ est \`a
valeurs
dans la sous-cat\'egorie $W$ de $D$ si pour tout objet $X$ de $C$,
$u(X)$ est un morphisme de $W$ (de m\^eme pour une cha\^{\i}ne de
transformations
naturelles).

\begin{lemme}
\label{nonadjoint}
Soit $C$ et $C'$ deux
cat\'egories munies de  sous-cat\'egories $W$ et $W'$. Soient
$$
F: C\rightarrow C' ,\;\;\; G: C' \rightarrow C
$$
deux foncteurs avec $F(W)\subset W'$ et $G(W')\subset W$.
Supposons qu'il existe  des cha\^{\i}nes de transformations naturelles $u:
GF\leftrightarrow
1_{C}$
et $v: FG \leftrightarrow 1_{C'}$
\`a valeurs respectivement dans $W$ et $W'$.
Alors $F$ induit une \'equivalence de cat\'egories
simpliciales $$
L(C,W)\stackrel{\cong}{\rightarrow}L(C',W').
$$
\end{lemme}
{\em Preuve:}  L'existence d'une cha\^{\i}ne de transformations
naturelles \`a valeurs dans $W$
implique que le foncteur induit  par
$GF$ sur $L(C,W)$ est une \'equivalence: il r\'esulte alors des
Propositions 3.3 et
3.5 de \cite{DwyerKan2} que le morphisme $L(C,W)(X,Y)\rightarrow
L(C,W)(GFX,GFY)$ est
une \'equivalence, et il est facile de voir que $GF:L(C,W)\rightarrow
L(C,W)$ est
essentiellement surjectif. De m\`eme le foncteur induit par $FG$ est une
\'equivalence. Par cons\'equent $F$ et $G$ induisent des \'equivalences.
A la limite, cette derni\`ere \'etape peut \^etre vue en utilisant la
structure de cat\'egorie de mod\`eles ferm\'ee pour les $1$-pr\'ecats de Segal:
un foncteur est une \'equivalence si et seulement s'il induit un
isomorphisme dans la cat\'egorie homotopique, et, dans celle-ci, on peut
utiliser le fait que l'existence d'inverses \`a droite et \`a gauche implique
l'inversibilit\'e; enfin, pour conclure, on observe qu'un foncteur entre
cat\'egories
simpliciales est une \'equivalence si et seulement si le foncteur correspondant
entre $1$-cat\'egories de Segal est une \'equivalence.
\eop

L'observation suivante est utile pour traiter des sous-cat\'egories
d\'efinies par des conditions homotopiquement invariantes, par
exemple des cat\'egories de complexes \`a cohomologie support\'ee dans un
intervalle.

\begin{proposition}
\label{stabilite}
Soit $C$ une cat\'egorie et $W$ une sous-cat\'egorie de $C$ admettant un
calcul de
fractions homotopique
(cf \cite{GabrielZisman}).
Soit $B\subset C$ une sous-cat\'egorie pleine d\'efinie
par une propri\'et\'e $W$-invariante, i.e. l'image inverse d'un sous-ensemble
$B'$ de l'ensemble de classes d'\'equivalence d'objets de $L(C,W)$.
Alors $(B, W\cap B)$ admet un calcul de fractions homotopique et $L(B, W\cap
B)\rightarrow L(C,W)$ est homotopiquement pleinement fid\`ele, avec pour image
essentielle la sous-cat\'egorie pleine de $L(C,W)$ image inverse du m\^eme
sous-ensemble de l'ensemble des classes d'\'equivalence.
\end{proposition}
C'est une cons\'equence directe de la d\'efinition de Dwyer-Kan
(\cite{DwyerKan2} Definition 6.1).
\eop

L'application principale de la construction de Dwyer-Kan cf \cite{DwyerKan2}
fournit, \`a partir d'une
cat\'egorie de mod\`eles ferm\'ee $M$ et de sa
sous-cat\'egorie $W$ des \'equivalences faibles, une cat\'egorie
simpliciale $L(M,W)$.

\begin{convention}
\label{convenir}
Si $M$ est une cat\'egorie de mod\`eles ferm\'ee
(ou une sous-cat\'egorie d'ob\-jets fibrants, cofibrants etc.) on notera
$L(M)$ la localis\'ee de Dwyer-Kan par rapport \`a la sous-cat\'egorie des
\'equivalences faibles. Cette convention s'\'etendra \`a toute cat\'egorie dans
laquelle il y a une notion naturelle ``d'\'equivalence''. En cas de confusion
possible, e.g. quand il s'agit de comparer les \'equivalences pour la topologie
grossi\`ere et les $\Gg$-\'equivalences, on remettra $W$ dans la  notation par
exemple $L(nSePCh, W^{\rm gro})$ ou $L(nSePCh, W^{\Gg})$.
\end{convention}

Signalons maintenant quelques m\'ethodes de calcul des ``complexes de
fonctions'', i.e. des ensembles simpliciaux de morphismes dans $L(M)$ o\`u $M$
est une cat\'egorie de mod\`eles ferm\'ee.
Il y a bien s\^ur la d\'efinition directe de \cite{DwyerKan1} \`a l'aide d'une
r\'esolution de $M$ par des cat\'egories libres simpliciales;
nous renvoyons le
lecteur \`a \cite{DwyerKan1} pour plus de d\'etails.

La deuxi\`eme m\'ethode est celle des
``hamacs'', de Dwyer-Kan \cite{DwyerKan2}. Cette m\'ethode convient pour toute
localisation mais avec des hamacs de longueur arbitraire. Comme il est
remarqu\'e
dans \cite{DwyerKan2}, dans le cas d'une cat\'egorie de mod\`eles ferm\'ee, on
peut se restreindre aux hamacs de longueur $3$. Pour $x,y\in M$ on note
$$
{\rm ham}^3(M; x,y)
$$
la cat\'egorie de diagrammes de la forme
$$
x\leftarrow a\rightarrow b \leftarrow y
$$
o\`u la premi\`ere et la derni\`ere fl\`eche sont des \'equivalences faibles
i.e. dans $W$. Les morphismes dans ${\rm ham}^3(M; x,y)$
sont les diagrammes de la forme
$$
\begin{array}{ccccccc}
x&\leftarrow &a&\rightarrow &b& \leftarrow &y\\
\downarrow &&\downarrow &&\downarrow && \downarrow \\
x&\leftarrow &a&\rightarrow &b& \leftarrow &y
\end{array}
$$
o\`u les fl\`eches verticales aux deux bouts sont les identit\'es de $x$ et
de $y$.
Dwyer et Kan donnent dans \cite{DwyerKan2} une \'equivalence faible
d'ensembles simpliciaux naturelle entre le complexe de fonctions entre $x$ et
$y$, et le nerf de la cat\'egorie des hamacs:
$$
L(M)_{1/}(x,y) \cong \nu  {\rm ham}^3(M; x,y).
$$

La derni\`ere m\'ethode est celle des ``complexes de fonctions homotopiques''
(``homotopy function complexes''), de Reedy \cite{Reedy}, et de Dwyer,
Hirschhorn, Kan \cite{DHK} \cite{Hirschhorn} \cite{DwyerKan3}. Ceci est une
g\'en\'eralisation ou affaiblissement de la m\'ethode de Quillen  \cite{Quillen}
de consid\'erer une structure simpliciale sur la cat\'egorie de mod\`eles
ferm\'ee.
Dans la m\'ethode des ``function complexes'',
on n'exige pas autant de fonctorialit\'e, ce qui donne
plus de flexibilit\'e (au prix du fait, comme remarqu\'e dans \cite{DHK}, que la
composition des fonctions n'est plus directement accessible---mais ceci ne pose
pas de probl\`eme car on dispose d\'ej\`a de la cat\'egorie $L(M)$).

Soit $x,y\in M$. Une {\em r\'esolution simpliciale fibrante} de $y$ est un objet
simplicial de $M$, ${\bf y}\in M^{\Delta ^o}$ muni d'un morphisme
$c^{\ast}y\rightarrow {\bf y}$ o\`u $c^{\ast} y$ est l'objet constant \`a
valeurs $y$, tel que $y\rightarrow {\bf y}(p)$ soit une \'equivalence
faible pour tout $p\in \Delta$, et tel que ${\bf y}$ soit fibrant pour la {\em
structure de Reedy} de  $M^{\Delta ^o}$ (voir \cite{Reedy} \cite{Hirschhorn}
\cite{DHK} \cite{JardineGoerssBook} ainsi que notre discussion plus bas).
Dualement, une  {\em r\'esolution cosimpliciale cofibrante} de $x$ est un objet
cosimplicial de $M$, ${\bf x}\in M^{\Delta}$ muni d'un
morphisme ${\bf x}\rightarrow
c^{\ast}x$ tel que ${\bf x}(p)\rightarrow x$ soit une \'equivalence faible pour
tout $p\in \Delta$, et tel que ${\bf x}$ soit cofibrant pour la structure de
Reedy sur $M^{\Delta}$.

Avec ces r\'esolutions, on obtient (voir {\em loc cit.}) une
\'equivalence naturelle d'ensembles simpliciaux
$$
L(M)_{1/}(x,y) \cong d^{\ast} M_{1/}({\bf x}, {\bf y})
$$
o\`u le terme de droite est la diagonale de l'ensemble bisimplicial de
morphismes entre les deux r\'esolutions. Si $x$ est cofibrant, on a
$$
L(M)_{1/}(x,y)\cong M_{1/}(x, {\bf y})
$$
tandis que si $y$ est fibrant on a
$$
L(M)_{1/}(x,y)\cong M_{1/}({\bf x}, y).
$$
Ceci est une  g\'en\'eralisation directe de la notion de
calcul de complexes de fonctions en alg\`ebre homologique ( \`a l'aide
r\'esolutions
projectives et injectives de complexes). Ce cadre tr\`es satisfaisant est (\`a
notre connaissance) d\^u \`a  Reedy \cite{Reedy}, et \`a Dwyer,
Hirschhorn, Kan \cite{DHK} \cite{Hirschhorn} \cite{DwyerKan3}; et c'est ce \`a
quoi pensait Quillen en intitulant ``alg\`ebre homotopique'' son livre
\cite{Quillen}. Plus loin on retrouve les origines dans la th\'eorie des
hyper-recouvrements de
Verdier \cite{SGA4b}.

Si $M$ est une
cat\'egorie de mod\`eles ferm\'ee simpliciale (c'est le cas principal
envisag\'e par Quillen \cite{Quillen}), on obtient des
r\'esolutions \`a partir de la structure simpliciale, donc la localis\'ee
$L(M)$ est
\'equivalente \`a
la cat\'egorie simpliciale $M^{\rm spl}_{cf}$ des objets cofibrants
et fibrants.

Dans l'introduction de \cite{DwyerKan3}, Dwyer-Kan indiquent sans
entrer dans les d\'etails que $L(M)$ ``capture toute l'information  homotopique
sup\'erieure implicite dans $M$''   que Quillen \cite{Quillen} avait
recherch\'ee, et partiellement reconstitu\'ee dans la cat\'egorie homotopique de
$M$---pour le cas non simplicial. En fait on peut pr\'eciser leur affirmation.
En effet, le lemme suivant et son corollaire autorisent
les constructions de Gabriel-Zisman qu'on a rappel\'ee au \S 7 plus haut
et, dans le cas point\'e, celles-ci induisent bien
sur $Ho(M)=\tau _{\leq 1}
L(M,W)$
la ``structure triangul\'ee'' d\'efinie par Quillen.

\begin{lemme}
\label{prodfibhom}
Si $M$ est une cat\'egorie de mod\`eles ferm\'ee, alors $L(M)$ admet des
produits fibr\'es homotopiques et des coproduits homotopiques au sens du \S 7
ci-dessus.
\end{lemme}
{\em Preuve:}
On applique la m\'ethode des r\'esolutions rappel\'ee
ci-dessus. Soit
$$
x\rightarrow y\leftarrow z
$$
un diagramme entre objets  fibrants de $M$, avec les deux fl\`eches fibrantes
(tout diagramme pour un produit fibr\'e dans $L(M)$ est \'equivalent \`a un
diagramme de cette forme). Pour $u\in M$, on choisit une r\'esolution
cosimpliciale cofibrante ${\bf u}\rightarrow u$. Alors le carr\'e
$$
\begin{array}{ccc}
L(M)_{1/}(u, x\times _yz)& \rightarrow & L(M)_{1/}(u, z)\\
\downarrow && \downarrow \\
L(M)_{1/}(u, x)& \rightarrow & L(M)_{1/}(u, y)
\end{array}
$$
est \'equivalent au
carr\'e  cart\'esien d'ensembles simpliciaux
$$
\begin{array}{ccc}
M_{1/}({\bf u}, x\times _yz)& \rightarrow & M_{1/}({\bf u}, z)\\
\downarrow && \downarrow \\
M_{1/}({\bf u}, x)& \rightarrow & M_{1/}({\bf u}, y) .
\end{array}
$$
Les morphismes en bas et \`a droite sont des fibrations de Kan
(on peut le v\'erifier en utilisant le fait que ${\bf u}$ est cofibrant de
Reedy et que les morphismes $x\rightarrow y$
et $z\rightarrow y$ sont fibrants). Donc ce carr\'e est aussi
homotopiquement cart\'esien, donc le carr\'e des complexes de fonctions est
homotopiquement cart\'esien. Ceci montre que $x\times_yz$ est un produit fibr\'e
homotopique de $x$ et $z$ au-dessus de $y$.

La d\'emonstration pour les coproduits homotopiques est duale.
\eop

\begin{corollaire}
\label{gzpourM}
Si $M$ est une cat\'egorie de mod\`eles ferm\'ee,
la $2$-cat\'egorie stricte $\tau _{\leq 2} L(M)$ v\'erifie les
propri\'et\'es (A), (B), (C) de Gabriel-Zisman et (dans  le cas point\'e) leur
construction  des op\'erations $\Omega$ et $\Sigma$
s'applique, et donne donc la ``structure triangul\'ee'' sur $Ho\, M=\tau _{\leq
1} L(M)$.
\end{corollaire}
La v\'erification est laiss\'ee aux soins du lecteur.
\eop

Pour une version plus g\'en\'erale de ce lemme, qui permet de calculer les
limites index\'ees par une petite cat\'egorie, voir la proposition \ref{calclim}
ci-dessous (pour les colimites, voir \ref{calccolim}).

On revient maintenant \`a des consid\'erations g\'en\'erales sur la
localisation.
On va d\'efinir la localisation pour une $n$-cat\'egorie de Segal. Soit donc
$A$ une $n$-pr\'ecat de Segal et $f\in  A_{1^i}$ une $i$-fl\`eche, $1\leq i\leq
n$. On rappelle qu'avec les notations de \cite{limits} on a
$$
A_{1^i} = Mor _{nSePC}(h(1^i),A)
=Mor _{nSePC}(i\Upsilon (\ast ),A)
$$
o\`u $h(1^i)$ est la $n$-pr\'ecat de Segal repr\'esent\'ee par $1^i\in \Theta
^{n+1}$
et o\`u
$i\Upsilon$ d\'esigne l'op\'eration $\Upsilon$ de \cite{limits} it\'er\'ee $i$
fois. On a l'\'egalit\'e $\Upsilon (\ast )=  I$, o\`u $I$
est la $1$-cat\'egorie ayant deux objets
$0,1$ et un morphisme $0\rightarrow 1$.
On peut \'ecrire $i\Upsilon (\ast ) = (i-1)\Upsilon (I)$. Notre
$i$-fl\`eche $f$
correspond donc \`a un morphisme
$$
f: (i-1)\Upsilon (I)\rightarrow A.
$$
Soit $\overline{I}$ la $1$-cat\'egorie ayant deux
objets $0,1$ et un isomorphisme $0\cong 1$. On a $I\subset \overline{I}$,
et donc une cofibration
$$
(i-1)\Upsilon (I)\hookrightarrow
(i-1)\Upsilon (\overline{I})
$$
qu'on combine avec $f$ pour d\'efinir
$$
nSeL(A, f) := A \cup ^{ (i-1)\Upsilon (I)}(i-1)\Upsilon
(\overline{I}).
$$
C'est une $n$-pr\'ecat de Segal (m\^eme si $A$ \'etait une $n$-cat\'egorie de
Segal, ce n'est pas en g\'en\'eral une
$n$-cat\'egorie de Segal). Il faut appliquer l'op\'eration $SeCat$ (ou un
remplacement fibrant) pour obtenir une $n$-cat\'egorie de Segal
$SeCat(nSeL(A, f))$.

Si $B$ est une $n$-cat\'egorie de Segal fibrante et $u:A\rightarrow B$ un
morphisme, alors $u$ s'\'etend en un morphisme
$$
u':  nSeL(A, f) \rightarrow B
$$
si et seulement si l'image $u(f)$ est inversible
\`a homotopie pr\`es dans $B$, i.e. si elle est inversible comme $i$-fl\`eche
dans $\tau
_{\leq i}(B)$.

On peut pr\'eciser un peu plus cette propri\'et\'e universelle; mais il
nous faut
faire appel aux  notations de la section \S 11 ci-dessous. Le lecteur est
invit\'e
\`a lire le \S 11 (voir aussi \cite{nCAT} et
\cite{limits}) avant d'aborder ce qui suit
jusqu'\`a la fin de la d\'emonstration de la prochaine proposition (sigle
$\oslash$).

Pour $B$ fibrant, le morphisme
$$
\underline{Hom}(nSeL(A, f), B)\rightarrow
\underline{Hom}(A, B)
$$
est pleinement fid\`ele avec pour image essentielle la
sous-$n$-cat\'egorie de Segal pleine de $\underline{Hom}(A, B)$ dont les objets
sont les $u$ tels que $u(f)$ soit inversible \`a homotopie pr\`es. Ce r\'esultat
(pour la preuve voir la proposition ci-dessous) signifie que
$SeCat(nSeL(A, f))$ poss\`ede la propri\'et\'e universelle qui en fait une
localisation de $A$.

En particulier on note que si $f$ \'etait d\'ej\`a inversible \`a
homotopie pr\`es dans $SeCat(A)$, alors le morphisme
$A\rightarrow nSeL(A, f)$ est une \'equivalence faible.

On \'etend maintenant cette construction au cas des
sous-ensembles de $i$-fl\`eches. Soit $A$ une $n$-pr\'ecat de
Segal et soit $W=\{ W^i \}$ un syst\`eme de sous-ensembles de $i$-fl\`eches
$W^i\subset A_{1^i}$, $1\leq i \leq n$. On d\'efinit alors $nSeL(A,
W)$ comme le coproduit de $A$ avec un exemplaire de $(i-1)\Upsilon
(\overline{I})$ recoll\'e le long de  $(i-1)\Upsilon (I)$,  pour chaque $f\in
W^i$. Si on veut avoir une $n$-cat\'egorie de Segal on peut ensuite
regarder$SeCat(nSeL(A, W))$. On a la m\^eme propri\'et\'e universelle, qu'on
\'enonce sous la forme suivante.

\begin{proposition}
\label{propuniSeL}
Soit $A$ une $n$-cat\'egorie de Segal et
$W=\{ W^i \}$ un syst\`eme de sous-ensembles de $i$-fl\`eches
$W^i\subset A_{1^i}$. Alors $A\rightarrow nSeL(A, W)$ poss\`ede la
propri\'et\'e
universelle suivante. Si $B$ est une
$n$-cat\'egorie de Segal fibrante alors
$$
\underline{Hom}(nSeL(A, W), B)\rightarrow \underline{Hom}(A, B)
$$
est pleinement fid\`ele avec pour image essentielle
la sous-cat\'egorie pleine des $u: A\rightarrow B$ tels que, pour tout
$f$ dans l'un des $W^i$,
$u(f)$ soit inversible
\`a homotopie pr\`es.
\end{proposition}
{\em Preuve:}
Au vu des propri\'et\'es du $\underline{Hom}$ interne cf \ref{interne}
ci-dessous, il suffit de voir cette propri\'et\'e pour l'inclusion
$$
(i-1)\Upsilon (I) \rightarrow (i-1)\Upsilon (\overline{I}).
$$
Un morphisme
$$
(i-1)\Upsilon (I)
$$
est la m\^eme chose qu'un morphisme
$$
I\rightarrow B_{1^{i-1}/}.
$$
Comme $B_{1^{i-1}/}$ est fibrant, un tel morphisme s'\'etend \`a $\overline{I}$
si et seulement si la $i$-fl\`eche en question est inversible \`a \'equivalence
pr\`es. Il reste \`a prouver que le morphisme
$$
\underline{Hom}((i-1)\Upsilon (\overline{I}),B)_{1/}(f,g)
\rightarrow
\underline{Hom}((i-1)\Upsilon (I),B)_{1/}(f,g)
$$
est une \'equivalence; on va montrer que c'est une fibration triviale en
montrant la
propri\'et\'e de rel\`evement pour toute cofibration triviale
$E\rightarrow E'$.

Pour cela il faut montrer que le morphisme
$$
(i-1)\Upsilon (I) \times \Upsilon (E) \cup ^{
(i-1)\Upsilon (I)\times \{ 0,1\}}
(i-1)\Upsilon (\overline{I})\times \{ 0,1\}
\rightarrow
(i-1)\Upsilon (\overline{I})\times \Upsilon (E)
$$
est une cofibration triviale. L'argument est le m\^eme que dans la
d\'emonstration de (\cite{limits} Theorem 2.5.1) et on ne le reproduit pas ici.

Avec ce fait, on obtient que pour toute cofibration $E\rightarrow E'$,
la cofibration
$$
(i-1)\Upsilon (\overline{I})\times \Upsilon (E) \cup ^{
(i-1)\Upsilon (I)\times \Upsilon (E)}
(i-1)\Upsilon (I)\times \Upsilon (E')
$$
$$
\rightarrow (i-1)\Upsilon (\overline{I})\times \Upsilon (E')
$$
est triviale, ce qui donne la propri\'et\'e de rel\`evement voulue
(en notant qu'un morphisme
$$
E\rightarrow \underline{Hom}((i-1)\Upsilon (\overline{I}),B)_{1/}(f,g)
$$
s'identifie \`a un morphisme
$$
(i-1)\Upsilon (\overline{I})\times \Upsilon (E)\rightarrow B
$$
se restreignant en $f,g$ sur $(i-1)\Upsilon (\overline{I})\times \{ 0,1\}$
et de m\^eme pour les autres morphismes du carr\'e dans le diagramme
de rel\`evement).
\eop

\noindent
$\oslash$

Revenons maintenant au cas des $1$-cat\'egories de Segal. On a d\'efini une
localisation $1SeL(A,W)$ par tout sous-ensemble de $1$-fl\`eches $W\subset A_1$.
Dans ce cas on n'a pas besoin de la notation $\Upsilon$ puisqu' on
attache simplement
des exemplaires de $\overline{I}$ le long des morphismes $I\rightarrow A$
correspondant aux fl\`eches de $W$. Cette construction est compatible
avec
celle de Dwyer-Kan au sens suivant.

\begin{proposition}
\label{dklocsegloc}
Soit $C$ une cat\'egorie et $W\subset C$ une sous-cat\'egorie.
Alors il y a une
\'equivalence de $1$-cat\'egories de Segal (compatible avec
l'identit\'e de
$C$)
$$
L(C,W) \stackrel{\cong}{\rightarrow} 1SeL(C,W)'
$$
o\`u $1SeL(C,W)'$ est un remplacement fibrant de
$1SeL(C,W)$.
\end{proposition}
{\em Preuve:}
Rappelons (\cite{DwyerKan1}) qu'une {\em cat\'egorie libre} est une cat\'egorie
librement engendr\'ee par un ensemble de fl\`eches qu'on appellera les {\em
g\'en\'erateurs}. Si on fixe un ensemble $O$ d'objets, on a la notion de {\em
produit libre} d'une collection de $O$-cat\'egories
(i.e. de cat\'egories dont $O$ est
l'ensemble d'objets). Si on fixe un ensemble de g\'en\'erateurs $\{
f_{\alpha}\}$ avec $f_{\alpha}$ une fl\`eche de $x_{\alpha}\in O$ vers
$y_{\alpha}\in O$, notons par $\langle f_{\alpha}\rangle$ la $O$-cat\'egorie
libre engendr\'ee par $f_{\alpha}$. La $O$-cat\'egorie libre engendr\'ee par
$\{ f_{\alpha}\}$ est alors le produit libre des $O$-cat\'egories
$\langle f_{\alpha}\rangle$.

On va prouver d'abord le r\'esultat suivant: soit $C$ une cat\'egorie libre
(avec ensemble $O$ d'objets) et
soient $x,y\in O$. On va rajouter une fl\`eche $f$ de source
$x$ et but $y$. Soit
$$
C':= C\ast \langle f\rangle
$$
la $O$-cat\'egorie libre engendr\'ee par les g\'en\'erateurs de $C$ plus $f$.
D'autre part la fl\`eche $f$ de $C'$ correspond \`a un morphisme
$I\rightarrow C'$. On va prouver que le morphisme de $1$-pr\'ecats de Segal
$$
C\cup ^{\{ x,y\} }I\rightarrow C'
$$
est une \'equivalence faible.

Pour prouver ce r\'esultat, on montre que $C'$
s'obtient \`a partir de $C\cup ^{\{ x,y\} } I$ par une suite de cofibrations
triviales explicites. Pour tout $\ell$, soit $F_{\ell}C'$ le sous-ensemble
simplicial de $C'$ d\'etermin\'e par les simplexes comportant au total,
dans toutes
leurs ar\^etes principales, au plus  $\ell$ des g\'en\'erateurs de $C'$ (i.e.
des g\'en\'erateurs de $C$ plus $f$). Ici on compte, sur une m\^eme ar\^ete (qui
correspond \`a une fl\`eche de $C'$), le nombre des g\'en\'erateurs (avec
leur multiplicit\'e) qui entre en jeu dans la fl\`eche en question. Soit
$F_{\ell}C$ l'intersection de $F_{\ell}C'$ avec $C$,  c'est-\`a-dire le
sous-ensemble simplicial de $C$ des simplexes comportant au plus
$\ell$ g\'en\'erateurs de $C$ au sens ci-dessus. On a l'\'egalit\'e
$$
F_1C' = F_1C \cup ^{\{ x,y\} } I.
$$
On va montrer que, pour $\ell \geq 2$, $F_{\ell}C'$ est obtenu  \`a partir de
$F_{\ell -1}C'\cup ^{F_{\ell -1}C}F_{\ell}C$ par coproduit avec une
collection de cofibrations triviales.  Pour cela,  on appellera {\em
$\ell$-simplexe
\'el\'ementaire} tout \'el\'ement  $u\in C'_{\ell}$ dont les ar\^etes
principales sont des g\'en\'erateurs de $C'$. Rappelons aussi (cf.
\cite{nCAT}) que $h(\ell )$ d\'esigne l'ensemble simplicial repr\'esent\'e
par $\ell$
(consid\'er\'e comme $1$-pr\'ecat de Segal constant en la deuxi\`eme variable
simpliciale).
Le morphisme
$$
i: h(\ell -1) \cup ^{h(\ell -2)}h(\ell -1)\rightarrow h(\ell )
$$
correspondant \`a l'inclusion de la premi\`ere et de la derni\`ere face, est une
cofibration triviale de $1$-pr\'ecats de Segal.
\footnote{
Il est facile de voir par r\'ecurrence que le morphisme
$$
h(m -1) \cup ^{\ast} h(1) \rightarrow h(m )
$$
est une cofibration triviale, et on voit que $i$ est une cofibration
triviale
en combinant ce r\'esultat pour $m=\ell$ et $m=\ell -1$
 avec la propri\'et\'e ``trois pour le prix de deux''.}
Si $u$ est un $\ell$-simplexe \'el\'ementaire vu comme
morphisme
de $h(\ell )$ vers $C'$, on a soit $u: h(\ell )\rightarrow C$, soit
$$
u^{-1} (F_{\ell -1}C' \cup ^{F_{\ell -1}C} C) =
i\left( h(\ell -1) \cup ^{h(\ell -2)}h(\ell -1)\right) .
$$
Dans le premier cas on n'a pas besoin d'ajouter $u$; dans le deuxi\`eme cas on
peut ajouter \`a $F_{\ell -1} C'\cup ^{F_{\ell -1}C} F_{\ell}C$ un exemplaire de
$h(\ell )$ recoll\'e le long de  $ h(\ell -1) \cup ^{h(\ell -2)}h(\ell -1)$,
correspondant \`a $u$. Appelons provisoirement $B$ le r\'esultat de tous ces
coproduits, un pour chaque $\ell$-simplexe \'el\'ementaire non contenu
dans $C$. Il y a un morphisme \'evident $B\rightarrow C'$ et il est facile
de voir
que ce morphisme est injectif, car un \'el\'ement (i.e. $k$-simplexe)
introduit
par la cofibration associ\'ee \`a $u$ connait tous les g\'en\'erateurs
apparaissant
dans $u$, en particulier il ne peut provenir d'une
cofibration associ\'ee \`a un $u'\neq u$. D'autre part l'image de ce
morphisme est
visiblement $F_{\ell} C'$. Donc $B=F_{\ell} C'$ et on a prouv\'e que
$$
F_{\ell -1} C'\cup ^{F_{\ell -1}C} F_{\ell}C\rightarrow
F_{\ell} C'
$$
est une cofibration triviale. L'argument montre aussi que
$$
F_{\ell -1} C'\cup ^{F_{\ell -1}C} C\rightarrow
F_{\ell } C'\cup ^{F_{\ell }C} C
$$
est une cofibration triviale. La composition de ces
cofibrations triviales donne la cofibration triviale
$$
F_1C'\cup ^{F_1C} C\rightarrow C',
$$
c'est-\`a-dire le r\'esultat cherch\'e.

Maintenant on continue la d\'emonstration de la proposition. Notons
$A\ast B$ le produit libre de deux $O$-cat\'egories, et $A\cup ^OB$ le
coproduit en tant que $O$-pr\'ecats de Segal (i.e. $1$-pr\'ecats de Segal
ayant $O$
comme ensemble d'objets). Le r\'esultat ci-dessus implique, par
r\'ecurrence (plus
\'eventuellement passage \`a une colimite filtrante qui ne pose pas de
probl\`eme), que si $C$ est la cat\'egorie libre engendr\'ee par $\{
f_{\alpha}\}$
alors
$$
\coprod _{\alpha}^{O} \langle f_{\alpha} \rangle \rightarrow C
$$
est une \'equivalence faible de $1$-pr\'ecats de Segal. Il s'ensuit que si
$A$ et
$B$ sont des $O$-cat\'egories libres, alors le morphisme
$$
A\cup ^OB\rightarrow A\ast B
$$
est une \'equivalence faible de $1$-pr\'ecats de Segal.

Si $A$ est une $O$-cat\'egorie simpliciale on peut la consid\'erer comme
$1$-cat\'egorie de Segal (en utilisant la variable simpliciale comme deuxi\`eme
variable simpliciale pour la cat\'egorie de Segal). Le r\'esultat ci-dessus
implique directement que si $A$ et $B$ sont des $O$-cat\'egories simpliciales,
libres \`a chaque \'etage, alors le morphisme
$$
A\cup ^OB\rightarrow A\ast B
$$
est une \'equivalence faible. Maintenant, le coproduit (dans une cmf) d'objets
cofibrants (et pour \ref{SeCmf} tout objet est cofibrant) pr\'eserve le type
d'\'equivalence. Dwyer-Kan montrent que le produit libre des
$O$-cat\'egories simpliciales
quelconques preserve l'\'equivalence
(\cite{DwyerKan1} Proposition 2.7). En utilisant des r\'esolutions simpliciales
libres (\cite{DwyerKan1}) on obtient que pour toute paire de cat\'egories
simpliciales $A,B$ ayant l'ensemble $O$ pour objets,
le morphisme
$$
A\cup ^OB\rightarrow A\ast B
$$
est une \'equivalence faible de $1$-pr\'ecats de Segal.

Revenons maintenant \`a la d\'efinition de la localisation de Dwyer-Kan
\cite{DwyerKan1}. Si $C$ est une cat\'egorie et $W$ une sous-cat\'egorie,
Dwyer-Kan choisissent une r\'esolution libre de $C$ de la forme
$A\ast B$ avec $A$ et $B$ des cat\'egories simpliciales libres \`a
chaque niveau
(et partageant le m\^eme ensemble d'objets $O$),
et de sorte que $B$ soit une r\'esolution de $W$. Notons $\overline{B}$
la cat\'egorie simpliciale obtenue en inversant les fl\`eches de chaque
niveau
de $B$. Par d\'efinition on a
$$
L(C,W):= A\ast \overline{B}
$$
tandis que
$$
1SeL(A\ast B, B)= (A\ast B)\cup ^{B}1SeL(B,B).
$$
Observons en outre que $1SeL(B,B)\rightarrow \overline{B}$ est une \'equivalence
faible (pour cela on peut supposer que $B$ est une cat\'egorie libre, auquel cas
$\overline{B}$ est son compl\'et\'e $0$-groupique et $1SeL(B,B)$ aussi).

On
a une suite de fl\`eches $$
A\cup ^O \overline{B} = (A\cup ^OB)\cup ^B\overline{B}\rightarrow
(A\ast B)\cup ^{B}\overline{B}\rightarrow A\ast \overline{B},
$$
et d'apr\`es les remarques ci-dessus,
la premi\`ere fl\`eche et la compos\'ee sont
des \'equivalences faibles de $1$-pr\'ecats de Segal. D'apr\`es
``trois pour le prix de
deux'' on obtient que la fl\`eche
$$
1SeL(A\ast B, B)= (A\ast B)\cup ^{B}\overline{B}\rightarrow A\ast
\overline{B}=L(C,W)
$$
est une \'equivalence. D'autre part la localisation $1SeL$
pr\'eserve les \'equivalences donc la fl\`eche
$$
1SeL(A\ast B, B)\rightarrow 1SeL(C,W)
$$
est une \'equivalence faible donc, en choisissant un remplacement fibrant
$1SeL(L,W)'$ on obtient une \'equivalence (essentiellement bien d\'efinie)
de $1$-cat\'egories de Segal
$$
L(C,W)\stackrel{\cong}{\rightarrow} 1SeL(C,W)'.
$$
Ceci termine la d\'emonstration.
\eop

\begin{corollaire}
\label{dklocseglocspl}
Soit $C$ une cat\'egorie simpliciale et $W\subset C$ une sous-cat\'egorie
simpliciale. On peut supposer que $W$ est satur\'ee au sens que $W(x,y)$ est
une r\'eunion de composantes connexes de $C(x,y)$ (cf \cite{DwyerKan1}).
L'ensemble $W_{1,0}$ des sommets
de l'ensemble simplicial des fl\^eches de $W$ est un sous-ensemble de
$1$-fl\`eches de la $1$-cat\'egorie de Segal correspondante \`a $C$. On a alors
une
\'equivalence de $1$-cat\'egories de Segal (compatible avec les injections de
$C$)
$$
L(C,W) \stackrel{\cong}{\rightarrow} 1SeL(C,W_{1,0})'.
$$
\end{corollaire}
{\em Preuve:}
On suppose $C$ libre comme dans \cite{DwyerKan1}, et on applique la
proposition \`a
chaque niveau. On peut noter que d'apr\`es \cite{DwyerKan1}, $L(C,W)$ est
\'equivalent \`a $L(C, W_{1,0})$.
\eop

\begin{corollaire}
\label{propunivM}
Soit $C$ une cat\'egorie simpliciale et $W\subset C$ une sous-cat\'egorie
simpliciale. Si $B$ est une $1$-cat\'egorie de Segal fibrante, alors
le morphisme de $1$-cat\'egories de Segal
$$
\underline{Hom}(L(C,W), B)\rightarrow \underline{Hom}(C,B)
$$
est pleinement fid\`ele avec pour image essentielle la
sous-cat\'egorie des morphismes $C\rightarrow B$ qui envoient
les \'el\'ements de $W$ sur des
morphismes inversibles (\`a \'equivalence pr\`es) dans $B$.
\end{corollaire}
{\em Preuve:}
Il suffit
d'appliquer la proposition \ref{propuniSeL} et \ref{dklocseglocspl} (ou juste
\ref{dklocsegloc} s'il s'agit d'une $1$-cat\'egorie $C$ non simpliciale).
\eop

Ce corollaire devrait avoir une variante formul\'ee en termes de
cat\'egories simpliciales avec les cat\'egories simpliciales
$Coh(\cdot , \cdot )$ de  morphismes de Cordier et Porter, voir les probl\`emes
\ref{probleme1}--\ref{probleme2} ci-dessus.

{\em Remarque:} En toute rigueur, il n'existe pas forc\'ement de morphisme
$C\rightarrow L(C,W)$. En revanche, si $\widehat{C}\rightarrow C$
est la r\'esolution
simpliciale standard de $C$ (qui est une \'equivalence de cat\'egories
simpliciales) on a par construction \cite{DwyerKan1} un morphisme
$\widehat{C}\rightarrow L(C,W)$. Nous allons, dans le reste du papier,
n\'egliger
ce point technique (on l'a d\'ej\`a ignor\'e dans l'\'enonc\'e du corollaire
pr\'ec\'edent) et {\em pr\'etendre} qu'on a un morphisme $C\rightarrow L(C,W)$.
Cette difficult\'e ne concerne pas $1SeL(C,W)$
qui re\c{c}oit par construction
un morphisme de source $C$.

\bigskip

Le corollaire suivant est l'analogue pour les cat\'egories de Segal ayant
plusieurs objets, d'un fait bien connu pour les mono\"{\i}des
topologiques: tout mono\"{\i}de faible (i.e. objet de la cat\'egorie
sous-jacente \`a
une machine de d\'ela\c{c}age; par exemple une cat\'egorie de Segal avec un seul
objet, pour la machine de Segal) est \'equivalent \`a un mono\"{\i}de
topologique
strict. Voir le travail de
Fiedorowicz \cite{Fiedorowicz}, ou la proposition 1.12 de Dunn
\cite{Dunn}.

En fait ce corollaire se trouve 
dans \cite{DKS}, avec pour la premi\`ere fois
la notion de $1$-cat\'egorie de Segal (la pr\'esente remarque 
a \'et\'e ajout\'ee dans la version 2).

\begin{corollaire}
\label{catsimpl}
Si $A$ est une $1$-cat\'egorie de Segal fibrante alors il existe une cat\'egorie
simpliciale $C$ et une \'equivalence $C\rightarrow A$.
\end{corollaire}
{\em Esquisse de d\'emonstration:}
En consid\'erant $A$ comme un pr\'efaisceau simplicial sur la cat\'egorie
$\Delta$, on peut appliquer la ``subdivision barycentrique''
au-dessus de chaque objet de $\Delta$ (voir le lemme \ref{bd} ci-dessous) pour
obtenir un pr\'efaisceau de cat\'egories sur $\Delta$. Ensuite on
applique l'op\'eration d'int\'egration de Grothendieck (cf \S 16 ci-dessous), en
utilisant un analogue du lemme \ref{bd} pour les pr\'efaisceaux de cat\'egories
au-dessus de $\Delta$ au lieu des ensembles simpliciaux.  Ceci est l'analogue
pour les cat\'egories de Segal, de
la platification (``flattening'') utilis\'ee dans
Dwyer-Kan (\cite{DwyerKanDiags} ``A delocalization theorem'' 2.5). On obtient
ainsi une $1$-cat\'egorie $\beta (A)$ munie d'un sous-ensemble $\delta (A)$ de
$1$-fl\`eches: on y inclut toutes les fl\`eches ``verticales'' i.e.
des cat\'egories-fibres, ainsi que les fl\`eches ``horizontales'' qui
entrent dans le $\delta (A)$ du \ref{bd}. Comme dans \ref{bd},
on a
$$
A\cong 1SeL(\beta (A), \delta (A))'
$$
(\'equivalence dans la cat\'egorie homotopique de $1SeCat$).
Alors, d'apr\`es la proposition pr\'ec\'edente, on a aussi
$A \cong L(\beta (A), \delta (A))$. Et comme $A$ est fibrante,
cette \'equivalence peut \^etre r\'ealis\'ee comme une \'equivalence
$L(\beta (A), \delta (A))\rightarrow A$.
\eop

Cette d\'emonstration montre aussi qu'il existe un ensemble simplicial $X$
v\'erifiant
$$
SeCat(X) \cong A.
$$
En effet, $\beta (A)$ est une $1$-cat\'egorie qui
correspond \`a un ensemble simplicial (son nerf); et la localisation
$1SeL(\beta (A), W)$ s'obtient en recollant des exemplaires de
$\overline{I}$ le long des morphismes $I\rightarrow \beta (A)$ de $W$.
Le coproduit reste dans les ensembles simpliciaux (en effet, le
r\'esultat reste
constant dans la deuxi\`eme variable simpliciale des $1$-pr\'ecats de Segal).
On a donc bien un ensemble simplicial
$$
X:= 1SeL(\beta (A), W)
$$
v\'erifiant $A\cong SeCat(X)$.

\subnumero{Adjoints homotopiques}

La propri\'et\'e d'adjonction se transmet aux localis\'es de Dwyer-Kan
et devient une propri\'et\'e d'adjonction homotopique.  Soient $A$ et $B$ deux
$n$-cat\'egories de Segal (qu'on peut supposer fibrantes), et
$$
F: A\rightarrow B,\;\;\;\; G: B\rightarrow A
$$
deux morphismes. Soit aussi
$$
u\in \underline{Hom}(A,A)_{1/}(1_A, GF)
$$
une transformation naturelle. Alors
pour tout $x\in A_0$ et tout $y\in B_0$, $u$ induit un morphisme
(essentiellement
bien d\'efini)
$$
B_{1/}(Fx,y)\rightarrow A_{1/}(x,Gy).
$$
On dira que $F$ et $G$ sont {\em homotopiquement adjoints via $u$} si
ce morphisme est une \'equivalence de $n-1$-cat\'egories de Segal pour tout
$x$ et tout $y$. On dira que $F$ est {\em l'adjoint \`a gauche} et $G$ {\em
l'adjoint \`a
droite}. On pourrait aussi construire la d\'efinition \`a partir
d'une transformation
naturelle
$$
v\in \underline{Hom}(B,B)_{1/}(FG,1_B).
$$
Et pour bien faire, il faudrait v\'erifier
que les deux d\'efinitions co\"{\i}ncident, que l'adjoint
$(G,u)$ d'un morphisme $F$ est essentiellement unique s'il existe, etc.

Voir Cordier-Porter \cite{CordierPorter} dans le cas $n=1$.

De m\^eme on pourrait envisager de d\'efinir la notion de paire adjointe de
$1$-fl\`eches dans une $n$-cat\'egorie de Segal quelconque (et r\'ecup\'erer la
notion ci-dessus avec la $n+1$-cat\'egorie $nSeCAT$). Ceci pourrait correspondre
\`a ce
que Baez-Dolan appellent  une ``duale'' dans une $n$-cat\'egorie
\cite{BaezDolan}.

\begin{lemme}
\label{adjointh}
Soit $F:M\rightarrow N$ un foncteur de Quillen \`a gauche entre cat\'egories de
mod\`eles ferm\'ees, et soit $G$ son adjoint \`a droite.
Alors
$$
L(F): L(M)\leftrightarrow L(N): L(G)
$$
sont homotopiquement adjoints.
\end{lemme}
{\em Preuve:}
D'abord on doit noter que le calcul des $L(M)_{1/}(x,y)$ par r\'esolutions est
compatible avec la composition dans la mesure du possible: si par exemple
$y\rightarrow {\bf y}$ est une r\'esolution simpliciale fibrante, et si
$x\rightarrow x'$ est un morphisme entre objets cofibrants, alors
le morphisme induit
$$
M_{1/}(x',{\bf y})\rightarrow M_{1/}(x,{\bf y})
$$
co\"{\i}ncide avec le morphisme de composition
$$
L(M)_{1/}(x',y)\rightarrow L(M)_{1/}(x,y)
$$
via l'\'equivalence de \cite{DwyerKan2} qui a \'et\'e rappel\'ee plus haut.

On applique ceci au morphisme d'adjonction $x\rightarrow GFx$.
On obtient que si $x$ est cofibrant et si $y\rightarrow {\bf y}$ est une
r\'esolution simpliciale fibrante, alors le morphisme
$$
L(N)_{1/}(Fx,y)\rightarrow L(M)_{1/}(x,Gy)
$$
induit par le morphisme d'adjonction, co\"{\i}ncide (\`a \'equivalence pr\`es)
avec le morphisme
$$
N_{1/}(Fx, {\bf y})\rightarrow M_{1/}(x,G{\bf y})
$$
(ici le fait que $F$ et $G$ sont des foncteurs de Quillen (\`a gauche et \`a
droite) implique que $Fx$ est cofibrant, et que $Gy\rightarrow G{\bf y}$ est une
r\'esolution simpliciale fibrante). Or ce dernier morphisme est un isomorphisme,
donc le morphisme d'adjonction sur les localis\'es est une \'equivalence.
\eop

Pour une autre propri\'et\'e d'une cat\'egorie de mod\`ele ferm\'ee
$M$ qui peut
\^etre ``amplifi\'ee''
en propri\'et\'e de $L(M)$ voir la proposition
\ref{hominterne} sur les $Hom$ internes homotopiques.
On termine en posant le probl\`eme analogue pour les ``foncteurs
d\'eriv\'es''.

\begin{probleme}
\label{derive}
D\'efinir la notion de {\em foncteur d\'eriv\'e homotopique} dans le cadre des
cat\'egories simpliciales ou $1$-cat\'egories de Segal (ou \'eventuellement les
$n$-cat\'egories de Segal). Prouver que si $M$ et $N$ sont des cat\'egories de
mod\`eles ferm\'ees et si $F: M\rightarrow N$ est un foncteur de Quillen (\`a
gauche, disons), alors le foncteur compos\'e
$$
L(M)\cong L(M_c)\rightarrow L(N_c)\cong L(N)
$$
est le foncteur d\'eriv\'e homotopique de $M\rightarrow L(N)$ le long de
$M\rightarrow L(M)$. Donner un \'enonc\'e d'unicit\'e essentielle pour le
foncteur d\'eriv\'e homotopique, s'il existe.
\end{probleme}

La solution de
ce probl\`eme permettrait de voir que le $1$-pr\'echamp de Segal
$L(\mm )$ associ\'e \`a un pr\'efaisceau de Quillen \`a gauche $\mm$ (voir
\S 17)
est essentiellement bien d\'efini ind\'ependamment du choix des
sous-cat\'egories
de cofibrations dans les $\mm (y)$.

\numero{Premi\`ere d\'efinition de champ}
\label{notionchamppage}

On consid\`ere \`a nouveau la structure de cat\'egorie de
mod\`eles ferm\'ee du th\'eor\`eme \ref{cmf}, relatif \`a un site $\Xx$ avec
topologie de Grothendieck $\Gg$.

Soit $A$ un $n$-pr\'echamp de Segal.
Soit $A\rightarrow A'$ une cofibration triviale vers un $n$-pr\'echamp de Segal
$\Gg$-fibrant.  On dira que $A$ est un {\em $n$-champ de Segal} (ou bien un {\em
$\infty$-champ} ou m\^eme juste un {\em champ}) si: \newline
---$A(X)$ est une $n$-cat\'egorie de Segal pour tout $X\in \Xx$; et
\newline
---le morphisme $A(X)\rightarrow A'(X)$ est une \'equivalence de
$n$-cat\'egories de Segal pour tout $X\in \Xx$.

On voit facilement (voir l'argument du corollaire \ref{unique} ci-dessous) que
cette condition est ind\'ependante du choix de $A\rightarrow A'$.

Pour $n=0$ les pr\'efaisceaux simpliciaux qui sont des $0$-champs de Segal sont
pr\'e\-ci\-s\'e\-ment ceux
qui v\'erifient
la condition de {\em descente cohomologique} de Thomason (voir
\cite{TravauxThomason}), ou encore ceux qui sont {\em
flasques
par rapport \`a tout objet $X$ du site $\Xx$} au sens de Jardine \cite{Jardine}
(ce sont aussi les {\em faisceaux homotopiques} de
\cite{kobe}
ou les {\em faisceaux flexibles} de \cite{flexible}). La terminologie
``descente cohomologique''
de Thomason est inspir\'ee de la descente cohomologique de Deligne-Saint Donat \cite{SGA4b}
et \cite{HodgeIII}.

Il est int\'eressant de constater que les lemmes suivants, analogues des
r\'esultats classiques sur les faisceaux, sont des cons\'equences enti\`erement
formelles de la th\'eorie de Quillen \cite{Quillen}
(d'ailleurs on peut
d\'erouler la th\'eorie des faisceaux d'ensembles sur un site, comme application
essentiellement facile de ces arguments standards en th\'eorie des
cmf, voir le travail de Rezk \cite{Rezk}).

\begin{lemme}
\label{morphentrechamps}
Si $f:A\rightarrow B$ est un morphisme de $n$-champs de Segal qui est une
$\Gg$-\'equivalence faible alors $f$ est aussi une \'equivalence faible
objet-par-objet (i.e. pour la topologie grossi\`ere).
\end{lemme}
{\em Preuve:}
Par d\'efinition des $n$-champs de Segal on peut supposer que $A$ et $B$ sont
$\Gg$-fibrants. On peut factoriser le morphisme donn\'e en
$$
A\rightarrow B' \rightarrow B
$$
o\`u le premier morphisme est une cofibration $\Gg$-triviale et le deuxi\`eme
morphisme est une fibration $\Gg$-triviale. Une fibration $\Gg$-triviale
poss\`ede la propri\'et\'e de rel\`evement vis-\`a-vis de toutes les
cofibrations; en
particulier elle est aussi triviale pour la topologie grossi\`ere. D'autre part,
comme $A$ est $\Gg$-fibrant, il existe une r\'etraction $B'\rightarrow A$
(induisant l'identit\'e sur $A$). En particulier le morphisme $A\rightarrow B'$
admet un inverse \`a gauche. On a ainsi montr\'e que tout morphisme
$f:A\rightarrow B$
entre $n$-champs de Segal (pour $\Gg$) admet un inverse $g$ \`a gauche,
dans la
cat\'egorie homotopique des $n$-pr\'echamps de Segal pour l'\'equivalence de la
topologie grossi\`ere. Ceci s'applique en
particulier \`a cet inverse $g:B\rightarrow A$ (c'est ici
qu'on voit pourquoi on a suppos\'e que $B$ est un
champ). On obtient ainsi un inverse \`a gauche $h$ de $g$; or $f$ est inverse
\`a droite de $g$ donc $g$ est inversible, ce qui implique que $f$ est
inversible. Ceci permet de conclure puisque, d'apr\`es Quillen
\cite{Quillen} les seuls morphismes
dont la classe soit inversible dans la cat\'egorie homotopique
sont les
\'equivalences (ici, pour la topologie grossi\`ere).
\eop

\begin{lemme}
\label{grofibchampimplfib}
Si $A$ est un $n$-pr\'echamp de Segal alors $A$ est $\Gg$-fibrant si et
seulement
si $A$ est un champ, et fibrant pour la topologie grossi\`ere.
\end{lemme}
{\em Preuve:}
Si $A$ est $\Gg$-fibrant alors $A$ est
fibrant pour la topologie grossi\`ere, et par d\'efinition c'est un champ.

Supposons donc que $A$ est fibrant pour la topologie grossi\`ere et
que c'est un champ.
Soit $A\hookrightarrow A'$ un remplacement $\Gg$-fibrant.
Par d\'efinition c'est une \'equivalence faible pour la topologie grossi\`ere.
Comme $A$ est fibrant pour la topologie grossi\`ere, il existe une r\'etraction
$r:A'\rightarrow A$. Montrons alors que $A$ poss\`ede la propri\'et\'e de
rel\`evement vis-\`a-vis de toute cofibration $\Gg$-triviale $X\rightarrow Y$.
Comme $A$ est $\Gg$-fibrant, tout morphisme $X\rightarrow A$ s'\'etend
en un morphisme $Y\rightarrow A'$. En composant avec la r\'etraction $r$
on obtient l'extension $Y\rightarrow A$ d\'esir\'ee.
\eop

{\em Remarque:} On a un r\'esultat similaire (et m\^eme plus joli) pour la
structure de HBKQ (Th\'eor\`eme \ref{cmfHirschho}): si $A$ est un
$n$-pr\'echamp de Segal alors $A$ est $\Gg$-fibrant de HBKQ, si et
seulement
si les $A(X)$ sont fibrants pour tout $X$, et $A$ est un champ. En particulier,
on obtient la m\^eme notion de ``champ'' en utilisant
\ref{cmfHirschho} au lieu de \ref{cmf}.

\begin{lemme}
\label{stuff}
Si $f: A\rightarrow B$ est une
$\Gg$-fibration et si $B$ est un $n$-champ de Segal, alors $A$ est un
$n$-champ de Segal.

Si $f: A\rightarrow B$ est un morphisme entre $n$-champs
de Segal et si $f$ est fibrant pour la topologie grossi\`ere alors
$f$ est $\Gg$-fibrant.
\end{lemme}
{\em Preuve:}
Pour la premi\`ere partie, on commence par le cas ponctuel, pour lequel on va
seulement donner une esquisse de d\'emonstration:
si $f:A\rightarrow B$ est une fibration de $n$-pr\'ecats de Segal, et si $B$ est
une $n$-cat\'egorie de Segal, alors $A$ est une $n$-cat\'egorie de Segal.
On remarque d'abord que les morphismes
$A_{p/}(x_0,\ldots , x_p)\rightarrow B_{p/}(fx_0,\ldots ,fx_p)$ sont des
fibrations de $n-1$-pr\'ecats de Segal, donc par r\'ecurrence sur $n$ on peut
supposer que les $A_{p/}$ sont des $n-1$-cat\'egories de Segal. Maintenant
on consid\'ere le diagramme $$
\begin{array}{ccc}
A{p/}&\rightarrow &A_{1/}\times _{A_0}\ldots , \times _{A_0}A_{1/}\\
\downarrow && \downarrow \\
B{p/}&\rightarrow &B_{1/}\times _{B_0}\ldots , \times _{B_0}B_{1/}.
\end{array}
$$
Soit $E\rightarrow A_{1/}\times _{A_0}\ldots , \times _{A_0}A_{1/}$
un morphisme de $n-1$-pr\'ecats de Segal,
ce qui correspond \`a un morphisme $U^p(E)\rightarrow A$
(o\`u $U^p(\, )$ est le focteur adjoint \`a gauche de
$A\mapsto A_{1/}\times _{A_0}\ldots , \times
_{A_0}A_{1/}$, qui associe donc \`a toute $n-1$-pr\'ecat de Segal
une $n$-pr\'ecat
de Segal). On a une inclusion $U^p(E)\rightarrow [p]E$ o\`u $[p]$ est
adjoint \`a gauche de $A\mapsto A_{p/}$.
Cette inclusion est une cofibration triviale
de $n$-pr\'ecats de Segal. Maintenant, le fait que $B$ soit une $n$-cat\'egorie
de Segal implique qu'il existe un prolongement de $U^p(E)\rightarrow B$
en un morphisme $([p]E)'\rightarrow B$ o\`u $([p]E)'$ est \`equivalent \`a
$[p]E$ objet-par-objet au-dessus de $\Delta$. Le fait que $A\rightarrow B$
soit fibrant (combin\'e au fait que $U^p(E)\rightarrow ([p]E)'$ soit une
cofibration triviale) implique qu'il existe \'egalement une extension
$([p]E)'\rightarrow A$. Ceci implique que le morphisme $E\rightarrow
A_{1/}\times _{A_0}\ldots , \times _{A_0}A_{1/}$ se rel\`eve \`a homotopie
pr\`es en un morphisme $E\rightarrow A_{p/}$. On pourrait donner aussi une
version relative de cette discussion: \'etant donn\'e $F\rightarrow A_{p/}$,
$F\hookrightarrow E$ et $E\rightarrow A_{1/}\times _{A_0}\ldots , \times
_{A_0}A_{1/}$, il existe \`a homotopie pr\`es une extension en $E\rightarrow
A_{p/}$. On obtient que l'application de Segal est une \'equivalence, donc $A$
est une $n$-cat\'egorie de Segal.

On traite maintenant la premi\`ere partie pour $f:A\rightarrow B$ au-dessus de
$\Xx$. Soit $B\hookrightarrow B'$ un remplacement $\Gg$-fibrant
(qui est une \'equivalence objet-par-objet), et soit
$$
A\hookrightarrow A'\rightarrow B'
$$
une factorisation o\`u la deuxi\`eme fl\`eche est une fibration pour la
topologie
grossi\`ere et la premi\`ere est une \'equivalence objet-par-objet.
On va montrer que la deuxi\`eme fl\`eche est une $\Gg$-fibration. Soit
$E\rightarrow
E'$ une cofibration $\Gg$-triviale avec un diagramme
$$
\begin{array}{ccc}
E&\rightarrow & A'\\
\downarrow && \downarrow \\
E'&\rightarrow & B'.
\end{array}
$$
On peut supposer que $E\rightarrow A'$ et $E'\rightarrow B'$ sont des fibrations
pour la topologie grossi\`ere. Avec cette reduction, on a par propr\'et\'e
cf \ref{proper} que
$E\times _{A'}A\rightarrow E$ et $E'\times _{B'}B\rightarrow E'$ sont des
\'equivalences objet-par-objet (on utilise le fait que les valeurs  $B(X)$ et
$B'(X)$ sont des $n$-cat\'egories de Segal; de m\^eme pour $A(X)$, et $A'(X)$,
pour cela on utilise la version ponctuelle mentionn\'ee au d\'ebut).
D'autre part,
on  a un diagramme
$$
\begin{array}{ccc}
E\times _{A'}A&\rightarrow & A\\
\downarrow && \downarrow \\
E'\times _{B'}B&\rightarrow & B,
\end{array}
$$
o\`u la fl\`eche verticale de gauche est une cofibration $\Gg$-triviale.
Comme par hypoth\`ese la fl\`eche de droite est une $\Gg$-fibration,
on obtient un rel\`evement $E'\times _{B'}B\rightarrow A$.
Maintenant
$$
(E'\times _{B'}B)\cup ^{E\times _{A'}A}E\rightarrow E'
$$
est une \'equivalence pour la topologie grossi\`ere.
On la factorise en
$$
(E'\times _{B'}B)\cup ^{E\times _{A'}A}E\rightarrow E''\rightarrow E'
$$
o\`u la premi\`ere fl\`eche est une cofibration triviale pour la topologie
grossi\`ere et la deuxi\`eme est une fibration triviale pour la
topologie grossi\`ere. On a un morphisme provenant de notre rel\`evement
$$
(E'\times _{B'}B) \cup ^{E\times _{A'}A}E \rightarrow A'.
$$
Comme $A'\rightarrow B'$ est une fibration pour la
topologie grossi\`ere, il existe une extension $E''\rightarrow A'$.
Maintenant, comme $E''\rightarrow E'$ est une fibration triviale pour
la topologie grossi\`ere, il existe une section $E'\rightarrow E''$
qui co\"{\i}ncide sur $E$ avec le morphisme donn\'e $E\rightarrow E''$.
On obtient donc par composition, un morphisme $E'\rightarrow A'$ qui fournit le
rel\`evement cherch\'e.

Pour la deuxi\`eme partie du lemme,
on a un morphisme entre $n$-champs de Segal $f: A\rightarrow B$ qui est
fibrant pour la topologie grossi\`ere. On le factorise en
$$
A\stackrel{g}{\rightarrow} A'\stackrel{h}{\rightarrow} B
$$
avec $g$ une $\Gg$-\'equivalence et $h$ une $\Gg$-fibration. D'apr\`es la
premi\`ere partie,
$A'$ est un champ. Par le lemme pr\'ec\'edent, $g$ est une \'equivalence pour la
topologie grossi\`ere. Si $E\rightarrow E'$ est une cofibration $\Gg$-triviale
avec un diagramme
$$
\begin{array}{ccc}
E&\rightarrow & A\\
\downarrow && \downarrow \\
E'&\rightarrow & B
\end{array}
$$
alors il existe un rel\`evement interm\'ediaire $E'\rightarrow A'$. Comme
plus haut, on peut supposer que le morphisme $E'\rightarrow A'$ est une
fibration
pour la topologie grossi\`ere. Dans ce cas $E'\times _{A'} A\rightarrow E'$
est une \'equivalence objet-par-objet (ici on utilise encore le fait  que les
$A(X)$ sont des $n$-cat\'egories de Segal, cf. le d\'ebut de la
d\'emonstration),
donc c'est une cofibration triviale pour la topologie grossi\`ere.
On a l'inclusion
$E\subset E'\times _{A'}A$. On applique alors la propri\'et\'e
de rel\`evement de $f$ (qui est une
fibration pour la topologie grossi`ere) au diagramme
$$
\begin{array}{ccc}
E'\times _{A'}A&\rightarrow & A\\
\downarrow && \downarrow \\
E'&\rightarrow & B,
\end{array}
$$
pour obtenir un rel\`evement $E'\rightarrow A$ qui r\'epond \`a la question.
\eop

\subnumero{Le champ associ\'e}

Si $A$ est un $n$-pr\'echamp de Segal sur $\Xx$ on appellera {\em champ
associ\'e \`a $A$} toute \'e\-qui\-va\-len\-ce faible
$A\rightarrow A'$ o\`u $A'$ est un $n$-champ de Segal. On note qu'un champ
associ\'e existe toujours, il suffit de prendre pour $A'$ un remplacement
fibrant de
$A$.

\begin{corollaire}
\label{unique}
Soit $A$ un $n$-pr\'echamp et soient $A\rightarrow A'$ et $A\rightarrow A''$
deux champs associ\'es \`a $A$. Supposons que $A''$ est fibrant pour la
topologie grossi\`ere. Alors il existe une \'equivalence (pour la topologie
grossi\`ere) $A'\rightarrow A''$ rendant commutatif le triangle \`a homotopie
pr\`es. Si de plus $A\rightarrow A'$ est une cofibration, on peut supposer
le triangle
strictement commutatif.
\end{corollaire}
{\em Preuve:}
Si $A\rightarrow A'$ est une cofibration, alors c'est une cofibration
triviale et
le fait que $A''$ soit fibrant implique l'existence de l'extension
cherch\'ee $A'\rightarrow
A''$.  Si $A\rightarrow A'$ n'est pas une cofibration, on factorise
$$
A\rightarrow B \rightarrow A'
$$
avec \`a gauche une cofibration triviale et \`a droite une fibration
triviale. Il existe donc une section $A'\rightarrow B$ ainsi qu'une extension
du morphisme $A\rightarrow A''$ en $B\rightarrow A''$. Ceci donne un morphisme
$A'\rightarrow A''$ et on v\'erifie (en suivant les d\'efinitions de
\cite{Quillen})
que ce morphisme rend le triangle commutatif \`a homotopie pr\`es.
\eop

{\em Am\'elioration:} \,\, Avec le $\underline{Hom}$ interne des
$n$-pr\'echamps de
Segal fibrants qui sera d\'efini au \S 11 ci-dessous, on pourra formuler un
meilleur \'enonc\'e d'unicit\'e du ``champ associ\'e'' via la propri\'et\'e
universelle suivante. Si $A\rightarrow A'$ est
un champ associ\'e au pr\'echamp $A$ alors pour tout $n$-pr\'echamp de Segal $B$
fibrant pour la topologie grossi\`ere et qui est un champ, le morphisme
induit
$$
\underline{Hom}(A',B)\rightarrow \underline{Hom}(A,B)
$$
est une \'equivalence de $n$-champs de Segal (fibrants).
Voir \S 13.

\subnumero{Produits fibr\'es}

Soient $A,B,C$ des $n$-pr\'ecats de Segal. Si $B\rightarrow A$
et $C\rightarrow A$ sont des morphismes, on dispose d'un {\em produit fibr\'e
homotopique} $B\times ^h_AC$, essentiellement bien d\'efini: on peut le
d\'efinir comme $B'\times _{A'}C'$ o\`u
$$
\begin{array}{ccccc}
B&\rightarrow & A & \leftarrow & C \\
\downarrow && \downarrow && \downarrow \\
B'&\rightarrow & A' & \leftarrow & C'
\end{array}
$$
est un diagramme dans lequel $A',B'$ et $C'$ sont des $n$-cat\'egories de
Segal avec
$B'\rightarrow A'$ fibrant, et  o\`u les fl\`eches verticales sont des
\'equivalences
faibles. Il est montr\'e dans \cite{limits} Lemma 6.1.2, que ce produit est une
limite dans la $n+1$-cat\'egorie de Segal $nSeCAT$ (\cite{limits} traite le cas
des $n$-cat\'egories mais le cas des $n$-cat\'egories de Segal est identique).

Le remplacement de $A,B,C$ par
$A',B',C'$ peut \'etre choisi fonctoriellement. D'o\`u, si $B\rightarrow
A\leftarrow C$ est un diagramme de $n$-pr\'echamps de Segal sur une cat\'egorie
$\Xx$, on peut consid\'erer leur {\em produit fibr\'e homotopique
objet-par-objet}
$B\times ^h_AC$. C'est un $n$-pr\'echamp
essentiellement bien d\'efini sur
$\Xx$ et (avec les notations de structure interne qu'on introduira au \S 11
ci-dessous) c'est la limite dans $nSeCHAMP (\Xx ^{\rm gro})$ du
diagramme
donn\'e.

Le lemme suivant est l'analogue de l'\'enonc\'e de th\'eorie des
faisceaux selon lequel
le noyau d'un morphisme de faisceaux de groupes (ou le noyau d'une double
fl\`eche s'il
s'agit de faisceaux d'ensembles) est encore un faisceau.

\begin{lemme}
\label{fiprodchamp}
Soit $\Xx , \Gg$ un site et $B\rightarrow A\leftarrow C$ un
diagramme de $n$-champs de Segal sur $\Xx$. Alors
le produit fibr\'e homotopique objet-par-objet $B\times ^h_AC$ est un $n$-champ
de Segal (et c'est aussi la limite du diagramme dans $nSeCHAMP(\Xx ^{\Gg})$,
voir la
notation au \S 11).
\end{lemme}
{\em Preuve:}
On peut choisir un remplacement du diagramme
$$
\begin{array}{ccccc}
B&\rightarrow & A & \leftarrow & C \\
\downarrow && \downarrow && \downarrow \\
B'&\rightarrow & A' & \leftarrow & C'
\end{array}
$$
avec $A', B', C'$ $\Gg$-fibrants, et o\`u les fl\`eches verticales sont des
\'equivalences faibles pour $\Gg$, et $B'\rightarrow A'$ est une
$\Gg$-fibration. D'apr\`es le lemme \ref{morphentrechamps} les fl\`eches
verticales sont des \'equivalences faibles  objet-par-objet, et on a
donc un remplacement ad\'equat pour la
construction du produit fibr\'e homotopique
$$
B\times ^h_AC = B'\times _{A'}C'
$$
qui est $\Gg$-fibrant, en particulier c'est un $n$-champ de Segal. Toute
autre construction de $B\times ^h_AC$ donnera un $n$-pr\'echamp
\'equivalent objet-par-objet au pr\'ec\'edent; et par d\'efinition, ce sera
toujours un
$n$-champ de Segal.
\eop

Le corollaire suivant r\'epond \`a la question originellement pos\'e par C. Rezk
dans \cite{Rezk} (o\`u il donne la d\'emonstration pour les
$0$-pr\'echamps): le passage au ``champ associ\'e'' est compatible
avec les produits fibr\'es homotopiques.

\begin{corollaire}
\label{rezkQ}
Soit $B\rightarrow A \leftarrow C$ un diagramme de $n$-pr\'echamps de Segal.
On suppose que
les valeurs $A(X)$, $B(X)$ et $C(X)$ sont des $n$-cat\'egories de Segal.
Soit $B'\rightarrow A'\leftarrow C'$ le diagramme obtenu en passant
aux $\Gg$-champs associ\'es. Alors
$$
B\times ^h_AC\rightarrow B'\times ^h_{A'}C'
$$
est un $\Gg$-champ associ\'e; ici les produits fibr\'es homotopiques des deux
cot\'es sont pris objet-par-objet au-dessus de $\Xx$.
\end{corollaire}
{\em Preuve:}
La fl\`eche en question est une $\Gg$-\'equivalence faible d'apr\`es le
corollaire
\ref{cordereedy}; et d'apr\`es le lemme pr\'ec\'edent,
le but est un $\Gg$-champ. Il
s'ensuit par d\'efinition que la fl\`eche est un ``champ associ\'ee''.
\eop

\numero{Crit\`eres pour qu'un pr\'echamp soit un champ}
\label{caracterisationpage}

On commence par introduire plusieurs notions
\'equivalentes {\em d'effectivit\'e des donn\'ees de
descente}.
Par la suite, nous utiliserons seulement la version (ii)
ci-dessous; les autres
sont l\`a pour \'eclairer la d\'efinition---le
lecteur peut donc sauter la d\'emonstration.

L'int\'er\^et des conditions (vii) et (viii) ci-dessous est leur caract\`ere
\'el\'ementaire: en effet, elles ne font r\'ef\'erence
ni \`a la structure de mod\`eles \ref{cmf}, ni
\`a l'op\'eration $SeCat$. Cependant, sans la structure de cmf ces conditions
seraient quasiment inutilisables. On les a incluses pour
faciliter \'eventuellement une exposition sans
d\'emonstrations de la notion de champ.

\begin{lemmedefinition}
\label{lemmedef}
Soit $A$ un $n$-pr\'echamp de Segal sur un site $\Xx$. Soit $X\in \Xx$ et
$\Bb \subset \Xx /X$ un crible couvrant $X$. On dira que {\em les donn\'ees
de descente pour $A$ sur $\Bb$ sont effectives} si l'une des
conditions \'equivalentes suivantes est satisfaite.
\newline
(i)\,  Tout morphisme de $\ast _{\Bb}$ dans $A|_{\Bb}$ dans la cat\'egorie
homotopique  $Ho(nSePCh(\Bb ^{\rm gro}))$ des $n$-pr\'echamps de Segal par
rapport
aux \'equivalences objet-par-objet, s'\'etend en un morphisme de $\ast$ vers
$A(X)$ dans $Ho(nSePC)$. Autrement dit, l'application naturelle
$$
[\ast , A(X)]\rightarrow [\ast _{\Bb} , A|_{\Bb}]_{\Bb ^{\rm gro}}
$$
est surjective;
\newline
(ii)\, Pour tout remplacement fibrant $A\rightarrow A'$ pour la topologie
grossi\`ere, tout morphisme $\ast _{\Bb}\rightarrow A'$ s'\'etend en un
morphisme $\ast _{\Xx /X}\rightarrow A'$;
\newline
(ii')\, Il existe un remplacement fibrant $A\rightarrow A'$ pour la topologie
grossi\`ere, tel que tout morphisme $\ast _{\Bb}\rightarrow A'$ s'\'etende en un
morphisme $\ast _{\Xx /X}\rightarrow A'$;
\newline
(iii)\, Pour tout remplacement fibrant $A\rightarrow A'$ pour la topologie
grossi\`ere, le morphisme de $n$-cat\'egories de Segal
$$
A'(X)=\Gamma (\Xx /X , A'|_{\Xx /X})\rightarrow \Gamma (\Bb , A'|_{\Bb})
$$
est essentiellement surjectif;
\newline
(iii')\, Il existe un remplacement fibrant $A\rightarrow A'$ pour la topologie
grossi\`ere, tel que le morphisme de $n$-cat\'egories de Segal
$$
A'(X)=\Gamma (\Xx /X , A'|_{\Xx /X})\rightarrow \Gamma (\Bb , A'|_{\Bb})
$$
soit essentiellement surjectif;
\newline
(iv)\, Le morphisme
$$
A(X)\rightarrow \lim _{\leftarrow , \Bb ^o} A|_{\Bb }
$$
est essentiellement surjectif;
\newline
(v)\, Pour tout remplacement fibrant $A\rightarrow A'$ pour la structure
HBKQ (\ref{cmfHirschho})
i.e. une \'equivalence objet-par-objet avec chaque $A'(X)$
fibrant, et pour tout morphisme
$$
\eta : D_{\Bb}\rightarrow A'
$$
il existe une homotopie $h: D_{\Bb}\times \overline{I}\rightarrow A'$
avec $h(0)=\eta$ et telle que $h(1)$ provienne d'un morphisme $\ast\rightarrow
A(X)$. Ici $D_{\Bb}\rightarrow \ast _{\Bb}$ est le remplacement HBKQ-cofibrant
construit \`a la fin de \S 5.
\newline
(v')\, Il existe un remplacement fibrant $A\rightarrow A'$ pour la
structure HBKQ
(\ref{cmfHirschho}) tel que pour tout morphisme
$$
\eta : D_{\Bb}\rightarrow A'
$$
il existe une homotopie $h: D_{\Bb}\times \overline{I}\rightarrow A'$
avec $h(0)=\eta$ et telle que $h(1)$ provienne d'un morphisme $\ast\rightarrow
A(X)$.

Si on suppose, de plus, que $A(Y)$ est une $n$-cat\'egorie de Segal pour tout
$Y\in \Xx$, alors ces conditions sont \'equivalentes aux conditions suivantes:
\newline
(vi)\, Pour tout diagramme de pr\'efaisceaux de $n$-cat\'egories de Segal
$$
\ast _{\Bb } \leftarrow C \rightarrow A
$$
o\`u le premier morphisme est une \'equivalence objet-par-objet, il existe $\ast
\rightarrow A(X)$ tel que le compos\'e
$$
C\rightarrow \ast _{\Bb} \rightarrow \ast _{\Xx /X} \rightarrow A
$$
soit homotope (au sens de Quillen) au morphisme $C\rightarrow A$ de d\'epart;
\newline
(vii) \, Tout diagramme de pr\'efaisceaux de $n$-cat\'egories de Segal
$$
\ast _{\Bb } \leftarrow C \rightarrow A
$$
o\`u le premier morphisme est une \'equivalence objet-par-objet se
compl\`ete en un
diagramme
$$
\begin{array}{ccccc}
\ast _{\Bb }& \leftarrow &C &\rightarrow &A\\
\downarrow &&\downarrow && \parallel \\
\ast _{\Xx /X }& \leftarrow &C' &\rightarrow &A
\end{array}
$$
o\`u le morphisme $C'\rightarrow \ast _{\Xx /X}$ est une \'equivalence
objet-par-objet; et
\newline
(viii)\, Soit $nSeCat^{\Bb^o}$ la cat\'egorie des pr\'efaisceaux de
$n$-cat\'egories de Segal sur $\Bb$; et soit $Ho( nSeCat^{\Bb^o})$
sa localis\'ee de Gabriel-Zisman, obtenue en inversant les \'equivalences
objet-par-objet. Alors tout morphisme
$\ast _{\Bb}\rightarrow A$ dans $Ho( nSeCat^{\Bb^o})$
provient d'un morphisme $\ast \rightarrow A(X)$.
\end{lemmedefinition}
{\em Preuve:}
\newline
(i)$\Leftrightarrow$(ii)($\Leftrightarrow$(ii'))\,\,
D'apr\`es Quillen \cite{Quillen} tout
\'el\'ement de
$[\ast _{\Bb}, A|_{\Bb}]_{\Bb ^{\rm gro}}$ provient d'un morphisme
$f:\ast _{\Bb}\rightarrow A'|_{\Bb }$ (remarquer que $A'|_{\Bb }$ est
fibrant pour
la topologie grossi\`ere voir \S 4, et que $\ast_{\Bb}$ est automatiquement
cofibrant). Si $f$ s'\'etend \`a $\ast _{\Xx /X}$ alors la classe de $f$
provient \'evidemment de $[\ast , A(X)]$. Dans l'autre sens, si $f$ correspond
\`a un \'el\'ement de $[\ast _{\Bb}, A|_{\Bb}]_{\Bb ^{\rm gro}}$ provenant de
$[\ast , A(X)]$ alors il existe un morphisme $g:\ast _{\Xx /X}\rightarrow
A'$ tel
que $g|_{\Bb}$ soit homotope (\`a gauche) \`a $f$ au sens de Quillen
\cite{Quillen}. Comme $\ast _{\Bb} \times \overline{I}$ est un
objet-cylindre pour
$\ast_{\Bb}$, il existe un morphisme $$
h: \ast _{\Bb} \times \overline{I}\rightarrow A'
$$
avec $h(0)= f$ et $h(1)= g|_{\Bb}$. Maintenant on a un morphisme
$$
(\ast _{\Bb} \times \overline{I})\cup^{\ast _{\Bb} \times \{ 1\}}
(\ast _{\Xx /X} \times \{ 1\}) \rightarrow A'
$$
qui se prolonge le long de la cofibration triviale  (pour la topologie
grossi\`ere)
$$
(\ast _{\Bb} \times \overline{I})\cup^{\ast _{\Bb} \times \{ 1\}}
(\ast _{\Xx /X} \times \{ 1\}) \hookrightarrow
\ast _{\Xx /X} \times \overline{I},
$$
en un morphisme dont la restriction \`a $\ast _{\Xx /X} \times \{ 0\}$ fournit
l'\'extension de $f$ cherch\'ee.
\newline
(ii)$\Leftrightarrow$(iii)($\Leftrightarrow$(iii')) \, est imm\'ediat au vu de
la d\'efinition de $\Gamma$.
\newline
(iii)$\Leftrightarrow$(iv) \, d\'ecoule de la proposition \ref{calclim}
ci-dessous. \newline
(i)$\Leftrightarrow$(v)($\Leftrightarrow$(v')) On raisonne comme pour
(i)$\Leftrightarrow$(ii)  mais avec la structure HBKQ; cf. \S 5. Notons que
la notion de classe d'homotopie d'applications est ind\'ependante du choix de
la structure de cmf et ne d\'epend que des \'equivalences faibles; en
particulier dans (i) il s'agit des m\^emes ensembles dans la structure de
\ref{cmf} ou dans la structure HBKQ \ref{cmfHirschho}.

On suppose maintenant que $A$ est un pr\'efaisceau de $n$-cat\'egories de Segal.
On a l'\'equivalence $Ho(nSeCat^{\Bb})\cong Ho(nSePCh(\Bb ^{\rm gro})_f)$ car
il y a une \'equivalence naturelle ``remplacement fibrant''
de $nSeCat^{\Bb}$ vers $nSePCh(\Bb ^{\rm gro})_f$. D'autre part on a
l'\'equivance $Ho(nSePCh(\Bb ^{\rm gro})_f) \cong
Ho(nSePCh(\Bb ^{\rm gro}))$
(\cite{Quillen}). Ceci donne (i)$\Leftrightarrow$(viii).
Pour les conditions (vi) et (vii), on choisit un remplacement fibrant
$A\rightarrow A'$ pour la topologie grossi\`ere. Tout morphisme
$\ast _{\Bb}\rightarrow A'$ admet donc une factorisation
$$
\ast _{\Bb}\rightarrow E \rightarrow A'
$$
o\`u la premi\`ere fl\`eche est une cofibration triviale pour la topologie
grossi\`ere et la deuxi\`eme une fibration pour la topologie
grossi\`ere. Il existe une (unique) r\'etraction $E\rightarrow \ast _{\Bb}$
(on peut dire que $\ast _{\Bb}$ est fibrant pour la topologie grossi\`ere mais
en fait il suffit de consid\'erer les composantes connexes de $C$). Posons
$$
C:= E\times _{A'} A.
$$
Le fait que les $A(Y)$ soient des $n$-cat\'egories de Segal permet
d'appliquer \ref{proper} pour obtenir que $C\rightarrow E$ est une \'equivalence
objet-par-objet. Notons d'apr\`es \ref{stuff} que les $C(Y)$ sont des
$n$-cat\'egories de Segal. Ils'ensuit que tout \'el\'ement de $[\ast _{\Bb},
A|_{\Bb }]_{\Bb ^{\rm gro}}$ provient d'un diagramme (de pr\'efaisceaux de
$n$-cat\'egories de Segal)
$$
\ast _{\Bb}\leftarrow C
\rightarrow A
$$
o\`u le premier morphisme est une \'equivalence objet-par-objet.
On laisse au lecteur de compl\'eter la d\'emonstration que (vi) et (vii) sont
\'equivalentes aux autres conditions, en exprimant de la m\^eme fa\c{c}on la
condition qu'un tel diagramme donne lieu \`a un \'el\'ement de
$[\ast _{\Bb},
A|_{\Bb }]_{\Bb ^{\rm gro}}$ qui provient de $[\ast , A(X)]$.
\eop

\subnumero{Crit\`ere principal}

On a la caract\'erisation suivante des champs, qui est analogue \`a la
d\'efinition de
$1$-champ dans Giraud \cite{GiraudThese} et Laumon-Moret-Bailly
\cite{LaumonMB}.
Les conditions (a) et (b) sont les analogues des deux conditions dans la
d\'efinition de faisceau. Notons aussi que la condition (a) nous am\`ene
r\'ecursivement vers la condition (b).

\begin{proposition}
\label{critere}
Soit $n\geq 1$.
Soit $A$ un $n$-pr\'echamp de Segal sur $\Xx$ dont les valeurs sont
des $n$-cat\'egories de Segal.
Alors $A$ est un $n$-champ de Segal si et seulement si:
\newline
(a)\,\, pour tout $X\in \Xx$ et tout $x,y\in A_0(X)$ le $(n-1)$-pr\'echamp
de Segal $A_{1/}(x,y)$ sur $\Xx /X$ est un $(n-1)$-champ de Segal; et
\newline
(b) \,\, pour tout $X\in \Xx$ et tout crible $\Bb \subset \Xx /X$ dans un
ensemble de cribles qui engendre la topologie, les donn\'ees de descente
pour $A$
par rapport \`a $\Bb$ sont effectives, voir la d\'efinition
\ref{lemmedef}
ci-dessus.
\end{proposition}
{\em Preuve:}
Les conditions (a) et (b) sont \'evidemment n\'ecessaires (pour (b) on
utilise le
fait que $A'$ est $\Gg$-fibrant en vertu du  Lemme \ref{grofibchampimplfib}).
Supposons que (a) et (b) sont satisfaites et soit $A'\rightarrow A''$ une
cofibration $\Gg$-triviale vers un $n$-pr\'echamp de Segal $\Gg$-fibrant. Pour
tout $X\in \Xx$ et $x,y\in A_0(X)$, le morphisme
$$
A_{1/} (x,y)\rightarrow         A''_{1/}(x,y)
$$
est une cofibration $\Gg_X$-triviale de $(n-1)$-pr\'echamps de Segal sur
$\Xx /X$ dont le but est
$\Gg _X$-fibrant (Corollaire \ref{Gfibrestrict}).  La condition (a) 
implique que
ce morphisme est  une \'equivalence pour la
topologie grossi\`ere, en particulier
$$
A(X)_{1/} (x,y)\rightarrow      A''(X)_{1/}(x,y)
$$
est une \'equivalence. Il en est de m\^eme pour le morphisme
$$
A'(X)_{1/} (x',y')\rightarrow  A''(X)_{1/}(x',y')
$$
pour
tous $x',y'\in A'_0(X)$ car $x'$ et $y'$ sont \'equivalents \`a des objets
provenant
de $A_0(X)$. Ceci montre que $A'\rightarrow A''$ est pleinement fid\`ele
objet-par-objet. Il reste \`a voir que $A'\rightarrow A''$ est essentiellement
surjectif objet-par-objet. Pour cela, on peut choisir une factorisation
$$
A'\rightarrow B \rightarrow A''
$$
o\`u le premier morphisme est une cofibration triviale pour la topologie
grossi\`ere et le deuxi\`eme une fibration pour la topologie
grossi\`ere. Quitte \`a remplacer $A'$ par $B$, on peut dor\'enavant
supposer que
$A'\rightarrow A''$ est fibrant pour la topologie grossi\`ere.

Soit $X\in \Xx$ et
$u\in A''_0(X)$, qu'on consid\`ere comme morphisme
$$
u: \ast _X \rightarrow A'' |_{\Xx /X}.
$$
On pose
$$
C:= (A'|_{\Xx /X})\times _{A'' |_{\Xx /X}}\ast _X .
$$
Le morphisme
$C\rightarrow \ast _X$ est une fibration pour la topologie grossi\`ere
et aussi une \'e\-qui\-va\-len\-ce pour la topologie $\Gg$. De plus, comme
$A'\rightarrow A''$ est pleinement fid\`ele, le morphisme
$C\rightarrow \ast _X$ est pleinement fid\`ele. En particulier, $C$ est
$0$-tronqu\'e, et \'equivalent objet-par-objet \`a un pr\'efaisceau de
sous-ensembles
de $\ast _X$. Comme notre morphisme est une \'equivalence pour la topologie
$\Gg$,
il existe un crible $\Bb \subset \Xx /X$ de la topologie (et m\^eme de
l'ensemble g\'en\'erateur de cribles mentionn\'e
dans (b)) tel que $\ast _{\Bb}$ soit
contenu dans ce sous-pr\'efaisceau de $\ast _X$. En particulier, le morphisme
$$
C\times _{\ast _X}\ast _{\Bb} \rightarrow \ast _{\Bb}
$$
est une \'equivalence objet-par-objet et, comme c'est aussi une fibration pour
la topologie grossi\`ere, il existe une section $\ast _{\Bb }\rightarrow C$.
On obtient ainsi un morphisme
$$
\ast _{\Bb} \rightarrow A'|_{\Xx /X}
$$
qui, gr\^ace \`a (b), s'\'etend en un morphisme
$$
u':\ast _{X} \rightarrow A'|_{\Xx /X}.
$$
On obtient ainsi un \'el\'ement $u'\in A'_0(X)$. Avec le m\^eme type
d'argument on voit que l'image de $u'$ dans $A''_0(X)$ est \'equivalente
\`a $u$, ce qui prouve l'essentielle surjectivit\'e.
\eop

\subnumero{Crit\`ere pour les $0$-champs de Segal}

On a un r\'esultat analogue \`a \ref{critere} pour les $0$-pr\'echamps de
Segal i.e.
les pr\'efaisceaux simpliciaux. Pour l'obtenir, on utilise le lemme suivant.

\begin{lemme}
\label{pimseg}
Soit $m,n\geq 0$ et soit $A$ un pr\'efaisceau de $n$-cat\'egories de Segal.
Alors $A$ est un $n$-champ de Segal si et seulement si $\Pi _{m,Se}\circ A$ est
un $n+m$-champ de Segal.
\end{lemme}
{\em Preuve:}
D'abord on montre que si $A$ est un $n$-champ de Segal, alors
$\Pi _{m,Se}\circ A$ est un $n+m$-champ de Segal. En raisonnant par r\'ecurrence
sur $m$ et en appliquant le crit\`ere \ref{critere}, on se ram\`ene
au cas $n=0$ et $m=1$. Soit donc $A$ un $0$-champ de Segal. On peut
supposer que $A$ est $\Gg$-fibrant. On va prouver les conditions
(a) et (b) de \ref{critere}
pour $B=\Pi _{1,Se}\circ A$. Notons que $A$ et $\Re _{\geq 0}B$ sont
\'equivalents
(objet-par-objet). Si $X$ est un objet de
$\Xx$ et $\Bb$ un crible
couvrant $X$, tout morphisme $\ast _{\Bb}\rightarrow B$ induit
un morphisme $\ast
_{\Bb}\rightarrow \Re _{\geq 0}B$ et donc une donn\'ee de descente pour $A$
(i.e. un \'el\'ement de $[\ast _{\Bb}, A|_{\Bb }]_{\Bb ^{\rm gro}}$).
Maintenant,
l'inclusion $\ast _{\Bb}\rightarrow \ast _{\Xx /X}$ est une $\Gg$-\'equivalence
et le fait que $A$ soit $\Gg$-fibrant implique donc l'existence d'une extension
de notre donn\'ee de descente en un \'el\'ement de $[\ast , A(X)]= [\ast ,
B(X)]$. Ceci prouve
la condition (b) de \ref{critere} pour $B$. D'autre part, soient
$a$ et $b$ dans $B_0(X)$. On a l'\'equivalence
$$
B_{1/}(a,b) \cong Path ^{a,b}(A):= \ast \times _{A|_{\Xx /X}}\ast ,
$$
et ce
dernier pr\'echamp
est un $0$-champ de Segal d'apr\`es \ref{fiprodchamp}. Ceci donne
la condition (a) de
\ref{critere}, et donc, par \ref{critere}, $B$ est un $1$-champ de Segal.

Pour l'autre direction, supposons que
$\Pi _{m,Se}\circ A$ est un $n+m$-champ de Segal.
Choisissons un remplacement $\Gg$-fibrant
$A\rightarrow A'$. On obtient un morphisme
$$
\Pi _{m,Se}\circ A\rightarrow \Pi _{m,Se}\circ A',
$$
dont le but est un $n+m$-champ de Segal d'apr\`es le premier paragraphe.
Par hypoth\`ese la source est un $n+m$-champ de Segal. Comme ce morphisme est
une $\Gg$-\'equivalence faible, \ref{morphentrechamps} implique que c'est
une \'equivalence objet-par-objet. Ceci implique \`a son tour que $A\rightarrow
A'$ est une \'equivalene objet-par-objet, et donc $A$ est un $n$-champ de Segal.
\eop

On donne maintenant une version du crit\`ere adapt\'ee aux
$0$-pr\'echamps de Segal
(ou pr\'efaisceaux simpliciaux). Rappelons ici qu'un pr\'efaisceau
simpliciale $U$
est un $0$-champ de Segal, si et seulement si $U$ est {\em flasque par rapport
aux objets du site} au sens de Jardine \cite{Jardine}, ou {\em de
descente cohomologique} dans la terminologie de Thomason
\cite{Thomason} \cite{TravauxThomason}.

\begin{corollaire}
\label{0critere}
Si $U$ est un pr\'efaisceau simplicial au-dessus de $\Xx$, alors $U$ est un
$0$-champ de Segal si et seulement si
\newline
(a)\,\,  pour tout $X\in \Xx$ et tous $x,y\in U_0(X)$ le pr\'efaisceau
simplicial \newline
$Y\mapsto Path ^{x|_Y,y|_Y}(U(Y))$ sur $\Xx /X$ est un
$0$-champ de Segal; et
\newline
(b)\,\, pour tout $X\in \Xx$ et tout crible $\Bb \subset \Xx /X$ (dans une
famille de cribles qui engendre la topologie), les donn\'ees de descente
pour $U$
sur $\Bb$ sont effectives.
\end{corollaire}
{\em Preuve:}
On observe, comme dans la d\'emonstration pr\'ec\'edente, que si $A:= \Pi
_{1,Se}\circ U$ alors on a
$$
\left( Y\mapsto Path ^{x|_Y,y|_Y}(U(Y))\right) \cong A_{1/}(x,y)
$$
et d'autre part, l'ensemble des donn\'ees de descente (resp. effectives)
pour $U$ est isomorphe \`a l'ensemble des donn\'ees de descente (resp.
effectives) pour $A$. On conclut \`a l'aide de \ref{critere} et \ref{pimseg}.
\eop

Dans la condition (a) de ce corollaire, $Path ^{x|_Y,y|_Y}(U(Y))$ d\'esigne
l'ensemble simplicial des chemins entre $x|_Y$ et
$y|_Y$ dans un remplacement de Kan de $U(Y)$.

Si $U$ est un $0$-champ de Segal $k$-tronqu\'e (i.e. les $U(X)$ ont
leurs groupes d'homotopie nuls en degr\'e $>k$) alors le crit\`ere de
\ref{0critere} devient r\'ecursif puisque dans la condition $(a)'$, $Path
^{x|_Y,y|_Y}(U(Y))$ est $k-1$-tronqu\'e. On peut convertir cette
remarque en
une version de la Proposition \ref{critere} adapt\'ee aux $k$-champs (tout
court et
non de
Segal).

\medskip

\subnumero{Compatibilit\'e avec la restriction aux $\Xx /X$}

On revient aux $n$-champs de Segal; on peut enfin compl\'eter
les r\'esultats du \S
4 concernant la topologie $\Gg$.

\begin{lemme}
\label{grest}
Soit $\Xx$ un site dont on note $\Gg$ la topologie de Grothendieck. Un
$n$-pr\'echamp de
Segal $A$ sur $\Xx$ est un champ si et seulement si $A|_{\Xx /X}$ est un champ
pour tout $X\in \Xx$. La restriction des $n$-pr\'echamps de Segal
$$
p^{\ast}: A\mapsto A|_{\Xx /X}
$$
pr\'eserve les $\Gg$-fibrations et les $\Gg$-fibrations triviales; pour la
topologie
$\Gg$, $p^{\ast}$ est un foncteur de Quillen \`a la fois \`a gauche et \`a
droite.
\end{lemme}
{\em Preuve:}
Il est \'evident d'apr\`es nos crit\`eres (\ref{critere} et \ref{0critere}) que
$A$ est un champ si et seulement si $A|_{\Xx /X}$ est un champ pour tout $X\in
\Xx$. En particulier la restriction $p^{\ast}$ envoie les objets $\Gg$-fibrants
sur des objets $\Gg$-fibrants (cf le lemme \ref{grofibchampimplfib}). Soit $f:
A\rightarrow B$ une $\Gg$-fibration entre objets $\Gg$-fibrants. Alors $f|_{\Xx
/X}$ est une fibration pour la topologie grossi\`ere (d'apr\`es \ref{shriek})
entre objets $\Gg$-fibrants. D'apr\`es le lemme \ref{stuff}, on a que $f|_{\Xx
/X}$ est $\Gg$-fibrant.

Maintenant on sait que $p^{\ast}$ transforme les fibrations entre objets
fibrants
en fibrations. Or, une cofibration dans une cmf est triviale si et seulement si
elle v\'erifie la propri\'et\'e de rel\`evement pour toutes les fibrations entre
objets fibrants. Ceci implique que l'adjoint $p_!$ de
$p^{\ast}$ transforme les cofibrations triviales en cofibrations triviales (et
on savait d\'ej\`a que $p_!$ pr\'eserve les cofibrations). Donc $p_!$ est un
foncteur de Quillen \`a gauche, et par suite $p^{\ast}$ est un foncteur de
Quillen \`a droite.
\eop

\subnumero{A propos de $\emptyset$}
Habituellement dans un site $\Xx$ on a un objet initial $\iota$
(g\'en\'eralement l'objet ``vide'' $\iota = \emptyset$), et
les champs sur $\Xx$ doivent satisfaire la condition $A(\iota )\cong \ast$,
du moins lorsque la topologie $\Gg$ admet
le crible vide $\Bb
=\emptyset \subset \Xx /\iota$ comme crible couvrant $\iota$ (ce qui sera le
cas dans les exemples g\'eom\'etriques).  En effet,
sous cette hypoth\`ese, si $A$ est $\Gg$-fibrant,
alors comme  $E'\times \ast _{\iota} \rightarrow E\times \ast _{\iota}$ est une
cofibration triviale, pour toute cofibration de $n$-pr\'ecats de Segal
$E'\subset
E$, tout morphisme $E'\rightarrow A(\iota )$ s'\'etend en
$E\rightarrow A(\iota )$ ce qui entraine que $A(\iota )\rightarrow \ast$ est
une fibration triviale.

Voyons ce que \c{c}a donne  avec nos crit\`eres: puisque  $\Bb
=\emptyset \subset \Xx /\iota$ est un crible couvrant $\iota$  ($\Bb
$ ne contient m\^eme pas $\iota$), le pr\'echamp $\ast _{\Bb}$ est
\'egal au vide
$\emptyset$ au-dessus de $\Xx /\iota$ donc il existe un unique morphisme $\ast
_{\Bb}\rightarrow A'|_{\Xx /\iota})$. La condition (b) de nos crit\`eres
implique
donc que $A'(\iota )$ est non-vide. D'autre part, la condition (a) de nos
crit\`eres
implique (par r\'ecurrence sur $n$) que pour chaque paire d'objets de $A(\iota
)$ la
$n-1$-cat\'egorie de morphismes correspondante est $\ast$. D'o\`u $A(\iota
)=\ast$.
Pour $n=0$ le m\^eme argument (en rempla\c{c}ant \`a chaque fois $A$ par son
$\Pi _{1,Se}\circ A$) montre par r\'ecurrence que tous les groupes
d'homotopie de
l'ensemble simplicial $A(\iota )$ s'annulent;
c'est le cas initial de la r\'ecurrence
pour $n>0$.
En revanche, pour ce qu'on appelle la ``topologie grossi\`ere'', on n'admet pas
le crible vide
$\emptyset \rightarrow \Xx /\iota$ comme crible couvrant. Dans ce cas on n'a
pas la condition $A(\iota )\cong \ast$.

\subnumero{D\'evissage}

Pour $n>0$ on peut d\'evisser la condition pour \^etre un champ en termes
des $k<n$,
gr\^ace au corollaire suivant.  Rappelons encore qu'un pr\'efaisceau simplicial
$U$ est un $0$-champ de Segal, si et seulement si $U$ est flasque par
rapport aux objets du site (Jardine \cite{Jardine}), ou
de descente cohomologique (Thomason \cite{Thomason}
\cite{TravauxThomason}). Cette condition est donc bien connue.

On rappelle (cf. \S 2) que si $A$ est une $n$-cat\'egorie de Segal, pour
$k\leq n$, on
dispose de l' {\em int\'erieur $k$-groupique} $A^{int, k}$ qui est une
$k$-cat\'egorie de Segal. Essentiellement cette construction consiste \`a
ne garder que les $i$-fl\`eches qui sont inversibles \`a \'equivalence
pr\`es ($i>k$) et \`a consid\'erer le r\'esultat comme une
$k$-cat\'egorie de Segal. Pour $k=0$, $A^{int, 0}$ est donc un ensemble
simplicial avec $\Pi _{n,Se}\circ A^{int, 0}$ \'equivalent au plus
grand $n$-groupo\"{\i}de de Segal contenu dans $A$. Ces constructions
s'\'etendent de mani\`ere \'evidente \`a des pr\'efaisceaux de $n$-cat\'egories
de Segal.

Le corollaire suivant permet de d\'evisser la condition pour qu'un
$n$-pr\'echamp de
Segal $A$ soit un $n$-champ de Segal, en termes de la condition  pour
qu'un pr\'efaisceau simplicial soit un $0$-champ de Segal.

\begin{corollaire}
\label{devissage}
Soit $A$ un $n$-pr\'echamp de Segal ($n\geq 1$) tel que $A(X)$ soit une
$n$-cat\'egorie
de Segal pour tout $X\in \Xx$. Alors $A$ est un $n$-champ de Segal
si et seulement si:
\newline
(a) \,  pour tout $X\in \Xx$ et toute paire d'objets $x,y\in A_0(X)$, le
$n-1$-pr\'echamp de Segal $A_{1/}(x,y)$ est un $n-1$-champ de Segal au-dessus de
$\Xx /X$; et
\newline
(b) \, $A^{int, 0}$ est un $0$-champ de Segal.
\end{corollaire}
{\em Preuve:}
On a
$$
[\ast _{\Bb}, A|_{\Bb}]_{\Bb ^{\rm gro}} =
[\ast _{\Bb}, A^{int, 0}|_{\Bb}]_{\Bb ^{\rm gro}} ,
$$
et on applique \ref{critere} et \ref{0critere}.
\eop

\begin{corollaire}
Si $A$ est un $n$-champ de Segal alors pour $0\leq k < n$, $A^{int,k}$ est un
$k$-champ de Segal.
\end{corollaire}
\eop

En fait, le lecteur pourra se convaincre que le pr\'esent d\'evissage nous
ram\`ene
\`a un nombre infini de conditions d'effectivit\'e des donn\'ees de
descente, une pour les $i$-fl\`eches pour chaque $i\geq 0$. Pour \'enoncer ceci
plus pr\'ecisement, introduisons des notations. Rappelons d'abord une notation
de Tamsamani \cite{Tamsamani}: l'ensemble de $i$-fl\`eches d'une
$n$-cat\'egorie $A$, not\'e $Fl^i(A)$, est l'ensemble $A_{1^i}$
o\`u $1^i = (1,\ldots , 1)\in \Theta ^n$ est la classe de $(1,\ldots ,
1,0,\ldots , 0)\in \Delta
^n$. On a les applications $s$ (``source'') et $b$ (``but'') de
$Fl^i(A)$ vers $Fl^{i-1}(A)$. Dans notre situation, si $A$ est une
$n$-cat\'egorie de Segal, on conviendra de dire que les
$i$-fl\`eches pour $i>n$ sont les $i$-fl\`eches de $\Pi _{m,Se}\circ A$
pour $m> i-n$.  De fa\c{c}on plus syst\'ematique on pourrait d\'efinir
(au vu de notre d\'efinition de $\tau _{\leq m}$ voir \S 2)
$$
Fl^i(A):= Fl^i(\tau _{\leq i+1}A).
$$
Soit $i\geq 1$ et soient $x,y\in Fl^{i-1}(A)$ avec $s(x)=s(x)$ et $b(x)=b(y)$.
On d\'efinit alors la $(n-i)$-cat\'egorie de Segal
$\underline{Fl}^i(A; x,y)$ comme
la sous-$(n-i)$-cat\'egorie de Segal de $A_{1^i/}$ form\'ee des objets
dont la source est $x$ et le but est $y$ (ici, si $n-i$ est n\'egatif, on le
remplace par $0$).
\footnote{
Cette inf\'elicit\'e devrait dispara\^{\i}tre \`a terme quand on pourra
parler des
{\em $\infty$-cat\'egories}.}

Si $A$ est un $n$-pr\'echamp de Segal dont les valeurs $A(X)$ sont des
$n$-cat\'egories de Segal, alors pour tout $X\in \Xx$ et $x,y\in
Fl^{i-1}(A(X))$ avec $s(x)=s(y)$ et $b(x)=b(y)$, on dispose du $n-i$-pr\'echamp
de Segal $\underline{Fl}^i(A|_{\Xx /X}; x,y)$ qui \`a
$Y\rightarrow X$ associe la $n-i$-cat\'egorie
$\underline{Fl}^i(A(Y); x|_Y,y|_Y)$.

\begin{proposition}
\label{conditiond}
Soit $A$ un $n$-pr\'echamp de Segal sur $\Xx$ tel que les $A(X)$ soient des
$n$-cat\'egories de Segal. Alors $A$ est un $n$-champ de Segal si et seulement
si pour tout $X\in \Xx$ et tout crible $\Bb \subset \Xx /X$ (dans une famille
de cribles qui engendre la topologie), on a:
\newline
---les donn\'ees de descente pour $A$ par rapport \`a $\Bb$ sont effectives;
et
\newline
---pour tout $1\leq i < \infty$, pour tout $x,y\in Fl^{i-1}(A(X))$ avec
$s(x)=s(y)$ et $b(x)=b(y)$, les donn\'ees de descente pour
$\underline{Fl}^i(A|_{\Xx /X}, x,y)$ sur $\Bb$ sont effectives.
\end{proposition}
La d\'emonstration est laiss\'ee en exercice au lecteur. {\em
Indication:} combiner le lemme \ref{pourloc}, avec la remarque qu'il suffit,
pour obtenir toutes les cofibrations triviales, de consid\'erer des suites
ad\'equates compos\'ees de cofibrations triviales pour la topologie grossi\`ere,
ainsi que de cofibrations de la forme $cof(U\hookrightarrow U'; \Bb \subset \Xx
/X)$ o\`u $U'= h(1^i)$ repr\'esente les $i$-fl\`eches, et $U$ est son
''bord'' qui
repr\'esente les couples de $i-1$-fl\`eches avec m\^eme source et but.
\eop

\subnumero{Le cas des $1$-pr\'echamps de Segal et pr\'efaisceaux de
cat\'egories simpliciales}

Dans la derni\`ere partie du papier, on s'int\'eressera surtout au cas $n=1$;
pour cette raison, on recopie \ref{devissage} pour les $1$-cat\'egories de
Segal ainsi que pour les pr\'efaisceaux de cat\'egories simpliciales. On
rappelle
encore que les $0$-champs de Segal sont exactement les pr\'efaisceaux
simpliciaux
flasques par rapport aux objets du site \cite{Jardine} ou encore v\'erifiant
la condition de
descente cohomologique de \cite{Thomason} \cite{TravauxThomason}.

\begin{corollaire}
Soit $A$ un $1$-pr\'echamp de Segal tel que $A(X)$ soit une cat\'egorie de
Segal pour tout $X\in \Xx$. Alors $A$ est un $1$-champ de Segal
($\infty$-champ) si et seulement si:
\newline
(a) \, pour tout $X\in \Xx$ et toute paire d'objets $x,y\in A_0(X)$, le
pr\'efaisceau simplicial $A_{1/}(x,y)$ au-dessus de $\Xx /X$ est un $0$-champ
de Segal; et
\newline
(b) \, le pr\'efaisceau simplicial $A^{int, 0}$ est un $0$-champ
de Segal.
\end{corollaire}
\eop

Si $B$ est une cat\'egorie simpliciale, on peut calculer l'ensemble simplicial
$B^{int, 0}$ de la fa\c{c}on suivante: soit $WB$ la sous-cat\'egorie
simpliciale o\`u on ne garde que les fl\`eches qui sont inversibles \`a
homotopie pr\`es; soit $\nu WB$ son ``nerf'', i.e. l'ensemble
bisimplicial associ\'e (c'est aussi la $1$-pr\'ecat de Segal associ\'ee).
Alors, si
$d(\nu WB)$ est l'ensemble simplicial diagonal correspondant,
on a $d(\nu WB)\cong B^{int, 0}$.
Cette discussion s'\'etend de
fa\c{c}on \'evidente aux pr\'efaisceaux de cat\'egories simpliciales (qui
donnent lieu \`a des $1$-pr\'echamps de Segal qu'on note avec la m\^eme
lettre).

\begin{corollaire}
Soit $B$ un pr\'efaisceau de cat\'egories simpliciales. Alors $B$ est un
$1$-champ de Segal si et seulement si:
\newline
(a) \, pour tout $X\in \Xx$ et toute paire d'objets $x,y\in Ob (B)(X)$, le
pr\'efaisceau simplicial $Mor_B(x,y)$ au-dessus de $\Xx /X$ est un $0$-champ de
Segal, et
\newline (b) \, le pr\'efaisceau simplicial
$d(\nu WB)= B^{int,0}$ est un $0$-champ de
Segal.
\end{corollaire}
\eop

\subnumero{Une variante du crit\`ere}


On donne ici un \'enonc\'e tr\`es l\'eg\`erement
diff\'erent de \ref{critere}.
On remarque d'abord que si $A'$ est un $n$-pr\'echamp de Segal fibrant
pour la topologie grossi\`ere, alors $A'|_{\Xx /X}$ est aussi fibrant pour la
topologie grossi\`ere, ainsi que $A'|_{\Bb}$ pour tout crible $\Bb$
couvrant $\Xx /X$ (voir \S 4). Par ailleurs, on observe que le morphisme
$$
\Gamma (\Xx /X , A'|_{\Xx /X}) \rightarrow A'(X)
$$
est un isomorphisme de $n$-cat\'egories de Segal, car $X$ est l'objet final
de la
cat\'egorie $\Xx /X$ et $\Gamma$ est simplement le foncteur des sections
globales pour
les pr\'efaisceaux.

L'\'enonc\'e suivant n'est probablement pas vrai pour le cas des $n$-pr\'echamps
de Segal m\^eme pour $n=0$; cela tient au fait que
le lemme \ref{pourloc} (corrig\'e dans {\tt v3}) n'est vrai que pour 
le cas $n$-tronqu\'e. Il semblerait qu'on obtient un contre-exemple
\`a \ref{newcritere} en consid\'erant un espace topologique $X$ ayant 
deux ouverts  
non-vides $U$ et $V$
et aucun autre \`a part les ouverts de $U\cap V$ introduits par la condition
r\'ecurrente $U\cap V \cong X$. 

\begin{proposition}
\label{newcritere}
Soit $A$ un $n$-pr\'echamp (non de Segal) sur $\Xx$ dont les valeurs sont
des $n$-cat\'egories. Fixons une cofibration triviale pour la
topologie grossi\`ere $A\rightarrow A'$ avec $A'$ fibrant pour la topologie
grossi\`ere.
Alors $A$ est un $n$-champ si et seulement si:
\newline
(c)\,\,
pour tout $X\in \Xx$ et tout crible $\Bb \subset \Xx /X$ de la topologie,
le morphisme
$$
A'(X)=\Gamma (\Xx /X, A'|_{\Xx /X}) \rightarrow \Gamma (\Bb , A'|_{\Bb})
$$
est une \'equivalence de $n$-cat\'egories.
\end{proposition}
{\em Preuve:}
Supposons que $A$ est un $n$-champ. Quitte \`a remplacer
$A$ par $A'$, on peut
aussi supposer que $A$ est fibrant pour la topologie grossi\`ere, et
donc m\^eme
$\Gg$-fibrant (voir \ref{stuff}). On va montrer que, dans ce cas, le
morphisme
de (c) est une fibration triviale; pour cela on montrera qu'il poss\`ede la
propri\'et\'e de rel\`evement par rapport \`a toute cofibration $E\rightarrow
E'$ de $n$-pr\'ecats.

On commence par observer que si un morphisme $U\rightarrow
V$ de $n$-pr\'echamps au-dessus de $\Xx /X$ est une \'equivalence
au-dessus de chaque objet d'un crible $\Bb \subset \Xx /X$ couvrant $X$,
alors c'est
une $\Gg$-\'equivalence faible. En effet, c'est facile \`a voir pour $n=0$
et pour $n>0$ on peut proc\'eder par r\'ecurrence sur $n$ en utilisant
directement
la d\'efinition de $\Gg$-\'equivalence faible (\S 3). On applique cela au
morphisme
$$
j:E'_{\Bb} \cup ^{E_{\Bb}} E_{\Xx /X} \rightarrow E'_{\Xx /X},
$$
o\`u $E\hookrightarrow E'$ est la cofibration triviale de $n$-pr\'ecats de
Segal donn\'ee, et $E_{\Bb}$ est le
$n$-pr\'echamp obtenu en appliquant $p_{!}$ au $n$-pr\'echamp
constant \`a valeurs $E$ sur $\Bb$, $p$ d\'esignant
l'inclusion $\Bb \rightarrow \Xx /X$.
On obtient ainsi que $j$ est une cofibration $\Gg$-triviale. La propri\'et\'e de
rel\`evement de $A\rightarrow \ast$ pour le morphisme $j$ se traduit en la
propri\'et\'e de rel\`evement du morphisme
$$
\Gamma (\Xx /X, A|_{\Xx /X}) \rightarrow \Gamma (\Bb , A|_{\Bb})
$$
par rapport \`a la cofibration $E\rightarrow E'$. Ceci prouve que le morphisme
en question est une fibration triviale, ce qui donne la condition (c) pour $A$.

Pour l'autre sens, supposons que $A$ satisfait la condition (c).
On note que si
$A$ est fibrant pour la topologie grossi\`ere, la condition (c) implique
que $A$
poss\`ede la
propri\'et\'e de rel\`evement pour toute cofibration de la forme
$cof(U\hookrightarrow U', \Bb \subset \Xx /X)$ (voir la d\'efinition juste avant
la proposition \ref{naqualocaliser}).  D'apr\`es le
lemme \ref{pourloc}, ceci implique
que $A$ poss\`ede la propri\'et\'e de rel\`evement pour toute cofibration
$\Gg$-triviale, i.e. que $A$ est $\Gg$-fibrant. La condition (c)
implique donc que $A$
est un $n$-champ.
\eop

{\em Remarques:} On indique ici un argument alternatif de r\'ecurrence sur $n$
qui permet facilement de r\'eduire la proposition au cas $n=0$. Supposons que
$n\geq 1$ et qu'on a d\'ej\`a d\'emontr\'e la proposition pour les
$n-1$-pr\'echamps de Segal. Si $\Bb$ est un crible couvrant $\Xx /X$,
alors
pour tous $x,y\in
A'_0(X)$ on a
$$
\Gamma (\Bb  , A'_{1/}(x,y)|_{\Bb})= \Gamma (\Bb , A'|_{\Bb})_{1/}(x,y).
$$
Au vu de la d\'efinition de l'\'equivalence entre $n$-cat\'egories de Segal,
la condition (c) pour $A'$ est manifestement \'equivalente \`a la
conjonction de la condition
(b) de \ref{critere} avec la condition (c) pour les $A'_{1/}(x,y)$. Par
hypoth\`ese de r\'ecurrence, cette derni\`ere est \'equivalente \`a la condition
(a) de \ref{critere}. On obtient donc bien que la condition (c) est
\'equivalente aux
conditions (a) plus (b) de \ref{critere}.
Pour $n=0$ i.e. pour les pr\'efaisceaux d'ensembles il y a un argument facile.

Dans le cas des pr\'efaisceaux simpliciaux, 
le fait que la condition (c) est n\'ecessaire
a fait l'objet d'une note de Toen \cite{ToenDescente}.
Pour tout $n$ la condition est n\'ecessaire pour les $n$-pr\'echamps de Segal.

Il est probable qu'en rempla\c{c}ant le lemme \ref{pourloc} par
la version  en termes d'hyper-recouvrements
d\'emontr\'ee r\'ecemment par Dugger \cite{Dugger}, 
on obtiendrait une version n\'e\-ces\-saire et suffisante de ce
crit\`ere (c) en termes d' hyper-recouvrements. Nous laissons 
ouverte cette int\'eressante question.

\subnumero{Les protochamps}

Pour les $1$-champs ordinaires (non de Segal), il appara\^{\i}t chez
Giraud \cite{GiraudThese} et Laumon-Moret-Bailly \cite{LaumonMB}
une notion de ``pr\'echamp'': en leur sens, un pr\'echamp
est un de nos pr\'echamps qui satisfait en outre la condition (a) de
\ref{critere}. Dans le cas des pr\'efaisceaux, ceci correspond \`a exiger la
premi\`ere des deux conditions de la d\'efinition habituelle des faisceaux
(Houzel appelait ca des semi-faisceaux mais ``semi-champ'' ne sonne pas bien) et
la terminologie de \cite{GiraudThese} \cite{LaumonMB} n'est donc pas compatible
avec la terminologie habituelle concernant les pr\'efaisceaux;
c'est pour cela que nous avons pr\'ef\'er\'e garder le terme
``pr\'echamp'' pour des objets d\'epourvus de relation avec la topologie $\Gg$.
Cependant, la notion de pr\'echamp satisfaisant la condition (a) de
\ref{critere} est utile, et nous adoptons la d\'efinition
suivante:

\begin{definition}
Soit $n\geq 1$.
Un $n$-pr\'echamp de Segal $A$ sera appel\'e {\em protochamp} si les $A(X)$ sont
des $n$-cat\'egories de Segal, et si $A$ satisfait la condition (a) de
\ref{critere}, \`a savoir que pour tout $X\in \Xx$ et tous $x,y\in A(X)_0$, le
$n-1$-pr\'echamp de Segal $A_{1/}(x,y)$ est un $n-1$-champ de Segal sur
$\Xx /X$.
\end{definition}

\begin{lemme}
Soit $A$ un $n$-pr\'echamp de Segal dont les valeurs $A(X)$ sont des
$n$-ca\-t\'e\-go\-ries de Segal. Soit $f:A\rightarrow A'$ un champ
associ\'e (\S 9). Alors $A$ est un protochamp si et seulement si le
morphisme $f$
est pleinement fid\`ele (objet-par-objet).
\end{lemme}
{\em Preuve:}
Si $f$ est pleinement fid\`ele alors, comme les $A'_{1/}(fx,fy)$  sont des
$n-1$-champs de Segal, il en est de m\^eme des $A_{1/}(x,y)$. Dans l'autre sens,
supposons que $A$ est un protochamp, et soient $x,y\in A(X)_0$.
Alors le morphisme $A_{1/}(x,y)\rightarrow A'_{1/}(fx,fy)$ est une
$\Gg_X$-\'equivalence faible entre $n-1$-champs de Segal sur $\Xx /X$,
donc par \ref{morphentrechamps}, c'est une \'equivalence objet-par-objet.
Donc $f$ est pleinement fid\`ele.
\eop

La prochaine proposition est l'analogue d'un r\'esultat de \cite{flexible}. On
laisse la d\'e\-mon\-stra\-tion en ``exercice'' au lecteur.

\begin{proposition}
\label{exo}
Soit $F$ un $n$-protochamp de Segal sur $\Xx$, et d\'efinissons un
$n$-pr\'e\-champ de Segal $G$ par
$$
G(X):= {\rm colim} _{\Bb \subset \Xx /X}\left( \lim _{\leftarrow , \Bb ^o}
F|_{\Bb} \right)
$$
(ici on utilise les limites et colimites de \cite{limits}, voir \S 14 pour
comment les calculer; la colimite est prise sur l'ensemble filtrant des cribles
couvrant $X$). Alors $F\rightarrow G$ est le champ associ\'e.
\end{proposition}
\eop

Rappelons qu'on peut consid\'erer une $n$-cat\'egorie comme $n$-cat\'egorie de
Segal en
consid\'erant les ensembles de $n$-morphismes
comme des ensembles simpliciaux $0$-tronqu\'es. Les $n$-cat\'egories
de Segal qui correspondent ainsi aux $n$-cat\'egories sont exactement celles
qui sont {\em $n$-tronqu\'es} i.e. pour lesquelles
on a $A\cong A^{int, n}\cong \tau
_{\leq n}A$. Le lemme suivant indique qu'on peut utiliser cette correspondance
pour traduire les r\'esultats sur les $n$-pr\'echamps de Segal en
r\'esultats sur les $n$-pr\'echamps (non de Segal).
Dans ce lemme, l'hypoth\`ese selon laquelle les $A(X)$
doivent \^etre des $n$-cat\'gories de Segal est bien n\'ecessaire
(car l'op\'eration $SeCat$ ne conserve pas la propri\'et\'e d'\^etre
tronqu\'e).

\begin{corollaire}
\label{preservetronque}
Si $A$ est un pr\'efaisceau de $n$-cat\'egories de Segal qui est $n$-tronqu\'e
(i.e. $A\cong \tau _{\leq n}A$) et si $A\rightarrow A'$ est un champ associ\'e,
alors $A'$ est aussi $n$-tronqu\'e.
\end{corollaire}
{\em Preuve:} Ce serait une cons\'equence imm\'ediate de
la proposition \ref{exo}.
Comme on n'a pas fourni de preuve de
cette proposition, on indique ici une d\'emonstration alternative du corollaire.
On peut supposer qu'on a d\'evelopp\'e parall\`element \`a notre discussion des
$n$-pr\'echamps de Segal, une discussion analogue pour les $n$-pr\'echamps
non de Segal.
Par exemple,  pour un $n$-pr\'echamp non de Segal $T$, on peut d\'efinir le
$n$-champ associ\'e $f:T\rightarrow T'$. Ce morphisme est une
$\Gg$-\'equivalence
faible de $n$-pr\'echamps.  Notons $Ind(T)$ (resp. $Ind(T')$) les
$n$-pr\'echamps de Segal induits.
Supposons que les $T(X)$ sont des
$n$-cat\'egories.
Le fait  que $f$ soit une $\Gg$-\'equivalence de
$n$-pr\'echamps implique que $Ind(f)$ est une $\Gg$-\'equivalence de
$n$-pr\'echamps de Segal. D'autre part, d'apr\`es \ref{critere} et son analogue
pour les $n$-pr\'echamps non de Segal, $Ind(A')$ est un $n$-champ
de Segal. En effet, avec la condition (a) on se r\'eduit par r\'ecurrence au cas
o\`u $n=0$; pour celui-ci, la question est de savoir si un faisceau
d'ensembles est un
$0$-champ de Segal, ce qui est le cas (on peut le voir en appliquant encore
deux fois \ref{0critere}). Ceci d\'emontre que $Ind(A')$ est un $n$-champ de
Segal. Alors par d\'efinition c'est le $n$-champ de Segal associ\'e \`a
$Ind(A)$. Dans la situation du corollaire, on part d'un $n$-pr\'echamp
de Segal $A$ qui est objet-par-objet $n$-tronqu\'e. Il s'ensuit que
$A$ est \'equivalent \`a
un pr\'echamp $Ind(T)$ o\`u $T$ est un $n$-pr\'echamp non de Segal
(et dont les valeurs sont des $n$-cat\'egories) et la discussion
pr\'ec\'edente prouve le
corollaire.
\eop

{\em Exercice:} Si $A$ est un $n$-pr\'echamp non de Segal (ou de Segal mais
$n$-tronqu\'e) alors on obtient le champ associ\'e  $A'$  \`a partir de $A$ par
application $n+2$ fois de la construction de la proposition \ref{exo}.

\numero{Cat\'egories de mod\`eles internes}

\label{cmipage}

Nous voulons construire le $n+1$-pr\'echamp de Segal de modules des
$n$-champs de Segal. Pour ceci on utilise la notion suivante de {\em cat\'egorie
de mod\`eles ferm\'ee interne}. Soit $M$ une cat\'egorie de mod\`eles ferm\'ee.
On dira que $M$ est {\em interne} si elle v\'erifie les axiomes suivants:
\newline
IM(a)---$M$ admet un $Hom$ interne, i.e. le foncteur $X\mapsto M_1(X\times A,
B)$ est repr\'esentable par un objet $\underline{Hom} (A,B)$;
\newline
IM(b)---le produit direct commute aux colimites finies, en particulier $A\times
\emptyset = \emptyset$ pour l'objet initial $\emptyset$;
\newline
IM(c)---si $A\rightarrow A'$ et $B\rightarrow B'$ sont des cofibrations alors
le morphisme naturel
$$
{\rm colim} \left(
\begin{array}{ccc}
A\times B & \rightarrow & A'\times B \\
\downarrow &&  \\
A\times B' & &
\end{array}
\right) \rightarrow A'\times B'
$$
est aussi une cofibration; et
\newline
IM(d)--- l'\'equivalence faible est stable par produit direct.

Ces axiomes sont tr\`es proches des axiomes d'une {\em cat\'egorie de
mod\`eles ferm\'ee mo\-no\^{\i}\-da\-le}
au  sens de Hovey \cite{HoveyMonMod}. Essentiellement, on
demande que le produit direct d\'efinisse une structure mono\^{\i}dale
au sens
de \cite{HoveyMonMod}, et  on demande en outre que le $\underline{Hom}$
interne existe. Cette derni\`ere condition est plus ou moins
\'equivalente \`a
l'existence d'une certaine colimite, mais cette
colimite n'est peut-\^etre pas
``petite'' (i.e. index\'ee par un ensemble et non une classe)
et nous ne savons
pas si la simple condition d'existence des (petites) colimites suffit pour
garantir
IM(a).

La proposition suivante redonne les axiomes IM1-IM4 de
\cite{nCAT} \S 11, et nous consid\'erons d\'esormais
les axiomes IM(a)--IM(d) comme constituant la version
``officielle'' de la notion de cmf interne.

\begin{proposition}
\label{interne}
Soit $M$ une cat\'egorie de mod\`eles ferm\'ee interne.
\newline
(1) \,\, Si $A$ est cofibrant et $B$ fibrant alors $\underline{Hom}(A,B)$ est
fibrant; \newline
(2)\,\, Si $A\rightarrow A'$ est une cofibration (resp. cofibration triviale)
d'objets cofibrants, et si $B'\rightarrow B$ est une fibration (resp. fibration
triviale) entre objets fibrants, alors le morphisme
$$
\underline{Hom}(A',B')\rightarrow
\underline{Hom}(A',B) \times _{\underline{Hom}(A,B) }\underline{Hom}(A,B')
$$
est une fibration (resp. fibration triviale);
\newline
(3)\,\, Le $\underline{Hom}$ interne transforme tout coproduit
cofibrant dans
son premier argument
(resp. produit fibr\'e  fibrant dans son second argument) en un produit
fibr\'e fibrant.
\end{proposition}
La preuve est facile et laiss\'ee au lecteur (les arguments se trouvent au \S
7 de \cite{nCAT}).
\eop

On a une version ``homotopique'' de la notion de $Hom$ interne. Soit $A$ une
cat\'egorie simpliciale, et soient $x,y\in A$. On suppose que $A$ admet des
produits directs homotopiques; soit $h\in A$ et $f: h\times x\rightarrow y$ un
morphisme (i.e. on choisit un objet $h\times x$ muni de morphismes vers $h$
et $x$
qui repr\'esente le produit direct). On dit que $(h,f)$ est un $Hom$ interne
homotopique
de $x$ \`a $y$, si le morphisme (essentiellement bien d\'efini)
$$
A_{1/}(u,h)\rightarrow A_{1/}(u\times x, y)
$$
est une \'equivalence (ici encore on choisit un objet $u\times x$).

\begin{proposition}
\label{hominterne}
Soit $M$ une cat\'egorie de mod\`eles ferm\'ee interne. Alors la cat\'egorie
simpliciale $L(M)$ admet un $Hom$ interne au sens homotopique ci-dessus.
\end{proposition}
{\em Preuve:}
On utilise la technique des r\'esolutions cosimpliciales cofibrantes
pour calculer les complexes de fonctions
(cf \cite{Hirschhorn} \cite{DHK} \cite{DwyerKan3}
\cite{JardineGoerssBook} \cite{Reedy}).
Fixons $x,y\in M$ avec $x$ cofibrant et $y$ fibrant, et notons
$h:=\underline{Hom}(x,y)$ le $\underline{Hom}$ interne de $M$ (IM(a)). On a un
morphisme $h\times x\rightarrow y$. Ici et plus loin il convient de
remarquer que
le produit direct d'objets de $M$ est aussi le produit direct homotopique
dans $L(M)$, d'apr\`es IM(d). Fixons $u\in M$ et choisissons une
r\'esolution cosimpliciale cofibrante de Reedy ${\bf u}\rightarrow u$
(cf \cite{Hirschhorn} \cite{DHK} \cite{DwyerKan3} \cite{JardineGoerssBook}
\cite{Reedy}). On a alors l'\'egalit\'e entre ensembles simpliciaux
$$
M_{1/} ({\bf u}, h)= M_{1/}(x\times {\bf u}, y).
$$
Or, comme $h$ est fibrant (Proposition \ref{interne} (1)), on a aussi
(d'apr\`es {\em loc cit.})
$$
M_{1/} ({\bf u}, h)\cong L(M)_{1/}(u,h).
$$
D'autre part, les axiomes IM(b) et IM(c) impliquent que $x\times {\bf u}$ est
une r\'esolution cosimpliciale cofibrante de Reedy pour $x\times u$, donc
(encore d'apr\`es {\em loc cit.}) on a l'\'equivalence
$$
M_{1/} (x\times {\bf u}, y)\cong L(M)_{1/}(x\times u, y).
$$
Ces deux \'equivalences impliquent par d\'efinition que $h$ est un $Hom$ interne
homotopique de $x$ \`a $y$.
\eop

\subnumero{La cat\'egorie $M_f$-enrichie associ\'ee}

Soit $M$ une cat\'egorie de mod\`eles ferm\'ee interne. Notons $\ast$
son objet final. On obtient un foncteur $\sigma : M\rightarrow Sets$
en posant $\sigma (X):= M_1(\ast , X)$. Cette construction commute
aux produits directs.
On a la formule tautologique
$$
M_1(A,B) = \sigma \underline{Hom} (A,B).
$$
Le $\underline{Hom}$ interne permet de munir $M_{cf}$ d'une structure de
cat\'egorie
enrichie dans $M_f$, via le foncteur $\sigma$. Autrement dit, pour
tous $A,B\in M_{cf}$ on a $\underline{Hom} (A,B)\in M_f$ et
$M_{cf}(A,B)=\sigma (\underline{Hom}
(A,B))$;
pour tous $A,B,C$ on a une composition
$$
\underline{Hom} (A,B)\times \underline{Hom} (B,C)\rightarrow
\underline{Hom}(A,C)
$$
qui est associative (nous n\'egligeons ici les probl\`emes li\'es aux
contraintes
d'associativit\'e qui sont des isomorphismes naturels). Cette
composition est transform\'ee par
$\sigma$ en la composition de $M_{cf}$.

Si on sait en outre que les $\underline{Hom} (A,B)$ sont cofibrants (c'est par
exemple
le cas quand les
cofibrations de $M$ sont les injections), alors $M_{cf}$ est
auto-enrichie par son $Hom$ interne. En particulier, la structure
$INT(M_{cf})$ ne depend que de la cat\'egorie $M_{cf}$ dans ce cas.

On obtient aussi ce qu'on pourrait appeler une {\em $M_f$-cat\'egorie de Segal
stricte} et qu'on notera par $INT(M_{cf})$,
en posant
$$
INT(M_{cf})_{p/}(A_0,\ldots , A_p):=
\underline{Hom}(A_0,A_1)\times \ldots \times \underline{Hom}(A_{p-1}, A_p).
$$
Cet objet consiste plus pr\'ecisement en un ensemble
$INT(M_{cf})_0$ d'objets, et pour toute suite $A_0,\ldots , A_p$ d'objets,
un objet $INT(M_{cf})_{p/}(A_0,\ldots , A_p)$ de $M_f$ (avec les morphismes
habituels de fonctorialit\'e) tel que les morphismes de Segal soient des
isomorphismes. Ceci est \'equivalent bien s\^ur \`a la donn\'ee d'une
cat\'egorie $M_f$-enrichie, et nous allons utiliser la notation
$INT(M_{cf})$ indiff\'eremment pour la cat\'egorie $M_f$-enrichie ou pour la
$M_f$-cat\'egorie de Segal stricte.

\subnumero{Structure interne pour les $n$-cat\'egories de Segal}

Dans le cas o\`u $M$ est une cat\'egorie de pr\'efaisceaux d'ensembles sur
une certaine
cat\'egorie, et si les cofibrations sont les morphismes injectifs
(ou m\^eme seulement injectifs au-dessus d'un sous-ensemble d'objets dans
la cat\'egorie de base) alors les axiomes IM(a), IM(b) et IM(c) sont
automatiquement v\'erifi\'es.
Tel est le cas de la cat\'egorie de mod\`eles ferm\'ee des $n$-pr\'ecats de
\cite{nCAT} ou sa variante pour les $n$-pr\'ecats de Segal d\'efinie plus haut
dans le pr\'esent travail \ref{cmf}. Tel est aussi le cas pour la cat\'egorie de
mod\`eles ferm\'ee
des $n$-pr\'echamps de Segal \ref{SeCmf}, car un $n$-pr\'echamp de Segal n'est
rien d'autre qu'un
pr\'efaisceau d'ensembles sur $\Theta ^{n+1} \times \Xx$ (et les
cofibrations sont les injections de pr\'efaisceaux).
La seule difficult\'e dans ces exemples consiste \`a v\'erifier
IM(d).

Ceci a \'et\'e fait
pour les $n$-pr\'ecats dans \cite{nCAT} Theorem 5.1. Comme pour l'existence de
la structure de cat\'egorie de mod\`eles ferm\'ee, on peut soit recopier la
d\'emonstration, soit utiliser ce r\'esultat via $A\mapsto \Pi _m \circ A$, pour
obtenir le m\^eme r\'esultat pour les $n$-pr\'ecats de Segal. On obtient:

\begin{lemme}
La cmf $nSePC$ des $n$-pr\'ecats de Segal de \ref{cmf} est
interne.
\end{lemme}
\eop

Si on applique la discussion pr\'ec\'edente \`a  la cat\'egorie de mod\`eles
ferm\'ee $nSePC$ des $n$-pr\'ecats de Segal on obtient une $1$-cat\'egorie
stricte $INT(nSePC_{cf})$
enrichie sur $nSePC_f$; et en
prenant le nerf de celle-ci on obtient une $nSePC_f$-cat\'egorie de Segal, i.e.
une $n+1$-cat\'egorie de Segal qu'on note
$$
nSeCAT:= INT(nSePC_{cf}).
$$
Insistons bien sur le fait que les objets de
$nSeCAT$ sont les $n$-cat\'egories de Segal fibrantes, et pour toute
suite d'objets $A_0,\ldots , A_p$ on a
$$
nSeCAT_{p/}(A_0,\ldots , A_p):= \underline{Hom}(A_0, A_1)\times \ldots
\times \underline{Hom}(A_{p-1}, A_p).
$$
La m\^eme construction pour la $n+1$-cat\'egorie $nCAT$ des $n$-cat\'egories
fibrantes, a \'et\'e donn\'ee dans \cite{nCAT} \S 7.

On observe que $nSeCAT$ (resp. $nCAT$) n'est pas fibrante, et
on utilisera souvent un
remplacement fibrant qu'on notera $nSeCAT'$ (resp. $nCAT'$).

{\em Remarque ensembliste:}
Les objets de $(nSeCAT)_0$ forment une classe, et
techniquement parlant, $nSeCAT$ n'est
pas une $n+1$-cat\'egorie. Pour contourner ce probl\`eme, on convient de fixer
un ensemble (tr\`es grand!) et de ne retenir comme objet de $nSeCAT$
que les $n$-cat\'egories dont
l'ensemble sous-jacent (i.e. r\'eunion des $A_M$) est contenu dans cet ensemble.
Avec cette convention, le remplacement fibrant $nSeCAT'$ existe. Le lecteur
v\'erifiera que nos constructions ne nous \'eloignent pas de ce cadre: quand on
parle de limites ou colimites, par exemple,
elles sont toujours prises suivant un ensemble d'indices
dont on peut borner \`a l'avance la
puissance; il suffit donc de choisir notre ensemble de d\'epart
suffisamment grand pour que toutes les limites et colimites qu'on va prendre,
existent; de fa\c{c}on similaire, quand on utilise l'argument du petit objet,
on connait toujours \`a l'avance une borne transfinie du nombre
d'\'etapes n\'ecessaires.

\subnumero{Structure interne pour les $n$-champs de Segal}

On va maintenant adapter aux $n$-champs de Segal la construction vue pour les
$n$-cat\'egories de Segal.

\begin{lemme}
La cmf $nSePCh$ des $n$-pr\'echamps de Segal sur un site $\Xx$, est interne.
\end{lemme}
{\em Preuve:} comme on l'a remarqu\'e plus haut, la seule
difficult\'e concerne IM(d)
qu'on prouve par r\'ecurrence sur $n$. Pour $n=0$ la propri\'et\'e en question
r\'esulte facilement de la d\'efinition de l'\'equivalence d'Illusie (en notant
que le passage au faisceau associ\'e est
compatible avec le produit direct
des pr\'efaisceaux d'ensembles, voir par exemple \cite{Rezk}). On suppose donc
$n\geq 1$ et le r\'esultat connu pour $n-1$.  Soit $f:A\rightarrow A'$ une
$\Gg$-\'equivalence de $n$-pr\'echamps de Segal, et soit $B$ un
$n$-pr\'echamp de
Segal.  On peut supposer (gr\^ace au r\'esultat correspondant pour les
$n$-pr\'ecats
de Segal) que les $A(X)$, $A'(X)$ et $B(X)$ sont des $n$-cat\'egories de Segal.
On a
$$
\tau _{\leq 0} (A\times B)= \tau _{\leq 0}(A)\times \tau _{\leq 0} (B)
$$
(et de m\^eme pour $A'$)
donc le morphisme
$$
\tau _{\leq 0} (A\times B)\rightarrow \tau _{\leq 0} (A'\times B)
$$
induit une surjection sur les faisceaux d'ensembles associ\'es. D'autre part,
pour $(x,a)$ et $(y,b)$ dans $(A\times B)_0(X)$, on a
$$
(A\times B)_{1/}((x,a), (y,b))= A_{1/}(x,y)\times B_{1/}(a,b)
$$
et de m\^eme pour $A'$; on obtient donc (\`a partir du r\'esultat pour
$n-1$) que le morphisme induit par $f$
$$
(A\times B)_{1/}((x,a), (y,b))\rightarrow
(A'\times B)_{1/}((fx,a), (fy,b))
$$
est une $\Gg$-\'equivalence. Par d\'efinition (voir \S 3) cela veut dire que
$f$ est une $\Gg$-\'equivalence. Ceci donne IM(d).
\eop

Maintenant,
\`a partir de $M=nSePCh$, on peut construire le $n+1$-pr\'echamp de Segal
des $n$-champs de Segal fibrants, not\'e
$$
nSe\underline{CHAMP}(\Xx ):= INT(nSePCh_{cf}).
$$
Les objets de $nSe\underline{CHAMP}(\Xx )$ au-dessus de $X\in \Xx$ sont les
$n$-pr\'echamps de Segal fibrants sur $\Xx /X$. Si $A_0,\ldots , A_p$ sont de
tels objets, alors on d\'efinit sur $\Xx /X$ le $n$-pr\'echamp de Segal
$$
nSe\underline{CHAMP}(\Xx )_{p/}(A_0,\ldots , A_p)=
\underline{Hom}(A_0,A_1)\times \ldots \times \underline{Hom}(A_{p-1}, A_p).
$$

En appliquant le foncteur des sections globales $\Gamma$ \`a
$nSe\underline{CHAMP}(\Xx )$, on obtient {\em la
$n+1$-cat\'egorie de Segal des $n$-champs de Segal sur $\Xx$},
qu'on pourrait noter
$$
nSeCHAMP(\Xx ).
$$
En fait, si on analyse la phrase pr\'ec\'edente, on constate que
les objets de $nSeCHAMP(\Xx )$
sont les $n$-pr\'echamps de Segal $A$ tels que $A|_{\Xx /X}$ soit
fibrant pour tout $X\in \Xx$. Si $\Xx$ n'a pas d'objet final, cette condition
peut \^etre diff\'erente que la condition (plus forte) que $A$ soit fibrant
(noter le contre-exemple dans \S 6),
et nous allons l\'eg\`erement modifier la d\'efinition pour ne garder dans
$nSeCHAMP(\Xx )$ que les seuls objets qui sont fibrants sur $\Xx$.

En conclusion, on pr\'ef\`ere noter
$nSeCHAMP(\Xx )$
la $n+1$-cat\'egorie dont les objets
sont les $n$-pr\'echamps de Segal $\Gg$-fibrants sur le site
$\Xx$. Si $A_0,\ldots , A_p$ sont de
tels objets, alors on  prend comme pr\'ec\'edemment pour
$nSeCHAMP(\Xx )_{p/}(A_0,\ldots , A_p)$
le produit de $n$-pr\'ecats
$\Gamma \underline{Hom}(A_0,A_1)\times \ldots \times
\Gamma \underline{Hom}(A_{p-1}, A_p).$

Il convient de rappeler, pour comprendre cette d\'efinition, que
$\Gamma \underline{Hom}(A,B)$ peut \^etre d\'efini comme le repr\'esentant du
foncteur
$nSePC\rightarrow Ens$
$$
E\mapsto Hom _{nSePCh}(A\times \underline{E},B),
$$
o\`u $\underline{E}$ est le $n$-pr\'echamp de Segal constant \`a valeurs $E$.

Il est imm\'ediat que pour tout $X\in \Xx$, on a
$$
nSe\underline{CHAMP}(\Xx )(X) = nSeCHAMP(\Xx /X),
$$
et ceci pourrait servir de d\'efinition de $nSe\underline{CHAMP}(\Xx )$.

Enfin, on peut faire ici la m\^eme ``remarque ensembliste''
que plus haut (qui explique
pourquoi on se permet de prendre des remplacements fibrants
$nSe\underline{CHAMP}(\Xx )'$).

\subnumero{Calcul avec la structure de type HBKQ}

Nous indiquons ici comment utiliser la structure du th\'eor\`eme
\ref{cmfHirschho} pour calculer $\Gamma \underline{Hom}$ et
$\underline{Hom}$ (ce dernier dans le cas o\`u $\Xx$ admet des produits
fibr\'es).

\begin{lemme}
\label{hbkqinterne}
Si le site $\Xx$ admet des produits directs, alors la cmf de type HBKQ de
\ref{cmfHirschho} est interne.
\end{lemme}
{\em Preuve:} Il suffit de remarquer
que le produit de deux additions libres de cellules au-dessus de $X$
et de $Y$, est une addition libre au-dessus de $X\times Y$.
\eop

{\em Remarque:} Si on veut appliquer ce lemme sur les sites $\Xx /X$ alors il
faut que ceux-ci admettent des produits directs, i.e. que $\Xx$ admet des
produits fibr\'es.

Au vu de ce lemme, on pourrait, dans le cas o\`u $\Xx$ admet des produits
directs, d\'efinir une version HBKQ de $nSeCHAMP(\Xx )$ et
$nSe\underline{CHAMP}(\Xx )$. On signale maintenant qu'on obtiendrait
essentiellement les m\^emes objets que plus haut.

On commencera par montrer une propri\'et\'e valide m\^eme si $\Xx$ n'admet
pas tous les produits directs.

Soit $E$ est une $n$-pr\'ecat de Segal, et notons
$\underline{E}$ le $n$-pr\'echamp de Segal constant sur $\Xx$ de valeur $E$.
Si $A\rightarrow A'$ est la cofibration engendr\'ee librement par l'addition
d'une cellule au-dessus de $X\in \Xx$ (i.e. engendr\'ee par une cofibration de
$n$-pr\'ecats de Segal $A(X)\rightarrow C$) alors $A\times
\underline{E}\rightarrow A'\times \underline{E}$ est aussi engendr\'ee
librement par addition de la cellule $A(X)\times E \rightarrow C\times E$.
Il s'ensuit (du fait que le produit direct avec $\underline{E}$ pr\'eserve
les compositions transfinies et les r\'etractions) que si
$$
A\rightarrow A'
$$
est une cofibration pour la structure de HBKQ, alors
il en est de m\^eme de
$$
A\times
\underline{E}\rightarrow A'\times \underline{E}.
$$
Soit maintenant $A$ cofibrant et $B$ fibrant pour la structure de
HBKQ (Th\'eor\`eme \ref{cmfHirschho}). La $n$-pr\'ecat de
Segal
$\Gamma \underline{Hom}(A,B)$ repr\'esente le foncteur $nSePC\rightarrow Ens$
$$
E\mapsto Hom _{nSePCh}(A\times \underline{E}, B).
$$
Le fait que les $A\times \underline{E}$ soient cofibrants (pour
\ref{cmfHirschho})
implique que si $B\rightarrow B'$ est une \'equivalence entre deux
$n$-pr\'echamps de Segal fibrants pour \ref{cmfHirschho}, alors le morphisme
induit
$$
\Gamma \underline{Hom}(A,B)\rightarrow \Gamma \underline{Hom}(A,B')
$$
est une \'equivalence de $n$-pr\'ecats de Segal (les deux sont des
$n$-pr\'ecats de Segal fibrantes). En particulier
on peut choisir $B'$
fibrant par rapport \`a la structure du th\'eor\`eme \ref{cmf},
auquel cas on obtient la
``bonne'' $n$-cat\'egorie de Segal des morphismes de $A$ vers $B$. La conclusion
est le lemme suivant, qui dit que $A$ et $B$ peuvent etre utilis\'es pour
calculer les $n$-cat\'egories de Segal de morphismes dans $nSeCHAMP (\Xx )$.

Notons avant d'\'enoncer le lemme, que $B$ est fibrant pour la structure
\ref{cmfHirschho} si et seulement si pour tout $X\in \Xx$, $B(X)$ est fibrant.
Ceci est donc une condition ponctuelle (v\'erifi\'ee par exemple par tout
foncteur $\Xx ^o \rightarrow nSeCAT$).

\begin{lemme}
\label{bonhomHirschho}
Soient $A$ cofibrant et $B$ fibrant pour la structure de
HBKQ (Th\'eor\`eme \ref{cmfHirschho}). Soient $A'$ et $B'$
leurs remplacements fibrants pour la structure de \ref{cmf}. Alors
il y a une \'equivalence de $n$-cat\'egories de Segal
$$
nSeCHAMP(\Xx )_{1/}(A',B')=\Gamma \underline{Hom}(A',B')\cong
\Gamma \underline{Hom}(A,B).
$$
\end{lemme}
\eop

Signalons qu'on peut calculer de la m\^eme fa\c{c}on les $n$-pr\'echamps
de Segal de morphismes dans $nSe\underline{CHAMP}(\Xx )$, {\em \`a condition
que $\Xx$ admette des produits fibr\'es}; gr\^ace au lemme
\ref{hbkqinterne}:

\begin{lemme}
\label{bonhomHirschho2}
Supposons que le site $\Xx$ admet des produits fibr\'es.
Soient $A$ cofibrant et $B$ fibrant pour la structure de
HBKQ (Th\'eor\`eme \ref{cmfHirschho}). Soient $A'$ et $B'$
leurs remplacements fibrants pour la structure de \ref{cmf}. Alors
il y a une \'equivalence de $n$-cat\'egories de Segal
$$
nSe\underline{CHAMP}(\Xx )_{1/}(A',B')=\underline{Hom}(A',B')\cong
\underline{Hom}(A,B).
$$
\end{lemme}
\eop

Le lemme \ref{bonhomHirschho} nous aide \`a expliciter la notion de ``donn\'ee
de descente''. Soit $X\in \Xx$ et soit $\Bb \subset \Xx /X$ un crible
couvrant $X$. Soit
$$
{\cal D}_{\Bb} \rightarrow \ast _{\Bb}
$$
le morphisme de $n$-pr\'echamps de Segal (qui peuvent \^etre consid\'er\'es
comme
au-dessus de $\Xx /X$ ou $\Bb$) construit au \S 5.
Soit $A$ un $n$-pr\'echamp de Segal sur $\Xx$ fibrant pour la structure
\ref{cmfHirschho}, i.e. un pr\'efaisceau sur $\Xx$ \`a valeurs dans la
cat\'egorie
$nSePC_f$ des $n$-cat\'egories de Segal fibrantes.
Soit $A\rightarrow A'$ un remplacement fibrant pour la structure de \ref{cmf}
pour la topologie grossi\`ere. Alors on a l'\'equivalence
$$
\Gamma \underline{Hom} (\ast _{\Bb}, A'|_{\Bb})
\cong \Gamma \underline{Hom}({\cal D}_{\Bb}, A|_{\Bb}).
$$
En particulier en prenant $\tau _{\leq 0}$,
on obtient
le corollaire suivant.

\begin{corollaire}
\label{calculDescente}
Soit $A$ un $n$-pr\'echamp de Segal sur $\Xx$ fibrant pour la structure
\ref{cmfHirschho}, i.e. un pr\'efaisceau sur $\Xx$ \`a valeurs dans la
cat\'egorie
$nSePC_f$ des $n$-cat\'egories de Segal fibrantes.
Soit $X\in \Xx$ et  $\Bb \subset \Xx /X$ un crible couvrant $X$, et soit
$$
{\cal D}_{\Bb} \rightarrow \ast _{\Bb}
$$
la ``r\'esolution'' construite au \S 5. Alors l'ensemble
``des donn\'ees de descente''
$$
[\ast _{\Bb}, A|_{\Bb}]_{\Bb ^{\rm
gro}}
$$
se calcule comme l'ensemble de morphismes (dans $nSePCh$) de ${\cal
D}_{\Bb}$ vers
$A$, modulo la relation qui
identifie deux morphismes s'ils sont li\'es par un
morphisme ${\cal D}_{\Bb} \times \underline{\overline{I}}\rightarrow A$.
\end{corollaire}
{\em Preuve:}
Soit $A\rightarrow A'$ un remplacement fibrant pour la structure de \ref{cmf}
pour $\Xx ^{\rm gro}$.
L'ensemble des donn\'ees de descente se calcule commme
$$
[\ast _{\Bb}, A|_{\Bb}]_{\Bb ^{\rm
gro}}
=\tau _{\leq 0}\Gamma \underline{Hom} (\ast _{\Bb}, A'|_{\Bb}).
$$
Par l'\'equivalence ci-dessus ceci est isomorphe \`a
$$
\tau _{\leq 0}\Gamma \underline{Hom} ({\cal D}_{\Bb}, A|_{\Bb}).
$$
Comme
$\Gamma \underline{Hom} ({\cal D}_{\Bb}, A|_{\Bb})$ est une $n$-pr\'ecat de
Segal
fibrante, sa troncation $\tau _{\leq 0}$ se calcule comme l'ensemble des
objets
modulo la relation qui identifie
deux objets li\'es par un morphisme de source $\overline{I}$.
\eop

\subnumero{Relation avec la structure simpliciale}

D'apr\`es Dwyer-Kan \cite{DwyerKan3}, la ``structure simpliciale'' d'une
cat\'egorie de mod\`eles ferm\'ee est bien d\'efinie \`a homotopie pr\`es, et
on voudrait la comparer \`a la structure interne quand cette derni\`ere existe.
Il serait possible de formuler cette question en toute g\'en\'eralit\'e mais
par commodit\'e (et du fait que cela suffit pour nos applications) nous allons
seulement regarder le cas d'une cat\'egorie de mod\`eles ferm\'ee interne avec
structure simpliciale compatible. On commence par la
d\'efinition: si $M$ est une cat\'egorie de mod\`eles ferm\'ee interne, une
{\em structure simpliciale compatible} est un foncteur $R: \Delta \rightarrow M$
tel que la cat\'egorie simpliciale $M_{\ast}$
d\'efinie par
$$
Hom _{M_k}(x,y):= Hom _M(R(k) \times x, y)
$$
soit une cat\'egorie de mod\`eles simpliciale au sens de Quillen.
On peut aussi \'ecrire
$$
Hom_{M_{\ast}}(x,y) = \left( k\mapsto Hom _M(R(k)\times x, y) \right) .
$$

Le foncteur
$R$ s'\'etend par passage aux limites finies
en un foncteur (qu'on notera encore $R$) de la cat\'egorie
des ensembles simpliciaux n'ayant qu'un nombre fini de simplexes
non-d\'eg\'en\'er\'es
vers $M$. Dans les notations introduites par Quillen pour la structure
simpliciale, on a
$$
K\otimes x = R(K)\times x,
$$
et
$$
Hom (K, y) = \underline{Hom}(R(K), y)
$$
($\underline{Hom}$ est le $Hom$
interne de $M$). On peut aussi prendre ces formules comme
d\'efinition alternative de la structure simpliciale associ\'ee \`a $R$
(ce point de vue correspond \`a
la d\'efinition de cat\'egorie de mod\`eles simpliciale
donn\'ee par Goerss-Jardine \cite{JardineGoerssBook}).

Si $(M,R)$ est une cat\'egorie de mod\`eles ferm\'ee interne avec structure
simpliciale compatible, on obtient un foncteur $\Xi : M\rightarrow EnsSpl$
de $M$ vers les ensembles simpliciaux, en posant pour $x\in M$,
$$
\Xi (x)_p := Hom _M(R(p), x),
$$
ou en termes de foncteurs repr\'esentables,
$$
Hom _{EnsSpl}(K, \Xi (x)):= Hom _M(R(K), x).
$$
En d'autres termes $\Xi (x)$ est l'ensemble simplicial
$Hom_{M_{\ast}}(\ast , x)$ o\`u $\ast$ est l'objet final de $M$.
Le foncteur $\Xi$ est compatible aux produits directs. De ce fait,
\`a partir d'une cat\'egorie $M$-enrichie $C$, on peut fabriquer la
cat\'egorie $\Xi \circ  C$ qui est
enrichie sur les ensembles simpliciaux.

\begin{lemme}
\label{splinterne}
Soit $(M,R)$ une cat\'egorie de mod\`eles ferm\'ee interne avec structure
simpliciale compatible, et soit $INT(M_{cf})$ la cat\'egorie $M_f$-enrichie
construite plus haut. D'autre part, soit $M_{cf,\ast}$ la cat\'egorie
simpliciale des objets cofibrants et fibrants avec sa structure simpliciale
d\'etermin\'ee par $R$. Alors on a un isomorphisme de cat\'egories
simpliciales
$$
M_{cf, \ast} = \Xi \circ  INT(M_{cf}).
$$
Par cons\'equent, si $L(M)$ est la localis\'ee de Dwyer-Kan de $M$ par rapport
aux \'equivalences faibles, on a une \'equivalence de cat\'egories simpliciales
$$
L(M)\cong \Xi \circ INT(M_{cf}).
$$
\end{lemme}
{\em Preuve:}
Les objets de $M_{cf, \ast}$ et de $INT(M_{cf})$ sont les objets cofibrants et
fibrants de $M$. Pour deux tels objets $x,y$, on a par d\'efinition
$$
INT(M_{cf})_{1/}(x,y):= \underline{Hom}(x,y),
$$
et donc
$$
(\Xi \circ INT(M_{cf}))_{1/}(x,y):= \Xi (INT(M_{cf})_{1/}(x,y))
= \left( p\mapsto Hom (R(p), \underline{Hom}(x,y)) \right)
$$
$$
= \left( p\mapsto Hom (R(p)\times x,y) \right)
= Hom _{M_{cf, \ast}}(x,y).
$$
Ceci donne l'isomorphisme souhait\'e. La deuxi\`eme partie est cons\'equence
de l'\'equivalence $L(M) \cong M_{cf, \ast}$ voir \cite{DwyerKan2}.
\eop

On va appliquer tout ceci aux $n$-cat\'egories et $n$-champs de Segal.
On v\'erifie facilement qu'on peut munir la
cat\'egorie de mod\`eles ferm\`ee int\`erne $nSePC$  des $n$-pr\'ecats de Segal
d'une structure simpliciale compatible en prenant pour $R(p)$
la cat'egorie $\overline{I}^{(p)}$
ayant $p+1$ objets $0',\ldots , p'$
et un seul isomorphisme entre chaque paire d'objets. Soit
$\Xi : nSePC \rightarrow EnsSpl$ le foncteur correspondant.

\begin{lemme}
\label{sigma}
Si $A$ est une $n$-cat\'egorie de Segal fibrante (et automatiquement
cofibrante) alors il y a une \'equivalence naturelle d'ensembles simpliciaux
$$
\Xi (A) \stackrel{\cong}{\rightarrow} A^{int, 0}.
$$
\end{lemme}
{\em Preuve:}
On a un morphisme de $n$-pr\'ecats de Segal
$$
R\Xi (A) \rightarrow A.
$$
En notant que $SeCat(R\Xi (A))$ est un $n$-groupo\"{\i}de de Segal (comme
colimite homotopique de $n$-groupo\"{\i}des de Segal), et
en utilisant la notation $A^{int, 0'}$ du \S 2, on obtient l'\'equivalence
$$
R\Xi (A)\stackrel{\cong}{\rightarrow}
SeCat(R\Xi (A))^{int, 0'} \rightarrow A^{int, 0'}
$$
(il faut choisir une r\'etraction $SeCat(A)\rightarrow A$,
ce qui possible parce que $A$ est
fibrante). Le compos\'e est le morphisme ci-dessus. Si $\Re$ est la
r\'ealisation comme
ensemble simplicial, on a par d\'efinition $A^{int, 0}:= \Re A^{int, 0'}$,
d'o\`u un morphisme $$
\Re (R\Xi (A)) \rightarrow A^{int, 0}.
$$
Si on note $I^{(p)}$ la cat\'egorie avec $p+1$ objets $0',\ldots , p'$ et
un morphisme $i'\rightarrow j'$ pour tout $i\leq j$, alors
l'ensemble simplicial
$$
\Re (I^{(p)})= h(p)
$$
est le $p$-simplexe standard, et l'inclusion
$I^{(p)} \subset \overline{I}^{(p)}$ fournit donc un morphisme
$$
h(p)\rightarrow \Re (R(p))
$$
qui s'\'etend par passage aux
colimites en un morphisme $K\rightarrow \Re (R(K))$ pour tout
ensemble simplicial $K$. En particulier, on a un morphisme
$$
\Xi (A)\rightarrow \Re(R\Xi (A)),
$$
lequel compos\'e avec le morphisme ci-dessus donne
$$
\Xi (A)\rightarrow A^{int, 0}.
$$
Il nous faut que c'est une \'equivalence.

La construction $\Xi$ respecte les \'equivalences entre objets fibrants.
En utilisant l'\'egalit\'e $\Xi (A) = \Xi
(A^{int, 0'})$, on peut---quitte \`a remplacer $A$ par $A^{int, 0'}$--
supposer que $A$ est un $n$-groupo\"{\i}de de Segal. Ce
dernier correspond \`a un
espace topologique d'apr\`es le travail de
Tamsamani \cite{Tamsamani} (les techniques de Tamsamani
s'appliquent directement au cadre ``de Segal'' car il travaille de fait
dans ce cadre, et n'op\`ere la troncation qu'\`a la fin de l'argument).
On peut donc
supposer qu'il existe un espace
topologique $Y$ avec $A= \Pi _{n,Se}(Y)$. On note que ce dernier
est automatiquement fibrant. On a une \'equivalence faible d'espaces $|\Re
(A)|\rightarrow Y$ (encore gr\^ace \`a \cite{Tamsamani}).
Par ailleurs, l'ensemble simplicial
$\Xi (A)$ est une variante du ``complexe
singulier'' de $Y$, correspondant
au foncteur $\Delta \rightarrow Top$ qui \`a $p\in
\Delta$ associe $\Re (R(p))$. Il s'ensuit que $| \Xi (A)| \rightarrow Y$ est
une \'equivalence faible.  On obtient que le morphisme d'espaces
topologiques
$$
| \Xi (A)|\rightarrow |\Re (A)|
$$
est une \'equivalence faible, ce qui implique que
$$
\Xi (A)\rightarrow \Re (A) = A^{int, 0}
$$
est une \'equivalence d'ensembles simpliciaux.
\eop

\begin{theoreme}
\label{intereqloc}
Il y a une \'equivalence naturelle de $1$-cat\'egories de Segal
$$
nSeCAT ^{int,1} \cong L(nSePC ).
$$
De plus,
si $\Xx$ est un site, il y a une \'equivalence naturelle (et fonctorielle en
$\Xx$)
$$
nSeCHAMP (\Xx )^{int, 1} \cong L(nSePCh(\Xx )).
$$
Ici les $L(\, )$ sont les localis\'ees de Dwyer-Kan des cat\'egories de
mod\`eles par rapport aux \'equivalences faibles.
\end{theoreme}
{\em Preuve:}
D'abord $nSeCAT ^{int,1}$ s'obtient \`a partir de $nSeCAT$ en prenant
l'int\'erieur $(\, )^{int, 0}$ des $n$-cat\'egories de Segal $nSeCAT _{1/}(A,B)$
qui sont fibrantes.
D'apr\`es   le lemme \ref{sigma}, on obtient l'int\'erieur d'un objet
en appliquant la construction $\Xi$. On a donc une \'equivalence
$$
nSeCAT ^{int,1}\cong \Xi \circ nSeCAT,
$$
et le lemme \ref{splinterne} donne l'\'equivalence avec
$L(nSePC)$.

Le foncteur $\underline{R}: \Delta \rightarrow nSePCh$ qui \`a $p$
associe le pr\'efaisceau constant de valeur $R(p)\in nSePC$,
d\'efinit une structure
simpliciale compatible sur la cat\'egorie de mod\`eles ferm\'ee interne
$nSePCh$. Si on note $\Xi _{\underline{R}}$ la construction $\Xi$
associ\'ee \`a cette structure (et $\Xi _R$ celle associ\'ee \`a la
structure $R$ sur $nSePC$), on a, pour tout $n$-pr\'echamp de Segal fibrant $A$,
$$
\Xi _{\underline{R}}(A) = \Xi _R (\Gamma (\Xx , A)).
$$
Rappelons qu'on a pos\'e
$$
nSe\underline{CHAMP} := INT (nSePCh_{cf}),
$$
et
$$
nSeCHAMP _{1/} (A, B) := \Gamma (\Xx , nSe\underline{CHAMP}_{1/}(A,B)).
$$
On obtient l'\'egalit\'e
$$
(\Xi _{\underline{R}}\circ nSe\underline{CHAMP})_{1/}(A,B)
= (\Xi _R \circ nSeCHAMP )_{1/}(A,B),
$$
autrement dit
$$
\Xi _{\underline{R}}\circ INT (nSePCh_{cf} )=
\Xi _{\underline{R}}\circ nSe\underline{CHAMP}=
\Xi _R \circ nSeCHAMP.
$$
Le lemme \ref{splinterne} donne alors l'\'equivalence
$$
L(nSePCh )\cong \Xi _R \circ nSeCHAMP.
$$
Comme pour la premi\`ere partie du th\'eor\`eme, en
utilisant que les $nSeCHAMP
_{1/}(A,B)$ sont fibrants et en leur appliquant le lemme \ref{sigma}, on
obtient l'\'equivalence
$$
\Xi _R \circ nSeCHAMP\cong  nSeCHAMP ^{int, 1}.
$$
Ces deux \'equivalences donnent la deuxi\`eme partie du th\'eor\`eme.
\eop

\subnumero{Le $\underline{Hom}$ interne relatif}

On pr\'esente ici une extension souvent utile de la notion de
$\underline{Hom}$
interne: il s'agit du $\underline{Hom}$ ``fibre par fibre'' pour un morphisme
$f:A\rightarrow B$, qu'on notera $\underline{Hom}(A/B, C)\rightarrow B$. Il
n'existera que sous certaines conditions sur le morphisme $f$.

Soit $M$ une cat\'egorie de pr\'efaisceaux (i.e. on suppose qu'il
existe une autre cat\'egorie $\Yy$ telle que $M$ soit \'equivalente \`a la
cat\'egorie des pr\'efaisceaux d'ensembles sur $\Yy$). Soient
$A,B,C\in M$ avec
une fl\`eche $f:A\rightarrow B$. Posons, pour $E\in M$,
$$
H(E):= \{ E\rightarrow B, \, \, E\times _BA\rightarrow C\} .
$$
Nous affirmons que ce foncteur est repr\'esent\'e par un objet
$\underline{Hom}(A/B,
C)$ de $M$. En effet, $H(-)$ transforme les colimites en limites,
ce qui suffit
pour garantir sa repr\'esentabilit\'e.

Supposons de plus que $M$ est une cmf interne. On dira qu'une fibration
$f:A\rightarrow B$ est {\em compatible aux changements de base (ccb)}
si, pour
tout diagramme
$$
B'\stackrel{a}{\rightarrow} B'' \rightarrow B
$$
o\`u $a$ est une \'equivalence faible, on a que
$$
A\times _BB'\rightarrow A\times _BB''
$$
est aussi une \'equivalence faible. On note que si $M$ est propre
(e.g. pour $M=0SePC$
ou $M=0SePCh$) alors toute fibration est compatible aux changements de
base---c'est la d\'efinition de ``propre''.

\begin{lemme}
\label{ccb1}
Si
$$
\begin{array}{ccc}
A&\stackrel{f}{\rightarrow}& B\\
\downarrow && \downarrow \\
A'&\stackrel{f'}{\rightarrow}& B'
\end{array}
$$
est un diagramme avec $f$ et $f'$ des fibrations et les fl\`eches
verticales des
\'equivalences, alors:
\newline
---si $f'$ est ccb alors $f$ estr ccb;
\newline
---si $B\rightarrow B'$ est une fibration et $f$ est ccb alors $f'$
est ccb; et
\newline
---si $B$ et $B'$ sont fibrants et $f$ est ccb alors $f'$ est ccb.
\end{lemme}
{\em Preuve:}
Laiss\'ee au lecteur (dans le troisi\`eme cas on peut choisir
deux fibrations $B\leftarrow B'' \rightarrow B'$ et
appliquer le deuxi\`eme cas).
\eop

La proposition suivante \'etablit la relation entre la condition ``ccb''
et le $\underline{Hom}$ interne relatif: pour que ce dernier ait un sens
homotopique il faut que le morphisme $f:A\rightarrow B$ soit ccb.

\begin{proposition}
\label{hominternerelatif}
Soit $M$ une cmf interne qui est une cat\'egorie de pr\'efaisceaux.
On suppose
en outre que $cof \subset M$ est la cat\'egorie des injections des
pr\'efaisceaux. Soit
$f:A\rightarrow B$ une fibration compatible aux changements de base. Alors le
foncteur
$$
C\mapsto \underline{Hom}(A/B, C)
$$
transforme fibrations en
fibrations, en particulier si $C$ est fibrant alors le morphisme structurel
$\underline{Hom}(A/B, C)\rightarrow B$ est une fibration. Cette construction
est---\`a \'equivalence pr\`es---invariante par \'equivalence en les
variables $A$ et $C$ (ces invariances sont soumises \`a la condition que la
fl\`eche $f$ reste une fibration ccb et que $C$ reste fibrant).
Si
$$
\begin{array}{ccc}
A' &\rightarrow & A \\
\downarrow && \downarrow \\
B'&\rightarrow &B
\end{array}
$$
est un diagramme cart\'esien
alors on a la formule
$$
\underline{Hom}(A'/B', C) = \underline{Hom}(A/B,C)\times _BB'.
$$
Si les fl\`eches verticales sont des fibrations ccb et le morphisme
$B\rightarrow
B'$ est soit une fibration triviale soit une \'equivalence entre deux
objets fibrants, alors le morphisme induit
$$
\underline{Hom}(A'/B', C) \rightarrow  \underline{Hom}(A/B,C)
$$
est une \'equivalence.
\end{proposition}
{\em Preuve:} Soit $g:C'\rightarrow C$ une fibration; montrons que
la fl\`eche
$\underline{Hom}(A/B,g)$ satisfait la propri\'et\'e de rel\`evement pour une
cofibration triviale $E\rightarrow E'$. Cette propri\'et\'e revient \`a la
propri\'et\'e de rel\`evement pour $g$ par rapport au morphisme
$$
j:A\times _BE\rightarrow A\times _BE' .
$$
Le fait que $cof$ est compatible aux produits fibr\'es implique que $j$
est une cofibration; et la propri\'et\'e ccb pour $f$ implique que $j$
est une
\'equivalence faible; d'o\`u la propri\'et\'e de rel\`evement voulue.

On note de la m\^eme fa\c{c}on que si $g:C'\rightarrow C$ est une fibration
triviale alors $g$ induit une fibration triviale $\underline{Hom}(A/B,g)$.
Ceci
implique l'invariance de la construction par
\'equivalence en l'argument fibrant
$C$ (en effet toute \'equivalence entre dans un
triangle---commutatif \`a une r\'etraction pr\`es---o\`u les deux autres
cot\'es
sont des fibrations triviales).

Soit $A\stackrel{\cong}{\hookrightarrow}A'\rightarrow B$ un diagramme
o\`u les
deux morphismes $f:A\rightarrow B$ et $f':A'\rightarrow B$ sont des
fibrations
ccb. Dans ce cas, pour $C$ fibrant, le morphisme
$$
\underline{Hom}(A'/B,C)\rightarrow \underline{Hom}(A/B, C)
$$
est une fibration triviale: pour une cofibration $E\rightarrow E'$ il
s'agit de voir que
$$
A\times _B E'
\cup ^{A\times _BE} A'\times _B E \rightarrow A'\times _BE'
$$
est une cofibration triviale, ce qui est le cas gr\^ace au fait que
$$
A\times _BE\rightarrow A'\times _BE \;\;\; et \;\;\;
A\times _BE' \rightarrow A'\times _BE'
$$
sont des cofibrations triviales (pour cela on n'a pas besoin de la
propri\'et\'e
ccb). Il s'ensuit que la formation de $\underline{Hom}(A/B, C)$ est
invariante
\`a
\'equivalence pr\`es par \'equivalence en la variable $A$ (\`a condition que
$f$ reste fibrant et ccb et que $C$ soit fibrant).

Enfin, la formule en cas de changement de base $B\rightarrow B'$
est une cons\'equence imm\'ediate de la d\'efinition du $\underline{Hom}$
interne relatif
comme foncteur repr\'esentable. Si le morphisme $B'\rightarrow B$
est une fibration triviale,
le produit fibr\'e avec celui-ci induit encore une
fibration triviale. Si cette fl\`eche est une \'equivalence entre objets
fibrants on peut se r\'eduire (en introduisant un $B"$ comme au d\'ebut de la
d\'emonstration) au cas d'une fibration triviale, ce qui donne l'\'enonc\'e
d'invariance par changement de $B$.
\eop

{\em Question:} Si $f:A\rightarrow B$ est ccb, en est-il de m\^eme pour
$\underline{Hom}(A/B,C)\rightarrow B$?

\begin{lemme}
\label{propunivHominterne}
Soit $M$ une cmf comme dans la proposition ci-dessus. Soit $f:A\rightarrow B$
une fibration compatible aux changements de base et soit $C$ un objet fibrant.
Pour tout morphisme $B'\rightarrow B$, on a
$$
Sect (B', \underline{Hom}(A/B, C)\times _BB') =
\underline{Hom}(A\times _BB',C).
$$
Ici $Sect(U,V)$ d\'esigne,
pour un
morphisme donn\'e $V\rightarrow U$,
la fibre de
$\underline{Hom}(U,V)$ au-dessus de $1_U\in \underline{Hom}(U,U)$.
\end{lemme}
{\em Preuve:}
Utiliser la propri\'et\'e d\'efinissant $\underline{Hom}(A/B, C)$.
Un morphisme
$$
E\rightarrow Sect (B', \underline{Hom}(A/B, C)\times _BB')
$$
est un morphisme $E\times B'\rightarrow \underline{Hom}(A/B, C)$
dont le compos\'e vers $B$ se factorise \`a travers $B'$. Ou encore un
morphisme
$$
E\times (B'\times _BA)= (E\times B')\times _B A\rightarrow C,
$$
i.e. un morphisme $E\rightarrow \underline{Hom}(A\times _BB',C)$.
\eop

Si la r\'eponse \`a la question pr\'ec\'edente est affirmative,
ou bien si $B$
et $B'$ sont fibrants, ou bien si toute fibration de base $B$
est ccb (cf le lemme ci-dessous par exemple), alors la propri\'et\'e
universelle
de ce lemme est aussi une propri\'et\'e universelle homotopique car
les morphismes rentrent en jeu dans les produits
fibr\'es sont fibrants,
les produits fibr\'es en question sont invariants sous \'equivalence,
et les sections sont prises le long d'un morphisme fibrant.

On applique maintenant cette discussion \`a nos exemples $M=nSePC$ et
$M=nSePCh$. En fait il suffira de traiter l'exemple des $n$-pr\'echamps de
Segal $M=nSePCh$ car on r\'ecup\`ere le cas des $n$-pr\'ecats
de Segal ($M=nSePC$)
comme le cas ``ponctuel'' i.e. le cas o\`u le site sous-jacent est
$\Xx =\ast$.

On note d'abord que si $n=0$ (i.e. dans le cas des pr\'efaisceaux
simpliciaux)
alors $M$ est propre et toute fibration $f$ est ccb. Pour $n>0$ on
obtient que
les objets qui sont essentiellement des $0$-pr\'echamps de Segal, i.e. les
objets $0$-groupiques ou $n$-groupo\"{\i}des de Segal, se comportent
comme des
pr\'efaisceaux simpliciaux. On obtient en outre une petite am\'elioration,
\`a savoir qu'il
suffit que la base $B$ soit $0$-groupique.

\begin{lemme}
\label{cestccb}
Dans la cmf $M=nSePCh$ des $n$-pr\'echamps de Segal \ref{SeCmf}, une
fibration
$f:A\rightarrow B$ de base $B$ telle que chaque $B(X)$ soit une
$n$-cat\'egorie de
Segal $0$-groupique (i.e. un $n$-groupo\"{\i}de de Segal), est compatible aux
changements de base.
\end{lemme}
{\em Preuve:}
On peut supposer que $B$ est fibrant (\`a cause de l'invariance par
changement
de $A$).  Au vu de l'\'enonc\'e partiel de propr\'et\'e \ref{proper},
il suffit
de voir la compatibilit\'e aux changements de base de la forme
$$
B'\rightarrow SeCat(B')\rightarrow B.
$$
Pour cela il suffit de voir la compatibilit\'e sous les m\^emes
changements de
base  au-dessus de chaque objet $X\in \Xx$. On est donc ramen\'e au cas
ponctuel
$\Xx =\ast$. Rappelons que dans ce cas l'op\'eration $SeCat$ peut \^etre vue
comme compos\'ee transfinie de coproduits avec des cofibrations standards
not\'ees $\Sigma \rightarrow h(M)$ dans \cite{nCAT} (ici $h(M)$ est la
$n$-pr\'ecat de Segal repr\'esent\'e par $M\in \Theta^{n+1}$). On peut donc
r\'eduire
au cas d'un changement de base de la forme
$$
\Sigma \rightarrow h(M) \rightarrow B.
$$
Soit $\overline{h(M)}$ le ``compl\'et\'e $0$-groupique'' de $h(M)$,
obtenu par
exemple en prenant le $\Pi _{n,Se}$ de la r\'ealisation topologique
de $h(M)$.
Le fait que $B$ soit $0$-groupique et fibrant implique l'existence d'une
factorisation de la forme
$$
\Sigma \rightarrow h(M) \rightarrow \overline{h(M)} \rightarrow B.
$$
On peut donc se r\'eduire au cas $B= \overline{h(M)}$. Maintenant le
morphisme
$B\rightarrow \ast$ est une fibration triviale donc, si
$$
\begin{array}{ccc}
A&\rightarrow &A'\\
\downarrow &&\downarrow \\
B&\rightarrow &\ast
\end{array}
$$
est une factorisation de $A\rightarrow \ast$ avec $A'$ fibrant,
la propri\'et\'e
ccb pour $A\rightarrow B$ est \'equivalente \'a cette propri\'et\'e pour
$A'\rightarrow \ast$. Or la propri\'et\'e ccb pour $A'\rightarrow \ast$
r\'esulte
de la stabilit\'e des \'equivalences faibles par produit
direct, qui n'est autre que IM(d) pour $M=nSePC$ (on rappelle qu'on a fait
r\'ef\'erence, pour sa d\'emonstration, \`a la m\^eme propri\'et\'e pour les
$n$-pr\'ecats non de Segal, qui \`a son tour \'etait l'une des \'etapes
principales de \cite{nCAT}).
\eop

En cons\'equence de ce lemme, pour toute fibration $f:A\rightarrow B$
de $n$-pr\'echamps de Segal avec $B$ $0$-groupique, et pour
tout
$n$-pr\'echamp de Segal fibrant $C$, le morphisme
$$
\underline{Hom}(A/B,C)\rightarrow B
$$
est une fibration poss\'edant
la propri\'et\'e \ref{propunivHominterne} qui est une propri\'et\'e
universelle homotopique.

Le contre-exemple \`a la propret\'e donn\'e dans \cite{nCAT} est le
morphisme
$$
a_{02}: A=I\rightarrow I^{(2)}=B
$$
o\`u $B$ est la cat\'egorie avec trois objets $0,1,2$ et des fl\^eches
$0\rightarrow 1\rightarrow 2$, et $a_{02}$ est l'inclusion de la fl\^eche
$0\rightarrow 2$. Ce morphisme est une fibration mais n'est pas compatible aux
changements de base. Cet exemple montre qu'il n'existe pas toujours un
$\underline{Hom}$ interne relatif de la forme $\underline{Hom}(A/B,C)$
satisfaisant \`a la propri\'et\'e universelle homotopique voulue, est
obstru\'ee.

\numero{La famille universelle}

\label{universellepage}

Soit pour le moment $\Xx$ une cat\'egorie munie de la topologie grossi\`ere.
On a un morphisme de $n+1$-pr\'echamps
$$
\Psi : nSeCHAMP(\Xx ) \times \Xx ^o \rightarrow nSeCAT.
$$
Ce morphisme peut \^etre vu comme la {\em famille universelle} de
$n$-cat\'egories
au-dessus de $nSeCHAMP(\Xx ) \times \Xx$. On d\'efinit $\Psi$ en posant
$$
\Psi (A,X):= A(X);
$$
et pour $f:X\rightarrow Y$ dans $\Xx $ on utilise le morphisme naturel
$$
\Gamma \underline{Hom}(A,B) \times \{ f\} \rightarrow \underline{Hom}(A(Y),
B(X))
$$
pour d\'efinir l'action de $\Psi$ sur les morphismes. A gauche,
il s'agit du
$\underline{Hom}$ interne des $n$-pr\'echamps de Segal; et, \`a droite,
il s'agit du
$\underline{Hom}$ interne des $n$-pr\'ecats de Segal.

Soit $nSeCAT'$ le remplacement fibrant de $nSeCAT$. Par composition on
obtient
un morphisme
$$
nSeCHAMP(\Xx ) \times \Xx ^o\rightarrow nSeCAT',
$$
qui correspond, par d\'efinition du $\underline{Hom}$ interne (des
$n+1$-pr\'ecats de Segal), \`a un morphisme
$$
\Phi : nSeCHAMP(\Xx ) \rightarrow \underline{Hom}(\Xx ^o, nSeCAT').
$$

Le r\'esultat suivant fait le
lien entre ces d\'efinitions et les notions de champs propos\'ees dans
\cite{limits}.

\begin{theoreme}
\label{correlation}
Soit $nSeCAT'$ un remplacement fibrant de $nSeCAT$.
Si la topologie de $\Xx$ est
grossi\`ere alors le morphisme
$\Phi$ induit une \'equivalence entre
$nSeCHAMP(\Xx )$ et la
$n+1$-cat\'egorie de Segal $\underline{Hom}(\Xx ^o, nSeCAT')$.

Si $\Xx$ est muni d'une autre topologie
$\Gg$ alors $\Phi$ induit une \'equivalence entre \newline
$nSeCHAMP(\Xx )$ et la sous-cat\'egorie
pleine de
$\underline{Hom}(\Xx ^o, nSeCAT')$ form\'ee des morphismes
$F:\Xx ^o\rightarrow nSeCAT'$ qui correspondent \`a des $n$-champs de Segal.

Pour le cas des $n$-champs non de Segal cette sous-cat\'egorie peut \^etre
identifi\'ee comme celle des morphismes qui 
satisfont \`a la condition de descente
de \cite{limits} 6.3 \`a savoir que la fl\`eche
$$
\lim _{\leftarrow} F|_{\Xx /X} \rightarrow
\lim _{\leftarrow} F|_{\Bb}
$$
est une \'equivalence pour tout crible $\Bb \subset \Xx /X$ de $\Gg$.
\end{theoreme}

Dans cette section (Proposition \ref{pleinfidele} ci-dessous) nous allons
seulement
d\'emontrer que $\Phi$ est pleinement fid\`ele. La surjectivit\'e essentielle de
$\Phi$ dans le premier cas du th\'eor\`eme
sera le th\'eor\`eme \ref{strictif3}
plus bas. Le deuxi\`eme cas r\'esultera alors
(une fois admis le th\'eor\`eme \ref{strictif3}) du Corollaire
\ref{critereaveclim}.

On note d'abord que si $\Gg$ d\'esigne la
topologie de $\Xx$,
alors les $n$-pr\'echamps $\Gg$-fibrants sont en particulier fibrants pour la
topologie grossi\`ere, tandis que la d\'efinition des morphismes dans
$nSeCHAMP(\Xx ^{\Gg})$
ne fait pas intervenir $\Gg$. Donc $nSeCHAMP(\Xx ^{\Gg} )$ est une
sous-cat\'egorie
pleine de $nSeCHAMP(\Xx ^{\rm gro})$ o\`u $\Xx ^{\rm gro}$ est la cat\'egorie
$\Xx$ avec
la topologie grossi\`ere. Ainsi la propri\'et\'e que $\Phi$ soit pleinement
fid\`ele  se r\'eduit au cas de la topologie grossi\`ere
et nous pouvons ignorer la topologie $\Gg$
pour le reste de cette
section.

Pour insister sur le fait qu'on ne consid\`ere pas la topologie, on change la
notation pour $\Xx$ en $\Yy$: on consid\`ere donc une cat\'egorie $\Yy$ avec sa
topologie grossi\`ere.

Rappelons que les objets de $nSeCHAMP(\Yy )$ sont les $n$-pr\'echamps de Segal
fibrants sur $\Xx$, et que si $A$ et $B$ sont deux tels objets, on a,
$$
nSeCHAMP (\Yy )_{1/}(A,B) := \Gamma \underline{Hom}(A,B).
$$
Cette $n$-pr\'ecat de Segal fibrante repr\'esente le foncteur qui \`a une
$n$-pr\'ecat de Segal $E$ associe l'ensemble $Hom (A\times \underline{E}, B)$.
D'autre part si $A$ est fibrant alors les valeurs $A(X)$ sont fibrantes, i.e.
objets de
$nSeCAT$. Le $n$-pr\'echamp de Segal fibrant $A$ induit donc un morphisme
$\Yy ^o \rightarrow  nSeCAT$ qu'on peut composer avec le
remplacement fibrant
pour obtenir un morphisme $\Yy ^o \rightarrow  nSeCAT'$.
Si $E$ est une $n$-pr\'ecat de Segal $E$, se donner un
morphisme
$$
E\rightarrow \underline{Hom}(\Yy ^o, nSeCAT')_{1/}(A,B)
$$
revient \`a se donner un morphisme
$$
F:\Yy ^o\times \Upsilon (E) \rightarrow nSeCAT'
$$
avec
$$
F|_{\Yy ^o \times \{ 0\}} = A,\;\;\;\;
F|_{\Yy ^o \times \{ 1\} } = B.
$$
(Voir \cite{limits} pour la notation $\Upsilon$.)
Or, la
donn\'ee d'un tel
morphisme strict (i.e. \`a valeurs dans $nSeCAT$ au lieu de $nSeCAT'$)
$$
F: \Yy ^o\times \Upsilon (E) \rightarrow nSeCAT
$$
avec
$$
F|_{\Yy ^o \times \{ 0\}} = A,\;\;\;\;
F|_{\Yy ^o \times \{ 1\}} = B
$$
revient \`a la donn\'ee, pour tout $y\in \Yy$, d'un morphisme
$E\times A(y)\rightarrow B(y)$ telle que pour $y\rightarrow z$ dans $\Yy$, ces
morphismes soient strictement compatibles avec les restrictions dans $A$ et $B$.
Autrement dit, on obtient l'\'egalit\'e
$$
\underline{Hom}(\Yy ^o, nSeCAT)_{1/}(A,B) = \Gamma \underline{Hom}(A,B).
$$
Par suite l'inclusion $nSeCAT \rightarrow nSeCAT'$ induit un morphisme
$$
\Gamma \underline{Hom}(A,B) \rightarrow
\underline{Hom}(\Yy ^o, nSeCAT')_{1/}(A,B),
$$
et la pleine fidelit\'e du foncteur
$\Phi$ est \'equivalente \`a la proposition suivante.

\begin{proposition}
\label{pleinfidele}
Soient $A$ et $B$ deux $n$-pr\'echamps de Segal fibrants sur $\Yy$. Alors le
morphisme
$$
\Gamma \underline{Hom}(A,B) \rightarrow
\underline{Hom}(\Yy ^o, nSeCAT')_{1/}(A,B).
$$
est une \'equivalence de $n$-cat\'egories de Segal.
\end{proposition}
{\em Preuve:} On introduit d'abord quelques notations.
Soit $Y$ une $n+1$-cat\'egorie de Segal  et $E$ une $n$-pr\'ecat de Segal. On
d\'efinit la $n+1$-pr\'ecat de Segal $Y^{\otimes E}$  qui a les m\^emes objets
que $Y$ en posant (pour $p\geq 1$)
$$
Y^{\otimes E}_{p/}(y_0,\ldots , y_p):=
Y_{p/}(y_0,\ldots , y_p) \times E.
$$
On s'int\'eresse \`a cette construction surtout quand $Y$ est
notre $1$-cat\'egorie de base $\Yy$, ou bien encore pour $Y=\Yy \times I$.

Si $\Zz \subset \Yy$ est une inclusion de $1$-cat\'egories ayant les m\^emes
objets, si $C'$ est une $n+1$-cat\'egorie fibrante, et si $A : \Zz \rightarrow
C'$ est un
morphisme, on introduit la $n$-cat\'egorie de Segal fibrante $Diag (\Yy , \Zz ,
A; C')$, dont les objets sont les prolongements de $A$ \`a $\Yy$, et qui
est caract\'eris\'ee par la formule
$$
Hom _{nSePC}(E,Diag (\Yy , \Zz , A;C'))=\{ f: \Yy ^{\otimes E} \rightarrow C',
\;\;\; f|_{\Zz ^{\otimes E}} = A|_{\Zz ^{\otimes E}} \} .
$$
En effet, on v\'erifie que le membre de droite d\'efinit un foncteur
en $E$ qui
transforme colimites en limites. Ce foncteur est donc
repr\'esentable par une $n$-pr\'ecat de Segal qu'on note
$Diag (\Yy , \Zz ,
A;C')$.
Comme $C'$ est
fibrante, notre foncteur a la propri\'et\'e d'extension pour les
cofibrations triviales $E\subset E'$, et donc $Diag (\Yy , \Zz , A;C')$ est
fibrante.

Soient $\nu (\Zz )\subset \nu (\Yy )$ et $\nu (C')$ les nerfs de ces $1$- ou
$n+1$-cat\'egories de Segal, consid\'er\'ees comme $n$-pr\'echamps de Segal
au-dessus de $\Delta$. On note que $\nu (C' )$ est fibrant au-dessus de
$\Delta$. On peut v\'erifier
que la $n$-cat\'egorie de Segal $Diag (\Yy , \Zz , A;C')$  est la fibre du
morphisme
$$
\underline{Hom}(\nu (\Yy ), \nu (C' ))
\rightarrow
\underline{Hom}(\nu (\Zz ), \nu (C' ))
$$
au-dessus de $A$. Il s'agit ici des $\underline{Hom}$ internes de
$n$-pr\'echamps de
Segal au-dessus de $\Delta$.

Soient maintenant $A,B: \Yy \rightarrow C'$ et notons $A\sqcup B$ le morphisme
correspondant
$$
A\sqcup B:\Yy \times \{ 0,1\} = \Yy \sqcup \Yy \rightarrow C'.
$$
On a la formule
$$
\underline{Hom}(\Yy , C')_{1/}(A,B) = Diag (\Yy \times I, \Yy \times \{ 0,1\} ,
A\sqcup B ; C').
$$
Pour prouver la pleine fidelit\'e, nous allons appliquer cette formule
avec
$C'=nSeCAT'$.

On va d'abord modifier l\'eg\`erement les choses.
Dans la situation du paragraphe
pr\'ec\'edent, soit $\Uu$ la  cat\'egorie dont les objets sont les paires de
suites $(y_0,\ldots , y_p; z_0, \ldots , z_q)$ munies de morphismes
$y_i\rightarrow y_{i+1}$,
$z_i\rightarrow z_{i+1}$ et aussi $y_p\rightarrow z_0$; et o\`u les morphismes
sont les inclusions et les d\'eg\'enerescences de suites.
Autrement dit, $\Uu$ est la {\em cat\'egorie des simplexes de $\Yy \times I$}
qui ne sont ni dans $\Yy \times \{ 0\}$ ni dans $\Yy \times
\{ 1\}$. Remarquons qu'on dispose d'un foncteur d'oubli
$$
(y_0,\ldots , y_p; z_0, \ldots , z_q)\mapsto
[\{ 0,\ldots , p\} , \{ 0,\ldots ,
q\} ]
$$
de
$\Uu$  vers $\Delta \times \Delta$.
On peut d\'efinir
un $n$-pr\'echamp $F = F(A,B; C')$ sur
$\Uu$ en prenant pour
$F(y_0,\ldots , y_p,z_0,\ldots ,z_q)$
la fibre au-dessus de
$(f,g)$ du morphisme
$$
C'_{p+q+1/}(A(y_0),\ldots , A(y_p),B(z_0),\ldots , B(z_q))
\rightarrow
$$
$$
C'_{p/}(A(y_0),\ldots , A(y_p))\times
C'_{q/}(B(z_0),\ldots , B(z_q))
$$
o\`u $f$ (resp. $g$) est l'image de l' \'el\'ement $(y_0,\ldots , y_p)$
(resp. $(z_0,\ldots , z_q)$) du nerf de $\Yy$ par le morphisme $A$ (resp. $B$).

Observons ici la propri\'et\'e suivante:

$(\ast ):$ pour tout $(y_0,\ldots , z_q)$,
le morphisme de restriction
$$
F(y_0,\ldots , y_p,z_0,\ldots ,z_q)\rightarrow F(y_p,z_0)
= C'_{1/}(A(y_p), B(z_0))
$$
correspondant \`a la fl\`eche
$$
(y_p; z_0)\rightarrow (y_0,\ldots ,y_p; z_0,ldots , z_q)
$$
de $\Uu$, est
une \'equivalence de $n$-cat\'egories de Segal.

Pour le voir, il faut utiliser
l'hypoth\`ese que $C'$ est fibrant, ce qui fait que le produit fibr\'e dans
la d\'efinition de $F$ est aussi un produit fibr\'e homotopique.

Nous allons finir la d\'emonstration en cours en admettant
que $F$ est fibrant au-dessus de $\Uu$, et nous consacrerons la fin de
la pr\'esente section \`a ce probl\`eme.
Observons maintenant qu'on a
$$
Diag (\Yy \times I, \Yy \times \{ 0,1\} ,
A\sqcup B ; C') = \Gamma (\Uu , F).
$$
En effet, si $E$ est une
$n$-pr\'ecat de Segal, la donn\'ee d'un morphisme
$E\rightarrow \Gamma (\Uu ,F)$
est \'equivalente \`a la donn\'ee, pour tout objet
$(y_0,\ldots , y_p,z_0,\ldots ,z_q)$ de $\Uu$, d'un morphisme
$$
E\rightarrow
C'_{p+q+1/}(A(y_0),\ldots , A(y_p),B(z_0),\ldots , B(z_q))
$$
compatible aux restrictions et
d\'eg\'en\'erescences de suites,
et compatible aux morphismes de transition de $A$ et $B$
sur les suites $(y_0,\ldots , y_p)$ et $(z_0,\ldots , z_q)$. Une
telle donn\'ee est encore \'equivalente \`a celle d'un
morphisme $(\Yy \times I)^{\otimes E}\rightarrow C'$
compatible avec $A\sqcup B$, d'o\`u la formule.

Soit $Arr(\Yy )$ la cat\'egorie dont les objets sont les fl\`eches
$x\rightarrow y$ de $\Yy$ et dont les morphismes de
$x\rightarrow y$ vers $x'\rightarrow y'$ sont les diagrammes commutatifs
$$
\begin{array}{ccc}
x & \rightarrow & x' \\
\downarrow && \downarrow \\
y& \leftarrow & y' .
\end{array}
$$
On a un foncteur $\varphi : \Uu \rightarrow Arr (\Yy )$ d\'efini par
$$
\varphi : (y_0,\ldots , y_p,z_0,\ldots , z_q)\mapsto (y_p\rightarrow z_0).
$$

On fixe dor\'enavant $C'=nSeCAT'$, et on supprime la r\'ef\'erence \`a
$C'$ dans
$Diag$. Par ailleurs on pose $F:=F(A,B; nSeCAT')$,
o\`u $A$ et $B$ sont d\'esormais
des foncteurs stricts $\Yy \rightarrow nSeCAT$,
i.e. des $n$-pr\'echamps de Segal sur $\Yy ^o$ avec
$B$
fibrant en tant que $n$-pr\'echamp sur $\Yy^o$.

Soit $\overline{hom}(A,B)$ le $n$-pr\'echamp sur $Arr (\Yy )$ qui \`a
$y\rightarrow z$ associe la $n$-cat\'egorie \linebreak
$\underline{Hom}(A(y),B(z))$ (compte-tenu de notre d\'efinition des
morphismes de $Arr(\Yy )$ c'
est bien un foncteur contravariant sur $Arr(\Yy )$). On a choisi
la notation $\overline{hom}$ pour sugg\'erer
qu'il s'agit l\`a d'une sorte de  ``$Hom$ externe''.

Si on avait d\'efini $F$ avec $nSeCAT$ au lieu de $nSeCAT'$
on aurait eu
$$
F(A,B, nSeCAT) = \varphi ^{\ast} (\overline{hom}(A,B)).
$$
Au lieu de cela, avec
$F:= F(A,B; nSeCAT')$
et l'inclusion
$nSeCAT\rightarrow nSeCAT'$ on a un morphisme
$$
\varphi ^{\ast} (\overline{hom}(A,B))\rightarrow F,
$$
qui est, objet-par-objet au dessus de $\Uu$, une \'equivalence de
$n$-cat\'egories de Segal (gr\^ace \`a l'\'equivalence $(\ast )$
vue plus haut). Par adjonction on obtient un morphisme
$$
\overline{hom}(A,B)\rightarrow \varphi _{\ast}(F).
$$
D'autre part on a les formules
$$
\underline{Hom}(\Yy , nSeCAT')_{1/}(A,B)=\Gamma (\Uu ,F) = \Gamma (Arr(\Yy ),
\varphi _{\ast} (F))
$$
et
$$
\Gamma \underline{Hom}(A,B)= \Gamma (   Arr(\Yy ), \overline{hom}(A,B)),
$$
qui induisent un diagramme commutatif
$$
\begin{array}{ccc}
\Gamma \underline{Hom}(A,B)&= &\Gamma ( Arr(\Yy ), \overline{hom}(A,B))\\
\downarrow && \downarrow \\
\underline{Hom}(\Yy , nSeCAT')_{1/}(A,B)&=& \Gamma (Arr(\Yy ),
\varphi _{\ast} (F)).
\end{array}
$$
Pour obtenir la proposition
i.e. que la fl\`eche verticale de gauche est une \'equivalence, il suffit de
montrer que la fl\`eche verticale de droite en est une. Pour ceci, il
suffirait  (compte-tenu du fait que $\varphi _{\ast} (F)$ est fibrant
sur $Arr(\Yy )$
cf \S 4)
de montrer les deux choses suivantes:
\newline
(i)\,\, le morphisme $\overline{hom}(A,B)\rightarrow \varphi _{\ast}(F)$ est
une \'equivalence objet-par-objet sur $Arr(\Yy )$; et
\newline
(ii)\,\, $\overline{hom}(A,B)$ est fibrant sur $Arr(\Yy )$.

En fait nous ne savons pas si (b) est vrai mais on va en montrer une
variante plus
faible:
\newline
(ii)'\,\, pour tout $A$ il existe une \'equivalence $A'\rightarrow A$ tel que
$\overline{hom}(A',B)$ soit fibrant sur $Arr(\Yy )$.

Les conditions (i)$+$(ii)' suffisent pour prouver le r\'esultat
voulu. En effet, elles entrainent
que dans
le diagramme $$
\begin{array}{ccc}
\Gamma \underline{Hom}(A,B)&\rightarrow &\underline{Hom}(\Yy ,
nSeCAT')_{1/}(A,B)\\
\downarrow & & \downarrow \\
\Gamma \underline{Hom}(A',B)&\rightarrow &\underline{Hom}(\Yy ,
nSeCAT')_{1/}(A',B)
\end{array}
$$
la fl\`eche du bas est une \'equivalence; or on sait que les deux fl\`eches
verticales sont des \'equivalences et par cons\'equent
il en est de m\^eme de la fl\`eche du
haut.

Pour prouver (i) fixons $\alpha : y\rightarrow z$ dans $Arr(\Yy )$
et notons
$\varphi /\alpha $ la cat\'egorie des $(y_0,\ldots , z_q)$ dans $\Uu $
munis d'un
morphisme
de $(y_p\rightarrow z_0)$ vers $\alpha$ dans $Arr(\Yy )$.
Par d\'efinition on a
$$
\varphi _{\ast}(F)(\alpha ) = \lim _{\leftarrow , (y,z)\in \varphi /\alpha }
F(y,z).
$$
Or on va voir que la cat\'egorie $\varphi /\alpha $ s'obtient par le m\^eme
proc\'ed\'e que
$\Uu $ mais avec d'autres arguments pour $Diag$. Soit $\Zz$ la cat\'egorie qui
contient la r\'eunion disjointe de $\Yy /y$ et $z/\Yy$ et qui, en plus des
fl\`eches de ces
cat\'egories, a une fl\`eche et une seule partant de chaque objet de $\Yy
/y$ vers
chaque objet de $z/\Yy$. On a  un morphisme de $1$-pr\'ecats de Segal
$$
\Yy /y \sqcup ^{\{ y\} } I \sqcup ^{\{ z \} } z/\Yy \rightarrow \Zz
$$
et on peut v\'erifier que ce morphisme est une \'equivalence faible.
En particulier on a la cofibration
$$
\Yy /y \sqcup z/\Yy \hookrightarrow \Zz ,
$$
et $A|_{\Yy /y}\sqcup B|_{z/\Zz}$ fournit un morphisme
$$
\Yy /y \sqcup z/\Yy\rightarrow nSeCAT'.
$$
Maintenant, par un argument
analogue \`a celui donn\'e plus haut pour $\Uu$,
$Diag (\Zz ,\Yy /y \sqcup z/\Yy, A|_{\Yy /y}\sqcup B|_{z/\Zz})$
est \'egale \`a
$\Gamma (\varphi /\alpha , F|_{\varphi /\alpha }).$
Par d\'efinition des limites, on a
$$
\lim _{(y,z)\in \varphi /\alpha }
F(y,z) = \Gamma (\varphi /\alpha , F|_{\varphi /\alpha }).
$$
D'o\`u l'\'egalit\'e
$$
\varphi _{\ast}(F)(\alpha )=
Diag (\Zz ,\Yy /y \sqcup z/\Yy, A|_{\Yy /y}\sqcup B|_{z/\Zz}).
$$

Par ailleurs, l'\'equivalence faible
$$
\Yy /y \sqcup ^{\{ y\} } I \sqcup ^{\{ z \} } z/\Yy \rightarrow \Zz
$$
nous donne
$$
Diag (\Zz ,\Yy /y \sqcup z/\Yy, A|_{\Yy /y}\sqcup B|_{z/\Zz})
\stackrel{\cong}{\rightarrow}
Diag(I, \{ 0,1\} , A(y)\sqcup B(z)).
$$
Ce dernier terme s'identifie \`a $\overline{hom}(A,B)(\alpha )$ (de
\facon
compatible avec les morphismes de source $\overline{hom}(A,B)(\alpha
)$).
Ceci prouve que le morphisme
$$
\overline{hom}(A,B)(\alpha )\rightarrow
\varphi _{\ast}(F)(\alpha )
$$
est une \'equivalence de $n$-cat\'egories de Segal, c'est-\`a-dire (i).

Pour prouver (ii)', on va montrer l'\'enonc\'e
\newline
(ii)'' le foncteur
$$
A\mapsto \overline{hom}(A,B)
$$
transforme les cofibrations de la structure de type
HBKQ, en fibrations.
\newline
Avec ce r\'esultat, pour $A$ quelconque,
on pourra choisir un remplacement HBKQ-cofibrant
$A'\rightarrow A$, et $\overline{hom}(A',B)$ sera fibrant ce
qui donne (ii)'. Pour montrer (ii)'' il suffit de consid\'erer une cellule
\'el\'ementaire, i.e. une cofibration de la forme
$$
h(X) \times N \rightarrow h(X) \times N'
$$
o\`u $N\rightarrow N'$ est une cofibration de $n$-pr\'ecats de Segal, et
o\`u $h(X)$ est le foncteur repr\'esent\'e par $X$ \`a savoir
$h(X)(Y)=\{ X\rightarrow Y\}$. On doit montrer que
$$
\overline{hom}( h(X) \times N',B)\rightarrow
\overline{hom}( h(X) \times N,B)
$$
est une fibration de $n$-pr\'echamps de Segal au-dessus de $Arr(\Yy )$.
Autrement dit, on doit montrer que si $U\rightarrow U'$ est une cofibration
triviale de $n$-pr\'echamps de Segal au-dessus de $Arr(\Yy )$,
et si on a des morphismes
$$
U\rightarrow \overline{hom}( h(X) \times N',B)
$$
et
$$
U'\rightarrow \overline{hom}( h(X) \times N,B)
$$
qui donnent le m\^eme morphisme
$U\rightarrow \overline{hom}( h(X) \times N,B)$, alors ces morphismes
se factorisent \`a travers
un m\^eme rel\`evement: $U'\rightarrow \overline{hom}( h(X) \times N',B)$.

De fa\c{c}on g\'en\'erale, un morphisme
$U\rightarrow \overline{hom}(A,B)$ est la donn\'ee,
pour tout objet $x\rightarrow y$ de $Arr (\Yy )$, d'un morphisme
$$
U(x\rightarrow y)\times A(x) \rightarrow B(y),
$$
telle que, pour tout morphisme
$(x\rightarrow y)\rightarrow (x'\rightarrow y')$ de $Arr (\Yy
)$ (cf ci-dessus), les deux morphismes induits
$$
U(x'\rightarrow y')\times A(x) \rightarrow B(y),
$$
soient \'egaux.

Notons $U(X\rightarrow ?)$ le $n$-pr\'echamp de Segal image inverse de $U$ par
le morphisme $(X/\Yy )^o \rightarrow Arr (\Yy )$; c'est un $n$-pr\'echamp de
Segal sur $(X/\Yy )^o$. Il ressort de la description du paragraphe
pr\'ec\'edent qu'un morphisme
$$
U\rightarrow \overline{hom}( h(X) \times N,B)
$$
n'est rien d'autre qu'un morphisme de $n$-pr\'echamps de Segal sur
$(X/\Yy )^o$
$$
U(X\rightarrow ?)\times N \rightarrow B|_{(X/\Yy )^o}.
$$
Notre probl\`eme devient donc: \'etant donn\'es deux morphismes
$$
U(X\rightarrow ?)\times N' \rightarrow B|_{(X/\Yy )^o}
$$
et
$$
U'(X\rightarrow ?)\times N \rightarrow B|_{(X/\Yy )^o}
$$
qui co\"{\i}ncident sur $U(X\rightarrow ?)\times N$,
admettent-ils un rel\`evement commun
\`a $U'(X\rightarrow ?)\times N'$?

Le fait que $B$ soit fibrant sur $\Yy ^o$  implique que
$B|_{(X/\Yy )^o}$ est fibrant (gr\^ace
\`a l'egalit\'e $(X/\Yy )^o=\Yy ^o /X$ et
au corollaire \ref{fibrestrict}). Etant donn\'ee une cofibration quelconque
$N\rightarrow N'$ et une cofibration triviale $U\rightarrow U'$, le
morphisme associ\'e $U(X\rightarrow ?)\rightarrow U'(X\rightarrow ?)$
est aussi une cofibration triviale et on
obtient une cofibration triviale
$$
U(X\rightarrow ?)\times N'\cup ^{U(X\rightarrow ?)\times N}
U'(X\rightarrow ?)\times N \rightarrow U'(X\rightarrow ?)\times N',
$$
ce qui (du fait que
$B|_{(X/\Yy )^o}$ est fibrant) donne la propri\'et\'e de rel\`evement voulue.
\eop

\subnumero{Analyse de $F(A,B;C')$ sur $\Uu$}
Dans la d\'emonstration ci-dessus on a laiss\'e de cot\'e le
probl\`eme de
montrer que $F(A,B;C')$ est fibrant sur $\Uu$, que nous
traitons maintenant. Pour cela nous utilisons la
notion de {\em
cat\'egories de Reedy}
pour laquelle nous ne donnerons (au \S 17) qu'une
pr\'esentation tr\`es
br\`eve: le lecteur devra donc se reporter aux
r\'ef\'erences (Reedy \cite{Reedy}, Bousfield-Kan \cite{BousfieldKan} p. 274,
Dwyer-Kan \cite{DwyerKan3},
Jardine-Goerss \cite{JardineGoerssBook}, Dwyer-Hirschhorn-Kan \cite{DHK},
Hirschhorn \cite{Hirschhorn} Chapter 16). On ne donne ci-dessous
qu'un argument rapide et
cette section s'adresse donc plut\^ot au lecteur
intr\'epide qui conna\^{\i}t d\'ej\`a ces r\'ef\'erences.

La cat\'egorie $\Uu$ est une cat\'egorie de Reedy pour le degr\'e
d\'efini par
$$
deg (y_0,\ldots , y_p; z_0,\ldots , z_q):= p+q;
$$
les morphismes ``directs''
sont
les applications ``faces'' qui correspondent \`a des inclusions de suites, et
les morphismes ``inverses'' sont les morphismes de d\'eg\'en\'er\'escence
qui correspondent \`a des surjections dans $\Delta \times \Delta$.
La cat\'egorie
oppos\'e $\Uu^o$ est une cat\'egorie de Reedy avec la m\^eme fonction de
degr\'e: l'oppos\'e d'un morphisme direct est inverse
et vice-versa.

En particulier, pour toute cmf $M$ la cat\'egorie $M^{\Uu
^o}$ des pr\'efaisceaux d'objets de $M$ au-dessus de $\Uu$, admet une
structure de
cmf qu'on appelle ``structure de Reedy'', voir \cite{Hirschhorn}.  La
cat\'egorie de Reedy $\Uu$ (c'est essentiellement l'exemple 16.1.10 de
\cite{Hirschhorn}) a la propri\'et\'e suivante, analogue d'une
propri\'et\'e de $\Delta$ (voir \cite{Hirschhorn} Corollary 16.4.7):
tout
pr\'efaisceau d'ensembles est cofibrant pour la structure de Reedy, et en fait
toute injection de pr\'efaisceaux d'ensembles sur $\Uu$ est une cofibration de
Reedy (on laisse la d\'emonstration au lecteur). Il s'ensuit la m\^eme
propri\'et\'e pour les $\Uu$-diagrammes dans $M$ si $M$ est une cmf de
pr\'efaisceaux d'ensembles dont les cofibrations sont les injections.
Tel est en particulier le cas pour $M=nSePC$.  \footnote{
On pourrait aussi s'int\'eresser \`a la cmf $nPC$ des $n$-pr\'ecats non de
Segal,
d\'efinie dans \cite{nCAT}.
Dans ce cas, les cofibrations sont les morphismes qui sont
injectifs sur les objets de $\Theta ^ n$ de longueur non-maximale; donc notre
remarque s'applique et tout morphisme de $\Uu$-diagrammes dans $nPC$ qui est
une cofibration objet-par-objet est une cofibration de Reedy.}

La cat\'egorie sous-jacente \`a la cmf de Reedy $nSePC^{\Uu}$ est
exactement
$nSePCh(\Uu )$; et les \'equivalences faibles (objet-par-objet) sont les
m\^emes. La propri\'et\'e \'enonc\'ee au paragraphe pr\'ec\'ecent
dit que les cofibrations sont les m\^emes; donc (\cite{Quillen}) les fibrations
sont les m\^emes et en fait {\em la structure de Reedy sur $nSePC^{\Uu}$
s'identifie \`a la structure de cmf de \ref{cmf} sur $nSePCh(\Uu )$}.

En particulier, pour v\'erifier que $F=F(A,B; C')$ est fibrant
dans $nSePCh(\Uu )$
il suffit de v\'erifier que $F$ est fibrant pour la structure de Reedy.
Rappelons ce que cela veut dire (\cite{Hirschhorn} Definition 16.3.2 (3)).
Pour un objet $(y,z)= (y_0,\ldots , y_p; z_0,\ldots , z_q)$ de $\Uu$,
on fabrique un (gigantesque) produit fibr\'e dont les termes principaux
(correspondant \`a des sous-suites de $(y,z)$ de longueur $p+q-1$)
sont les
$$
C'_{p+q-1/}(A(y_0),\ldots , \widehat{A(y_i)},\ldots , A(y_p),
B(z_0),\ldots , B(z_q))
$$
et
$$
C'_{p+q-1/}(A(y_0),\ldots , \ldots , A(y_p),
B(z_0),\ldots ,\widehat{B(z_j)},\ldots , B(z_q)),
$$
qu'on multiplie entre eux
au-dessus des produits fibr\'es analogues pour les
sous-suites de longueur $p+q-2$ (plus pr\'ecis\'ement, l'objet en question
s'exprime comme une limite sur la ``cat\'egorie appariante''
(``matching
category'') de $(y,z)\in \Uu $). Notons $Match((y,z); C')$
cette (gigantesque) $n$-pr\'ecat de Segal. On a un morphisme naturel
$$
C'_{p+q/}(A(y_0),\ldots ,  A(y_p),
B(z_0),\ldots , B(z_q))\rightarrow
Match ((y,z); C').
$$
Dire que l'objet qui nous int\'eresse, \`a savoir
$$
(y,z)\mapsto C'_{p+q/}(A(y_0),\ldots ,  A(y_p),
B(z_0),\ldots , B(z_q)),
$$
est fibrant pour la structure de Reedy,
revient exactement \`a dire que le
morphisme pr\'ec\'edent est une fibration.

Notons $O$ l'ensemble
des objets de $C'$ et $\Delta O$ la cat\'egorie des suites d'objets dans
$O$, i.e.
des $(x_0,\ldots , x_k)$ avec $x_i\in O$ (c'est la cat\'egorie des simplexes
de la cat\'egorie ayant $O$
pour ensemble d'objets et un isomorphisme entre chaque paire d'objets).
L'application $$
\nu _O(C'): (x_0,\ldots , x_k)\mapsto C'_{k/}(x_0,\ldots , x_k)
$$
d\'efinit un $n$-pr\'echamp de Segal au-dessus de $\Delta O$ et il est facile de
voir que si $C'$ est fibrante, alors $\nu _O(C')$ est fibrant.
D'autre part $\Delta O$ est aussi une cat\'egorie de Reedy et (comme pour $\Uu$)
la structure de Reedy co\"{\i}ncide avec la structure de \ref{cmf} sur
$nSePCh(\Delta O)$. Donc $\nu _O(C')$ est fibrant pour la structure de Reedy.
L'objet $Match ((y,z); C')$ est essentiellement le m\^eme
que l'objet appariant pour $\nu _O(C')$
(du moins si $p>0$ et
$q>0$---si $p=0$ ou $q=0$ il y a une face manquante tandis que si $p=q=0$ tout
est trivial), donc le fait que  $\nu
_O(C')$ soit fibrant implique que le morphisme
$$
C'_{p+q/}(A(y_0),\ldots ,  A(y_p),
B(z_0),\ldots , B(z_q))\rightarrow
Match ((y,z); C')
$$
est fibrant.

Rappelons que $F$ est d\'efini par le diagramme cart\'esien
suivant (o\`u les notations sont abr\'eg\'ees):
$$
\begin{array}{ccc}
F(y,z) & \rightarrow & C'_{p+q/}(Ay,
Bz)\\
\downarrow && \downarrow \\
\ast & \rightarrow & C'_{p/}(Ay)\times
C'_{q/}(Bz).
\end{array}
$$
Pour montrer que
$F$ est fibrant il nous faut voir que le morphisme de $n$-pr\'echamps de
Segal au-dessus de $\Uu$
$$
C'_{p+q/}(Ay, Bz)\rightarrow C'_{p/}(Ay)\times
C'_{q/}(Bz)
$$
est une fibration de Reedy; pour cela il convient
de remplacer l'objet appariant
$$
Match((y,z);C')
$$
ci-dessus, par {\em l'objet appariant relatif}
(``relative matching objet'') (\cite{Hirschhorn} Definition 16.2.22)
qu'on notera (\`a d\'efaut d'une meilleure notation)
$$
Match((y,z);
C'(Ay,Bz)\rightarrow C'(Ay)\times C'(Bz)).
$$
On renvoie le lecteur \`a
\cite{Hirschhorn} Definition 16.2.22 pour la d\'efinition que nous
n'\'ecrivons pas plus explicitement (ce serait trop long).
On remarque seulement que l'application appariante relative
$$
C'(Ay,Bz)\rightarrow Match((y,z);
C'(Ay,Bz)\rightarrow C'(Ay)\times C'(Bz))
$$
est identique \`a l'application
$$
C'(Ay,Bz)\rightarrow Match((y,z);C')
$$
pour $p>0$ et $q>0$; tandis que pour $p=0$ ou $q=0$, dans l'objet appariant
relatif on retrouve exactement la face manquante ci-dessus et l'application
$$
C'(Ay,Bz)\rightarrow Match((y,z);
C'(Ay,Bz)\rightarrow C'(Ay)\times C'(Bz))
$$
est \'egale \`a l'application appariante en $(Ay,Bz)$ pour $\nu _O(C')$ sur
$\Delta O$. Dans tous les cas on obtient que l'application appariante relative
pour $C'(Ay,Bz)\rightarrow C'(Ay)\times C'(Bz)$ est une fibration. Ceci prouve
que ce morphisme est une fibration pour
la structure de Reedy, donc une fibration pour la
structure \ref{cmf} sur $nSePCh(\Uu )$, et donc que sa fibre $F$ est
un objet fibrant.

\numero{Le champ associ\'e \`a un pr\'echamp}

\label{associepage}

La construction du faisceau associ\'e \`a un
pr\'efaisceau joue un r\^ole central dans la th\'eorie des faisceaux.
Le faisceau associ\'e est d\'efini par une propri\'et\'e universelle.
On commence par faire la m\^eme chose pour le champ associ\'e \`a un pr\'echamp.
Il faut se rappeler,
en lisant cette section, que l'analogue de la distinction
faisceau/pr\'efaisceau, est la distinction entre les $n$-champs de Segal pour la
topologie $\Gg$, et les $n$-champs de Segal pour la topologie grossi\`ere.
Ces derniers sont
simplement les pr\'efaisceaux de $n$-cat\'egories de Segal (la seule
des conditions du \S 9 qui entre en jeu pour la topologie
grossi\`ere est la condition que chaque fibre soit une $n$-cat\'egorie).

\begin{lemme}
\label{propuniv1}
Soit $A$ un $n$-pr\'echamp de Segal et $A\rightarrow A'$ un choix de champ
associ\'e. Si $B$ est un $n$-pr\'echamp de Segal $\Gg$-fibrant, alors le
morphisme
$$
\underline{Hom}(A', B)\rightarrow \underline{Hom}(A,B)
$$
est une \'equivalence.
\end{lemme}
{\em Preuve:} Cela r\'esulte
de la proposition \ref{interne} pour la cat\'egorie de
mod\`eles interne \newline
$nSePSh(\Xx ^{\Gg}$, puisque
$A\rightarrow A'$ est une $\Gg$-\'equivalence faible (par d\'efinition).
\eop

On rappelle que les objets de $nSeCHAMP (\Xx ^{\rm gro})$
sont les $n$-champs de Segal
fibrants pour la topologie grossi\`ere, et que
$$
nSeCHAMP (\Xx ^{\Gg })\subset nSeCHAMP (\Xx ^{\rm gro})
$$
est la sous-cat\'egorie enrichie pleine des objets qui sont $\Gg$-fibrants,
qui est
(par \ref{grofibchampimplfib}) aussi la sous-cat\'egorie enrichie pleine
des objets
qui sont des champs. Dans ce contexte,
si $A\in nSeCHAMP (\Xx ^{\rm gro})_0$ est un objet,
un {\em champ associ\'e \`a $A$} est un
remplacement $\Gg$-fibrant $A\rightarrow A'$ (on exige donc que $A'$ soit
$\Gg$-fibrant, et non pas seulement un champ, car on veut que $A'$ soit dans
la sous-cat\'egorie ci-dessus).
Autrement dit, un tel champ associ\'e
est un morphisme de $A$ vers $A'$ dans
$nSeCHAMP (\Xx ^{\rm gro})$, i.e. un objet de la $n$-cat\'egorie de
Segal $nSeCHAMP (\Xx ^{\rm gro})_{1/}(A,A')$.
Dans ce cadre on a la propri\'et\'e universelle suivante:

\begin{lemme}
\label{propuniv2}
Soit $u:A\rightarrow A'$ un champ associ\'e au sens du paragraphe
pr\'ec\'edent.
Pour tout $B\in nSeCHAMP (\Xx ^{\Gg})$, la composition avec
$u$ induit une \'equivalence
$$
nSeCHAMP (\Xx ^{\Gg})_{1/}(A',B)
\stackrel{\cong}{\rightarrow}
nSeCHAMP (\Xx ^{\rm gro})_{1/}(A,B).
$$
\end{lemme}
{\em Preuve:}
Au vu des d\'efinitions des
cat\'egories $nPC$-enrichies $nSeCHAMP(\ldots )$
et de l'\'egalit\'e
$$
nSeCHAMP (\Xx ^{\Gg})_{1/}(A',B)=
nSeCHAMP (\Xx ^{\rm gro})_{1/}(A',B),
$$
c'est une simple transcription du lemme \ref{propuniv1}.
\eop

Consid\'erons maintenant des
familles faibles de $n$-cat\'egories de Segal au-dessus de $\Xx ^o$,
qu'on voit comme morphismes
$\Xx ^o\rightarrow nSeCAT'$ ou comme objets de
$\underline{Hom}(\Xx ^o, nSeCAT')$.
Rappelons qu'on a l'\'equivalence
du th\'eor\`eme \ref{correlation}
$$
nSeCHAMP(\Xx ^{\rm gro}) \cong \underline{Hom}(\Xx ^o, nSeCAT').
$$
Si $F:\Xx ^o\rightarrow nSeCAT'$ on appellera {\em champ associ\'e \`a
$F$} tout couple $(F^{\rm ch}, f)$ o\`u
$F^{\rm ch}:\Xx ^o\rightarrow nSeCAT'$ est un champ i.e. satisfait la condition
du deuxi\`eme paragraphe de \ref{correlation}, et o\`u
$f: F \rightarrow F^{\rm ch}$ est, dans la $n+1$-cat\'egorie de Segal
$\underline{Hom}(\Xx ^o, nSeCAT')$, un
$1$-morphisme
\'equivalent, via l'\'equivalence de \ref{correlation} rappel\'ee
ci-dessus, \`a une $\Gg$-\'equivalence faible $A\rightarrow A'$.
Dans ces conditions $F$ est \'equivalent \`a $A$ et $F^{\rm ch}$ \'equivalent
\`a $A'$ (en tant qu'objets de $\underline{Hom}(\Xx ^o, nSeCAT')$).
En outre $A\rightarrow A'$ est un ``champ associ\'e'' au sens du
paragraphe pr\'ec\'edent. Ceci nous
permet, \`a partir de \ref{propuniv2}, de d\'emontrer la
propri\'et\'e universelle du lemme suivant.

Avant de l'\'enoncer, rappelons qu'un $1$-morphisme $f$ de
$\underline{Hom}(\Xx ^o, nSeCAT')$ est un morphisme
$$
I=\Upsilon (\ast )\rightarrow \underline{Hom}(\Xx ^o, nSeCAT'),
$$
autrement dit, un morphisme
$$
f:I\times \Xx ^o\rightarrow nSeCAT'.
$$
Sa source est la restriction $f|_{\{ 0\} \times \Xx ^o}$ et son but
est la restriction $f|_{\{ 1\} \times \Xx ^o}$. Si $f$ est un morphisme de
source $F$ et but $F^{\rm ch}$ et si $G:\Xx ^o\rightarrow nSeCAT'$ est une
autre famille, on peut construire
un morphisme ``composition avec $f$'' essentiellement
bien d\'efini
$$
-\circ f : \underline{Hom}(\Xx ^o, nSeCAT')_{1/}(F^{\rm ch}, G)
\rightarrow \underline{Hom}(\Xx ^o, nSeCAT')_{1/}(F, G).
$$
Pour en donner une d\'efinition pr\'ecise, on note
$$
\underline{Hom}(\Xx ^o, nSeCAT')_{2/}(F,F^{\rm ch}, G;f)
$$
la fibre de
$$
r_{01}:\underline{Hom}(\Xx ^o, nSeCAT')_{2/}(F,F^{\rm ch}, G)
\rightarrow \underline{Hom}(\Xx ^o, nSeCAT')_{1/}(F,F^{\rm ch})
$$
au-dessus de $f$; cette fibre vient avec deux morphismes
(le premier desquels est une \'equivalence)
$$
r_{12}:\underline{Hom}(\Xx ^o, nSeCAT')_{2/}(F,F^{\rm ch}, G;f)
\stackrel{\cong}{\rightarrow}
\underline{Hom}(\Xx ^o, nSeCAT')_{1/}(F^{\rm ch}, G)
$$
et
$$
r_{02}:\underline{Hom}(\Xx ^o, nSeCAT')_{2/}(F,F^{\rm ch}, G;f)
\rightarrow
\underline{Hom}(\Xx ^o, nSeCAT')_{1/}(F, G).
$$
On obtient le morphisme $-\circ f$ en choisissant un inverse de $r_{12}$
\`a \'equivalence pr\`es (ceci est possible car toutes les $n$-cat\'egories de
Segal consid\'er\'ees sont fibrantes) et en le composant avec $r_{02}$.

\begin{lemme}
\label{propuniv3}
Soit $F\in \underline{Hom}(\Xx ^o, nSeCAT')_0$ et soit $(F^{\rm ch}, f)$
un champ associ\'e au sens pr\'ec\'edemment d\'efini. Alors pour tout champ
$G$, i.e. objet de $\underline{Hom}(\Xx ^o, nSeCAT')_0$ v\'erifiant la
propri\'et\'e du deuxi\`eme paragraphe de \ref{correlation}, le morphisme
``composition avec $f$'' defini ci-dessus est une \'equivalence de
$n$-cat\'egories de Segal
$$
-\circ f : \underline{Hom}(\Xx ^o, nSeCAT')_{1/}(F^{\rm ch}, G)
\stackrel{\cong}{\rightarrow} \underline{Hom}(\Xx ^o, nSeCAT')_{1/}(F, G).
$$
\end{lemme}
{\em  Preuve:}
Par d\'efinition, $f$ est \'equivalent \`a un morphisme $u: A\rightarrow A'$
dans la $n+1$-cat\'egorie de Segal $nSeCHAMP(\Xx ^{\rm gro})$ qui est un champ
associ\'e pour la topologie $\Gg$, i.e. une $\Gg$-\'equivalence faible vers un
objet $\Gg$-fibrant. On peut aussi supposer par \ref{correlation} que $G$
provient d'un objet de $nSeCHAMP(\Xx ^{\Gg})$. Le morphisme ``composition avec
$f$'' est alors \'equivalent via l'\'equivalence de \ref{correlation}
$$
\underline{Hom}(\Xx ^o, nSeCAT')\cong
nSeCHAMP(\Xx ^{\rm gro}),
$$
au morphisme ``composition avec $u$'' (ce dernier \'etant bien d\'efini car
$nSeCHAMP(\Xx ^{\rm gro})$ est une cat\'egorie enrichie sur $nSeCat$). Le lemme
\ref{propuniv2} permet de conclure.
\eop

Si $(F^{\rm ch,1}, f^1)$ et $(F^{\rm ch,2}, f^2)$ sont deux champs associ\'es
\`a $F$, soit
$$
\underline{Hom}(\Xx ^o, nSeCAT')_{2/}(F,F^{\rm ch,1}, F^{\rm ch,2}; f^1, f^2)
:= (r_{01}, r_{02})^{-1}(f^1,f^2)
$$
la fibre de
$$
(r_{01}, r_{02}): \underline{Hom}(\Xx ^o, nSeCAT')_{2/}(F,F^{\rm ch,1}, F^{\rm
ch,2}) \rightarrow
$$
$$
\underline{Hom}(\Xx ^o, nSeCAT')_{1/}(F,F^{\rm ch,1})
\times
\underline{Hom}(\Xx ^o, nSeCAT')_{1/}(F, F^{\rm ch,2})
$$
au-dessus de $(f^1,f^2)$. On voit gr\^ace au lemme \ref{propuniv3} que
cette fibre $(r_{01}, r_{02})^{-1}(f^1,f^2)$ est contractile i.e. \'equivalente
(en tant que $n$-cat\'egorie de Segal) \`a $\ast$. Pour cela
il faut observer que le
morphisme $(r_{01}, r_{02})$ est une fibration de $n$-pr\'ecats de Segal.
On obtient donc un morphisme essentiellement bien d\'efini entre
$(F^{\rm ch,1}, f^1)$ et $(F^{\rm ch,2}, f^2)$.
On peut montrer la m\^eme chose dans l'autre sens
(et aussi la m\^eme chose  pour les diagrammes de dimension deux
entre $F^{\rm ch,1}$, $F^{\rm ch,2}$, $F^{\rm ch,1}$ ou
$F^{\rm ch,2}$, $F^{\rm ch,1}$, $F^{\rm ch,2}$) ce qui montre que le morphisme
qu'on a construit
est une \'equivalence. Ceci donne l'unicit\'e essentielle du ``champ associ\'e
\`a $F$''.

On veut maintenant recoller ces champs associ\'es en un foncteur
``faible''
$$
nSeCHAMP(\Xx ^{\rm gro})\rightarrow nSeCHAMP(\Xx ^{\Gg}),
$$
autrement  dit un foncteur vers le remplacement fibrant
$nSeCHAMP(\Xx ^{\Gg})'$. En fait, on va construire un diagramme
$$
\underline{Hom}(\Xx ^o, nSeCAT')\leftarrow C\rightarrow
\underline{Hom}(\Xx ^o, nSeCAT'),
$$
o\`u le premier morphisme est une \'equivalence, et le deuxi\`eme
prend ses valeurs dans la
sous-$n+1$-cat'egorie de Segal pleine des champs. Ceci donnera  (par
\ref{correlation} et les arguments
habituels) un vrai morphisme,
essentiellement bien d\'efini, vers le remplacement fibrant voulu.
On va utiliser librement les notations $Diag$ cf \S
12, et $\Upsilon ^k$ cf \cite{limits}.

On d\'efinit $C$ de la fa\c{c}on suivante. Les objets de $C$ sont les objets de
\newline
$nSeCHAMP(\Xx ^{\rm gro})$, i.e. les $n$-pr\'echamps de Segal fibrants
pour la topologie grossi\`ere. On choisit pour chaque objet $F\in C_0$ un
morphisme
$\eta _F: F\rightarrow F'$ vers un remplacement $\Gg$-fibrant.
Ensuite si $F_0,\ldots ,
F_p$ sont des objets, on d\'efinit (avec les notations du \S 12)
$$
C_{p/}(F_0,\ldots , F_p):= Diag(\Xx ^o\times I^{(k)}\times I,
\Xx ^o\times \{ 0,\ldots , p\}
\times I;  \eta _{F_0}\sqcup \ldots \sqcup \eta _{F_p}; nSeCAT').
$$
C'est la $n$-cat\'egorie de Segal des diagrammes dans $nSeCAT'$ qui ont la forme
d'un rectangle $p\times 1$, se restreignant \`a $\{ i\} \times I$,
pour tout $i$,
en le morphisme $\eta _{F_i}: F_i\rightarrow F'_i$ vu comme
morphisme $I\rightarrow nSeCAT'$.

En divisant---\`a cofibration triviale pr\`es---le
rectangle en coproduit de carr\'es $1\times 1$ et en utilisant le fait que
$nSeCAT'$
est fibrante, on obtient que l'application de Segal
$$
C_{p/}(F_0,\ldots , F_p)\rightarrow
C_{1/}(F_0,F_1)\times \ldots \times C_{1/}(F_{p-1}, F_p)
$$
est une \'equivalence; donc $C$ est une $n+1$-cat\'egorie de Segal.
Ensuite, en divisant un carr\'e $1\times 1$ en deux triangles,
on voit qu'un morphisme
$$
E\rightarrow Diag(\Xx ^o\times I\times I,
\Xx ^o\times \{ 0,1\}\times I; \eta _{F_0}\sqcup \eta _{F_1};
nSeCAT')
$$
n'est rien d'autre qu'une paire de morphismes
$$
u: \Xx ^o\times \Upsilon ^2(\ast ,E)\rightarrow nSeCAT', \;\;\;
v: \Xx ^o\times \Upsilon  ^2(E,\ast )\rightarrow nSeCAT',
$$
avec
$$
r_{01}(u)=\eta _{F_0},\;\; r_{12}(v)= \eta _{F_1}, \;\; r_{02}(u)=r_{02}(v).
$$
Si on note par exemple $\underline{Hom}(\Xx ^o,nSeCAT')_{2/}(F_0, F'_0, F'_1;
\eta _{F_0},-)$ la fibre de
$$
r_{01}:\underline{Hom}(\Xx ^o,nSeCAT')_{2/}(F_0, F'_0, F'_1)\rightarrow
\underline{Hom}(\Xx ^o,nSeCAT')_{1/}(F_0, F'_0)
$$
au-dessus de $\eta _{F_0}$, alors on obtient
$$
Diag(\Xx ^o\times I\times I,
\Xx ^o\times \{ 0,1\}\times I; \eta _{F_0}\sqcup \eta _{F_1};
nSeCAT') =
$$
$$
\underline{Hom}(\Xx ^o,nSeCAT')_{2/}(F_0, F'_0, F'_1; \eta _{F_0},-)\times _{
\underline{Hom}(\Xx ^o,nSeCAT')_{1/}(F_0,  F'_1)}
$$
$$
\underline{Hom}(\Xx ^o,nSeCAT')_{2/}(F_0, F_1, F'_1; -,\eta _{F_1}).
$$

La propri\'et\'e universelle pour le morphisme $\eta _{F_0}$, et le fait que
$F'_1$ soit $\Gg$-fibrant (i.e. un champ), impliquent que
le morphisme
$$
u\mapsto r_{02}(u)
$$
est une fibration triviale de $n$-pr\'ecats de Segal,
$$
r_{02}: \underline{Hom}(\Xx ^o,nSeCAT')_{2/}(F_0, F'_0, F'_1; \eta _{F_0},-)
\rightarrow \underline{Hom}(\Xx ^o,nSeCAT')_{1/}(F_0,  F'_1).
$$
D'autre part, le fait que
$$
\Upsilon (E)\sqcup ^{\{ 1\} }\Upsilon (\ast )
\rightarrow \Upsilon ^2(E,\ast )
$$
soit une cofibration triviale implique que le morphisme
$$
\underline{Hom}(\Xx ^o,nSeCAT')_{2/}(F_0, F_1, F'_1; -,\eta _{F_1})
\rightarrow
\underline{Hom}(\Xx ^o,nSeCAT')_{1/}(F_0, F_1)
$$
est une \'equivalence. On obtient que le morphisme
$$
C_{1/}(F_0,F_1)=
Diag(\Xx ^o\times I\times I, \Xx ^o \times \{ 0,1\}\times I; \eta _{F_0}\sqcup
\eta _{F_1}; nSeCAT')
$$
$$
\rightarrow \underline{Hom}(\Xx ^o,nSeCAT')_{1/}(F_0,F_1)
$$
est une \'equivalence.

La restriction des diagrammes de $I^{(p)}\times I$ sur
$I^{(p)}\times \{ 0\}$ fournit un morphisme
$$
C_{p/}(F_0,\ldots , F_p)\rightarrow\underline{Hom}(\Xx
^o,nSeCAT')_{p/}(F_0,\ldots , F_p),
$$
donc un morphisme de $n+1$-cat\'egories de Segal
$$
R_0: C\rightarrow \underline{Hom}(\Xx ^o,nSeCAT').
$$
On  a montr\'e que ce morphisme est pleinement fid\`ele;
et il est essentiellement surjectif par construction gr\^ace \`a
\ref{correlation}.

D'autre part,  la restriction des diagrammes de $I^{(p)}\times I$ sur
$I^{(p)}\times \{ 1\}$ fournit un morphisme
$$
C_{p/}(F_0,\ldots , F_p)\rightarrow\underline{Hom}(\Xx
^o,nSeCAT')_{p/}(F_0,\ldots , F_p),
$$
d'o\`u un morphisme
$$
R_1: C\rightarrow \underline{Hom}(\Xx ^o,nSeCAT').
$$
Ces deux morphismes forment le diagramme cherch\'e:
$$
\underline{Hom}(\Xx ^o,nSeCAT')
\stackrel{R_0}{\leftarrow}C
\stackrel{R_1}{\rightarrow}\underline{Hom}(\Xx ^o,nSeCAT').
$$

Observons qu'on a $R_0(F)=F$ et $R_1(F)=F'$.
On dira que ce diagramme est le ``foncteur champ associ\'e'',
encore not\'e ${\bf ch}$
et on \'ecrira aussi $F'= {\bf ch}(F)$.

Il reste maintenant \`a construire une transformation naturelle
entre l'identit\'e et le foncteur qu'on a construit ci-dessus. Ce
doit \^etre un diagramme
$$
nSeCHAMP(\Xx ^{\rm gro})\stackrel{T_0}{\leftarrow} C\times I
\stackrel{T_1}{\rightarrow}
nSeCHAMP(\Xx ^{\rm gro}).
$$
On prend $T_0:= R_0\circ pr_1$ et on veut construire un morphisme $T_1$
qui donne $R_0$ sur $C\times \{ 0\}$ et $R_1$ sur $C\times \{ 1\}$.
Pour cela,
rappelons que les objets de $C\times I$ sont les couples $(F,\epsilon
)$ o\`u $F$ est un objet de $C$ et $\epsilon$ est $0$ ou $1$. On
observe que
$$
(C\times I)_{p/}((F_0,\epsilon _0),\ldots , (F_p, \epsilon _p))
$$
est vide sauf si $\epsilon _0,\ldots , \epsilon _p)$ est une suite
non-d\'ecroissante, auquel cas, c'est tout simplement
$C_{p/}(F_0,\ldots , F_p)$. A toute suite non-d\'ecroissante $\epsilon
_0,\ldots , \epsilon _p$ on associe le morphisme
$$
i_{\epsilon} :I^{(p)}\rightarrow I^{(p)}\times I
$$
d\'efini par $ i_{\epsilon}(k):= (k, \epsilon _k)$.

L'image inverse par ce morphisme
fournit donc pour chaque suite $\epsilon$ un morphisme
$$
i_{\epsilon}^{\ast}: C_{p/}(F_0,\ldots , F_p)\rightarrow \underline{Hom}(\Xx
^o,nSeCAT')_{p/}(F_0^{(\epsilon _0)},\ldots , F_p^{(\epsilon _p)} ).
$$
Ici $F^{(0)}$ et $F^{(1)}$ d\'esignent
respectivement $F$ et $F'$.

Pour $\epsilon = (0,\ldots , 0)$,
on retrouve le morphisme de restriction qu'on a
utilis\'e pour d\'efinir $R_0$, tandis que
pour $\epsilon =(1,\ldots , 1)$ on retrouve celui qu'on a
utilis\'e pour d\'efinir $R_1$. Pour toute suite
d'objets $(F_0,\epsilon _0),\ldots , (F_p, \epsilon _p)$ on
d\'efinit \`a l'aide
de $i_{\epsilon}^{\ast}$ le morphisme
$$
(C\times I)_{p/}((F_0,\epsilon _0),\ldots , (F_p, \epsilon _p))\rightarrow
\underline{Hom}(\Xx
^o,nSeCAT')_{p/}(F_0^{(\epsilon _0)},\ldots , F_p^{(\epsilon _p)} );
$$
ce qui permet de d\'efinir le morphisme cherch\'e
$$
T_1: C\times I \rightarrow \underline{Hom}(\Xx
^o,nSeCAT').
$$

On
peut v\'erifier que $T_1|_{\{ F\} \times I}$ est bien
le morphisme $F\rightarrow F'$.
Le diagramme $(T_0,T_1)$ constitue donc bien une transformation naturelle
$$
1\rightarrow {\bf ch}
$$
entre un pr\'echamp et son champ associ\'e.

\begin{theoreme}
\label{cestadjoint}
Avec cette transformation naturelle, le foncteur
``champ associ\'e'' ${\bf ch}$ est homotopiquement adjoint \`a l'inclusion
$$
nSeCHAMP(\Xx ^{\Gg})\subset nSeCHAMP(\Xx ^{\rm gro}).
$$
\end{theoreme}
{\em Preuve:}
L'\'equivalence necessaire est une cons\'equence de la propri\'et\'e universelle
de ${\bf ch}(F)$, \ref{propuniv2} par exemple.
\eop

\subnumero{Calcul des morphismes dans le champ associ\'e}

On a la formule suivante:

\begin{lemme}
\label{formule}
Pour un $n$-prechamp de Segal $F$ dont les
valeurs sont des $n$-cat\'egories de Segal ($n\geq 1$), et pour $x,y\in F(X)$,
il y a une \'equivalence objet-par-objet au-dessus de $\Xx /X$,
$$
({\bf ch}(F))_{1/}(x,y) \cong {\bf ch}(F_{1/}(x,y)).
$$
\end{lemme}
{\em Preuve:}
Notons que
$({\bf ch}(F))_{1/}(x,y)$ est un $n-1$-champ de Segal par \ref{critere}.
Puisque $F\rightarrow {\bf ch}(F)$ est une $\Gg$-\'equivalence faible, par
d\'efinition (\S 3)
$$
F_{1/}(x,y)\rightarrow ({\bf ch}(F))_{1/}(x,y)
$$
est une $\Gg$-\'equivalence faible. Alors par d\'efinition (\S 9)
$({\bf ch}(F))_{1/}(x,y)$ est le $n-1$-champ de Segal associ\'e \`a
$F_{1/}(x,y)$.
\eop

\subnumero{Version \'el\'ementaire}
Il semble difficile d'obtenir le foncteur ${\bf ch}$ directement en utilisant la
structure de cat\'egorie de  mod\`eles interne, car le foncteur ``remplacement
$\Gg$-fibrant'' n'est pas strictement compatible aux produits directs.
On peut cependant remarquer que la th\'eorie de la localisation de Dwyer-Kan
donne imm\'ediatement
$$
L(nSePCh , W^{\rm gro}) \rightarrow L(nSePCh, W^{\Gg}),
$$
o\`u $W^{\rm gro}\subset nSePCh$ est la sous-cat\'egorie
des \'equivalences faibles  pour la topologie grossi\`ere, et
$W^{\Gg}\subset W^{\rm gro}$ est sa sous-cat\'egorie des $\Gg$-\'equivalences
faibles.
Le fait que $W^{\rm gro}$ est contenu dans $W^{\Gg}$ donne directement le
morphisme ci-dessus.
Maintenant le r\'esultat du th\'eor\`eme
\ref{intereqloc}
montre que le morphisme pr\'ec\'edent est \'equivalent \`a un morphisme
$$
nSeCHAMP(\Xx ^{\rm gro})^{int, 1}\rightarrow
nSeCHAMP'(\Xx ^{\Gg})^{int, 1}
$$
(on a pris un remplacement fibrant not\'e $nSeCHAMP'$ \`a l'arriv\'ee).
Ce  morphisme est la restriction de ${\bf ch}$ \`a l'int\'erieur $1$-groupique.
Dans le cas $n=0$ i.e. des pr\'efaiceaux simpliciaux, on a
$$
0SeCHAMP(\Xx ^{\rm gro})^{int, 1}=0SeCHAMP(\Xx ^{\rm gro})
$$
et
$$
0SeCHAMP(\Xx ^{\Gg})^{int, 1}=0SeCHAMP(\Xx ^{\Gg})
$$
et on obtient ainsi (i.e. avec la localis\'ee de Dwyer-Kan) une d\'efinition
\'el\'ementaire du foncteur ${\bf ch}$.

\numero{Limites}

\label{limitespage}

Dans cette section, on fait quelques rappels concernant
les limites. Dans \cite{limits}
on d\'efinit les notions de limite et colimite dans une $n$-cat\'egorie.
On peut donner les m\^emes
d\'efinitions {\em mutatis mutandis} dans une $n$-cat\'egorie de
Segal. On renvoie le lecteur \`a \cite{limits}, 3.1 et 3.2 pour ces
d\'efinitions. En transcrivant les d\'emonstrations de \cite{limits}
dans le cadre des $n$-cat\'egories de Segal, on obtient:

\begin{theoreme}
La $n+1$-cat\'egorie de Segal $nSeCAT'$ admet toutes les limites et colimites
index\'ees par des (petites) $n+1$-cat\'egories de Segal.
\end{theoreme}

Rappelons comment calculer concr\`etement une limite (c'est la construction
utilis\'ee dans \cite{limits} 4.1). Si $\Yy$ est une $n+1$-cat\'egorie
munie d'un morphisme $F: \Yy \rightarrow nSeCAT'$, alors
on a (\cite{limits} 4.1.2)
$$
\lim _{\leftarrow , \Yy} F
= \underline{Hom}(\Yy , nSeCAT')_{1/}(\ast _{\Yy}, F).
$$

D'autre part, on dispose des notions de
limite et colimite homotopique ($holim$ et $hocolim$) pour les
diagrammes (index\'es par une petite $1$-cat\'egorie) dans une cmf
essentiellement arbitraire. Ces notions
sont dues \`a Bousfield et Kan \cite{BousfieldKan},
voir aussi 
Vogt \cite{Vogt}, Edwards et Hastings \cite{EH},
Cordier et Porter \cite{CordierPorter2},
Hirschhorn \cite{Hirschhorn}, 
Dwyer-Hirschhorn-Kan \cite{DHK}.

Normalement nous utiliserons la notation $holim$ pour les
constructions de type Bous\-field-Kan, et la notation $\lim$ pour les
limites selon \cite{limits}. Heureusement, les deux notions co\"{\i}ncident
pour les limites de $n$-cat\'egories de Segal index\'ees par des
$1$-cat\'egories. C'est essentiellement une cons\'equence de la proposition
suivante. On laisse au lecteur le soin
de faire la comparaison analogue avec le $holim$ de
\cite{Hirschhorn} pour la cmf $nSePC$.

\begin{proposition}
\label{calclim}
Supposons que $\Yy$ est une cat\'egorie et $F$ un $n$-pr\'echamp de Segal
au-dessus de $\Yy$, qu'on peut aussi voir comme morphisme
$\Yy ^o\rightarrow nSeCAT$. Soit $F\rightarrow F'$ une \'equivalence faible
objet-par-objet vers un $n$-pr\'echamp de Segal fibrant pour la topologie
grossi\`ere. Alors on a des \'equivalences de $n$-cat\'egories de Segal
$$
\Gamma (\Yy , F')
\stackrel{\cong}{\rightarrow}
$$
$$
\lim _{\leftarrow , \Yy ^o} F' \cong
\lim _{\leftarrow , \Yy ^o} F .
$$
Ici les limites
sont celles d\'efinies suivant \cite{limits}.
\end{proposition}
{\em Preuve:}
On a
l'\'egalit\'e
$\Gamma (\Yy , F')=\Gamma \underline{Hom}(\ast _{\Yy}, F')$ car
les deux membres
repr\'esentent le m\^eme foncteur sur la cat\'egorie
des $n$-pr\'ecats de Segal.
En particulier, on a
$$
\Gamma (\Yy , F')= nSeCHAMP (\Yy )_{1/}(\ast, F').
$$
D'autre part, d'apr\`es (\cite{limits} 4.0.1 et 4.1.2), on a
$$
\lim _{\leftarrow , \Yy ^o} F' = \underline{Hom}(\Yy ^o, nSeCAT')_{1/} (\ast
_{\Yy}, F').
$$
Maintenant le fait que le foncteur $\Phi$ du Th\'eor\`eme
\ref{correlation} soit pleinement fid\`ele (cf Proposition \ref{pleinfidele})
signifie \'exactement que
$$
nSeCHAMP (\Yy )_{1/}(\ast, F')\rightarrow
\underline{Hom}(\Yy ^o, nSeCAT')_{1/} (\ast
_{\Yy}, F')
$$
est une \'equivalence. Enfin l'\'equivalence entre les limites de $F$ et $F'$
r\'esulte des propri\'et\'es d'invariance de la notion de limite
\cite{limits}.
\eop

On a une caract\'erisation similaire pour les colimites, qui utilise la
structure de type HBKQ (\S 5). On l'\'enonce ci-dessous
bien qu'on n'en ait pas besoin. Cette proposition
peut aussi \^etre vue comme un \'enonc\'e de compatibilit\'e entre la
d\'efinition de
$hocolim$ de Bousfield-Kan,
et la construction de \cite{limits}.

\begin{proposition}
\label{calccolim}
Soient $\Yy$ une cat\'egorie et $F$ un $n$-pr\'echamp de Segal
au-dessus de $\Yy$, vu aussi comme morphisme
$\Yy ^o\rightarrow nSeCAT$. Soit $F'\rightarrow F$ une \'equivalence faible
objet-par-objet avec $F'$ cofibrant pour la
structure de cmf de type HBKQ, pour la topologie
grossi\`ere. Alors on a des \'equivalences de $n$-cat\'egories de Segal
$$
colim_{nSePC}(F')
\stackrel{\cong}{\rightarrow}
$$
$$
\lim _{\rightarrow , \Yy ^o} F' \stackrel{\cong}{\rightarrow}
\lim _{\rightarrow , \Yy ^o} F
$$
o\`u $colim_{nSePC}(F')$ est la colimite standard du foncteur $\Yy \rightarrow
nSePC$ de $\Yy$ vers la cat\'egorie des $n$-pr\'ecats de Segal,
tandis que les autres
colimites
sont les colimites d\'efinies suivant \cite{limits}.
\end{proposition}
\eop

Gr\^ace \`a
la Proposition \ref{calclim} on obtient une version de nos crit\`eres
qui utilise les limites de \cite{limits}. Cette version de la d\'efinition des
champs a \'et\'e sugg\'er\'ee
dans \cite{limits}.

\begin{corollaire}
\label{critereaveclim}
Soit $A$ un $n$-pr\'echamp de Segal sur $\Xx$ dont les valeurs sont
des $n+1$-cat\'egories de Segal, qui correspond \`a un morphisme
$\Xx \rightarrow nSeCAT$.
Alors pour que $A$ soit un $n$-champ de Segal il faut et il suffit que:
\newline
(a)\,\, pour tout $X\in \Xx$ et tout $x,y\in A_0(X)$ le $n-1$-pr\'echamp
de Segal $A_{1/}(x,y)$ sur $\Xx /X$ soit un $n-1$-champ de Segal
(resp. $Path^{x,y}A$ soit un $0$-champ de Segal, dans le cas $n=0$); et
\newline
(b)\,\, pour tout $X\in \Xx$ et tout crible $\Bb \subset \Xx /X$ de la
topologie,
le morphisme
$$
A(X)\rightarrow \lim _{\leftarrow , \Bb ^o} A|_{\Bb}
$$
soit essentiellement surjectif.

Ou encore, pour qu'un $n$-pr\'echamp (non de Segal)
 $A$ soit un $n$-champ il faut et il suffit que:
\newline
(c)\,\, pour tout $X\in \Xx$ et tout crible $\Bb \subset \Xx /X$ de la
topologie,
le morphisme
$$
A(X)\rightarrow \lim _{\leftarrow , \Bb ^o} A|_{\Bb}
$$
soit une \'equivalence de $n$-cat\'egories de Segal.
\end{corollaire}
{\em Preuve:}
Au vu de la proposition \ref{calclim}, la condition (b) est \'equivalente \`a
celle de \ref{critere} ou \ref{0critere} (les conditions (a) sont
identiques). De
m\^eme la condition (c) est \'equivalente \`a celle de \ref{newcritere}
(qui ne marche que pour les $n$-pr\'echamps non de Segal d'apr\`es la correction
{\tt v3}).
\eop

{\em D\'emonstration de la deuxi\`eme partie du th\'eor\`eme
\ref{correlation}:} si on suppose connue l'es\-sen\-ti\-el\-le surjectivit\'e de
$\Phi$ pour la premi\`ere partie de \ref{correlation} (qui sera d\'emontr\'ee en
\ref{strictif3} ci-dessous), la deuxi\`eme partie de \ref{correlation}
d\'ecoule maintenant de la condition (c) du corollaire ci-dessus.

\bigskip

\numero{Un peu plus sur la condition de descente}

\label{peupluspage}

On voudrait avoir une version un peu plus explicite des limites qui entrent
dans la caract\'erisation \ref{critereaveclim}.

On suppose ici que le site $\Xx$ admet des produits fibr\'es.

On dira que le site $\Xx$ admet {\em suffisamment de sommes disjointes
compatibles aux produits fibr\'es}
s'il existe un ordinal $\beta$ tel que les cribles correspondant aux familles
couvrantes de taille strictement plus petite que
$\beta$ engendrent la topologie (i.e. sont cofinaux
parmi tous les cribles), et si les sommes disjointes \`a
moins de $\beta$ \'el\'ements
existent dans $\Xx$ et sont compatibles aux produits fibr\'es (i.e. qu'on a la
formule de distributivit\'e pour le produit fibr\'e de deux sommes disjointes
au-dessus d'un objet de $X$).

Par exemple si $\Xx$ est quasi-compact, on peut prendre
$\beta = \omega$ et la condition est qu'il existe des sommes disjointes finies.
Pour le site des espaces paracompacts on prendrait $\beta = \omega +1$.

Par ailleurs, on dira que {\em les sommes disjointes sont couvrantes} si
pour toute famille $\Uu$ dont la somme disjointe $X$ existe, la famille $\Uu$
est une famille couvrant $X$.

\begin{remarque}
\label{casuffit}
Dans les sites consid\'er\'es ici,
il suffit pour la condition (b) de \ref{critere} ou \ref{critereaveclim}
de prendre en compte les cribles engendr\'es par les familles couvrantes \`a
moins de $\beta$ \'el\'ements.
\end{remarque}
En effet, ces cribles engendrent la topologie.

Supposons maintenant que $\Uu = \{ U_{\alpha}\rightarrow X\}$ est une famille
couvrante $X\in \Xx$, ayant moins de $\beta$ \'el\'ements. On pose
$$
\amalg \Uu :=
\amalg _{\alpha} U_{\alpha}.
$$
C'est un objet de $\Xx /X$. Cependant, le crible engendr\'e par cet objet
est plus grand que le crible engendr\'e par $\Uu$ et en fait les cribles
engendr\'es par les $\amalg \Uu$ ne forment pas en g\'en\'eral
une famille cofinale.

On a donc besoin de la d\'efinition suivante.

\begin{definition}
\label{compsomdis}
Soit $A$ un $n$-pr\'echamp de Segal sur $\Xx$
dont les valeurs $A(U)$ sont fibrantes.
On dira que {\em $A$ est compatible aux sommes disjointes} si pour toute famille
$\Uu = \{ U_{\alpha} \}$ telle que $\amalg \Uu$ existe, le morphisme naturel
$$
A(\amalg \Uu )\rightarrow \prod _{\alpha} A(U_{\alpha})
$$
est une \'equivalence faible.
\end{definition}

Dans notre cadre o\`u l'ordinal $\beta$ est fix\'e on dira juste ``compatible
aux sommes disjointes'' sans pr\'eciser qu'il s'agit seulement des sommes
disjointes \`a moins de $\beta$ \'el\'ements.

On note en particulier que la somme disjointe de la famille vide (i.e. la
famille index\'ee  par $\emptyset$) est un objet initial $\iota$ du site
$\Xx$ (g\'en\'eralement l'objet vide $\iota = \emptyset$ si $\Xx$ est un site
d'objets g\'eom\'etriques), et la  condition que $A$ soit compatible aux sommes
disjointes entra\*{\i}ne $A(\iota )= \ast$.

On peut comprendre la condition \ref{compsomdis} \`a l'aide de la topologie
suivante: on appelle {\em topologie des composantes connexes} la topologie
engendr\'ee par les familles  $\Uu \rightarrow X$ v\'erifiant $X=\amalg \Uu$
(cette condition a un sens sans qu'on suppose l'existence de toutes les sommes
disjointes: elle signifie que $X$ est une somme disjointe des
\'el\'ements de $\Uu$). L'existence de produits fibr\'es et
la compatibilit\'e de
ceux-ci avec les sommes disjointes garantit qu'on d\'efinit bien
ainsi une topologie.

\begin{remarque}
\label{topcompcon}
Soit $A$ un $n$-pr\'echamp de Segal sur $\Xx$ dont les valeurs $A(U)$
sont fibrantes.
Alors $A$ est compatible aux sommes disjointes si et seulement si $A$ est un
champ pour la topologie des composantes connexes.
\end{remarque}

En effet, la condition \ref{compsomdis} est la m\^eme que la condition (c) de
\ref{critereaveclim}: la limite de $A$ sur le crible engendr\'e par $\Uu$
est \'equivalente \`a $\prod _{\alpha} A(U_{\alpha})$ et, d'apr\`es la remarque
\ref{casuffit}, il suffit de consid\'erer les cribles engendr\'es par les
familles couvrantes.

On va maintenant
montrer que si $A$ est compatible aux sommes disjointes, alors dans la
condition (b) de
\ref{critere} ou \ref{critereaveclim}, il suffit de consid\'erer les cribles
engendr\'es par une famille couvrante avec $1$ \'el\'ement.
Plus pr\'ecisement on montre:

\begin{lemme}
\label{unseul}
Supposons que $\Xx$ admet suffisamment de sommes disjointes.
Soit $A$ un $n$-pr\'echamp de Segal sur $\Xx$ fibrant objet-par-objet.
Supposons que $A$ est compatible aux sommes disjointes (\ref{compsomdis})
et satisfait la condition (a) de \ref{critereaveclim}, ainsi que la
condition (b) pour les cribles engendr\'es par les familles couvrantes \`a
$1$ \'el\'ement. Alors $A$ est un $n$-champ de Segal.
\end{lemme}
{\em Preuve:}
On peut supposer $A$ fibrant pour la topologie grossi\`ere. Par la remarque
\ref{topcompcon}, $A$ est aussi fibrant pour la topologie des composantes
connexes. Soit $\Bb$ un crible engendr\'e par une famille couvrante $\Uu
\rightarrow X$. Soit $\Bb '$ le crible engendr\'e par $\amalg \Uu$; on a $\Bb
\subset \Bb '$. Le morphisme
$$
\ast _{\Bb}\rightarrow \ast _{\Bb '}
$$
est une cofibration triviale pour la topologie des composantes connexes.
Donc si $f:\ast _{\Bb}\rightarrow A$ est un morphisme, il en existe une
extension en $f':\ast _{\Bb '}\rightarrow A$. On dispose de la condition (b) de
\ref{critereaveclim} pour $\Bb '$, qui
implique celle
de  \ref{critere}, et il s'ensuit que ce morphisme $f'$ s'\'etend en un
morphisme $\ast _{\Xx /X}\rightarrow A$. On obtient ainsi la condition (b) pour
tout crible $\Bb$ engendr\'e par une famille couvrante, ce qui suffit
(d'apr\`es la remarque \ref{casuffit}).
\eop

\bigskip

On suppose toujours que $\Xx$ admet des
produits directs et suffisamment de sommes
disjointes. On fixe dor\'enavant $X\in \Xx$ et une famille couvrante $\Uu$ qui
consiste en un \'el\'ement $U\rightarrow X$. On note $\Bb$ le crible
engendr\'e par $\Uu$. On d\'efinit un foncteur
$$
\rho (\Uu ): \Delta ^o \rightarrow \Xx
$$
par
$$
\rho (\Uu )(p):= U\times \ldots \times U \;\;\;\; p+1\;\;
\mbox{fois}.
$$

On va encore re\'ecrire la condition (b) de \ref{critereaveclim}. On suppose
que $A$ est un $n$-pr\'echamp de Segal, et on a $\rho (\Uu
)^{\ast }A$ qui est un $n$-pr\'echamp de Segal sur $\Delta ^0$, i.e. une
$n$-pr\'ecat
de Segal cosimpliciale ou encore un foncteur de $\Delta$ vers les
$n$-pr\'ecats de
Segal. D'autre part, soit
$$
t: \rho (\Uu )\rightarrow \underline{X}
$$
la transformation naturelle de foncteurs $\Delta ^o\rightarrow \Xx$ vers le
foncteur constant \`a valeurs $X$, qui provient des projections $\Uu \times
_X\ldots \times _X\Uu \rightarrow X$. Si on note $c^{\ast} A(X)$ la
$n$-pr\'ecat de Segal cosimpliciale constante (i.e. le foncteur
constant de $\Delta$ dans
les $n$-pr\'ecats de Segal et de valeur $A(X)$) alors $t$ induit
un morphisme
$$
t^{\ast} : c^{\ast}A(X) \rightarrow  \rho (\Uu )^{\ast }A
$$
de $n$-pr\'ecats de Segal cosimpliciales.

En particulier, $t$ induit un morphisme
$$
t^{\ast} : A(X) \rightarrow \lim _{\leftarrow , \Delta} \rho (\Uu )^{\ast }A.
$$
Ici la limite est celle suivant
\cite{limits}. D'apr\`es la construction de
\cite{limits} rappel\'ee au \S 14 ci-haut, on peut \'ecrire
$$
\lim _{\leftarrow , \Delta} \rho (\Uu )^{\ast }A
= \underline{Hom}(\Delta , nSeCAT')_{1/}(\ast , \rho (\Uu )^{\ast }A).
$$
Ici $nSeCAT'$ d\'esigne un remplacement fibrant de $nSeCAT$.

\begin{corollaire}
Soit $\rho (\Uu )^{\ast }A\rightarrow B$ une \'equivalence faible
objet-par-objet
vers un $n$-pr\'echamp de Segal au-dessus de $\Delta ^o$ fibrant pour la
topologie
grossi\`ere. Alors
\newline
$\lim _{\leftarrow , \Delta} \rho (\Uu )^{\ast }A$
est \'equivalente \`a la $n$-pr\'ecat de
Segal $\Gamma (\Delta ^o, B)$ des sections globales du pr\'efaisceau $B$.
\end{corollaire}
{\em Preuve:}
C'est une cons\'equence directe de la proposition \ref{calclim}.
\eop

Nous avons appris le lemme suivant dans le livre de Dwyer-Hirschhorn-Kan
\cite{DHK}, mais il est d\^u \`a Bousfield-Kan (\cite{BousfieldKan} Ch. XI \S
9).
Nous en donnons une version un peu plus faible que ce qu'on trouve dans
\cite{BousfieldKan} \cite{DHK}. On dira qu'un foncteur $f:\Yy \rightarrow \Zz$
est  {\em fortement initial}, si pour tout objet $Z\in \Zz$, la cat\'egorie
$f/Z$ est filtrante (\`a titre de comparaison, dans {\em loc cit.} $f$ est
dit {\em
initial} si les $f/Z$ sont \`a nerf contractile).

\begin{lemme}
Si $f:\Yy \rightarrow \Zz$ est un foncteur de $1$-cat\'egories qui est fortement
initial, alors pour tout foncteur $G: \Zz \rightarrow nSeCAT'$ le morphisme
naturel  $$ \lim _{\leftarrow , \Zz } G \stackrel{\cong}{\rightarrow}
\lim _{\leftarrow , \Yy } f^{\ast}G
$$
est une
\'equivalence de $n$-cat\'egories de Segal.
\end{lemme}
{\em Preuve:}
D'apr\`es \ref{strictif3}
\footnote{Il n'y a pas ici
de circularit\'e dans la mesure o\`u si on utilise ce lemme
avant la d\'emonstration de \ref{strictif3}, ce sera sur un $G$ qui est d\'ej\`a
strict.}
on peut supposer  que
$G$ provient d'un  $n$-pr\'echamp de Segal fibrant pour la topologie
grossi\`ere de $\Zz$.

Le remont\'e $f^{\ast}G$ est alors fibrant pour la topologie grossi\`ere de
$\Yy$.  En effet, il suffit de montrer que $f_!$ pr\'eserve les cofibrations
triviales. On a la formule
$$
f_!(U)(Z)= \lim _{\rightarrow , f/Z} U|_{f/Z},
$$
et comme par hypoth\`ese $f/Z$ est filtrante, cette colimite pr\'eserve les
cofibrations et les cofibrations triviales.

On prouve le lemme d'abord pour $n=0$ i.e. pour les $0$-cat\'egories de Segal ou
ensembles simpliciaux. Dans ce cas, d'apr\`es la Proposition \ref{calclim}, les
limites en question sont \'egales \`a $\Gamma (\Zz , G)$ et
$\Gamma (\Yy , f^{\ast}G)$ respectivement. Or celles-ci sont \'equivalentes aux
$holim_{\Zz}(G)$ et
$holim_{\Yy}(f^{\ast}G)$ de Bousfield-Kan \cite{BousfieldKan}.
La propri\'et\'e en question pour ces $holim$
(\cite{BousfieldKan} Ch. XI th\'eor\`eme 9.2) donne
l'\'enonc\'e.

Pour $n$ quelconque on peut dire exactement la m\^eme chose en utilisant
l'extension de Dwyer-Hirschhorn-Kan du $holim$ de Bousfield-Kan, \`a une
cat\'egorie de mod\`eles plus g\'en\'erale en l'occurrence la cat\'egorie des
$n$-pr\'ecats de Segal.

Une autre fa\c{c}on de compl\'eter l'argument pour $n$ quelconque
est de noter qu'on a
$$
\Gamma (\Zz , G)^{int, 0}=
\Gamma (\Zz , G^{int, 0})
$$
et de m\^eme
$$
\Gamma (\Yy , f^{\ast}G)^{int, 0}=
\Gamma (\Yy , (f^{\ast}G)^{int, 0});
$$
il s'ensuit que
$$
\Gamma (\Zz , G)^{int, 0}\rightarrow
\Gamma (\Yy , f^{\ast}G)^{int, 0}
$$
est une \'equivalence. En particulier le foncteur en question est
essentiellement surjectif.

D'autre part pour $x,y\in \Gamma (\Zz , G)$ on a
$$
\Gamma (\Zz , G)_{1/}(x,y) =
\Gamma (\Zz , G_{1/}(x,y))
$$
et $G_{1/}(x,y)$ est un $n-1$-pr\'echamp fibrant sur $\Zz$ (car
la propri\'et\'e de
rel\`evement pour $G_{1/}(x,y)$
vis-\`a-vis d'une cofibration triviale $E\rightarrow E'$,
est impliqu\'ee par la m\^eme propri\'et\'e
pour $G$ vis-\`a-vis de la cofibration triviale $\Upsilon
(E)\rightarrow \Upsilon (E')$); et
de m\'eme pour $(f^{\ast}G )_{1/}(f^{\ast}x, f^{\ast}y)$.  Par r\'ecurrence,
on peut
supposer connu le lemme pour les $n-1$-pr\'ecats de Segal,
d'o\`u on obtient que
$$
\Gamma (\Zz , G_{1/}(x,y))\rightarrow \Gamma (\Yy , (f^{\ast}G)_{1/}(f^{\ast}x,
f^{\ast}y))
$$
est une \'equivalence, ce qui entra\^{\i}ne que
$$
\Gamma (\Zz , G)_{1/}(x,y)\rightarrow
\Gamma (\Yy , f^{\ast}G)_{1/}(f^{\ast}x, f^{\ast}y)
$$
est une \'equivalence et le foncteur en question est pleinement fid\`ele.
\eop

{\em Remarque--Question:}
Bousfield-Kan et Dwyer-Hirschhorn-Kan donnent le m\^eme \'e\-non\-c\'e pour les
colimites
dans une cat\'egorie de mod\`eles ferm\'ee $M$, avec seulement l'hypoth\`ese que
$f$ est initial i.e. les $f/Z$ ont leur
nerf contractile. Peut-on obtenir le m\^eme
\'enonc\'e ici? Dans la preuve ci-dessus, nous avons utilis\'e la condition
``fortement initial'' pour dire que $f_!$ pr\'eserve les cofibrations et
les cofibrations triviales.

\begin{corollaire}
\label{limsurdelta}
Supposons que $\Xx$ admet des produits fibr\'es.
Soit $A$ un $n$-pr\'echamp de Segal, et soit $\Bb \subset
\Xx /X$
le crible associ\'e \`a une famille couvrante $\Uu \rightarrow X$.
Alors on a une \'equivalence naturelle
$$
\lim _{\leftarrow , \Bb ^o} A|_{\Bb} \stackrel{\cong}{\rightarrow}
\lim _{\leftarrow , \Delta  } \rho (\Uu )^{\ast}A.
$$
\end{corollaire}
{\em Preuve:}
D'apr\`es le lemme, il suffit de v\'erifier que le morphisme
$$
\rho (\Uu ): \Delta \rightarrow \Bb ^o
$$
est initial. Pour $Y\in \Bb ^o$, la cat\'egorie $\rho (\Uu )/Y$
est \'egale \`a la cat\'egorie
$$
\{ (n,f): \;\;\; n\in \Delta , \;\;\; f: Y\rightarrow P_n(\Uu )\} .
$$
Un morphisme $f: Y\rightarrow P_n(\Uu )$ est la m\^eme chose qu'un
$n+1$-uplet de morphismes $Y\rightarrow \Uu$ au-dessus de $X$. Soit $H:= Hom
(Y,\Uu )$ l'ensemble des morphismes $Y\rightarrow \Uu$ au-dessus de $X$,
et adoptons la notation $P_n(H)$ pour le produit
de $n+1$ copies de $H$ (ce qui d\'efinit un
ensemble simplicial $P_{\cdot}(H)$ qui est le nerf de la
cat\'egorie ayant $H$
pour ensemble d'objets et un isomorphisme entre chaque paire d'objets). La
cat\'egorie $\rho (\Uu )/Y$ est donc \'egale \`a
$$
\{ (n,h): \;\;\; n\in \Delta , \;\;\; h\in P_n(H)\} ,
$$
autrement dit c'est la cat\'egorie totale associ\'ee \`a
l'ensemble simplicial
$P_{\cdot}(H)$,
qui est contractile car $H$ est non-vide (la condition $H\neq \emptyset$
est \'equivalente \`a la condition $Y\in \Bb$, c'est la d\'efinition du crible
$\Bb$ associ\'e \`a $\Uu$). Par le lemme, on obtient l'\'equivalence entre
limites cherch\'ee.
\eop

\begin{proposition}
\label{critere7}
Supposons que $\Xx$ admet des produits fibr\'es, et suffisamment de
\linebreak
sommes disjointes compatibles aux produits fibr\'es. Soit $A$ un
$n$-pr\'echamp  de Segal
tel que $A(X)$ soit une $n$-cat\'egorie de Segal pour tout $X\in \Xx$.
Alors $A$ est un champ si et seulement si:
\newline
(a) \,  pour tout $X\in \Xx$ et toute paire d'objets $x,y\in A_0(X)$, le
$n-1$-pr\'echamp de Segal $A_{1/}(x,y)$ est un $n-1$-champ de Segal au-dessus de
$\Xx /X$;
et
\newline
(b)\, $A$ est compatible aux sommes disjointes (d\'efinition
\ref{compsomdis}), et pour tout $X$ et toute famille couvrante $\Uu \rightarrow
X$ \`a un \'el\'ement, le morphisme
$$
t^{\ast} : A(X) \rightarrow \lim _{\leftarrow , \Delta} \rho (\Uu )^{\ast }A
$$
d\'efini ci-dessus, est essentiellement surjectif.

La condition (b) peut \^etre
remplac\'ee par la m\^eme condition pour n'importe lequel des
$A^{int, k}$.

Si $A$ est un $n$-pr\'echamp (non de Segal) alors
le fait que  $A$ soit un champ est aussi \'equivalent \`a la condition
\newline
(c)\,\,  $A$ est compatible aux sommes disjointes (d\'efinition
\ref{compsomdis}), et pour tout $X$ et toute famille couvrante $\Uu
\rightarrow
X$ \`a un \'el\'ement, le
morphisme
$$
t^{\ast} : A(X) \rightarrow \lim _{\leftarrow , \Delta} \rho (\Uu )^{\ast }A
$$
est une \'equivalence de $n$-cat\'egories de Segal.
\end{proposition}
{\em Preuve:}
Combiner les Corollaires \ref{critereaveclim} et \ref{limsurdelta}.
\eop

Si on remplace la condition (b) par la m\^eme pour le pr\'efaisceau simplicial
$A^{int, 0}$, alors la limite
en question est la m\^eme que la $holim$ de Bousfield-Kan. On
peut donc re\'ecrire:
\begin{corollaire}
Supposons que $\Xx$ admet des produits fibr\'es, et suffisamment de sommes
disjointes compatibles aux produits fibr\'es. Soit $A$ un $n$-pr\'echamp de
Segal tel que $A(X)$ soit une $n$-cat\'egorie
de Segal pour tout $X\in \Xx$. Alors $A$ est un champ si et seulement si:
\newline
(a) \,  pour tout $X\in \Xx$ et toute paire d'objets $x,y\in A_0(X)$, le
$n-1$-pr\'echamp de Segal $A_{1/}(x,y)$ est un $n-1$-champ de Segal au-dessus de
$\Xx /X$; et
\newline
(b)\, $A$ est compatible aux sommes disjointes (d\'efinition
\ref{compsomdis}), et pour tout $X$ et toute famille couvrante $\Uu \rightarrow
X$ \`a un \'el\'ement, le morphisme
$$
t^{\ast} : A^{int, 0}(X) \rightarrow holim _{\leftarrow , \Delta} \rho
(\Uu )^{\ast }A^{int ,0}
$$
induit une surjection sur les $\pi _0$.

Si $A$ est un champ alors le morphisme de (b) est une \'equivalence faible.
\end{corollaire}
\eop

Ce crit\`ere permet---par r\'ecurrence sur $n$--- d'appr\'ehender la condition
pour
\^etre un champ,  en termes de pr\'efaisceaux simpliciaux uniquement, et
dans ce cadre la
notion est d\'ej\`a bien connue.

\numero{La construction de Grothendieck}

\label{gropage}

On indique ici comment on peut construire une ``$n$-cat\'egorie fibr\'ee''
au-dessus
d'une $1$-cat\'egorie $Y$, \`a partir d'un pr\'efaisceau de $n$-cat\'egories
au-dessus de $Y^o$; c'est l'analogue de la construction de SGA1 \cite{SGA1}.

Il serait int\'eressant d'avoir une d\'efinition de
{\em $n$-cat\'egorie fibr\'ee} et une construction allant dans
l'autre sens (pour obtenir l'analogue de toute la th\'eorie de SGA1);
mais nous ne traitons pas cette question ici.

Soit $Y$ une $1$-cat\'egorie et
$A$ un $n$-pr\'echamp de Segal au-dessus de $Y^o$.
(cette covariance
est mieux adapt\'ee pour
certains aspects de la discussion
qui suit).
\footnote{En fait il y a un choix \`a faire quant \`a la direction des
fl\`eches, comme la distinction de SGA 1 entre ``cat\'egorie fibr\'ee'' et
``cat\'egorie cofibr\'ee'' et nous ne consid\'erons qu'une des
deux possibilit\'es.
L'autre s'en d\'eduit par conjugaison avec l'op\'eration
``cat\'egorie oppos\'ee''.}
On consid\`ere $A$ comme un foncteur covariant $\Yy
\rightarrow nSePC_f$ (i.e. chaque $A(y)$ est suppos\'e fibrant) et on va
d\'efinir une $n$-pr\'ecat de Segal munie d'un foncteur vers $Y$
$$
\int _YA \rightarrow Y.
$$
Ceci g\'en\'eralise la {\em
construction de Grothendieck} (cf \cite{ThomasonLimits} et \cite{DwyerKan2}) qui
concerne le cas o\`u $A$ est un pr\'efaisceau de $1$-cat\'egories.

Les objets de $\int _YA$ sont les couples $(y,a)$ avec $y\in Y$ et
$a\in A(y)$. Pour une suite $(y_0,a_0),\ldots , (y_p,a_p)$
et une suite de fl\`eches $\varphi _i : y_{i-1}\rightarrow y_i$ de $Y$ on
construira la $n-1$-cat\'egorie
$$
(\int _YA)_{p/}((y_0,a_0),\ldots , (y_p,a_p); \varphi _1,\ldots , \varphi _p)
$$
qui sera l'image inverse de
$$
(\varphi _1,\ldots , \varphi _p)\in Y_{p}(y_0,\ldots , y_p).
$$
On fait cette construction par r\'ecurrence sur $p$; en supposant que c'est fait
pour $p-1$ on pose $(\ast )$
$$
(\int _YA)_{p/}((y_0,a_0),\ldots , (y_p,a_p); \varphi _1,\ldots , \varphi _p)
:=
$$
$$
(\int _YA)_{p-1/}((y_0,a_0),\ldots , (y_{p-1},a_{p-1});
\varphi _1,\ldots , \varphi
_{p-1})
$$
$$
\times _{A(y_{p})_{p-1/}(\phi _0(a_0), \ldots , \phi _{p-1}(a_{p-1}))}
A(y_{p})_{p/}(\phi _0(a_0), \ldots , \phi _{p}(a_{p}))
$$
o\`u les
$$
\phi _i : A(y_i)\rightarrow A(y_p)
$$
sont les morphismes compos\'es des $res_{\varphi _p}\cdots res _{\varphi
_{i+1}}$. Ici par ailleurs, on a par construction un morphisme
$$
(\int _YA)_{p/}((y_0,a_0),\ldots , (y_p,a_p); \varphi _1,\ldots , \varphi _p)
\rightarrow
A(y_{p})_{p/}(\phi _0(a_0), \ldots , \phi _{p}(a_{p}))
$$
et le premier morphisme dans le produit fibr\'e $(\ast )$ est ce morphisme pour
$p-1$, compos\'e avec $res_{\varphi _p}$ pour passer de $A(y_{p-1})_{p-1/}$ \`a
$A(y_p)_{p-1/}$.

Maintenant on pose
$$
(\int _YA)_{p/}((y_0,a_0),\ldots , (y_p,a_p))
:=
\coprod _{\varphi _1,\ldots , \varphi _p}(\int _YA)_{p/}((y_0,a_0),\ldots ,
(y_p,a_p); \varphi _1,\ldots , \varphi _p),
$$
o\`u le coproduit est pris sur tous les
$$
(\varphi _1,\ldots , \varphi _p)\in Y_{p}(y_0,\ldots , y_p).
$$
Ceci d\'efinit bien la $n$-pr\'ecat de Segal $\int _YA$ avec une application
tautologique vers $Y$.

On v\'erifie que $\int _YA$ est une $n$-cat\'egorie de Segal avec
$$
(\int _YA)_{1/}((y_0,a_0), (y_1,a_1))=\coprod _{\varphi _1: y_0\rightarrow y_1}
A(y_1)_{1/}(\varphi _1(a_0), a_1).
$$
Le fait que $\int _YA$ est une $n$-cat\'egorie de Segal utilise le fait que
chaque $A(y_p)$ est fibrant, et plus particuli\`erement le fait qu'en
cons\'equence le deuxi\`eme morphisme dans le produit fibr\'e $(\ast )$ est
fibrant.

Si $A$ est un pr\'efaisceau de $1$-cat\'egories sur $Y^o$ (les $A(y)$ sont
dans ce cas fibrants) alors $\int _AY$ est la cat\'egorie dont les objets sont
les couples $(y,a)$ et les morphismes de $(y,a)$ vers $(z,b)$ sont les couples
$(f,r)$ o\`u $f: y\rightarrow z$ dans $Y$ et $r: f(a)\rightarrow b$ dans $A(z)$.
C'est la {\em construction de Grothendieck} \cite{SGA1} \cite{ThomasonLimits}
\cite{DwyerKan2}.

Une premi\`ere utilisation de cette construction est la {\em subdivision
barycentrique}, suivant le point de vue de Thomason \cite{ThomasonLimits}.
Soit $X$ un ensemble simplicial et posons
$$
B:= \int _{\Delta ^o} X.
$$
C'est une cat\'egorie, qui est aussi la {\em cat\'egorie des simplexes de $X$}
not\'ee $\Delta X$ dans \cite{Hirschhorn} \cite{DHK}.
Son nerf $\nu B$ s'explicite de la fa\c{c}on suivante: les \'el\'ements de
$\nu B_p$ sont les
$$
(n,\sigma , x)= (n_0,\ldots , n_p; \sigma _1,\ldots , \sigma _p; x)
$$
o\`u $n_i\in \Delta$, $\sigma _i: n_{i-1}\rightarrow n_i$ dans $\Delta$,
et $x\in X(n_p)$. On note
$$
g(n,\sigma ): p\rightarrow n_p
$$
le morphisme qui \`a $i'\in p$ associe l'image
$$
g(n,\sigma )(i'):= \sigma _p\cdots \sigma _{i+1}(n'_i)
$$
o\`u $n'_i\in n_i$ est le dernier \'el\'ement (on utilise ici le ``prime''
pour distinguer entre les \'el\'ements des ensembles ordonn\'es, et les
ensembles
qui sont eux-m\^emes objets de $\Delta$). On  a un morphisme
d'ensembles simpliciaux
$$
\psi : \nu B \rightarrow X
$$
d\'efini par
$$
\psi (n,\sigma , x):= g(n,\sigma )^{\ast} (x)\in X(p).
$$
D'autre part, soit $\Delta ^r$ la sous-cat\'egorie de $\Delta$ o\`u les
morphismes sont seulement ceux qui envoient le dernier \'el\'ement sur le
dernier \'el\'ement. On pose
$$
D:=\int _{\Delta ^r} X|_{(\Delta ^r)^o} \subset B.
$$
Les \'el\'ements de $\nu D_p$ sont les $(n,\sigma , x)$ tels que
$\sigma _i(n'_{i-1})=n'_i$. En particulier, le morphisme
$g(n,\sigma )$ se factorise \`a travers la projection $p\rightarrow 0$
via l'inclusion $0\rightarrow n_p$ correspondant au dernier \'el\'ement $n'_p$.
On a la factorisation
$$
\nu D\rightarrow X_0 \rightarrow X.
$$
Par cons\'equent le morphisme $\psi$ s'\'etend directement en un morphisme
$$
L\psi : 1SeL(B,D)\rightarrow X
$$
car tout morphisme $I\rightarrow D$ s'envoie dans $X$ sur un objet $\{ x\}$,
et s'\'etend donc automatiquement en $\overline{I}\rightarrow \{ x\}$.

Le lemme suivant doit en r\'ealit\'e \^etre d\^u \`a Thomason
\cite{ThomasonLimits}
(en tout cas,
il montre que le morphisme sur les compl\'et\'es groupiques, est une
\'equivalence).

\begin{lemme}
\label{bd}
Pour tout ensemble simplicial $X$, le morphisme
$L\psi : 1SeL(B,D)\rightarrow X$ ci-dessus est une \'equivalence faible de
$1$-pr\'ecats de Segal.
\end{lemme}
{\em Preuve:}
La formation de $B$, $D$, $1SeL(B,D)$ et des morphismes $\psi$ et $L\psi$
est compatible aux colimites. Comme tout ensemble simplicial est colimite
d'ensembles simpliciaux de la forme $h(m)$ (le pr\'efaisceau repr\'esent\'e par
$m$, i.e. le $m$-simplexe standard), il suffit de prouver le lemme pour
$X=h(m)$. Dans ce cas $B=\Delta /m$ et $D$ est la sous-cat\'egorie des objets de
$\Delta ^r$ au-dessus (par rapport \`a $\Delta$) de $m$, not\'ee $i/m$ o\`u $i:
\Delta ^r\rightarrow \Delta $ est l'inclusion. Soit $I^{(m)}$ la cat\'egorie
avec pour objets $0'',\ldots , m''$ et un morphisme $i''\rightarrow j''$ pour
$i\leq j$. On a une \'equivalence faible de $1$-pr\'ecats de Segal
$h(m)\rightarrow I^{(m)}$.
On a un foncteur $a: I^{(m)}\rightarrow B$ qui \`a $i''$ associe l'objet
$i\rightarrow m$ o\`u ce morphisme envoie $j'\in i$ sur $j'\in m$
(ici encore, on note avec un ``prime'' les \'el\'ements des $p\in \Delta$).
D'autre part on a un foncteur $b: B\rightarrow I^{(m)}$ qui \`a un objet
$p\rightarrow m$ de $B$, associe $i''\in I^{(m)}$  o\`u (pour le m\^eme
$i$) le dernier \'el\'ement dans l'image de $p\rightarrow m$ est $i'$.
On note que $b$ se factorise \`a travers le morphisme de $1$-pr\'ecats de Segal
$B\rightarrow h(m)$ via l'inclusion $h(m)\rightarrow I^{(m)}$. Aussi, $b$
envoie les morphismes de $D$ sur les identit\'es de $I^{(m)}$. Pour prouver le
lemme, il suffit donc de prouver que $b$ induit une \'equivalence entre
$1SeL(B,D)$ et $I^{(m)}$. Or, d'apr\`es  la proposition \ref{dklocsegloc}, il
suffit de prouver que $b$ induit une \'equivalence
entre $L(B,D)$ et $I^{(m)}$. Pour ceci, on applique le Corollaire 3.6 de
\cite{DwyerKan2} sous la forme du lemme \ref{nonadjoint} ci-dessus:  on a que
$ba$ est l'identit\'e de $I^{(m)}$, et il y a une transformation
naturelle $\eta : 1\rightarrow ab$ dont la valeur sur $p\rightarrow m$
(objet de $\Delta /m$) est le morphisme $p\rightarrow i \rightarrow m$ o\`u $i'$
est le dernier \'el\'ement dans l'image de $p\rightarrow m$ et  $i\rightarrow
m$ est $a(i'')$, l'injection qui envoie $j'$ sur $j'$ pour $j'\in i$.
On note que les $\eta _{p\rightarrow m}$ sont des morphismes dans $D$, donc
cette transformation naturelle induit une d\'eformation-r\'etraction de
$L(B,D)$ sur
$I^{(m)}$.
\eop

Lorsqu'il sera question d'appliquer cette construction \`a plusieurs ensembles
simpliciaux $X$, on notera $\beta (X):= B$ et $\delta (X):= D\subset \beta (X)$.
Le lemme dit qu'on a une \'equivalence faible
$$
1SeL(\beta (X), \delta (X))\stackrel{\cong}{\rightarrow} X.
$$

On donne maintenant une variante qui peut \^etre utile.  Cette
variante entre dans le cadre
plus g\'en\'eral de la th\'eorie des ``cat\'egories-test'' envisag\'ee par
Grothendieck dans \cite{Grothendieck}. Ici on  utilise un exemple qui provient
des travaux de Thomason \cite{ThomasonLimits} \cite{ThomasonCatMod}; le
deuxi\`eme auteur remercie G. Maltsiniotis et A. Brugui\`eres de lui avoir
expliqu\'e cette id\'ee.

Pour $p\in \Delta$ soit ${\rm poset}(p)$
l'ensemble partiellement ordonn\'e des hyperfaces dans le $p$-simplexe standard
(incluant le $p$-simplexe lui-m\^eme). On consid\`ere ${\rm poset}(p)$ comme une
cat\'egorie. Pour ce qu'on en fera, la direction des fl\`eches n'est pas
importante; on peut fixer un choix en disant que le $p$-simplexe est l'objet
initial de  ${\rm poset}(p)$. On obtient un foncteur
$$
{\rm poset}: \Delta \rightarrow Cat.
$$
On d\'efinit aussi la sous-cat\'egorie ${\rm Dposet} (p)\subset
{\rm poset}(p)$ constitu\'ee des hyper-faces contenant le premier sommet. On a
d'autre part le foncteur $p\mapsto I^{(p)}$ o\`u $I^{(p)}$ est la cat\'egorie
avec $p+1$ objets $0',\ldots , p'$ et un morphisme $i'\rightarrow j'$ pour
$i\leq j$. On notera ${\rm std}: \Delta \rightarrow Cat$ ce foncteur.
Le foncteur
$$
{\rm poset}(p)\rightarrow {\rm std}(p)= I^{(p)}
$$
qui envoie
un hyper-face du $p$-simplexe sur son premier sommet consid\'er\'e
comme objet de $I^{(p)}$, donne une transformation naturelle
$$
{\rm poset} \rightarrow {\rm std}.
$$
de foncteurs $\Delta
\rightarrow Cat$.
Cette transformation naturelle envoie la sous-cat\'egorie $ {\rm Dposet}$
sur les identit\'es dans ${\rm std}$ et le morphisme de $1$-pr\'ecats de Segal
$$
1SeL({\rm poset}(p), {\rm Dposet}(p) )\rightarrow {\rm std}(p)=I^{(p)}
$$
est une \'equivalence faible.

Si $X$ est un ensemble simplicial alors $X$ peut \^etre vu comme le coproduit
de ses simplexes, et on peut \'ecrire
$$
X= \coprod _p\coprod _{x\in X_p}\nu I^{(p)} /\sim .
$$
Si on effectue le m\^eme coproduit, mais avec les $\nu {\rm poset}(p)$
(resp. $\nu {\rm Dposet}(p)$) comme
facteurs, on obtient un ensemble simplicial
$$
\beta ^{\rm poset} (X):= \coprod _p \coprod _{x\in X_p}{\rm poset}(p) /\sim
$$
(resp.
$$
\delta ^{\rm poset} (X):= \coprod _p \coprod _{x\in X_p} {\rm Dposet}(p)
/\sim  )
$$
et le morphisme
$$
1SeL(\beta ^{\rm poset} (X), \delta ^{\rm poset} (X))\rightarrow
X
$$
est une \'equivalence faible de $1$-pr\'ecats de Segal. Par contre, on peut
remarquer que $\beta ^{\rm poset} (X)$ est le nerf de la cat\'egorie
sous-jacente \`a un ensemble partiellement ordonn\'e, et
$\delta ^{\rm poset} (X)$ une sous-cat\'egorie. Ceci donne la variante suivante
du lemme \ref{bd}.

\begin{lemme}
\label{bdvariant}
Pour tout ensemble simplicial $X$, il y a une cat\'egorie
$\beta ^{\rm poset} (X)$ sous-jacente \`a un ensemble partiellement ordonn\'e,
avec une sous-cat\'egorie $\delta ^{\rm poset} (X)$ et une
\'equivalence faible de $1$-pr\'ecats de Segal
$$
1SeL(\beta ^{\rm poset} (X), \delta ^{\rm poset} (X))
\stackrel{\cong}{\rightarrow}
X.
$$
\end{lemme}
\eop

On appliquera cela au nerf $X$ d'une $1$-cat\'egorie $Y$. Dans ce cas
le morphisme
$$
\beta ^{\rm poset} (X)\rightarrow Y
$$
est un morphisme de
cat\'egories et le lemme implique que le morphisme
$$
L(\beta ^{\rm poset} (X), \delta ^{\rm poset} (X))
\stackrel{\cong}{\rightarrow}
Y
$$
est une \'equivalence de cat\'egories simpliciales.

\subnumero{Sections de l'int\'egrale}

Si $F\rightarrow Y$ est un morphisme (suppos\'e fibrant) d'une
$n$-cat\'egorie de Segal vers une $1$-cat\'egorie $Y$, on note $Sect
(Y,F)$ la
$n$-cat\'egorie de Segal des sections. Si $F$ est une $1$-cat\'egorie alors
$Sect(Y,F)$ est une $1$-cat\'egorie. Les objets sont les sections $\sigma :
Y\rightarrow F$ et les morphismes entre sections sont les transformations
naturelles qui se projettent sur la transformation naturelle
identit\'e du
foncteur identit\'e de $Y$.

Si $A$ est un $n$-pr\'echamp de Segal au-dessus de $Y$ on dispose donc
de la
$n$-pr\'ecat
de Segal
$$
Sect(Y, \int _YA).
$$
Celle-ci pourrait ne pas avoir les propri\'et\'es escompt\'ees. On
commence donc par choisir
un remplacement fibrant (au-dessus de $Y$), de la forme
$$
\int _YA\rightarrow \int '_YA \rightarrow Y
$$
o\`u la  premi\`ere fl\`eche est une \'equivalence et la deuxi\`eme est
une fibration; et on consid\`ere plut\^ot $Sect(Y,\int '_YA)$
qui est une $n$-cat\'egorie de Segal fibrante.

\begin{lemme}
\label{inchange}
Avec les notations pr\'ec\'edentes, pour tout $y\in Y$, soit $A'(y)$ la fibre de
$\int '_YA$ (i.e. le produit fibr\'e de ce dernier avec $\ast = \{ y\}$
au-dessus de $Y$). Alors le morphisme
$$
g:\int _YA\rightarrow \int '_YA
$$
induit une \'equivalence $g(y): A(y)\stackrel{\cong}{\rightarrow} A'(y)$.
\end{lemme}
{\em Preuve:}
Le fait que $g$ est (par hypoth\`ese) une \'equivalence implique
directement que le morphisme $g(y):A(y)\rightarrow A'(y)$ est pleinement
fid\`ele. On montre que $g(y)$ est essentiellement surjectif.
Si $(y,a)\in A'(y)$ alors il y a une \'equivalence dans $\int '_YA$ entre
$(y,a)$ et un objet $(z,b)\in A(z)$ pour un autre $z\in Y$. Cet \'equivalence
qu'on notera $I\rightarrow \int '_YA$ s'\'etend en $\eta :
\overline{I}\rightarrow \int '_YA$ qui se projette donc vers un morphisme
$\psi :
\overline{I}\rightarrow Y$ qui correspond \`a un isomorphisme entre $y$ et $z$.
Maintenant on peut relever ceci en un morphisme ``horizontal'' (voir les
descriptifs dans \ref{casi} et \ref{gammasect} ci-dessous)
$h:\overline{I}\rightarrow \int _YA$. Soit maintenant  $\overline{I}^{(2)}$ la
cat\'egorie avec pour objets $0,1,2$ et un isomorphisme entre chaque paire
d'objets,
contenant des sous-cat\'egories \'equivalentes \`a $\overline{I}$ qu'on notera
$\overline{I}^{01}$, $\overline{I}^{12}$, $\overline{I}^{02}$. Toutes ces
cat\'egories sont contractiles i.e. \'equivalentes \`a $\ast$. On a d'autre part
une projection
$$
p_{0,12}: \overline{I}^{(2)}\rightarrow \overline{I}
$$
qui envoie $0$ sur $0$ et $1,2$ sur $1$. Avec ceci on obtient un morphisme
$$
\psi \circ p_{0,12}: \overline{I}^{(2)}\rightarrow Y,
$$
avec rel\`evement partiel sur la r\'eunion $\overline{I}^{01}\cup
\overline{I}^{02}$ en un morphisme
$$
\eta \sqcup ^{(z,b)}h: \overline{I}^{01}\cup
\overline{I}^{02}\rightarrow \int '_YA.
$$
Ce dernier envoie $0$ sur $(z,b)$, envoie $\overline{I}^{01}$ via $\eta$, et
envoie $\overline{I}^{02}$ via $h$.  Comme l'inclusion
$$
\overline{I}^{01}\cup
\overline{I}^{02}\subset
\overline{I}^{(2)}
$$
est une \'equivalence (donc une cofibration triviale), il y a un rel\`evement
et extension en un morphisme
$$
H: \overline{I}^{(2)}\rightarrow \int '_YA.
$$
La restriction $H^{12}$ de $H$ \`a $\overline{I}^{12}$ reste au-dessus de la
sous-cat\'egorie $\{ y\} \subset Y$ (contenant l'objet $y$ et l'identit\'e
comme seul morphisme), donc on peut \'ecrire
$$
H^{12}: \overline{I}^{12}\rightarrow A'(y).
$$
Aussi $H^{12}$ envoie
le sommet $1$ sur $\eta (1)= (y,a)$, et envoie le sommet $2$ sur un objet de
$A(y)$. On obtient donc que $(y,a)$ est \'equivalent \`a un objet de $A(y)$ ce
qui prouve l'essentielle surjectivit\'e de $g(y)$.
\eop

\begin{corollaire}
La formation de $\int '_YA$ commute au changement de base: si $f:Z\rightarrow Y$
alors $(\int '_YA)\times _YZ$ est un choix possible pour $\int '_Z(A|_Z)$.
\end{corollaire}
{\em Preuve:}
Le morphisme $(\int '_YA)\times _YZ\rightarrow Z$ est fibrant, et il ne
reste qu'\`a
prouver que
$$
\int _Z(A|_Z)\rightarrow  (\int '_YA)\times _YZ
$$
est une \'equivalence. Au vu du lemme, ce morphisme est essentiellement
surjectif.
Si $\varphi : u\rightarrow v$ est une fl\`eche de $Z$ (qui donne une fl\`eche
$f\varphi : fu\rightarrow fv$ de $Y$), alors, pour $a\in A(fu)$ et $b\in A(fv)$,
\begin{eqnarray*}
(\int _ZA|_Z)_{1/}((u,a), (v,b); \varphi ) &=& A(fv)_{1/}(res_{f\varphi}(a),
b)\\
&= &(\int _YA)_{1/}((fu,a), (fv,b); f\varphi ) \\
& = & ((\int _YA)\times _YZ)_{1/}(((fu,a), u), ((fv,b), v); \varphi )\\
& \stackrel{\cong}{\rightarrow}&
((\int' _YA)\times _YZ)_{1/}(((fu,a), u), ((fv,b), v); \varphi ).
\end{eqnarray*}
La d\'erni\`ere \'equivalence
provient du fait que $\int _YA\rightarrow \int '_YA$ est une
\'equivalence. Ceci montre que le morphisme
$$
\int _Z(A|_Z)\rightarrow  (\int '_YA)\times _YZ
$$
est pleinement fid\`ele.
\eop

On retourne maintenant \`a la consid\'eration des sections de $\int '_YA$.
Un objet de \linebreak
$Sect(Y,\int '_YA)$ est un morphisme $\sigma : Y\rightarrow
\int '_YA$. Pour tout morphisme $f:y\rightarrow z$ dans $Y$ on obtient un
morphisme $\sigma ^r(f): res_f(\sigma (y))\rightarrow \sigma (z)$, bien d\'efini
\`a homotopie pr\`es comme morphisme de $A(z)$. On dira que $\sigma$ est une
$eq$-section si les $\sigma ^r(f)$ d\'efinis ci-dessus sont des \'equivalences.
On note $Sect^{\rm eq}(Y,\int '_YA)$ la sous-$n$-pr\'ecat pleine de
$Sect(Y,\int '_YA)$  dont les objets sont les
$eq$-sections.

Si $A$ est un pr\'efaisceau de $1$-cat\'egories, la section $\sigma$
est d\'etermin\'ee par les morphismes $\sigma ^r(f)$ (qui sont bien d\'efinis
dans ce cas).

Avant de traiter le cas g\'en\'eral des sections au-dessus d'une cat\'egorie
$Y$, on traite le cas (en quelque sorte irr\'eductible) de la cat\'egorie $Y=I$
avec deux objets $0,1$ et un morphisme $0\rightarrow 1$.

\begin{lemme}
\label{casi}
Soit $A$ un $n$-pr\'echamp de Segal au-dessus de $I^o$ dont les valeurs
$A(0)$ et
$A(1)$ sont fibrantes. Soit $\int '_IA$ le remplacement fibrant pour $\int  _IA$
au-dessus de $Y$ et notons $A'(0)$ la fibre de ce dernier au-dessus de $0\in
I$. Alors le morphisme d'\'evaluation
$$
ev(0): Sect ^{\rm eq}(I, \int '_IA)\rightarrow A'(0)
$$
est une \'equivalence.
\end{lemme}
{\em Preuve:}
On rappelle d'abord qu'un $n$-pr\'echamp de Segal $A$ au-dessus de $I^o$
consiste
en un triplet $(A(0), A(1), r: A(0)\rightarrow A(1))$.
Soit $a\in A(0)_0$, et
soit $r(a)$ son image dans $A(1)$. On va construire une section
``horizontale''
$\sigma :I\rightarrow\int _IA$ qui envoie $0$ sur $a$ et $1$ sur $r(a)$. Un
simplexe dans le nerf $\nu I$ est de la forme $(0,\ldots , 0, 1, \ldots
1)$ (disons, avec $p$ fois $0$ et $q$ fois $1$).
Rappelons que
$$
A(0)_{p+q/}((0, a),\ldots , (0,a),(1,r(a)),\ldots , (1,r(a)))
$$
$$
= A(0)_{p/}(a,\ldots , a)\times _{A(1)_{p/}(r(a),\ldots , r(a))}
A(1)_{p+q/}(r(a),\ldots , r(a)).
$$
On sp\'ecifie $\sigma$ en disant que $\sigma (0,\ldots, 0, 1,\ldots , 1)$
est la paire des deux simplexes d\'eg\'en\'er\'es dans ce produit fibr\'e.
Noter que $\sigma$ est une $eq$-section.

Cette construction montre que tout objet de $A'(0)$ qui provient d'un objet de
$A(0)$, est dans l'image essentielle du morphisme $ev(0)$. Au vu du lemme
\ref{inchange}, ceci montre que $ev(0)$ est essentiellement surjectif.

Il reste maintenant \`a montrer que $ev(0)$ est pleinement fid\`ele.
L'argument ci-dessus donne en fait un morphisme $A(0)\rightarrow
Sect ^{\rm eq}(I, \int '_IA)$ qui scinde $ev(0)$.
On dira qu'un morphisme provenant de ce scindage est ``horizontal''.
Tout morphisme $f:(0,u)\rightarrow (1,v)$ dans $\int '_IA$
(i.e. objet  $f\in (\int '_IA)_{1/}((0,u),(1,v))_0$) s'\'ecrit \`a
\'equivalence pr\`es comme le ``compos\'e'' $gh$ (ce
compos\'e \'etant d\'efini \`a
\'equivalence pr\`es) o\`u $h$ est horizontal et $g$ un morphisme de $A'(1)$.
Si $a: I\rightarrow \int '_IA$ alors l'image $a(\varphi )$ est bien d\'efinie
comme morphisme $f\in (\int '_IA)_{1/}((0,u),(1,v))_0$. D'apr\`es la
d\'efinition de $eq$-section, $a$ est une $eq$-section si et seulement si
$f=gh$ avec $g$ une \'equivalence dans $A'(1)$. Il s'ensuit
(nous laissons au lecteur le soin de d\'etailler ce point) que toute
$eq$-section est \'equivalente (dans $Sect(I, \int '_IA)$) \`a une section
horizontale.

Pour montrer que $ev(0)$ est pleinement fid\`ele on peut alors consid\'erer
deux sections horizontales $a,b: I\rightarrow \int _IA$, et on voudrait
montrer que le morphisme
$$
Sect ^{\rm eq}(I, \int '_IA)_{1/}(a,b)\rightarrow A'(1)_{1/}(a(1), b(1))
$$
est une \'equivalence.

On a l'\'egalit\'e suivante
(formule par foncteur repr\'esentable
de $n-1$-pr\'ecats de Segal)
$$
Sect ^{\rm eq}(I, \int '_IA)_{1/}(a,b)= \left(
E\mapsto Hom ^{a,b/I}(I\times \Upsilon (E), \int '_IA)
\right)
$$
o\`u $Hom ^{a,b/I}$ est l'ensemble des morphismes qui se restreignent en $a$
et $b$ sur $I \times \{ 0\}$ et $I\times \{ 1\}$ et qui se projettent vers
la deuxi\`eme projection $I\times \Upsilon (E)\rightarrow I$.
Le morphisme d'evaluation correspond \`a la restriction de morphismes $I\times
\Upsilon (E)\rightarrow \int '_IA$ sur des morphismes $\{ 0\} \times \Upsilon
(E) \rightarrow A'(0)$. On voudrait montrer que cette restriction est une
\'equivalence.  Divisons le carr\'e $I\times \Upsilon (E)$ en deux triangles
not\'es $\Upsilon ^2(E,\ast )$ et $\Upsilon ^2(\ast , E)$ recoll\'es le long de
la diagonale $\Upsilon (E)$. Le foncteur (en $E$) des morphismes sur le premier
triangle $\Upsilon ^2(E, \ast )\rightarrow \int '_IA$ est \'equivalent, par
restriction \`a l'ar\^ete $01$, \`a $A'(0)_{1/}(a(0), b(0))$.
Il suffit donc de voir que l'application de restriction entre le foncteur
en $E$ des
morphismes avec source le carr\'e, vers le foncteur en $E$ des morphismes
avec source le premier
triangle, est une \'equivalence.
Maintenant cette question ne concerne que le
deuxi\`eme triangle: il suffit de voir que le morphisme (de restriction sur
l'ar\^ete $02$)
$$
r_{02}:\left(
E\mapsto Hom ^{a,b(1)/I}(\Upsilon ^2(\ast , E), \int '_IA)
\right) \rightarrow
$$
$$
\left(
E\mapsto Hom ^{a(0),b(1)/I}( \Upsilon (E), \int '_IA)
\right)
$$
est une \'equivalence de $n-1$-pr\'ecats de Segal. Ici
$Hom ^{a,b(1)/I}$ est l'ensemble des morphismes qui donnent $a$ sur $\Upsilon
(\ast )$ (ar\^ete $01$),qui donnent $b(1)$ sur le sommet $2$, et qui sont
compatibles avec la projection vers $I$; et
$Hom ^{a(0),b(1)/I}$ est l'ensemble des morphismes qui donnent  $a(0)$ sur le
premier sommet $0$ et $b(1)$ sur le deuxi\`eme sommet (qu'on note ici $2$, et
qui sont compatibles avec la projection vers $I$. Mais
le l'arriv\'ee
du morphisme
$r_{02}$ est juste ce qu'on a not\'e
$$
(\int '_IA)_{1/}((0,a(0)), (1,b(1)); \varphi )
$$
ci-dessus, o\`u  $\varphi : 0\rightarrow 1$ est le morphisme de $I$. Et
\`a cause de l'\'equivalence
faible $\Upsilon ^2(\ast , E)\cong \Upsilon (\ast )\sqcup ^{\{ 1\} } \Upsilon
(E)$,
le depart
du morphisme $r_{02}$ ci-dessus est \'equivalent via
$r_{12}$ \`a
$$
\left(
E\mapsto Hom ^{a(1),b(1)}(\Upsilon ( E), \int '_IA)
\right) = A'(1)_{1/}(a(1), b(1)).
$$
Rappelons qu'on a suppos\'e que $a$ correspond \`a un morphisme horizontal
$I\rightarrow \int _IA$, i.e. $a(1)=ra(0)$ pour $r$ le morphisme structurel de
$A$. On dispose donc du m\^eme morphisme que $r^{02}$, qu'on appellera
$\tilde{r}^{02}$, qu'on aurait pu fabriquer avec $\int _IA$ au lieu de $\int
'_IA$, ce qui est  (via les m\^emes identifications que ci-dessus mais pour
$\int _IA$) le morphisme de composition avec $a(\varphi ):a(0)\rightarrow a(1)$
$$
A(1)_{1/}(a(1), b(1))\rightarrow
(\int _IA)_{1/}((0,a(0)), (1,b(1)); \varphi ).
$$
Mais ce morphisme est par construction un isomorphisme. Le fait que
$\int _IA\rightarrow \int '_IA$ est une \'equivalence de cat\'egories
implique que dans le diagramme
$$
\begin{array}{ccccc}
A(1)_{1/}(a(1), b(1))&\stackrel{=}{\leftarrow} &
{\rm source(\tilde{r}^{02})} &\stackrel{\tilde{r}^{02}}{\rightarrow}&
(\int _IA)_{1/}((0,a(0)), (1,b(1)); \varphi )\\
\downarrow && \downarrow && \downarrow \\
A'(1)_{1/}(a(1), b(1)) &
\stackrel{\cong}{\leftarrow} &
{\rm source(r^{02})} &\stackrel{r^{02}}{\rightarrow}&
(\int '_IA)_{1/}((0,a(0)), (1,b(1)); \varphi )
\end{array}
$$
les fl\`eches verticales sont des \'equivalences. Le fait que le morphisme en
haut \`a droite soit un isomorphisme (par construction---et c'est ici
l'endroit essentiel o\`u on utilise le fait que $a$ soit horizontale) implique
que $r^{02}$ est une \'equivalence.  Ceci termine la d\'emonstration.
\eop

\begin{proposition}
\label{gammasect}
Soit $A$ un $n$-pr\'echamp de Segal sur $Y^o$. Alors il y a un morphisme
naturel de
$n$-pr\'ecats de Segal
$$
\Gamma (Y^o,A)\rightarrow Sect (Y,\int _YA),
$$
et si $A$ est fibrant alors le morphisme compos\'e
$$
\Gamma (Y^o,A)\rightarrow Sect (Y,\int '_YA)
$$
est pleinement fid\`ele avec pour image essentielle
$Sect^{\rm eq} (Y,\int '_YA)$.
\end{proposition}
{\em Preuve:}
D'abord on d\'efinit le morphisme. Au niveau des objets, c'est facile: une
section  $\sigma$ de
$A$ au-dessus de $Y^o$ donne directement un morphisme
$Y\rightarrow \int _YA$ qui envoie $y$ sur $(y,\sigma (y))$.
Le simplexe
$(y_0,\ldots , y_p; \varphi _1,\ldots , \varphi _p)$ va sur l'\'el\'ement de
$$
(\int _YA)_{p/}((y_0,\sigma (y_0)),\ldots , (y_p, \sigma (y_p));
\varphi _1,\ldots , \varphi _p)
$$
dont les deuxi\`emes projections dans les produits fibr\'es $(\ast )$
(par la r\'ecurrence il y en a une pour chaque $0<p\leq p$) sont les
simplexes d\'eg\'en\'er\'es dans $A_{q/}(\sigma (y_q), \ldots , \sigma (y_q))$.

Pour sp\'ecifier le morphisme plus pr\'ecisement, on utilise la d\'efinition de
$\Gamma (Y^o, A)$ par propri\'et\'e universelle: pour une $n$-pr\'ecat de Segal
$E$, un morphisme $E\rightarrow \Gamma (Y^o , A)$ est la m\^eme chose qu'un
morphisme $\underline{E}\rightarrow A$ o\`u $\underline{E}$ est le
$n$-pr\'echamp
de Segal constant sur $Y^o$ \`a valeurs $E$. On a (avec la m\^eme d\'efinition
que ci-dessus m\^eme si $E$ n'est pas fibrant)
$$
\int _Y\underline{E} = Y\times E.
$$
Par fonctorialit\'e de la construction $\int _Y$, un morphisme
$\underline{E}\rightarrow A$ fournit un morphisme
$$
Y\times E =\int _Y\underline{E} \rightarrow \int _YA
$$
dont la composition avec la projection sur $Y$ est la premi\`ere projection
$Y\times E\rightarrow Y$. Ceci donne, par d\'efinition de $Sect (Y, -)$, un
morphisme
$$
E\rightarrow Sect (Y, \int _YA)
$$
et on a d\'efini le morphisme $\Gamma (Y^o, A)\rightarrow Sect (Y, \int _YA)$.

Pour montrer l'\'equivalence,
on voudrait d\'ecomposer $Y$ en cellules et traiter chaque cellule
individuellement. La premi\`ere partie de l'argument consiste donc \`a r\'eduire
la proposition au cas $Y=I^{(p)}$.  Dans la deuxi\`eme partie de l'argument,
on traitera ce cas particulier (et pour ceci on r\'eduit encore facilement
au cas
$Y=I$). Pour la r\'eduction de la premi\`ere partie, on a besoin de la
subdivision
barycentrique d\'efinie ci-dessus. Cet argument sera en quelque sorte un
pr\'ecurseur de l'utilisation qu'on fera des {\em cat\'egories de Reedy} au \S
18 ci-dessous; mais pour le moment nous n'avons pas besoin de ce formalisme.

Pour une cat\'egorie $Y$ on notera (avec un l\'eger abus de notation)
$$
\beta (Y):= \beta (\nu Y)\stackrel{\psi}{\rightarrow} Y
$$
la subdivision barycentrique du nerf de $Y$, avec foncteur $\psi$ vers $Y$.
(Aussi on rappelle que, par d\'efinition, il y a un foncteur structurel
$\beta (Y)\rightarrow \Delta ^o$.) Un
objet de la cat\'egorie $\beta (Y)$ est une suite $(y_0,\ldots , y_p;
f_1,\ldots , f_p)$ composable de fl\`eches $f_i: y_{i-1}\rightarrow y_i$ de $Y$.
Un morphisme de $\beta (Y)$ se d\'eduit d'un morphisme de $\Delta ^o$, i.e. d'un
morphisme $q\rightarrow p$ mais allant dans l'autre sens.
Si $u: q\rightarrow p$ est d\'efinie par la suite $u_0,\ldots , u_q\in \{
0,\ldots , p\}$ alors les morphismes correspondants sont de la forme
$$
(y_0,\ldots , y_p;
f_1,\ldots , f_p)\rightarrow (y_{u_0}, \ldots , y_{u_q};
g_1, \ldots , g_q)
$$
o\`u $g_i = f_{u_i} \cdots f_{u_{i-1}+1}$ (si $u_i=u_{i-1}$ alors $g_i$ est
l'identit\'e). Le foncteur vers $Y$ est
$$
(y_0,\ldots , y_p;
f_1,\ldots , f_p)\mapsto y_0.
$$
On note qu'en fait on utilisera souvent l'oppos\'e $\beta (Y)^o$ et que les
morphismes dans $\beta (Y)^o$, qui vont dans l'autre sens que ci-dessus,
peuvent \^etre consid\'er\'es comme des ``inclusions (ou d\'eg\'en\'erescences)
de suites''.

On dira
qu'une suite $(y_0,\ldots , y_p; f_1,\ldots , f_p)$ (qui correspond \`a un
objet de $\beta (Y)$) est {\em non-d\'eg\'en\'er\'ee} si $f_i\neq 1_{x_{i-1}}$
pour tout $i$. A toute suite d\'eg\'en\'er\'ee, on associe la suite
non-d\'eg\'en\'er\'ee obtenue en enlevant toutes les identit\'es. Soit $\beta
_m(Y)\subset \beta (Y)$ la sous-cat\'egorie pleine constitu\'ee
par les objets dont la
suite  non-d\'eg\'en\'er\'ee correspondante est de longueur $\leq m$.  D'autre
part, pour chaque suite non-d\'eg\'en\'er\'ee  $(p; y,f)=(y_0,\ldots , y_p;
f_1,\ldots , f_p)$, on consid\`ere le morphisme
$$
\eta (p; y,f): I^{(p)}\rightarrow Y
$$
o\`u $I^{(p)}$ est la cat\'egorie avec $p+1$ objets $0',\ldots , p'$ et un
morphisme $i'\rightarrow j'$ pour tout couple $i\leq j$. Le morphisme
$\eta (p; y,f)$ envoie $i'$ sur $y_i$ et la fl\`eche $(i-1)'\rightarrow i'$ sur
$f_i$. Ce morphisme induit
$$
\beta (I^{(p)})\rightarrow \beta (Y),\;\;\;\; \mbox{et} \; \forall m,\;\;
\beta _m (I^{(p)})\rightarrow \beta _m(Y).
$$
On note l'\'egalit\'e
$$
\beta _p(I^{(p)})=\beta (I^{(p)}).
$$
On a la formule
$$
\beta _{m+1}(Y)= \beta _{m}(Y) \cup ^{\coprod _{(m+1; y,f)}\beta _m(I^{(m+1)})}
\coprod _{(m+1; y,f)}\beta (I^{(m+1)}).
$$
Le coproduit est pris sur toutes les suites non-d\'eg\'en\'er\'ees de longueur
$m+1$, not\'ees $(m+1; y,f)$.

Cette formule est vraie aussi bien avec le coproduit de cat\'egories qu'avec
le coproduit de leurs nerfs (on pourrait \'ecrire la deuxi\`eme version:
$$
\nu \beta _{m+1}(Y)= \nu \beta _{m}(Y) \cup ^{\coprod _{(m+1; y,f)}\nu \beta
_m(I^{(m+1)})}  \coprod _{(m+1; y,f)}\nu \beta (I^{(m+1)}).
$$
pour \^etre plus pr\'ecis).
Pour voir ceci on note (en \'ecrivant explicitement ce qu'est un
$p$-simplexe du nerf de $\beta (Y)= \beta (\nu Y)$, comme  l'a fait ci-dessus
pour $\beta (X)$) qu'un $p$-simplexe dans le nerf de $\beta _{m+1}(Y)$ est ou
bien dans $\nu \beta _{m}(Y)_p$, ou bien dans l'un des
$\nu \beta (I^{(m+1)})$.

D'apr\`es cette deuxi\`eme interpr\'etation,
la formule d\'ecompose $\beta _{m+1}(Y)$
en un coproduit en tant que $1$-cat\'egorie de Segal. Les sous-cat\'egories
$\delta$
s'y int\'egrent de fa\c{c}on compatible (en
ce sens qu'on a la m\^eme d\'ecomposition
pour $\nu \delta _{m+1}(Y)$), et on a
$$
1SeL(\beta _{m+1}(Y), \delta  _{m+1}(Y)) =
1SeL(\beta _{m}(Y),\delta _m(Y))
$$
$$
\cup ^{\coprod _{(m+1; y,f)}1SeL(\beta
_m(I^{(m+1)}),\delta _m(I^{(m+1)}))}  \coprod _{(m+1; y,f)} 1SeL(\beta
(I^{(m+1)}),\delta (I^{(m+1)})).
$$

Si $U\rightarrow Y$ est un morphisme fibrant de $n$-pr\'ecats de Segal, alors
$\psi$ induit une \'equivalence
$$
\psi ^{\ast}: Sect (Y, U)\rightarrow Sect (1SeL(\beta (Y), \delta (Y)), \psi
^{\ast}(U)).
$$
Ceci est d\^u au fait que $\psi $ induit une \'equivalence faible entre
$1SeL(B,D)$
et $Y$ (ce principe est le sujet de \cite{DwyerKanDiags}).  Si $A$ est un
$n$-pr\'echamp de Segal sur $Y$ on obtient l'\'equivalence $$
\psi ^{\ast}: Sect (Y, \int '_YA)\stackrel{\cong}{\rightarrow}
Sect ^{\delta (Y)-eq}(\beta (Y), \int '_{\beta (Y)}\psi ^{\ast}(A)).
$$
Ici $Sect ^{\delta (Y)-eq}(\beta (Y), \int '_{\beta (Y)}\psi ^{\ast}(A))$
est la sous-$n$-cat\'egorie de Segal
pleine des sections qui, restreintes \`a $\delta (Y)$, sont des $eq$-sections.
Pour voir ceci il faudrait utiliser le th\'eor\`eme 2.5.1 de \cite{limits}.

Une section sur $Y$ est une $eq$-section si et seulement si sa restriction sur
$\beta (Y)$ est une $eq$-section. Donc $\psi$ induit une \'equivalence
$$
\psi ^{\ast}: Sect ^{eq}(Y, \int '_YA)\stackrel{\cong}{\rightarrow}
Sect ^{eq}(\beta (Y), \int '_{\beta (Y)}\psi ^{\ast}(A)).
$$

D'autre part $\psi$ induit un {\em isomorphisme}
$$
\psi ^{\ast}: \Gamma (Y^o,A) \stackrel{\cong}{\rightarrow}\Gamma (\beta (Y) ^o,
\psi ^{\ast}A).
$$
Ceci car les morphismes de $\delta (Y)$ agissent comme l'identit\'e dans le
pr\'efaisceau $\psi ^{\ast}(A)$, et donc toute section de
$\psi ^{\ast}(A)$ sur $\beta (Y)^o$ est forc\'ement constante sur $\delta
(Y)^o$.
Les sections de $\psi ^{\ast}(A)$ descendent donc (uniquement) sur la
cat\'egorie
quotient de $\beta (Y)^o$ par $\delta (Y)^o$, qui est $Y^o$.

Malheureusement, le remont\'e $\psi ^{\ast}(A)$ n'est pas en g\'en\'eral
fibrant
sur $\beta (Y)^o$. On contournera ce probl\`eme en utilisant \ref{calclim},
qui donne une \'equivalence naturelle
$$
\Gamma (Y^o, A)\stackrel{\cong}{\rightarrow}
\lim _{\leftarrow , Y} A.
$$
En utilisant la description
$$
\lim _{\leftarrow , Y} A= \underline{Hom}(Y, nSeCAT')_{1/}(\ast , A)
$$
on obtient (\`a cause du lemme \ref{bd} ci-dessus ainsi que le Theorem 2.5.1
de \cite{limits}) que le morphisme
$$
\lim _{\leftarrow , Y} A\rightarrow \lim _{\leftarrow , \beta (Y)} \psi ^{\ast}A
$$
est une \'equivalence. Ceci implique que le morphisme
$$
\Gamma (\beta (Y)^o, \psi ^{\ast} A) \rightarrow
\lim _{\leftarrow , \beta (Y)} \psi ^{\ast}A
$$
est une \'equivalence.

Soit $\psi ^{\ast} (A) \rightarrow F$ un remplacement fibrant au-dessus de
$\beta (Y)$, et regardons le diagramme
$$
\begin{array}{ccc}
\Gamma (\beta (Y)^o, \psi ^{\ast} A) &\rightarrow &
\lim _{\leftarrow , \beta (Y)} \psi ^{\ast}A \\
\downarrow && \downarrow \\
\Gamma (\beta (Y)^o, F)&\rightarrow &
\lim _{\leftarrow , \beta (Y)}F.
\end{array}
$$
On sait, comme ci-dessus, que la fl\`eche du bas est une \'equivalence.
D'autre part, comme le morphisme $\psi ^{\ast} (A) \rightarrow F$
est une \'equivalence dans $\underline{Hom}(Y, nSeCAT')$, on obtient que
la fl\`eche verticale \`a droite est une \'equivalence. Enfin, on a vu
pr\'ec\'ed\'emment que la fl\`eche du haut est une \'equivalence. Donc
$$
\Gamma (\beta (Y)^o, \psi ^{\ast} A)
\stackrel{\cong}{\rightarrow}
\Gamma (\beta (Y)^o, F).
$$
Le fait que $\psi ^{\ast} (A) \rightarrow F$ soit une \'equivalence
objet-par-objet au-dessus de $\beta (Y)$ implique que le morphisme induit
$$
\int '_{\beta (Y)}\psi ^{\ast}(A) \rightarrow
\int '_{\beta (Y)}F
$$
est une \'equivalence, et donc que
$$
Sect ^{eq}(\beta (Y), \int '_{\beta (Y)}\psi ^{\ast}(A) )
\rightarrow
Sect ^{eq}(\beta (Y), \int '_{\beta (Y)}F )
$$
est une \'equivalence. En somme, il suffit de prouver que le morphisme
$$
\Gamma (\beta (Y)^o, F)\rightarrow Sect ^{eq}(\beta (Y), \int '_{\beta (Y)}F )
$$
est une \'equivalence.

Au vu de cette r\'eduction, on pourra dire qu'on s'est ramen\'e \`a prouver la
proposition pour des cat\'egories de la forme $\beta (Y)$. Le m\^eme argument
tient avec $\beta^{\rm poset}(Y)$ \`a la place de $\beta (Y)$ (en utilisant le
lemme \ref{bdvariant} \`a la place de \ref{bd}). Ceci veut dire qu'on peut
dor\'enavant supposer que $Y$ est d\'ej\`a de la forme  $Y= \beta^{\rm
poset}(Z)$
pour une cat\'egorie $Z$, i.e. que $Y$ est la cat\'egorie associ\'ee \`a un
ensemble partiellement ordonn\'e (cette id\'ee d'utiliser la subdivision
barycentrique deux fois provient de Thomason \cite{ThomasonCatMod}).

Maintenant, pour $F$ on peut utiliser la d\'ecomposition de $\beta _{m+1}(Y)$ en
coproduit: on a
$$
\Gamma (\beta _{m+1}(Y)^o, F|_{\beta _{m+1}(Y)})=
$$
$$
\Gamma (\beta _m(Y)^o, F|_{\beta _{m}(Y)^o}) \times
_{\prod _{(m+1; y,f)}\Gamma (\beta _m(I^{(m+1)})^o, F|_{\beta
_{m}(I^{(m+1)})^o})}
\prod _{(m+1; y,f)}\Gamma (\beta (I^{(m+1)})^o, F|_{\beta (I^{(m+1)})^o}).
$$

Soit $i_m: \beta_m (Y)^o\rightarrow \beta (Y)^o$ l'inclusion. On note que
$i_m$ est un crible, i.e.  que si
$a\rightarrow b$ est un morphisme dans $\beta (Y)^o$ (i.e. la suite
non-d\'eg\'en\'er\'ee qui correspond \`a $a$ est incluse dans celle qui
correspond \`a $b$) et si $b\in \beta _m(Y)^o$ alors $a\in \beta _m(Y)$. On
obtient, d'apr\`es le corollaire \ref{crible}, que $i_m^{\ast}$ est un
foncteur de
Quillen \`a droite, i.e. pr\'eserve les fibrations. Donc $F|_{\beta
_m(Y)^o}$ est
fibrant. Aussi, on peut \'ecrire
$$
\Gamma (\beta (Y)^o, F) =
\underline{Hom}(\ast_{\beta (Y)^o}, F),
$$
tandis que
$$
\Gamma (\beta _m(Y)^o, F|_{\beta _m(Y)^o}) = \underline{Hom}(i_{m,!}\ast_{\beta
_m(Y)^o}, F|_{\beta (Y)^o})
$$
(ceci au vu de la description explicite de $i_{m,!}$ cf la remarque
apr\`es \ref{crible}). Le morphisme
$$
\Gamma (\beta (Y)^o, F)
\rightarrow
\Gamma (\beta _m(Y)^o, F|_{\beta _m(Y)^o})
$$
est fibrant, car il
provient de la cofibration
$$
i_{m,!}\ast_{\beta
_m(Y)^o}\hookrightarrow \ast _{\beta (Y)^o}
$$
(le fait que l'avant-dernier morphisme est fibrant sera utilis\'e
directement sur $Y$ mais aussi sur les $I^{(m+1)}$).

On applique maintenant \ref{shriek} \`a un morphisme
$$
p: \beta (I^{(m+1)})^o\rightarrow \beta (Y)^o
$$
provenant d'une suite non-d\'eg\'en\'er\'ee $(m+1; y,f)$ dans $Y$. C'est ici
qu'on utilise l'hy\-po\-th\`e\-se que $Y$ est la cat\'egorie associ\'ee \`a un
ensemble $S$ partiellement ordonn\'e: dans ce cas, les objets de $\beta
(Y)$ sont les suites d\'ecroissantes
pour l'ordre dans $S$. Donc, pour le morphisme
$p$ associ\'e \`a une suite non-d\'eg\'en\'er\'ee $(m+1; y,f)$ dans $Y$ (i.e.
une suite strictement d\'ecroissante de $S$), on a que la cat\'egorie des objets
$b\in \beta (I^{(m+1)})$ au-dessous d'un objet  $a\in \beta (Y)$ est, ou bien
vide (si $a$ ne correspond pas \`a une sous-suite d'objets de la suite $(m+1;
y,f)$), ou bien admet un objet initial (si $a$ correspond \`a une sous-suite
d'objets de la suite $(m+1; y,f)$, l'objet initial \'etant la sous-suite
correspondante d'objets de $I^{(m+1)}$). Donc dans les deux cas, la
colimite de
\ref{shriek} pr\'eserve les cofibrations et cofibrations triviales, i.e. $p_!$
est un foncteur de Quillen \`a gauche et donc $p^{\ast}$ pr\'eserve les objets
fibrants.
En particulier,
$$
F|_{\beta (I^{(m+1)})^o}
$$
est une $n$-pr\'echamp de Segal fibrant au-dessus de $\beta (I^{(m+1)})^o$.
Par les
arguments ant\'erieurs, on obtient que
$$
F|_{\beta _m(I^{(m+1)})^o}
$$
est aussi fibrant.

Les morphismes
$$
\Gamma (\beta (I^{(m+1)})^o,
F|_{\beta (I^{(m+1)})^o}) \rightarrow
\Gamma (\beta _m(I^{(m+1)})^o, F|_{\beta _m(I^{(m+1)})^o})
$$
sont fibrants (d'apr\`es une des remarques ci-dessus,
qu'on avait apliqu\'eee \`a $Y$, et
qu'on applique ici \`a $I^{(m+1)}$). Donc la d\'ecomposition de $\Gamma (\beta
_{m+1}(Y)^o, F|_{\beta _{m+1}(Y)^o})$ en produit fibr\'e est
aussi un produit fibr\'e
homotopique.

On a exactement la m\^eme d\'ecomposition (qu'on ne re\'ecrit pas) pour
$$
Sect ^{eq}(\beta
_{m+1}(Y), \int '_{\beta _{m+1}(Y)} F|_{\beta _{m+1}(Y)} ),
$$
et cette d\'ecomposition est aussi un
produit fibr\'e homotopique. En effet, la d\'ecomposition de
$\beta _{m+1}(Y)$ en coproduit est une d\'ecomposition en coproduit de
$1$-pr\'ecats de Segal, d'o\`u la d\'ecomposition pour les sections; et
le morphisme sur les sections qui correspond \`a $\beta
_m(I^{(m+1)})\hookrightarrow \beta (I^{(m+1)})$, est une fibration d'o\`u
il d\'ecoule que
le produit fibr\'e dans la d\'ecomposition est un produit fibr\'e homotopique.

On remarque maintenant qu'on a
$$
\Gamma (\beta (Y)^o, F)
= \lim _{\leftarrow , m}
\Gamma (\beta _m(Y)^o, F|_{\beta _m(Y)^o}),
$$
et cette limite est une limite dans laquelle tous les morphismes sont fibrants;
c'est donc aussi une limite homotopique. De plus, on a
$$
Sect ^{eq}(\beta (Y), \int '_{\beta (Y)} F)=
\lim _{\leftarrow , m}
Sect ^{eq}(\beta
_{m}(Y), \int '_{\beta _{m}(Y)} F|_{\beta _{m}(Y)} ),
$$
et cette limite est encore une limite o\`u les morphismes sont fibrants, donc
c'est une limite homotopique.

A cause de ces derni\`eres remarques sur les limites,
il suffit de prouver que les morphismes
$$
\Gamma (\beta _m(Y)^o, F|_{\beta _m(Y)^o})\rightarrow
Sect ^{eq}(\beta _m(Y), \int
'_{\beta _m(Y)}F |_{\beta _m(Y)})
$$
sont des \'equivalences. On va le faire par r\'ecurrence sur $m$.
Nous pouvons supposer que c'est connu pour $m$, pour toutes
les donn\'ees initiales $(Y,A)$ et
pour tout choix de $F$. On peut en particulier appliquer
l'hypoth\`ese de r\'ecurrence pour obtenir que le morphisme
$$
\Gamma (\beta _m(I^{(m+1)})^o, F|_{\beta _m(I^{(m+1)})^o})
\rightarrow Sect ^{eq}(\beta _m(I^{(m+1)}), \int
'_{\beta _m(I^{(m+1)})}F |_{\beta _m(I^{(m+1)})})
$$
est une \'equivalence
(car, en r\'ealit\'e, $F|_{\beta (I^{(m+1)})^o}$ aurait pu
aussi bien provenir de la couple
$(I^{(m+1)}, A|_{I^{(m+1), o}})$). Au vu de nos d\'ecompositions des deux
cot\'es  du morphisme
$$
\Gamma (\beta _{m+1}(Y)^o, F|_{\beta _{m+1}(Y)^o})\rightarrow
Sect ^{eq}(\beta _{m+1}(Y),
\int '_{\beta _{m+1}(Y)}F |_{\beta _{m+1}(Y)})
$$
en produits fibr\'es homotopiques, et au vu de l'hypoth\`ese de r\'ecurrence, il
suffit, pour obtenir que c'est une \'equivalence, de prouver que le
morphisme
$$
\Gamma (\beta (I^{(m+1)})^o, F|_{\beta (I^{(m+1)})^o})\rightarrow
Sect ^{eq}(\beta (I^{(m+1)}),
\int '_{\beta (I^{(m+1)})}F |_{\beta (I^{(m+1)})})
$$
est une \'equivalence.

En faisant marche-arri\`ere \`a travers notre argument
ci-dessus, on voit qu'il suffit de prouver que le morphisme
$$
\Gamma (I^{(m+1),o}, A|_{I^{(m+1),o}})\rightarrow
Sect ^{eq}(I^{(m+1)},
\int '_{I^{(m+1)}}A |_{I^{(m+1)}})
$$
est une \'equivalence.
En somme, nous nous sommes maintenant r\'eduits \`a prouver la proposition pour
le cas $Y= I^{(m+1)}$.

Pour simplifier les notations, on pose $p:= m+1$ et
on va
prouver la proposition dans le cas $Y:= I^{(p)}$. Pour deux objets
$i$ et $j$ de $Y$, on a
une fl\`eche $i\rightarrow j$ pour $i\leq
j$. En particulier, $0$ est l'objet initial de $Y$. Si $A$ est un $n$-pr\'echamp
de Segal sur $Y^o$, on a
$$
\Gamma (Y, A^o) =A(0).
$$
Il s'agit donc (au vu du lemme \ref{inchange}) de prouver que le morphisme
$$
Sect ^{eq}(Y, \int '_YA)\rightarrow A'(0)
$$
est une \'equivalence (ici $A'(0)$ est la fibre de $\int '_YA$ au-dessus de
$0\in Y$, qui est \'equivalente \`a $A(0)$). En \'ecrivant l'\'equivalence
faible de $1$-pr\'ecats de Segal
$$
I\cup ^{\{ 1\}} I \cup ^{\{ 2\} } \cup \ldots \cup ^{\{ (p-1)\} } I
\rightarrow Y
$$
on voit qu'il suffit de traiter le cas $Y=I$, i.e. de prouver que
$$
Sect ^{eq} (I, \int '_IA)\rightarrow A'(0)
$$
est une \'equivalence.
Or ceci est le r\'esultat du lemme \ref{casi}, ce qui termine
la d\'emonstration de la proposition
\ref{gammasect}. \eop

Cet \'enonc\'e, combin\'e avec la proposition \ref{calclim}
indique
que si $A$ est un foncteur strict d'une $1$-cat\'egorie $Y$ vers $nSeCat$
on peut calculer sa limite \`a l'aide des sections de l'int\'egrale:
$$
\lim _{\leftarrow , Y}A \cong Sect ^{\rm eq}(Y, \int '_YA).
$$
Il y a un r\'esultat ``dual'' (mais qui ne s'en d\'eduit pas par dualit\'e)
qui permet le calcul des colimites: c'est le th\'eor\`eme de Thomason
\cite{ThomasonLimits}
dans le cas des pr\'efaisceaux simpliciaux ($n=0$). Pour $n$ quelconque,
nous le formulons comme une conjecture.

\begin{conjecture}
Soit $A$ un pr\'efaisceau de $n$-cat\'egories de Segal au-dessus de $Y^o$.
Soit $horiz$ l'ensemble des fl\`eches ``horizontales'' i.e. de $y$ \`a
$res_f(y)$
correspondant \`a l'identit\'e de $res_f(y)$, dans $\int _YA$. Alors
on a
$$
\lim _{\rightarrow , Y} A \cong nSeL(\int _YA, horiz )
$$
i.e. la colimite de $A$ sur $Y$ se calcule en localisant $\int _YA$ en les
$1$-fl\`eches horizontales.
\end{conjecture}

On termine cette section par la formulation d'un probl\`eme. Dans
SGA1, la th\'eorie des champs et de la descente \'etait envisag\'ee dans le
cadre
des {\em cat\'egories fibr\'ees} \cite{SGA1} \cite{GiraudThese}, mais nous
n'avons repris ci-dessus que la moiti\'e facile de cela, i.e. la construction
d'une ``$n$-cat\'egorie de Segal fibr\'ee'' \`a partir d'un pr\'efaisceau strict
de $n$-cat\'egories de Segal.

\begin{probleme}
\label{problemefibree}
D\'efinir une notion de {\em $n$-cat\'egorie de Segal fibr\'ee}
$F\rightarrow \Yy$ au-dessus d'une $1$-cat\'egorie $\Yy$ (ou m\^eme ---plus
difficile--- au-dessus d'une autre $n$-cat\'egorie de Segal), et \'etablir
la correspondance avec nos $n$-champs de Segal. Par exemple, d\'efinir la
$n+1$-cat\'egorie de Segal $nSeCATFIB /\Yy $ des $n$-cat\'egories de
Segal fibr\'ees au-dessus de $Y$ et donner une \'equivalence $nSeCATFIB
/\Yy \cong nSeCHAMP (\Yy ^{\rm gro})$.

Cette \'equivalence devrait passer par
une construction $A\mapsto \int _{\Yy} A$ pour tout morphisme faible $A$ i.e.
morphisme $A: \Yy ^o\rightarrow nSeCAT'$.
\end{probleme}

\numero{Pr\'efaisceaux de Quillen}

\label{quipage}

C'est en \'etudiant l'exemple du champ des modules des complexes (cf
\S 21 ci-dessous) que nous avons ressenti la n\'ecessit\'e de consid\'erer des
``pr\'efaisceaux
de cmf''.
Un point de vue qui se r\'epand actuellement est que les ``bons''
morphismes entre
cat\'egories de mod\`eles ferm\'ees, sont les {\em foncteurs de Quillen \`a
gauche} (notion qu'on a rappel\'ee au \S 4)---cf par exemple le livre de
Dwyer-Hirschhorn-Kan \cite{DHK}, ou \'egalement
\cite{HoveyBook} o\`u Hovey d\'efinit la {\em $2$-cat\'egorie des
cat\'egories de mod\`eles ferm\'ees} dont les objets sont les cat\'egories de
mod\`eles ferm\'ees et les morphismes sont les foncteurs de Quillen \`a
gauche.
On est donc conduit naturellement vers la notion suivante.

Soit $\Xx$ une cat\'egorie. Un {\em pr\'efaisceau de
Quillen \`a gauche} (resp. {\em a droite}) sur $\Xx$ est un pr\'efaisceau
de cat\'egories
${\bf M}$, o\`u chaque ${\bf M}_X$ est muni d'une structure de cat\'egorie de
mod\`eles ferm\'ee, telle que les foncteurs de restriction soient des
foncteurs de
Quillen \`a gauche (resp. \`a droite).
\footnote{
Il serait peut-etre pr\'ef\'erable de parler de {\em
$1$-champ} de Quillen car les valeurs sont des $1$-cat\'egories et on pourrait
envisager un probl\`eme de coh\'erence. En fait, puisque la notion de
cat\'egorie
de mod\`eles ferm\'ee est stable par \'equivalence de la cat\'egorie
sous-jacente, on peut toujours strictifier et ne parler que de pr\'efaisceaux de
cat\'egories.  Il pourrait aussi y avoir un probl\`eme de coh\'erence pour le
choix des adjoints des foncteurs de restriction. Nous ignorons ces probl\`emes,
consid\'erant que s'ils existent, ils n'ont d'inter\^et ni g\'eom\'etrique ni
topologique (de m\^eme que nous ignorons les subtilit\'es des questions de
cardinalit\'e).}
On note que les foncteurs de restriction ne pr\'eservent pas forcement
toutes les
structures (elles pr\'eservent les cofibrations et cofibrations triviales,
mais pas obligatoirement
les fibrations ni les \'equivalences faibles qui ne sont pas des cofibrations).

Si on fait un rapprochement entre les notions de topos et de cmf, les
notion de pr\'efaisceau de
Quillen \`a gauche ou \`a droite sont analogues aux notions de
``cat\'egorie bifibr\'ee en topos''
et ``cat\'egorie bifibr\'ee en duaux de topos'' de SGA 4 (\cite{SGA4b}
expos\'e Vbis, d\'ef. 1.2.1,
1.2.2).

Si ${\bf M}$ est un pr\'efaisceau de Quillen
\`a gauche, soit ${\bf M}^{\ast}$ le pr\'efaisceau qui \`a chaque $x\in \Xx$
associe la cat\'egorie de mod\`eles ferm\'ee ${\bf M}^{\ast}(x)$ duale
de ${\bf M} (x)$. Alors ${\bf M}^{\ast}$ est un pr\'efaisceau de Quillen \`a
droite. Ceci nous permet de ne nous int\'eresser qu'aux pr\'efaisceaux de
Quillen \`a gauche. Un peu plus int\'eressant: soit ${\bf M}^{\dag}$ le
pr\'efaisceau sur la cat\'egorie oppos\'ee $\Xx ^o$, avec les m\^emes
cat\'egories de mod\`eles ferm\'ees ${\bf M}(x)$ comme valeurs, mais obtenu en
utilisant les adjoints \`a droite comme morphismes de restriction.  Alors ${\bf
M}^{\dag}$ est aussi un pr\'efaisceau de Quillen (\`a droite).

Soit ${\bf M}$ un pr\'efaisceau de Quillen au-dessus de
$\Xx$. Les objets cofibrants d\'efinissent
un sous-pr\'efaisceau de cat\'egories ${\bf M}_c$
et son sous-pr\'efaisceau $W{\bf M}_c$ o\`u les morphismes sont les
\'equivalences faibles entre objets cofibrants (qui sont bien
pr\'eserv\'ees par
les restrictions). On d\'efinit un
pr\'efaisceau de cat\'egories simpliciales $L({\bf M})$ au-dessus de
$\Xx$
en prenant pour
$
L({\bf M})(X)$
la cat\'egorie simpliciale localis\'ee de Dwyer-Kan
$$
L(\mm )(X):= L({\bf M}_c(X), W{\bf M}_c(X)).
$$
Comme d'habitude, on consid\`ere ce
pr\'efaisceau  comme un $1$-pr\'echamp de Segal sur $\Xx$.

Pour des raisons de simplification des notations ult\'erieures, il est plus
commode de prendre le
point de vue ``covariant''. On consid\`ere dor\'enavant une cat\'egorie $Y$ et
un pr\'efaisceau de
Quillen \`a gauche $\mm$ sur $Y^o$, c'est donc une famille de cmf $\mm (y)$ avec
des foncteurs de Quillen \`a gauche de restriction, qui sont covariants.

Soit $\mm$ un pr\'efaisceau de Quillen \`a gauche sur $Y^o$; si $f: y\rightarrow
z$ est un morphisme de $Y$ notons
$$
res_f: \mm (y)\rightarrow \mm (z)
$$
le foncteur de Quillen \`a gauche de restriction; et notons $res_f^{\ast}$
son adjoint \`a droite. On dispose de la
``cat\'egorie fibr\'ee'' (cf la section pr\'ec\'edente)
$$
\int _Y\mm \rightarrow Y
$$
et de la cat\'egorie des sections
$$
Sect (Y,\int _Y\mm ).
$$
Un objet $\sigma$ de $Sect (Y,\int _Y\mm )$ consiste en une collection
d'objets $\sigma (y)\in
\mm (y)$ pour $y\in Y$ avec des morphismes
$$
\sigma (f): res_f\sigma (y)\rightarrow \sigma (z).
$$
qu'on peut interpr\^eter par adjonction:
$$
\sigma (f): \sigma (y)\rightarrow res_f^{\ast}\sigma (z).
$$
Ces morphismes sont soumis \`a une contrainte d'associativit\'e.
Un morphisme $a:\sigma \rightarrow \tau$ dans $Sect (Y,\int _Y\mm )$
est une collection de morphismes
$$
a(y): \sigma (y)\rightarrow \tau (y)\;\;\;\; \mbox{dans}\;\;\; \mm (y)
$$
soumis \`a la condition de naturalit\'e habituelle.
On dira qu'un morphisme $a$ est une {\em \'equivalence faible} si, pour
chaque $y$, $a(y)$ est une  \'equivalence faible dans $\mm (y)$.

On peut observer que si $\mm = \underline{M}$ est le pr\'efaisceau
constant
prenant pour
valeur la cat\'egorie de mod\`eles ferm\'ee fixe $M$ (qui est bien un
pr\'efaisceau de Quillen \`a
gauche), alors on a
$$
\int _Y\underline{M} = Y\times M
$$
et
$$
Sect (Y,\int _Y\underline{M}) = M^Y.
$$
Rappelons qu'on dispose de plusieurs structures de
cat\'egories de mod\`eles ferm\'ees
sur $M^Y$
sous diverses hypoth\`eses concernant $M$ et $Y$.
\newline
(I)\,\, D'abord il y  a la  structure engendr\'ee par les
cofibrations de Hirschhorn, de type Bousfield-Kan (on l'avait appel\'ee de
``type HBKQ'' au \S 5);
celle-ci existe d\`es que $M$ est
engendr\'ee par les cofibrations (Hirschhorn \cite{Hirschhorn}). Les
fibrations sont
les fibrations objet-par-objet.
\newline
(II)\,\, Ensuite il y a la structure que Dwyer-Hirschhorn-Kan qualifient
de ``unusual'' mais qui est la structure habituellement utilis\'ee en
$K$-th\'eorie, et apparemment due \`a K. Brown, Joyal et Jardine dans le cas de
diagrammes d'ensembles simpliciaux (Hirschhorn l'appelle
``structure de
Heller''); ici les cofibrations sont les cofibrations objet-par-objet.
L'existence de cette structure semble co\"{\i}ncider \`a peu pr\`es avec celle
de la structure (I) de type HBKQ, mais nous n'avons pas trouv\'e
dans la litt\'erature un crit\`ere g\'en\'eral pour cette existence (cf
cependant le livre de
Goerss-Jardine \cite{JardineGoerssBook}).
\newline
(III)\,\, Enfin, si $Y$ est une cat\'egorie de Reedy, il y a une structure
de mod\`eles sur $M^Y$ sans hypoth\`ese suppl\'ementaire sur $M$ (tout
au plus, si la fonction
degr\'e est index\'ee par un ordinal plus grand que $\omega$, il faut
supposer que  les limites de cette taille existent dans $M$). Cette structure
est en g\'en\'eral diff\'erente des deux autres (pour la d\'efinition,
voir la d\'emonstration du th\'eor\`eme ci-dessous). Pour cette {\em structure
de Reedy} voir \cite{Reedy}, \cite{Hirschhorn}, \cite{DHK},
\cite{JardineGoerssBook}. La partie essentielle de la notion de fibration
de Reedy appara\^{\i}t aussi
dans la notion d'hyper-recouvrement expos\'ee par Verdier dans SGA 4
(\cite{SGA4b}, expos\'e V, \S 7, d\'efinition 7.3.1.1 (H3)).

Si $\mm$ est un pr\'efaisceau de Quillen \`a gauche sur $Y^o$,
on a les structures analogues sur $Sect (Y,\int _Y\mm )$:

\begin{theoreme}
\label{sectionscmf}
Soit $\mm$ un pr\'efaisceau de Quillen \`a gauche sur $Y^o$. Alors on
dispose des
structures suivantes de cat\'egorie de mod\`eles ferm\'ee sur $Sect (Y,\int
_Y\mm )$ (les \'equivalences faibles sont toujours les morphismes qui sont
objet-par-objet des \'equivalences faibles):
\newline
(I)\,\, La structure de type HBKQ o\`u les fibrations sont
les morphisms $a:\sigma \rightarrow \tau $ tel que $a(y)$ est une fibration
pour tout $y\in Y$: elle existe si chaque $\mm (y)$ est engendr\'e par
cofibrations;
\newline
(II)\,\, La structure de type Brown-Jardine-Heller o\`u les cofibrations sont
les morphismes $a:\sigma \rightarrow \tau$ tel que $a(y)$ soit une cofibration
pour tout $y\in Y$,
\`a condition que cette classe de cofibrations (resp. les cofibrations
$\cap$ les \'equivalences
faibles) admette un ensemble g\'en\'erateur $I$ (resp. $J$) permettant
l'argument du petit objet
(voir avant \ref{dhklemme}),
et qu'en outre chaque $\mm (y)$ soit engendr\'e par cofibrations; et
\newline
(III)\,\, La structure de type Reedy (cf la d\'efinition
ci-dessous),  si
$Y$ est une cat\'egorie de Reedy.
\end{theoreme}
Comparer avec
la proposition 1.2.12 de SGA 4 \cite{SGA4b} expos\'e Vbis.

\noindent
{\em Preuve:}
Il peut \^etre utile de revoir la preuve de l'\'enonce \ref{dhklemme}
ci-dessus, mais pour (I)
nous conseillons au lecteur de se reporter directement aux r\'ef\'erences
\cite{Hirschhorn} \cite{DHK} au lieu de s'en tenir \`a \ref{dhklemme}, car nous
utiliserons ici directement les arguments du type \cite{DHK}.

On note d'abord que l'ensemble $W$ des \'equivalences faibles satisfait
automatiquement la clot\^ure sous r\'etractes et la propri\'et\'e ``trois
pour le
prix de deux''.

Dans (I) et (II) la condition que chaque $\mm (y)$ soit engendr\'e par
cofibrations (et
donc admette des petites limites et colimites) implique que
$Sect (Y,\int
_Y\mm )$ admet des petites limites et colimites.

Pour (I) la preuve est la m\^eme que dans \cite{Hirschhorn},
\cite{DHK}.
Soit $\sigma$ une section et $i: \sigma (y)\rightarrow u$ une cofibration dans
$\mm (y)$. On obtient alors la {\em cofibration librement engendr\'ee par $i$},
$Lib(i): \sigma \rightarrow \nu$, d\'efinie par
$$
\nu (z) = \coprod _{f: y\rightarrow z} res_f(u) \sqcup ^{\coprod_{f}
res_f(\sigma (y))}\sigma (z).
$$
Soient  $I$ (resp. $J$) l'ensemble des cofibrations
de la forme $Lib(i)$ o\`u $i\in I(y)$, $y\in Y$ (resp. $i\in J(y)$, $y\in
Y$).
Les $I$-injectifs (resp.
$J$-injectifs) sont juste les morphismes qui sont objet-par-objet des fibrations
(resp. fibrations triviales); et on v\'erifie facilement que le crit\`ere de
reconnaissance  (``Recognition lemma'' \cite{Hirschhorn} 13.3.1, \cite{DHK} 8.1)
est satisfait. \footnote{Essentiellement la
seule chose \`a v\'erifier---comme c'est indiqu\'e dans
(\cite{DHK}, preuve de 48.7: la structure de cat\'egorie de mod\`eles sur
$S-Cat$, p. 61), voir aussi notre lemme \ref{dhklemme}---est que la colimite de
$J$-cofibrations r\'egulie\`eres est encore dans $W$.  Dans le pr\'esent cas
ceci est
vrai ponctuellement au-dessus de chaque objet de $Y$.}

l'ensemble
contiennent les
ponctuellement une
reguliere,
$J$ est

Pour (II) on s'appuie sur le lemme \ref{dhklemme}. Rappelons que la
condition que chaque $\mm (y)$
est engendr\'e par cofibrations nous donne les propri\'et\'es de
\ref{dhklemme} (0)-(7) pour chaque
$\mm (y)$.  Comme  cons\'equence on obtient imm\'ediatement les
propri\'et\'es (0), (1), (2),
(6) et (7) pour $Sect (Y,\int
_Y\mm )$.  On peut remarquer que les cofibrations de la structure (I)
ci-dessus sont des
cofibrations ici aussi, donc la propri\'et\'e de rel\`evement pour ces
cofibrations  implique qu'un
morphisme est objet-par-objet une \'equivalence faible.  Ceci donne (3).
Enfin (4) et (5) forment
la premi\`ere partie de l'hypoth\`ese ici, donc par le lemme \ref{dhklemme}
on obtient (II).

{\em Remarque:} Tout comme pour l'application du lemme \ref{dhklemme}, on
peut dire que la partie
de l'hypoth\`ese (II) qui correspond \`a \ref{dhklemme} (4) et (5) serait
cons\'equence du
caract\`ere ensemblistement raisonnable des notions de cofibration et
d'\'equivalence faible,
dans tous les exemples qu'on va rencontrer. Voir Jardine \cite{Jardine} pour la
technique n\'ecessaire pour
donner une preuve rigoureuse de ce type de condition. Au vu de cette
remarque, on pourra
consid\'erer que l'hypoth\`ese de (II) se r\'eduit \`a demander (comme pour
(I)) simplement que les
$\mm (y)$ soient engendr\'es par cofibrations.

\subnumero{Les cat\'egories de Reedy}

Avant d'aborder la partie (III) du th\'eor\`eme \ref{sectionscmf},
on fait quelques rappels sur les ``cat\'egories de Reedy''. Les
r\'ef\'erences sont
\cite{Reedy}, \cite{BousfieldKan}, \cite{Hirschhorn}, \cite{DHK},
\cite{JardineGoerssBook}, voir aussi \cite{DwyerKan3}.

Le jeu entre les ``latch'' et les ``match'' qui suivra, n'est qu'une fa\c{c}on
syst\'ematique et avec une notation abstraite de prendre en compte le fait
qu'il y a des morphismes de deg\'en\'er\'escence dans $\Delta$ qui vont dans le
sens oppos\'e des morphismes ``face''. On peut probablement, dans toutes les
applications, consid\'erer au lieu d'objets simpliciaux, des objets
param\'etris\'ees par la sous-cat\'egorie $\direct{\Delta}$ engendr\'ee par les
morphismes face. Dans ce cas on n'aurait qu'\`a consid\'erer la partie
``Latch'' de la discussion ci-dessous et beaucoup d'arguments deviendraient
plus facile.

On rappelle d'abord la d\'efinition: une {\em cat\'egorie de Reedy}
est une cat\'egorie $Y$ munie de deux sous-cat\'egories (avec les m\^emes objets
que $Y$) appel\'ees respectivement {\em sous-cat\'egorie directe} et {\em
sous-cat\'egorie inverse}
$$
\direct{Y}\subset Y \;\;\; \mbox{et} \;\;\; \inverse{Y}\subset Y,
$$
tel qu'il existe une fonction {\em degr\'e} (vers un ordinal qui sera le plus
souvent $\omega$)  telle que les morphismes de $\direct{Y}$ sauf les
identit\'es augmentent strictement le degr\'e, les morphismes de $\inverse{Y}$
sauf les identit\'es diminuent strictement le degr\'e, et tout morphisme $f$ se
factorise uniquement
$$
f= \direct{f}\inverse{f} \;\;\; \mbox{avec} \;\;\; \direct{f}\in \direct{Y},
\;\; \inverse{f}\in \inverse{Y}.
$$
On pr\'ecise ici qu'on \'ecrit les compositions dans le sens habituel, donc
dans cette formule la source de $\direct{f}$ (qui est le but de
$\inverse{f}$) est de degr\'e $\leq$ le minimum des degr\'es de la source et du
but de $f$; et en cas d'\'egalit\'e $f$ est ou bien directe ou bien inverse.

Les exemples classiques sont la cat\'egorie $\Delta$ et son oppos\'e $\Delta
^o$, \cite{BousfieldKan}, \cite{Reedy}.
L'exemple que nous utiliserons est la cat\'egorie $\Delta K$ des simplexes d'un
ensemble simplicial $K$ (ou son oppos\'e $\Delta ^oK$). Cet exemple
appara\^{\i}t
dans \cite{Hirschhorn}, \cite{DHK}, etc. Plus particuli\`erement on va utiliser
le cas de $\Delta \nu \Yy$ o\`u $\nu \Zz$ est le nerf d'une
$1$-cat\'egorie $\Zz$. C'est la cat\'egorie not\'ee $\beta (\Zz )$ au \S
16 ci-dessus.

Pour un objet $y\in Y$ on d\'efinit les cat\'egories
{\em rel\^achantes} (``latching'')
\footnote{
Le mot ``rel\^achant'' n'est pas une traduction stricte de
``latching'' mais comporte la m\^eme sonorit\'e et semble
bien correspondre \`a la situation.}
$Latch(y)$
et {\em appariantes} (``matching'') $Match(y)$.
En fait on d\'efinit
$$
Latch(y)+\{ y\} := \direct{Y}/y
$$
(les objets de $\direct{Y}$ au-dessus de $y$)
et on pose
$$
Latch(y):= \direct{Y}/y -\{ y\} \subset Latch(y)+\{ y\} .
$$
Aussi
$$
\{ y\} + Match(y):= y/\inverse{Y}
$$
(les objets de $\inverse{Y}$ au-dessous de $y$)
et
$$
Match(y):= y/\inverse{Y} -\{ y\} \subset \{ y\} + Match(y).
$$
On note
$$
Latch(y)+ \{ y\} + Match(y)
$$
la cat\'egorie obtenue en attachant $Latch(y)+\{ y\} $ \`a
$\{ y\} + Match(y) $ le long de
l'objet $y$, avec exactement un morphisme entre les objets de $Latch(y)$ et
$Match(y)$. On note
$$
Latch(y)+Match(y)\subset Latch(y)+ \{ y\} + Match(y)
$$
la sous-cat\'egorie obtenue en enlevant $y$.
Le morphisme de $1$-pr\'ecats de Segal
$$
(Latch(y)+\{ y\} ) \cup ^{\{ y\} } (\{ y\} + Match(y))\rightarrow
Latch(y)+ \{ y\} + Match(y)
$$
est une \'equivalence faible (on laisse la preuve, qui est similaire \`a celle
de \ref{reedydecomp} ci-dessous, aux soins du lecteur). En particulier, pour une
$n$-cat\'egorie de Segal fibrante $C'$ un morphisme
$$
Latch(y)+ \{ y\} + Match(y)\rightarrow C'
$$
est essentiellement la m\^eme chose que deux morphismes
$$
Latch(y)+\{ y\} \rightarrow C', \;\;\; \{ y\} + Match(y) \rightarrow C'
$$
qui envoient $y$ sur le m\^eme objet de $C'$.

Pour un degr\'e $k$  fix\'e, on note $F^kY\subset Y$ la sous-cat\'egorie
des objets de degr\'e $\leq k$. Si $y$ est de degr\'e $k$ alors on a un
morphisme
$$
Latch(y)+Match(y) \rightarrow F^{k-1}Y.
$$
On a que $Y$ est la colimite filtrante des $F^kY$.

La description des diagrammes sur une cat\'egorie de Reedy de \cite{DHK},
\cite{JardineGoerssBook}, \cite{BousfieldKan},
\cite{DwyerE2}, \cite{BousfieldFriedlander}, \cite{Reedy}, \cite{DwyerKan3},
et particuli\`erement Hirschhorn (\cite{Hirschhorn} Theorem 16.2.12) o\`u cela
appara\^{\i}t tr\`es clairement, revient essentiellement au lemme suivant.

\begin{lemme}
\label{reedydecomp}
Pour tout $d$, le morphisme de $1$-pr\'ecats de Segal
$$
F^{d-1}Y \cup ^{\coprod _{deg(y)=d}Latch(y)+Match(y)}\coprod_{deg(y)=d}
Latch(y)+ \{ y\} + Match(y)
$$
$$
\rightarrow F^dY
$$
est une \'equivalence faible.
\end{lemme}

Ce lemme
g\'en\'eralise \'eventuellement les r\'ef\'erences ci-dessus en prenant en
compte (via la
notion de $1$-pr\'ecat de Segal) les homotopies sup\'erieures. Nous aurons
besoin
de cet aspect dans les prochains chapitres. Pour la partie
(III) du th\'eor\`eme \ref{sectionscmf}, la version
(\cite{Hirschhorn} Theorem 16.2.12)---qui est l'\'enonc\'e obtenu \`a
partir de notre lemme en
appliquant la troncation $\tau_{\leq 1}$---suffirait. Au vu de cela, nous
repoussons la
d\'emonstration jusqu'\`a la fin du chapitre.

Le lemme \ref{reedydecomp} dit que pour une $n$-cat\'egorie de Segal
fibrante $C'$, se donner un morphisme
$f_k:F^kY\rightarrow C'$ revient \`a se donner un morphisme
$f_{k-1}: F^{k-1}Y\rightarrow C'$ et pour tout objet $y$, des extensions des
morphismes
$$
f_{k-1}|_{Latch(y)} \;\; \mbox{en}\;\; Latch(y)+\{ y\} \rightarrow C',
$$
et
$$
f_{k-1}|_{Match(y)} \;\; \mbox{en}\;\; Match(y)+\{ y\} \rightarrow C',
$$
prenant la m\^eme valeur sur $y$,
et telles que le morphisme induit (essentiellement bien d\'efini)
$$
Latch(y)+ \{ y\} + Match(y)\rightarrow C'
$$
ait $f_{k-1}$ pour restriction \`a $Latch(y)+Match(y)$.
Cette interpr\'etation provient de Hirsch\-horn \cite{Hirschhorn}.

\subnumero{La structure (III) du th\'eor\`eme \ref{sectionscmf}}

On doit dire quelles sont les fibrations et cofibrations.
Soit $\sigma$ une section et $y\in Y$.
Soient $latch(\sigma , y)$ et $match(\sigma , y)$ dans $\mm (y)$ les {\em objets
rel\^achants} et {\em appariants} de $\sigma$. Ces objets sont d\'efinis
respectivement
comme la colimite sur $Latch(y)$ du foncteur
$$
res_{Latch(y)\rightarrow y} \circ \sigma : Latch(y)\rightarrow \mm (y)
$$
o\`u
$$
res_{Latch(y)\rightarrow y} : \int _{Latch(y)}\mm |_{Latch(y)} \rightarrow \mm
(y)
$$
est le morphisme donn\'e par les restrictions $res_f$ pour $f: z\rightarrow y$;
et la limite sur $Match(y)$ du foncteur
$$
res_{y\rightarrow Match(y)} ^{\ast}\circ \sigma : Match(y)\rightarrow \mm
(y)
$$
o\`u
$$
res_{y\rightarrow Match(y)}^{\ast} : \int _{Match(y)}\mm |_{Match(y)}
\rightarrow \mm (y)
$$
est le morphisme donn\'e par les adjoints des restrictions
$res_g^{\ast}$ pour $g: y\rightarrow w$.
Pour d\'efinir $res_{y\rightarrow Match(y)}^{\ast}$  on note qu'on a, pour
$$
y\stackrel{g}{\rightarrow} w \stackrel{h}{\rightarrow} x,
$$
une transformation naturelle
$$
res_g^{\ast} \rightarrow res _{hg}^{\ast} \circ res_h,
$$
qui provient du morphisme d'adjonction
$$
1\rightarrow res_h^{\ast} \circ res_h.
$$
On a une factorisation
$$
latch(\sigma ,y)\rightarrow \sigma (y) \rightarrow match(\sigma , y).
$$
On dira (en suivant \cite{Reedy} \cite{Hirschhorn} \cite{DHK} {\em et al.}) que
$\sigma$ est {\em cofibrant (de Reedy)} si les morphismes  $latch(\sigma
,y)\rightarrow \sigma (y)$  sont des cofibrations dans $\mm (y)$, et que
$\sigma$
est  {\em fibrant (de Reedy)} si
les morphismes $\sigma (y)\rightarrow match(\sigma ,y)$ sont des fibrations dans
$\mm (y)$. Plus g\'en\'eralement pour $a:\sigma \rightarrow \tau$
on obtient l'objet rel\^achant relatif
$$
latch(a,y):= latch(\tau , y) \cup ^{latch(\sigma , y)}\sigma (y)
$$
avec morphisme
$$
latch (a,y)\rightarrow \tau (y);
$$
on dira que $a$ est une {\em cofibration (de Reedy)} si ce morphisme est une
cofibration dans $\mm (y)$. Dualement on obtient l'objet appariant relatif
$$
match (a,y):= match (\sigma ,y)\times _{match(\tau , y)}\tau (y)
$$
avec morphisme
$$
\sigma (y)\rightarrow match(a,y);
$$
et on dira que $a$ est une {\em fibration (de Reedy)} si ce morphisme est une
fibration dans $\mm (y)$.

La preuve que ces trois classes de morphismes forment une structure de
cat\'egorie
de mod\`eles ferm\'ee est la m\^eme que dans le cas du pr\'efaisceau
constant $\mm =\underline{M}$, cf \cite{Reedy} \cite{Hirschhorn} {\em et al.};
nous laissons au lecteur le soin de faire cette v\'erification.

Ceci termine la d\'emonstration de (III) et donc du th\'eor\`eme
\ref{sectionscmf}.
\eop

\subnumero{Preuve du lemme \ref{reedydecomp}}

On garde les notations du lemme \`a d\'emontrer.

On peut supposer que $Y=F^dY$ et on pose $Z:= F^{d-1}Y$. Notons
$A^0$ le coproduit de l'\'enonc\'e. C'est un ensemble simplicial
(sous-ensemble de $\nu Y$), qu'on consid\`ere comme $1$-pr\'ecat de Segal. On va
obtenir $\nu Y$ comme limite d'une suite (\'eventuellement transfinie si le
nombre
d'objets de $Y$ est plus que d\'enombrable)
$$
\ldots A^i\subset A^{i+1}\ldots
$$
de sous-ensembles simpliciaux de $\nu Y$. A chaque fois on rajoutera un
``simplexe r\'egulier'' (d\'efinition ci-dessous), et on va prouver que les
$A^i\subset A^{i+1}$ sont des cofibrations triviales de $1$-pr\'ecats de Segal.

Rappelons que les \'el\'ements de $\nu (Y)_p$ sont les suites
composables
$$
(y,f) =(y_0,\ldots , y_p; f_1,\ldots , f_p)
$$
avec $f_i:
y_{i-1}\rightarrow y_i$ dans $Y$. On dira qu'un tel \'el\'ement est un {\em
simplexe r\'egulier} si aucun des $f_i$ n'est l'identit\'e, et si pour tout
$1\leq i \leq p$, ou bien $f_i \in \direct{Y}$, ou bien $f_i \in \inverse{Y}$.
Pour un simplexe r\'egulier, on appellera {\em pic} tout sommet i.e. objet $y_i$
avec $f_{i-1}\in \direct{Y}$ et $f_i\in \inverse{Y}$, ce qui \'equivaut
\`a la condition $deg(y_{i-1}) < deg(y_i) > deg(y_{i+1}$. On appellera {\em
vall\'ee} tout objet $y_i$ tel que $f_{i-1}\in \inverse{Y}$ et $f_i\in
\direct{Y}$,, ce qui \'equivaut \`a la condition $deg(y_{i-1}) > deg(y_i) <
deg(y_{i+1}$. La condition de Reedy implique que si $y_i$ est une vall\'ee,
alors le simplexe en question est le seul simplexe r\'egulier de longueur $\leq
p$ contenant la face o\`u on enl\`eve $y_i$.

Pour un \'el\'ement $(y,f)\in \nu (Y)_p$ on obtient un morphisme
$(y,f):h(p)\rightarrow Y$ o\`u $h(p)$ est l'ensemble simplicial repr\'esent\'e
par $p$, i.e. le $p$-simplexe standard.
On appellera ``image d'un simplexe $(y,f)$'' l'image
de ce morphisme.
L'image contient
tous les simplexes ``hyper-faces'' de $(y,f)$, qui sont tous les simplexes
obtenus en enlevant certains des $y_i$; ainsi que tous leurs
d\'eg\'en\'er\'escences.

On choisit un ordre total $\prec$ sur l'ensemble des simplexes r\'eguliers de
$\nu Y$ non contenus dans $A^0$, assujetti aux conditions suivantes: si
$$
(y_0,\ldots , y_p; f_1,\ldots , f_p) \prec
(y'_0,\ldots , y'_q; f'_1,\ldots , f'_q)
$$
alors
$$
p\leq q
$$
et, si $p=q$ alors
$$
\sum _{i=0}^p deg (y_i) \leq \sum _{i=0}^p deg(y'_i).
$$
Maintenant, on part de $A^0$ et on ajoute les images de $h(p)$ pour les
simplexes r\'eguliers $(y,f): h(p)\rightarrow \nu (Y)$, un par un, suivant
l'ordre total $\prec$ (aux ordinaux limites on ins\`ere un objet qui est la
colimite de ce qui vient avant, puis on recommence en ajoutant le prochain
simplexe). On obtient une suite $A^i\subset \nu (Y)$. On va montrer qu'une
cofibration $A^i\rightarrow A^{i+1}$ qui correspond \`a l'addition d'un
simplexe r\'egulier $(y_0,\ldots , y_p; f_1,\ldots , f_p)$, est une
\'equivalence faible de $1$-pr\'ecats de Segal. Ceci donnera le lemme
(il y a un argument de colimite filtrante \`a faire pour les ordinaux limites et
\`a la fin du processus, que nous laissons au lecteur).

On va prouver, par r\'ecurrence sur $i$, l'\'enonc\'e suivant:

$(\ast )$ un simplexe d\'eg\'en\'er\'e
est dans $A^i$ si et seulement si le simplexe non-d\'eg\'en\'er\'e
correspondant  est dans $A^i$.

Soit $A^i\rightarrow A^{i+1}$ une cofibration qui correspond \`a
l'addition d'un simplexe r\'egulier $(y_0,\ldots , y_p; f_1,\ldots , f_p)$, et
soit $h(p)\rightarrow \nu (Y)$ le morphisme correspondant \`a ce simplexe.
Si on pose
$$
U:= h(p)\times _{A^i} \nu (Y) \subset h(p),
$$
alors on a
$$
A^{i+1} = h(p)\cup ^U A^i.
$$
Il s'agit donc de prouver que $U\rightarrow h(p)$ est une cofibration triviale
de $1$-pr\'ecats de Segal. On note d'abord que par l'hypoth\`ese $(\ast )$ pour
$A^i$, un \'el\'ement d\'eg\'en\'er\'e de $h(p)_q$ est dans $U$ si et seulement
si l'\'el\'ement non-d\'eg\'en\'er\'e correspondant est dans $U$. Il s'agit donc
de comprendre quelles hyper-faces de $h(p)$ sont dans $U$. On a la
caract\'erisation suivante:

$(\ast \ast )$ \, {\em Une hyper-face de $h(p)$ qui correspond \`a un
sous-ensemble $J\subset \{ 0,\ldots , p\}$, appartient \`a $U$ si et seulement
s'il existe $i\in \{ 0,\ldots , p\}$ avec $i\not \in J$ et  tel que $y_i$ ne
soit pas une
vall\'ee de $(y_0,\ldots , y_p; f_1,\ldots , f_p)$.}

Prouvons ceci. On consid\`ere l'hyper-face correspondant \`a
$$
J= \{j_0,\ldots ,
j_q\} \subset \{ 0,\ldots , p\}
$$
(inclusion stricte). Ceci correspond \`a
l'\'el\'ement
$$
J^{\ast}(y,f):=(y_{j_0}, \ldots , y_{j_q}, f'_1,\ldots , f'_q)\in \nu Y
$$
o\`u
les $f'_j$ sont compositions des $f_i$. Soit $(z,g)=(z_0,\ldots , z_r;
g_1,\ldots , g_r)$ le simplexe r\'egulier obtenu en rempla\c{c}ant dans
$J^{\ast}(y,f)$ toute fl\`eche  $f'_k: y_{j_{k-1}}\rightarrow y_{j_k}$ qui n'est
pas dans $\direct{Y}$ ni dans $\inverse{Y}$, par la paire de fl\`eches
$$
y_{j_{k-1}}\stackrel{g'}{\rightarrow}z\stackrel{g''}{\rightarrow} y_{j_k}
$$
avec $g'\in \inverse{Y}$ et $g''\in \direct{Y}$.
On laisse inchang\'ees les
fl\`eches $f'_k$ qui sont soit dans $\direct{Y}$ soit dans
$\inverse{Y}$.
Aussi on contracte toutes les identit\'es $f'_k$ qui apparaissent dans
$J^{\ast}(y,f)$. On note que $J^{\ast}(y,f)$ est dans l'image de $(z,g)$.

Une fl\`eche $f'_k$ qui est rempla\c{c}\'ee (ou contract\'ee) dans le
proc\'ed\'e
ci-dessus, provient forc\'ement d'une composition d'au moins deux fl\`eches
$f_i$.  Donc la longueur $r$ de $(z,g)$ est $\leq p$. D'autre part, si $r=p$
alors
$$
\sum _{i=0}^p deg (z_i) \leq \sum _{i=0}^p deg (y_i).
$$
En effet, si $f'_k=f_{i+1}f_i = g_{i+1}g_i$ alors $y_i$ est ou bien un pic
ou bien une vall\'ee: dans le cas contraire, $f'_k$ est ou bien dans
$\direct{Y}$
ou bien dans $\inverse{Y}$ et n'est pas rempla\c{c}\'e dans $(z,g)$ (et $f'_k$
n'est pas non plus une identit\'e donc elle n'est pas non plus contract\'ee).
Si $y_i$ est une vall\'ee alors $z_i=y_i$; si $y_i$ est un pic, dans ce cas
$z_i$ est une vall\'ee de $(z,g)$ et $deg(z_i) < deg (y_i)$.

D'apr\`es le paragraphe pr\'ec\'edent, si le compl\'ementaire de $J$ dans $\{
0,\ldots , p\}$ contient un pic $y_i$ alors l'une des deux
in\'egalit\'es sera stricte et (d'apr\`es ce qu'on a impos\'e \`a l'ordre
$\prec$) on a $(z,g)\prec (y,f)$ et donc $J^{\ast}(y,f)\in A^i$. Si le
compl\'ementaire
de $J$ contient un $y_i$ qui n'est ni pic ni vall\'ee, alors comme l'une des
fl\`eches $f'_k$ qui est une composition d'au moins deux fl\`eches $f_i$, n'est
pas rempla\c{c}\'ee par deux fl\`eches $g$, on a $r<p$ et de nouveau
$J^{\ast}(y,f)\in A^i$.
Ceci prouve une des directions de l'\'enonc\'e $(\ast \ast )$
(il faut noter que, dans l'argument ci-dessus, on peut tomber, sans le
dire, sur un
simplexe qui est dans $A^0$ au lieu d'\^etre ajout\'e dans la suite organis\'ee
par $\prec$).

Par ailleurs, s'il y a une contraction d'un $f'_k= 1$ alors la longueur
d\'ecro\*{\i}t: $ r<p$ et $J^{\ast}(y,f)\in A^i$.  On voit que tout nouveau
simplexe
de la forme $J^{\ast}(y,f)$ dans $A^{i+1}$ qui est d\'eg\'en\'er\'e appartient
d\'ej\`a \`a $A^{i}$; et les autres simplexes d\'eg\'en\'er\'es
proviennent des d\'eg\'en\'er\'escences de simplexes de la forme
$J^{\ast}(y,f)$. Ceci prouve $(\ast )$ pour $A^{i+1}$.

Supposons pour compl\'eter la preuve de $(\ast \ast )$ que le compl\'ementaire
de $J$ consiste enti\`erement de vall\'ees de $(y,f)$. Dans ce cas, pour tout
$i\not \in J$, on a $i-1\in J$ et $i+1\in J$, et l'un des morphismes dans
$J^{\ast}(y,f)$ est
$$
f'_k=f_{i+1}f_i: y_{i-1}\rightarrow y_{i+1}.
$$
On note que $f'_k$ n'est pas dans $\direct{Y}$ ni dans $\inverse{Y}$.
Donc, si le simplexe $J^{\ast}(y,f)$ provient de l'image d'un simplexe
r\'egulier
$(w,h)$ alors la fl\`eche $f'_k$ est forc\'ement d\'ecompos\'ee dans $(w,h)$
comme produit d'au moins deux fl\`eches $h$. Il s'ensuit que la longueur de
$(w,h)$
est au moins $p$. En cas d'\'egalit\'e, chacune des fl\`eches $f'_k$ sera
d\'ecompos\'ee en un produit de la forme $h_{i+1} h_i$. Ici il y a deux
options: soit $w_i$ est une vall\'ee, auquel cas on a $w_i=y_i$ par
l'unicit\'e de
la d\'ecomposition dans l'axiome des cat\'egories de Reedy; soit $w_i$ est
un pic, auquel cas on a $deg (w_i)> deg(y_i)$. Donc, finalement, il y a
trois cas:
soit $(w,h)$ est de longueur $>p$; soit $(w,h)$ est de longueur $p$ mais
$$
\sum _{i=0}^p deg(w_i)>  \sum _{i=0}^p deg(y_i);
$$
soit $(w,h)=(y,f)$. Dans tous les cas, $(w,h)$ ne pr\'ec\`ede pas $(y,f)$ pour
l'ordre $\prec$ et donc $J^{\ast}(y,f)$ n'appartient pas \`a $A^i$.

En fait pour cette derni\`ere phrase il faut justifier aussi que
$J^{\ast}(y,f)$ n'appartient pas \`a $A^0$. On note que $(y,f)$ a
au moins un pic de degr\'e $d$ (i.e. dans $Y$ mais pas dans $Z$), et ce pic
persiste dans
$J^{\ast}(y,f)$. Les seuls simplexes de cette forme dans $A^0$ sont
les simplexes r\'eguliers avec un seul pic (i.e. les simplexes ajout\'es dans
le coproduit de l'\'enonc\'e du lemme), or en enlevant des vall\'ees on a rendu
$J^{\ast}(y,f)$ irr\'egulier et on a donc $J^{\ast}(y,f)\not \in A^0$.
Ceci
prouve $(\ast \ast )$.

L'\'enonc\'e $(\ast \ast )$ dit que $U$ est la r\'eunion de toutes les
faces de $h(p)$ obtenues en enlevant des sommets autre que des vall\'ees.
Ceci implique
que le morphisme
$$
U\rightarrow
h(p)
$$
est une \'equivalence faible de $1$-pr\'ecats de Segal, par le lemme suivant
(qu'on isole car l'\'enonc\'e pourrait avoir un int\'er\^et propre). Ce qui
terminera la
d\'emonstration du
lemme \ref{reedydecomp}.
\eop

\begin{lemme}
\label{faces}
Notons $h(p)$ l'ensemble simplicial ``$p$-simplexe standard'', et soit
$f:U\hookrightarrow h(p)$ une inclusion telle que $U$ soit r\'eunion d'un
certain
nombre de faces de $h(p)$. Supposons que $U$ contient la premi\`ere et la
derni\`ere face, et ne contient pas toutes les autres; alors le morphisme
$f$ est une \'equivalence faible de $1$-pr\'ecats de Segal, i.e. le
morphisme induit
$$
SeCat(U)\rightarrow SeCat(h(p))\cong I^{(p)}
$$
est une \'equivalence.
\end{lemme}
{\em Preuve:}
On raisonne par r\'ecurrence sur $p$; supposons donc
l'\'enonc\'e connu pour $p-1$.
Ecrivons
$$
U= V^0\cup V^p \cup V^{b_1}\cup \ldots \cup V^{b_k}
$$
comme r\'eunion (non-disjointe!) de faces o\`u $V^j$ est la face (ensemble
simplicial isomorphe \`a $h(p-1)$) obtenue en enlevant le $j^{\rm ieme}$
sommet. On fixe
$a\in \{ 1,\ldots , p-1\}$ avec $b_i\neq a$, i.e. la face $V^a$
n'appara\^{\i}t pas
dans $U$. Pour $0\leq j\leq k$ posons
$$
U^j:= V^0\cup V^p \cup V^{b_1}\cup \ldots \cup V^{b_j}\subset U.
$$
En particulier $U^0= V^0\cup V^p$, et $U^k=U$. On note d'abord que
$U^0\rightarrow h(p)$ est une \'equivalence (d'ici la fin de la preuve, ceci
voudra dire que c'est une \'equivalence faible de $1$-pr\'ecats de Segal).
Pour voir cela, soit $W(p)\subset h(p)$ la r\'eunion des ar\^etes principales
i.e. des ar\^etes de la forme $\{ i-1, i\}$. Le morphisme
$W(p)\rightarrow h(p)$ est une \'equivalence (en effet, $h(p)$ est obtenu \`a
partir de $W$ par application une fois de l'op\'eration $Gen[p]$ voir \S 2;
voir aussi la discussion dans \cite{limits} \S 2.4.8).
Maintenant $V^0\cap W(p)$ est de la forme $W(p-1)\subset V^0=h(p-1)$.
Donc on peut exprimer $Q(p):=V^0\cup W(p)$ comme coproduit:
$$
Q(p):=V^0\cup W(p) = V^0 \cup ^{W(p-1)}W(p)
$$
o\`u le premier morphisme est une cofibration triviale (de $1$-pr\'ecats de
Segal). Donc le morphisme $W(p)\rightarrow Q(p)$ est une cofibration
triviale et on obtient que $Q(p)\rightarrow h(p)$ est une
\'equivalence (et c'est une cofibration donc c'est une cofibration triviale).
Notons que $Q(p)$ est l'unioni de la premi\`ere face avec la derni\`ere ar\^ete.
Maintenant $Q(p) \cap V^p = Q(p-1)\subset V^p=h(p-1)$. Donc (par le r\'esultat
qu'on vient de d\'emontrer mais pour $Q(p-1)$) le deuxi\`eme  morphisme dans le
coproduit
$$
U^0= V^0\cup V^p = Q(p)\cup V^p = Q(p) \cup ^{Q(p-1)}V^p
$$
est une cofibration triviale. Il s'ensuit que $Q(p)\rightarrow U^0$ est une
cofibration triviale et $U^0\rightarrow h(p)$ est une \'equivalence.

Montrons maintenant par r\'ecurrence sur $j$ que $U^j\rightarrow h(p)$ est
une \'equivalence. On vient de le d\'emontrer pour $j=0$, donc on peut supposer
que $j\geq 1$ et que c'est connu pour $j-1$. On a
$$
U^j = U^{j-1} \cup ^{(U^{j-1}\cap V^{b_j})}V^{b_j}.
$$
Pour prouver que $U^j\rightarrow h(p)$ est une \'equivalence il suffit de
prouver que
$$
(U^{j-1}\cap V^{b_j})\rightarrow V^{b_j}
$$
est une \'equivalence. Or, $V^{b_j} = h(p-1)$ et on note
$0,\ldots , \widehat{b_j}, \ldots , p$
les sommets de ce simplexe. Le sous-ensemble simplicial
$(U^{j-1}\cap V^{b_j})$ est une r\'eunion de faces de ce $h(p-1)$. Cette
r\'eunion
de faces contient la premi\`ere et la derni\`ere face (car $V^0\cup V^p\subset
U^{j-1}$), et ne contient pas la face o\`u l'on enl\`eve le sommet $a$. Donc, le
lemme pour $p-1$ s'applique et on obtient que
$(U^{j-1}\cap V^{b_j})\rightarrow V^{b_j}$ est une \'equivalence comme voulu.
\eop

\numero{Strictification}

\label{strictpage}

La question de l'essentielle surjectivit\'e du morphisme $\Phi$ du
th\'eor\`eme \ref{correlation}
pose le probl\`eme de ``strictifier'' un foncteur large $A: Y \rightarrow
nSeCAT'$ en un foncteur strict $Y \rightarrow nSeCat$. Dans la litt\'erature,
les r\'eponses \`a ce type de questions s'appellent des th\'eor\`emes de
``strictification'', ou de ``coh\'erence''.
Voir la discussion du \S 4 de Baez-Dolan \cite{BaezDolan};
ils y font r\'ef\'erence \`a Gordon-Power-Street \cite{Gordon-Power-Street},
et au th\'eor\`eme de coh\'erence de Mac Lane \cite{MacLane2}.
Voir aussi Dunn \cite{Dunn97}. 
Voir Cordier et Porter \cite{CordierPorter2} \cite{CordierPorter3}
qui parlent de ``rectification'' des diagrammes, et Segal
(\cite{SegalTopology} Prop. B1).

Un th\'eor\`eme tr\`es proche de ce qu'on va faire dans ce chapitre est
``A realization theorem 2.4'' de Dwyer-Kan \cite{DwyerKanDiags}.
Dans leur th\'eor\`eme (qui concerne les pr\'efaisceaux simpliciaux i.e.
$0$-champs de Segal), un ``foncteur large'' un foncteur
$Free _{\cdot}(Y)\rightarrow EnsSpl$ o\`u $Free_{\cdot}(Y)$ est la
r\'esolution standard par
des cat\'egories libres. Ils montrent que tout foncteur
large en ce sens, est \'equivalent \`a un foncteur strict $Y\rightarrow
EnsSpl$.

Le th\'eor\`eme-type de ce genre est le r\'esultat de SGA 1 \cite{SGA1}
qui dit qu'une cat\'egorie fibr\'ee est \'equivalente \`a une cat\'egorie
fibr\'ee scind\'ee. Notre m\'ethode de d\'emonstration
pour \ref{strictif3} sera
bas\'ee sur la m\'ethode de SGA 1, mais avec un apport des techniques de
``cat\'egories de Reedy'' issues de la th\'eorie de l'homotopie \cite{Reedy}
\cite{BousfieldKan} \cite{Hirschhorn} \cite{DHK} \cite{DwyerKanDiags}.
Cette m\'ethode a d\'ej\`a \'et\'e utilis\'ee par le deuxi\`eme auteur
pour obtenir un r\'esultat de strictification
des pr\'efaisceaux faibles d'espaces
topologiques (``pr\'efaisceaux flexibles'') dans \cite{flexible},
qui se trouve \^etre essentiellement \'equivalent \`a celui de
\cite{DwyerKanDiags}.

On peut envisager un cadre un peu plus g\'en\'eral.
Soit $M$ une cat\'egorie de mod\`eles ferm\'ee, et soit $Y$ une
$1$-cat\'egorie. Soit $L(M)'$ un remplacement fibrant pour la cat\'egorie de
Segal $L(M)$ (qui est la localis\'ee de Dwyer-Kan de $M$ par
rapport aux \'equivalences
faibles). Soit
$F: Y \rightarrow L(M)'$ un foncteur; on voudrait trouver un foncteur
$G: Y \rightarrow M$ tel
que
$p\circ G \sim F$ o\`u $p: M\rightarrow L(M)'$ est le morphisme
tautologique. Pour la
cat\'egorie
de mod\`eles ferm\'ee $M=nSePC$, ceci donnerait la strictification pour des
foncteurs $A:Y \rightarrow nSeCAT'$ au vu
du th\'eor\`eme \ref{intereqloc}.

On va g\'en\'eraliser encore plus en consid\'erant, au lieu d'une seule
cat\'egorie de mod\`eles ferm\'ee $M$, un pr\'efaisceau de Quillen \`a gauche
$\mm$ au-dessus de $Y ^o$. On retrouvera le r\'esultat du paragraphe
pr\'ec\'edent en prenant le pr\'efaisceau de Quillen constant $\mm =
\underline{M}$.

On rappelle que $\int '_YL(\mm )$ est un remplacement fibrant pour $\int
_YL(\mm )$ au-dessus de $Y$
(i.e. le morphisme vers $Y$ devient un morphisme fibrant de
$1$-cat\'egories de Segal).

\begin{definition}
\label{sestrictifient}
Soit $Y$ une petite cat\'egorie et $\mm$ un pr\'efaisceau de Quillen \`a gauche
sur $Y^o$. On dira que {\em les sections faibles de $\mm$ au-dessus de $Y$ se
strictifient} si le morphisme de $1$-pr\'ecats de Segal
$$
Sect (Y, \int _Y\mm )\rightarrow Sect (Y, \int '_YL(\mm ))
$$
est essentiellement surjectif.
\end{definition}

\begin{theoreme}
\label{strictif1A}
Si $\mm $ est un pr\'efaisceau de Quillen \`a gauche sur une cat\'egorie de
Reedy
$Y$ alors les sections faibles de $\mm$ au-dessus de $Y$ se
strictifient au sens de la d\'efinition \ref{sestrictifient}.
\end{theoreme}

On va donner la d\'emonstration plus bas.

On devrait avoir le m\^eme r\'esultat sur une cat\'egorie quelconque, \`a
condition que les $\mm (y)$ admettent suffisamment de limites. Nous avons une
esquisse de d\'emonstration mais il y a une partie que nous n'avons pas voulu
faire, donc nous laissons l'\'enonc\'e sous forme ``conjecture''.

\begin{conjecture}
\label{strictif1B}
Soit $Y$ une petite  cat\'egorie et $\mm$ un pr\'efaisceau de Quillen \`a gauche
sur $Y$. Supposons que $\mm$ satisfait l'hypoth\`ese (o) qui sera introduite
au \S 19 (qui dit en particulier que les $\mm (y)$ admettent des
petites limites arbitraires).  Alors les sections faibles de $\mm$
au-dessus de $Y$ se
strictifient au sens de la d\'efinition \ref{sestrictifient}.
\end{conjecture}

Nous proposerons plus bas une esquisse de d\'emonstration pour cette
conjecture.

\begin{corollaire}
\label{strictif2A}
Soit $Y$ une $1$-cat\'egorie de Reedy.
Soit $M$ une cat\'egorie de mod\`eles ferm\'ee.
Alors le morphisme (de $1$-cat\'egories de Segal)
$$
\underline{Hom} (Y,M)\rightarrow \underline{Hom}(Y,L(M)')
$$
est essentiellement surjectif.
\end{corollaire}
{\em Preuve:} On applique le th\'eor\`eme \ref{strictif1A} au pr\'efaisceau de
Quillen constant $\underline{M}$.
\eop

La conjecture \ref{strictif1B} donnerait le m\^eme corollaire pour toute
cat\'egorie $Y$, \`a condition d'avoir l'hypoth\`ese (o) du \S 19 pour
$\underline{M}$.
En fait, on s'int\'eresse \`a l'\'enonc\'e en question pour $M=nSePC$: le
th\'eor\`eme suivant terminera enfin la d\'emonstration du Th\'eor\`eme
\ref{correlation}.

\begin{theoreme}
\label{strictif3}
Soit $Y$ une petite $1$-cat\'egorie (non de Segal). Les sections faibles du
pr\'efaisceau constant  $\underline{nSePC}$ au-dessus de $Y$ se strictifient
au sens de la d\'efinition \ref{sestrictifient}.

En particulier, tout
objet de $\underline{Hom}(Y, nSeCAT')$ est \'equivalent (dans cette
$n+1$-cat\'egorie de Segal) \`a un objet provenant d'un pr\'efaisceau de
$n$-cat\'egories de Segal sur $Y$, i.e. le morphisme $\Phi$ du Th\'eor\`eme
\ref{correlation} est essentiellement surjectif.
\end{theoreme}

{\em D\'emonstration de \ref{strictif3} utilisant la conjecture
\ref{strictif1B}:}
On applique \ref{strictif1B} au pr\'efaisceau constant $\underline{M}$ pour
la cmf
$M=nSePC$ des $n$-pr\'ecats de Segal. Il est facile de voir que
$\underline{M}$ satisfera
l'hypoth\`ese (o). Le morphisme tautologique
$$
M_f \rightarrow nSeCAT'
$$
envoie les \'equivalences sur
des \'equivalences, et s'\'etend donc (par la propri\'et\'e
universelle de la localis\'ee) en un morphisme
$$
L(M_f) \rightarrow nSeCAT'.
$$
Le morphisme $L(M_f)\rightarrow L(M)$ est une cofibration triviale,
tout comme le morphisme vers un remplacement fibrant
$L(M)\rightarrow L(M)'$, et on obtient un
morphisme
$$
L(M)' \rightarrow nSeCAT'.
$$
D'apr\`es le th\'eor\`eme \ref{intereqloc} ce
morphisme induit une \'equivalence
$$
L(M)' \stackrel{\cong}{\rightarrow} (nSeCAT')^{int, 1}.
$$
Ce morphisme est par choix strictement compatible avec le morphisme de
source $M_f$. Le fait que $Y$ est une $1$-cat\'egorie (en particulier
$1$-groupique) entra\^{\i}ne l'\'egalit\'e
$$
\underline{Hom}(Y, nSeCAT')^{int, 1} =
\underline{Hom}(Y, (nSeCAT')^{int, 1}).
$$
Un objet  $f\in \underline{Hom}(Y, nSeCAT')$ est donc \'equivalent (dans cette
$n+1$-cat\'egorie de Segal) \'a un morphisme provenant de
$g: Y\rightarrow L(M)'$. En appliquant la version du corollaire
\ref{strictif2A}
correspondant \`a \ref{strictif1B}, on obtient que $g$ est
\'equivalent (dans $\underline{Hom}(Y, L(M)')$) \`a un morphisme provenant de
$h: Y\rightarrow M$; en composant avec un remplacement fibrant fonctoriel qui
existe pour $M=nSePC$, on obtient que $h$  peut \^etre choisi comme
morphisme $h:
Y\rightarrow M_f$. La projection de $h$ dans $\underline{Hom}(Y, nSeCAT')$ est
donc \'equivalente \`a la projection de $g$, qui \`a son tour est \'equivalente
\`a $f$.

Cet \'enonc\'e est exactement l'essentielle surjectivit\'e de $\Phi$.
\eop

Du fait que nous ne pouvons proposer qu'un esquisse de
d\'emonstration pour \ref{strictif1B}, nous allons
maintenant donner une d\'emonstration
compl\`ete du th\'eor\`eme \ref{strictif3}.

\subnumero{D\'emonstration du th\'eor\`eme \ref{strictif1A}}
On commence par la d\'emonstration du th\'eor\`eme \ref{strictif1A},
en gardant les notations $\mm$ et $Y$. Notre m\'ethode de d\'emonstration
consiste en
une analyse pr\'ecise des diagrammes index\'es par une cat\'egorie de Reedy,
utilisant les id\'ees et techniques d\'evelopp\'ees justement
dans ce but par Reedy
\cite{Reedy}, Bousfield-Kan \cite{BousfieldKan}, Hirschhorn \cite{Hirschhorn},
Dwyer-Hirschhorn-Kan \cite{DHK}. Il faut consid\'erer que notre d\'emonstration
est juste une application des techniques de ces r\'ef\'erences. On suit plus
particuli\`erement le point de vue de \cite{Hirschhorn}.

Une partie de l'analyse a d\'ej\`a \'et\'e effectu\'ee dans le lemme
\ref{reedydecomp} ci-dessus, qu'on va utiliser maintenant. Soit $F^kY$ la
filtration de $Y$ par
les sous-cat\'egories pleines form\'ees des objets de degr\'e $\leq k$. On
suppose donn\'ee une section $\sigma$ de $\int '_YL(\mm )$ et on voudrait
construire une section $u$ de $\int _Y\mm$ qui est \'equivalente \`a $\sigma$.
On va proc\'eder par r\'ecurrence sur le degr\'e,
et d\'efinir $u$ sur les
$F^kY$.
Dans la r\'ecurrence on choisira de plus $u$ fibrant et cofibrant pour la
structure de Reedy.

On suppose
donc qu'on a trouv\'e une section convenable $v$ sur $F^kY$ et on veut
l'\'etendre en $u$ sur $F^{k+1}Y$.
D'apr\`es \ref{reedydecomp}, $F^{k+1}Y$ se d\'ecompose (\`a \'equivalence faible
pr\`es) comme coproduit de $F^kY$
avec des cat\'egories de la forme
$$
Y(y):= Latch (y) + \{ y\} + Match (y)
$$
pour des objets $y\in Y$ de degr\'e $k+1$.
Posons  $F^kY(y):= F^kY\cap Y(y)$ d'o\`u
$$
F^kY(y)= Latch (y)  + Match (y).
$$
On va \'etendre $v|_{F^kY(y)}$ en une section $u_y$ sur $Y(y)$,
\'equivalente \`a $\sigma |_{Y(y)}$. Ceci suffira pour trouver $u$ et
l'\'equivalence entre $u$ et $\sigma |_{F^{k+1}Y}$ au vu de \ref{reedydecomp}.
D'autre part, la condition que $u$ soit cofibrant et fibrant pour la structure
de Reedy, s'exprime en termes des restrictions $u|_{Y(y)}$ i.e. $u$ est
cofibrant et fibrant si et seulement si (1) $v= u|_{F^kY}$ l'est, et
(2) pour tout $y$ de degr\'e  $k+1$, $u_y= u|_{Y(y)}$ l'est.

En somme, il suffit
de r\'esoudre le probl\`eme d'extension de $v$ en $u$, pour $F^kY(y) \subset
Y(y)$.
On s'est donc ramen\'e \`a la situation suivante: $Y=Y(y)$ et $Z:= F^kY(y)$;
on  a une section $\sigma$ de $\int '_YL(\mm )$ sur $Y$, et une section
$v$ de $\int _Z\mm |_Z$ qui est \'equivalente \`a $\sigma |_Z$; et on
cherche une section $u$ sur $Y$ qui \'etend $v$ et qui est \'equivalente
(en tant que section de $\int '_YL(\mm )$) \`a $\sigma$.
Rappelons les notations suivantes:
$$
Y= Latch (y) + \{ y\} + Match (y);\;\;\;
Z= Latch (y) + Match (y).
$$

On a un morphisme
$$
R: \int _{Latch(y)}L(\mm )\rightarrow L(\mm (y))
$$
correspondant aux
morphismes de restriction $res_f: \mm (z)\rightarrow \mm (y)$
pour tout $f:z\rightarrow y$ dans $Z$. Ceci donne
(par extension le long d'une cofibration triviale)
$$
R': \int '_{Latch(y)}L(\mm )\rightarrow L(\mm (y))' .
$$
Notre section $\sigma |_{Latch (y)}$ donne un morphisme
$$
R'\circ \sigma |_{Latch(y)} : Latch (y) \rightarrow L(\mm (y))'.
$$
La section $\sigma$ au-dessus de $Latch(y)+ \{ y\}$ correspond
(via une version homotopique
adapt\'ee aux cat\'egories de Segal de la description des
diagrammes de Reedy de \cite{Hirschhorn} \cite{DHK}---g\'en\'eralisation que
nous utilisons ici sans d\'emonstration) \`a un morphisme
$$
hocolim (R'\circ \sigma |_{Latch(y)})\rightarrow \sigma (y)
$$
dans $L(\mm (y))'$. D'autre part, on a aussi le morphisme
$$
{\bf r}: \int _{Latch(y)}\mm \rightarrow \mm (y)
$$
d'o\`u
$$
{\bf r}\circ v|_{Latch (y)}: Latch(y)\rightarrow \mm (y).
$$
Le fait que $v|_{Latch(y)}$ soit cofibrant de Reedy, coupl\'e avec la condition
que les foncteurs de restriction soient des foncteurs de Quillen \`a gauche,
implique que ${\bf r}\circ v|_{Latch (y)}$ est un $Latch(y)$-diagramme de
$\mm (y)$ lui-m\^eme cofibrant de Reedy. Ceci implique qu'il existe une
\'equivalence naturelle dans $L(\mm (y))'$:
$$
colim ({\bf r}\circ v|_{Latch (y)}) \cong
hocolim (R'\circ \sigma |_{Latch(y)})
$$
(en fait, c'est la d\'efinition de $hocolim$ voir \cite{BousfieldKan}
\cite{Hirschhorn} \cite{DHK}).

Pour avoir une extension de $v$ \`a une section $u_L$
au-dessus de $Latch (y)+ \{ y\}$ il suffit (d'apr\`es la description de
\cite{Hirschhorn} \cite{DHK}) d'avoir un objet $u(y)\in \mm (y)$ et un
morphisme
$$
colim ({\bf r}\circ v|_{Latch (y)})\rightarrow u(y).
$$
En outre, $u_L$ est cofibrant de Reedy (ainsi qu'une \'eventuelle extension
de $u_L$ \`a $u$ sur $Z$) si et seulement si le morphisme pr\'ec\'edent est une
cofibration. La section $u_L$ est  \'equivalente
(en tant que section des localis\'es)
\`a $\sigma |_{Latch(y)+ \{ y\}}$ s'il existe une \'equivalence dans
$L(\mm (y))'$, $u(y)\cong \sigma (y)$, et une homotopie de  commutativit\'e pour
le diagramme
$$
\begin{array}{ccc}
colim ({\bf r}\circ v|_{Latch (y)})
& \cong & hocolim (R'\circ \sigma |_{Latch(y)})\\
\downarrow && \downarrow \\
u(y) & \cong & \sigma (y).
\end{array}
$$

Dualement (encore que l'argument soit un peu diff\'erent) on \'etudie les
extensions de
$v$ en une section $u_M$ au-dessus de $Match (y)+ \{ y\}$. D'abord, notons
$$
{\bf t}: \mm (y)\times Y \rightarrow \int _{Match (y)}\mm
$$
le morphisme donn\'e par les restrictions. Pour
$$
y\stackrel{f}{\rightarrow}z\stackrel{g}{\rightarrow}w
$$
avec compos\'e $h:= gf$, on obtient, en utilisant les adjonctions, un
morphisme (pour $a\in \mm (z)$)
$$
res_f^{\ast} a \rightarrow res _h^{\ast} (res_ga).
$$
Ceci donne un morphisme
$$
{\bf t}^{\ast} : \int _{Match (y)}\mm \rightarrow
\mm (y).
$$
On obtient \'egalement sur les localis\'es
$$
T^{\ast} : \int '_{Match (y)}L(\mm )\rightarrow
L\mm (y))'.
$$
En composant avec $\sigma$ on obtient
$$
T^{\ast}\circ \sigma |_{Match(y)}: Match(y)\rightarrow L(\mm (y)).
$$
La restriction de $\sigma$ sur $\{ y\} + Match(y)$ correspond
(de nouveau via
la version homotopique pour cat\'egories de Segal, de la description des
diagrammes de Reedy de \cite{Hirschhorn} \cite{DHK} que nous ne d\'emontrons pas
ici) \`a un morphisme
$$
\sigma (y) \rightarrow
holim (T^{\ast} \circ \sigma |_{Match(y)})
$$
dans $L(\mm (y))'$.

De fa\c{c}on similaire, une extension $u_M$ de $v$ \`a $\{ y\} + Match(y)$
correspond \`a un objet $u(y)\in \mm (y)$ et un morphisme
$$
u(y)\rightarrow lim ({\bf t}^{\ast} \circ v|_{Match(y)}).
$$
Ici la limite existe car l'indexation est finie.
La section $u_M$ ainsi que son \'eventuelle extension $u$ \`a tout $Z$, est
fibrante de Reedy si et seulement si ce morphisme est une fibration dans $\mm
(y)$.

D'autre part, comme plus haut on a une \'equivalence dans $L(\mm (y))'$
$$
lim ({\bf t}^{\ast} \circ v|_{Match(y)})\cong
holim (T^{\ast} \circ \sigma |_{Match(y)})
$$
car $v|_{Match(y)}$ est fibrant de Reedy par hypoth\`ese de r\'ecurrence, et les
adjoints $res _f^{\ast}$ sont des foncteurs de Quillen \`a droite, donc
${\bf t}^{\ast} \circ v|_{Match(y)}$ est un $Match(y)$-diagramme de $\mm (y)$
qui est
fibrant de Reedy. Une extension de l'\'equivalence de $v$ avec $\sigma |_Z$ en
une \'equivalence entre $u_M$ et $\sigma |_{\{ y\} + Match (y)}$ correspond
\`a la donn\'ee d'une \'equivalence $u(y)\cong \sigma (y)$ et d'une homotopie de
commutativit\'e pour le diagramme
$$
\begin{array}{ccc}
u(y) & \cong & \sigma (y) \\
\downarrow && \downarrow \\
lim ({\bf t}^{\ast} \circ v|_{Match(y)})& \cong &
holim (T^{\ast} \circ \sigma |_{Match(y)}).
\end{array}
$$

Soient $x\in Latch(y)$ et $z\in Match(y)$
et notons par $f: x\rightarrow y$ et $g: y\rightarrow z$ les morphismes.
Alors---
par la propri\'et\'e d'adjoint homotopique (Lemme \ref{adjointh}) des foncteurs
induits par $res_g$ et $res_g^{\ast}$ entre les cat\'egories simpliciales $L(\mm
(y))$ et $L(\mm (z))$--- le morphisme
$$
res_{gf}(\sigma (x))\rightarrow \sigma (z)
$$
est \'egal \`a un morphisme
$$
res_f(\sigma (x)) \rightarrow res_g^{\ast} (\sigma (z)).
$$
Quand $x$ varie dans $Latch (y)$ et $z$ dans $Match (y)$ on obtient le
r\'esultat suivant.
Notons $Latch(y)\sqcup Match(y)$ la r\'eunion disjointe de ces deux
cat\'egories; elle
est contenu dans $Z$ mais dans $Z$ il y a, en plus, un (unique)
morphisme de chaque objet de $Latch(y)$ vers chaque objet de $Match(y)$.
Une extension de $\sigma |_{Latch(y)\sqcup Match(y)}$ \`a $\sigma |_Z$
correspond \`a un morphisme
$$
hocolim (R'\circ \sigma |_{Latch(y)})\rightarrow
holim (T^{\ast} \circ \sigma |_{Match(y)}).
$$
Si $\sigma |_Z$ est \'equivalente \`a la section induite par $v$, l'application
ci-dessus est induite par
$$
colim ({\bf r}\circ v|_{Latch (y)})\rightarrow
lim ({\bf t}^{\ast} \circ v|_{Match(y)}).
$$

L'extension de $\sigma$ de $Z$ \`a $Y=Z\cup \{ y\}$ est d\'etermin\'ee
par une factorisation
$$
hocolim (R'\circ \sigma |_{Latch(y)})\rightarrow \sigma (y) \rightarrow
holim (T^{\ast} \circ \sigma |_{Match(y)}).
$$
Celle-ce \'etant donn\'ee, pour \'etendre $v$ en une section $u$ sur $Y$,
il faut
trouver une factorisation
$$
colim ({\bf r}\circ v|_{Latch (y)})\rightarrow u(y)\rightarrow
lim ({\bf t}^{\ast} \circ v|_{Match(y)})
$$
qui induise la factorisation correspondant \`a $\sigma$.

On note par ailleurs que  $colim ({\bf r}\circ v|_{Latch (y)})$
est cofibrant (c'est la colimite d'un diagramme qui est cofibrant pour la
structure de Reedy, car ${\bf r}$ est un foncteur de Quillen \`a gauche);
et $lim ({\bf t}^{\ast} \circ v|_{Match(y)})$ est fibrant---c'est la limite
d'un diagramme qui est fibrant pour la structure de Reedy, car ${\bf t}^{\ast}$
est un foncteur de Quillen \`a droite.

On est donc ramen\'e au probl\`eme suivant: \'etant donn\'e un morphisme
$$
h:a\rightarrow b
$$
dans une cat\'egorie de mod\`eles ferm\'ee $M=\mm (y)$, et une factorisation
$$
\overline{a} \rightarrow w \rightarrow \overline{b}
$$
dans $L(M)$ o\`u $\overline{a}$ et $\overline{b}$ sont les images de $a,b\in M$
dans $L(M)$, on veut la relever en une factorisation
$$
a\rightarrow c \rightarrow b
$$
dans $M$ avec le premier morphisme cofibrant et le deuxi\`eme
fibrant. On pourra supposer que $a$ est cofibrant et $b$ fibrant.

On peut aussi supposer qu'on
a $w=\overline{c}$ dans $L(M)$ avec $c$ fibrant et
cofibrant. Dans ce cas, les morphismes de $L(M)$ proviennent de morphismes
$f: a\rightarrow c$ et $g: c\rightarrow b$. On peut supposer de plus que $g$
est une fibration. La factorisation dans $L(M)$ est une homotopie entre $gf$ et
$h$, en tant que morphismes de $a$ vers $b$. L'homotopie peut \^etre
consid\'er\'ee
comme Quillen l'a d\'efinie, c'est-\`a-dire comme un diagramme
$$
i_0,i_1 :a\tworightarrows \alpha \stackrel{p}{\rightarrow} a ,\;\;\; k:\alpha
\rightarrow b
$$
avec $pi_0=pi_1=1_a$ et $ki_0=gf$, $ki_1=h$ (et $p$ une fibration triviale, et
$i_0\sqcup i_1$ une cofibration). Le morphisme  $i_0:a\rightarrow \alpha$
est une
cofibration triviale. Le fait qu'on a suppos\'e
que $g$ est une fibration implique qu'il
existe un rel\`evement $w: \alpha \rightarrow c$ avec $gw=k$
et $wi_0= f$. Si on pose $f':= wi_1$, on obtient une homotopie (qu'on note
encore $w$) entre $f$ et $f'$;
on a aussi $gf'=h$ et l'homotopie $w$ compos\'ee avec
$g$ est l'homotopie $k$. Autrement dit, la factorisation $h=gf'$ avec
l'identit\'e comme homotopie, est \'equivalente \`a la factorisation qui
correspond \`a l'homotopie $k$ entre $gf$ et $h$. On obtient le rel\`evement de
notre factorisation dans $L(M)$ en une factorisation
$$
a\stackrel{f'}{\rightarrow}c\stackrel{g}{\rightarrow}b
$$
dans $M$, ce qui r\'esout le probl\`eme.

Pour compl\'eter la d\'emonstration du th\'eor\`eme \ref{strictif1A},
appliquons la solution du
probl\`eme ci-dessus \`a la situation pr\'ec\'edente avec
$$
a= colim ({\bf r}\circ v|_{Latch (y)}),\;\;\;\;\;\;  b =
lim ({\bf t}^{\ast} \circ v|_{Match(y)});
$$
$$
\overline{a} = hocolim (R'\circ \sigma |_{Latch(y)}), \;\;\;\;\;\;
\overline{b}=  holim
(T^{\ast} \circ \sigma |_{Match(y)});
$$
et $w = \sigma (y)$. On trouve la factorisation d\'esir\'ee avec $c=u(y)$.
\eop

\subnumero{Le cas des $1$-champs}

Avant d'aborder l'esquisse de d\'emonstration pour \ref{strictif1B}
et la d\'emonstration du th\'eo\-r\`e\-me \ref{strictif3}, signalons
l'origine de
l'id\'ee en faisant la comparaison avec le cas des $1$-champs.

On peut montrer que la $2$-cat\'egorie des $1$-champs au-dessus de $\Xx$
est \'equivalente (via la comparaison de Tamsamani \cite{Tamsamani} entre les
$2$-cat\'egories suivant sa d\'efinition, et les $2$-cat\'egories de Benabou)
\`a la $2$-cat\'egorie $1CHAMP(\Xx )$ des $1$-champs (non de Segal) au-dessus de
$\Xx$. Si on adopte la d\'efinition de $1$-champ qui utilise la notion de
cat\'egorie fibr\'ee \cite{SGA1}, il s'agit pour l'essentiel de montrer qu'une
cat\'egorie fibr\'ee est \'equivalente \`a une cat\'egorie fibr\'ee
scind\'ee, i.e. une qui provient d'un pr\'efaisceau de cat\'egories.

Ce r\'esultat est un th\'eor\`eme de SGA 1 \cite{SGA1}. L'id\'ee de la
d\'emonstration est que si $\Ff \rightarrow Y$ est une cat\'egorie fibr\'ee
au-dessus d'une cat\'egorie $Y$, alors le foncteur
$$
A:y\mapsto Sect^{eq} (Y/y, \Ff \times _Y(Y/y))
$$
(o\`u $Sect^{eq}$ sont les sections ``cart\'esiennes'' i.e. qui envoient les
fl\`eches
d'en bas en fl\`eches cart\'esiennes en haut) sera un foncteur strict de
$Y$ vers
la $1$-cat\'egorie
$1Cat$ des $1$-cat\'egories, tel qu'on ait
un e\'equivalence $\int _YA \cong \Ff $ de
cat\'egories fibr\'ees
au-dessus de $Y$.

L'observation-cl\'e est que cette propri\'et\'e de strictification provient
du fait que la base $Y$ est une $1$-cat\'egorie stricte.

Cette id\'ee a \'et\'e reprise dans \cite{flexible} pour donner un r\'esultat de
strictification des ``faisceaux flexibles'' d'espaces topologiques, ce qui est
l'analogue du th\'eor\`eme \ref{strictif3} pour le cas $n=0$.

\subnumero{Id\'ee pour preuve directe de \ref{strictif3}}

On esquisse ici une id\'ee pour une preuve directe de la strictification
pour les foncteurs $\Xx ^o \rightarrow nSeCAT'$ (Th\'eor\`eme \ref{strictif3}),
en suivant le cas des $1$-champs mentionn\'e ci-dessus.

Supposons maintenant que $M= nSePC$ est la cmf
des $n$-pr\'ecats de Segal. Un foncteur $Y ^o\rightarrow M$ est un
$n$-pr\'ecat de
Segal au-dessus de $Y$. D'autre part, on a \ref{intereqloc}
$$
L(M) \cong nSeCAT ^{int, 1}.
$$
Comme $Y$ est d\'ej\`a $1$-groupique, un foncteur
$Y ^o \rightarrow L(M) '$ est donc \`a \'equivalence pr\`es la m\^eme chose
qu'un foncteur $Y ^o\rightarrow nSeCAT '$. Notre question devient donc: est-ce
que tout foncteur  $Y ^o\rightarrow nSeCAT'$ est \'equivalent \`a
une $n$-pr\'ecat de Segal sur $Y$?

L'id\'ee de base est d'utiliser le fait que $Y$ est une $1$-cat\'egorie qui est
forcement stricte.
Soit $F: Y ^o \rightarrow nSeCAT'$ un morphisme. Pour
$y\in Y$, on pose
$$
G(y):= \Gamma (Y/y, F|_{Y/y}):= \underline{Hom}((Y/y)^o,
nSeCAT')_{1/}(\ast _y,  F|_{Y/y}).
$$
Ceci varie fonctoriellement en $y$ et les valeurs sont des $n$-cat\'egories
de Segal fibrantes. Le probl\`eme est de montrer que $G$ est \'equivalent
\`a $F$.

Nous ne savons pas actuellement trouver directement un morphisme
entre
$F$ et $G$ (et c'est pour cela que ce paragraphe n'est qu'une esquisse).
Cependant, notons pour plus tard, que si $F$ provenait d\'ej\`a
d'une section stricte, alors il y
aurait une \'equivalence entre $F$ et $G$ donn\'ee
par le morphisme
$$
F(y)\rightarrow \Gamma (Y/y, F|_{Y/y}).
$$

Nous allons contourner le probl\`eme de trouver un morphisme entre $F$ et $G$
en nous appuyant sur le th\'eor\`eme \ref{strictif1A},
et en rempla\c{c}ant $Y$ par une
cat\'egorie de Reedy $B$ munie d'une sous-cat\'egorie $D\subset B$ telle que
la localis\'ee $L(B,D)$ soit \'equivalente \`a $Y$. L'id\'ee que les diagrammes
sur $Y$ sont les m\^emes que les diagrammes sur $B$ qui sont des \'equivalences
dans la direction de $D$, provient du papier de Dwyer-Kan \cite{DwyerKanDiags}.

\subnumero{Esquisse de d\'emonstration de la conjecture \ref{strictif1B}:}
Soit $Y$ une  $1$-cat\'egorie et $\mm$ un pr\'efaisceau de Quillen \`a
gauche sur
$Y$ qui satisfait les conditions (o) du \S 19 ci-dessous. On aura besoin de la
cat\'egorie de Reedy suivante. Soit  $\nu Y$ le nerf de $Y$, et soit
$B=\int _{\Delta} \nu Y$ la cat\'egorie introduite pour le lemme \ref{bd}.
On note que $B$ est aussi la cat\'egorie des simplexes de $\nu Y$, not\'ee
$B=\Delta \nu Y$ dans \cite{Hirschhorn} \cite{DHK}; d'apr\`es {\em loc cit.}
$B$ est une cat\'egorie de Reedy (la structure de Reedy \'etant induite par le
foncteur $B\rightarrow \Delta$). Le morphisme d'ensembles simpliciaux
$$
\nu B\rightarrow \nu Y
$$
provient d'un foncteur $B\rightarrow Y$, et en plus ce foncteur envoie
les morphismes de $D\subset B$ sur les identit\'es de $Y$. D'apr\`es le Lemme
\ref{bd}, ceci induit une \'equivalence de $1$-cat\'egories simpliciales ou de
Segal
$$
L(B,D)\cong 1SeL(B,D) \stackrel{\cong}{\rightarrow} Y.
$$

Consid\'erons le diagramme suivant:
$$
\begin{array}{ccc}
Sect (Y, \int _Y\mm ) & \rightarrow & Sect (Y, \int '_YL(\mm ))\\
\downarrow && \downarrow \\
Sect (B, \int _B \varphi ^{\ast} \mm )&\rightarrow & Sect (B, \int '_B \varphi
^{\ast}L(\mm )).
\end{array}
$$
A partir du r\'esultat \ref{strictif1A}, on obtient que le morphisme du bas
est essentiellement surjectif. D'autre part, l'\'equivalence entre $Y$ et
$L(B,D)$ implique que le morphisme vertical de droite est pleinement fid\`ele
avec pour image essentielle la classe
des sections $\sigma : B\rightarrow \varphi ^{\ast}L(\mm )'$ telles que
$\sigma (g)$
soit une \'equivalence pour tout morphisme $g$ de $D$
(cette id\'ee provient de Dwyer-Kan \cite{DwyerKanDiags}). On appellera cette
condition $D-eq$. Si une section $\sigma $ est $D-eq$ et si $\sigma$ provient
d'une section $u: B\rightarrow \int _B\varphi ^{\ast}\mm$ alors $u$ est aussi
$D-eq$ dans le sens que les morphismes de restriction $res_g(u)$ sont des
\'equivalences faibles
pour $g\in D$. Il nous suffit de d\'emontrer que si
$u\in Sect (B, \int _B\varphi ^{\ast} \mm )$ est une section $D-eq$ en ce
sens, alors $u$ provient (\`a \'equivalence faible pour la structure de
cat\'egorie de mod\`eles pr\`es) d'une section $v\in Sect (Y, \int _Y\mm )$.

Le foncteur
$$
p^{\ast}: Sect (Y, \int _Y\mm )\rightarrow Sect (B, \int _Bp^{\ast}\mm )
$$
admet un adjoint \`a droite
$$
p_{\ast}: Sect (B, \int _Bp^{\ast}\mm )
\rightarrow  Sect (Y, \int _Y\mm ).
$$
On pr\'etend que si $u\in Sect ^{D-eq}(B, \int _Bp^{\ast}\mm )$
est fibrant pour la structure de type (II),
alors on a l'\'equivalence
$p^{\ast}p_{\ast} u\cong u$. Pour ceci on envisage une d\'emonstration qui
suivrait les lignes de l'argument du \S 19 ci-dessous (en particulier c'est
pour cela qu'on a besoin de l'hypoth\`ese (o) dont l'essentiel est l'existence
de limites dans les $\mm (y)$). Ce r\'esultat serait en quelque sorte une
version
relative de l'argument du \S 19. Cependant, pour la pr\'esente version du papier
nous n'avons pas v\'erifi\'e les d\'etails de cet argument et
c'est pour cela que la
d\'emonstration de \ref{strictif1B} comporte une lacune.

En admettant ce fait, on obtient que tout
$u\in Sect ^{D-eq}(B, \int _Bp^{\ast}\mm )$ est \'equivalent \`a
une section de la forme $p^{\ast} v$ pour
$v\in Sect (Y, \int _Y\mm )$, ce qui termine la d\'emonstration d'apr\`es les
remarques pr\'ec\'edentes.
\eop

\subnumero{Preuve du Th\'eor\`eme \ref{strictif3}}

Du fait qu'il y a une lacune dans notre esquisse de d\'emonstration de la
conjecture \ref{strictif1B}, et que l'id\'ee pour la preuve directe donn\'ee
auparavant n'a pas abouti non plus, on donne enfin une
d\'emonstration compl\`ete en combinant ces deux id\'ees. On devait appliquer
\ref{strictif1B} au cas du pr\'efaisceau de Quillen constant $\mm =
\underline{M}$
\`a valeurs $M=nSePC$, la cat\'egorie de mod\`eles de $n$-pr\'ecats de Segal.
On note par exemple que $\int _Y\mm = Y\times M$ etc.  Avec ces notations, on
revient au diagramme
$$
\begin{array}{ccc}
Sect (Y, Y\times M ) & \rightarrow & Sect (Y, Y\times L(M )')\\
\downarrow && \downarrow \\
Sect (B, B\times M )&\rightarrow & Sect (B, B\times L(M )').
\end{array}
$$
Etant donn\'e $\sigma : Y\rightarrow L(M)'$,
on peut appliquer la construction
du paragraphe avant le pr\'ec\'edent, pour obtenir une section $G_{\sigma} :
Y\rightarrow
M$, essentiellement par $G_{\sigma}(y)= \lim _{Y/y} \sigma |_{Y/y}$. Ceci est
compatible avec la restriction \`a $B$ (via le morphisme $\psi :
B\rightarrow Y$) en ce sens qu'on a un morphisme
$$
r:G_{\sigma} \circ  \psi \rightarrow G_{\sigma \circ \psi }.
$$
Le fait que $\psi$ induise une \'equivalence $L(B,D)\cong Y$ implique
que pour $b\in B$, on a
$$
\lim _{Y/\psi (b)}\sigma |_{Y/\psi (b)}\stackrel{\cong}{\rightarrow}
\lim _{B/b}(\sigma \circ \psi )|_{B/b},
$$
i.e. que le morphisme $r$ est une \'equivalence (objet-par-objet au-dessus de
$B$). D'autre part, par le th\'eor\`eme \ref{strictif1A}, $\sigma
\circ \psi $ est \'equivalente \`a une section provenant de $u: B\rightarrow M$.
Pour les sections strictes, on a le morphisme d\'esir\'e au paragraphe
avant le pr\'ec\'edent. On a donc une \'equivalence
$$
G_u \cong u,
$$
ce qui donne
$$
G_{\sigma} \circ \psi \cong G_{\sigma \circ \psi} \cong \sigma \circ \psi .
$$
En posant $v:= G_{\sigma} : Y\rightarrow M$, on a r\'esolu (pour ce cas) le
probl\`eme  qui a donn\'e naissance \`a la lacune dans la d\'emonstration de la
conjecture \ref{strictif1B}. On finit la d\'emonstration comme avant: le
fait que
$\psi$ induise une \'equivalence $L(B,D)\stackrel{\cong}{\rightarrow} Y$
implique que le morphisme $Sect (Y,Y\times L(M)')\rightarrow Sect(B, B\times
L(M)')$ est pleinement fid\`ele (cf \cite{DwyerKanDiags}).
Donc, l'\'equivalence $G_{\sigma}\circ \psi\cong \sigma \circ \psi$ provient
d'une \'equivalence $G_{\sigma} \cong \sigma$.
Ceci donne enfin une d\'emonstration compl\`ete du th\'eor\`eme \ref{strictif3}.
\eop

Le th\'eor\`eme \ref{strictif3} donne l'essentielle surjectivit\'e qui
manquait
pour terminer la d\'e\-mon\-stra\-tion du th\'eor\`eme \ref{correlation}.

{\em Remarque--Exercice:}
\newline
Dans la d\'emonstration de \ref{strictif3} on
aurait pu
prendre $B=\beta ^{\rm poset}(Y)$, voir \ref{bdvariant}. Dans ce cas,
$B$ est la cat\'egorie sous-jacente \`a un ensemble partiellement ordonn\'e;
de ce fait on n'aurait besoin du th\'eor\`eme \ref{strictif1A} que pour ce
type de cat\'egorie qui est une cat\'egorie de Reedy directe. Cela
simplifierait beaucoup la d\'emonstration de \ref{strictif1A}, en particulier
les objets appariants $Match$ serait triviaux et on n'aurait  pas besoin de
consid\'erer les morphismes de connexion $Latch\rightarrow Match$. On laisse au
lecteur l'exercice de r\'ediger une d\'emonstration de \ref{strictif3}
(incluant la partie de \ref{strictif1A} dont on aurait besoin) sur la base de
cette remarque.

\subnumero{Une \'equivalence}

On donne maintenant une am\'elioration du r\'esultat pr\'ec\'edent
concernant
l'essentielle surjectivit\'e du morphisme \ref{sestrictifient}; ce morphisme
devient une
\'equivalence si l'on localise (\`a la Dwyer-Kan) la cat\'egorie de mod\`eles
\`a la source du morphisme \ref{sestrictifient}. On ne donne l'\'enonc\'e
que pour le cas des
cat\'egories de Reedy; on pense qu'il reste vrai pour une cat\'egorie
de base quelconque, sous
des hypoth\`eses du type (o) du \S 19.

\begin{theoreme}
\label{uneequivalence}
Soit $Y$ une cat\'egorie de Reedy, et soit $\mm$ un pr\'efaisceau de Quillen \`a
gauche sur $Y$. Supposons que chaque $\mm (y)$ admet des ``factorisations
fonctorielles'' pour la structure de Quillen. Supposons la m\^eme chose pour
$Sect (Z,\int _Z\mm |_Z)$ pour toute cat\'egorie de Reedy $Z$ avec foncteur
vers $Y$ (ce foncteur ne respectant pas forc\'ement la structure de Reedy).
Alors le
morphisme
$$
L(Sect (Y,\int _Y\mm ))\rightarrow Sect (Y,\int '_YL(\mm ))
$$
est une \'equivalence.
\end{theoreme}

\begin{corollaire}
\label{pqcalclim}
Soit $Y$ une cat\'egorie de Reedy et soit $\mm$ un pr\'efaisceau de Quillen \`a
gauche sur $Y$. Alors on a une \'equivalence
$$
\lim _{\leftarrow , Y^o}L(\mm )\cong L(Sect ^{\rm eq}(Y,\int _Y\mm )).
$$
\end{corollaire}
{\em Preuve:}
Dans l'\'equivalence du th\'eor\`eme, une section $\sigma \in Sect (Y,\int _Y\mm
)$ va sur une $eq$-section de $\int 'L(\mm )$ si et seulement si les morphismes
de transition $res_f\sigma (y)\rightarrow \sigma (z)$ (pour $f:y\rightarrow z$
dans $Y$) sont des \'equivalences faibles; i.e. si et seulement si $\sigma$ est
une $eq$-section de $\int _Y\mm$. On a donc (voir \ref{stabilite}) une
\'equivalence:
$$
L(Sect ^{\rm eq}(Y,\int _Y\mm ))
\stackrel{\cong}{\rightarrow} Sect ^{\rm eq}(Y,\int '_YL(\mm )).
$$
D'autre part, on a par la proposition \ref{gammasect} que pour le remplacement
fibrant $L(\mm )'$ de $L(\mm )$,
$$
\Gamma (Y, L(\mm )')\rightarrow
Sect ^{\rm eq}(Y,\int '_YL(\mm )')
$$
est une \'equivalence; le morphisme
$$
Sect ^{\rm eq}(Y,\int '_YL(\mm ))\rightarrow Sect ^{\rm eq}(Y,\int '_YL(\mm )')
$$
est une \'equivalence; et d'apr\`es la proposition \ref{calclim},
$$
\lim _{\leftarrow , Y^o}L(\mm ) \cong \Gamma (Y, L(\mm )').
$$
D'o\`u l'\'enonc\'e.
\eop

Avant de faire la d\'emonstration du th\'eor\`eme \ref{uneequivalence}
en g\'en\'eral, nous allons en traiter quelques cas particuliers. On isole
ces lemmes car ils montrent bien, par des
exemples, ce qui se passe dans \ref{uneequivalence}.

Pour ces lemmes, on utilisera la technique des ``homotopy function complexes''
d\'e\-ve\-lop\-p\'ee dans Hirschhorn \cite{Hirschhorn}, Dwyer-Hirschhorn-Kan
\cite{DHK} et qui a ses origines dans le travail sur la localisation
\cite{DwyerKan3} ainsi que dans les travaux de Reedy \cite{Reedy}.
On rappelle bri\`evement de
quoi il s'agit.
Ce rappel chevauche partiellement celui du \S 8.

Ici $M$ est une cat\'egorie de mod\`eles ferm\'ee, et on veut
calculer les types d'homotopie des ensembles simplicaux $Hom$ dans la
localis\'ee de Dwyer-Kan $L(M)$, i.e. pour $x,y\in M$ on veut calculer
$L(M)_{1/}(x,y)$. Pour ceci, on introduit (dans les r\'ef\'erences ci-dessus)
la notion de {\em r\'esolution simpliciale fibrante} $y\rightarrow {\bf y}$.
Ici ${\bf y}$ est un objet simplicial de $M$, i.e. un objet ${\bf y}\in
M^{\Delta ^o}$, avec morphisme $c^{\ast} y\rightarrow {\bf y}$ o\`u
$c^{\ast} y$ est l'objet simplicial constant \`a valeurs $y$. On demande que
pour $p\in \Delta$, $y\rightarrow {\bf y}(p)$ soit une \'equivalence faible;
et que ${\bf y}$ soit fibrant pour la structure de cat\'egorie de mod\`eles
ferm\'ee de Reedy (i.e. ce qu'on appelle ici ``de
type (III)'') sur $M^{\Delta ^o}$.
Dans ce cas, l'ensemble simplicial
$$
M_{1/}(x, {\bf y}):= \; p \mapsto M_{1/}(x, {\bf y}(p))
$$
est un ensemble simplicial
qui est naturellement \'equivalent \`a $L(M)_{1/}(x,y)$
\cite{DwyerKan3}.
\footnote{
Techniquement le
morphime naturel ici est un morphisme
$$
M_{1/}(x, {\bf y})\rightarrow L^H(M)
$$
dans la localis\'ee {\em par hamacs} voir \cite{DwyerKan3} 4.4 et 7.2.
On note que l'argument de \cite{DwyerKan3} 7.2, qu'ils donnent pour
le cas d'une paire de r\'esolutions ${\bf x}$ ${\bf y}$ (cosimpliciale et
simpliciale respectivement), marche \'egalement pour fournir le morphisme qu'on
cherche ici dans le cas d'une seule r\'esolution. Dans ce qui suit on utilisera
la notation $L(M)$ mais pour \^etre techniquement correct (i.e. pour avoir les
morphismes qu'on dit) il faudrait lire $L^H(M)$. }
La notation $M_{1/}(a,b):= Hom _M(a,b)$
d\'esigne l'ensemble des morphismes pour la cat\'egorie $M$.

On note maintenant que l'ensemble simplicial $M_{1/}(x, {\bf y})$ est
fibrant i.e. de Kan. En effet, la condition d'extension de Kan pour
$M_{1/}(x, {\bf y})$ devient une condition de rel\`evement pour des morphismes
de $x$ vers une fl\`eche de la forme ${\bf y}(p)\rightarrow {\bf w}$ o\`u
${\bf w}$ est un produit fibr\'e (r\'ecursif) de composants de ${\bf y}$ (ce
produit fibr\'e d\'epend de quelle condition d'extension i.e. quelle
``corne'' on
regarde;  nous laissons au lecteur d'\'ecrire pr\'ecis\'ement l'expression
pour le
produit fibr\'e en question). La condition que ${\bf y}$ est fibrante pour la
structure de Reedy implique que ces morphismes ${\bf y}(p)\rightarrow {\bf w}$
sont fibrants; et la condition que  tous les ${\bf y}(p)$ sont \'equivalents \`a
$y$ implique (au vu de la forme du produit fibr\'e qui correspond \`a une
``corne'') que ${\bf y}(p)\rightarrow {\bf w}$ est une \'equivalence faible.
Donc (puisque $x$ est cofibrant) tout morphisme de $x$ dans  ${\bf w}$ se
rel\`eve en un morphisme $x\rightarrow {\bf y}(p)$, ce qui donne la condition
d'extension de Kan.

Le m\^eme argument montre qu'une cofibration $x\rightarrow x'$ induit une
fibration de Kan
$$
M_{1/}(x', {\bf y})\rightarrow M_{1/}(x, {\bf y}).
$$

Le premier lemme concerne le produit direct: c'est le cas du th\'eor\`eme
\ref{uneequivalence} o\`u la cat\'egorie de base $Y$ est discr\`ete, \'egale
\`a un ensemble. Malgr\'e
les apparences, ce cas n'est pas totalement \'evident \`a
cause du fait qu'un produit infini d'ensembles simpliciaux peut ne pas avoir le
bon type d'homotopie, si les facteurs ne sont pas fibrants (de Kan). Ce fait a
\'et\'e remarqu\'e notamment par Jardine dans \cite{JardineBool} (qui donne un
contre-exemple).

\begin{lemma}
\label{equivalence1}
Soit $\{ M_i\}$ une collection de cat\'egories de mod\`eles
ferm\'ees (index\'ee par un ensemble). On munit le produit $\prod _i M_i$
de la structure de cat\'egorie de mod\`eles ferm\'ee produit: les fibrations,
cofibrations et \'equivalences faibles sont les morphismes qui le sont
par rapport \`a
chaque variable.  Soient $L(M_i)'$ des remplacements fibrants (en tant que
cat\'egories de Segal) des localis\'ees $L(M_i)$ de Dwyer-Kan.
Alors le morphisme naturel
$$
L(\prod _i M_i)\rightarrow \prod _i L(M_i)'
$$
(induit par sa restriction \`a $\prod _i M_i$) est une \'equivalence.
\end{lemma}
{\em Preuve:}
On note d'abord que l'ensemble des objets est le m\^eme des deux cot\'es;
il s'agit donc de prouver que le morphisme est (homotopiquement) pleinement
fid\`ele. Soient $(y_i)$ et $(z_i)$ deux objets du produit $\prod _i M_i$.
On choisit pour chaque $i$ une r\'esolution simpliciale fibrante $z_i
\rightarrow
{\bf z}_i$.  On obtient ainsi une r\'esolution simpliciale fibrante
$$
(z_i)\rightarrow ({\bf z}_i)
$$
pour le produit. On a donc
$$
L(\prod _i M_i) _{1/}((y_i), (z_i))\cong
(\prod _i M_i)_{1/}((y_i), ({\bf z}_i)) = \prod _i \left(
(M_i)_{1/}(y_i,{\bf z}_i)\right) .
$$
D'autre part, pour chaque $i$ le morphisme
$$
L(M_i)_{1/}(y_i, z_i)\rightarrow L(M_i)'_{1/}(y_i,z_i)
$$
est une \'equivalence d'ensembles simpliciaux avec le deuxi\`eme terme fibrant.
On a  pour chaque $i$
$$
L(M_i)_{1/}(y_i, z_i)\cong (M_i)_{1/}(y_i,{\bf z}_i),
$$
et comme on l'a remarqu\'e ci-dessus, les
$(M_i)_{1/}(y_i,{\bf z}_i)$ sont fibrants. Le produit direct d'ensembles
simpliciaux fibrants conserve les \'equivalences d'homotopie
\cite{JardineBool}, donc on a
$$
\prod _i \left(L(M_i)'_{1/}(y_i, z_i)\right) \cong \prod _i \left(
(M_i)_{1/}(y_i,{\bf z}_i)\right) .
$$
Notons qu'on a
$$
\prod _i \left(L(M_i)'_{1/}(y_i, z_i)\right) =
(\prod _i L(M_i)')_{1/}(y_i, z_i);
$$
on obtient l'\'equivalence
$$
(\prod _i L(M_i)')_{1/}(y_i, z_i)\cong
L(\prod _i M_i) _{1/}((y_i), (z_i)).
$$
Cette \'equivalence est homotope au morphisme
$$
L(\prod _i M_i) _{1/}((y_i), (z_i))\rightarrow
(\prod _i L(M_i)')_{1/}(y_i, z_i)
$$
du lemme, qui est (par d\'efinition) celui induit par le morphisme
$$
(\prod _i M_i) _{1/}((y_i), (z_i)\rightarrow
(\prod _i L(M_i)')_{1/}(y_i, z_i),
$$
qui envoie les \'equivalences faibles sur des \'equivalences.
\eop

\begin{lemma}
\label{equivalence2}
Soit $M$ une cat\'egorie de mod\`eles ferm\'ee. Soit $I$ la cat\'egorie
avec objets $0,1$ et un morphisme $0\rightarrow 1$, et soit $M^I$ la
cat\'egorie des $I$-diagrammes dans $M$ (qui admet une structure de cmf de
Reedy, par exemple). Soit $L(M)$ un remplacement fibrant (en tant que
cat\'egorie de Segal) pour $L(M)$. Alors le morphisme naturel
$$
L(M^I)\rightarrow \underline{Hom}(I, L(M)')
$$
(induit par sa restriction sur $M^I$) est une \'equivalence. Il s'agit ici du
$\underline{Hom}$ interne pour les $1$-pr\'ecats de Segal.
\end{lemma}
{\em Preuve:}
Une fl\`eche de $L(M)'$ est une classe d'homotopie de morphismes de $M$,
donc tout objet de $\underline{Hom}(I, L(M)')$ est \'equivalent \`a un objet
provenant de $M^I$. Il s'agit donc de prouver que le morphisme en question est
homotopiquement pleinement fid\`ele. Soient $y$ et $z$ des $I$-diagrammes de
$M$, i.e. $y= (f: y_0\rightarrow y_1)$ avec $y_i\in M$ et
$z= (g: z_0\rightarrow z_1)$.
On suppose que $y$ est cofibrant pour celle
des deux structures de Reedy possibles pour laquelle
$f$ est une cofibration.
Si ${\bf z}$ est une r\'esolution simpliciale
fibrante de $z$ alors on a
$$
L(M^I)_{1/}(y,z) \cong M^I_{1/}(y, {\bf z}).
$$
On peut \'ecrire ${\bf z} = ({\bf g}: {\bf z}_0\rightarrow {\bf z}_1)$.

D'autre part,
$$
\underline{Hom}(I, L(M)')_{1/}(y,z)
$$
repr\'esente le foncteur d'ensembles simpliciaux
$$
K\mapsto Hom ^{y,z}(I\times \Upsilon (K), L(M)')
$$
o\`u le terme de droite est le sous-ensemble de morphismes
$r:I\times \Upsilon (K)\rightarrow L(M)'$ avec $r|_{I\times \{ 0\} } = y$
et $r|_{I\times \{ 1\} } = z$. Un morphisme
$$
I\times \Upsilon (K)\rightarrow L(M)'
$$
\'equivaut (par division du carr\'e en deux triangles---le lecteur
est invit\'e \`a
dessiner un carr\'e avec d'un cot\'e $f: y_0\rightarrow y_1$ et de l'autre $g:
z_0\rightarrow z_1$ et avec des fl\`eches \'etiquet\'ees $K$ entre $y_0$ et
$z_0$
et $y_1$ et $z_1$; divis\'e en deux triangles par une fl\`eche \'etiquet\'ee $K$
de $y_0$ vers $z_1$) \`a deux morphismes
$$
\Upsilon ^2(\ast , K)\rightarrow
L(M)' $$ et
$$
\Upsilon ^2(K,\ast )\rightarrow L(M)'
$$
avec la m\^eme restriction \`a l'$\Upsilon (K)$ diagonal. Les
morphismes de restriction des ensembles simpliciaux de ces diagrammes, sur
l'ensemble simplicial des diagrammes pour la diagonale $\Upsilon (K)$, sont
des fibrations. On obtient une formule avec produit fibr\'e homotopique,
i.e. un carr\'e homotopiquement-cart\'esien
$$
\begin{array}{ccc}
\underline{Hom}(I, L(M)')_{1/}(y,z)&\rightarrow & L(M)_{1/}(y_0, z_0)\\
\downarrow && \downarrow \\
L(M)_{1/}(y_1, z_1) & \rightarrow &  L(M)_{1/}(y_0, z_1).
\end{array}
$$
Il y a deux structures de Reedy possibles sur $M^I$ (on peut
choisir $deg(0)>deg(1)$ ou l'inverse);  on prendra la structure ($deg(0) <
deg (1)$) pour laquelle les objets cofibrants sont les cofibrations
$y_0\rightarrow y_1$ entre $y_i$ cofibrants. Dans ce cas et pour une
r\'esolution simpliciale fibrante ${\bf z}$ (dont le composant ${\bf z}_1$
est en
particulier aussi une r\'esolution simpliciale fibrante), le morphisme de
composition avec $f$
$$
M_{1/}(y_1, {\bf z}_1)\rightarrow M_{1/}(y_0, {\bf z}_1)
$$
est une fibration d'ensembles simpliciaux (cf la discussion avant ces lemmes).
Il s'ensuit que  le carr\'e cart\'esien d'ensembles simpliciaux
$$
\begin{array}{ccc}
M^I_{1/}(y,{\bf z}) & \rightarrow &  M_{1/}(y_0,{\bf z}_0)\\
\downarrow && \downarrow \\
M_{1/}(y_1,{\bf z}_1)& \rightarrow & M_{1/}(y_0,{\bf z}_1)
\end{array}
$$
est aussi homotopiquement-cart\'esien. Les trois coins
du bas et de droite co\"{\i}ncident
\`a homotopie pr\`es avec les trois coins correspondants dans le carr\'e
ci-dessus pour les localis\'ees. On en d\'eduit l'\'equivalence
$$
M^I_{1/}(y,{\bf z})\cong \underline{Hom}(I, L(M)')_{1/}(y,z)
$$
qui donne l'\'enonc\'e voulu.
\eop

\subnumero{D\'emonstration du th\'eor\`eme \ref{uneequivalence}}
On note que le morphisme en question est essentiellement surjectif par
\ref{strictif1A}. Il s'agit donc de prouver qu'il est pleinement fid\`ele.

Pour le moment on va supposer que le th\'eor\`eme est connu pour les
cat\'egories {\em directes} ou {\em inverses}, i.e.
les cat\'egories avec fonction
``degr\'e'' o\`u toutes les fl\`eches sauf l'identit\'e sont strictement
monotones pour le degr\'e. Ceci sera justifi\'e ult\'erieurement (cf point (1)
ci-dessous). En fait, on fera une premi\`ere passe de notre d\'emonstration pour
r\'egler ce cas, ensuite la deuxi\`eme passe que nous d\'ecrivons
maintenant---la
raison pour cette contorsion \'etant que pour le cas direct nous avons besoin
d'exactement les m\^emes arguments que pour le cas g\'en\'eral, avec m\^eme
quelques simplifications mais pas suffisamment pour justifier, pour nous qui
sommes
paresseux, de r\'ep\'eter deux fois l'argument; donc, on commence par la version
la plus compliqu\'ee (cette manipulation sera mieux expliqu\'ee en (1) \`a la
fin de la preuve). Nous avons d\'ecouvert cette technique
qui consiste \`a consid\'erer d'abord le cas direct et
ensuite le cas de Reedy dans \cite{Hirschhorn} o\`u
Hirschhorn o\`u prend soin de traiter le cas ``direct'' ind\'ependamment
d'abord.

On proc\'edera par r\'ecurrence sur la longueur de la fonction degr\'e sur $Y$.
On suppose que tous les objets de $Y$ sont de degr\'e $\leq k$, et que le
th\'eor\`eme est d\'emontr\'e pour les cat\'egories de Reedy avec fonction
degr\'e de longueur $\leq k-1$. On pose $Z:= F^{k-1}Y$, en particulier cette
hypoth\`ese s'applique \`a $Z$.

On notera $y_i$ les objets de $Y$ de degr\'e $k$ (l'indice $i$ est dans un
ensemble d'indices que nous ne mettons pas dans les notations). Pour chacun de
ces objets, on note
$$
Y_i:= Latch(y_i)+\{ y _i\} +Match (y_i),
$$
et $Z_i := Y_i \cap Z= Latch(y_i)+Match(y_i)$. On note que $Y_i$ et $Z_i$ sont
des cat\'egories directes (ou alors inverses), donc on peut supposer que le
th\'eor\`eme s'y applique.

Soit $[0,1,2]$ la cat\'egorie associ\'ee \`a
l'ensemble ordonn\'e $0 < 1 < 2$ et
$[0,2]$ sa sous-cat\'egorie pleine avec objets $0$ et $2$.

Soit $\mm (y_i) ^{[0,1,2]}$
(resp. $\mm (y_i) ^{[0,2]}$)
la cat\'egorie de $[0,1,2]$-diagrammes (resp. $[0,2]$-diagrammes)
dans $\mm (y_i)$. On notera les diagrammes $u=(u_0\rightarrow u_1\rightarrow
u_2)$ (ou $u=(u_0\rightarrow u_2)$).

On a le carr\'e cart\'esien suivant de cat\'egories:
$$
\begin{array}{ccc}
Sect (Y_i,\int _{Y_i}\mm ) & \rightarrow & \mm (y_i) ^{[0,1,2]}\\
\downarrow && \downarrow \\
Sect (Z_i,\int _{Z_i}\mm ) & \rightarrow & \mm (y_i) ^{[0,2]}.
\end{array}
$$
Le morphisme du haut envoie une section $\sigma$ vers le diagramme
$$
latch(\sigma ,y_i) \rightarrow \sigma (y_i) \rightarrow match (\sigma , y_i).
$$
Le morphisme du bas envoie une section $\tau$ au-dessus de $Z$, vers le
diagramme
$$
latch(\tau  ,y_i) \rightarrow  match (\tau , y_i).
$$

On notera $Sect (Y_i,\int _{Y_i}\mm )_{c,f}$ la sous-cat\'egorie des sections
dont la partie directe (i.e. la restriction de la section sur $Latch(y_i)+\{
y_i\}$) est cofibrante, et la partie inverse (i.e. la restriction de la section
sur $\{ y_i\} +Match(y_i)$) est fibrante. De m\`eme pour
$Sect (Z_i,\int _{Z_i}\mm )_{c,f}$. On notera $\mm (y_i) ^{[0,1,2]}_{c,f}$
la cat\'egorie des diagrammes $u$ dans lesquels
$u_0$ est cofibrant, $u_0\rightarrow u_1$ est
une cofibration, $u_2$ est fibrant et $u_1\rightarrow u_2$
est une fibration. On
notera $\mm (y_i) ^{[0,2]}_{c,f}$ la cat\'egorie des diagrammes $v$ avec
$v_0$
cofibrant et $v_2$ fibrant.

Le carr\'e ci-dessus donne aussi un carr\'e cart\'esien
$$
\begin{array}{ccc}
Sect (Y_i,\int _{Y_i}\mm ) _{c,f}& \rightarrow & \mm (y_i) ^{[0,1,2]} _{c,f}\\
\downarrow && \downarrow \\
Sect (Z_i,\int _{Z_i}\mm )  _{c,f}& \rightarrow & \mm (y_i) ^{[0,2]} _{c,f}.
\end{array}
$$
Pour chacune de ces cat\'egories on a une notion \'evidente d'\'equivalence
faible---qu'on se dispensera de mettre dans la notation des localis\'ees.
L'avantage
de ce deuxi\`eme carr\'e est que les morphismes horizontaux respectent les
\'equivalences faibles (ce qui n'est pas le cas pour le premier carr\'e,
car les
morphismes horizontaux comportent des limites et colimites).

Maintenant on a aussi le carr\'e cart\'esien de cat\'egories
$$
\begin{array}{ccc}
Sect (Y,\int _{Y}\mm )& \rightarrow & \prod _iSect (Y_i,\int _{Y_i}\mm )
\\ \downarrow && \downarrow \\
Sect (Z,\int _{Z}\mm ) & \rightarrow &
\prod _i Sect (Z_i,\int _{Z_i}\mm ).
\end{array}
$$
Si on note $Sect (Y,\int _{Y}\mm )_{c,f}$ la sous-cat\'egorie des objets
cofibrants et fibrants pour la structure de Reedy
(avec la notation analogue pour
$Z$) alors on obtient le carr\'e cart\'esien
$$
\begin{array}{ccc}
Sect (Y,\int _{Y}\mm ) _{c,f}& \rightarrow & \prod _iSect (Y_i,\int _{Y_i}\mm )
_{c,f}
\\ \downarrow && \downarrow \\
Sect (Z,\int _{Z}\mm )  _{c,f}& \rightarrow &
\prod _i Sect (Z_i,\int _{Z_i}\mm )  _{c,f}.
\end{array}
$$
Ici encore, les morphismes horizontaux respectent les \'equivalences faibles.
Il en d\'ecoule le carr\'e cart\'esien
$$
\begin{array}{ccc}
Sect (Y,\int _{Y}\mm ) _{c,f}& \rightarrow & \prod _i \mm (y_i) ^{[0,1,2]}
_{c,f}\\ \downarrow && \downarrow \\
Sect (Z,\int _{Z}\mm )  _{c,f}& \rightarrow & \prod _i\mm (y_i) ^{[0,2]} _{c,f}.
\end{array}
$$

Soit $N$ l'une des cat\'egories de mod\`eles dans le carr\'e cart\'esien
pr\'ec\'edent (soit l'un des facteurs d'un produit, soit le produit), et soit
$N_{c,f}$ sa sous-cat\'egorie d\'efinie ci-dessus. Sous l'hypoth\`ese des
factorisations fonctorielles, il existe un foncteur $\zeta : N\rightarrow
N_{c,f}$ et une transformation naturelle $t_u: \zeta
(u)\stackrel{\cong}{\rightarrow} u$. Ce foncteur respecte les \'equivalences
faibles et les $t_u$ sont des \'equivalences faibles. En particulier, il
\'etablit une \'equivalence $L(N_{c,f})\cong L(N)$.

Soient $\sigma , \tau \in Sect (Y, \int _Y\mm )_{c,f}$, et notons
$\sigma |_Z$, $\tau |_Z$, $\sigma ^i_{012}$, $\tau ^i_{012}$, et $\sigma
^i_{02}$, $\tau ^i_{02}$ leurs images respectivement dans $Sect (Z,\int _Z\mm
|_Z)_{c,f}$, $\mm (y_i)^{[0,1,2]}_{c,f}$, et $\mm (y_i)^{[0,2]}_{c,f}$. On
pr\'etend que le carr\'e  $(\ast )$
$$
\begin{array}{ccc}
L(Sect (Y, \int _Y\mm )_{c,f})_{1/}(\sigma , \tau )
&\rightarrow & \prod _i L(\mm (y_i)^{[0,1,2]}_{c,f})_{1/}(\sigma
^i_{012},\tau ^i
_{012})\\ \downarrow && \downarrow \\
L(Sect (Z, \int _Z\mm |_Z)_{c,f})_{1/}(\sigma |_Z, \tau |_Z)
& \rightarrow & \prod _i L(\mm (y_i)^{[0,2]}_{c,f})_{1/}(\sigma ^i_{02},\tau
^i_{02}) \end{array}
$$
est homotopiquement-cart\'esien.

Supposons qu'on a montr\'e que $(\ast )$ est homotopiquement-cart\'esien,
et terminons la d\'emonstration du th\'eor\`eme.
On a un carr\'e homotopiquement-cart\'esien
$$
\begin{array}{ccc}
Sect (Y, \int _Y L(\mm ))_{1/}(\sigma , \tau )
&\rightarrow & \prod _i L(\mm (y_i))^{[0,1,2]}_{1/}(\sigma ^i_{012},\tau ^i
_{012})\\ \downarrow && \downarrow \\
Sect (Z, \int _Z L(\mm )|_Z)_{1/}(\sigma |_Z, \tau |_Z)
& \rightarrow & \prod _i L(\mm (y_i))^{[0,2]}_{1/}(\sigma ^i_{02},\tau
^i_{02}), \end{array}
$$
o\`u pour les produits directs on prend d'abord des remplacements de Kan des
ensembles simpliciaux (i.e. on prend les produits directs homotopiques). Pour
justifier ceci, on fait une version ``homotopique'' pour les $L(\mm )$ de la
discussion des diagrammes sur les cat\'egories de Reedy, en utilisant la notion
d'adjoint homotopique du \S 8. Voir aussi dans la d\'emonstration du
th\'eor\`eme \ref{strictif1A}. Nous ne donnons pas plus de d\'etails
sur ce point.

Le morphisme du carr\'e $(\ast )$ dans celui-ci
induit une \'equivalence sur les trois termes
en bas et \`a droite. En effet, pour $Z$ on suppose
connu le th\'eor\`eme par r\'ecurrence, et \`a droite il s'agit de diagrammes
index\'es par des cat\'egories directes $[0,1,2]$ et $[0,2]$---que nous
supposons
d\'ej\`a trait\'ees.
Par le lemme
\ref{equivalence1}, les produits directs dans le carr\'e $(\ast )$ sont des
produits directs homotopiques.

Le fait que les deux carr\'es soient
homotopiquement-cart\'esiens implique alors que le morphisme
$$
L(Sect (Y, \int _Y\mm )_{c,f})_{1/}(\sigma , \tau )
\rightarrow Sect (Y, \int _Y L(\mm ))_{1/}(\sigma , \tau )
$$
est une \'equivalence, ce qui d\'emontre l'\'etape de r\'ecurrence.

Pour finir la d\'emonstration, nous devons: (1) justifier l'hypoth\`ese que
le th\'eor\`eme est connu pour les cat\'egories directes (de longueur finie,
c'est ce qu'on utilise dans la d\'emonstration ci-dessus); (2) montrer que le
carr\'e $(\ast )$ ci-dessus est homotopiquement-cart\'esien; et (3) justifier le
passage \`a la limite sur
le $k$ des $F^kY$ pour une cat\'egorie de Reedy $Y$ avec
fonction degr\'e  non-born\'ee.

Pour (1), supposons maintenant qu'on veut d\'emontrer le th\'eor\`eme pour une
cat\'egorie directe $Y$ (de longueur finie). On proc\`ede de la m\^eme fa\c{c}on
par r\'ecurrence sur la longueur de $Y$, et on recopie les \'etapes de la
d\'emonstration ci-dessus. Le seul changement est que les cat\'egories
$Match(y_i)$ sont vides et donc les termes $match (\sigma , y_i)$
n'apparaissent pas. En particulier, on peut rempla\c{c}er
la cat\'egorie des diagrammes
$\mm (y_i)^{[0,1,2]}$ de longueur $3$, par une cat\'egorie de
diagrammes $\mm (y_i)^{[0,1]}$ de longueur $2$; et la cat\'egorie des diagrammes
$\mm (y_i)^{[0,1]}$ devient juste $\mm (y_i)$. Ceci veut dire qu'\`a la fin
on est ramen\'e \`a consid\'erer le cas de $\mm (y_i)^{[0,1]}$,
autrement dit le cas des diagrammes index\'es par $Y= I=[0,1]$ \`a valeurs dans
une cat\'egorie de mod\`eles ferm\'ee fixe $M=\mm (y_i)$. Ce cas a \'et\'e
trait\'e par le lemme \ref{equivalence2}.

Pour (2), on va encore utiliser la m\'ethode des r\'esolutions. D'abord, on peut
supposer que $\sigma$ est fibrant et cofibrant, ainsi que $\tau$. On choisit
une r\'esolution simpliciale $\tau \rightarrow {\bf t}$ et on peut
supposer que ${\bf t}$ est fibrant et cofibrant en tant qu'objet simplicial
i.e.
$$
{\bf t}\in Sect (Y,\int _Y\mm )^{\Delta ^o}.
$$
Notons ici qu'il s'agit d'it\'erer deux fois l'op\'eration
qui consiste \`a ``prendre la
structure de Reedy'': d'abord on a pris la structure de Reedy sur
$Sect (Y,\int _Y\mm )$ et ensuite, par rapport \`a cette structure, on a pris
la structure de Reedy sur les diagrammes simpliciaux l\`a-dedans.

En particulier,  pour tout $p\in \Delta$, ${\bf t}(p)\in Sect (Y,\int _Y\mm )$
est \`a la fois fibrant et cofibrant, donc on a
$$
{\bf t}\in [Sect (Y,\int _Y\mm )_{c,f}]^{\Delta ^o}.
$$
Les images ${\bf t} |_Z$,  ${\bf t} ^i_{012}$, et  ${\bf t} ^i_{02}$,
qui sont des objets simpliciaux dans les autres cmf $N$ apparaissant ci-dessus,
sont en fait des objets simpliciaux des $N_{c,f}$.

Il faut fixer une structure
de cmf sur $\mm (y_i)^{[0,1,2]}$ et sur $\mm (y_i)^{[0,2]}$. Pour cela, on
utilisera une
structure de Reedy sur $[0,1,2]$ pour laquelle l'objet $1$ est de degr\'e
maximal; donc la fl\`eche $01$ est directe, et $12$ est inverse. Maintenant il
y a un choix \`a faire quant \`a la fl\`eche $02$: nous allons choisir de dire
que c'est une fl\`eche directe, ce qui correspond \`a dire $deg (0)<deg(2)$
(il faut faire ``pencher'' la fonction degr\'e d'un cot\'e ou de l'autre).
La structure correspondante sur $\mm (y_i)^{[0,2]}$ consiste \`a dire
encore que la
fl\`eche $02$ est directe. Avec ce choix, un objet de $\mm (y_i)^{[0,2]}$ est
fibrant si et seulement si ses composantes sont fibrantes, tandis qu'un
objet est
cofibrant si et seulement si ses composantes sont fibrantes et la fl\`eche $02$
est une cofibration. Un objet de $\mm (y_i)^{[0,1,2]}$ est fibrant si et
seulement si
ses composantes ainsi que la fl\`eche $12$ le sont; et un objet est cofibrant si
et seulement si ses composantes ainsi que les fl\`eches $01$ et $02$ le sont.

On a toujours que $\tau ^i_{012}$ et ${\bf t}^i_{012}$ sont des objets fibrants
(de m\^eme pour $\tau ^i_{02}$ et ${\bf t}^i_{02}$). Par contre, $\sigma
^i_{012}$ et $\sigma ^i_{02}$ ne sont plus cofibrantes car on ne peut pas
garantir que la fl\`eche $\sigma ^i_0\rightarrow \sigma ^i_2$ soit une
cofibration.
Pour arranger cela on choisit un remplacement cofibrant
$\tilde{\sigma}^i_{012}\rightarrow  \sigma ^i_{012}$ de la mani\`ere
suivante: on
choisit un diagramme
$$
\begin{array}{ccc}
\sigma ^i_1 & \rightarrow & \tilde{\sigma}^i_2\\
\uparrow & \nearrow & \downarrow \\
\sigma ^i_0 & \rightarrow & \sigma ^i_2
\end{array}
$$
avec une fibration triviale \`a droite, et
une cofibration
en diabonale. Le morphisme de gauche reste une cofibration et on ne dit
rien du morphisme du haut. Le compos\'e $\sigma ^i_1\rightarrow \sigma ^i_2$
devrait \^etre celui qu'on a d\'ej\`a. On pose
$$
\tilde{\sigma}^i_{012}:= [\sigma ^i_0\rightarrow \sigma ^i_1 \rightarrow
\tilde{\sigma}^i_2].
$$
Le carr\'e d'ensembles simpliciaux
$$
\begin{array}{ccc}
\mm (y_i)^{[0,1,2]}(\sigma ^i_{012},{\bf t} ^i
_{012})
&\rightarrow & \mm (y_i)^{[0,1,2]}(\tilde{\sigma} ^i_{012},{\bf t} ^i
_{012})\\ \downarrow && \downarrow \\
\mm (y_i)^{[0,2]}_{1/}(\sigma ^i_{02},{\bf t}
^i_{02})
& \rightarrow & \mm (y_i)^{[0,2]}_{1/}(\tilde{\sigma} ^i_{02},{\bf t}
^i_{02})
\end{array}
$$
est cart\'esien.  De m\^eme le carr\'e
$$
\begin{array}{ccc}
Sect (Y, \int _Y\mm )_{1/}(\sigma , {\bf t} )
&\rightarrow & \prod _i \mm (y)^{[0,1,2]}(\sigma ^i_{012},{\bf t} ^i
_{012})\\ \downarrow && \downarrow \\
Sect (Z, \int _Z\mm |_Z)_{1/}(\sigma |_Z, {\bf t} |_Z)
& \rightarrow & \prod _i \mm (y)^{[0,2]}_{1/}(\sigma ^i_{02},{\bf t}
^i_{02}) \end{array}
$$
est cart\'esien et en composant avec le produit sur $i$ des carr\'es
pr\'ec\'edents
on obtient le carr\'e cart\'esien $(\ast \ast )$ d'ensembles simpliciaux
$$
\begin{array}{ccc}
Sect (Y, \int _Y\mm )_{1/}(\sigma , {\bf t} )
&\rightarrow & \prod _i \mm (y_i)^{[0,1,2]}(\tilde{\sigma} ^i_{012},{\bf t} ^i
_{012})\\ \downarrow && \downarrow \\
Sect (Z, \int _Z\mm |_Z)_{1/}(\sigma |_Z, {\bf t} |_Z)
& \rightarrow & \prod _i \mm (y_i)^{[0,2]}_{1/}(\tilde{\sigma} ^i_{02},{\bf t}
^i_{02}).
\end{array}
$$
Ce carr\'e est \'equivalent au carr\'e $(\ast )$
(pour trouver le morphisme entre carr\'es, notons de fa\c{c}on g\'en\'erale que
si $(C,W)$ est une cat\'egorie avec sous-cat\'egorie ``d'\'equivalences'' et si
$x,y$
en sont des objets, avec un objet simplicial augment\'e $y\rightarrow {\bf y}$
tel que pour tout $p$ le morphisme $y\rightarrow {\bf y}(p)$ soit dans $W$,
alors on obtient un morphisme naturel d'ensembles simpliciaux
$C_{1/}(x,{\bf y})\rightarrow L^H(C,W)_{1/}(x,y)$ o\`u $L^H$ est la
localisation par hamacs \cite{DwyerKan2} \cite{DwyerKan3}; on applique ceci aux
cat\'egories $C=N_{c,f}$ qui apparaissent dans $(\ast )$).

Donc pour (2) il suffit de voir que $(\ast \ast )$ est
homotopiquement-cart\'esien,
donc il suffit de voir que pour tout $i$ le morphisme
$$
\mm (y_i)^{[0,1,2]}(\tilde{\sigma} ^i_{012},{\bf t} ^i
_{012})\rightarrow
\mm (y_i)^{[0,2]}_{1/}(\tilde{\sigma} ^i_{02},{\bf t}
^i_{02})
$$
est une fibration de Kan.

Pour ceci
on va simplifier la notation: on prend une cfm $M$ et deux diagrammes
$x,y\in M^{[0,1,2]}$ plus une r\'esolution simpliciale fibrante $y\rightarrow
{\bf y}$. Pour ceci on munit (comme avant) $M^{[0,1,2]}$ d'une structure de cmf
telle que les objets fibrants soient les diagrammes d'objets fibrants avec le
deuxi\`eme morphisme $12$ fibrant, et les objets cofibrants sont les diagrammes
d'objets cofibrants avec les morphismes $01$ et $02$ cofibrants. Pareillement on
munit $M^{[0,2]}$ de la structure pour laquelle les objets cofibrants sont
les diagrammes d'objets cofibrants avec une cofibration en $02$.

Notons $i$ le morphisme $ [0,2]\hookrightarrow [0,1,2]$; on a le foncteur
$$
i^{\ast}: M^{[0,1,2]}\rightarrow M^{[0,2]}
$$
et son adjoint \`a gauche $i_!$. Cet adjoint a la description concr\`ete
suivante. Si $u= (u_0\rightarrow u_2)\in M^{[0,2]}$ alors
$$
i_!(u)= [u_0 \stackrel{=}{\rightarrow} u_0 \rightarrow u_2]\in M^{[0,1,2]}.
$$
Si $x$ est un objet cofibrant de $M^{[0,1,2]}$ alors
le morphisme d'adjonction $i_!i^{\ast}x\rightarrow x$ est
$$
[x_0 \stackrel{=}{\rightarrow} x_0 \rightarrow x_2]
\rightarrow [x_0\rightarrow x_1\rightarrow x_2]
$$
qui est une cofibration.

On a
$$
M^{[0,2]}_{1/}(i^{\ast} x, i^{\ast} {\bf y})=
M^{[0,1,2]}_{1/}(i_!i^{\ast} x, {\bf y}).
$$
Le morphisme en question est donc le morphisme induit par
$i_!i^{\ast}x\rightarrow x$,
$$
M^{[0,1,2]}_{1/}(x, {\bf y})
\rightarrow
M^{[0,1,2]}_{1/}(i_!i^{\ast} x, {\bf y}).
$$
Si $x$ est cofibrant alors la cofibration $i_!i^{\ast}x\rightarrow x$
induit ici une fibration de Kan d'ensembles simpliciaux,
ce qui montre que $(\ast \ast )$ est homotopiquement-cart\'esien,
d'o\`u la m\^eme chose
pour $(\ast )$ ce qui donne (2).

Au passage, l'argument ci-dessus montre que si $\sigma$ et $\tau$ (resp. ${\bf
t}$) sont des sections cofibrantes et fibrantes (resp.
est une r\'esolution simpliciale
fibrante et cofibrante) alors le morphisme
$$
Sect (Y, \int _Y\mm )_{1/}(\sigma , {\bf t} )\rightarrow
Sect (Z, \int _Z\mm |_Z)_{1/}(\sigma , {\bf t} )
$$
est une fibration de Kan.

Pour (3) on fixe maintenant une cat\'egorie de Reedy $Y$ avec fonction
degr\'e non-n\'ecessairement born\'ee. On fixe un remplacement fibrant
$$
\int '_YL(\mm )\rightarrow Y
$$
ce qui d\'etermine des remplacements fibrants par restriction
$$
\int '_{F^kY} L(\mm )|_{F^kY}\rightarrow F^kY.
$$
Dans ces conditions, $Sect (Y, \int '_YL(\mm ))$ est la limite de la suite
$$
\ldots \rightarrow Sect (F^kY, \int '_{F^kY}L(\mm )|_{F^kY})
\rightarrow Sect (F^{k-1}Y, \int '_{F^{k-1}Y}L(\mm
)|_{F^{k-1}Y})\rightarrow \ldots ,
$$
et les morphismes de transition dans cette suite sont des fibrations
de $1$-cat\'egories de Segal. D'autre part, fixons des sections $\sigma$ et
$\tau$ dans $Sect (Y, \int _Y\mm )$, avec $\sigma$ cofibrante, et fixons une
r\'esolution simpliciale fibrante $\tau \rightarrow {\bf t}$. On a que
$$
Sect (Y, \int _Y\mm )_{1/}(\sigma , {\bf t})
$$
est la limite de la suite d'ensembles simpliciaux
$$
\ldots \rightarrow Sect (F^kY, \int _{F^kY}\mm |_{F^kY})_{1/}(\sigma
|_{F^kY}, {\bf t}|_{F^kY})
\rightarrow \ldots .
$$
Par l'argument donn\'e pour (2), les morphismes de transition
de cette suite sont des
fibrations de Kan. Cette suite est, terme \`a terme \'equivalente
(par le m\^eme argument) \`a la suite
$$
\ldots \rightarrow Sect (F^kY, \int '_{F^kY}L(\mm )|_{F^kY})_{1/}(\sigma
|_{F^kY}, \tau |_{F^kY}
) \rightarrow \ldots ,
$$
et ces \'equivalences etant fonctorielles on obtient une \'equivalence
de diagrammes entre les deux suites (en fait, pour bien avoir un morphisme ici
il conviendrait d'utiliser la localisation par hamacs $L^H$ au lieu de $L$). Le
fait que les morphismes de transition dans les deux cas soient des fibrations de
Kan implique que le morphisme sur les limites est une \'equivalence
$$
Sect (Y, \int _Y\mm )_{1/}(\sigma , {\bf t})\stackrel{\cong}{\rightarrow}
Sect (Y, \int '_YL(\mm ))_{1/}(\sigma , \tau ).
$$
Or
$$
Sect (Y, \int _Y\mm )_{1/}(\sigma , {\bf t}) \cong
L(Sect (Y, \int _Y\mm ))_{1/}(\sigma , \tau )
$$
d'o\`u la partie (3).

Ceci termine la d\'emonstration du th\'eor\`eme \ref{uneequivalence}.
\eop

\numero{La descente pour les pr\'efaisceaux de Quillen \`a gauche}

\label{quidescentepage}

On va prouver un th\'eor\`eme de descente pour un pr\'efaisceau de Quillen
\`a gauche, qui donnera des conditions pour que $L(\mm )$ soit un champ.
Nous avons trouv\'e cete d\'emonstration \`a partir du cas des complexes.
Le probl\`eme principal pour ce r\'esultat est qu'on doit partir d'une
section faible de $L(\mm )$ qui peut \^etre, {\em a priori}, assez ``sauvage''.
Dans les chapitres pr\'ec\'edents nous
avons trouv\'e une s\'erie de r\'esultats qui permettent de ``domestiquer'' une
telle section, i.e. de
se ramener au cas d'une section de $\int \mm$ au lieu d'une
section de $\int 'L(\mm )$.

Le travail maintenant consiste \`a d\'emontrer qu'une
section de $\int \mm$ descend sur l'objet de base. Cette derni\`ere partie
est en
r\'ealit\'e d\'ej\`a bien connue, par exemple on sait prendre le complexe simple
associ\'e \`a un complexe cosimplicial, voir Deligne \cite{HodgeIII}.
On conjecture que c'est ce probl\`eme qui a \'et\'e abord\'e
dans les ``papiers secrets'' de Deligne auxquels Illusie
fait r\'ef\'erence dans
\cite{Illusie}.
C'est aussi le sujet de l'expos\'e de B. Saint-Donat sur la {\em descente
cohomologique} dans SGA 4 (\cite{SGA4}). La m\'ethode dans tous les
cas est de commencer
avec une section $\sigma$ dans $Sect ^{\rm eq} (\Delta , \int _{\Delta} \mm
)$ (dans SGA 4 cet
ensemble de sections s'appelle $\Gamma ^{\rm cocart}(\, \, )$) et de la
``descendre'' via
l'augmentation, en appliquant le foncteur ``image directe'' que nous noterons
$\varphi _{\ast}$.
Le probl\`eme est de v\'erifier que cela r\'epond bien au probl\`eme de
descente, i.e. que
$\varphi ^{\ast} (\varphi _{\ast} \sigma )$ est \'equivalente \`a $\sigma$.

Voir la fin de ce chapitre pour une comparaison
plus d\'etaill\'ee entre notre discussion et celle de SGA 4.

Au d\'ebut de ce num\'ero on va \'etudier un pr\'efaisceau de Quillen \`a gauche
$\mm$ sur $(\Delta ^+)^o$. Si $p\rightarrow q$ est un morphisme de $\Delta^+$
alors le foncteur de restriction (i.e. le foncteur de Quillen \`a gauche) va
dans le sens $\mm (p)\rightarrow \mm (q)$. Ici  $\Delta ^+$ est la cat\'egorie
$\Delta$ augment\'ee par un objet initial qu'on notera $\iota$. Sauf en cas de
confusion possible, on \'evitera la notation $\mm |_{\Delta ^o}$ et
on d\'esignera aussi
cette restriction par $\mm$.

On fera l'hypoth\`ese (o) suivante:
\newline
$(o)'$\,\, Que  chaque $\mm (y)$  admet des limites et colimites
petites arbitraires, et factorisations fonctorielles;
\newline
$(o)''$\,\,
que tout objet de $\mm (y)$ est cofibrant; et
\newline
$(o)'''$\,\, que $Sect (\Delta ^+, \int
_{\Delta ^+}\mm )$ et $Sect (\Delta , \int _{\Delta} \mm )$ admettent des
structures de cmf de type (II)
(du Th\'eor\`eme \ref{sectionscmf}) engendr\'ees par
cofibrations.

L'hypoth\`ese $(o)''$ que tout objet est cofibrant est l\`a par pure
commodit\'e et \'evidemment inessentielle. La structure de type (II)
de l'hypoth\`ese $(o)'''$ pourrait probablement \^etre remplac\'ee, dans
notre argument, par  une structure de Reedy. Par contre, l'hypoth\`ese $(o)'$
(qui par ailleurs est devenue standard
cf \cite{Hirschhorn} \cite{DHK} \cite{HoveyBook})
et en particulier la fermeture de $\mm (\iota )$ par petites limites---par
exemple, par limites index\'ees par $\Delta$---est essentielle et le
th\'eor\`eme
\ref{dansleschamps} ci-dessous ne serait probablement plus vrai sans
cette condition.

On munit $Sect (\Delta , \int _{\Delta} \mm )$ de sa
structure de cmf de type (II) de l'hypoth\`ese $(o)''$.
Le foncteur ``section
constant''
$$
\varphi ^{\ast} : \mm (\iota )\rightarrow Sect (\Delta , \int _{\Delta} \mm )
$$
est alors un foncteur de Quillen \`a gauche, et son adjoint \`a droite
qu'on notera $\varphi _{\ast}$ (qui existe d'apr\`es l'hypoth\`ese $(o)'$)
est un foncteur de Quillen \`a droite.

On peut d\'ecomposer $\varphi _{\ast}$ de la mani\`ere suivante. L'adjoint des
restrictions fournit un morphisme
$$
r_{\ast}:Sect (\Delta , \int _{\Delta} \mm )\rightarrow \mm (\iota )^{\Delta}.
$$
En composant
ensuite avec le foncteur $lim : \mm (\iota )^{\Delta} \rightarrow \mm (\iota
)$  on obtient
$$
\varphi _{\ast} : Sect (\Delta , \int _{\Delta} \mm )\rightarrow \mm (\iota ),
$$
$$
\varphi _{\ast} (\sigma ) := \lim _{\Delta} r_{\ast}\sigma .
$$

On note qu'un objet fibrant pour la structure de type (II) est aussi fibrant
pour la structure de Reedy (qui existe automatiquement sur les sections).
En outre $r$ transforme les objets fibrants en objets fibrants
(c'est un foncteur de Quillen \`a droite car compos\'e d'adjoints de
restrictions).  Donc, pour $\sigma \in Sect (\Delta , \int _{\Delta} \mm )_f$,
$\varphi _{\ast} (\sigma )$ est \'equivalent \`a  $holim _{\Delta}r_{\ast}\sigma $.

On veut montrer que $\varphi ^{\ast}$ induit une \'equivalence
$$
L(\varphi ^{\ast} ): L(\mm (\iota ))\cong
L(Sect ^{\rm eq}(\Delta , \int _{\Delta} \mm )).
$$
On commence par l'observation suivante. On garde pour ce lemme l'hypoth\`ese de
commodit\'e que tous les objets sont cofibrants, mais elle n'est
certainement pas essentielle.

\begin{lemme}
\label{dpq1}
Soient $N$ et $M$ des cat\'egories de mod\`eles ferm\'ees et
$\varphi _{\ast} : N\rightarrow M$, $\varphi ^{\ast} : M\rightarrow N$ une
paire de foncteurs
adjoints, o\`u $\varphi ^{\ast}$ est un foncteur de Quillen \`a gauche et
$\varphi _{\ast}$ un
foncteur de Quillen \`a droite. Supposons que tous les objets de $M$ et
de $N$ sont cofibrants. Supposons que pour tout objet $x\in M$, le morphisme
d'adjonction $x\rightarrow \varphi _{\ast} (\varphi ^{\ast} x)'$ est une
\'equivalence (o\`u $(\varphi ^{\ast}
x)'$ est le remplacement fibrant de $\varphi ^{\ast} x$).
\newline Alors le foncteur
$$
L(\varphi ^{\ast} ) : L(M)\rightarrow L(N)
$$
est pleinement fid\`ele, avec image essentielle
constitu\'e des objets  fibrants $y\in N$
tels que $\varphi ^{\ast} (\varphi _{\ast} y)\rightarrow y$ soit une
\'equivalence, et avec
$L(\varphi _{\ast} )$ pour inverse sur cette image.
\end{lemme}
{\em Preuve:}
Soit $N'\subset N_f$ la sous-cat\'egorie des objets
fibrants $y$ tels que $\varphi ^{\ast} (\varphi _{\ast} y)\rightarrow y$
soit une \'equivalence. On
note que $N'$ est stable par \'equivalence dans $N_f$ (i.e. si $z\in N_f$
et $y\in N'$ avec $z\cong y$ alors $z\in N'$). Par le principe \ref{stabilite},
$(N',W_{N_f}\cap N')$ admet un calcul de fractions homotopique, et la
localis\'ee
$L(N',W_{N_f}\cap N')$ est une sous-cat\'egorie simpliciale pleine de
$L(N_f,W_{N_f})\cong L(N,W_N)$. Par  \cite{DwyerKan2} Corollary 3.6, les
foncteurs $\varphi _{\ast}$ et $\varphi ^{\ast}$ \'etablissent une
\'equivalence entre $L(M, W_M)$
et $L(N', W_N\cap N')$.
\eop

Retournons \`a la situation pr\'ec\'edente. On va appliquer le lemme \`a
la paire de
foncteurs $\varphi _{\ast} , \varphi ^{\ast}$ d\'efinis pr\'ec\'edemment.
Si $\sigma \in
Sect (\Delta , \int _{\Delta} \mm )$ est \'equivalente \`a $\varphi ^{\ast}
(x)$ pour $x\in
\mm (\iota )$ alors les morphismes $\sigma ^r(f)$ sont des \'equivalences, i.e.
$\sigma$ est dans $Sect ^{\rm eq} (\Delta , \int _{\Delta} \mm )$.

On va conserver l'hypoth\`ese (o) ci-dessus.
En plus, on va supposer qu'on a les deux propri\'et\'es
suivantes: \newline
(i)\,\, un morphisme $a$ de $\mm (\iota )$ est une \'equivalence
si et seulement si $\varphi ^{\ast} (a)$ est une \'e\-qui\-va\-len\-ce; et
\newline
(ii)\,\, pour tout $\sigma \in Sect ^{\rm eq} (\Delta , \int _{\Delta} \mm )_f$,
le morphisme d'adjonction $\varphi ^{\ast} \varphi _{\ast} (\sigma )
\rightarrow \sigma$
est une \'equivalence.

On peut remarquer que la propri\'et\'e (i) est \'equivalente
\`a la condition:
\newline
(i)'\,\, un morphisme $a$ de $\mm (\iota )$ est
une \'equivalence si et seulement si
sa restriction \`a $\mm (0)$ en est une.

Alors pour $x\in \mm (\iota )$ la propri\'et\'e (ii) s'applique \`a
$\varphi ^{\ast} (x)'$
(le remplacement fibrant de $\varphi ^{\ast} (x)$),
donc le morphisme
$\varphi ^{\ast} \varphi _{\ast} (\varphi ^{\ast} (x)')\rightarrow \varphi
^{\ast} (x)'$ est une \'equivalence. Le
compos\'e
$$
\varphi ^{\ast} (x)\rightarrow \varphi ^{\ast}
 \varphi _{\ast} (\varphi ^{\ast} (x)') \rightarrow \varphi ^{\ast} (x)'
$$
est \'egal \`a l'\'equivalence faible $\varphi ^{\ast} (x)\rightarrow
\varphi ^{\ast} (x)'$ du
remplacement fibrant. Le premier morphisme est donc une \'equivalence faible,
mais ce morphisme est $\varphi ^{\ast} (a)$ o\`u $a$ est le morphisme
d'adjonction
$a:x\rightarrow \varphi _{\ast} \varphi ^{\ast} (x)$. L'hypoth\`ese du lemme
est donc v\'erifi\'ee; il s'ensuit que
$L(\varphi ^{\ast} )$ est pleinement fid\`ele. D'autre part,
la condition (ii)
ci-dessus
identifie l'image de $L(\varphi ^{\ast} )$ comme la localis\'ee de  $Sect
^{\rm eq} (\Delta
, \int _{\Delta} \mm )$.

On note comme plus haut que (i) est \'equivalente \`a la condition
selon laquelle
un morphisme $a$ de $\mm (\iota )$ est une \'equivalence
si et seulement si sa restriction \`a $\mm (0)$ en est une.
On a obtenu le

\begin{corollaire}
\label{dpq2}
Soit $\mm$ un pr\'efaisceau de Quillen \`a gauche sur $(\Delta ^+)^o$
qui satisfait les conditions ci-dessus qu'on rappelle:
\newline
(o)\,\, chaque $\mm (y)$ admet des petites limites et colimites
arbitraires, et des factorisations fonctorielles;
tout objet de $\mm (y)$ est cofibrant; et $Sect (\Delta ^+, \int
_{\Delta ^+}\mm )$ et $Sect (\Delta , \int _{\Delta} \mm )$ admettent des
structures de cmf de type (II) (du
Th\'eor\`eme \ref{sectionscmf}) engendr\'ees par
cofibrations.
\newline
(i)\,\, un morphisme $a$ de $\mm (\iota )$ est une \'equivalence
si et seulement si sa restriction \`a $\mm (0)$ en est une; et
\newline
(ii)\,\, pour tout $\sigma \in Sect ^{\rm eq} (\Delta , \int _{\Delta} \mm )_f$,
le morphisme d'adjonction $\varphi ^{\ast} \varphi _{\ast} (\sigma )
\rightarrow \sigma$
est une \'equivalence.
\newline
Alors le morphisme
$$
L(\mm (\iota ))\rightarrow
L(Sect ^{\rm eq} (\Delta , \int _{\Delta} \mm ))
$$
est une \'equivalence de cat\'egories simpliciales.
\end{corollaire}
\eop

\bigskip

On aborde maintenant le cas d'un pr\'efaisceau de Quillen
sur un site. Soit $\Xx$ un site admettant des produits fibr\'es et suffisamment
de sommes disjointes (voir \S 15) et soit $\mm$ un pr\'efaisceau de Quillen \`a
gauche sur $\Xx$.
Si $f:X\rightarrow Y$ est un
morphisme de $\Xx$ alors on note
$$
f^{\ast}: \mm (Y)\rightarrow \mm (X)
$$
la restriction $res_f$, et $f_{\ast}$ son adjoint \`a droite.

On fait l'hypoth\`ese que $L(\mm )$ est compatible aux sommes disjointes
\ref{compsomdis}. On pr\'esente d'abord un crit\`ere qui permet de v\'erifier
cette hypoth\`ese.
Si $\Uu = \{ U_{\alpha}\}$ est une famille d'objets de $\Xx$, on pose
$$
\mm (\Uu ):= \prod _{\alpha} \mm (U_{\alpha}).
$$
Ceci a une structure de cmf.
Le morphisme de restriction fournit un morphisme
$$
r_{\Uu}^{\ast}:\mm (\amalg \Uu )\rightarrow \mm (\Uu ),
$$
qui est un
foncteur de Quillen \`a gauche,
dont on note $r_{\Uu , \ast}$ l'adjoint \`a droite.

\begin{corollaire}
\label{compsomcritere}
On suppose que, pour toute famille $\Uu$ d'objets de $\Xx$ (de taille $<\beta$),
et pour tout $u$ dans $\mm (\Uu )$, le morphisme d'adjonction
$$
r_{\Uu}^{\ast} r_{\Uu , \ast} (u)\rightarrow u
$$
est une \'equivalence faible. On suppose en outre que, pour toute famille $\Uu$
comme plus haut, un
morphisme $a$ de $\mm (\amalg \Uu )$ est une \'equivalence faible
si et seulement si $r_{\Uu}^{\ast} (a)$ en est une. Alors
$L(\mm )$ est compatible aux sommes disjointes.
\end{corollaire}
{\em Preuve:} Ceci est une caract\'erisation bien connue des {\em \'equivalences
de Quillen}. Par ailleurs, c'est un corollaire du lemme \ref{dpq1}.
\eop

L'\'enonc\'e du th\'eor\`eme suivant, d\'ej\`a corrig\'e dans la 
version 2 du papier par l'addition de  l'hypoth\`ese (4), est corrig\'e dans la
pr\'esente version 3 par l'addition de 
l'hypoth\`ese (5).

\begin{theoreme}
\label{dansleschamps}
Supposons
que $\Xx$ est un site qui admet des produits fibr\'es et suffisamment de sommes
disjointes compatibles aux produits fibr\'es.
Soit
$\mm$ un pr\'efaisceau de Quillen sur $\Xx$.  On suppose que $\mm$ satisfait les
trois propri\'et\'es suivantes: \newline
(0)\,\, chaque $\mm (y)$ admet des petites limites et colimites
arbitraires, et des factorisations fonctorielles;
tout objet de $\mm (y)$ est cofibrant; et
pour tout foncteur $\Yy \rightarrow \Xx$ de $1$-cat\'egories,
$Sect (\Yy , \int _{\Yy }\mm |_{\Yy})$ admet une structure de cmf de
type (II) engendr\'ee par cofibrations.
\newline
(1)\,\, Pour tout $X$ de $\Xx$ et
toute famille
$\Uu =\{ U_{\alpha}\stackrel{p_{\alpha}}{\rightarrow}
X\}$ couvrant $X$ (ou v\'erifiant $X=\amalg \Uu$)
un morphisme $a$ de $\mm (X)$ est une \'equivalence faible
si et seulement si sa restriction  $p_{\alpha}^{\ast}(a)$ \`a $\mm
(U_{\alpha})$  est une \'equivalence faible pour tout $\alpha$;
\newline
(2)\,\, si $X= \amalg \Uu$ est la somme disjointe d'une famille $\Uu =\{
U_{\alpha}\}$ alors le morphisme d'adjonction
$$
r_{\Uu}^{\ast} r_{\Uu , \ast}u\rightarrow u
$$
est une \'equivalence faible pour tout $u\in  \mm (\Uu ):= \prod _{\alpha} \mm
(U_{\alpha})$; 
\newline
(3)\,\, le pr\'efaisceau $\mm$ est {\em cart\'esien} en ce sens que pour
tout diagramme cart\'esien d'objets de $\Xx$
$$
\begin{array}{ccc}
Y\times _XZ & \stackrel{q}{\rightarrow} & Z \\
{\scriptstyle p}\downarrow && \downarrow {\scriptstyle r}\\
Y& \stackrel{s}{\rightarrow} & X
\end{array}
$$
et pour tout $u$ dans $\mm (Y)$.
l'application naturelle
$$
r^{\ast} s_{\ast} (u) \stackrel{\cong}{\rightarrow} q_{\ast} p^{\ast}(u).
$$
est une \'equivalence
faible; 
\newline
(4)\,\, les morphismes de restriction du pr\'echamp localis\'e $L(\mm )$ 
pr\'eservent les limites  homotopiques (on note que ces morphismes pr\'eservent
automatiquement les colimites homotopiques car 
$\mm$ est de Quillen \`a gauche); et
\newline
(5) \,\, le pr\'echamp de Segal $L(\mm )$ est  un protochamp, c'est-\`a-dire
que pour tout $x\in \Xx$ et tous $m,m' \in L(\mm )(x)$ le pr\'efaisceau
simplicial $L(\mm )_{1/}(m,m')$  est un champ sur $\Xx / x$.

Alors le pr\'efaisceau de cat\'egories simpliciales $L(\mm )$ est un champ (i.e.
$1$-champ de Segal) sur $\Xx$.
\end{theoreme}
{\em Preuve:}
On montre d'abord \`a l'aide du corollaire \ref{compsomcritere}
que $L(\mm )$ est compatible aux sommes disjointes. Si $X=\amalg U_{\alpha}$
est la somme disjointe de la famille $\Uu = \{ U_{\alpha}\}$ alors
on a le foncteur
$$
r_{\Uu}^{\ast}: \mm (X)\rightarrow \prod _{\alpha} \mm (U_{\alpha})
$$
et son adjoint \`a droite $r_{\Uu}^{\ast}$. La condition (1) plus le fait
que $\Uu$ est une famille couvrant $X$ impliquent qu'une fl\`eche $a$ de
$\mm (X)$ est une \'equivalence faible si et seulement si $r^{\ast}_{\Uu}(a)$
est une \'equivalence faible. La condition (2) donne l'autre  partie de
l'hypoth\`ese du corollaire \ref{compsomcritere},
et on obtient donc que $L(\mm )$
est compatible aux sommes disjointes.

On va appliquer \ref{critere7} \`a $L(\mm )$. On peut donc fixer $X\in \Xx$ et
une famille couvrant $\Uu$ \`a un seul \'el\'ement
$p:U\rightarrow X$.
On obtient le foncteur
$\rho ^+\Uu: \Delta ^+ \rightarrow \Xx
^o$ d\'efini par
$$
\rho ^+(\Uu )(p) := U \times _X \ldots \times _XU \;\;\; (p+1 \;\;
\mbox{fois}).
$$
En particulier
on a $\rho ^+ (\Uu )(\iota ) = X$. On a le pr\'efaisceau de Quillen \`a
gauche $\rho ^+ (\Uu )^{\ast} \mm$ sur $(\Delta ^+)^o$. On appellera
$\rho (\Uu )^{\ast} \mm$ sa restriction \`a $\Delta$. Par le corollaire
\ref{pqcalclim} (c'est ici qu'on strictifie nos donn\'ees de descente) on a
$$
\lim _{\leftarrow , \Delta} \rho (\Uu )^{\ast}L(\mm )
\cong L(Sect ^{\rm eq}(\Delta , \int _{\Delta}\rho (\Uu )^{\ast}\mm )).
$$
Pour prouver que $L(\mm )$ est un champ, il suffit, d'apr\`es la proposition
\ref{critere7} ((a) + (b)), de prouver (b) que le morphisme de restriction
$$
L(\mm (X))\rightarrow \lim _{\leftarrow , \Delta} \rho (\Uu )^{\ast}L(\mm )
$$
est essentiellement surjectif. 
En effet la condition (a) de \ref{critere7}
est l'hypoth\`ese (5) du pr\'esent \'enonc\'e.
Donc il suffit de prouver que
$$
L(\mm (X))\rightarrow
L(Sect ^{\rm eq}(\Delta , \int _{\Delta}\rho (\Uu )^{\ast}\mm ))
$$
est une \'equivalence.
Comme on a $M(X)=\rho (\Uu )^{\ast}\mm (\iota )$, la question ne concerne
que la restriction $\rho ^+(\Uu )^{\ast}\mm$ \`a $\Delta ^+$, et on peut
appliquer le corollaire \ref{dpq2} ci-dessus.
On munit
$$
Sect (\Delta , \int _{\Delta}\rho (\Uu )^{\ast}\mm )
$$
de sa structure de cmf de type (II) garantie par l'hypoth\`ese (0).

Les hypoth\`eses (0) et  (1)  impliquent
imm\'ediatement les conditions (o) et (i) du corollaire \ref{dpq2}.

Pour l'hypoth\`ese (ii) de \ref{dpq2}, soit
$$
\sigma \in
Sect ^{\rm eq}(\Delta , \int _{\Delta}\rho (\Uu )^{\ast}\mm )_f
$$
(c'est notre donn\'ee de descente). On va la descendre en $\varphi _{\ast}
(\sigma )$.
On va montrer que le morphisme d'adjonction
$$
\varphi ^{\ast} \varphi _{\ast} (\sigma )\rightarrow \sigma
$$
est une \'equivalence, apr\`es quoi \ref{dpq2} s'appliquera pour donner le
r\'esultat voulu.

Au-dessus de $(\Delta ^+\times \Delta ^+)^o$ on a le pr\'efaisceau
de Quillen \`a gauche $\nn$ d\'efini par
$$
\nn (a,b) := \mm (\rho ^+ (\Uu )(a) \times _X \rho ^+ (\Uu )(b)).
$$
Il y a un foncteur
$$
\gamma : \Delta ^+ \times \Delta ^+ \rightarrow
\Delta ^+
$$
de concat\'enation d'ensembles ordonn\'es ($\iota$ est l'ensemble vide),
et on a
$$
\nn = \gamma ^{\ast} (\rho (\Uu )^{\ast}\mm ).
$$
Avec la section $\sigma$ on obtient trois sections de $\int _{\Delta \times
\Delta} \nn$, not\'ees
$$
p_1^{\ast}\sigma , \;\;\;
p_2^{\ast}\sigma , \;\;\;
\gamma ^{\ast}\sigma .
$$
La troisi\`eme $\gamma ^{\ast}\sigma$ est juste la remont\'ee de $\sigma$ via
$\gamma$; et par exemple $p_1^{\ast}(\sigma )(a,b)$ est l'image de
$\sigma (\rho (\Uu )(a))$ par le morphisme
$$
\mm (\rho (\Uu )(a))\rightarrow \nn (a,b) = \mm (\rho (\Uu )(a)\times
_X\rho (\Uu
) (b))
$$
de restriction via la premi\`ere projection du produit.
On a des morphismes
$$
p_1^{\ast}\sigma \rightarrow \gamma ^{\ast} \sigma \leftarrow p_2^{\ast} \sigma
$$
et notre hypoth\`ese que $\sigma$ est une $eq$-section implique que ces
morphismes sont des \'e\-qui\-va\-len\-ces faibles (objet par objet au-dessus de
$\Delta
\times \Delta$).

On munit
$Sect (\Delta \times \Delta ,
\int _{\Delta \times \Delta}\nn )$
d'une structure de cmf qu'on appellera {\em de type (I-II)} (ce qui veut dire
{\em grosso modo} ``type (I) en la premi\`ere variable et type (II) en la
deuxi\`eme'') de la fa\c{c}on suivante. Pour $p\in \Delta$ on munit
$$
Sect (\{ p\} \times \Delta   ,
\int _{\{ p\} \times \Delta   }\nn |_{\{ p\} \times \Delta })
$$
de sa structure de type (II) engendr\'ee par cofibrations (d'apr\`es
l'hypoth\`ese (0)). Ces cmf forment un pr\'efaisceau de Quillen \`a gauche sur
$\Delta$ dont chaque valeur est engendr\'ee par cofibrations. La cat\'egorie des
sections de ce pr\'efaisceau de Quillen---qu'on munit de sa structure de type
(I) du th\'eor\`eme \ref{sectionscmf}---est exactement $Sect (\Delta \times
\Delta , \int _{\Delta \times \Delta}\nn )$.

Dans cette structure de type (I-II), les cofibrations sont les morphismes qui,
restreints \`a chaque $\Delta \times \{ p\}$, sont des cofibrations pour la
structure de type (I) (i.e. de type HBKQ). Par contre, les
fibrations sont les morphismes qui, restreints \`a chaque $\{ p\} \times
\Delta$, sont des fibrations pour la structure de type (II) (i.e. de type
Jardine-Brown-Heller).

Avec ces structures, le foncteur (``remont\'e par la premi\`ere projection'')
$$
p_1^{\ast} : Sect (\Delta , \int _{\Delta} \rho (\Uu )^{\ast} \mm )
\rightarrow
Sect (\Delta \times \Delta ,
\int _{\Delta \times \Delta}\nn )
$$
pr\'eserve les cofibrations, i.e. c'est un foncteur de Quillen \`a gauche.
Son adjoint \`a droite, qui existe d'apr\`es (0), et qu'on notera
$$
p _{1,\ast}: Sect (\Delta \times \Delta ,
\int _{\Delta \times \Delta}\nn )\rightarrow
Sect (\Delta , \int _{\Delta} \rho (\Uu )^{\ast} \mm )
$$
est un foncteur de Quillen \`a droite. C'est
l'analogue de $\varphi _{\ast}$ pour la situation relative \`a la
premi\`ere projection.


On va prouver l'assertion suivante: 
le morphisme
naturel
$$
\varphi ^{\ast} \varphi _{\ast} (\sigma )\stackrel{\cong}{\rightarrow}
p _{1,\ast}p_2^{\ast} (\sigma )
$$
est une \'equivalence. Pour cela, on va utiliser
 les conditions (3) et (4) du th\'eor\`eme.
La phrase correspondante \`a cette assertion dans la premi\`ere version
du papier \'etait incorrecte; c'est pour corriger l'argument qu'il 
a fallu rajouter
l'hypoth\`ese (4).  Pour la preuve, on commence par rappeler que 
$$
\varphi _{\ast} ( \sigma ) = \lim \, r_{\ast} (\sigma )
$$
o\`u 
$$
r_{\ast} (\sigma ) \in \mm (X) ^{\Delta} .
$$
La limite dont il s'agit est homotopique (cf la discussion au d\'ebut de ce 
chapitre).

On a une fl\^eche 
$$
\ell : \varphi ^{\ast} \lim _{\Delta}\, r_{\ast} (\sigma ) \rightarrow 
holim _{\Delta} \varphi ^{\ast} r_{\ast} (\sigma )
$$ 
et la condition (4) du th\'eor\`eme signifie que $\ell $ 
est  une \'equivalence. Ici, 
$$ 
\varphi ^{\ast}r_{\ast} (\sigma ) \in 
Sect (\Delta , \int _{\Delta}\rho (\Uu )^{\ast} \mm )^{\Delta},
$$
et l'occurence de $\Delta$ suivant laquelle on prend le $holim$ est celle qui
figure en exposant
(et qui correspond \`a la  variable not\'ee $b$ plus haut).

D'autre part, la condition (3) du th\'eor\`eme signifie que la 
fl\^eche naturelle
$$
\varphi ^{\ast}r_{\ast} (\sigma ) \rightarrow 
r_{1,\ast} p_2^{\ast} (\sigma )
$$
est une \`equivalence objet-par-objet au-dessus de 
$\Delta \times \Delta$, o\`u
la notation $r_{1,\ast}$ d\'esigne le morphisme induit par
les adjoints des restrictions
$$
r_{1,\ast} : 
Sect (\Delta \times \Delta , \int _{\Delta \times \Delta}\nn )
\rightarrow
Sect (\Delta , \rho (\Uu )^{\ast} \mm )^{\Delta}.
$$
Pour voir cela, fixons $(a,b)\in \Delta \times \Delta$. Alors la fl\`eche 
en question au point $(a,b)$ est la fl\`eche naturelle entre les deux
compos\'es dans le diagramme suivant:
$$
\begin{array}{ccccc}
\nn (a,b) = & \mm (\rho (\Uu )(a) \times _X \rho (\Uu )(b) )
& \stackrel{r_{1,\ast}}{\rightarrow} & \mm (\rho (\Uu )(a)) & 
= \rho (\Uu )^{\ast} \mm (a) \\
& {\scriptstyle p_2^{\ast}}\uparrow  && \uparrow {\scriptstyle \varphi ^{\ast}}
& \\
\rho (\Uu )^{\ast}\mm (b) = & \mm (\rho (\Uu )(b)) &
\stackrel{r_{\ast}}{\rightarrow} & \mm (X) & = \rho (\Uu )^{\ast}(\iota ).
\end{array}
$$
La condition (3) du th\'eor\`eme dit exactement que cette fl\`eche est une
\'equivalence.

On montre plus bas (comme dans la version initiale) 
que $p _{1,\ast}p_2^{\ast} (\sigma )$ est la
limite homotopique (le long de $\Delta$) de 
$r_{1,\ast} p_2^{\ast} (\sigma )$.  Donc l'\'equivalence donn\'ee ci-dessus
par la condition (3) implique que 
$$
p _{1,\ast}p_2^{\ast} (\sigma )\cong 
holim \varphi ^{\ast} r_{\ast} (\sigma ).
$$
La fl\^eche en question dans l'assertion est donc \'equivalente \`a
la fl\^eche $\ell$.

D'autre part, les $holim$ dans 
$Sect (\Delta , \rho (\Uu )^{\ast} \mm )$ se calculent objet-par-objet
au-dessus de $\Delta$. On peut voir cela en utilisant la structure (I)
du th\'eor\`eme \ref{sectionscmf}. L'hypoth\`ese (4) implique que
$\varphi ^{\ast}$ compos\'ee avec n'importe quelle restriction sur un objet
de $\Delta$, commute aux limites homotopiques. Il ensuit que $\varphi ^{\ast}$ 
commute
aux limites homotopiques. 
On conclut que $\ell$ est une \'equivalence, ce qui donne l'assertion.

Pour le reste de la pr\'esente d\'emonstration on reprend 
la premi\`ere version du papier. Plus loin, dans
les sections 20 et 21, nous aurons \`a v\'erifier l'hypoth\`ese
(4) du th\'eor\`eme lors de son utilisation. 


Soit ${\bf f}$ le foncteur ``remplacement fibrant'' pour la structure
de type (I-II) sur les sections de $\nn$. Le foncteur
$p _{1,\ast} \circ {\bf f}$ est invariant par \'equivalence donc on
obtient des \'equivalences
$$
p _{1,\ast} {\bf f}(p_1^{\ast} \sigma )\rightarrow
p _{1,\ast} {\bf f}(\gamma ^{\ast}\sigma ) \leftarrow
p _{1,\ast} {\bf f}(p_2^{\ast} \sigma ).
$$
Sous l'hypoth\`ese que $\sigma$ est fibrant (pour la structure de
type (II)), on
a que $p_2^{\ast}\sigma$ est fibrant pour la structure de type (I-II),
et donc
$$
p _{1,\ast}p_2^{\ast} (\sigma )
\rightarrow
p _{1,\ast} {\bf f}p_2^{\ast} (\sigma )
$$
est une \'equivalence. 
Notons que cela implique que $p _{1,\ast}p_2^{\ast} (\sigma )$ est une
limite homotopique de $r_{1,\ast}p_2^{\ast} (\sigma )$.
D'autre part on a le morphisme d'adjonction
$$
\varphi ^{\ast} \varphi _{\ast} (\sigma )
\rightarrow \sigma
$$
et encore
$$
\sigma \rightarrow p _{1,\ast} {\bf f} p_1^{\ast} \sigma .
$$
Leur compos\'e est le  morphisme
de gauche dans le carr\'e $(\ast )$
$$
\begin{array}{ccc}
\varphi ^{\ast} \varphi _{\ast} (\sigma )&\rightarrow &
p _{1,\ast} p_2^{\ast} (\sigma )\\
\downarrow && \downarrow \\
p _{1,\ast} p_1^{\ast} (\sigma )&\rightarrow &
p _{1,\ast} \gamma ^{\ast}(\sigma ).
\end{array}
$$
Les autres morphismes se d\'eduisent des morphismes ant\'erieurs.
C'est un fait pas tout-\`a-fait \'evident que ce carr\'e commute.
Pour le prouver, posons $\eta :=  \varphi _{\ast} (\sigma )$
et consid\'erons
d'abord le carr\'e $(\ast \ast )$
$$
\begin{array}{ccc}
t^{\ast} \eta & \rightarrow & p _2^{\ast} \sigma \\
\downarrow && \downarrow \\
p_1^{\ast} \sigma & \rightarrow & \gamma ^{\ast} \sigma ,
\end{array}
$$
o\`u
$$
t^{\ast} \eta := p_1^{\ast} \varphi ^{\ast} \eta = p_2^{\ast} \varphi
^{\ast} \eta .
$$
Ce carr\'e $(\ast \ast )$  commute. Pour le voir, notons $P_a$ le
produit fibr\'e de $a+1$ exemplaires de $U$ au-dessus de $X$ correspondant \`a
$a\in \Delta$; et notons $[x\rightarrow y]^{\ast} $ et $[x\rightarrow
y]_{\ast}$ les foncteurs de restriction et leurs adjoints. Notons qu'on a
$$
P_{\gamma (a,b)} = P_a\times _XP_b
$$
et
$$
\eta = \lim _{\leftarrow , a} [P_a\rightarrow X]_{\ast}(\sigma (a)).
$$
Notons $
f_a: \eta \rightarrow [P_a\rightarrow X]_{\ast}(\sigma (a)).
$ les morphismes structurels correspondants.
Les morphismes de transition dans la limite proviennent par exemple
des morphismes structurels de $\sigma$ comme
$$
[P_a\times _XP_b\rightarrow P_a]^{\ast} (\sigma (a)) \rightarrow
\sigma (\gamma (a,b)),
$$
qui donne
$$
[P_a\rightarrow X]_{\ast} (\sigma (a))\rightarrow
[P_a\times P_b\rightarrow X]_{\ast}
[P_a\times _XP_b\rightarrow P_a]^{\ast} (\sigma (a))
$$
$$
\rightarrow
[P_a\times P_b\rightarrow X]_{\ast}\sigma (\gamma (a,b)).
$$
On obtient que le compos\'e
$$
\eta \stackrel{f_a}{\rightarrow}
[P_a\rightarrow X]_{\ast} (\sigma (a))\rightarrow
\rightarrow
[P_a\times P_b\rightarrow X]_{\ast}\sigma (\gamma (a,b))
$$
est \'egal au morphisme $f_{\gamma (a,b)}$.
Le compos\'e en question est l'adjoint de
la composante en $(a,b)$ du compos\'e
inf\'erieur dans le carr\'e $(\ast \ast )$. Il s'ensuit que ce compos\'e
est \'egal au morphisme diagonal naturel
$$
t^{\ast} \eta\rightarrow \gamma ^{\ast} \sigma
$$
dont la composante en $(a,b)$ est l'adjoint de $f_{\gamma (a,b)}$.
Par sym\'etrie, le carr\'e $(\ast \ast )$ commute.

En appliquant au carr\'e $(\ast \ast )$ l'op\'eration
$p_{1,\ast }$ on obtient un carr\'e qu'on peut composer avec  le carr\'e
de naturalit\'e des morphismes d'adjonction
$$
\begin{array}{ccc}
\varphi ^{\ast} \eta & \rightarrow & p_{1,\ast } p_1^{\ast} \varphi ^{\ast}
\eta \\
\downarrow && \downarrow \\
\sigma &\rightarrow & p_{1,\ast } p_1^{\ast} \sigma
\end{array}
$$
pour montrer que le carr\'e $(\ast )$ commute.

En appliquant l'op\'eration ``remplacement fibrant'' au milieu sur la ligne en
bas on obtient  (par naturalit\'e de la transformation naturelle $u\rightarrow
{\bf f}u$) le carr\'e $(\ast )'$
$$
\begin{array}{ccc}
\varphi ^{\ast} \eta &\rightarrow &
p _{1,\ast} p_2^{\ast} (\sigma )\\
\downarrow && \downarrow \\
p _{1,\ast} {\bf f} p_1^{\ast} (\sigma )&\rightarrow &
p _{1,\ast} {\bf f} \gamma ^{\ast}(\sigma )
\end{array}
$$
qui commute. Or, ici, on sait d\'ej\`a que les morphismes en haut, \`a droite
et en bas sont des \'equivalences. Par ``trois pour le prix de deux'', le
morphisme \`a gauche est une \'equivalence. Rappelons que c'est le
compos\'e
$$
\varphi ^{\ast} \varphi _{\ast} (\sigma ) =\varphi ^{\ast} \eta
\rightarrow \sigma \rightarrow p _{1,\ast} {\bf f} p_1^{\ast} (\sigma ).
$$
En somme, pour prouver que notre morphisme d'adjonction
$$
\varphi ^{\ast} \varphi _{\ast} (\sigma ) \rightarrow \sigma
$$
est une \'equivalence, il suffit de
prouver que  le morphisme
$$
\sigma \rightarrow p _{1,\ast} {\bf f} p_1^{\ast} (\sigma )
$$
en est une.

Ce dernier fait rel\`eve de l'alg\`ebre homologique ou homotopique plus ou moins
standard. On indique bri\`evement
l'argument, en faisant r\'ef\'erence \`a Illusie
\cite{Illusie}
\footnote{
Le lecteur prendra garde que la signification de
notre exposant $+$ (e.g. $\Delta
^+$) diff\`ere de celle de \cite{Illusie}---c'est
plut\^ot exactement l'inverse.}
(qui lui-m\^eme fait r\'ef\'erence \`a des ``papiers
secrets'' de Deligne---mais le lecteur pourra prouver le r\'esultat en question
tout seul). Posons $$
P_m:= \rho (\Uu )(m) = \Uu \times _X\ldots \times _X \Uu .
$$
On utilisera les notations $[x\rightarrow y]^{\ast}$ et
$[x\rightarrow y]_{\ast}$ pour le foncteur de restriction et son adjoint
entre $\mm (y)$ et $\mm (x)$,
correspondant \`a un morphisme $x\rightarrow y$. D'autre
part on travaille dans $\Xx /X$ donc on peut supprimer l'indice $X$ pour le
produit fibr\'e au-dessus de $X$.
Avec ces notations, pour $(a,b)\in \Delta \times \Delta$ on a
$$
(p_1^{\ast} \sigma )(a,b)=
[P_a \times P_b \rightarrow P_a]^{\ast} \sigma (a).
$$
On obtient l'objet cosimplicial de $\mm (P_a)$, avec indice
(co)simplicial $b$
$$
C_{a,b}:= [P_a \times P_b \rightarrow P_a]_{\ast}
[P_a \times P_b \rightarrow P_a]^{\ast}\sigma (a).
$$
On a
$$
(p _{1,\ast} {\bf f}p_1^{\ast} \sigma )(a)\cong
holim _{b\in \Delta} C_{a,b}.
$$
Notons $\epsilon : P_a\times P_0 \rightarrow P_a$
la projection. On peut d\'ecomposer
la projection $P_a\times P_b\rightarrow P_a$ en produit de $b+1$ exemplaires de
$\epsilon$. On note $\epsilon ^{\ast}$ et $\epsilon _{\ast}$ les
fonctorialit\'es pour $\mm$ par rapport \`a $\epsilon$. En utilisant la formule
donn\'ee par l'hypoth\`ese (2) du pr\'esent th\'eor\`eme, on a
$$
C_{a,b} = (\epsilon _{\ast} \epsilon ^{\ast} )^{b+1}\sigma (a).
$$
En particulier, cet objet cosimplicial de $\mm (P_a)$ est la {\em r\'esolution
cosimpliciale standard de $\sigma (a)$ fournie par le couple de foncteurs
adjoints $\epsilon _{\ast}, \epsilon ^{\ast}$}---cf Illusie \cite{Illusie}
\S 1.5. Le th\'eor\`eme 1.5.3 de {\em loc cit.} implique que
l'objet cosimplicial augment\'e
$$
\epsilon ^{\ast} \sigma (a) \rightarrow \epsilon ^{\ast} C_{a,-}
$$
de $\mm (P_a\times P_0)$ est ``homotopiquement trivial''.
D'apr\`es \cite{Illusie} (1.1.5 et le descriptif pour 1.5.3) ceci veut dire
qu'il y a un morphisme
$$
\epsilon ^{\ast} C_{a,-}\rightarrow \epsilon ^{\ast} \sigma (a)
$$
tel que le compos\'e
$$
\epsilon ^{\ast} C_{a,-}\rightarrow
\epsilon ^{\ast} C_{a,-}
$$
soit homotope (au sens de \cite{Illusie} 1.1.5) \`a l'identit\'e
(l'autre compos\'e \'etant lui-m\^eme l'i\-den\-ti\-t\'e).

Il y a une section
$$
s:P_a\rightarrow P_a\times P_0
$$
(diagonale en derni\`ere variable), et en appliquant $s^{\ast}$ \`a l'objet
cosimplicial augment\'e ci-dessus, on obtient (avec la formule
$s^{\ast} \epsilon ^{\ast} = 1$) que
l'objet cosimplicial augment\'e
$$
\sigma (a) \rightarrow C_{a,-}
$$
de $\mm (P_a)$ est homotopiquement trivial.
Il en r\'esulte (car $holim$ transforme les homotopies de \cite{Illusie} 1.1.5
en  homotopies \`a la Quillen \cite{Quillen} pour la structure de cmf de $\mm
(P_a)$) que le morphisme
$$
holim _{b\in \Delta}\sigma (a)\rightarrow  holim _{b\in \Delta}
C_{a,b}
$$
admet un inverse \`a homotopie (de Quillen) pr\`es. Ceci implique
que c'est une \'equivalence faible. On note que
$$
\sigma (a)\rightarrow holim _{b\in \Delta}\sigma (a)
$$
est une \'equivalence. Pour ceci on renvoie au r\'esultat de Hirschhorn
(\cite{Hirschhorn}
Theorem 19.5.1) qui implique que la $holim$ d'un diagramme constant (dans
n'importe quelle cmf qui satisfait (0) par exemple), prise sur une cat\'egorie
dont le nerf est contractile (tel est le cas pour $\Delta$), est \'equivalente
\`a la valeur prise sur le diagramme.

On conclut qu'on a une \'equivalence
$$
\sigma (a)\stackrel{\cong}{\rightarrow} holim _{b\in \Delta}
C_{a,b}
$$
ce qui termine la d\'emonstration.
\eop

{\em Remarque:}
On pourrait penser que si $M$ est une cat\'egorie de mod\`eles ferm\'ee fixe
admettant des limites et colimites arbitraires, alors le pr\'efaisceau constant
$\underline{M}$ \`a valeurs dans $M$ v\'erifie les hypoth\`eses du
th\'eor\`eme pour une topologie quelconque $\Gg$. Pourtant $L(\underline{M})$,
qui est
le $1$-pr\'echamp de Segal constant \`a valeurs $L(M)$, n'est pas en g\'en\'eral
un champ. La raison en est que $L(\mm )$ n'est pas compatible aux sommes
disjointes.

{\bf Exercice:} Re\'ecrire l'\'enonc\'e et la d\'emonstration du th\'eor\`eme
\ref{dansleschamps} avec des structures de type Reedy au lieu de type (II);
sans l'hypoth\`ese que tous les objets des $\mm (y)$ soient cofibrants; et avec
les morphismes
des hypoth\`eses (2) et (3) suppos\'es \^etre
des \'equivalences faibles seulement
quand il s'agit d'objets cofibrants et fibrants.

On pr\'ecise que les auteurs n'ont pas fait cet exercice.

\subnumero{Comparaison avec SGA 4}
On compare ici notre r\'esultat du pr\'esent chapitre avec celui de
l'expos\'e Vbis de SGA 4 \cite{SGA4b} de Saint-Donat (d'apr\`es des
``notes
succinctes'' de Deligne).

La premi\`ere remarque \'a faire est que nous consid\'erons la situation
``instable'' via nos
cmf quelconques, tandis que \cite{SGA4b} est consacr\'e \`a la descente des
complexes ce qui
correspond \`a l'homotopie stable. N\'eanmoins, le cadre est tr\`es proche,
et l'application qui
nous a motiv\'es est justement l'exemple des complexes. Pour faire cette
transcription on note que
dans \cite{SGA4b} il s'agit d'une famille de topos annel\'es index\'ee
par la cat\'egorie de
base; on peut retrouver le cadre ``famille de cmf'' en prenant la famille des
cmf de complexes de
modules sur les anneaux structurels dans les topos fibres (on n'entre pas
dans les d\'etails d'une
d\'emonstration potentielle de la pr\'esente assertion). L'objet
$D^+(\underline{\Gamma}(E), A)$
de \cite{SGA4b} Vbis 2.2.7 est l'analogue de notre $Sect (\Delta , \int
_{\Delta}\mm )$; et
la sous-cat\'egorie pleine dont il s'agit dans
\cite{SGA4b} Vbis 2.2.7 (des sections dont les cohomologies sont
cocart\'esiennes) correspond \`a notre $Sect ^{\rm eq}(\Delta , \int
_{\Delta}\mm )$.

Dans le cadre de la descente des complexes (i.e. pour l'application au \S
21) on peut tr\`es bien
substituer l'argument de \cite{SGA4b} \`a notre argument du pr\'esent
chapitre. Cependant, il est
\`a noter que notre th\'eor\`eme \ref{dansleschamps} comporte aussi un
volet ``strictification des
donn\'ees de descente'' concr\'etis\'e dans l'utilisation du corollaire
\ref{pqcalclim}.
Ce volet n'a pas de contre-partie dans \cite{SGA4b} o\`u on part
d'une ``donn\'ee de
descente'' qui est d\'ej\`a un complexe dans le topos des sections,
c'est-\`a-dire d'une famille stricte
de complexes. On ne disposait pas \`a l'\'epoque du langage des
$\infty$-cat\'egories qui nous permet d'introduire
la notion de donn\'ee de descente faible.

Quant \'a la m\'ethode
pour descendre une donn\'ee stricte, l'argument de
\cite{SGA4b} pour le cas
des complexes est un d\'evissage par troncation du complexe, permettant de
traiter les objets de
cohomologie un \`a la fois. En partant d'une section $\sigma$
cocart\'esienne, on applique le
foncteur ``image directe'' que nous appelons $\varphi _{\ast}(\sigma )$
(dans la notation de
\cite{SGA4b} la donn\'ee de descente est $F^{\cdot}$ et son image directe est
${\bf R}^+ (\overline{\theta}_{\ast})(F^{\cdot})$). Il s'agit de montrer
que $\varphi
^{\ast}\varphi_{\ast} \sigma  $ est \'equivalent \`a $\sigma$; ce qui se
fait sur une page pour le
cas des complexes dans \cite{SGA4b} Vbis Prop. 2.2.7.

Cette m\'ethode serait aussi envisageable pour le cas ``instable''
pour peu qu'on dispose
d'une th\'eorie suffisante de la ``tour de Postnikov''.

La majeure partie de \cite{SGA4b} est consacr\'ee \`a la recherche de conditions
garantissant que
l'image inverse ${\bf L}^+(\overline{\theta}^{\ast})$ (que nous appelons
$L(\varphi ^{\ast})$) est
pleinement fid\`ele. Notre lemme \ref{dpq1} est l'analogue de la remarque
(\cite{SGA4b} Vbis Prop. 2.2.7).
Ensuite (pour notre corollaire \ref{dpq2},
dont la d\'emonstration se trouve
apr\`es \ref{dpq1})
nous avons utilis\'e la condition \ref{dpq2} (ii), qui
correspond \`a l'effectivit\'e
des donn\'ees de descente.
Cela rend plus facile l'argument de notre \ref{dpq2}; et on peut le faire
parce que de toutes fa\c{c}ons
nous avons besoin d'un argument sp\'ecifique pour obtenir l'effectivit\'e
de
\ref{dpq2} (ii). L'\'enonc\'e de \ref{dpq2} est \`a comparer \`a l'\'enonc\'e de
\cite{SGA4b} Vbis Prop. 2.2.7.

Regardons enfin les diverses hypoth\`eses dans notre th\'eor\`eme
\ref{dansleschamps}.
L'hypoth\`ese sur l'existence de suffisamment de sommes disjointes, est \`a
rapprocher de
\cite{SGA4b} Vbis 3.0.0. On peut comparer la condition (3) avec
\cite{SGA4b} Vbis 3.2.
Ici on peut remarquer que Saint-Donat utilise une famille quelconque de
changements de base pour
``tester'', tandis que notre strat\'egie est de tester avec le
changement de base qu'on
utilise pour la descente.

La technique bisimpliciale de la d\'emonstration de \ref{dansleschamps}
s'apparente \`a la
technique bisimpliciale pour la pleine fid\'elit\'e de SGA 4 (\cite{SGA4b}
Vbis \S 2.3, voir aussi
3.3.1(a)).

En somme, la strat\'egie globale, qui consiste \`a
descendre $\sigma$ en prenant
$\varphi _{\ast}(\sigma )$,
est commune, mais le probl\`eme de prouver que cela r\'epond bien \`a la
question est pour nous plus
compliqu\'e car nous ne disposons pas d'argument de r\'ecurrence
``\`a la Postnikov''. Cela n\'ecessite
alors un argument sp\'ecifique (preuve de \ref{dansleschamps}); et avec
cela on peut contourner le
probl\`eme de la pleine fid\'elit\'e de $L(\varphi ^{\ast})$ (\ref{dpq2}).
Naturellement, les
arguments de \cite{SGA4b} pour la pleine fid\'elit\'e, se retrouvent dans
notre argument
sp\'ecifique pour l'effectivit\'e de
\ref{dansleschamps}.

\numero{Exemple: la descente pour les $n$-champs}

\label{exemplepage}

Le premier exemple d'application du th\'eor\`eme \ref{dansleschamps}
est la descente pour les $n$-champs, analogue du fait que le pr\'echamp des
faisceaux est un champ. Ce th\'eor\`eme de descente est classique pour $n=1$.
Pour $n=2$, il a \'et\'e \'enonc\'e par Breen dans \cite{Breen}, sans
d\'emonstration. Breen l'utilise pour donner l'une des deux directions de sa
description des $2$-champs en termes de cocycles.

\begin{theoreme}
\label{descente}
Le $n+1$-pr\'echamp de Segal $nSe\underline{CHAMP}(\Xx )$ est un $n+1$-champ
de Segal (i.e. $\infty$-champ) au-dessus de $\Xx$.
\end{theoreme}

Ce th\'eor\`eme signifie que les $n$-champs de Segal se recollent, i.e.
notamment que les donn\'ees de descente pour les $n$-champs de Segal, sont
effectives. Cependant, comme nous n'avons pas une tr\`es bonne description
concr\`ete de ce que c'est exactement
qu'une donn\'ee de descente (voir la fin de
\S 5 pour le mieux qu'on  puisse dire actuellement), nous ne pouvons pas donner
une description enti\`erement en termes de cocycles, \`a la Breen
\cite{Breen}, de ce que veut dire cette effectivit\'e des donn\'ees de
descente.

\subnumero{Preuve utilisant le th\'eor\`eme \ref{dansleschamps}}
Supposons que $\Xx$ admet des produits fibr\'es et a suffisamment de sommes
disjointes compatibles aux produits fibr\'es. Dans  ce cas,
on peut d\'eduire le pr\'esent th\'eor\`eme \ref{descente} directement du
th\'eor\`eme \ref{dansleschamps}.
En effet, on peut v\'erifier que la famille des cmf
$$
X\mapsto nSePCh(\Xx /X)
$$
est un pr\'efaisceau de Quillen \`a gauche, avec compatibilit\'e aux produits
directs. L'adjoint de la restriction pour $s: Y\rightarrow X$ et $A\in
nSePCh(\Xx /Y)$ est d\'efini par
$$
s_{\ast}(A)(U\rightarrow X):= A(Y\times _XU\rightarrow Y).
$$
La restriction pr\'eserve \'evidemment les cofibrations et
les cofibrations triviales, et
la compatibilit\'e aux produits directs (i.e. la condition (3) de
\ref{dansleschamps}) r\'esulte d'un calcul facile.

Pour la condition (2),
il suffit d'observer que si
$i: U\hookrightarrow X$ est une composante
d'une somme disjointe, et si $A$ est n'importe quel type de pr\'efaisceau sur
$U$ alors on a
$i^{\ast}i_{\ast}A=A$.

La condition (1) est une
cons\'equence imm\'ediate du caract\`ere local de la
notion de $\Gg$-\'equivalence faible entre $n$-pr\'echamps de Segal.

Les
cat\'egories de mod\`eles fibres admettent des limites et colimites
arbitraires, et des
factorisations fonctorielles. Pour les structures de type (II), voir le
th\'eor\`eme
\ref{sectionscmf}; on renvoie comme d'habitude \`a \cite{Jardine} pour la
technique n\'ecessaire pour
donner une d\'emonstration rigoureuse de la premi\`ere partie de
l'hypoth\`ese de \ref{sectionscmf}
(II). On obtient ainsi la condition (0) du th\'eor\`eme \ref{dansleschamps}.


Il faut v\'erifier la condition (4)
(ajout\'ee dans version 2 du papier) du th\'eor\`eme \ref{dansleschamps}. 
Soit donc $f:X\rightarrow
Y$ un morphisme dans $\Xx$. Il s'agit de voir que
$$
f^{\ast} :L(nSePCh , W)(Y) \rightarrow L(nSePCh , W) (X)
$$
pr\'eserve les limites (homotopiques). Ici $f^{\ast}$ est juste la restriction
\`a $\Xx /X$ des $n$-champs de Segal sur $\Xx /Y$. Les limites des
$n$-champs de Segal se calculent objet-par-objet, donc la restriction pr\'eserve
ces limites. Le fait que les limites se calculent objet-par-objet,
est une g\'en\'eralisation du lemme \ref{fiprodchamp} qui concerne le cas du produit
fibr\'e. Pour le cas g\'en\'eral on peut utiliser la commutation des limites
(\cite{Vogt2}, \cite{limits} 3.4.11) ainsi que le crit\`ere \ref{critereaveclim} (c) pour voir qu'une 
limite (en tant que pr\'echamp)
de champs est encore un champ. Les limites de pr\'echamps se calculent
objet-par-objet d'apr\`es \cite{limits} 3.4.4. 



Il ne restera plus qu'\`a v\'erifier l'hypoth\`ese (5) de
\ref{dansleschamps} (nouveau dans {\tt v3})
pour prouver que $L(nSePCh, W)$
est un $1$-champ de Segal. On fera cela en m\^eme temps que la preuve par
le crit\`ere \ref{critere} que $nSe\underline{CHAMP}$ est un champ.

Rappelons que par \ref{intereqloc}, $L(nSePCh, W)$
est \'equivalent (pour la topologie
grossi\`ere) \`a l'int\'erieur $1$-groupique
$$
L(nSePCh, W) \cong nSe\underline{CHAMP}^{int,1}.
$$
La partie (a) du crit\`ere \ref{critere} pour $nSe\underline{CHAMP}$
est automatique  car, par d\'efinition, les
objets de $nSe\underline{CHAMP}(X)$ sont les $n$-champs de Segal $\Gg$-fibrants,
et si $A$ et $B$ sont $\Gg$-fibrants sur $\Xx /X$
alors 
$$
nSe\underline{CHAMP}_{1/}(A,B) := 
\underline{Hom}(A,B)
$$ 
est fibrant---en particulier c'est un champ.
Il s'ensuit que l'int\'erieur $0$-groupique de 
$nSe\underline{CHAMP}_{1/}(A,B)$ est un champ. L'objet des morphismes dans
l'int\'erieur $1$-groupique est aussi l'int\'erieur $0$-groupique de
\'objet des morphismes, donc on obtient que 
l'int\'erieur $1$-groupique
$nSe\underline{CHAMP}^{int,1}$
(\'equivalent \`a $L(nSePCh, W)$)
satisfait la partie (a) du crit\`ere \ref{critere}. Or cette condition
est exactement la condition (5) du th\'eor\`eme \ref{dansleschamps}.

On a donc fini de v\'erifier les hypoth\`eses de \ref{dansleschamps},
et d'apr\`es ce th\'eor\`eme, le localis\'e $L(nSePCh, W)$
est un $1$-champ de Segal. 
Par ce r\'esultat, on obtient la partie (b) du crit\`ere de \ref{critere} 
(voir aussi \ref{devissage}) pour
$nSe\underline{CHAMP}$. 

La proposition \ref{critere}  
s'applique maintenant pour conclure que
$nSe\underline{CHAMP}$ est
un $n+1$-champ de Segal.
\eop


\subnumero{Preuve directe}
On donne maintenant une preuve directe du th\'eor\`eme \ref{descente},
qui n'utilise pas
le th\'eor\`eme \ref{dansleschamps}. Cette d\'emonstration semble
donner plus d'indications que la premi\`ere sur la mani\`ere dont les
$n$-champs de Segal se recollent. De plus, elle s'applique
sans autre hypoth\`ese sur
$\Xx$ et s'adapte {\em mutatis mutandis} au cas des $n$-champs non de
Segal.

Pour la suite on fixera
un site $\Xx$ qu'on ne mentionnera
plus dans les notations. En cas de besoin on notera
$\Gg$ la
topologie sur $\Xx$ (par opposition \`a la topologie grossi\`ere).

On utilisera la construction $\Upsilon$ (ou ses variantes $\Upsilon ^k$) de
\cite{limits}. Elle se g\'en\'eralise im\-m\'e\-dia\-te\-ment
aux $n$-pr\'ecats de
Segal et, par fonctorialit\'e, aux pr\'efaisceaux de
$n$-pr\'ecats de Segal, i.e. aux $n$-pr\'echamps de Segal, en conservant ses
propri\'et\'es universelles, dont on aura quelquefois besoin.

Notons
$$
nSe\underline{CHAMP}\rightarrow nSe\underline{CHAMP}'
\rightarrow nSe\underline{CHAMP}''
$$
deux cofibrations: la premi\`ere est une cofibration triviale pour la topologie
grossi\`ere (en particulier une \'equivalence objet-par-objet) vers un
mod\`ele
fibrant pour la topologie grossi\`ere;  la deuxi\`eme est une cofibration
triviale
pour la topologie $\Gg$ vers un mod\`ele fibrant pour la topologie $\Gg$. Leur
compos\'e est \'egalement une cofibration triviale pour la topologie
$\Gg$ vers
un mod\`ele fibrant pour la topologie $\Gg$.

Par construction les $nSe\underline{CHAMP}(X)$ sont d\'ej\`a des
$n+1$-cat\'egories de Segal. La condition que $nSe\underline{CHAMP}$ soit
un $n+1$-champ de Segal est donc \'equivalente \`a la condition que le morphisme
$$
nSe\underline{CHAMP}'\rightarrow nSe\underline{CHAMP}''
$$
soit une \'equivalence pour la topologie grossi\`ere, i.e. une \'equivalence
objet-par-objet.

On va utiliser le crit\`ere \ref{critere}.
Pour $A$ et $B$ dans $nSe\underline{CHAMP}(X)$
on  a par construction
$$
nSe\underline{CHAMP}(X)_{1/} (A,B) = \underline{Hom}(A,B).
$$
En particulier, comme $A$ et $B$ sont cofibrants et fibrants (par
la d\'efinition m\^eme de $nSe\underline{CHAMP}$) ceci est un $n$-champ de
Segal.
La premi\`ere hypoth\`ese de \ref{critere} est donc v\'erifi\'ee.

Il s'ensuit que, pour tout $X$, le morphisme
$$
nSe\underline{CHAMP}(X)
\rightarrow
nSe\underline{CHAMP}''(X)
$$
est pleinement fid\`ele. En effet, le morphisme
$$
nSe\underline{CHAMP}
\rightarrow
nSe\underline{CHAMP}''
$$
est une $\Gg$-\'equivalence faible de $n+1$-pr\'echamps de
Segal dont les valeurs sont des $n+1$-cat\'egories de Segal; par d\'efinition
pour tout $X$ et tous $A,B\in nSe\underline{CHAMP}'(X)_0$, le morphisme
$$
nSe\underline{CHAMP}_{1/}(A,B)
\rightarrow
nSe\underline{CHAMP}''_{1/}(A,B)
$$
est une \'equivalence faible de $n$-pr\'echamps de Segal dont la
source et le but
sont d\'ej\`a des $n$-champs de Segal; par \ref{morphentrechamps}
c'est une \'equivalence faible pour la topologie grossi\`ere.
Ce r\'esultat admet la g\'en\'eralisation suivante:

\begin{lemme}
\label{hompleinfidele}
Si $A$ est un $n+1$-pr\'echamp de Segal alors le morphisme de
$n+1$-pr\'echamps de
Segal sur $\Xx$
$$
\underline{Hom}(A, nSe\underline{CHAMP}')\rightarrow
\underline{Hom}(A,nSe\underline{CHAMP}'')
$$
est pleinement fid\`ele pour la topologie grossi\`ere.
\end{lemme}
{\em Preuve:}
Soit $Z$ le sous-pr\'echamp plein de $nSe\underline{CHAMP}''$
constitu\'e avec les objets qui sont dans l'image de $nSe\underline{CHAMP}'$.
Alors le morphisme
$nSe\underline{CHAMP}'\rightarrow Z$ est une \'equivalence faible pour la
topologie grossi\`ere. Donc il en est de m\^eme du morphisme
$$
\underline{Hom}(A, nSe\underline{CHAMP}')\rightarrow
\underline{Hom}(A, Z).
$$
Mais le morphisme
$$
\underline{Hom}(A, Z)\rightarrow
\underline{Hom}(A, nSe\underline{CHAMP}'')
$$
est pleinement fid\`ele pour la topologie grossi\`ere: si $a,b$ sont
deux objets de $\underline{Hom}(A, Z)(X)$ alors un morphisme
$$
U\rightarrow \underline{Hom}(A, Z)_{1/}(a,b)
$$
correspond \`a un morphisme
$$
A\times \Upsilon (U)\rightarrow Z
$$
qui donne $a$ (resp. $b$) sur $A\times 0$ (resp. $A\times 1$);
comme tous les objets de   $A\times \Upsilon (U)$ sont soit dans
$A\times 0$ soit dans $A\times 1$, ceci
correspond encore \`a un morphisme
$$
A\times \Upsilon (U)\rightarrow nSe\underline{CHAMP}''
$$
qui donne $a$ (resp. $b$) sur $A\times 0$ (resp. $A\times 1$),
autrement dit \`a un morphisme
$$
U\rightarrow \underline{Hom}(A, nSe\underline{CHAMP}'')_{1/}(a,b).
$$
En fait, on a l'\'egalit\'e
$$
\underline{Hom}(A, Z)_{1/}(a,b)  =
\underline{Hom}(A, nSe\underline{CHAMP}'')_{1/}(a,b).
$$
Par composition, on obtient que
$$
\underline{Hom}(A, nSe\underline{CHAMP}')\rightarrow
\underline{Hom}(A, nSe\underline{CHAMP}'')
$$
est pleinement fid\`ele pour la topologie grossi\`ere.
\eop

On poursuit la d\'emonstration du th\'eor\`eme \ref{descente}.
Au vu du crit\`ere de la proposition \ref{critere},
il s'agit maintenant de prouver que pour tout crible
couvrant $\Bb \subset \Xx /X$, tout
morphisme
$$
f:\ast _{\Bb} \rightarrow nSe\underline{CHAMP}'|_{\Xx /X}
$$
s'\'etend en un morphisme d\'efini sur $\ast _X$. Comme la formation de
$nSe\underline{CHAMP}$ est compatible aux changements de base, et comme
le remplacement fibrant pour la topologie grossi\`ere peut aussi \^etre
choisi compatible aux restrictions (lemme \ref{shriek}), on peut se placer
dor\'enavant sur le site $\Xx /X$ et supprimer dans nos notations la
r\'ef\'erence \`a la ``restriction \`a $\Xx
/X$''.

Soit $$
E:= \underline{Hom}(\ast _{\Bb},nSe\underline{CHAMP}')_{1/}(\ast , f).
$$
Ici $\underline{Hom}$ est le $Hom$ interne des $n+1$-pr\'echamps de Segal sur
$\Xx /X$. En particulier,  $E$ est un $n$-pr\'echamp de Segal sur $\Xx /X$.

Le morphisme
$$
\underline{Hom}(\ast _{\Bb},nSe\underline{CHAMP}')_{1/}(\ast , f)
\rightarrow
\underline{Hom}(\ast _{\Bb},nSe\underline{CHAMP}'')_{1/}(\ast , f)
$$
est une \'equivalence faible de $n$-pr\'echamps
pour la topologie grossi\`ere sur
$\Xx /X$ d'apr\`es le lemme \ref{hompleinfidele}. Comme le but de ce
morphisme est
$\Gg$-fibrant, c'est un $n$-champ de Segal donc
sa source
---qui est $E$---est
\'egalement un $n$-champ de Segal. Ce dernier \'etant aussi fibrant pour la
topologie grossi\`ere, on obtient que $E$ est $\Gg$-fibrant par
\ref{grofibchampimplfib}.

L'id\'ee de la d\'emonstration est que $E$ r\'esout le probl\`eme de
descente qui est pos\'e. Il faut
encore un peu de technique pour obtenir un morphisme
entre $f$ et le morphisme constant $cst(E)$, et
achever la preuve du th\'eor\`eme.

On d\'efinit d'abord un $n+1$-pr\'echamp de Segal $\Upsilon (E)$ en
appliquant la
d\'efinition de \cite{limits} objet-par-objet, i.e.
$$
\Upsilon (E)(Y):= \Upsilon (E(Y))
$$
(l'extension des constructions de \cite{limits} objet-par-objet sera utilis\'ee
ci-dessous sans autre commentaire).

Par la propri\'et\'e
universelle qui d\'efinit $\Upsilon$ on obtient un morphisme $$
\Upsilon (E)\rightarrow \underline{Hom}(\ast _{\Bb},nSe\underline{CHAMP}')
$$
ou encore
$$
h: \ast _{\Bb} \times \Upsilon (E) \rightarrow nSe\underline{CHAMP}'.
$$
On a
$$
h|_{\ast _{\Bb} \times 0} = \ast
$$
et
$$
h|_{\ast _{\Bb} \times 1} = f.
$$

On utilise maintenant les propri\'et\'es d'extension de morphismes d\'efinis
sur des parties de $\Upsilon ^k$---on
renvoie \`a \cite{limits} pour les d\'etails---,
en notant qu'elles s'\'etendent au cadre des pr\'echamps.
Rappelons les notations $01$, $12$ et $02$ pour les copies de $\Upsilon (A)$,
$\Upsilon (B)$ et $\Upsilon (A\times B)$ respectivement dans $\Upsilon ^2(A,B)$.

Par l'analogue de \cite{limits} 4.3.5, il existe un morphisme
$$
h^2: \ast _{\Bb} \times \Upsilon ^2(E,\ast )
$$
avec
$$
h^2|_{01} = 1_E,
$$
et
$$
h^2|_{02} = h.
$$
On rappelle que dans les notations de \cite{limits}
$1_E$ est une application
$$
\Upsilon (E) \rightarrow nSe\underline{CHAMP}'
$$
qui induit $\ast$ sur le premier sommet $0$; et
qui induit le morphisme constant
$$
cst(E) :\ast \rightarrow nSe\underline{CHAMP}'
$$
\`a valeur $E$, sur le deuxi\`eme sommet $1$.

La restriction de $h$ \`a l'ar\^ete $12$ fournit donc un morphisme
$$
g: \ast _{\Bb} \times I \rightarrow nSe\underline{CHAMP}',
$$
o\`u nous avons not\'e $I$ la $1$-cat\'egorie avec deux objets et une fl\`eche,
qui est aussi \'egale \`a $\Upsilon (\ast )$.
Les objets de $I$ seront appel\'es $0$ et $1$ bien que cette notation soit
incompatible avec la notation pour les sommets de $\Upsilon ^2(E,\ast )$
(l'objet $0$ de $I$ correspondant au sommet $1$,
et l'objet $1$ de $I$ correspondant au
sommet $2$).

L'image du premier objet $0$ de $I$ (donc l'image du sommet $1$ de
$\Upsilon ^2$)
est l'application constante $cst(E)$ compos\'ee avec la projection $\ast
_{\Bb}\rightarrow \ast$. L'image
du deuxi\`eme objet $1$ de $I$ est l'application $f$.

Noua affirmons que $g$ s'\'etend en un morphisme
$$
\overline{g}: \ast _{\Bb} \times \overline{I} \rightarrow
nSe\underline{CHAMP}'
$$
o\`u $\overline{I}$ est la $1$-cat\'egorie avec deux objets $0,1$ et un
isomorphisme entre eux (on a $I\subset \overline{I}$). On
observe d'abord la propri\'et\'e suivante:

$(\ast
)$ pour tout
$Y\in \Bb$, le morphisme induit par $g$
$$
E(Y)\rightarrow f_Y(\ast _{\Bb}(Y))
$$
est une \'equivalence faible de $n+1$-pr\'ecats de Segal.

En effet, on a
$$
E(Y)= \Gamma \underline{Hom}(\ast _{\Bb }|_{\Xx /Y},
nSe\underline{CHAMP}'|_{\Xx /Y})_{1/}(\ast , f),
$$
mais aussi $\ast _{\Bb }|_{\Xx /Y}=\ast _{\Xx /Y}$ et donc
$$
E(Y)= \Gamma (\Yy /Y, nSe\underline{CHAMP}')_{1/}(\ast , f|_{\Xx /Y})
$$
$$
\cong nSeCHAMP (Y)_{1/}(\ast , f_Y(\ast _{\Bb}(Y))) \cong
f_Y(\ast _{\Bb}(Y)).
$$

Notons $C$ le r\'esultat obtenu \`a
partir de $\ast _{\Bb} \times I$, en ajoutant librement
au-dessus de chaque objet de $\Bb$, une cellule correspondant \`a la cofibration
$I\hookrightarrow \overline{I}$. On a un morphisme de $C$ vers
$\ast_{\Bb} \times \overline{I}$. Ce morphisme est une \'equivalence, car par
\cite{limits} 2.5.1, le fait de faire (pour inverser la m\^eme fl\`eche)
plusieurs fois le coproduit avec la cofibration $I\rightarrow
\overline{I}$ a le m\^eme effet que de le faire une seule fois.
Maintenant, l'\'enonc\'e $(\ast )$ ci-dessus implique que $g$ s'\'etend en un
morphisme $C\rightarrow nSe\underline{CHAMP}'$. Par les arguments habituels,
cela implique que $g$ s'\'etend en $\overline{g}$ comme annonc\'e.

Gr\^ace \`a ce r\'esultat, on peut poser
$$
V:= (\ast _{\Bb}\times \overline{I}) \cup ^{\ast _{\Bb} \times 0}\ast _X.
$$
L'application $V\rightarrow \ast _X$ est une \'equivalence faible pour la
topologie grossi\`ere. L'application
$$
cst(E): \ast _X \rightarrow nSe\underline{CHAMP}'
$$
se recolle avec l'application
$$
\overline{g}: \ast _{\Bb} \times \overline{I} \rightarrow
nSe\underline{CHAMP}'
$$
car on a, par construction,
$$
\overline{g}|_{\ast _{\Bb} \times 0} = cst(E).
$$
On obtient une application
$$
V\rightarrow nSe\underline{CHAMP}'
$$
dont la restriction sur $\ast _{\Bb} \times 1$ est \'egale \`a $f$.
Soit
$$
V\rightarrow V' \rightarrow \ast _X
$$
une factorisation avec comme premier morphisme une cofibration triviale et
comme
deux\-i\`eme morphisme une fibration triviale (le
tout pour la topologie grossi\`ere).
L'application ci-dessus s'\'etend en une application
$$
\alpha : V'\rightarrow nSe\underline{CHAMP}'
$$
dont la restriction \`a
$\ast _{\Bb} \times 1\subset V \subset V'$ est encore \'egale \`a $f$.
Consid\'erons le diagramme
$$
\begin{array}{ccc}
\ast _{\Bb} & \rightarrow & V' \\
\downarrow && \downarrow \\
\ast _X & = & \ast _X.
\end{array}
$$
La fl\`eche du haut est l'inclusion
$$
i: \ast _{\Bb}=\ast _{\Bb} \times 1\subset V
\subset V'.
$$
On peut \'ecrire $\alpha \circ i = f$.
La fl\`eche de droite est une fibration triviale pour la topologie grossi\`ere,
en particulier elle poss\`ede la propri\'et\'e de rel\`evement pour toute
cofibration.
Donc il existe un rel\`evement $j:\ast _X\rightarrow V'$ \'egal \`a
$i$ sur $\ast _{\Bb}$.   Le morphisme
$$
\alpha \circ j: \ast _X \rightarrow nSe\underline{CHAMP}'
$$
est l'extension de $f$ que nous recherchons, ce qui compl\`ete la
d\'emonstration
du th\'eor\`eme \ref{descente}.
\eop

\begin{corollaire}
Les pr\'efaisceaux de cat\'egories simpliciales $X\mapsto L(nSePCh(\Xx /X))$
sont des $1$-champs de Segal. En particulier le pr\'efaisceau de cat\'egories
simpliciales
$$
X\mapsto L(PrefSpl(\Xx /X),W^{\rm Illusie})
$$
est un $1$-champ de Segal.
\end{corollaire}
{\em Preuve:}
Par le th\'eor\`eme \ref{descente}, $nSe\underline{CHAMP}$ est un $n+1$-champ de
Segal. Donc (cf \S 10) son int\'erieur $1$-groupique est un $1$-champ de Segal.
On applique \ref{intereqloc}.
\eop

Ce corollaire montre que les donn\'ees de descente homotopiques pour les
pr\'efaisceaux simpliciaux ou $n$-pr\'echamps de Segal, sont effectives, a
\'equivalence pr\`es.

Soit $nSeChamp(\Xx )$ la $1$-cat\'egorie stricte des pr\'efaisceaux de
$n$-cat\'egories de Segal qui sont des $n$-champs de Segal sur $\Xx$.
Ici, la notion de $\Gg$-\'equivalence co\"{\i}ncide avec celle d'\'equivalence
objet-par-objet. On obtient que {\em pour
$L(nSeChamp(\Xx ))$ les donn\'ees de descente sont effectives}, ce qui veut dire
que les $n$-champs de Segal se recollent, \`a \'equivalence objet-par-objet
pr\`es. Cet \'enonc\'e est
l'analogue du r\'esultat classique de recollement des faisceaux.

\subnumero{Enonc\'es pour les $n$-champs (non de Segal)}

Pour l'essentiel, tous les \'enonc\'es du
pr\'esent travail, sauf ceux concernant
directement les localis\'ees de Dwyer-Kan, sont aussi valables pour les
$n$-cat\'egories, $n$-champs etc., non de Segal. Sans entrer dans les details,
nous indiquons ici les faits principaux (entre autres, pour justifier le titre
choisi!).

D'abord, la cat\'egorie $nPC$ des $n$-pr\'ecats remplace $nSePC$, celle
des $n$-pr\'ecats de Segal, et on dispose de l'op\'eration $Cat$ (analogue \`a
$SeCat$), et la notion de $n$-cat\'egorie a la m\^eme d\'efinition---tout cela
se trouve dans
\cite{Tamsamani} et \cite{nCAT}.  Un {\em $n$-pr\'echamp sur $\Xx$} est un
pr\'efaisceau sur $\Xx$ \`a valeurs dans $nPC$; on note
$nPCh(\Xx )$ la cat\'egorie des
$n$-pr\'echamps.

Les cat\'egories $nPC$ et $nPCh(\Xx )$ ont des
structures de cmf analogues \`a celles de \ref{SeCmf} et \ref{cmf}, et
$nPCh(\Xx )$ a \'egalement la structure de type HBKQ de \ref{cmfHirschho}.

On dit (de la m\^eme fa\c{c}on  qu'au \S 9) qu'un $n$-pr\'echamp $A$ est
un {\em
$n$-champ} si le morphisme $A\rightarrow A'$ de remplacement fibrant pour la
topologie $\Gg$, est une \'equivalence objet-par-objet. Les propri\'et\'es du
\S 9 et l'autre d\'efinition du \S 10 s'adaptent {\em mutatis mutandis}.

Les cmf $nPC$ et $nPCh(\Xx )$ sont internes (\S 11; pour $nPC$ c'est dans
\cite{nCAT}), ce qui permet de d\'efinir les $n+1$-cat\'egories
$$
nCAT\;\;\; \mbox{et} \;\;\; nCHAMP(\Xx )
$$
ainsi que le $n+1$-champ $n\underline{CHAMP}(\Xx )$.

On d\'efinit une famille universelle et un morphisme $\Phi$ comme au \S 12.
Le th\'eor\`eme \ref{correlation} prend la forme suivante.

\begin{theoreme}
\label{correlationNdS}
Soit $nCAT'$ un remplacement fibrant de $nCAT$.
Si la topologie de $\Xx$ est
grossi\`ere alors $nCHAMP(\Xx )$ est \'equivalent via $\Phi$ \`a la
$n+1$-cat\'egorie $\underline{Hom}(\Xx ^o, nCAT')$.

Si $\Xx$ est muni d'une topologie quelconque
$\Gg$ alors $nCHAMP(\Xx )$ est \'equivalente \`a la sous-cat\'egorie
pleine de
$\underline{Hom}(\Xx ^o, nCAT')$ form\'ee des morphismes
$F:\Xx ^o\rightarrow nCAT'$ qui satisfont \`a la condition de descente
(\cite{limits} 6.3) qui dit que la fl\`eche
$$
\lim _{\leftarrow} F|_{\Xx /X} \rightarrow
\lim _{\leftarrow} F|_{\Bb}
$$
est une \'equivalence pour tout crible $\Bb \subset \Xx /X$ de $\Gg$.
\end{theoreme}
\eop

On a un foncteur ``champ associ\'e'' entre $n+1$-cat\'egories
$$
{\bf ch} : nCHAMP (\Xx ^{\rm gro})\rightarrow nCHAMP (\Xx ^{\Gg})
$$
avec transformation naturelle $F\rightarrow {\bf ch}(F)$,
et (th\'eor\`eme \ref{cestadjoint}) cette transformation
naturelle fait de ce foncteur l'adjoint homotopique de l'inclusion
$$
nCHAMP (\Xx ^{\Gg})\subset nCHAMP (\Xx ^{\rm gro}).
$$

Enfin, le th\'eor\`eme principal du pr\'esent chapitre \ref{descente}
devient:

\begin{theoreme}
\label{descenteNdS}
Le $n+1$-pr\'echamp  $n\underline{CHAMP}(\Xx )$ est un $n+1$-champ
au-dessus de $\Xx$.
\end{theoreme}

La preuve directe s'adapte {\em mutatis mutandis}. Si on veut obtenir une preuve
par application de \ref{dansleschamps}, il faut observer le fait suivant:

\begin{proposition}
\label{segalounon}
La $n+1$-cat\'egorie $nCAT$, consid\'er\'ee
par induction comme une $n+1$-cat\'egorie de Segal, est \'equivalente
\`a la sous-$n+1$-cat\'egorie de Segal pleine de $nSeCAT$ form\'ee
des objets qui sont $n$-tronqu\'es.

La $n+1$-cat\'egorie $nCHAMP(\Xx )$,
consid\'er\'ee de la m\^eme fa\c{c}on comme une $n+1$-cat\'egorie de Segal,
est \'equivalente \`a la sous-$n+1$-cat\'egorie de Segal
pleine de
$nSeCHAMP$ form\'ee par les $n$-champs de Segal qui sont $n$-tronqu\'es
objet-par-objet au-dessus de $\Xx$.

Le $n+1$-pr\'echamp $n\underline{CHAMP}(\Xx )$,
consid\'er\'e par induction comme un $n+1$-pr\'echamp de Segal,
est \'equivalent au sous-pr\'echamp de Segal plein de
$nSe\underline{CHAMP}(\Xx )$
form\'e des objets dont les valeurs (sur tout $X\in \Xx$) sont
des $n$-champs de
Segal (sur $\Xx/X$) qui sont $n$-tronqu\'es objet-par-objet.
\end{proposition}
\eop

On peut alors appliquer \ref{preservetronque}.
\eop

\numero{Exemple: la descente pour les complexes}

\label{complexepage}

Quand on parle des champs, l'exemple des complexes a toujours \'et\'e pr\'esent
aussi bien en g\'eom\'etrie alg\'ebrique qu'en topologie alg\'ebrique
(voir \cite{SGA4b}, \cite{SGA6}, \cite{GabrielZisman}, \cite{Quillen}, \cite{ResiduesDuality},
\cite{Illusie}, \cite{ThomasonTrobaugh},
\cite{BondalKapranov}, \cite{Hinich2}, \cite{HoveyBook}  par exemple).
C'est d'ailleurs cet exemple des complexes qui nous a entra\^{\i}n\'e dans
ce travail, et qui nous a guid\'e notamment dans la formulation et
d\'emonstration des r\'esultats des \S 17-\S 18-\S 19.

Soit $\Xx$ un site muni d'un faisceau d'anneaux $\Oo$.
On suppose que $\Xx$ admet des produits fibr\'es et a suffisamment de sommes
disjointes compatibles aux produits fibr\'es. Pour $X\in \Xx$, soit
$Cpx_{\Oo}(X)$ la cat\'egorie des complexes de faisceaux de
$\Oo$-modules sur $\Xx /X$. Soit
$$
qis(X)\subset Cpx_{\Oo}(X)
$$
la sous-cat\'egorie des fl\`eches qui sont des quasi-isomorphismes.

Notons $Cpx_{\Oo}^{[0, \infty )}(X)$ la sous-cat\'egorie des complexes
concentr\'es en degr\'es
$\geq 0$. D'apr\`es
Quillen \cite{Quillen}, $Cpx_{\Oo}^{[0,
\infty )}(X)$ a une structure de cat\'egorie de mod\`eles ferm\'ee o\`u les
\'e\-qui\-va\-len\-ces faibles sont les quasi-isomorphismes et les cofibrations
sont les injections.

R\'ecemment Hinich \cite{Hinich2} et Hovey \cite{HoveyBook} ont muni
la cat\'egorie $Cpx_{\Oo}(X)$ de tous les complexes (non-n\'ecessairement
born\'es) d'une structure de cmf dont les
\'equivalences faibles sont les quasi-isomorphismes et les cofibrations
sont les injections.

En fait, Quillen aussi bien que Hinich (et Hovey dans une version
pr\'ec\'edente de \cite{HoveyBook})
travaillent dans la situation duale, i.e. o\`u les fibrations sont les
surjections; on doit donc
appliquer leurs constructions
\`a la cat\'egorie ab\'elienne duale de celle des $\Oo$-modules.  Cette
dualit\'e pourrait
\'eventuellement
poser des probl\`emes de nature ensembliste. Hovey y fait allusion dans
son livre
\cite{HoveyBook}, et attribue \`a Grodal
la construction o\`u les cofibrations sont
les injections.
Dans la version la plus r\'ecente de \cite{HoveyBook}, Hovey donne la
construction de cette cmf qu'il
appelle la ``structure injective de cmf pour les complexes'', et il prouve
qu'elle est engendr\'ee
par cofibrations
(\cite{HoveyBook} Def. 2.3.13).
Le lecteur pourra aussi retrouver ce r\'esultat en appliquant le crit\`ere
\ref{dhklemme}
\`a la cat\'egorie $Cpx_{\Oo}(X)$ avec $W = qis(X)$ et $cof$ les injections
terme-\`a-terme de complexes. On v\'erifie les propri\'et\'es (4) et (5) de
\ref{dhklemme}
en utilisant encore la technique de Jardine \cite{Jardine} (les autres
propri\'et\'es sont
faciles).


On a donc bien une cmf $Cpx_{\Oo}(X)$ engendr\'ee par cofibrations, et on
v\'erifie que le foncteur
$$
X\mapsto Cpx_{\Oo}(X)
$$
est un pr\'efaisceau de Quillen \`a gauche sur $\Xx$.

On voudrait appliquer
\ref{dansleschamps} pour obtenir que ce pr\'efaisceau est un champ.
Malheureusement, dans la pr\'esente version {\tt v3} nous ne pouvons
pas faire cela car nous ne savons pas calculer directement les
pr\'efaisceaux simpliciaux de morphismes pour les complexes non-born\'es.

Nous nous restreignons donc au cas des complexes born\'es inf\'erieurement
(on retire donc  notre l'\'enonc\'e pour les complexes non-born\'es).

Soit $Cpx_{\Oo}^{\{ 0,\infty )}(X)$ la sous-cat\'egorie des complexes
concentr\'es
cohomologiquement en degr\'es $\geq 0$.
 Par \ref{stabilite} $L(Cpx_{\Oo}^{\{ 0,\infty )}(X))\subset
L(Cpx_{\Oo}(X))$ est une sous-cat\'egorie simpliciale pleine.
La descente pr\'eserve la condition de troncation cohomologique. Un argument
\`a base de troncation de complexes montre que
$$
L(Cpx_{\Oo}^{[0,\infty )}(X))\rightarrow L(Cpx_{\Oo}^{\{ 0,\infty )}(X))
$$
est une \'equivalence. On peut munir $Cpx_{\Oo}^{[0,\infty )}(X)$ d'une
structure de cmf engendr\'ee par cofibrations (c'est plus facile que le
resultat de Hovey mentionn\'e ci-dessus pour les complexes non-born\'es).

Maintenant on applique le th\'eor\`me
\ref{dansleschamps} pour obtenir le corollaire suivant,
que le localis\'e $L(Cpx_{\Oo}^{[0,\infty )})$
est un $1$-champ de
Segal.

Il est \`a noter
que ce corollaire est
essentiellement la ``descente cohomologique'' des complexes de SGA 4
\cite{SGA4b} expos\'e Vbis,
mais avec un volet
suppl\'ementaire ``strictification des donn\'ees de descente''
(voir \S 18) qui n'appara\^{\i}t pas dans SGA 4. Pour ce corollaire,
on pourrait probablement remplacer la partie de
la d\'emonstration de
\ref{dansleschamps} qui se trouve au \S 19, par les arguments de \cite{SGA4b}.


\begin{corollaire}
\label{cestunchamp}
Supposons que $\Xx$ admet des produits fibr\'es et suffisamment de sommes
disjointes. Alors le pr\'efaisceau de cat\'egories simpliciales
$$
L(Cpx_{\Oo}^{[0,\infty )}):= 
\left( X\mapsto  L(Cpx_{\Oo}^{[0,\infty )}(X), qis(X)) \right)
$$
est un $1$-champ de Segal sur $\Xx$ qu'on appelle le {\em $1$-champ de Segal de
modules des complexes sur $(X,\Oo )$}.
\end{corollaire}
{\em Preuve:} Appliquer \ref{dansleschamps}. La compatibilit\'e (3) avec les
produits est imm\'ediate, de m\^eme que  la compatibilit\'e (2) avec les
sommes disjointes.  La condition (1) est cons\'equence du caract\`ere
local de la
notion de quasi-isomorphisme. Pour la condition (0) on applique le th\'eor\`eme
\ref{sectionscmf}
(II), en disant comme d'habitude que la v\'erification de la premi\`ere
partie de l'hypoth\`ese
(celle qui correspond aux conditions (4) et (5) de \ref{dhklemme}) fait appel
aux techniques de
\cite{Jardine}.


Pour la condition (4) (ajout\'ee dans la version 2 du papier) on va prouver
que si $f: X\rightarrow Y$ est un morphisme dans $\Xx $ alors le foncteur
de restriction $f^{\ast}$ des complexes, qui est de Quillen \`a gauche,
est aussi de Quillen \`a droite. Cette assertion implique imm\'ediatement
la compatibilit\'e aux limites voulue. Pour prouver que $f^{\ast}$
est de Quillen \`a droite on construit son adjoint \'a gauche $f_{!}$ par
analogie avec la construction de \S 4. Comme il s'agit de complexes et non
de $n$-cat\'egories, on remplace la r\'eunion disjointe par la somme directe.
Pour un complexe $C^{\cdot}$ de faisceaux de $\Oo$-modules
sur $\Xx /X$, on pose:
$$
f_{!, 0}(C^{\cdot})(Z) := \bigoplus _{\varphi : Z\rightarrow X}
C^{\cdot} (Z).
$$
Ensuite on pose $f_!(C^{\cdot})$ \'egale au faisceau associ\'e \`a 
$f_{!,0}(C^{\cdot})$. (Si les objets repr\'esentables sont des faisceaux
sur $\Xx$ cet \'etape n'est pas n\'ecessaire car $f_{!,0}(C^{\cdot})$
est d\'ej\`a un complexe de faisceaux.)

Ce foncteur $f_!$ est l'adjoint \`a gauche de $f^{\ast}$ et 
le seul probl\`eme est
de prouver qu'il pr\'eserve les cofibrations (resp. les cofibrations 
triviales).
La pr\'eservation des cofibrations est imm\'ediate car les cofibrations sont 
juste
les injections (le foncteur $f_{!,0}$ pr\'eserve les injections trivialement
et le foncteur ``faisceau associ\'e'' envoi les injections de pr\'efaisceaux
sur des injections de faisceaux). 

Pour les cofibrations triviales il suffit de
prouver que le foncteur $f_{!}$ pr\'eserve les quasi-isomorphismes. 
On observe que le foncteur
$f_{!,0}$ pr\'eserve la structure exacte des complexes de pr\'efaisceaux, 
et donc
que les faisceaux de cohomologie de $f_{!}(C^{\cdot})$ sont  les faisceaux
associ\'es aux images par $f_{!,0}$ des pr\'efaisceaux de cohomologie de
$C^{\cdot}$. En utilisant la structure exacte il suffit de prouver que si $A$
est un pr\'efaisceau sur $\Xx /X$ dont le faisceau associ\'e est $0$
alors il en est de m\^eme de $f_{!,0}(A)$. L'hypoth\`ese sur $A$ veut dire
que toute section se trivialise sur un recouvrement. Une section $a$ de
$f_{!,0}(A)$ est une somme finie de sections de $A$, donc on peut prendre 
un raffinement commun des recouvrements qui trivialisent ces composantes,
pour trivialiser $a$.  Ceci compl\`ete la d\'emonstration de la condition (4),
ce qui permet d'appliquer \ref{dansleschamps}.


On s'attaque maintenant \`a la v\'erification de la condition (5)
de \ref{dansleschamps}.
On commence par identifier les pr\'efaisceaux simpliciaux de morphismes
dans $L(Cpx_{\Oo}^{[0, \infty )})$. Soit $DP$ la
construction de Dold-Puppe (\cite{DoldPuppe} cf \cite{Illusie}) qui
transforme un
pr\'efaisceau de complexes de groupes ab\'eliens en un pr\'efaisceau
de groupes ab\'eliens
simpliciaux; et notons $\tau ^{\leq 0} C$ la troncation ``intelligente''
(i.e. celle
qui pr\'eserve la cohomologie)  en degr\'es $\leq 0$ d'un complexe $C$.
Pour
$U,V\in Ob\,Cpx_{\Oo}^{[0, \infty )}(X)$, on a
$$
L(Cpx_{\Oo}^{[0, \infty )})_{1/}(U,V)
\cong DP(\tau ^{\leq 0} {\bf R} \underline{Hom} (U,V))
$$
o\`u pour $Y\in \Xx /X$,
$$
{\bf R} \underline{Hom} (U,V)(Y):=
Hom _{\Oo}(U|_{\Xx /Y}, V'|_{\Xx /Y})
$$
d\'esigne
le {\em pr\'efaisceau $\underline{Hom}$ interne} avec $V\rightarrow V'$ un
quasi-isomorphisme vers un complexe de faisceaux injectifs.
Plus g\'en\'eralement si $\Ff$ est une cat\'egorie ab\'elienne avec suffisamment
d'injectifs, si $Cpx^{[0, \infty )}(\Ff )$ est la cat\'egorie des complexes
en degr\'es $\geq 0$ d'objets de $\Ff$, alors on obtient la cat\'egorie
simpliciale $L(Cpx^{[0, \infty )}(\Ff ), qis)$. On a la m\^eme formule
$$
L(Cpx^{[0, \infty )}(\Ff ), qis)_{1/}(U,V)
\cong DP(\tau ^{\leq 0} {\bf R} \underline{Hom} (U,V)).
$$
En particulier, $L(Cpx^{[0, \infty )}(\Ff ), qis)$ est tr\`es proche
de la {\em cat\'egorie diff\'erentielle gradu\'ee} des complexes consid\'er\'ee
par Bondal-Kapranov \cite{BondalKapranov}. D'autre part la troncation
$$
\tau _{\leq 2} L(Cpx^{[0, \infty )}(\Ff ), qis),
$$
qui peut \^etre choisie comme
une $2$-cat\'egorie stricte, est celle qui a \'et\'e consid\'er\'ee par
Gabriel-Zisman \cite{GabrielZisman}.

Revenons \`a nos complexes de faisceaux.
On va prouver directement la condition (5) de \ref{dansleschamps}
(autrement dit la partie (a) du crit\`ere \ref{critere} ) \`a savoir que
$$
L(Cpx_{\Oo}^{[0, \infty )})_{1/}(U,V)
$$
est un champ. 
Soit 
$$
I^{\cdot}:= {\bf R} \underline{Hom} (U,V)
:= \underline{Hom} (U,V'); 
$$
c'est un complexe
de faisceaux de groupes ab\'eliens acycliques. Pour tout $k\geq 0$ on
d\'efinit le tronqu\'e $U_{\leq k}$ du complexe $U$ (cette
troncation pr\'eserve
les cohomologies). Ce tronqu\'e est un sous-complexe de $U$ et $U$ est
la colimite des $U_{\leq k}$. Les morphismes dans le syst\`eme sont
des injections i.e. cofibrations de complexes, donc la colimite est une
$hocolim$. On pose 
$$
I^{\cdot}_{k} :=\underline{Hom} (U_{\leq k},V').
$$
Alors $I^{\cdot }_{k}$ est un complexe concentr\'e en degr\'es $\geq -k$
et dont les composantes sont des faisceaux de groupes ab\'eliens
acycliques. 

Il s'ensuit que 
$DP(\tau ^{\leq 0} I^{\cdot}_{ k})$ admet un d\'evissage (fini) par
produits fibr\'es homotopiques successifs
avec des pr\'efaisceaux simpliciaux
de la forme $DP(I^j[j]) = K(I^j, j)$, o\`u les $I^k$ sont acycliques;
ces derniers sont des champs
d'apr\`es \cite{Jardine} \cite{KBrown} et au d\'epart ($j=0$)
on a $K(\ker (I^0_k\rightarrow I^1_k), 0)$ qui est un faisceau
d'ensembles, donc un champ. En appliquant le lemme \ref{fiprodchamp}
on obtient
que $DP(\tau ^{\leq 0} I^{\cdot}_{k})$ est un champ.

D'autre part, $I^{\cdot} = \lim _{\leftarrow , k} I^{\cdot }_{k}$.
On obtient que
$$
DP(\tau ^{\leq 0} I^{\cdot}) = \lim _{\leftarrow , k}
DP(\tau ^{\leq 0} I^{\cdot}_k)
$$ 
et cette limite est une limite homotopique objet-par-objet au-dessus
de $\Xx$: 
les morphismes de transition
sont des surjections de groupes abeliens simpliciaux 
(sauf en degr\'e z\'ero). Autrement dit
les morphismes dans le syst\`eme inverse sont des fibrations HBKQ
de pr\'efaisceaux simpliciaux. Si l'on remplace
ce syst\`eme inverse par un syst\`eme, \'equivalent objet-par-objet,
not\'e $\{ A_k\}$, et dont les objets sont fibrants et les
morphismes de transition des
fibrations pour la structure de Jardine, alors pour tout objet
$X$ de $\Xx$ on aura
$$
DP(\tau ^{\leq 0} I^{\cdot}) (X) = \lim _{\leftarrow , k}
DP(\tau ^{\leq 0} I^{\cdot}_k)(X)
\cong \lim _{\leftarrow , k} A_k(X).
$$
D'autre part chaque $A_k$ est un champ pour la topologie $\Gg$ du site
$\Xx$, donc d'apr\`es les lemmes
\ref {grofibchampimplfib} et
\ref{stuff} les  $A_k$ sont $\Gg$-fibrants
et les morphismes de transition sont des $\Gg$-fibrations.
Donc la limite $\lim _{\leftarrow , k} A_k$ est un champ,
ce qui prouve que 
$$
L(Cpx_{\Oo}^{[0, \infty )})_{1/}(U,V)
= DP(\tau ^{\leq 0} I^{\cdot})
$$
est un champ. 
Ceci termine la v\'erification de la condition (5) du th\'eor\`eme 
\ref{dansleschamps}, qui permet de conclure que 
$L(Cpx_{\Oo}^{[0,\infty )})$ est un $1$-champ de Segal.
\eop

On conjecture bien \'evidemment que l'hypoth\`ese sur le site $\Xx$ 
n'est pas
n\'ecessaire. D'autre part, l'affaiblissement (dans {\tt v3}) de ce corollaire 
d\^u \`a l'erreur
dans le lemme \ref{pourloc} nous pousse \`a poser la question suivante:
\newline
---{\em Question:} Le $1$-pr\'echamp de Segal 
localis\'e de Dwyer-Kan des complexes non-born\'es par rapport aux
quasi-isomorphismes $L(Cpx_{\Oo})$ est-il un champ?

On passe maintenant au cas des complexes born\'es
des deux cot\'es.

Pour $0\leq a < \infty$ on note
$Cpx_{\Oo}^{[0, a\} }(X)$ (resp. $Cpx_{\Oo}^{[0, a]}(X)$)
la sous-cat\'egorie des complexes \`a cohomologie nulle
en degr\'e $>a$ (resp. \`a
composantes nulles en degr\'e $>a$). En utilisant la propri\'et\'e
\ref{stabilite} on obtient que
$$
L(Cpx_{\Oo}^{[0, a\} }(X), qis(X))\rightarrow
L(Cpx_{\Oo}^{[0, \infty )}(X), qis(X))
$$
est une sous-cat\'egorie simpliciale pleine, et avec un argument de troncation
de complexes (utilisant par exemple \ref{stabilite}) on obtient que le
morphisme
$$
L(Cpx_{\Oo}^{[0, a] }(X), qis(X))\rightarrow
L(Cpx_{\Oo}^{[0, a\} }(X), qis(X))
$$
est une \'equivalence. On obtient une \'equivalence objet-par-objet
entre pr\'efaisceaux de cat\'egories simpliciales sur $\Xx$
$$
L(Cpx_{\Oo}^{[0, a] })\cong
L(Cpx_{\Oo}^{[0, a\} }).
$$
Ce sont des sous-$1$-pr\'echamps de Segal pleins de
$L(Cpx_{\Oo}^{[0, \infty )})$, en particulier
la condition (a) de \ref{critere} est
automatique; mais comme l'annulation de la cohomologie en
degr\'es  $>a$ est
une propri\'et\'e locale, la condition (b) de \ref{critere} est v\'erifi\'ee par
$L(Cpx_{\Oo}^{[0, a\} })$. On conclut que
$L(Cpx_{\Oo}^{[0, a] })$ est un $1$-champ de Segal.

Par translation, $L(Cpx _{\Oo}^{[a,b]})$ est un $1$-champ de Segal. Pour
$n=b-a$ les champs de morphismes ici sont des pr\'efaisceaux simpliciaux
$n$-tronqu\'es. En particulier (cf les notations de \S 2),
$$
L(Cpx _{\Oo}^{[a,b]})\cong \Re _{\geq 1}
\Pi _n \circ L(Cpx _{\Oo}^{[a,b]}),
$$
donc par abus de notation on peut identifier $L(Cpx _{\Oo}^{[a,b]})$
avec le $n+1$-champ (non de Segal)
$$
\tau _{\leq n+1}L(Cpx _{\Oo}^{[a,b]})=\Pi _n \circ L(Cpx _{\Oo}^{[a,b]}).
$$
Ce $n+1$-champ est $1$-groupique (en effet il provient d'une $1$-champ de
Segal $n+1$-tronqu\'e).

\subnumero{Complexes parfaits}

Les references principaux sont SGA 6 (\cite{SGA6}), Illusie \cite{Illusie} et
Thomason-Trobaugh \cite{ThomasonTrobaugh}, mais il y a de nombreux autres
travaux sur ce sujet.

Soit $CpxParf^{[a,b]}(X)$ la cat\'egorie des complexes parfaits de
$tor$-amplitude contenue dans l'intervalle $[a,b]$, et soit $W$ la
sous-cat\'egorie des quasi-isomorphismes. D'apr\`es la th\'eorie de Dwyer-Kan
\cite{DwyerKan3} cf \S 8, on obtient des $\infty$-cat\'egories
$L( CpxParf^{[a,b]}(X))$
qui sont en fait des
cat\'egories simpliciales strictes et que nous consid\'erons comme
$1$-cat\'egories de Segal, mais qui sont d'autre part $n+1$-tronqu\'es pour
$n=b-a$ et qu'on pourrait donc consid\'erer comme des $n+1$-cat\'egories. Quand
$X$ varie ceci fournit un $n+1$-pr\'echamp ${\cal A}^{[a,b]}$
au-dessus de $\Xx$:
$$
{\cal A}^{[a,b]} (X):= \Pi _n \circ L( CpxParf^{[a,b]}(X)).
$$

Comme corollaire de \ref{cestunchamp} on a:

\begin{proposition}
\label{intro1}
Le $n+1$-pr\'echamp ${\cal  A}^{[a,b]}$ est un $n+1$-champ ($n=b-a$).
\end{proposition}
\eop

Soit $CpxFib^{[a,b]}(X)$ la cat\'egorie des complexes de fibr\'es
(i.e. de faisceaux de $\Oo$-modules localement libres de rang fini) \`a support
dans l'intervalle $[a,b]$. On d\'efinit un nouveau $n+1$-pr\'echamp
$L(CpxFib^{[a,b]})$ en posant
$$
L(CpxFib^{[a,b]})(X):= L(CpxFib^{[a,b]}(X)).
$$
Ce $n+1$-pr\'echamp vient avec un morphisme
$$
L(CpxFib^{[a,b]})\rightarrow {\cal  A}^{[a,b]}=L(CpxParf^{[a,b]}).
$$

On introduit la sous-cat\'egorie $CpxProj^{[a,b]}(X)\subset CpxFib^{[a,b]}(X)$
des complexes de $\Oo (X)$-modules projectifs. Notons qu'\`a tout
complexes de $\Oo (X)$-modules
projectifs $M$ on peut associer le complexe de faisceaux sur $\Xx/X$
$$
\tilde{M}:Y\mapsto M_Y:= M\otimes _{\Oo
(X)}\Oo (Y).
$$
Sur $CpxProj^{[a,b]}(X)$, on dispose de la
structure de cmf de Quillen, dans laquelle on a
$$
L(CpxProj^{[a,b]}(X))_{1/}(U,V)
\cong DP(\tau ^{\leq 0} {\bf R} \underline{Hom} (U,V)).
$$
Il s'agit ici du ${\bf R} \underline{Hom}$ des complexes de
$\Oo (X)$-modules. En g\'en\'eral, il diff\`ere du
${\bf R} \underline{Hom}$ des faisceaux, i.e. la fl\`eche
$$
{\bf R} \underline{Hom} (U,V) \rightarrow {\bf R} \underline{Hom}
(\tilde{U},\tilde{V}) (X)
$$
n'est pas en g\'en\'eral une \'equivalence. Ceci tient \`a la non-trivialit\'e
cohomologique de $X$ (dans le site des sch\'emas on pourrait se restreindre
aux $X$ affines et on n'aurait pas ce probl\`eme). Cependant on peut voir
(admettons qu'il faudrait fournir plus de d\'etails sur ce point, ce qui sera
fait dans un travail ult\'erieur) que le morphisme de pr\'echamps sur $\Xx /X$:
$$
\left( Y\mapsto {\bf R} \underline{Hom} _{\Oo (Y)}(U_Y,V_Y)\right)
\rightarrow {\bf R} \underline{Hom}
(\tilde{U},\tilde{V})
$$
est une $\Gg$-\'equivalence faible. Or on a
$$
L(CpxParf^{[a,b]})_{1/}(\tilde{U}, \tilde{V})=
DP \tau ^{\leq 0}{\bf R} \underline{Hom}
(\tilde{U},\tilde{V})
$$
et en particulier, ce dernier est un champ. Donc c'est le champ associ\'e au
pr\'echamp
$$
L(CpxProj^{[a,b]})_{1/}(U,V) =
\left( Y\mapsto {\bf R} \underline{Hom} _{\Oo (Y)}(U_Y,V_Y)\right) .
$$
Par le lemme \ref{formule} on conclut que le morphisme
$$
({\bf ch}L(CpxProj^{[a,b]}))_{1/}(U,V) \rightarrow
L(CpxParf^{[a,b]})_{1/}(\tilde{U}, \tilde{V})
$$
est une \'equivalence, i.e. que
$$
{\bf ch}L(CpxProj^{[a,b]})\rightarrow L(CpxParf^{[a,b]})
$$
est pleinement fid\`ele. Le morphisme
$$
L(CpxProj^{[a,b]})\rightarrow L(CpxParf^{[a,b]})
$$
est clairement $\Gg$-essentiellement surjectif, donc on obtient que ce
morphisme induit une \'equivalence
$$
{\bf ch}L(CpxProj^{[a,b]})\stackrel{\cong}{\rightarrow} L(CpxParf^{[a,b]}).
$$
Il s'ensuit
(en utilisant la suite d'inclusions $CpxProj \subset CpxFib \subset CpxParf$
et le r\'esultat de la ligne pr\'ec\'edente) que ${\bf ch}L(CpxFib^{[a,b]})$
se r\'etracte sur  ${\bf ch}L(CpxProj^{[a,b]})$.

On voit d'autre part, \`a partir de la description de $L(CpxFib^{[a,b]})(X)$ par
``hamacs'' (Dwyer-Kan \cite{DwyerKan2}, cf \S 8), que toute $i$-fl\`eche dans
cette $n+1$-cat\'egorie provient localement (par rapport \`a la topologie
$\Gg$) d'une $i$-fl\`eche de $L(CpxProj^{[a,b]})$. Ceci implique que la
r\'etraction trouv\'ee est une \'equivalence
$$
{\bf ch}L(CpxProj^{[a,b]}) \cong {\bf ch}L(CpxFib^{[a,b]}).
$$
On a finalement obtenu la proposition suivante.

\begin{proposition}
${\cal  A}^{[a,b]}=L(CpxParf^{[a,b]})$ est le $n+1$-champ associ\'e au
$n+1$-pr\'echamp $L(CpxFib^{[a,b]})$.
\end{proposition}
\eop

On indique maintenant une d\'emonstration alternative du r\'esultat d'Illusie
\cite{Illusie}
sur l'existence des puissances sym\'etriques etc. des complexes parfaits.  Cette
d\'emonstration
donne vie \'a l'intuition selon laquelle ``un complexe parfait
\'etant localement un
complexe de fibr\'es
vectoriels, il suffit de prendre la puissance sym\'etrique de chacun de ces
complexes de fibr\'es et de
les recoller''.

\begin{corollaire}
Toute construction (telle
la puissance sym\'etrique) qui provient d'un foncteur
(naturel en $X\in \Xx$)
$$
CpxFib^{[a,b]}(X) \rightarrow CpxFib^{[a',b']}(X)
$$
et qui est compatible aux \'equivalences faibles,
s'\'etend en un morphisme de $n+1$-champs
$$
L(CpxParf^{[a,b]})\rightarrow L(CpxParf^{[a',b']}).
$$
Ce morphisme se tronque par la suite en un morphisme entre les cat\'egories
d\'eriv\'ees
$$
{\cal D}^{[a,b]}_{\rm parf}(X)\rightarrow {\cal D}^{[a',b']}_{\rm parf}(X).
$$
\end{corollaire}
{\em Preuve:} Appliquer la propri\'et\'e universelle du champ associ\'e \`a un
pr\'echamp (voir \S 13) pour obtenir le morphisme entre les champs
associ\'es. En
notant l'\'egalit\'e
$$
{\cal D}^{[a,b]}_{\rm parf}(X)=\tau _{\leq 1} \left(
L(CpxParf^{[a,b]})(X)\right) ,
$$
on obtient le morphisme entre les cat\'egories d\'eriv\'ees.
\eop

On se tourne maintenant vers le cas du site $\Xx = Sch$ des sch\'emas.
On dispose alors d'une notion de {\em $n$-champ localement g\'eom\'etrique} de
$n$-groupo\"{\i}des \cite{geometricN}, qui est l'analogue pour les
$n$-champs de la
notion de {\em $1$-champ alg\'ebrique} d'Artin \cite{ArtinInventiones}. La
notion s'\'etend aux $n$-champs non n\'ecessairement de groupo\"{\i}des. Sans
entrer dans les d\'etails de cette d\'efinition, on peut dire que si $A$ est
un $n$-champ (localement) g\'eom\'etrique alors son int\'erieur $A^{int,
0}$ est
un $n$-champ de groupo\"{\i}des (localement) g\'eom\'etrique au sens de
\cite{geometricN}. On \'enonce ici un r\'esultat qui fera l'objet d'un autre
travail.

\begin{theoreme}
\label{modparfait}
Sur le site des sch\'emas,
le $n+1$-champ ${\cal  A}^{[a,b]}$ est localement g\'eom\'etrique ($n=b-a$).
Il est couvert par les sous-champs ouverts ${\cal  A}^{s_a,\ldots ,s_b}$
obtenus en imposant \`a la cohomologie en tout point la majoration $h^i \leq
s_i$;
ces sous-champs sont g\'eom\'etriques.
\end{theoreme}

En fait, le lecteur pourra peut-\^etre trouver la d\'emonstration de ce
th\'eor\`eme en s'\-ap\-puy\-ant sur la ``th\'eorie homologique des
perturbations'' de Gugenheim {\em et al.} \cite{Gugenheim}
\cite{GugenheimLambe}. Il s'agit pour l'essentiel de prouver que le morphisme du
sch\'ema \'evident de param\`etres pour un complexe de fibr\'es triviaux (de
rangs donn\'es) vers le champ ${\cal  A}^{[a,b]}$ est lisse.

Le $1$-tronqu\'e $\tau _{\leq 1}{\cal  A}^{[a,b]}$ (dont les valeurs
sont les $1$-tronqu\'ees des $n+1$-cat\'egories ${\cal  A}^{[a,b]}(X)$
est le pr\'efaisceau de cat\'egories qui \`a $X$ associe la {\em cat\'egorie
d\'eriv\'ee} des complexes parfaits d'amplitude contenue dans $[a,b]$. Ni
le th\'eor\`eme \ref{intro1} ni
le th\'eor\`eme \ref{modparfait} ne restent vrais pour ce
tronqu\'e: ceci montre l'inter\^et de ``consid\'erer les homotopies
sup\'erieures''.

On peut utiliser les $L(Cpx _{\Oo})$ et ${\cal A}^{[a,b]}$ pour d\'efinir
les notions de complexe,
complexe parfait etc. sur des $n$-champs de Segal (en particulier, sur des
$1$-champs e.g. des
$1$-champs alg\'ebriques d'Artin). Si $B$ est un $n$-champ de Segal, on
dira qu'un {\em complexe de
$\Oo$-modules sur $B$} est un morphisme (vers le remplacement fibrant)
$$
B\rightarrow L(Cpx
_{\Oo})'.
$$
Un {\em complexe parfait d'amplitude contenue dans $[a,b]$} est un morphisme
$$
B\rightarrow ({\cal A}^{[a,b]})'.
$$
Les r\'esultats ci-dessus montrent que si $B$ est le $0$-champ repr\'esent\'e
par un objet du site
(e.g. un sch\'ema),
alors ces d\'efinitions co\"{\i}ncident avec les d\'efinitions habituelles.

\end{document}